\documentclass[reqno,11pt]{amsart}
\usepackage[utf8]{inputenc}
\usepackage{pifont}
\usepackage[T1]{fontenc}
\usepackage{textcomp}

\usepackage{pgffor}

\usepackage{stackengine}
\usepackage{esint}
\usepackage{pdfsync}
\usepackage[colorlinks=true,linkcolor=blue]{hyperref}%
\usepackage[T1]{fontenc}
\usepackage{pifont}
\usepackage[english]{babel}

\usepackage{ulem}
\usepackage{amsmath,esint,color}
\usepackage{amsthm,latexsym,epsfig,graphicx,marvosym, mathrsfs,bigints,stmaryrd,textgreek,upgreek,rsfso}
\usepackage{amsmath,stmaryrd}
\usepackage{relsize}
\usepackage{amsfonts}
\usepackage{amssymb}
\usepackage{fullpage}

\allowdisplaybreaks

\def\eps
{\varepsilon}

\def\N{{\mathbb N}}
\newcommand{\DDD}{\Delta\!\!\!\!\Delta}

\addcontentsline{toc}{section}{Preface}

\newcommand{\CC}{\mathbb{C}}
\newcommand{\NN}{\mathbb{N}}
\newcommand{\Z}{\mathbb{Z}}
\newcommand{\RR}{\mathbb{R}}
\newcommand{\T}{\mathbb{T}}

\newcommand{\C}{\mathbb{C}}
\newcommand{\tc}{\mathtt c}

\newcommand{\ii }{{\rm i} }

\newtheorem{definition}{Definition}[section]

\newtheorem{theorem}{Theorem}[section]
\newtheorem{proposition}{Proposition}[section]
\newtheorem{lemma}{Lemma}[section]
\newtheorem{remark}{Remark}[section]
\newtheorem{coro}{Corollary}[section]
\newtheorem{remarks}{Remarks}[section]
\numberwithin{equation}{section}

 \title{KAM theory for active scalar equations}
 \numberwithin{equation}{section}

 \author{Zineb Hassainia}
 
  \address{ NYUAD Research Institute, New York University Abu Dhabi, PO Box 129188, Abu Dhabi,   United Arab Emirates. }
\email{zh14@nyu.edu.}
 
  \author{ Taoufik Hmidi}

    \address{ NYUAD Research Institute, New York University Abu Dhabi, PO Box 129188, Abu Dhabi,   United Arab Emirates. 
    Univ Rennes, CNRS, IRMAR – UMR 6625, F-35000 Rennes, France}
\email{thmidi@univ-rennes1.fr.}
  
   \author{ Nader Masmoudi}   
   
   \address{ NYUAD Research Institute, New York University Abu Dhabi, PO Box 129188, Abu Dhabi,   United Arab Emirates.
 Courant Institute of Mathematical Sciences, New York University, 251 Mercer Street, New York, NY 10012, USA.}
\email{masmoudi@cims.nyu.edu.}

\begin{document}
\begin{abstract}
 In this paper, we establish  the existence of time quasi-periodic  solutions to generalized surface quasi-geostrophic equation $(\textnormal{gSQG})_\alpha$  in the patch form close to  Rankine vortices. We show that invariant tori survive when the order $\alpha$ of the singular operator  belongs to a Cantor set contained in $(0,\frac12)$ with almost full Lebesgue measure. The proof is based on several techniques from  KAM theory, pseudo-differential calculus  together with  Nash-Moser scheme  in the spirit of the recent works  \cite{Baldi-berti,BB13}. One key novelty here is a refined  Egorov type theorem  established through a new approach based on the kernel dynamics together with some hidden  T\"opliz structures.
\end{abstract}
\maketitle
\tableofcontents

\section{Introduction}

\quad In this paper we  are concerned  with the   generalized   surface quasi-geostrophic  equations {$($gSQG$)_\alpha$} described by 

\begin{equation}\label{Eq-F}
\left\{ \begin{array}{ll}
\partial_{t}\theta+u\cdot\nabla\theta=0,\quad(t,x)\in[0,T]\times\RR^2, &\\
u=-\nabla^\perp(-\Delta)^{-1+\frac{\alpha}{2}}\theta,\\
\theta_{|t=0}=\theta_0.
\end{array} \right.
\end{equation}
Here $u$ stands for  the velocity field, $\nabla^\perp=(-\partial_2,\partial_1)$  and $\alpha\in(0,2)$. 
The fractional Laplacian   $(-\Delta)^{-1+\frac{\alpha}{2}}$ is associated to Riesz potential according to the convolution law 
\begin{equation}\label{Integ1}
(-\Delta)^{-1+\frac{\alpha}{2}} \theta(x)=\frac{C_\alpha}{2\pi}\mathlarger{\int}_{\RR^2}\frac{\theta(y)}{\vert x-y\vert^\alpha}dy,\quad{\rm with}\quad C_\alpha\triangleq\frac{\Gamma(\alpha/2)}{2^{1-\alpha}\Gamma(\frac{2-\alpha}{2})}\cdot
\end{equation}

The limiting case $\alpha\to0$ corresponds to the 2D Euler equation, with $\theta$ representing the fluid vorticity,   the midpoint case $\alpha=1$ corresponds to the surface quasi-geostrophic equation (SQG), with $\theta$ denoting the temperature  in a rapidly rotating stratified fluid with uniform potential vorticity and  the case $\alpha=2$ produces stationary solutions. 

\smallskip

\quad While global regularity of smooth solutions to the 2D Euler equations has been known  since long ago \cite{Holder,Wolibner}, only local regularity persistence  is proved for the entire range  $\alpha\in(0,2)$, we refer for instance to  \cite{C-C-C-G-W,C-C-W,C-M-T}.
   The question of whether a finite-time singularity may develop from smooth initial datum remains open. However, certain blow-up scenarios have been ruled out in  \cite{Cord,C-F} and solutions with exponential growth (not excluding blow up in finite time) are shown to exist in \cite{H-K}. On the other hand,  
 examples of a non-trivial  global in time smooth solutions have been recently constructed  in \cite{ADDMW,CCG2,Grav-Smets} with different approaches related to bifurcation theory/assisted computer proof and variational principle.
 
 In the setting of weak solutions, it is known that  for Euler equations Yudovich solutions   exist  globally  in time and they are unique, see   \cite{Y1}. However, for $\alpha>0$  global existence   persists  but the uniqueness issue at the energy level  remains open, see for instance  \cite{Marchand,Resnick}. Recently,  non-unique weak solutions with negative Sobolev regularity were constructed in \cite{BSV}. An interesting sub-class which has been explored in various directions  during the past few decades is given by the vortex patches, that is, initial data taking the form of the characteristic function $\theta_0(x)=\chi_D$,  where $D$ is a bounded smooth domain.  In this framework, the contour dynamics  equation offers a suitable tool to construct  local solutions in the patch form with smooth boundary, see   \cite{C-C-C-G-W,Gancedo,R}. The uniqueness in this class  for $\alpha\in (0,1)$, was discussed in  \cite{ KYZ} and   the case $\alpha=1$ was analyzed in  \cite{Cord-Cord-Gan}.

  It is worthy to point out that while the boundary's regularity is globally preserved for $\alpha=0$ as proved in \cite{B-C,Ch}, numerical evidence \cite{C-F-M-R} suggests singularity formation in finite time for $\alpha>0$. In \cite{KRYZ} a new scenario with multiple patches with opposite signs colliding in finite time is established for a modified quasi-geostrophic equation in half plane. A finite-time singularity criterion have been recently provided in \cite{GP}. 
  
\subsection{Periodic vortex patches solutions}
\quad Due to the complexity of the motion and the deformation process that the vorticity undergoes, the dynamics of the boundary is hard to track and very little is known concerning the  nontrivial global vortex patches  solutions for the $($gSQG$)_\alpha$ equations.
Therefore, It is of important interest to look for the emergence of ordered structures for this Hamiltonian equation.  In our context ordered structures refers to relative equilibria corresponding to steady solutions relative to some translating or rotating frame of reference.

\smallskip

\quad Historically, relative equilibria was initiated  long time ago for  the two-dimensional  Euler equations with  the work Kirchhoff \cite{Kirc} who discovered that an elliptical patch rotates uniformly about it center of mass with constant angular velocity. Variational characterization as critical points of some functional energy was formulated by Kelvin. More implicit examples of uniformly rotating  vortex patches with higher symmetry,  called V-states, were numerically computed by Deem and Zabusky in \cite{DZ}.  Later, an  analytical   proof  was given by Burbea in \cite{B} and  based on the conformal mapping parametrization combined with some tools from  local bifurcation theory.
 Burbea's branches of solutions were extended to global ones in \cite{Hass-Mass-Wheel}. The  regularity and the convexity of the V-states  have been investigated in \cite{CCG4,Hass-Mass-Wheel,HMV}. 
 
 Similar research has been recently  carried out for the $($gSQG$)_\alpha$  equations.  The construction of simply connected V-states was established in \cite{HH} for $\alpha\in(0,1)$ 
 using Burbea's approach. In particular, the bifurcation  from the unit disc occurs at the angular velocities,
\begin{align}\label{Fre-per1}
\widetilde{\Omega}_m^\alpha\triangleq\frac{\Gamma(1-\alpha)}{2^{1-\alpha}\Gamma^2(1-\frac\alpha2)}\Big(\frac{\Gamma(1+\frac\alpha2)}{\Gamma(2-\frac\alpha2)}-\frac{\Gamma(m+\frac\alpha2)}{\Gamma(m+1-\frac\alpha2)}\Big),\quad m\in\N\setminus \{0,1\}.
\end{align}
The remaining case $\alpha\in[1,2)$  was  completely solved  by Castro, C\'ordoba and  G\'omez-Serrano in \cite{CCG} 
and the boundary regularity was discussed in \cite{CCG,CCG4}.

\quad We point out that a  rich investigation has been conducted  in the past few years  around  different  topological structures for the V-states and for various nonlinear transport equations. For instance, we mention  the existence results of rotating  
 multiply-connected patches \cite{DHH,HFMV,Gomez-Serrano,HM3,HR,Coralie}. It is shown in particular that rotating patches bifurcate from the annulus and two different  branches are assigned to  the same symmetry.   In addition, formation of small loops are discovered  when the two branches are close enough. Boundary effects on the emergence of the V-states  were analyzed through the disc  example in  \cite{DHHM}, with important numerical experiments putting in evidence the oscillation of the Burbea's curves. Note also that for Euler equations a second bifurcation of countable branches from the ellipses occurs but the shapes have in fact less symmetry and being at most two-folds.  The proof of the existence and analyticity of the boundary has been  investigated in \cite{CCG,HM2}.
 
  Bifurcation of V-states  for the  quasi-geostrophic shallow-water equation was accomplished in \cite{DHR}  supplemented with numerical experiments describing the imperfect bifurcation with respect to Rossby deformation length. Another delicate subject related to the construction of non homogeneous periodic solutions  around specific radial solutions has been developed in \cite{CCG2,GHS}. This offers an excellent starting point to explore  whether periodic solutions can be captured around generic radial profiles. The rigidity of radial symmetry properties  with respect to the angular velocity was explored  in a series of works \cite{Fraenkel,gomez2019symmetry,Hm}. Another connected subject is the desingularization of the point vortex system  to equilibria of vortex patches/smooth profiles type. This   was first  studied by Turkington  \cite{T}  for Euler equation using variational arguments. Following the same approach Turkington's result was extended for the $($gSQG$)_\alpha$ equations with $\alpha\in(0,1)$ in \cite{GGS}. We notice that this approach does not seem to be  efficient to describe neither the topological structure of the patches (for example whether they are connected or not) nor the regularity of the boundary.   In \cite{HM},  the second author and Mateu gave a direct proof showing  for $\alpha\in[0,1)$ the existence of co-rotating and counter-rotating pairs of simply connected smooth patches, using a desingularization of the contour dynamics equations and an application of the implicit function theorem. This approach sounds to be  flexible and robust and has been adapted  recently by different authors  to cover various interesting  point vortex  configurations associated to  multiple models. For instance, it was used by \cite{CQZZ} to extend the construction for  the case $\alpha\in[1,2)$ and in   \cite{HH2} for the desingularization of  the asymmetric pairs. 
  The same technique was used   to  desingularize a spatial periodic distribution called  
 Karman vortex street \cite{G-Kar} and a similar study was performed for  the  Thomson polygon \cite{Garcia}. Very recently, a generic system of rotating point vortices was analyzed in \cite{HW} and a weak condition on the validity of the desingularisation through the contour dynamics has been found.

\subsection{Quasi-periodic solutions}
It is worthy to point out that $($gSQG$)_\alpha$ is a reversible Hamiltonian system and as we shall see in Section \ref{Sec-Hamiltonian structure} this structure  persists at the level of the contour dynamics equation in the setting of the vortex patches. This system depends on  one degree of freedom given by the external   parameter $\alpha$ related to  the order of the nonlocal operator.  Then it is legitimate  to explore  whether quasi-periodic solutions constructed for the  linearized operator at the equilibrium state (obtained by  linear superposition of its eigenfunctions) could survive under small perturbation  for the nonlinear model when $\alpha$ is selected in a suitable Cantor set.\\  To fix the terminology,  a real-valued  function $f:\mathbb{R}\rightarrow\mathbb{R}$ is called quasi-periodic (shortened in QP) if there exists a multi-variable function $F:\mathbb{T}^{d}\rightarrow\mathbb{R}$, with $d\geqslant 1$ such that  
 $$
 \forall\, t\in\mathbb{R},\quad f(t)=F(\omega t)
 $$
 for some frequency vector $\omega\in\mathbb{R}^{d}$  which is non-resonant in the sense 
\begin{equation}\label{non-res}
\forall\, l\in\mathbb{Z}^{d}\backslash\{ 0\},\quad\omega\cdot l\neq 0,
\end{equation}
where we denote by $\mathbb{T}^{d}$ the flat torus of dimension $d$.  In  the case $d=1$, we recover   the definition of  periodic functions with frequency $\omega\in\mathbb{R}^{*}.$ The  persistence of invariant torus is a   relevant subject for nearly integrable  Hamiltonian systems in finite or infinite dimension spaces. In general, the construction stems from   KAM theory whose main query is to device a scheme allowing to avoid resonances that may destroy the invariant  torus (which may occur  even at the linear level). This is  an active area which  has been improved  and enriched through different important studies along the past few decades.

\smallskip

\quad  The main task   in this paper is to explore the existence of quasi-periodic solutions for \eqref{Eq-F}  in the setting of vortex patches. Even though several equations were subject to KAM studies,   it seems that the  vortex motion  which is an old topic in fluid dynamics   has escaped to these studies except for the periodic framework  where a lot of results have been obtained as we have  mentioned before. One of the main advantage  in getting QP solutions is the construction of non trivial global in time solutions around the stationary Rankine vortices for the singular \mbox{model \eqref{Eq-F}}.

\smallskip

\quad The literature related to KAM theory  is very abundant and  substantial progress connected  to various geometric and analytical aspects has been accomplished during the past  decades.  Here, we shall only focus on some specific contributions fitting with the main scope of the paper. This theory was initiated by Kolmogorov \cite{Kolmogorov}, Arnold \cite{Arnold} and Moser  \cite{moser}  who  proved in finite dimensional  space the  persistence of invariant tori for  small perturbation of integrable hamiltonian systems  under suitable non degeneracy and smoothness conditions.  Later, more development around lower dimensional elliptic/hyperbolic  invariant tori was carried out by R\"ussmann \cite{Russ}. We may also refer to Sevruyk \cite{sevryuk} who {constructed invariant tori for reversible systems}. The extension of KAM theory  to infinite dimensional Hamiltonian  systems with applications to PDE's like  the wave,  Schr\"odinger and Klein Gordon  equations  was implemented  by Kuksin \cite{Kuksin1}, Wayne \cite{wayne}, Bourgain \cite{Bourgain},  Kuksin-P\"oschel \cite{Kuksin-Poschel},  Eliasson-Kuksin \cite{Ell-Kuk }, Gr\'ebert-Kappeler \cite{Grebert-Kappeler},   Gr\'ebert-Paturel \cite{Grebert-Paturel}, Baldi-Berti-Haus-Montalto~\cite{Baldi-berti}. More recent studies can be found in  \cite{BBM-auto,bambusi-Berti,BB13,BertiMontalto,Berti-Bolle2013}.

\smallskip

\indent The complexity of the problem depends on the space dimension and on the structure of the equations. For example in the semi-linear case the nonlinearity can be seen as a bounded perturbation of the linear problem and this simplifies  a lot the problem of finding a right approximate inverse of the linearized operator around a state close to the equilibrium. However in the quasi-linear case where the nonlinearity is unbounded and has the same order as the linear part the situation turns to be much more tricky. This is the case for instance in  the water-waves equations where several results has been obtained in the past few years on the periodic and quasi-periodic settings \cite{Alazard-baldi,Baldi-berti,BFM21,BFM,BertiMontalto,Ioos}.  For similar discussion we refer to  \cite{Baldi-Montalto21} for the 3D incompressible Euler equations with quasi-periodic forcing and to \cite{BHM}  for  the vortex patches solutions to the 2D Euler equations. 

\smallskip

\quad Next we shall make general comments on  the  scheme commonly used to construct quasi-periodic  solutions for semi-linear or quasi-linear PDE's that was developed in particular  in the papers \cite{Baldi-berti,BB13,BertiMontalto} which is  robust and flexible and will  be adapted in our framework up to some important technical problems that have been settled differently. The first step is to write down   the dynamics according to  the action-angle coordinates  used to describe  the  tangential part  without affecting  the normal part which  lives in an infinite dimensional space. Notice that the tangential part is composed by a finite number of excited frequencies satisfying non-resonance condition that constitute the linear non-resonant torus.  In linearizing around a state near the equilibrium,  one finds an operator with variable coefficients that  should be inverted approximately up to different types of small errors, provided that  the external parameters belong to a suitable Cantor set defined through various  Diophantine conditions.
To do that it is, first, convenient  to look for an approximate inverse using  an intermediate  isotropic torus  built around the  initial one. It has the advantage to transform the linearized operator via symplectic change of coordinates into a  triangular system up to an error vanishing when  testing against an invariant torus. This perturbation can be later incorporated in Nash-Moser scheme.  Then the outcome is that  the Hamiltonian has a good normal form structure such that one can almost decouple the dynamics in the phase space in tangential and  normal modes. On the tangential part the system can be solved in a triangular way provided we can invert the linearized operator on the normal part up to a small coupling error term, and this is more or less finite dimensional KAM theory. Then, the analysis reduces to invert the linearized operator on the normal part, which can be viewed  as a small bounded/unbounded perturbation of a diagonal operator. The main ingredient to do that  is to conjugate in a suitable way the  linearized operator into  a diagonal one with constant coefficients that one could invert. This is the major step which turns out to be highly technical due to the resonances which emerge at all the scales and one should use the degree of freedom of the system  in order to avoid it during  the KAM scheme by making consecutive excisions ending with a Cantor like set. This allows  to build an approximate right  inverse to the  linearized operator with nice tame estimates that we can combine with  a frequency cut-off of  the Nash-Moser scheme. 

\subsection{Main result and ideas of the proof}
As we shall discuss later in Section \ref{sec contour}, the contour dynamics equation describing the  vortex patch motion  can be written in a more tractable way using polar coordinates. This description is meaningful at least for a  short time when the initial patch is sufficiently close to the equilibrium state given by Rankine vortex $\mathbf{1}_{\mathbb{D}},$ where $\mathbb{D}$ is the unit disc of the plane.  Thus the boundary $\partial D_t$ can  be parametrized as follows
\begin{equation}\label{BSh}
w(t,\theta)=R(t,\theta)e^{\ii \theta}\quad\mbox{ with }\quad R(t,\theta)=\left(1+2r(t,\theta)\right)^{\frac{1}{2}}.
\end{equation}
As we shall prove in \eqref{eq},  the function  $r$ satisfies the following  nonlinear transport equation   
\begin{equation}\label{equation introduction}
\partial_{t}r+F_{\alpha}[r]=0,
\end{equation}
with
\begin{align}\label{eqZ1}
F_\alpha[r](t,\theta)={\frac{ C_\alpha }{2\pi }\bigintsss_{0}^{2\pi} \frac{\partial^2_{\theta\eta} \big[\big(1+ 2r(t,\theta)\big)^{\frac12}\big(1+ 2r(t,\eta)\big)^{\frac12}\sin(\eta-\theta)\big]}{A_r^{\alpha/2}(\theta,\eta)} d\eta}
\end{align}
and
 \begin{equation*}
A_{r}(\theta,\eta)=\left|R(t,\theta)e^{\ii \theta}-R(t,\eta)e^{\ii \eta}\right|.
\end{equation*}
Next, we take a parameter $\Omega\neq0$ and look for the solutions in the form
\begin{equation}\label{ansatz1}
r(t, \theta)=\tilde r(t,\theta+\Omega t),
\end{equation} 
then the equation \eqref{equation introduction} is equivalent to 
\begin{equation}\label{equation-Master}
\partial_{t}\tilde r+\Omega\partial_{\theta} \tilde r
 +F_{\alpha}[\tilde r]=0.
\end{equation}
Notice that the introduction of  the parameter $\Omega$ sounds at this level  artificial  but it will be used later  to cancel the trivial degeneracy of the first eigenvalue associated  with  the linearized operator at the equilibrium state.
 In the quasi-periodic setting, we should find a non-resonant  vector \mbox{frequency  $\omega\in\mathbb{R}^{d}$}  such that the equation \eqref{equation-Master} admits a solution  in the form  $\tilde{r}(t,\theta)=\widehat{r}(\omega t,\theta)$ with $\widehat{r}:\mathbb{T}^d\times\mathbb{T}\to\mathbb{R}$ being a smooth $(2\pi)^{d+1}$-periodic function. Then we can easily check that  $\widehat{r}$ , still denoted in what follows by $r$, satisfies
\begin{equation*}
\omega\cdot\partial_{\varphi}r+\Omega\partial_{\theta}r+F_{\alpha}[r]=0.
\end{equation*}
The computation of  the linearized operator of \eqref{equation-Master} at a given state $r$ is described  in Proposition~ \ref{lin-eq-r}  and one gets
\begin{align*}
\partial_t h (t, \theta) & =  \partial_\theta\Big[-\big(\Omega+V_{r,\alpha}(t,\theta)\big)h(t, \theta) +\mathbb{K}_{r,\alpha}h(t,\theta)\Big]
\end{align*}
where $V_{r,\alpha}( t,\theta)  $ is the real  function
\begin{equation*}
V_{r,\alpha}(t,\theta)\triangleq  \frac{ C(\alpha) }{2\pi }\bigintsss_{0}^{2\pi} \frac{\partial_{\eta} \big[\big(1+ 2r(t,\eta)\big)^{\frac12}\sin(\eta-\theta)\big]}{\big(1+ 2r(t,\theta)\big)^{\frac12}A_r^{\alpha/2}(\theta,\eta)} d\eta
\end{equation*}
and the integral operator $\mathbb{K}_{r,\alpha}$ is defined by
\begin{equation*}
\mathbb{K}_{r,\alpha}h(t,\theta)\triangleq  \frac{ C(\alpha) }{2\pi }\bigintsss_{0}^{2\pi} \frac{h(t,\eta)}{A_r^{\alpha/2}(\theta,\eta)} d\eta.
\end{equation*}

At the equilibrium state $r\equiv 0$, we infer from  Proposition  \ref{linear-eq}  that  the linearized operator is a Fourier multiplier, since\begin{align}\label{FME}
h(\theta)=\sum_{j\in\mathbb{Z}\setminus \{0\}}h_{j}e^{\ii j\theta},\quad\partial_\theta\Big[-\big(\Omega+V_{0,\alpha}\big) +\mathbb{K}_{0,\alpha}\Big]h(\theta)=-\sum_{j\in\mathbb{Z}\setminus \{0\}}{\Omega_{j}(\alpha)}h_{j}e^{\ii j\theta}
\end{align}
where the frequency $\Omega_{j}$ is defined by
\begin{align}\label{Freq-equilibrium} 
\Omega_j(\alpha)&\triangleq j\left(\Omega+\tfrac{\Gamma(\frac12-\frac\alpha2)}{2\sqrt{\pi}\,\Gamma(1-\frac{\alpha}{2})}\bigg(\tfrac{  \Gamma(1+\frac\alpha2)}{ \Gamma(2-\frac\alpha2) }-\tfrac{  \Gamma(j+\frac\alpha2)}{ \Gamma(1+j-\frac\alpha2)} \bigg)\right)\\
\nonumber&
= j\left(\Omega+\widetilde{\Omega}_j^\alpha\right)\cdot
\end{align}
Here,  $\Gamma$ stands for  the usual Gamma function and $\widetilde{\Omega}_j^\alpha$ was introduced before in \eqref{Fre-per1}.
Consequently, the elements with zero average of the kernel of the linearized operator at the equilibrium state   are given by
$$
h(t,\theta)=\sum_{j\in\mathbb{Z}\setminus \{0\}}h_{j}e^{\ii(j\theta-{\Omega_{j}(\alpha)} t)},
$$
supplemented with a suitable decay property  of the Fourier coefficients to rend the sum meaningful.
Since the Hamiltonian system is reversible then we can use this property to filter trivial resonances coming from opposite frequencies by looking only for  the real solutions which are invariant by involution, that is,  $h(-t,-\theta)= h(t,\theta)$). Therefore we are led to consider only solutions  in the form
$$
h(t,\theta)=\sum_{j\in\mathbb{N}\setminus \{0\}}h_{j}\cos\big(j\theta-{\Omega_{j}(\alpha)} t\big).
$$
Then, by keeping only a finite number of  frequencies, this sum gives rise to  quasi-periodic solutions with non-resonant frequency provided that $\alpha$ belongs to a suitable set defined with Diophantine condition, for a precise statement  see Lemma \ref{lemma-Line-mes}. Our main result concerns the persistence of quasi-periodic solutions for the nonlinear model \eqref{equation-Master} when the perturbation is small enough. 

\begin{theorem}\label{main theorem}
Let $\Omega>0,\,0<\underline\alpha<\overline\alpha<\frac12$,  $\mathbb{S}\subset{\mathbb{N}\backslash\{0\}}$, with $\#\mathbb{S}=d\geqslant1.$
There exists  $\varepsilon_{0}\in(0,1)$ small enough  with the following properties:  For every amplitude ${\vec{\varepsilon}}=(\varepsilon_{j})_{j\in\mathbb{S}}\in(\mathbb{R}_{+}^{*})^{d}$ satisfying
$$
|\vec{\varepsilon}|\leqslant\varepsilon_{0}
$$ 
 there exists a Cantor-like set $\mathtt{C}_{\infty}\subset(\underline\alpha,\overline\alpha)$ with asymptotically full Lebesgue  measure as $\vec{\varepsilon}\rightarrow 0,$ i.e.,
$$\lim_{\vec{\varepsilon}\rightarrow 0}|\mathtt{C}_{\infty}|=\overline\alpha-\underline\alpha$$
such that for any $\alpha\in\mathtt{C}_{\infty}$, the  equation  \eqref{equation-Master} admits a time quasi-periodic solution 
 with diophantine frequency vector ${\omega}_{{\textnormal{pe}}}(\alpha,\vec{\varepsilon})\triangleq (\omega_{j}(\alpha,\vec{\varepsilon}))_{j\in\mathbb{S}}\in\mathbb{R}^{d}$ and taking  the form
$$\tilde r(t,\theta)=\sum_{j\in\mathbb{S}}\varepsilon_{j}\cos\big(j\theta+\omega_{j}(\alpha,\vec{\varepsilon})t\big)+\mathtt{p}\big({\omega}_{{\textnormal{pe}}}t,\theta\big)
$$
with
$$
{\mathtt{\omega}}_{{\textnormal{pe}}}(\alpha,\vec{\varepsilon})\underset{\vec{\varepsilon}\rightarrow 0}{\longrightarrow}(-{\Omega}_j(\alpha))_{j\in\mathbb{S}}
$$
and  ${\Omega}_j(\alpha)$ is  defined in \eqref{Freq-equilibrium} .  In addition, the perturbation $\mathtt{p}$ satisfies
$$\| \mathtt{p}\|_{H_{\textnormal{even}}^{{s}}(\mathbb{T}^{d+1},\mathbb{R})}\underset{\vec{\varepsilon}\rightarrow 0}{=}o(|\vec{\varepsilon}|),
$$
for $s$ large enough, where  the    Sobolev spaces $H_{\textnormal{even}}^{{s}}(\mathbb{T}^{d+1},\mathbb{R})$  are  defined in \eqref{FS-even}.
\end{theorem} 
Some remarks are in order.
\begin{remarks}
\begin{enumerate}
\item From this theorem,  \eqref{ansatz1} and \eqref{BSh} one gets that the boundary shape of the quasi-periodic patch can be parametrized  in polar coordinates as follows
\begin{equation*}
w(t,\theta)=R(t,\theta)e^{\ii \theta}\quad\mbox{ with }\quad R(t,\theta)=\left(1+2\tilde r(t,\theta+\Omega t)\right)^{\frac{1}{2}}
\end{equation*}
and $\tilde r$ is expanded as in the theorem. The time evolution of the shape is given  by small pulsation around the unit disc and the boundary is  localized in an annulus around the unit circle. 
\item We obtain global existence of non trivial solutions in the patch form for \eqref{Eq-F} where only a few results on global existence are known  in the periodic case, see \cite{ADDMW,CCG,CCG4,CCG2, DHHM,HH} and the references therein.   
\item Some technical problems are behind the limitation of $\overline\alpha>0$ and we cannot go to $\underline\alpha=0$. This is connected with the fact that in the fractional Laplacian \eqref{Integ1} associated to Riesz potential,  the constant $C(\alpha)$ blows up when $\alpha$ goes to zero. We know that this operator converges in a weak sense to $(-\Delta)^{-1}$ whose kernel is of logarithmic type and this change of behavior sounds to be the main obstruction in our result.   The second limitation concerns  the upper bound $\overline\alpha$ which should be smaller than $\tfrac12$. This is related first to the reducibility scheme for the fractional Laplacian part. Second, this limitation  is needed  in the proof of the refined Egorov theorem in Section \ref{Flows and Egorov theorem type}. We believe that the method developed here could handle  the remaining   case $\overline\alpha\in[\tfrac12,1)$ but the load  is very high due to several nonlocal commutators with positive order that should be reduced to constant coefficients one by one. In addition an adaptation of Egorov type theorem is also required
  \end{enumerate}
\end{remarks}
We shall now sketch the  main steps of the proof which will be implemented  following standard KAM scheme in the spirit  of the  preceding works \cite{BertiMontalto,bambusi-Berti} but with different substantial  variations as we shall discuss below. We basically use techniques from  KAM theory, pseudo-differential operators  combined with Nash Moser scheme. This will be done in   several steps which are detailed below.
\vspace{0.2cm}

\ding{202} {{\it Action-angle  reformulation.} We first notice that the equation \eqref{equation introduction} enjoys a Hamiltonian structure and from Proposition \ref{prop-hamilt}   we may write the contour dynamics equation in the form
\begin{equation}\label{HDCD}
\partial_{t}r=\partial_{\theta}\nabla H(r)
\end{equation}
and  the Hamiltonian $H$ can  expressed in terms of   the kinetic energy of the system and the angular momentum. Close to Rankine vortices, we can  write the PDE as a Hamiltonian perturbation of an  integrable system given by the linear dynamics at the equilibrium state. Indeed, according to  the expansion  \eqref{FME}  the linearized operator at  the equilibrium is given by a Fourier multiplier which gives rise to an integrable system and therefore  we can write \eqref{HDCD} in the form
$$
\partial_{t}r=\partial_{\theta}\mathrm{L}(\alpha)(r)+X_{P}(r)
$$
where $\mathrm{L}(\alpha)$ and  the perturbed Hamiltonian vector field  $X_{P}$ are defined by
$$
\mathrm{L}(\alpha)(r)= -\big(\Omega+V_{0,\alpha}\big) r +\mathbb{K}_{0,\alpha}r\quad\hbox{and}\quad X_{P}(r)=V_{0,\alpha}\partial_{\theta}r-\partial_{\theta}\mathbb{K}_{0,\alpha}r-F_{\alpha}[r].
$$
We point out that $V_{0,\alpha}$ is a constant and the operator $\mathbb{K}_{0,\alpha}$ is a Fourier multiplier.  Then we find it convenient to  rescale the solution size in the following way  $r\leadsto \varepsilon r$, with $\varepsilon$ a small parameter.  Thus we derive a new   equation that appears as a perturbation of the linearized equation at the equilibrium state.  Actually, we find that the new equation takes the form 
\begin{equation*}
\partial_{t}r=\partial_{\theta}\mathrm{L}(\alpha)(r)+\varepsilon X_{P_{\varepsilon}}(r)
\end{equation*}
where $X_{P_{\varepsilon}}$ is the Hamiltonian vector field defined by
$X_{P_{\varepsilon}}(r)\triangleq \varepsilon^{-2}X_{P}(\varepsilon r).$ Then, time quasi-periodic solutions oscillating at the  frequency $\omega\in \mathbb{R}^d$ are simply periodic solutions to the following equation
\begin{equation*}
\omega\cdot\partial_\varphi r=\partial_{\theta}\mathrm{L}(\alpha)(r)+\varepsilon X_{P_{\varepsilon}}(r).
\end{equation*}
Here we keep the same notation $r$ for the new profile which depends on the variables $(\varphi,\theta)\in\mathbb{T}^{d+1}$.
At this stage, we split the phase space into two parts: the tangential part given by  the  finite dimension space $\mathbb{H}_{\overline{\mathbb{S}}}$  and the normal one  $ \mathbb{H}^{\bot}_{\mathbb{S}_0}$, defined  in \eqref{decoacca}. The dynamics on the tangential space will be parametrized through  the action-angle variables $(I,\vartheta)$ leading to a new reformulation through the embedded torus. Indeed, by virtue of  \eqref{aacoordinates} we may decompose $r$ as follows
$$
r(\varphi,\theta)=\underbrace{v(\vartheta,I)(\varphi,\theta)}_{\in \mathbb{H}_{\overline{\mathbb{S}}}}+\underbrace{z(\varphi,\theta)}_{\in \mathbb{H}^{\bot}_{\mathbb{S}_0}},
$$
with
$$
v (\vartheta, I)(\varphi,\theta) \triangleq \sum_{j \in \overline{\mathbb{S}}}    
 \sqrt{\mathtt{a}_{j}^2+\tfrac{|j|}{2\pi}I_j(\varphi)}\,  e^{\ii (\vartheta_j(\varphi)+j\theta)}. 
$$
In this way, we view $r$ as an embedded torus
\begin{align}\label{embed-torus}
i:\begin{array}[t]{rcl}
\mathbb{T}^{d} & \rightarrow & \mathbb{T}^{d}\times\mathbb{R}^{d}\times \mathbb{H}^{\bot}_{\mathbb{S}_0}\\
\varphi & \mapsto & (\vartheta(\varphi),I(\varphi),z(\varphi)).
\end{array}
\end{align}
Therefore,  in the new coordinates system  the problem reduces  to finding an  invariant torus with non-resonant  frequency vector $\omega$ such that
\begin{equation}\label{Modi-Eq}
\omega\cdot\partial_{\varphi}i(\varphi)=X_{H_{\varepsilon}}(i(\varphi)),
\end{equation}
where the vector field $X_{H_{\varepsilon}}$ is associated to the new Hamiltonian $H_{\varepsilon}$ given in \eqref{cNP} by 
$$
H_{\varepsilon}=-\omega_{\textnormal{Eq}}(\alpha)\cdot I+\tfrac{1}{2}\langle\mathrm{L}(\alpha)z,z\rangle_{L^{2}(\mathbb{T})}+\varepsilon\mathcal{P}_{\varepsilon}\quad\hbox{and}\quad \omega_{\textnormal{Eq}}(\alpha)=\big(\Omega_j(\alpha)\big)_{j\in\mathbb{S}}\cdot
$$
We remind that the frequencies  $\Omega_j(\alpha)$ are defined in \eqref{Freq-equilibrium}. Instead of solving the  equation \eqref{Modi-Eq} we shall first  solve the relaxed problem 
$$
\omega\cdot\partial_{\varphi}i(\varphi)=X_{H_{\varepsilon}^{\tc}}(i(\varphi)),
$$
where the vector field $X_{H_{\varepsilon}^{\tc}}$ is associated to the modified Hamiltonian $H_{\varepsilon}^{\tc}$ given  in \eqref{H alpha} by 
$$H_{\varepsilon}^{\tc}=\tc\cdot I+\frac{1}{2}\langle\mathrm{L}(\alpha)z,z\rangle_{L^{2}(\mathbb{T})}+\varepsilon\,\mathcal{P}_{\varepsilon}.
$$
We emphasize that the advantage of proceeding in this way  is to get one degree of freedom with the vector $\tc$ that will be fixed  later to ensure   some compatibility assumptions in finding an approximate inverse of the linearized operator. At the end of Nash-Moser scheme, we will adjust the frequency $\omega$ in such a way that $\tc$  coincides with the equilibrium frequency $\omega_{\textnormal{Eq}}(\alpha)$ and  thereby  we get a solution to  the original  Hamiltonian equation. To find solutions to the relaxed problem  it suffices  to  construct zeros $(i,\tc)$ to  the nonlinear functional,
\begin{align}\label{PG1-C}
\nonumber\mathcal{F}(i,\tc,\omega,\alpha,\varepsilon)&\triangleq \omega\cdot\partial_{\varphi}i(\varphi)-X_{H_{\varepsilon}^{\tc}}(i(\varphi))\\
&=\left(\begin{array}{c}
\omega\cdot\partial_{\varphi}\vartheta(\varphi)-\tc-\varepsilon\partial_{I}\mathcal{P}_{\varepsilon}(i(\varphi))\\
\omega\cdot\partial_{\varphi}I(\varphi)+\varepsilon\partial_{\theta}\mathcal{P}_{\varepsilon}(i(\varphi))\\
\omega\cdot\partial_{\varphi}z(\varphi)-\partial_{\theta}\left(\mathrm{L}(\alpha)z(\varphi)+\varepsilon\nabla_{z}\mathcal{P}_{\varepsilon}(i(\varphi))\right)
\end{array}\right)
\end{align}
for given $(\omega, \alpha,  \varepsilon)\in \mathbb{R}^d\times[\underline\alpha,\overline\alpha] \times (0,1)$.
We remark that  the  flat  torus $i_{\textnormal{flat}}(\varphi)=(\varphi,0,0)$ is a trivial solution in the particular case $\varepsilon=0$,
$$
\mathcal{F}\big(i_{\textnormal{flat}},-\omega_{\textnormal{Eq}}(\alpha),-\omega_{\textnormal{Eq}}(\alpha),\alpha,0\big)=0.
$$
Then at  this point we are tempted  to apply the classical Implicit Function Theorem but unfortunately it  does not work because  the  linearized operator at the equilibrium state is not invertible due to the small divisors problem. We can remedy to this defect by  imposing suitable non-resonance  conditions on  the frequency $\omega$ and show that this linearized operator admits at least a right inverse but with loss of regularity. Then the main challenge is to extend this property for the linearized operator associated to any arbitrary small state near the flat torus which turns out to be no longer diagonal and admits variable coefficients affecting the main  part of its symbol.   This is considered as the main step  in the implementation of  Nash-Moser scheme.   Inverting this operator requires careful attention and delicate analysis of the resonances set and one needs  to set up some refined tools from toroidal pseudo-differential operators. 

\smallskip

\ding{203}{{ {\it  Approximate right inverse of the linearized operator.}}}
The  structure of the  linearized operator of  the functional $\mathcal{F}$ is  given by a complicated linear combination of the tangential and normal parts interacting through variables coefficients  and it  seems to be out of reach  to invert it in only one step. Then the approach  developed by Berti and Bolle in \cite{BB13}, which  is robust and  has been  performed  in several contexts through different papers \cite{Baldi-berti,BBM-auto, BertiMontalto}, consists in  linearizing first the functional around an isotropic torus sufficiently close to the original one and then proceed with a canonical conjugation through a symplectic change of coordinates leading   to a triangular   system (decoupling the  tangential and the normal parts) up to   small decaying errors, essentially of  "type $Z$" or highly decaying in frequency, that can be incorporated in Nash-Moser scheme. As a by-product of this formalism, to invert this latter triangular system it suffices to  get an approximate  right  inverse  for the linearized operator in the normal direction, denoted in what follows by $\widehat{\mathcal{L}}_{\omega}$.  We will see in Section \ref{sec:Approximate-inverse} that we can bypass the use of isotropic torus using exactly the same formalism. Actually, according to Proposition \ref{Prop-Conjugat}, we  can conjugate the linearized operator with the  transformation  described by \eqref{trasform-sympl}  and find at the end a triangular system with small errors mainly of  "type $Z$". The computations are done directly in a straightforward way and where only the Hamiltonian structure of the original system  sounds to be crucial to get the final triangular structure. Now, let us emphasize  that the transformation \eqref{trasform-sympl} is not  symplectic but  this does not matter because, first we are not interested in the persistence of the nonlinear Hamiltonian structure but simply concerned with its  linear level. Second, this transformation is almost symplectic up to errors of  "type $Z$". 
The main  advantage with this is  to require the invertibility only for   the linearized operator at the torus itself and not on a closer  isotropic one. By this way, we can avoid the accumulation of   different extra errors induced by the isotropic torus that one  encounters for example  in the estimates of the approximate inverse or in the multiple  Cantor sets elaborated along the different reduction steps where the coefficients should be computed at the isotropic torus. Therefore the outcome of this first step is to transform  the invertibility problem of the full operator  to simply invert the partial one resulting from the normal part and taking, according to  Proposition \ref{lemma-GS0} and \eqref{lin}, the form
$$\widehat{\mathcal{L}}_{\omega}=\Pi_{\mathbb{S}_0}^{\perp}\left(\mathcal{L}_{\varepsilon r,\lambda}-\varepsilon\partial_{\theta}\mathcal{R}\right)\Pi_{\mathbb{S}_0}^{\perp}\quad\hbox{with}\quad\mathcal{L}_{r,\lambda} h =\omega\cdot\partial_\varphi h +\partial_\theta\big( V_{r,\alpha}(\varphi,\theta) h  -\mathbb{K}_{r,\alpha} h\big)
$$
where $\Pi_{\mathbb{S}_0}^{\perp}$ is the orthogonal projection on the normal frequency phase space and  $\mathcal{R}$ is a smoothing  integral operator in the space variable $\theta$ with finite rank coming from the interaction between the tangential and the normal parts induced by the transformation \eqref{trasform-sympl}. By virtue of Lemma \ref{lemma-reste}, one has the asymptotic structure,
\begin{align*}
 \mathcal{L}_{r,\lambda}&=\omega\cdot\partial_\varphi+\partial_\theta\Big[ V_{r,\alpha}-\big(\mathscr{W}_{r,\alpha}\, |{\textnormal D}|^{\alpha-1}+|{\textnormal D}|^{\alpha-1}\mathscr{W}_{r,\alpha}\big)+\mathscr{R}_{r,\alpha}\Big]
\end{align*}
with $\mathscr{W}_{r,\alpha}$ being  a smooth function taking non vanishing constant value at the equilibrium state $r\equiv 0$ and the remainder $\mathscr{R}_{r,\alpha}$ is a pseudo-differential  operator of order $-2$ in the spatial variable. The operator $|{\textnormal D}|^{\alpha-1}$ is a modified fractional Laplacian of order $\alpha-1$ and whose kernel representation is detailed in \eqref{fract1}.
At the equilibrium state,  $\widehat{\mathcal{L}}_{\omega}$ coincides with  $\mathcal{L}_{0,\lambda}$ which is a  Fourier multiplier  that can be formally inverted with loss of algebraic regularity  by imposing  suitable  Diophantine conditions that can be guaranteed by  the parameters $(\omega,\alpha)$ which should be in a suitable  massive Cantor set.  However when  $r$ is taken  small,  the perturbation is propagated everywhere and affects all the positive order parts of the operator, that is, the transport part in terms of $V_{r,\alpha}$ and the nonlocal part (of order $\alpha$) through $\mathscr{W}_{r,\alpha}$. This situation is common to different models that we can encounter in the literature such as the   water waves \cite{Alazard-baldi,Baldi-berti,BertiMontalto}. Another observation that we want to stress concerns the order of the fractional Laplacian which depends  on the exterior  parameter $\alpha$. This is a new difficulty compared to  the equations subject to KAM studies, where the exterior parameter affects only the coefficients of the operator but not its order. This fact brings slightly more technical difficulties related to the functional calculus aspects.    Now the question, which is a central key point in KAM theory applied for PDE,  is how one could   invert the linearized operator with variable coefficients. For this aim, the intuitive idea that one could implement is to  diagonalize the operator   using different types of transformations in order  to get it  conjugated to a Fourier multiplier as  in the integrable case. Here in our case,  we distinguish three different reductions related to  the transport part, the fractional part and the remainder. Next, we intend to shed the light on the different techniques used to perform the reduction steps and isolate  the main technical difficulties where some of them are solved with a new approach.

\smallskip

\foreach \x in {\bf a} {%
  \textcircled{\x}
}  {\bf  Reduction of the transport part.}}  This procedure has been discussed recently throughout   several papers and consists in finding a suitable quasi-periodic symplectic change of  coordinates allowing to conjugate the transport part into a new one with constant coefficients \cite{BFM,FGMP19}.  The results related to this point can be found in Section \ref{sec-transport-1}. As we shall see in Proposition \ref{prop-chang} we may find an invertible  transformation
\begin{equation*}
\mathscr{B}h=\big(1+\partial_{\theta}\beta(\varphi,\theta)\big)h\big(\varphi,\theta+\beta(\varphi,\theta)\big)
\end{equation*}
such that for an arbitrary  $n\in\NN$ and if the parameter $\lambda$ is restricted to the  truncated set defined through the first order Melnikov condition (the notation $i$ below stands  for  the embeddings torus parametrization associated to $r$ according to \eqref{embed-torus})
$$\mathcal{O}_{\infty,n}^{\kappa,\tau_{1}}(i)=\bigcap_{(l,j)\in\mathbb{Z}^{d+1}\backslash\{0\}\atop |l|\leqslant N_{n}}\left\lbrace\lambda=(\omega,\alpha)\in \mathcal{O};\; \big|\omega\cdot l+jc(\lambda,i)\big|> 4\kappa^{\varrho}\tfrac{\langle j\rangle}{\langle l\rangle^{\tau_{1}}}\right\rbrace$$
we have 
\begin{align}\label{L0-op}
\mathcal{L}_{r,\lambda}^0\triangleq \mathscr{B}^{-1}\mathcal{L}_{\varepsilon r,\lambda}\mathscr{B}=\omega\cdot\partial_\varphi+c(\lambda,i)\partial_\theta-\partial_\theta\Big(\mu_{r,\lambda}|\textnormal{D}|^{\alpha-1}+\textnormal{D}|^{\alpha-1} \mu_{r,\lambda}\Big)+\partial_\theta\mathcal{R}_{r,\lambda}+\mathtt{E}_{n}^{0},\, 
\end{align}
with $N_n= N_0^{(\frac{3}{2})^n}, N_0\geqslant 2, \kappa,\varrho\in(0,1), \tau_1>d$, $ c(\lambda,i)$  a constant and $\mathtt{E}_{n}^{0}$ being  a linear operator satisfying in particular the decay estimate 
$$\|\mathtt{E}_{n}^{0}h\|_{s_0}^{q,\kappa}\lesssim\varepsilon\kappa^{-1}N_{0}^{\mu_{2}}N_{n+1}^{-\mu_{2}}\|h\|_{q,s_{0}+2}^{q,\kappa}.$$
The norms $\|\cdot \|_{s_0}^{q,\kappa}$ of weighted Sobolev spaces are given in \mbox{Definition \ref{Def-WS}.} As to the  number $\mu_2$, it is connected with the regularity of the torus associated to $r$ and can be taken arbitrary large allowing the remainder $\mathtt{E}_{n}^{0}$ to get sufficient frequency decay and thereby  evacuate it in the small errors  during the   Nash-Moser scheme implemented at the final stage for   the construction of the solutions to the nonlinear problem.
On the other hand, the open set $\mathcal{O}$ is without any importance at this level but later and due to several constraints it  should be an open set containing   the equilibrium frequencies, that is, 
$$
\mathcal{O}=(\underline\alpha,\overline\alpha)\times B(0,R),
$$
where the open  ball $B(0,R)$ with  radius $R$ contains  the equilibrium frequency vector set $\big\{{\omega}_{\textnormal{Eq}}(\alpha),\alpha\in[\underline\alpha,\overline\alpha]\big\}.$
The function $\mu_{r,\lambda}$ is not constant but it is close to  the constant $\mu_{0,\lambda}$ arising  from the equilibrium state. The operator $\partial_\theta\mathcal{R}_{\varepsilon r,\lambda}$ is not highly smoothing but  of order $-1$ and with small  size $\varepsilon.$ The construction of the change of coordinates follows the classical KAM scheme as in \cite{BFM,FGMP19} and consists in writing successive approximations implemented through solving the {\it homological equations}. This step requires to impose  Melnikov first order non-resonance condition for the transport part. The advantage of this scheme is to replace the first transport operator, after conjugation,  with a diagonal one up to a small quadratic error but still of order one. Then iterating this procedure allows to asymptotically get rid of the errors and reduce completely the operator to a diagonal one.  The cost of this scheme is directly reflected onto  the set of parameters which becomes smaller and smaller due to the excision imposed by Melnikov conditions. 
We point out that the final Cantor set $\mathcal{O}_{\infty,n}^{\kappa,\tau_{1}}(i)$ is constructed from the ultimate  coefficient $c(\lambda,i)$ (obtained as the limit of a suitable sequence) and it is  truncated in the time frequency. The advantage of manipulating  truncated Cantor sets, which are actually open sets, arises  at least at two different levels.  First,    during the Nash-Moser scheme where we need  to construct classical extensions  in the whole set of parameters by cut-off functions for the approximations. Second, it is crucial  when  we plan to   estimate  the  final Cantor set from which  nonlinear solutions emerge. There, we need some stability of the different Cantor sets generated throughout  Nash-Moser scheme. This is essential to show  that the final Cantor set is massive and  asymptotically with full Lebesgue measure.

\smallskip

The next step is to reduce the first  nonlocal part in \eqref{L0-op} which is a toroidal  pseudo-differential operator of order $\alpha$.

\smallskip

\foreach \x in {\bf b} {%
  \textcircled{\x}
}    {\bf  Reduction of the fractional Laplacian part for $\alpha<\frac12$}. This step will be explored in Section \ref{Sec-Local-P} where we  reduce the fractional Laplacian part in \eqref{L0-op}  to a Fourier multiplier through the use of infinite dimensional hyperbolic flows.  To be more precise, we show  in  Proposition \ref{prop-constant-coe} that we can construct a family of invertible operators $\Psi(\lambda)$ such that 
\mbox{when  $\lambda\in \mathcal{O}_{\infty,n}^{\kappa,\tau_{1}}(i)$}
\begin{align}\label{Red-Impo1}
\Psi^{-1}\mathcal{L}_{\varepsilon r,\lambda}\Psi\triangleq \mathcal{L}_{r,\lambda}^1\triangleq\omega\cdot\partial_\varphi+c(\lambda,i)\partial_\theta- \textnormal{m}(\lambda,i)\partial_\theta|\textnormal{D}|^{\alpha-1}+\mathcal{R}^1_{r,\lambda}+\mathtt{E}_{n}^{1}
\end{align}
with $m(\lambda,i)$ a constant close to  $m_{0,\lambda}$ obtained at the equilibrium state and  $\mathcal{R}^1_{r,\lambda}$ is a pseudo-differential operator of order $2\alpha-1$ in the spatial variable. The operator error $\mathtt{E}_{n}^{1}$ can be estimated in a similar way to $\mathtt{E}_{n}^{0}$ and will be included in the small remainders in the  Nash-Moser scheme. Let us now outline  the main ideas of the proof which will be done in the spirit of \cite{Baldi-berti}. 
Starting from the first reduction, we want to conjugate  the operator $\mathcal{L}_{r,\lambda}^0$ defined in \eqref{L0-op} to  a new one whose positive order part is a Fourier multiplier. For this aim, we use an infinite-dimensional flow $\Phi$ satisfying the autonomous pseudo-differential hyperbolic equation 
\begin{equation*}
\left\{ \begin{array}{ll}
  \partial_t \Phi =\partial_\theta\big(\rho \,|\textnormal  D|^{\alpha-1}+|\textnormal  D|^{\alpha-1}\rho)\Phi(t)\triangleq \mathbb{A}\Phi(t)  ,&\\ 
   \Phi(0)=\textnormal{Id}.
  \end{array}\right.
\end{equation*}
where the undetermined function $\rho$ will be fixed later according to a transport equation.
Then using Taylor expansion at the first order  for $\Phi(-1)\mathcal{L}_{r,\lambda}^0\Phi(1)$ we find that
\begin{align*}
\nonumber \Phi(-1)\mathcal{L}_{r,\lambda}^0\Phi(1)&=\mathcal{L}_{r,\lambda}^0+\big[\mathcal{L}_{r,\lambda}^0,\mathbb{A}_\rho \big]+\textnormal{l.o.t}\\
&=\omega\cdot\partial_\varphi+c(\lambda,i)\partial_\theta+\partial_\theta\Big(\widehat{\mu}_{r,\lambda}|\textnormal{D}|^{\alpha-1}+|\textnormal{D}|^{\alpha-1} \widehat{\mu}_{r,\lambda}\Big)+\mathcal{R},
\end{align*}
with $\mathcal{R}$ a remainder containing various terms that we shall comment later  and 
$$
\widehat{\mu}_{r,\lambda}\triangleq \big(\omega\cdot\partial_\varphi+c(\lambda,i)\partial_\theta\big)\rho-{\mu}_{r,\lambda}.
$$
Therefore the canonical  choice consists in solving  the following linear transport equation in  $\rho$
$$
\big(\omega\cdot\partial_\varphi+c(\lambda,i)\partial_\theta\big)\rho={\mu}_{r,\lambda}-\langle \mu_{r,\lambda}\rangle_{\varphi,\theta}
$$
where $\langle \mu_{r,\lambda}\rangle_{\varphi,\theta}$ stands for  the  average  of $\mu_{r,\lambda}$  in both variables. This equation can be solved in the periodic setting provided that $\lambda$ belongs to the Cantor set $ \mathcal{O}_{\infty,\infty}^{\gamma,\tau_{1}}(i)$.  Since we are interested in working with  the truncated Cantor set $\lambda\in \mathcal{O}_{\infty,n}^{\kappa,\tau_{1}}(i)$ then instead of solving the above transport equation, we simply solve the  modified one
$$
\big(\omega\cdot\partial_\varphi+c(\lambda,i)\partial_\theta\big)\rho=\Pi_{N_n}\big({\mu}_{r,\lambda}-\langle \mu_{r,\lambda}\rangle_{\varphi,\theta}\big),
$$
where $\Pi_{N}$ is a projector localizing in frequency at the modes range $|l|\leqslant N$, see Lemma \ref{L-Invert}. Then according to this lemma  we may  solve this equation with a loss of regularity but uniformly in $n.$ Notice that this frequency cut-off induces small errors that can be balanced into the error operator denoted by $\mathtt{E}_{n}^{1}$ which admits a fast decaying rate.  Let us now discuss some important technical  aspects  in estimating the remainder $\mathcal{R}.$ By referring  to  \eqref{LrYP0P} and \eqref{hatR},  one should deal with some remainders taking  the form $\Phi(-t)\mathcal{R}_0\Phi(t)$ with $\mathcal{R}_0$ is typically described through negative order commutators between pseudo-differential operators.
The delicate point is to find a suitable operator topology for which the next reduction on the remainder term, detailed in \foreach \x in {\bf c} {%
  \textcircled{\x}},  works and is compatible with the  KAM scheme. This step requires an adequate    topology satisfying   tame estimates in the scales of Sobolev spaces and obviously the bounded operator topology given in \eqref{strong-Top1} is not well adapted and seems to be  too weak for this purpose. Then one   common way   is to use a reinforced  topology on the symbols class  of order zero described through   \eqref{equivalent-norm} with  its weighted variant in  \eqref{Top-NormX}. This was for example used in the papers  \cite{Baldi-berti,BertiMontalto}. Now with this topology in mind we get the required functional tame estimates as it is indicated in Lemma \ref{comm-pseudo1}. However, and this is the main technical issue in this part,  it is not at all clear whether the new  remainder $\mathcal{R}$  still belongs to this topology with  tame estimates type.  To be more concrete, if we take an operator $\mathcal{R}_0$ whose symbol satisfies  \eqref{Top-NormX}, is it true that the conjugation by the flow  $\Phi(-t)\mathcal{R}_0\Phi(t)$  remains in the same class  with  suitable  tame estimates? We point out that the hyperbolic flow $\Phi(t)$  which acts continuously on the classical Sobolev spaces or their weighted version, see Proposition \ref{flowmap00}, does not  in general satisfy  the constraint \eqref{Top-NormX} which sounds to be too strong   for it. At a formal level, one has 
  $$
  \Phi(t)\sim \textnormal{Op}\left(e^{\ii\rho(\varphi,\theta)\xi|\xi|^{\alpha-1}}\right)
  $$
  then the differentiation in $\varphi$ or $ \theta$ will generate an unbounded operator  with a loss of $\alpha$ derivative.
  Therefore any attempt to use the law products turns to be  in vain and one should implement refined tools and this is one of our main technical contribution in this paper that will be investigated along  \mbox{Section \ref{Flows and Egorov theorem type}.}  Actually, in Theorem \ref{Prop-EgorV}, which can be understood as a refined Egorov type theorem, we shall  give a positive answer to our question  for a large class of pseudo-differential operators satisfying some weak constraints on the symbol structure, which is  large enough to include all the terms of $\mathcal{R}.$ The proof of this theorem is based on a new approach  that will be performed through  several steps based on the kernel dynamics combined with  pseudo-differential calculus and  some algebraic structures related  to T\"oplitz matrix operators introduced in Subsection \ref{Sect-Top-M}. In what follows, we intend to give some insights  on the proof of this theorem in the particular case of the norm $\interleave\Phi(-t)\mathcal{R}_0\Phi(t)\interleave_{0,s,0},$ for the definition  see \eqref{Def-Norm-M1}. Let   $\mathcal{K}_t$ denote  the kernel of the operator $\mathcal{R}_1(t)\triangleq \Phi(-t)\mathcal{R}_0\Phi(t)- \mathcal{R}_0$, then it  satisfies the following  transport equation, see  \eqref{kernel-dyna},
  \begin{equation*}
 \mathscr{L}\mathcal{K}_t(\varphi,\theta,\eta)\triangleq\partial_t\mathcal{K}_t(\varphi,\theta,\eta)-\big(\mathbb{A}_\theta+\mathbb{A}_\eta\big) \mathcal{K}_t(\varphi,\theta,\eta)+\big[\partial_\eta,\mathcal{T}_\eta \big]\mathcal{K}_t(\varphi,\theta,\eta)=K_0(\varphi,\theta,\eta),
\end{equation*}
with $K_0$  the kernel associated to $ [\mathbb{A},\mathcal{R}_0]$ and  $$
\mathbb{A}_\theta =\partial_\theta\big(\rho(\varphi,\theta)|\textnormal{D}_\theta|^{\alpha-1}+|\textnormal{D}_\theta|^{\alpha-1}\rho(\varphi,\cdot)\big)\triangleq \partial_\theta\mathcal{T}_\theta.
$$ 
According to the norm definition \eqref{Def-Norm-M1} and the relationship between the symbol and the kernel one has for any $n\in\N$ that
\begin{align*}
\interleave\mathcal{R}_1(t)\interleave_{0,n,0}
&\lesssim\sum_{i=0}^n\big(\|\partial_\chi^i\mathcal{K}_t\|_{L^2(\T^{d+2})}+\|\partial_\varphi^i\mathcal{K}_t\|_{L^2(\T^{d+2})}\big),\quad\partial_\chi\triangleq\partial_\theta+\partial_\eta
\end{align*}
Let us explain how to get the estimates of $\|\partial_\chi^i\mathcal{K}_t\|_{L^2(\T^{d+2})}$ stated in \eqref{Y-KL1}. We first establish the following vectorial equation: 
 \begin{equation*}
\mathcal{L}_n\mathcal{X}_n(t)=\mathcal{M}_n\mathcal{X}_n(t)+\mathcal{Y}_n\quad\hbox{with}\quad  \mathcal{L}_n\triangleq \mathscr{L}\textnormal{I}_{n+1},\quad \mathcal{X}_n(t)=\begin{pmatrix}
     \mathcal{K}_t &  \\
    \partial_\chi \mathcal{K}_t& \\
    ..&\\
    ..&\\
    \partial_\chi^n  \mathcal{K}_t
  \end{pmatrix},\quad\mathcal{Y}_n=\begin{pmatrix}
{K}_0 &  \\
    \partial_\chi{K}_0& \\
    ..&\\
    ..&\\
    \partial_\chi^n  {K}_0
  \end{pmatrix}
  \end{equation*}
with $\textnormal{I}_{n+1}$ is the identity matrix and $ \mathcal{M}_n$ being a lower T\"oplitz triangular  matrix operator
\begin{equation*}
 \mathcal{M}_n=\begin{pmatrix}
    0 &0&..&..&..&0  \\
 m_{1,1}&0&..&..&0&0 \\
   m_{2,2}& m_{1,2}&0&.. &..&0\\
    ..&..&..&..&0&0\\
 ..& ..&..&m_{1,n-1}& 0&0\\
 m_{n,n}& m_{n-1,n}&..&..& m_{1,n}&0
  \end{pmatrix}.
\end{equation*}
Notice that the entries of this matrix are themselves pseudo-differential operators with strictly  positive order $\alpha$ and therefore the matrix $ \mathcal{M}_n$ is not dissipative and we cannot at this stage implement energy estimates. For the structure of the entries we refer to \eqref{Bnk}. The key point is that the matrix $ \mathcal{M}_n$ is nilpotent of order $n+1$ and enjoying a rich  T\"oplitz structure. Thus  the suitable quantities to estimate are given by the iterated vectors $\big\{\mathcal{M}^k\,\mathcal{X}_n, \,\,k=0,.., n\big\}$ and one gets the equations
 \begin{align*}
\mathcal{L}_n\mathcal{M}_n^k\mathcal{X}_n(t)
&=-\sum_{j=0}^{k-1}\left(_j^k\right)\textnormal{Ad}_{\mathcal{M}_n}^{k-j}(\widehat{\mathcal{L}}_n)\, \mathcal{M}_n^{j}\mathcal{X}_n(t)+\mathcal{M}_n^{k+1}\mathcal{X}_n(t)+\mathcal{M}_n^{k}\mathcal{Y}_n.
\end{align*}
The estimate of the first term of the right-hand side is detailed in Lemma \ref{lem-com-it2} where we use in a crucial way the T\"oplitz structure not only for the matrix   $\mathcal{M}_n$ but also for the iterated matrices commutators $\textnormal{Ad}_{\mathcal{M}_n}^{k-j}(\widehat{\mathcal{L}}_n).$  Notice that we get in particular that the entries   of each matrix commutator are themselves commutators of scalar pseudo-differential operators of order $\alpha$, and therefore they are of negative order acting continuously in $L^2$. The tame estimates are very subtle and require more algebraic structure related to the number of vanishing lower sub-diagonals in the commutators combined with refined estimates on the  commutators as in Lemma  \ref{comm-pseudo}. As to the estimate of $\mathcal{M}_n^{k}\mathcal{Y}_n$, it is stated in Lemma \ref{V-13} and it is very involved due to the fact that $\mathcal{Y}_n$ is not smooth and it is singular at the diagonal line $\{\theta=\eta\}$. Then when  the T\"oplitz matrices $\mathcal{M}_n^k$ act on $\mathcal{Y}_n$  we should check that the singular contributions are cancelled  as for $\partial_\chi^k K_0$, in this latter  case the situation is easy since  the operator $\partial_\chi=\partial_\theta+\partial_\eta$ is local and the diagonal singularity belongs to  its kernel. The delicate point is that  the coefficients of $\mathcal{M}_n^k$ are very complicated and they are non local pseudo-differential operators and we should prove that the singularity is not affected by the iterated operators. The proof is based on suitable  recursive estimates governing the iterated kernels, see Lemma \ref{lemm-iter1}.

\smallskip

We want to precise that similar problem   occurs for gravity-capillary water waves as in \cite{BertiMontalto} where the authors proceed in a different way to deal with $\mathcal{R}_t\triangleq\Phi(-t)\mathcal{R}_0\Phi(t)$ based on the construction of an approximate solution to the  Heisenberg equation satisfied by $\mathcal{R}_t\triangleq\Phi(-t)\mathcal{R}_0\Phi(t)$, that is, $$
\partial_t \mathcal{R}_t=\big[\mathbb{A},\mathcal{R}_t\big].
$$
The approximation is given in  the spirit  of the  proof of Egorov theorem by expanding the symbol into a finite sum of symbols with decreasing order  associated to iterated commutators. Then each element of this sum belongs to the symbol class topology  \eqref{Top-NormX}, however  the error term operator is simply smoothing in the spatial variable and subject to a loss of regularity in the time variable due to the flow. Then to remedy to this lack of information on the error term the authors introduced \mbox{in \cite{BertiMontalto}} new "abstract" topology classes named  $\mathcal{D}^k$-tame and $\mathcal{D}^k$-modulo-tame operators, much weaker \mbox{than \eqref{Top-NormX}.} Compared to our case, we are able in this paper to check that the error operator  still  belongs to the  symbol class topology  \eqref{Top-NormX}, which has not only  an interest in itself but it brings  different technical simplifications  related to the use  the functional tools during the next steps.

  The outcome of this long and technical process is to achieve the reduction \eqref{Red-Impo1}, where the new remainder  $\mathcal{R}^1_{r,\lambda}$ is a pseudo-differential operator of negative order   with a tame estimate in the suitable strong topology over the symbol class. For a precise statement we refer to \mbox{Proposition \ref{prop-constant-coe}-(ii).}

\smallskip

\foreach \x in {\bf c} {%
  \textcircled{\x}
}   {\textbf{ KAM reduction of the remainder}}.
As we have seen  in  \foreach \x in {\bf b} {%
  \textcircled{\x},
} we obtain after  two reductions a new operator given by  \eqref{Red-Impo1} whose positive part is diagonal and the remainder is of order zero in all the variables and with small size. Then, in order to  diagonalize the new operator up to enough decaying small remainder we shall implement a KAM scheme as  for instance in \cite{Baldi-berti,Baldi-Berti-Montalto14,BertiMontalto}. Notice that  before proceeding in that way we should first implement the steps \foreach \x in {\bf a} {%
  \textcircled{\x}
} and \foreach \x in {\bf b} {%
  \textcircled{\x}
} with the  linearized operator restricted on the normal direction, given by \eqref{Norm-local-z}. This step is discussed in Section \ref{Norm-sec-1}, and roughly speaking  the localization into the normal direction induces new terms that can be described through finite rank operators with small sizes.
More precisely,    one gets according to Proposition \ref{projection in the normal directions} a suitable invertible operators $\Psi_{\perp}$ acting on the normal direction such that on 
  the  truncated Cantor set $\mathcal{O}_{\infty,n}^{\kappa,\tau_{1}}(i)$   we have
\begin{align*}
\Psi_{\perp}^{-1}\widehat{\mathcal{L}}_{\omega}\Psi_{\perp}&=\big(\omega\cdot\partial_\varphi+{c}(\lambda,i)\partial_\theta- m(\lambda,i)\partial_\theta|\textnormal{D}|^{\alpha-1}\big)\Pi_{\mathbb{S}_0}^{\perp}+\mathcal{R}^2_{r,\lambda}+\mathtt{E}_n^2\triangleq\mathscr{L}_0+\mathtt{E}_n^2, 
\end{align*}
where the orthogonal projector $\Pi_{\mathbb{S}_0}^\perp$  and the set  $\mathbb{S}_0$ are defined in \eqref{projectors-tan-normal} and \eqref{tangent-set2}, respectively.  The operator  $\mathtt{E}_n^2$ is similar to   $\mathtt{E}_n^1$ seen in  \foreach \x in {\bf b} {%
  \textcircled{\x}
} and describes different small  errors coming from  the time  truncation of the Cantor set $\mathcal{O}_{\infty,n}^{\kappa,\tau_{1}}(i)$.  
At this stage we get the following structure 
$$
\mathscr{L}_{0}=\omega\cdot\partial_{\varphi}+\mathscr{D}_{0}+\mathscr{R}_{0}
$$
where $\mathscr{D}_{0}$ is a diagonal operator and $\mathscr{R}_{0}$ is a {\it good} pseudo-differential operator of zero order and with small size,  and satisfying  in addition the  reversibility structures. Actually, $\mathscr{R}_{0}$ is of order $2\alpha-1$ in the spatial variable $\theta$. The KAM reduction result is stated in  Proposition \ref{reduction of the remainder term} and its ultimate target is   to eliminate the remainder $\mathscr{R}_{0}$ and transform it into a diagonal part  by suitable conjugation of $\mathscr{L}_{0}$. This will be developed in a standard way by constructing successive transformations through the KAM reduction  allowing to replace at each step the  remainder with a  smaller new one provided that we make the suitable parameters excision. This scheme works well if  we could solve in a reasonable way  the associated  {\it{homological equation}} which requires to eliminate resonances at each step by imposing {\it  the second order Melnikov condition} using the external   parameters. The final outcome can be summarized  as follows: there exists a reversible  invertible linear operator $\Phi_\infty$ defined in the full set of parameters such that on the truncated Cantor set 
\begin{align*}
\mathcal{O}_{{\infty,n}}^{\gamma,\tau_1,\tau_{2}}(i)\triangleq \Big\{\lambda=(\omega,\alpha)\in\mathcal{O}_{{\infty,n}}^{\gamma,\tau_{1}}(i);\;&\forall |l|\leqslant {N_n},\,j,j_{0}\in\mathbb{S}_{0}^{c},\,(l,j)\neq(0,j_{0}),\\
&\big|\omega\cdot l+\mu_{j}^{\infty}(\lambda,i)-\mu_{j_{0}}^{\infty}(\lambda,i)\big|>\tfrac{\kappa \langle j-j_{0}\rangle }{\langle l\rangle^{\tau_{2}}}\Big\}
\end{align*}
we have 
\begin{align}\label{Diagon-za}
\Phi_{\infty}^{-1}\mathscr{L}_{0}\Phi_{\infty}=\omega\cdot\partial_{\varphi}+\mathscr{D}_{\infty}+{\mathtt{E}^3_n}
\end{align}
where $\mathscr{D}_{\infty}=\left(\ii\mu_{j}^{\infty}(\lambda,i)\right)_{(l,j)\in\mathbb{Z}^{d}\times\mathbb{S}_{0}^{c}}$ is a diagonal reversible operator and the error operator ${\mathtt{E}^3_n}$ is similar to  ${\mathtt{E}^2_n}$ seen before. In addition, the  eigenvalues admit the following asymptotic expansion
\begin{align*}
\mu_{j}^{\infty}(\lambda,i)
&=\Omega_{j}(\alpha)+j\,r^1(\lambda,i)+j|j|^{\alpha-1}\,r^2(\lambda,i)+r_{j}^{\infty}(\lambda,i)
\end{align*}
where $\Omega_j$ is defined in \eqref{Freq-equilibrium}, $r^1, r^2$ and $r_{j}^{\infty}$ satisfies the frequency decay
$$
\max_{k=1,2}|r^k(\lambda,i)|\lesssim \varepsilon \kappa^{-1}\quad\hbox{and}\quad \sup_{j\in\mathbb{S}_{0}^{c}}|j|^{1-\epsilon-2\overline\alpha}| r_{j}^{\infty}(\lambda,i)|\leqslant C\varepsilon\kappa^{-2}.
$$
Actually,   we get more precise estimates related to the parameter differentiation of these perturbed coefficients. In fact, one  needs to control their derivatives $\partial^{k}_\lambda$ for any $0\leqslant k\leqslant q_0+1$, where $q_0$ is the index of  non-degeneracy  of the equilibrium frequencies constructed in Proposition \ref{lemma transversality}. More details will be given later on the measure estimates of the final Cantor set.

\smallskip

The final conclusion of this step is the construction of  an approximate inverse for the linearized operator restricted to the normal direction  as  stated in  Theorem  \ref{inversion of the linearized operator in the normal directions}. This  yields in turn  an approximate inverse for the full  linearized operator at an arbitrary state close to the equilibrium, according to Proposition \ref{thm:stima inverso approssimato}.

\smallskip

\ding{204 }{\it Nash-Moser scheme.\,}
This is the main  purpose of Section \ref{N-M-S1} where we 
 construct  solutions for the nonlinear function $\mathcal{F}$, defined in \eqref{PG1-C}, and   localized  around  a finite dimensional linear torus with a small   amplitude. We basically  follow   a modified Nash-Moser scheme as in the papers \cite{Baldi-berti,BCP,BertiMontalto} with slight variations. Notice that the Cantor sets used in getting an approximate inverse of  the linearized operator, see Theorem \ref{inversion of the linearized operator in the normal directions},  are truncated in the time frequency and associated to the final states.  We emphasize one time more that  working with truncated sets allows to  generate classical smooth  extensions, with a suitable frequency decay, to the whole set of parameters for the approximations defined a priori on the Cantor sets. Dealing with smooth extensions is useful later in measuring the final Cantor set where we shall implement perturbative arguments based on R\"ussmann techniques developed in \cite{Russ} and extended in \cite{bambusi-Berti}. 
 
\smallskip

 The conclusion of this step  is summarized in Corollary \ref{cor-ima12} where solutions to the nonlinear functional \eqref{PG1-C}  are constructed provided that the $(\omega,\alpha)$ belongs to a final Cantor set given by the intersection of all the intermediate Cantor sets needed during the Nash-Moser scheme. The ultimate  point to check, and which  achieves the proof of the main theorem,  is to show that the final Cantor set is massive. 
 
\smallskip

\ding{205}{\it  Non-degeneracy and Measure estimates.}
This is the main goal of Section \ref{Section 6.2} where we prove that the final Cantor set $\mathtt{C}_{\infty}^{\kappa,\varepsilon}$ described in  Corollary \ref{cor-ima12} is asymptotically with full Lebesgue measure. This set  is given by
\begin{equation*}
\mathtt{C}_{\infty}^{\kappa,\varepsilon}=\bigcap_{m\in\mathbb{N}}\mathtt{C}_{m}^{\kappa}\quad \mbox{ where }\quad \mathtt{C}_{m}^{\kappa}=\Big\{\alpha\in(\underline\alpha,\overline\alpha);\;\lambda(\alpha,\varepsilon)\in\mathcal{A}_m^{\kappa}\Big\},
\end{equation*}
where the intermediate sets $\mathcal{A}_m^{\kappa}$ are defined in Proposition \ref{Nash-Moser} and  constructed   along Nash-Moser scheme through non resonance conditions imposed at the approximate torus  $i_m$. To measure this set, we mainly use the stability of the Cantor sets combined with  the techniques developed in \cite{bambusi-Berti} and \cite{Russ}. One of the crucial key point is the transversality property stated in Lemma \ref{lemma R\"ussmann condition for the perturbed frequencies}. It will be first established for the  linear frequencies in Proposition \ref{lemma transversality}, using the analyticity of the eigenvalues and their asymptotic behavior combined with the poles structure of Gamma function. Then the  extension of the transversality assumption to the perturbed frequencies is done using perturbative arguments together with the asymptotic description of the approximate  eigenvalues detailed  in \eqref{asy-z1}, \eqref{uniform estimate rjinfty} and \eqref{uniform estimate r1}. Notice that the transversality is connected with the non-degeneracy of the eigenvalues in the sense of the Definition \ref{non-degeneracy}. For instance, we show that  the curve  $\alpha\in[0,1) \mapsto\big(\Omega_{j_1} (\alpha),...,\Omega_{j_d} (\alpha)\big)$ is  not contained in any vectorial plane, that is, if there exists a constant vector $c=(c_1,..,c_d)$ such that
$$
\forall \alpha\in[0,1),\quad \sum_{j=1}^d c_k \Omega_{j_k}(\alpha)=0
$$ 
then $c=0$. We point out that for most of the equations studied before  we may check this property by making Taylor expansion at zero leading to an invertible Vandermonde matrix type. However in the current  case the structure of the eigenvalues sounds at the first sight more complicated due to the Gamma quotient $\alpha\in\C\mapsto \frac{ \Gamma(j+\frac{\alpha}{2})}{\Gamma(1+j-\frac{\alpha}{2})}$ which has the advantage to be  a meromorphic function admitting countable sets of separated poles and zeroes altogether located at the real axis. For instance the zeroes are given by
$$
\mathcal{Z}_j\triangleq \big\{\alpha=2(n+j+1);\; n \in \mathbb{N}\big\}.
$$  
Then using the different list of zeroes we find successively $c_d=c_{d-1}=\cdots=c_1=0.$ \\
In Proposition \ref{lem-meas-es1}, we provide the following lower bound of the Lebesgue measure of the $\big|\mathtt{C}_{\infty}^{\kappa,\varepsilon}\big|$,
$$
\big|\mathtt{C}_{\infty}^{\kappa,\varepsilon}\big|\geqslant \overline\alpha-\underline\alpha-C\varepsilon^{\eta},
$$
for a suitable small number $\eta>0$. We notice that at this level we make the choice $\kappa=\varepsilon^a$ with $a$ being a sufficiently small number.

 \section{Hamiltonian reformulation}
 
This section is devoted to the contour dynamics equation governing the patch motion. We shall in particular reformulate the equations in the polar coordinates whose validity is guaranteed for any time for  periodic or quasi-periodic solutions near the unit circle. As we shall see, this can be transformed into  a Hamiltonian reformulation which is essential in the approach used for the construction  of  quasi-periodic solutions. 
\subsection{Contour dynamics equation}\label{sec contour}
We denote $\mathbb{D}$ the unit disc of $\mathbb{R}^{2}$ endowed with its usual Euclidean structure. Here and in the sequel, we identify $\mathbb{C}$ with $\mathbb{R}^{2}.$ In particular, the Euclidean structure of $\mathbb{R}^{2}$ is seen in the complex sense through the usual inner  product defined for all $z_{1}=a_{1}+\ii b_{1}\in\mathbb{C}$ and $z_{2}=a_{2}+\ii b_{2}\in\mathbb{C}$ by
\begin{equation}\label{definition of the scalar product on R2}
z_{1}\cdot z_{2}\triangleq \langle z_{1},z_{2}\rangle_{\mathbb{R}^{2}}=\mbox{Re}\left(z_{1}\overline{z_{2}}\right)=a_{1}a_{2}+b_{1}b_{2}.
\end{equation}
We point out that the Rankine vortex   $\mathbf{1}_{\mathbb{D}}$ (and actually any radial function) is a stationary solution to \eqref{Eq-F}. Then  to find    quasi-periodic vortex patches $t\mapsto\mathbf{1}_{D_{t}}$ around this trivial solution it is convenient to consider  the polar parametrization of the boundary given by:
\begin{align}\label{Param-curv}
w:\begin{array}[t]{rcl}
\mathbb{R}_{+}\times[0,2\pi] & \mapsto & \mathbb{C}\\
(t,\theta) & \mapsto & w(t,\theta)
\end{array}\quad\mbox{ with }\quad w(t,\theta)\triangleq \left(1+2r(t,\theta)\right)^{\frac{1}{2}}e^{\ii \theta}.
\end{align}
Here $r$ is the radial deformation of the patch which is small, namely $|r(t,\theta)|\ll1.$ Taking $r=0$ gives a parametrization of the  unit circle $\mathbb{T}.$ 
 To alleviate the notation we use sometimes \begin{equation}\label{definition of R}
R(t,\theta)\triangleq \left(1+2r(t,\theta)\right)^{\frac{1}{2}}.
\end{equation}
This particular form of $R$  is required later  for the Hamiltonian structure. Let us now write down the contour dynamics equation with the polar coordinates. It is known, see for instance \cite{HH}, that the particles on the boundary move with the flow and remain at the boundary  and therefore in the smooth case one has   
$$\left[\partial_{t}w(t,\theta)-\mathbf{v}(t,w(t,\theta))\right]\cdot\mathbf{n}(t,w(t,\theta))=0$$
where $\mathbf{n}(t,w(t,\theta))$ is the outward normal vector to the boundary $\partial D_{t}$ of $D_{t}$ at the point $w(t,\theta)$. 
 Since  one has, up to a real constant of renormalization, $\mathbf{n}(t,w(t,\theta))=-\ii\partial_{\theta}w(t,\theta)$, then  we get the complex formulation of the contour dynamics motion, 
\begin{equation}\label{complex vortex patch equation}
\mbox{Im}\left(\left[\partial_{t}w(t,\theta)-\mathbf{v}(t,w(t,\theta))\right]\overline{\partial_{\theta}w(t,\theta)}\right)=0.
\end{equation}
The velocity field is given by
\begin{align}\label{velocity}
\nonumber\mathbf{v}(t,w(t,\theta))&=\frac{C_\alpha}{2\pi}\int_{\partial D_{t}}\frac{1}{|w(t,\theta)-\xi|^\alpha} d\xi\\
&=\frac{ C(\alpha)}{2\pi}\int_{0}^{2\pi}\frac{ \partial_{\eta}w(t, \eta)}{\vert w(t,\theta)-w(t,\eta)\vert^\alpha} d\eta 
\end{align}
We are interested in patches that are perturbation from the unit disc, therefore, we consider the parametrization
\begin{equation}\label{z}
w(t,\theta)\triangleq \big(1+ 2r(t,\theta)\big)^{\frac12}e^{\ii\theta},
\end{equation}
Differentiating  $z(t,\theta)$ with respect to $t$ we get
\begin{equation}
\partial_t w(t,\theta)=\big(1+ 2r(t,\theta)\big)^{-\frac12}\partial_t r(t,\theta)e^{\ii\theta},
\end{equation}
and, with respect to the $\theta$,
\begin{equation}\label{part-z-theta}
  \partial_\theta w(t,\theta)= \Big(\big(1+ 2r(t,\theta)\big)^{-\frac12}\partial_\theta r(t,\theta)+\ii\big(1+ 2r(t,\theta)\big)^{\frac12}\Big)e^{\ii\theta}\cdot
\end{equation}
Combining the two last identities gives
\begin{align}\label{p1}
\textnormal{Im}\Big\{\partial_t w(t, \theta)\partial_\theta \overline{w(t,\theta)} \Big\}&=- \partial_tr(t,\theta).
\end{align}
and
\begin{align}\label{p3}
&
 \bigintsss_{0}^{2\pi}\frac{ {\rm Im} \big\{\partial_{\eta}w(t, \eta)\partial_\theta \overline{w(t, \theta) }\big\} }{\vert w(t, \theta)-w(t, \eta)\vert^\alpha} d\eta \,
 =\bigintsss_{0}^{2\pi} \frac{\partial^2_{\theta\eta} \big[\big(1+ 2r(t,\theta)\big)^{\frac12}\big(1+ 2r(t,\eta)\big)^{\frac12}\sin(\eta-\theta)\big]}{A_r^{\alpha/2}(\theta,\eta)} d\eta\cdot 
\end{align}
where
\begin{align}\label{A}
A_r(\theta,\eta)&=2\Big[1+ r(t,\theta)+r(t,\eta)-\big(1+ 2r(t,\eta)\big)^{\frac12}\big(1+ 2r(t,\theta)\big)^{\frac12}\cos(\eta-\theta)\Big].
\end{align}
Inserting the identities \eqref{p1}and \eqref{p3} into equation  \eqref{complex vortex patch equation} we get
\begin{align}\label{eq}
\partial_tr(t, \theta)&
+\underbrace{\frac{ C_\alpha }{2\pi }\bigintsss_{0}^{2\pi} \frac{\partial^2_{\theta\eta} \big[\big(1+ 2r(t,\theta)\big)^{\frac12}\big(1+ 2r(t,\eta)\big)^{\frac12}\sin(\eta-\theta)\big]}{A_r^{\alpha/2}(\theta,\eta)} d\eta}_{\triangleq F_\alpha[r](t,\theta )}
=0\cdot
\end{align}
we get  the vortex patch equation in the polar coordinates
\begin{equation}\label{equation0}
\partial_{t}r(t,\theta)
 +F_{\alpha}[r](t,\theta)=0.
\end{equation}
which is a nonlocal transport nonlinear PDE. We look for solutions 
  of the form
\begin{equation}
r(t, \theta)=\tilde r(t,\theta+\Omega t)
\end{equation} 
for some $\Omega>0$. Notice that the introduction of the angular velocity $\Omega$ is purely technical and   needed later  to circumvent  the degeneracy of the first frequency of the linearized operator at the equilibrium state.  
From \eqref{velocity} one gets 
\begin{equation}\label{u-rot}
F_{\alpha}[\tilde r](t,\theta{+\Omega t})= 
F_{\alpha}[ r](t,\theta) \, , 
\end{equation}
Thus, equation  \eqref{equation0} becomes 
\begin{equation}\label{equation}
\partial_{t}\tilde r(t,\theta)+\Omega\partial_{\theta} \tilde r(t,\theta)
 +F_{\alpha}[\tilde r](t,\theta)=0.
\end{equation}
Now, a time quasi-periodic solution of \eqref{equation} is nothing but  a solution in  the form  
$$
\tilde r(t,\theta)=\widehat{r}(\omega t,\theta),$$
where  $\widehat{r}:\,(\varphi,\theta)\in\mathbb{T}^{d+1}\mapsto \,\widehat{r}(\varphi,\theta)\in \mathbb{R}$ and  $\omega\in\mathbb{R}^{d}$ is a non-resonant vector frequency. Hence in this setting,  the equation \eqref{equation} becomes  
\begin{equation*}
\omega\cdot\partial_{\varphi}\widehat{r}(\varphi,\theta)+\Omega\partial_{\theta}\widehat{r}(\varphi,\theta)+F_{\alpha}[\widehat{r}](\varphi,\theta)=0.
\end{equation*}
In the sequel, we shall  alleviate the notation and denote $\widehat{r}$ simply by $r$ and the foregoing equation writes
\begin{equation}\label{equation with quasi-periodic ansatz}
\forall (\varphi,\theta)\in\mathbb{T}^{d+1},\quad\omega\cdot\partial_{\varphi}{r}(\varphi,\theta)+\Omega\partial_{\theta}{r}(\varphi,\theta)+F_{\alpha}[{r}](\varphi,\theta)=0.
\end{equation}

\subsection{Conservation laws}
This section is devoted to some conservation quantities that will be used to recover the Hamiltonian structure of the dynamical system \eqref{equation}. In our context, we shall explore  the circulation $C$,  the angular momentum $J$ and  the energy $E$ and write them in the patch setting using polar coordinates. 
Let $D(t)$ be a bounded simply connected region with smooth boundary $\partial D(t)$ and $\omega=\chi_{D(t)}$. We define the circulation $C$,  the angular momentum $J$ and  the energy $E$   as follows
\begin{equation}\label{J}
C(t)\triangleq \int_{D(t)}dA(z), \quad J(t)\triangleq -\frac12\int_{D(t)}|z|^2dA(z) \quad\hbox{and}\quad E(t)\triangleq \frac12\int_{D(t)} \psi (t, z) d A(z),
 \end{equation}
where $\psi$ is the stream function defined by 
\begin{equation}\label{psi}
\psi(z)\triangleq -\,\frac{ C(\alpha)}{2\pi}\bigintsss_{D(t)}\frac{1}{\vert z-\zeta\vert^\alpha} dA(\zeta)
\end{equation}
and $dA$ being the planar Lebesgue measure.
\begin{lemma} 
The following assertions hold true.
\begin{enumerate}
\item The quantities $C$, $J$ and $E$ are conserved in time.
\item In polar coordinates, they can be expressed in the form
\begin{align*}
 C(t)&=\pi+\int_{0}^{2\pi}r(t,\theta)d\theta,\\
 J(t)&=-\tfrac18\int_{0}^{2\pi}\big(1+2r(t,\theta)\big)^2 d\theta,\\
E(t)&=-\tfrac{ C(\alpha)}{16\pi \left(1-\frac\alpha2\right)^2}\bigintsss_{0}^{2\pi}\bigintsss_{0}^{2\pi}\frac{\partial^2_{\theta\eta}\big[(1+2r(t,\theta))^{\frac12}(1+2r(t,\eta))^{\frac12}\cos(\eta-\theta)\big]}{A_r^{\frac\alpha2-1}(\eta,\theta)}d\eta d\theta,
\end{align*}
where $A_r(\theta,\eta)$ is defined in \eqref{A}

\end{enumerate}
 
\end{lemma}
\begin{proof}
{\rm{(i)}} This result can be checked  by using any  flow associated to  the velocity field generated by the smooth patches ${\bf{1}}_{D_t}.$ 

\smallskip

{\rm{(ii)}}  We first recall    Green's formula written in the complex form,
\begin{equation}\label{green}
\forall z\in \mathbb{C},\quad \int_{D}\partial_{\overline{\zeta}}f(\zeta,\overline{\zeta})){dA(\zeta)}=\frac{1}{2 \ii}\int_{\partial D}f(\zeta,\overline{\zeta})d\zeta.
\end{equation}
In view of \eqref{J} and  \eqref{green} one has
  \begin{equation*}
  C(t)=\frac{1}{2\ii}\int_{\partial D(t)}\overline{\zeta} d\zeta,\qquad J(t)=-\frac{1}{8\ii}\int_{\partial D(t)}|\zeta|^2\overline{\zeta} d\zeta.
  \end{equation*}
  From \eqref{z} we get
   \begin{align*}
C(t)&= \frac{1}{2 \ii }\int_{0}^{2\pi}(1+2r(t,\theta))^{\frac12}e^{-\ii\theta} 
\partial_\theta\Big[(1+2r(t,\theta))^\frac12e^{\ii\theta}\Big] d\theta 
\\ &=\pi+\int_{0}^{2\pi}r(t,\theta)d\theta \, . 
\end{align*}
  and
  \begin{align}
J(t)&= -\frac{1}{8 \ii }\int_{0}^{2\pi}(1+2r(t,\theta))^{\frac32}e^{-\ii\theta} 
\partial_\theta\Big[(1+2r(t,\theta))^\frac12e^{\ii\theta}\Big] d\theta \nonumber
\\ &=-\frac{1}{8}\int_{0}^{2\pi}(1+2r(t,\theta))^{2}d\theta \, . \label{forJ}
\end{align}
 As to the energy expression, we shall first reformulate the stream function defined in \eqref{psi} in terms of polar coordinates. Applying  \eqref{green} yields 
 \begin{alignat*}{2}
\psi(z)& =&&\tfrac{\ii C(\alpha)}{4\pi\left(1-\frac\alpha2\right) }\int_{\partial D (t)}\frac{\overline{\zeta}-\overline{z}}{\vert \zeta-z \vert^\alpha} d\zeta.
\end{alignat*}
Let us now move to the energy $E$. According to \eqref{J} and  the preceding identity  we get
\begin{equation*}
 E(t)=\tfrac{\ii \,C(\alpha)}{8\pi\left(1-\frac\alpha2\right) }\bigintsss_{D(t)}\bigintsss_{\partial D (t)}\frac{\overline{\zeta}-\overline{z}}{\vert \zeta-z \vert^\alpha}d\zeta dA(z).
  \end{equation*}
Using once again  Green's formula \eqref{green}, then we deduce that
the contour integral
  \begin{alignat*}{2}
E(t)=\tfrac{ C(\alpha)}{16\pi\left(1-\frac\alpha2\right)\left(2-\frac\alpha2\right)}\bigintsss_{\partial D (t)}\bigintsss_{\partial D(t)}\frac{(\overline{\zeta}-\overline{z})^2}{\vert \zeta-z \vert^\alpha}dz  d\zeta.
\end{alignat*}
Using \eqref{z}   the last expression becomes 
\begin{align}\label{forE-r}
E(t)&=\tfrac{ C(\alpha)}{16\pi \left(1-\frac\alpha2\right)\left(2-\frac\alpha2\right)}\bigintsss_{0}^{2\pi}\bigintsss_{0}^{2\pi}\frac{\big( \overline{z(t,\eta)}- \overline{z(t,\theta)}\big)^2}{ |z(t,\eta)- z(t,\theta)|^\alpha}\partial_{\theta} z(t,\theta)\partial_{\eta} z(t,\eta)d\eta d\theta.
\end{align}
From the identity 
\begin{align*}
\frac{\big( \overline{z(t,\eta)}- \overline{z(t,\theta)}\big)^2}{ |z(t,\eta)- z(t,\theta)|^\alpha}\partial_{\theta} z(t,\theta)&=-\partial_{\theta} \overline{z(t,\theta)}|z(t,\eta)- z(t,\theta)|^{2-\alpha}\\ &-\frac{1}{\left(1-\frac\alpha2\right)}\big( \overline{z(t,\eta)}- \overline{z(t,\theta)}\big)\partial_{\theta}|z(t,\eta)- z(t,\theta)|^{2-\alpha}\cdot 
\end{align*}
We may write
\begin{align*}
E(t)&=-\tfrac{ C(\alpha)}{16\pi \left(1-\frac\alpha2\right)\left(2-\frac\alpha2\right)}\int_{0}^{2\pi}\int_{0}^{2\pi}\partial_{\theta} \overline{z(t,\theta)}\partial_{\eta} z(t,\eta)|z(t,\eta)- z(t,\theta)|^{2-\alpha}d\eta d\theta\\ 
&-\tfrac{ C(\alpha)}{8\pi \left(1-\frac\alpha2\right)^2\left(2-\frac\alpha2\right)}\int_{0}^{2\pi}\int_{0}^{2\pi}\big( \overline{z(t,\eta)}- \overline{z(t,\theta)}\big)\partial_{\theta}|z(t,\eta)- z(t,\theta)|^{2-\alpha}\partial_{\eta} z(t,\eta)d\eta d\theta.
\end{align*}
Integrating the last term by parts gives
\begin{align*}
E(t)&=-\tfrac{ C(\alpha)}{16\pi \left(1-\frac\alpha2\right)^2}\int_{0}^{2\pi}\int_{0}^{2\pi}|z(t,\eta)- z(t,\theta)|^{2-\alpha}\partial_{\theta} \overline{z(t,\theta)}\partial_{\eta} z(t,\eta)d\eta d\theta
\end{align*}
Since $E(t)$ is a real-valued function then 
\begin{align*}
E(t)&=-\tfrac{ C(\alpha)}{16\pi\left(1-\frac\alpha2\right)^2}\int_{0}^{2\pi}\int_{0}^{2\pi} |z(t,\eta)- z(t,\theta)|^{2-\alpha}\textnormal{Re}\big\{\partial_{\theta} \overline{z(t,\theta)}\partial_{\eta} z(t,\eta)\big\}d\eta d\theta.
\end{align*}
Then, substituting \eqref{z} in the  last expression  concludes the proof of the lemma.
\end{proof}

\subsection{Hamiltonian structure and reversibility}\label{Sec-Hamiltonian structure}

In this section we describe the Hamiltonian structure of the contour  dynamic  equation \eqref{eq}. We emphasize that the Hamiltonian formulation is crucial to approximately reduce the linearized operator to
a triangular form using the  "approximate inverse" approach  that will be developed in Section \ref{sec:Approximate-inverse}. 

\begin{proposition}\label{prop-hamilt}
The following assertions hold true.
\begin{enumerate}
\item The $L^2$-gradient of the angular momentum  $J$ is given by
\begin{align*}
\nabla J(r)=
-\frac12-r(\theta)\cdot
\end{align*}
\item The $L^2$-gradient of the pseudo energy $E$ 
is given by
\begin{align*}
\nabla E(r)=
\tfrac{ C(\alpha)}{4\pi \left(1-\frac\alpha2\right)}\bigintsss_{0}^{2\pi}\frac{1+2r(\eta)-(1+2r(\theta))^{\frac12}\partial_{\eta}\big[(1+2r(\eta))^{\frac12}\sin(\eta-\theta)\big]}{A_r^{\alpha/2}(\theta,\eta)}
d\eta,
\end{align*}
where  $A_r$ is given by \eqref{A}. 
\item The equation \eqref{eq} is Hamiltonian and takes  the form 
\begin{equation*}
{\partial_t} r={\partial_\theta}\nabla H(r) \quad \mbox{with}\qquad H\triangleq E +\Omega J.
\end{equation*}
Moreover, this equation  is reversible with respect to the involution
\begin{equation*}
( {\mathfrak S} r )(\theta) \triangleq  r (- \theta ) \, , 
\end{equation*}
namely the Hamiltonian vector field $ X_{H} = \partial_\theta \nabla H $ satisfies
\begin{equation*}
X_{H} \circ {\mathfrak S} = - {\mathfrak S} \circ X_{H} \, .
\end{equation*}
The reversibility property may be equivalently read
\begin{equation*}
\big(\partial_t r-X_{H}\big) \circ {\mathcal S} = - {\mathcal S} \circ \big(\partial_t r-X_{H}\big)  \, ,
\end{equation*}
with 
$$(\mathcal{S} r)(t,\theta)\triangleq r(-t,-\theta).$$
 Thus, if $r(t,\theta)$ is a solution of \eqref{eq} then $(\mathcal{S} r)(t,\theta)$ is also a solution.
\end{enumerate}
\end{proposition}
\begin{proof} For simplicity we shall  denote along this proof  $ r(t, \theta) = r(\theta)$.

\smallskip

${\bf {(i)}}$  Differentiating  \eqref{forJ} with respect to $r$ in the direction $h$ we get 
\begin{align*}
dJ(r)[h]&=- \frac12\int_{\T}h(\theta)\big(1+2r(\theta)\big)d\theta.
\end{align*}
Then ${\rm {(i)}}$ is immediate.

\smallskip

${\bf {(ii)}}$ Recall that the expression of the energy in polar coordinates is given by
\begin{align}\label{forE}
E&=-\tfrac{ C(\alpha)}{16\pi \left(1-\frac\alpha2\right)^2}\bigintsss_{0}^{2\pi}\bigintsss_{0}^{2\pi}\frac{\partial^2_{\theta\eta}\big[(1+2r(\theta))^{\frac12}(1+2r(\eta))^{\frac12}\cos(\eta-\theta)\big]}{A_r^{\frac\alpha2-1}(\eta,\theta)}d\eta d\theta.
\end{align}
 Differentiating  the last expression  with respect to $r$ in the direction $h$ gives
\begin{align*}
dE(r)[h]&=
-\tfrac{ C(\alpha)}{16\pi \left(1-\frac\alpha2\right)^2}\bigintsss_{0}^{2\pi}\bigintsss_{0}^{2\pi}\frac{\partial^2_{\theta\eta}\big[h(\theta)(1+2r(\theta))^{-\frac12}(1+2r(\eta))^{\frac12}\cos(\eta-\theta)\big]}{A_r^{\frac\alpha2-1}(\eta,\theta)}d\eta d\theta
\\ &\quad-\tfrac{ C(\alpha)}{16\pi \left(1-\frac\alpha2\right)^2}\bigintsss_{0}^{2\pi}\bigintsss_{0}^{2\pi}\frac{\partial^2_{\theta\eta}\big[h(\eta)(1+2r(\theta))^{\frac12}(1+2r(\eta))^{-\frac12}\cos(\eta-\theta)\big]}{A_r^{\frac\alpha2-1}(\eta,\theta)}d\eta d\theta\
\\ &\quad-\tfrac{ C(\alpha)}{16\pi \left(1-\frac\alpha2\right)}\bigintsss_{0}^{2\pi}\bigintsss_{0}^{2\pi}\frac{\partial^2_{\theta\eta}\big[(1+2r(\theta))^{\frac12}(1+2r(\eta))^{\frac12}\cos(\eta-\theta)\big]}{A_r^{\alpha/2}(\theta,\eta)}dA_r[h](\eta,\theta)d\eta d\theta,
\end{align*}
where $A_r$ is given by \eqref{A} and
\begin{equation*}
dA_r[h](\eta,\theta)=2h(\theta)+2h(\eta)- 2\bigg[h(\theta)\frac{\big(1+ 2r(\eta)\big)^{\frac12}}{\big(1+ 2r(\theta)\big)^{\frac12}}+h(\eta)\frac{\big(1+ 2r(\theta)\big)^{\frac12}}{\big(1+ 2r(\eta)\big)^{\frac12}}\bigg]\cos(\eta-\theta)\cdot
\end{equation*}
By symmetry arguments one may write
\begin{align*}
dE(r)[h]&=
-\tfrac{ C(\alpha)}{8\pi \left(1-\frac\alpha2\right)^2}\bigintsss_{0}^{2\pi}\bigintsss_{0}^{2\pi}\frac{\partial^2_{\theta\eta}\big[h(\theta)(1+2r(\theta))^{-\frac12}(1+2r(\eta))^{\frac12}\cos(\eta-\theta)\big]}{A_r^{\frac\alpha2-1}(\eta,\theta)}d\eta d\theta
\\ &\quad-\tfrac{ C(\alpha)}{4\pi \left(1-\frac\alpha2\right)}\bigintsss_{0}^{2\pi}\bigintsss_{0}^{2\pi}h(\theta)\frac{\partial^2_{\theta\eta}\big[(1+2r(\theta))^{\frac12}(1+2r(\eta))^{\frac12}\cos(\eta-\theta)\big]}{A_r^{\alpha/2}(\theta,\eta)}
\\ &\qquad\qquad\qquad\qquad\qquad\qquad\times\Big(1- {\big(1+ 2r(\eta)\big)^{\frac12}}{\big(1+ 2r(\theta)\big)^{-\frac12}}\cos(\eta-\theta)\Big)d\eta d\theta\cdot
\end{align*}
Integrating the first term by parts we get
\begin{align}\label{dE}
dE(r)[h]&=
\tfrac{ C(\alpha)}{8\pi \left(1-\frac\alpha2\right)}\bigintsss_{0}^{2\pi}\bigintsss_{0}^{2\pi}h(\theta)\frac{\partial_{\eta}\big[(1+2r(\eta))^{\frac12}\cos(\eta-\theta)\big]}{(1+2r(\theta))^{\frac12}}\frac{\partial_\theta A_r(\eta,\theta)}{A_r^{\alpha/2}(\theta,\eta)} d\eta d\theta
\\ &\quad-\tfrac{ C(\alpha)}{4\pi \left(1-\frac\alpha2\right)}\bigintsss_{0}^{2\pi}\bigintsss_{0}^{2\pi}h(\theta)\frac{\partial^2_{\theta\eta}\big[(1+2r(\theta))^{\frac12}(1+2r(\eta))^{\frac12}\cos(\eta-\theta)\big]}{A_r^{\alpha/2}(\theta,\eta)}
\notag\\ &\qquad\qquad\qquad\qquad\qquad\times\Big(1- {\big(1+ 2r(\eta)\big)^{\frac12}}{\big(1+ 2r(\theta)\big)^{-\frac12}}\cos(\eta-\theta)\Big)d\eta d\theta\cdot\notag
\end{align}
In view of \eqref{A} one has
\begin{align*}
\partial_\theta A_r(\eta,\theta)&=2\partial_\theta r(\theta)\Big(1-  \big(1+ 2r(\theta)\big)^{-\frac12}{\big(1+ 2r(\eta)\big)^{\frac12}}\cos(\eta-\theta)\Big)\\ &\quad-2\big(1+ 2r(\theta)\big)^{\frac12}{\big(1+ 2r(\eta)\big)^{\frac12}}\sin(\eta-\theta).
\end{align*}
On the other hand, we have the identity 
\begin{align*}
\partial^2_{\theta\eta}\big[(1+2r(\theta))^{\frac12}(1+2r(\eta))^{\frac12}\cos(\eta-\theta)\big]&=\partial_{\theta}r(\theta)(1+2r(\theta))^{-\frac12}\partial_{\eta}\big[(1+2r(\eta))^{\frac12}\cos(\eta-\theta)\big]\\ &\quad+(1+2r(\theta))^{\frac12}\partial_{\eta}\big[(1+2r(\eta))^{\frac12}\sin(\eta-\theta)\big].
\end{align*}
Inserting the two last identities into \eqref{dE} we find 
\begin{align*}
&dE(r)[h]=
\tfrac{ C(\alpha)}{4\pi \left(1-\frac\alpha2\right)}\bigintsss_{0}^{2\pi}\bigintsss_{0}^{2\pi}h(\theta)\frac{1+2r(\eta)-(1+2r(\theta))^{\frac12}\partial_{\eta}\big[(1+2r(\eta))^{\frac12}\sin(\eta-\theta)\big]}{A_r^{\alpha/2}(\theta,\eta)}d\eta d\theta
\cdot
\end{align*}
Consequently,
\begin{align}
\nabla E(r)=
\tfrac{ C(\alpha)}{4\pi \left(1-\frac\alpha2\right)}\bigintsss_{0}^{2\pi}\frac{1+2r(\eta)-(1+2r(\theta))^{\frac12}\partial_{\eta}\big[(1+2r(\eta))^{\frac12}\sin(\eta-\theta)\big]}{A_r^{\alpha/2}(\theta,\eta)}
d\eta.
\end{align}

\smallskip

${\bf {(iii)}}$  Since $H=E +\Omega J.$ then by   {\rm{(i)}} and {\rm{(ii)}} we have
\begin{align*}
\nabla H(r)=-\tfrac\Omega2\big(1+2 r(\theta)\big)+
\tfrac{ C(\alpha)}{4\pi \left(1-\frac\alpha2\right)}\bigintsss_{0}^{2\pi}\frac{1+2r(\eta)-(1+2r(\theta))^{\frac12}\partial_{\eta}\big[(1+2r(\eta))^{\frac12}\sin(\eta-\theta)\big]}{A_r^{\alpha/2}(\theta,\eta)}
d\eta.
\end{align*}
Differentiating the last expression with respect to $\theta$ gives
\begin{align*}
\partial_\theta\nabla H(r)&=-\Omega\,\partial_\theta r(\theta)
-\tfrac{ C(\alpha)}{4\pi \left(1-\frac\alpha2\right)}\bigintsss_{0}^{2\pi}\frac{\partial^2_{\theta\eta}\big[(1+2r(\theta))^{\frac12}(1+2r(\eta))^{\frac12}\sin(\eta-\theta)\big]}{A_r^{\alpha/2}(\theta,\eta)}
d\eta
\\ &\quad-\tfrac{\alpha C(\alpha)}{8\pi \left(1-\frac\alpha2\right)}\bigintsss_{0}^{2\pi}\frac{1+2r(\eta)-(1+2r(\theta))^{\frac12}\partial_{\eta}\big[(1+2r(\eta))^{\frac12}\sin(\eta-\theta)\big]}{A_r(\eta,\theta)^{\frac\alpha2+1}}\partial_\theta A_r(\eta,\theta)
d\eta\cdot
\end{align*}
Straightforward computations allow to get
\begin{align*}
&\Big(1+2r(\eta)-(1+2r(\theta))^{\frac12}\partial_{\eta}\big[(1+2r(\eta))^{\frac12}\sin(\eta-\theta)\big]\Big)\partial_\theta A_r(\eta,\theta)\\ &=-2\,\partial^2_{\theta\eta}\big[(1+2r(\theta))^{\frac12}(1+2r(\eta))^{\frac12}\sin(\eta-\theta)\big]A_r(\eta,\theta)\\ &+\Big(1+2r(\theta)+(1+2r(\eta))^{\frac12}\partial_{\theta}\big[(1+2r(\theta))^{\frac12}\sin(\eta-\theta)\big]\Big)\partial_{\eta} A_r(\eta,\theta).
\end{align*}
This leads to
\begin{align*}
\partial_\theta\nabla H(r)&=-\Omega\partial_\theta r(\theta)+
 \tfrac{(\alpha-1) C(\alpha)}{4\pi \left(1-\frac\alpha2\right)}\bigintsss_{0}^{2\pi}\frac{\,\partial^2_{\theta\eta}\big[(1+2r(\theta))^{\frac12}(1+2r(\eta))^{\frac12}\sin(\eta-\theta)\big]}{A_r^{\alpha/2}(\theta,\eta)}d\eta
\\ &\quad-\tfrac{\alpha C(\alpha)}{8\pi \left(1-\frac\alpha2\right)}\bigintsss_{0}^{2\pi}\frac{1+2r(\theta)+(1+2r(\eta))^{\frac12}\partial_{\theta}\big[(1+2r(\theta))^{\frac12}\sin(\eta-\theta)\big]}{A_r(\eta,\theta)^{\frac\alpha2+1}}\partial_{\eta} A_r(\eta,\theta)d\eta.
\end{align*}
Integrating by parts gives
\begin{align}\label{par-nabla-E}
&\partial_\theta\nabla H(r)=-\Omega\partial_\theta r(\theta)-
 \tfrac{ C(\alpha)}{2\pi}\bigintsss_{0}^{2\pi}\frac{\,\partial^2_{\theta\eta}\big[(1+2r(\theta))^{\frac12}(1+2r(\eta))^{\frac12}\sin(\eta-\theta)\big]}{A_r^{\alpha/2}(\theta,\eta)}d\eta.
\end{align}
Comparing  \eqref{eq}, \eqref{par-nabla-E}  concludes the proof of ${\rm{(iii)}}$.

\smallskip

${\bf {(iv)}}$
In view of \eqref{par-nabla-E} one has
\begin{align*}
({\mathfrak S}\circ X_H)(r)(\theta)&=-\Omega \,\partial_\theta r(-\theta)- \frac{ 1}{2\pi }\bigintsss_0^{2\pi}\frac{\partial^2_{\theta\eta} \big[\big(1+ 2r(-\theta)\big)^{\frac12}\big(1+ 2r(\eta)\big)^{\frac12}\sin(\eta+\theta)\big]}{A_r(-\theta,\eta)^{\alpha/2}} d\eta\cdot\notag
\end{align*}
Using the change of variable $\eta\mapsto-\eta$ in the last integral we get
\begin{align*}
({\mathfrak S}\circ X_H)(r)(-\theta)&=-\Omega \,\partial_\theta ({\mathfrak S} r)(\theta)+\frac{ 1}{2\pi }\bigintsss_{-2\pi}^0\frac{\partial^2_{\theta\eta} \big[\big(1+ 2({\mathfrak S}r)(\theta)\big)^{\frac12}\big(1+ 2({\mathfrak S}r)(\eta)\big)^{\frac12}\sin(\eta-\theta)\big]}{A^{\alpha/2}_{{\mathfrak S} r}(\theta,\eta)} d\eta\cdot
\end{align*}
where we have used the fact that $A_r(-\theta,-\eta)=A_{{\mathfrak S}r}(\theta,\eta)$.
Consequently 
\begin{align}
({\mathfrak S}\circ X_H)(r)(\theta)&= -X_H\circ ({\mathfrak S}r)(\theta).
\end{align}
This concludes the proof of the lemma.

\end{proof}

The symplectic form induced by the hamiltonian equation
\begin{equation}\label{HS}
{\partial_t} r={\partial_\theta}\nabla H(r) \quad \mbox{with}\qquad H\triangleq E +\Omega J.
\end{equation}
 is given by
\begin{equation}\label{sy2form}
\mathcal{W} ( r, h) \triangleq \int_0^{2\pi}  (\partial_\theta^{-1}r )h   d \theta,
\quad
 \partial_\theta^{-1} r= \sum_{j \in \Z \setminus \{0\}}\frac{1}{\ii j} r_j e^{\ii j \theta}, \quad \forall r,h\in L_0^2(\mathbb{T}).  
\end{equation}

The corresponding Poisson tensor is 
$ \mathcal{J}\triangleq  \partial_\theta $, and the Poisson bracket
$$
\{ F, G\} 
=\mathcal{W} (X_F, X_G) 
=\int_0^{2\pi} \nabla F (r) \partial_\theta \nabla G (r) \, d \theta  , 
$$
where $ \nabla F $, $ \nabla G $  denote the $ L^2 $-gradients of 
the functionals $ F, G : L^2_0(\mathbb{T}) \to \mathbb{R}  $.

Note that the Hamiltonian vector field $X_H=\partial_\theta\nabla H$, associated with the Hamiltonian $H$, is determined by 
 \begin{equation}\label{dhw}
 dH(r)[h]=(\nabla H(r),h)_{L^2(\mathbb{T})}=\mathcal{W} ( X_H(r),  \, h\, ),\quad \forall r,h\in L_0^2(\mathbb{T}).
  \end{equation}

\section{Linearization and frequencies structure}
In this section we compute the linear Hamiltonian PDE obtained from  linearizing  \eqref{HS}   at a given small state $r$ close to the equilibrium. At the latter one, we prove that the linearized operator acts as a Fourier multiplier with symbol related to the gamma functions. We shall also explore  some important structures on the equilibrium frequencies related to the  monotonicity,   the asymptotic behavior for large modes and the non-degeneracy/transversality properties. This last point turns out to be  crucial to establish the emergence of linear  quasi-periodic solutions for a massive set of exponents $\alpha\in(0,1).$  
\subsection{Linearized operator}
We intend to discuss the structure of the linearized operator around a general state close to the equilibrium  $r=0$. As a byproduct  we deduce that at the equilibrium state the linear operator  acts as a Fourier multiplier  with explicit eigenvalues related to Gamma function. The first main result reads as follows.
\begin{proposition}\label{lin-eq-r}
The linearized Hamiltonian equation  \eqref{HS} at  a small function  $(t,\theta)\mapsto r $ in the direction $(t,\theta)\mapsto  h( t, \theta ) $ is given by 
\begin{align*}
\partial_t h (t, \theta) & =  \partial_\theta\Big[-\big(\Omega+V_{r,\alpha}(t,\theta)\big)h(t, \theta) +\mathbb{K}_{r,\alpha}h(t,\theta)\Big]
\end{align*}
where $V_{r,\alpha}( t,\theta)  $ is the real  function
\begin{equation}\label{defV}
V_{r,\alpha}(t,\theta)\triangleq \tfrac{ C(\alpha) }{2\pi }{\big(1+ 2r(t,\theta)\big)^{-\frac12}}\bigintsss_{0}^{2\pi} \frac{\partial_{\eta} \big[\big(1+ 2r(t,\eta)\big)^{\frac12}\sin(\eta-\theta)\big]}{A_r^{\alpha/2}(\theta,\eta)} d\eta,
\end{equation}
 the integral operator $\mathbb{K}_{r,\alpha}$ is defined by
\begin{equation}\label{defK}
\mathbb{K}_{r,\alpha}h (t,\theta)\triangleq \tfrac{ C(\alpha) }{2\pi }\bigintsss_{0}^{2\pi} \frac{h(t,\eta)}{A_r^{\alpha/2}(\theta,\eta)} d\eta
\end{equation}
and   $ A_r(\theta, \eta)  $ is given by \eqref{A}.
\end{proposition}
\begin{proof} Throughout  the proof, and for the sake of simple notation, we remove the  time  dependency  from the functions. The computations of  the G\^ateaux derivative of the vector field $\nabla E$, given by Proposition  \ref{prop-hamilt}-${\rm (ii)}$, at the point $r$ in the direction $h$ are straightforward and standard and we shall only sketch the main lines. Notice that the functional  is smooth in a suitable functional setting and therefore its Frechet  differential  can  be recovered from   its G\^ateaux derivative.  From direct computations one gets\begin{align}\label{dnablaE}
 d\nabla E(r)h(\theta)&=\tfrac{ C(\alpha)}{4\pi \left(1-\frac\alpha2\right)}\Bigg[-\frac{h(\theta)}{\big(1+ 2r(t,\theta)\big)^\frac12}\bigintsss_{0}^{2\pi}\frac{\partial_{\eta}\big[(1+2r(\eta))^{\frac12}\sin(\eta-\theta)\big]}{A_r^{\alpha/2}(\theta,\eta)}d\eta
\\ &\qquad +\bigintsss_{0}^{2\pi}\frac{h(\eta)-(1+2r(\theta))^{\frac12}\partial_{\eta}\big[h(\eta)(1+2r(\eta))^{-\frac12}\sin(\eta-\theta)\big]}{A_r^{\alpha/2}(\theta,\eta)},
d\eta\notag\\ &\qquad-\frac{ \alpha}{2}\bigintsss_{0}^{2\pi}\frac{C_r(\eta,\theta)}{A_r(\eta,\theta)^{\frac\alpha2+1}}\notag
\big(h(\eta)B_r(\eta,\theta)+h(\theta)B_r(\theta,\eta)\big)d\eta\Bigg]\notag
 \end{align}
 with
\begin{align}\label{B}
B_r(\eta,\theta)&\triangleq 2-2\big(1+ 2r(t,\eta)\big)^{-\frac12}\big(1+ 2r(t,\theta)\big)^{\frac12}\cos(\eta-\theta)\big),\\
C_r(\eta,\theta)&\triangleq 1+2r(\eta)-(1+2r(\theta))^{\frac12}\partial_{\eta}\big[(1+2r(\eta))^{\frac12}\sin(\eta-\theta)\big]. \label{C}
\end{align}
From \eqref{A} and \eqref{B} one may observe that 
\begin{equation}\label{sum-BB}
\big(1+2r(\eta)\big)B_r(\eta,\theta)+\big(1 + 2 r(\theta)\big)B_r(\theta,\eta)=2A_r(\theta, \eta) \, . 
\end{equation}
Plugging the last identity into the last term of \eqref{dnablaE} yields
\begin{align*}
T_3(\theta) &\triangleq\frac{ \alpha}{2}\bigintsss_{0}^{2\pi}\frac{C_r(\eta,\theta)}{A_r(\eta,\theta)^{\frac\alpha2+1}}\notag
\big(h(\eta)B_r(\eta,\theta)+h(\theta)B_r(\theta,\eta)\big)
d\eta
\\ &= \frac{ \alpha }{2}\bigintsss_{0}^{2\pi}\frac{C_r(\eta,\theta)B_r(\eta,\theta)}{A_r(\eta,\theta)^{\frac\alpha2+1}}\bigg(h(\eta)-h(\theta)\frac{1+ 2r(t,\eta)}{1+ 2r(t,\theta)}\bigg)d\eta+\frac{ \alpha h(\theta)}{1+ 2r(t,\theta)}\bigintsss_{0}^{2\pi}\frac{C_r(\eta,\theta)}{A_r(\eta,\theta)^{\frac\alpha2}}
d\eta.
\end{align*}
In view of  \eqref{A}, \eqref{B}, \eqref{C} and \eqref{sum-BB}  one may easily check that
\begin{align}\label{partialR-new2}
&C_r(\eta,\theta) B_r(\eta,\theta)= 2A_r(\theta,\eta)-\big(1+2r(\eta)\big)^{-\frac12}\big(1+2r(\theta)\big)^{\frac12}\sin (\eta - \theta)\partial_{\eta} A_r(\theta,\eta).
\end{align}
Then, by \eqref{C} and \eqref{partialR-new2} we get
\begin{align*}
T_3(\theta)&= \alpha \bigintsss_{0}^{2\pi}\frac{h(\eta)}{A_r^{\alpha/2}(\theta,\eta)}d\eta
-\frac{\alpha h(\theta)}{\big(1+ 2r(t,\theta)\big)^{\frac12}}\bigintsss_{0}^{2\pi}\frac{\partial_{\eta}\big[(1+2r(\eta))^{\frac12}\sin(\eta-\theta)\big]}{A_r^{\alpha/2}(\theta,\eta)}
d\eta
\\ &\quad-\frac{ \alpha }{2}\bigintsss_{0}^{2\pi}\frac{\big(1+2r(\eta)\big)^{-\frac12}\big(1+2r(\theta)\big)^{\frac12}\sin (\eta - \theta)\partial_{\eta} A_r(\theta,\eta)}{A_r(\eta,\theta)^{\frac\alpha2+1}}\bigg(h(\eta)-h(\theta)\frac{1+ 2r(t,\eta)}{1+ 2r(t,\theta)}\bigg)
d\eta.
\end{align*}
Integrating the last term by parts implies
\begin{align*}
T_3(\theta)&=\alpha\bigintsss_{0}^{2\pi}\frac{h(\eta)}{A_r^{\alpha/2}(\theta,\eta)}d\eta+ \frac{(1-\alpha) h(\theta)}{\big(1+ 2r(t,\theta)\big)^{\frac12}}\bigintsss_{0}^{2\pi}\frac{\partial_{\eta}\big[(1+2r(\eta))^{\frac12}\sin(\eta-\theta)\big]}{A_r^{\alpha/2}(\theta,\eta)}
d\eta
\\ &\quad-\bigintsss_{0}^{2\pi}\frac{\big(1+2r(\theta)\big)^{\frac12}\partial_{\eta}\big[(h(\eta)\big(1+2r(\eta)\big)^{-\frac12}\sin (\eta - \theta)\big]}{A_r^{\alpha/2}(\theta,\eta)}
d\eta.
\end{align*}
Inserting the last formula into \eqref{dnablaE} gives
\begin{align*}
&d\nabla E(r)[h](\theta)=
\tfrac{ C(\alpha)}{2\pi }\bigintsss_{0}^{2\pi}\frac{h(\eta)}{A_r^{\alpha/2}(\theta,\eta)}
d\eta
-\frac{ C(\alpha)\,h(\theta)}{2\pi\big(1+ 2r(t,\theta)\big)^{\frac12}}\bigintsss_{0}^{2\pi}\frac{\partial_{\eta}\big[(1+2r(\eta))^{\frac12}\sin(\eta-\theta)\big]}{A_r^{\alpha/2}(\theta,\eta)}
d\eta.
\end{align*}
Combining this expression with Proposition \ref{prop-hamilt} gives the desired result. 
\end{proof}

Next, we shall move to the special case $r=0$ corresponding to the equilibrium state. We shall particularly show that the linear operator  acts as a Fourier multiplier allowing to get explicit form for the spectrum. More precisely, we establish the following result.
\begin{proposition}\label{linear-eq} The linearized  equation of
 \eqref{eq}  at the equilibrium state $r=0$  writes
\begin{align}\label{linearized-op}
\partial_t h & =\partial_\theta\, \mathrm{L}(\alpha)\, h=\partial_\theta \nabla H_{\mathrm{L}}(h) ,
\end{align}
where $H_{\mathrm{L}}(h)$ is the quadratic Hamiltonian
\begin{equation}\label{defQH}
H_{\mathrm{L}}(h)\triangleq\frac12\big(\mathrm{L}(\alpha)h,h\big)_{L^2}\, ,
\end{equation}
the self-adjoint operator  $\mathrm{L}(\alpha)$ is given by
\begin{equation}\label{defL}
\mathrm{L}(\alpha)\triangleq-\big(\Omega+V_{0,\alpha}\big)+\mathbb{K}_{0,\alpha},
\end{equation}
 the constant  $V_{0,\alpha}$ is given by
\begin{equation}\label{defTheta}
V_{0,\alpha}\triangleq \frac{ \Gamma(\frac12-\frac\alpha2)}{2\sqrt{\pi} \,\Gamma\left(1-\frac\alpha2\right)}\frac{\Gamma(1+\frac\alpha2)}{\Gamma\left(2-\frac\alpha2\right)},
\end{equation}
and the integral operator $\mathbb{K}_{0,\alpha}$ is diagonal:\, for any $h(t,\theta)=\displaystyle\sum_{j\in\mathbb{Z}}h_{j}(t)e^{\ii j\theta}$ ,
\begin{align}\label{k0}
\mathbb{K}_{0,\alpha}\,h(t,\theta)&= \tfrac{ \Gamma(\frac12-\frac\alpha2)}{2\sqrt{\pi} \,\Gamma\left(1-\frac\alpha2\right)}\sum_{j\in\mathbb{Z}}\frac{ \Gamma(j+\frac{\alpha}{2})}{\Gamma(1+j-\frac{\alpha}{2})}h_{j}(t)e^{\ii j\theta}\, \cdot
\end{align}
Moreover, the  reversible solutions to the linear equation \eqref{linearized-op} are given by
\begin{equation}\label{rev-sol}
h(t,\theta)=\sum_{j\in\mathbb{N}}h_j\,\cos\big(j\theta-\Omega_j(\alpha)t \big) ,\quad\textnormal{with}\quad h_j\in\mathbb{R},
\end{equation}
with 
\begin{equation}\label{omega}
\Omega_j(\alpha)\triangleq j\Omega+\frac{j\Gamma(\frac12-\frac\alpha2)}{2\sqrt{\pi}\,\Gamma(1-\frac{\alpha}{2})}\bigg(\frac{  \Gamma(1+\frac\alpha2)}{ \Gamma(2-\frac\alpha2) }-\frac{  \Gamma(j+\frac\alpha2)}{ \Gamma(1+j-\frac\alpha2)} \bigg)\cdot
\end{equation}

\end{proposition}
\begin{proof}
According to  Proposition \ref{lin-eq-r},  the linearized equation at  $r=0$ is given by
\begin{align}\label{linearized-op0}
\partial_t h  & =  \partial_\theta\Big[-\big(\Omega+V_{0,\alpha}\big)h+\mathbb{K}_{0,\alpha}(\alpha)h\Big],
\end{align}
where 
\begin{equation}\label{defV0}
V_{0,\alpha}= \frac{ C_\alpha }{2^{1+\alpha}\pi }\bigintsss_{0}^{2\pi} \frac{\cos(\eta-\theta)}{|\sin\big(\frac{\eta-\theta}{2}\big)|^\alpha} d\eta\quad\textnormal{and}\quad \mathbb{K}_{0,\alpha}h (t,\theta)= \frac{ C_\alpha }{2^{1+\alpha}\pi }\bigintsss_{0}^{2\pi} \frac{h(t,\eta)}{|\sin\big(\frac{\eta-\theta}{2}\big)|^\alpha} d\eta.
\end{equation}
We shall now recall the following identity, see for instance  \cite[p.449]{W}.
\begin{equation}\label{int}
\bigintsss_0^{\pi}\sin^x(\eta)e^{\ii y\eta}d\eta=\frac{\pi e^{\ii y\frac{\pi}{2}}\Gamma(x+1)}{2^x\Gamma(1+\frac{x+y}{2})\Gamma(1+\frac{x-y}{2})}\,\quad \forall x>-1.
\end{equation}
Apply this formula with $x = -\alpha$ and $y = 2$ yields,
\begin{align*}
V_{0,\alpha}&= -\frac{ C(\alpha) \Gamma(1-\alpha)}{ \Gamma\left(2-\frac\alpha2\right) \Gamma(-\alpha/2)}\cdot  
\end{align*}
Recalling $C_\alpha=\frac{\Gamma(\alpha/2)}{2^{1-\alpha}\Gamma(\frac{2-\alpha}{2})}$  and using the  identity $z\Gamma(z)=\Gamma(z+1)$  we obtain, for any $\alpha\in(0,1)$,
\begin{align*}
V_{0,\alpha} 
&=\frac{\Gamma(\frac{\alpha}{2}) \Gamma(1-\alpha)}{2^{1-\alpha}\Gamma(1-\frac{\alpha}{2})\Gamma\left(2-\frac\alpha2\right) \Gamma(-\frac\alpha2)} \\ 
&=\frac{\Gamma(1+\frac{\alpha}{2}) \Gamma(1-\alpha)}{2^{1-\alpha}\Gamma^2(1-\frac{\alpha}{2})\Gamma\left(2-\frac\alpha2\right) }\cdot
\end{align*}
Therefore  from  Legendre duplication formula (see \eqref{leg-dup}),
\begin{equation}\label{leg}
2^\alpha\sqrt{\pi}\,\Gamma(1-\alpha)=\Gamma\left(\tfrac12-\tfrac\alpha2\right)\Gamma\left(1-\tfrac\alpha2\right),
\end{equation}
 we obtain the expression of  $V_{0,\alpha}$ in \eqref{defTheta}.
Next, we shall study the action of the integral  operator $\mathbb{K}_{0,\alpha}$  on the scalar function $h$  that we expand in Fourier series as
\begin{equation}\label{rho6}
h(t,\theta)=\sum_{j\in\mathbb{Z}\setminus \{0\}}h_{j}(t){\bf{e}}_j(\theta)\quad\hbox{with}\quad {\bf{e}}_j(\theta)=e^{\ii j\theta}.
\end{equation}
By linearity, it suffices to evaluate the operator on the exponential basis. From \eqref{defV0} and \eqref{rho6} we write in  view of change of variables 
\begin{align*}
\mathbb{K}_{0,\alpha}\,{\bf{e}}_j(\theta)
&=\frac{ C(\alpha) }{2^{1+\alpha}\pi }{\bf{e}}_j(\theta)\bigintsss_{0}^{2\pi} \frac{e^{\ii j\eta}}{|\sin\big(\frac{\eta}{2}\big)|^{\alpha}} d\eta.
\end{align*}
Then, applying the formula \eqref{int} with $x=-\alpha$ and $y=2j$ we get
\begin{align*}
\mathbb{K}_{0,\alpha}\,{\bf{e}}_j(\theta) 
 &=  \frac{ (-1)^j\,C(\alpha)\Gamma(1-\alpha)}{\Gamma(1+j-\frac{\alpha}{2})\Gamma(1-j-\frac{\alpha}{2})}{\bf{e}}_j(\theta)
\\ & =\frac{ \Gamma(\frac\alpha2)}{2^{1-\alpha}\Gamma(1-\frac{\alpha}{2})}  \frac{ (-1)^j\Gamma(1-\alpha)}{\Gamma(1+j-\frac{\alpha}{2})\Gamma(1-j-\frac{\alpha}{2})}{\bf{e}}_j(\theta)\cdot\notag
\end{align*}
Using the identity, see \eqref{Pocc1},
\begin{equation}\label{Pocc16}
\frac{\Gamma(x+j)}{\Gamma(x)}=(-1)^j\frac{\Gamma(1-x)}{\Gamma(1-x-j)},
\end{equation}
we obtain in view of \eqref{leg}
\begin{align*}
\mathbb{K}_{0,\alpha}\,{\bf{e}}_j(\theta)&=\frac{ \Gamma(1-\alpha)}{2^{1-\alpha}\Gamma^2(1-\frac{\alpha}{2})}\frac{ \Gamma(j+\frac{\alpha}{2})}{\Gamma(1+j-\frac{\alpha}{2})}{\bf{e}}_j(\theta)\notag
\\
& =\frac{ \Gamma(\frac12-\frac\alpha2)}{2\sqrt{\pi}\,\Gamma(1-\frac{\alpha}{2})}\frac{ \Gamma(j+\frac{\alpha}{2})}{\Gamma(1+j-\frac{\alpha}{2})}{\bf{e}}_j(\theta)\cdot
\end{align*}
Finally, we obtain 
\begin{align}\label{expK0}
\mathbb{K}_{0,\alpha}h(t,\theta)& =\frac{ \Gamma(\frac12-\frac\alpha2)}{2\sqrt{\pi}\,\Gamma(1-\frac{\alpha}{2})}\sum_{j\in\mathbb{Z}}h_{j}(t)\frac{ \Gamma(j+\frac{\alpha}{2})}{\Gamma(1+j-\frac{\alpha}{2})}{\bf{e}}_j(\theta)\cdot
\end{align}
Inserting \eqref{defTheta} and  \eqref{expK0} into \eqref{linearized-op0} gives the Fourier expansion of  the solutions in the form
\begin{equation}\label{solrho}
h(t,\theta)=\sum_{j\in\mathbb{Z}}h_{j}(0)e^{-\ii (\Omega_j(\alpha)t-j\theta) },\, h_j(0)\in\mathbb{C}.
\end{equation}
where the frequencies $\Omega_j(\alpha)$ are defined in \eqref{omega}.
In view of \eqref{Pocc16} one has
\begin{equation}\label{2id}
\frac{(-1)^j\Gamma(1-\frac\alpha2)}{\Gamma\big(1-j-\frac\alpha2\big)}=\frac{\Gamma\big(j+\frac\alpha2\big)}{\Gamma\big(\frac\alpha2\big)}\quad\textnormal{and}\quad\frac{(-1)^j\Gamma(\frac\alpha2)}{\Gamma\big(-j+\frac\alpha2\big)}=\frac{\Gamma\big(j+1-\frac\alpha2\big)}{\Gamma\big(1-\frac\alpha2\big)}\cdot
\end{equation}
It follows that
\begin{equation}\label{sym-freq}
\forall\, j\in\Z,\quad \frac{  \Gamma(-j+\frac\alpha2)}{ \Gamma(1-j-\frac\alpha2)} =\frac{  \Gamma(j+\frac\alpha2)}{ \Gamma(1+j-\frac\alpha2)} \cdot
\end{equation}
Combining the last identity with the expression of $\Omega_j(\alpha)$, given by \eqref{omega}, we conclude that
\begin{align}\label{Freq-symm}
\Omega_{-j}(\alpha)=-\Omega_j(\alpha).
\end{align}
Therefore, every real-valued reversible solution to   \eqref{linearized-op0} has the form
$$
h(t,\theta)=\sum_{j\in\mathbb{N}}h_{j}\cos\big(\Omega_j(\alpha)t-j\theta \big),\quad\textnormal{with}\quad h_{j}\in\mathbb{R}.
$$
This concludes the proof of the lemma.
\end{proof}
\begin{remark}
In Fourier expansion $h(t,\theta)=\displaystyle\sum_{j\in\mathbb{Z}\setminus \{0\}}h_{j}(t)e^{\ii n\theta}$,  the quadratic Hamiltonian $H_{\rm L}$ in \eqref{defQH} writes 
\begin{equation}\label{defQH1}
H_{\mathrm{L}}(h)=-\pi \sum_{j\in\mathbb{Z}\setminus\{0\}}\tfrac{\Omega_j(\alpha)}{j}|h_j|^2=-2\pi \sum_{j\in\mathbb{N}}\tfrac{\Omega_j(\alpha)}{j}|h_j|^2, 
\end{equation}
and the self-adjoint operator $ \mathrm{L}(\alpha)$ reads
\begin{equation}
 \mathrm{L}(\alpha)h (t,\theta)=-\sum_{j\in\mathbb{Z}\setminus\{0\}} \tfrac{\Omega_j(\alpha)}{j} h_j(t) e^{\ii j\theta} .
\end{equation}
 The symplectic form writes
 \begin{equation}\label{sy2form-fourier2}
\mathcal{W} = \tfrac12\sum_{j \in \Z \setminus \{0\}}\frac{1}{\ii j}dr_j\wedge dr_{-j}=\sum_{j \geqslant 1}\frac{1}{\ii j}dr_j\wedge dr_{-j},
\end{equation}
the Hamiltonian vector field $X_H$ is
\begin{equation*}
[X_H(r)]_j={\ii j }\big(\partial_{ r_{-j}}H\big)(r),\quad \forall j\neq 0
\end{equation*}
 and the Poisson bracket,
\begin{equation*}
\{ F, G\} =\mathcal{W} (X_F, X_G) 
=- \sum_{j \in \Z \setminus \{0\}}\ii j\big(\partial_{ r_{-j}}F\big)(r)\big(\partial_{ r_{j}}G\big)(r).
\end{equation*}

\end{remark}

\subsection{Structure  of the linear frequencies}
The main  task in this section is to investigate some important structures of the equilibrium frequencies. In the first part we shall be concerned with their monotonicity and explore some  asymptotic behavior for large modes. However we shall discuss in the second part the non-degeneracy of these frequencies  through the so-called R\"ussmann conditions. This is  crucial  in the measure  of the final Cantor set giving rise to quasi-periodic solutions for the linear/nonlinear problem. Actually, this set appears as a perturbation of the Cantor set constructed from the equilibrium eigenvalues and therefore perturbative arguments based on their degeneracy  are  very useful and will be implemented in Section \ref{Section 6.2}.
\subsubsection{Monotonicity and asymptotic behavior}
In what follows, we intend to establish some basic properties related to the monotonicity and the asymptotic behavior for large modes of the spectrum  of the linearized operator at the equilibrium state. Notice that their explicit values are detailed  \mbox{in \eqref{omega}.} Our result reads as follows.
\begin{lemma}\label{lem-asym}
Let $\alpha\in [0,1)$ and $\Omega>0$. Then the following holds true.
\begin{enumerate}
\item  For any   $j\in\mathbb{Z}$,   
$
\Omega_{-j}(\alpha)=-\Omega_j(\alpha). 
$
\item  The sequence $\big(\tfrac{\Omega_j(\alpha)}{j}\big)_{j\in \mathbb{N}^*}$ is positive and   strictly increasing.
\item For any $j\in \mathbb{N}^*$, we have the expansion
\begin{align*}
{\Omega_j(\alpha)}&=
\left(\Omega+V_{0,\alpha}\right)j-W_{0,\alpha}\,j\,\mathtt{W}(j,\alpha)\\&=
 \left(\Omega+V_{0,\alpha}\right)j-W_{0,\alpha}\, j^{\alpha}+O\left({j^{\alpha-2}}\right).
\end{align*}
where $V_{0,\alpha}$ is given by \eqref{defTheta}
 and
\begin{equation*}\label{defTheta1}
W_{0,\alpha}\triangleq \frac{1}{2\sqrt\pi}\frac{ \Gamma(\frac12-\frac\alpha2)}{\Gamma(1-\frac{\alpha}{2})}, \quad \mathtt{W}(j,\alpha)\triangleq\frac{\Gamma\big(j+\frac{\alpha}{2}\big)}{\Gamma\big(j+1-\frac{\alpha}{2}\big)}\cdot
\end{equation*}

\item For all $j\in\mathbb{Z}$, we have 
\begin{equation*}
 |\Omega_j(\alpha)|\geqslant \Omega |j|.
\end{equation*}
\item Given $\overline\alpha\in(0,1)$, there exists $c>0$ such that 
\begin{equation*}
\forall \alpha\in[0,\overline\alpha],\quad \forall j,j'\in \mathbb{Z},\quad |\Omega_j(\alpha)\pm\Omega_{j'}(\alpha)|\geqslant c |j\pm j'|.
\end{equation*}
\item Given $\overline\alpha\in(0,1)$ and ${q}_0\in\N$, there exists $C>0$ such that
\begin{equation*}
\forall j,j'\in \mathbb{Z},\quad \max_{q\in\llbracket 0,{q}_0\rrbracket} \sup_{\alpha\in [0,\overline\alpha]}\big|\partial_\alpha^q\big(\Omega_j(\alpha)-\Omega_{j'}(\alpha)\big)\big|\leqslant C |j-j'|.
\end{equation*}
\end{enumerate}
\end{lemma} 
\begin{proof} 
${\bf(i)}$ It is proved in \eqref{Freq-symm}.

\smallskip

${\bf(ii)}$ Using the identity \eqref{Pocc1} together with \eqref{omega} allows to  find the alternative formula
\begin{equation}\label{disper00}
\Omega_j(\alpha)=j\Omega+ jV_{0,\alpha}\bigg(1-\frac{\big(1+\frac{\alpha}{2}\big)_{j-1}}{\big(2-\frac{\alpha}{2}\big)_{j-1}}\bigg),
\end{equation}
where 
 $(x)_j$ denotes Pokhhammer's symbol introduced in \eqref{Poch}.  Now, because  $x\mapsto (x)_{j-1}$ is strictly  increasing in the set $\RR_+$, we  conclude  that for all $\alpha\in[0,1)$ one has
\begin{equation}\label{bound-po}
0<\frac{\big(1+\frac{\alpha}{2}\big)_{j-1}}{\big(2-\frac{\alpha}{2}\big)_{j-1}}<1.
\end{equation} 
Since $x\mapsto \Gamma(x)$ is strictly positive, continuous and has a  minimum on  $\RR_+^*$ then we can easily show, from \eqref{defTheta}, that there exists $c_1>0$ such that 
\begin{equation}\label{inf-V0}
\inf_{\alpha\in(0,\overline\alpha)}V_{0,\alpha}\triangleq c_1>0.
\end{equation}
 Combining \eqref{disper00}, \eqref{bound-po} and \eqref{inf-V0} we obtain 
\begin{equation} \label{lowb}
\forall j\in \mathbb{N}^*,\quad\Omega_j(\alpha)> j\Omega>0.
\end{equation}
To prove that $j\mapsto \Omega_j(\alpha)/j$ is strictly increasing it suffices, according to \eqref{disper00}, to check that 
 the sequence $j\mapsto u_j= \frac{\big(1+\frac{\alpha}{2}\big)_{j-1}}{\big(2-\frac{\alpha}{2}\big)_{j-1}}$ is strictly decreasing. This follows from the obvious fact that 
$$
\forall \alpha\in[0,1),\quad  \frac{u_{j+1}}{u_j}=\frac{j+\frac\alpha2}{j+1-\frac\alpha2}<1.
$$
This ends the proof of ${\rm(ii)}$.

\smallskip

${\bf(iii)}$ The asymptotic behavior  follows from  \eqref{omega} and Lemma \ref{Stirling-formula}.

\smallskip

${\bf(iv)}$  It is an immediate consequence of ${\rm(i)}$ and  \eqref{lowb}.

\smallskip

${\bf(v)}$  By the oddness of $j\mapsto \Omega_j(\lambda)$ it is enough  to establish the estimate for $j,j'\in\mathbb{N}^*.$ We shall first focus on the estimate of the difference, that is,
\begin{equation}\label{diff}
\inf_{\alpha\in[0,\overline\alpha]}|\Omega_j(\alpha)-\Omega_{j^\prime}(\alpha)|\geqslant c|j-j^\prime|,
\end{equation}
for some constant $c>0$. The result is obvious when $j=j^\prime$, then we shall restrict  the discussion to $j\neq j^\prime$.  Without  loss of generality we can assume that   $j> j'\geqslant 1$.
According to {\rm (iii)} one may writes
\begin{align}\label{difference}
{\Omega_j(\alpha)}-{\Omega_{j'}(\alpha)}&=\frac{\Omega_j(\alpha)}{j}(j-j')-j'\Big(\frac{\Omega_{j'}(\alpha)}{j'}-\frac{\Omega_j(\alpha)}{j}\Big)\nonumber
\\&= \frac{\Omega_j(\alpha)}{j}(j-j')-j'\, W_{0,\alpha} \big(\mathtt{W}(j',\alpha)-\mathtt{W}(j,\alpha)\big)\cdot
\end{align}
 Then, by {\rm (iv)} we obtain
$$
{\Omega_j(\alpha)}-{\Omega_{j'}(\alpha)}\geqslant \Omega(j-j')
-j'\,\sup_{\alpha\in[0,\overline\alpha]} \Big(W_{0,\alpha} \big|\mathtt{W}(j',\alpha)-\mathtt{W}(j,\alpha)\big|\Big).
$$
  Using  Taylor formula combined with  the estimate in Lemma \ref{Stirling-formula}-{\rm (i)} we get, for \mbox{any $\alpha\in[0,\overline\alpha],$}
  \begin{align}\label{estm:dif}
 \big| \mathtt{W}(j',\alpha)-\mathtt{W}(j,\alpha)\big| &\leqslant C\Big|\int_{j'}^j\frac{dx}{x^{2-\alpha}}\Big|\nonumber \\ &\leqslant C|j^{\alpha-1}-(j')^{\alpha-1}|\nonumber \\ &\leqslant C \frac{|j-j'|}{j' j^{1-\alpha}}\cdot
  \end{align}
  It follows that 
\begin{align*}
\Omega_j(\alpha)-\Omega_{j^\prime}(\alpha)&\geqslant \Big(\Omega-\frac{C}{j^{1-\alpha}}\Big)(j-j')
\cdot
\end{align*}
 Therefore there exists $j_0$ independent of $\alpha\in[0,\overline\alpha]$ such that if   $j> j_0$ then we get \eqref{diff}. It remains to verify this property for $j'< j\leqslant  j_0$. Using the one-to-one property  of $j\mapsto \Omega_j(\alpha)$ combined with the continuity of $\alpha\in[0,\overline\alpha]\mapsto \Omega_j(\alpha)-\Omega_{j^\prime}(\alpha)$ we get for any $j> j^\prime\geqslant 1$
$$
\inf_{\alpha\in[0,\overline\alpha]}\big(\Omega_j(\alpha)-\Omega_{j^\prime}(\alpha)\big)\triangleq c_{jj^\prime}>0.
$$
Consequently 
$$
\inf_{j\neq j^\prime\in\llbracket 0,j_{0}\rrbracket\\
\atop \alpha\in[0,\overline\alpha]}|\Omega_j(\alpha)-\Omega_{j^\prime}(\alpha)|=\inf_{j\neq j^\prime\in\llbracket 0,j_{0}\rrbracket}c_{jj^\prime}>0
$$
which implies   \eqref{diff}.

\smallskip

			Let us now move to the estimate of $\Omega_{j}(\alpha)+\Omega_{j_{0}}(\alpha)$ for $j,j'\in\mathbb{N}^*.$ Since both quantities are positive then using the point \rm{(iv)} yields
				$$
				\forall\lambda\in[{0},\overline\alpha],\quad|\Omega_{j}(\alpha)+\Omega_{j_{0}}(\alpha)|=\Omega_{j}(\alpha)+\Omega_{j'}(\alpha)\geqslant \Omega(j+j')\geqslant c(j+j').
				$$
				This completes the proof of the desired estimate.

\smallskip

${\bf (vi)}$ From {\rm (iii)} we may write
\begin{align}\label{difference2}
{\Omega_j(\alpha)}-{\Omega_{j'}(\alpha)}&=\big(\Omega+V_{0,\alpha}\big)(j-j')- W_{0,\alpha} \big(j\mathtt{W}(j,\alpha)-j'\mathtt{W}(j',\alpha)\big).
\end{align}
 Let  ${q}_0\in\N$ and $q\in\llbracket 0,{q}_0\rrbracket$. Differentiating $q$ times the identity \eqref{difference} in $\alpha$ we find 
\begin{align}\label{difference-diff}
\partial_\alpha^q\big({\Omega_j(\alpha)}-{\Omega_{j'}(\alpha)}\big)&=\partial_\alpha^q\big(\Omega+V_{0,\alpha}\big)(j-j')- \big(\partial_\alpha^q W_{0,\alpha}\big) \big(j\mathtt{W}(j,\alpha)-j'\mathtt{W}(j',\alpha)\big)\nonumber\\ &\quad -  W_{0,\alpha} \,\partial_\alpha^q\big(j\mathtt{W}(j,\alpha)-j'\mathtt{W}(j',\alpha)\big)\cdot
\end{align}
By the  mean value theorem  combined with  the estimate in Lemma \ref{Stirling-formula}-{\rm (i)} we get,
\begin{equation}\label{ineq:fj}
 \sup_{\alpha\in [0,\overline\alpha]}\big|\partial_\alpha^q\big(j\mathtt{W}(j,\alpha)-j'\mathtt{W}(j',\alpha)\big)\big|\leqslant C {|j-j'|}.
\end{equation}
Combining \eqref{difference2} and  \eqref{ineq:fj} gives the estimate in {\rm (vi)}.

\end{proof}
\subsubsection{Non-degeneracy and transversality}
Next we shall explore the local structure of a given finite  set  of eigenvalues associated to the equilibrium state. In particular, we intend to establish their  non-degeneracy which is an essential property in measuring some suitable connected sets with Diophantine constraints. We first start with giving the following definition.
\begin{definition}\label{def-deg} 
Let  $N\in\mathbb{N}^*$. A function $f \triangleq  (f_1, \ldots , f_N ) : [\alpha_1,\alpha_2] \to \mathbb{R}^N$ is called non-degenerate if, for any vector $c \triangleq  (c_1,\ldots,c_N) \in  \mathbb{R}^N \backslash \{0\}$, the scalar function $f \cdot c = f_1c_1 + \cdots+ f_dc_N$ is not identically zero on the whole interval $[\alpha_1,\alpha_2]$. This means that the  curve of $f$ cannot be contained in an hyperplane.
\end{definition}
The first main goal is to prove the following result.
\begin{lemma}  \label{non-degeneracy}
Let $\overline\alpha\in(0,1)$, $N\in\mathbb{N}^*$ and  $1 \leqslant j_1 < j_2 < \cdots < j_N$ some fixed integers . Then 
the functions
\begin{align*}
[0,\overline\alpha]&\to \mathbb{R}^N,\qquad \\
\alpha &\mapsto\big(\Omega_{j_1} (\alpha),\ldots,\Omega_{j_N} (\alpha)\big),\\
[0,\overline\alpha]&\to \mathbb{R}^{N+1}\\
\alpha &\mapsto\big(\Omega+V_{0,\alpha},\Omega_{j_1} (\alpha),\ldots,\Omega_{j_N} (\alpha)\big)
\end{align*}
are non-degenerate in the sense of the Definition $\ref{def-deg},$  where $\Omega_j(\alpha)$ is given by \eqref{omega} and $V_{0,\alpha}$ is defined in \eqref{defTheta}. 
 \end{lemma}

\begin{proof}

According to Definition \ref{def-deg} one has to prove that, for all $c\in  \mathbb{R}^N \backslash \{0\}$, the function
\begin{align}\label{fn}
\alpha &\mapsto
g_1(\alpha)\triangleq \sum_{k=1}^N c_k  j_k\bigg[\frac{1}{2\sqrt{\pi}}\frac{\Gamma(\frac12-\frac\alpha2)}{\Gamma(1-\frac{\alpha}{2})}\bigg(\frac{  \Gamma(1+\frac\alpha2)}{ \Gamma(2-\frac\alpha2) }-\frac{  \Gamma(j_k+\frac\alpha2)}{ \Gamma(1+j_k-\frac\alpha2)} \bigg)+\Omega\bigg]
\end{align}
is not identically zero on the interval $[0,\overline\alpha]$. 
Suppose, by contradiction, that there exists $c\triangleq  (c_1,\ldots,c_N)\in \mathbb{R}^N\setminus \{0\}$ such that
\begin{equation}\label{comb0}
 \sum_{k=1}^N c_k  j_k\bigg[\frac{1}{2\sqrt{\pi}}\frac{\Gamma(\frac12-\frac\alpha2)}{\Gamma(1-\frac{\alpha}{2})}\bigg(\frac{  \Gamma(1+\frac\alpha2)}{ \Gamma(2-\frac\alpha2) }-\frac{  \Gamma(j_k+\frac\alpha2)}{ \Gamma(1+j_k-\frac\alpha2)} \bigg)+\,\Omega\bigg]
=0,\quad \forall \alpha\in [0,\overline\alpha].
\end{equation}
Since the gamma function has no real zeros but simple poles located at $-\mathbb{N}$ then, the function $\alpha\mapsto \frac{\Gamma(\frac12-\frac\alpha2)}{\Gamma(1-\frac{\alpha}{2})}$ can be extended to holomorphic function on  $\mathbb{C}\backslash \mathcal{Q}$ with
$$
\mathcal{Q}\triangleq \big\{\alpha=2n+1;\;  n \in \mathbb{N}\big\}
$$
and for all $j\in \mathbb{N}$,  
 the function $\alpha\mapsto \frac{  \Gamma(j+\frac\alpha2)}{ \Gamma(1+j-\frac\alpha2)} $ can be extended to holomorphic function on  $\mathbb{C}\backslash \mathcal{P}_j$ with 
$$
\mathcal{P}_j\triangleq \big\{\alpha=-2(n+j);\; n \in \mathbb{N}\big\}.
$$
Moreover, 
\begin{equation}
\mathcal{P}_{j+1}\subset \mathcal{P}_{j}\subset\cdots\subset \mathcal{P}_{1},\quad\textnormal{and}\quad \mathcal{P}_{j}\backslash \mathcal{P}_{j+1}=\{-2j\}.
\end{equation}

Therefore, the function $g_1(\alpha)$, defined in \eqref{fn}, can be extended to holomorphic function on $\mathbb{C}\backslash \mathcal{P}$, where
\begin{equation}\label{def:setP}
\mathcal{P}\triangleq \mathcal{P}_1\cup \mathcal{Q}=\big\{\alpha=-2(n+1);\; n \in \mathbb{N}\big\}\cup \big\{\alpha=2n+1;\;n \in \mathbb{N}\big\}.
\end{equation}
Thus, according to \eqref{comb0}, there exists $c\triangleq  (c_1,\ldots,c_N)\in \mathbb{R}^N\setminus \{0\}$ such that
\begin{equation}\label{comb}
 \sum_{k=1}^N c_k  j_k\bigg[\frac{1}{2\sqrt{\pi}}\frac{\Gamma(\frac12-\frac\alpha2)}{\Gamma(1-\frac{\alpha}{2})}\bigg(\frac{  \Gamma(1+\frac\alpha2)}{ \Gamma(2-\frac\alpha2) }-\frac{  \Gamma(j_k+\frac\alpha2)}{ \Gamma(1+j_k-\frac\alpha2)} \bigg)+\,\Omega\bigg]
=0,\quad \forall \alpha\in \mathbb{C}\backslash \mathcal{P}.
\end{equation}
Since $\Omega>0$ and the function $\alpha\mapsto 1/{\Gamma(1-\frac{\alpha}{2})}$ vanishes at $2$ then
 substituting $\alpha=2$ in  \eqref{comb} gives,
\begin{equation*}
\sum_{k=1}^N c_k  j_k 
=0.
\end{equation*}
It follows that
\begin{equation}\label{comb1}
g_1(\alpha)\triangleq
 \sum_{k=1}^N c_k  j_k\bigg[\bigg(\frac{  \Gamma(1+\frac\alpha2)}{ \Gamma(2-\frac\alpha2) }-\frac{  \Gamma(j_k+\frac\alpha2)}{ \Gamma(1+j_k-\frac\alpha2)} \bigg)\bigg]
=0,\quad \forall \alpha\in \mathbb{C}\backslash \mathcal{P}.
\end{equation}

For given $j\in\mathbb{N}$, the function $\alpha\mapsto 1/{ \Gamma(1+j-\frac\alpha2)} $  has zero set
$$
\mathcal{Z}_j\triangleq \big\{\alpha=2(n+j+1);\; n \in \mathbb{N}\big\}
$$
Moreover, one has 
\begin{equation}\label{inclusion-zeros}
\mathcal{Z}_{j}\subset \mathcal{Z}_{j-1},\quad\textnormal{and}\quad \mathcal{Z}_{j-1}\backslash \mathcal{Z}_{j}=\{2j\}.
\end{equation}
 Thus, substituting $\alpha=2j_N$ in  \eqref{comb1} we get
 \begin{align*}
g_1(2j_N)&= -c_{N} j_N\Gamma(2j_N)=0.
\end{align*}
It follows that
$$
c_{N} =0.
$$
Reproducing the same arguments for $c_{{N-1}},c_{{N-2}},\ldots,c_{{1}}$ respectively we get
$$
c_{{N-1}}=c_{{N-2}}=\cdots=c_{{1}}=0.
$$
This gives the contradiction. 
In order to conclude the proof of the lemma we have to prove that, for any  $N$, for any  ${1 \leqslant} n_1 < n_2 < \cdots < n_N$ the function $\alpha\in [0,\overline\alpha] \mapsto \big(\Omega+V_{0,\alpha},\Omega_{j_1} (\alpha),\ldots,\Omega_{j_N} (\alpha)\big)\in\mathbb{R}^{N+1}$  is non-degenerate according to Definition \ref{def-deg}.
Suppose, by contradiction, that there exists $c = (c_0,c_1,\ldots,c_N) \in \mathbb{R}^{N+1} \backslash\{0\}$ such that for all $\alpha\in \mathbb{C}\backslash \mathcal{P}$,
\begin{align*}
 c_0\Big(\Omega+\frac{2^{-1}}{\sqrt{\pi}}\frac{\Gamma(\frac12-\frac\alpha2)}{\Gamma(1-\frac{\alpha}{2})}\frac{  \Gamma(1+\frac\alpha2)}{ \Gamma(2-\frac\alpha2) }\Big) +\sum_{k=1}^N c_k  j_k\bigg[\frac{2^{-1}}{\sqrt{\pi}}\frac{\Gamma(\frac12-\frac\alpha2)}{\Gamma(1-\frac{\alpha}{2})}\bigg(\frac{  \Gamma(1+\frac\alpha2)}{ \Gamma(2-\frac\alpha2) }-\frac{  \Gamma(j_k+\frac\alpha2)}{ \Gamma(1+j_k-\frac\alpha2)} \bigg)+\,\Omega\bigg]=0.
\end{align*}
Substituting $\alpha=2$ in the last equation  gives
\begin{equation}\label{comb33}
c_0+\sum_{k=1}^N c_k  n_k =0.
\end{equation}
Therefore for all $\alpha\in \mathbb{C}\backslash \mathcal{P}$, where $\mathcal{P}$ is given by \eqref{def:setP}, one has
\begin{equation}\label{comb34}
 g_2(\alpha)\triangleq c_0\frac{  \Gamma(1+\frac\alpha2)}{ \Gamma(2-\frac\alpha2) }+\sum_{k=1}^N c_k  j_k\bigg(\frac{  \Gamma(1+\frac\alpha2)}{ \Gamma(2-\frac\alpha2) }-\frac{  \Gamma(j_k+\frac\alpha2)}{ \Gamma(1+j_k-\frac\alpha2)} \bigg)=0.
\end{equation}
It follows that
 \begin{align*}
g_2(2j_N)&=-c_{N} {j_N}\Gamma(2j_N).
\end{align*}
Then by \eqref{comb1} we conclude that
\begin{equation}\label{comb44cd}
c_{N} =0.
\end{equation}
Arguing as above we deduce we get
\begin{equation}\label{comb44}
c_{{N-1}}=c_{{N-2}}=\cdots=c_{{1}}=0.
\end{equation}
Finally, inserting  \eqref{comb44cd} and \eqref{comb44}  into \eqref{comb33} we find
\begin{equation}\label{comb55}
c_0=0.
\end{equation}
Then the vector $c$ is vanishing and this contradicts the assumption.\end{proof}

\begin{proposition}\label{lemma transversality}
Let $\overline\alpha\in(0,1)$.There exist $q_0\in\mathbb{N} $, $\rho_0>0$ such that, for all  $\alpha\in [0,\overline\alpha]$, the following assertions hold true.
\begin{enumerate}
\item For any $l\in \mathbb{Z}^d\setminus\{0\}$ we have
\begin{align*}
\max_{k\leqslant q_0}|\partial_\alpha^k\{{\omega}_{\textnormal{Eq}}(\alpha)\cdot l\}|&\geqslant \rho_0\langle l\rangle.
\end{align*}
 \item For any $ (j,l)\in ({\mathbb{N}}\times\Z^d)\setminus\{(0,0)\},$ we have
\begin{align*}
\max_{k\leqslant q_0}|\partial_\alpha^k\{{\omega}_{\textnormal{Eq}}(\alpha)\cdot l+jV_{0,\alpha}\}|&\geqslant \rho_0\langle l\rangle.
\end{align*}
\item For any $ (j,l)\in ({\mathbb{N}}\setminus{\mathbb{S}}_0)\times \Z^d,$ we have
\begin{align*}
\max_{k\leqslant q_0}|\partial_\alpha^k\{{\omega}_{\textnormal{Eq}}(\alpha)\cdot {l}+\Omega_j(\alpha)\}|&\geqslant \rho_0\langle l\rangle.
\end{align*}
\item For any $l\in\mathbb{Z}$, $j,j'\in\mathbb{N}\setminus {\mathbb{S}}_0$, with  $(j,l)\neq( j^\prime,0)$, we have 
\begin{align*}
\max_{k\leqslant q_0}|\partial_\alpha^k\{{\omega}_{\textnormal{Eq}}(\alpha)\cdot {l}+\Omega_j(\alpha)-\Omega_{j'}(\alpha)\}|\geqslant \rho_0\langle l\rangle
\end{align*}

\end{enumerate}
where ${\omega}_{\textnormal{Eq}}(\alpha)$ and   $\Omega_n(\alpha)$ are defined in \eqref{tan-nor} and \eqref{omega}.
\end{proposition}
\begin{proof}
${\bf{(i)}}$
Suppose, by contradiction, that for all $q_0\in \mathbb{N}$, $\rho_0>0$ there exist $\alpha\in [0,\overline\alpha]$ and $l\in\Z^d\backslash \{0\}$ such that
$$
\max_{k\leqslant q_0}|\partial_\alpha^k\{{\omega}_{\textnormal{Eq}}(\alpha)\cdot {l}\}|< \rho_0\langle l\rangle. 
$$
In particular, for all $n\in\N$, if $\rho_0=\frac{1}{n+1}$ and $q_0=n$ then there exist  $\alpha_n\in [0,\overline\alpha]$ and $l_n\in\Z^d\backslash \{0\}$ such that 
\begin{equation}\label{alphan}
\max_{k\leqslant n}\Big|\partial_\alpha^k\Big\{{\omega}_{\textnormal{Eq}}(\alpha_n)\cdot \tfrac{l_n}{\langle l_n\rangle}\Big\}\Big|< \tfrac{1}{n+1}. 
\end{equation}
The sequences $(\alpha_n)_n\subset [0,\overline\alpha]$ and $(c_n)_n\triangleq \big(\frac{l_n}{\langle l_n\rangle}\big)_n\subset \mathbb{R}^d\backslash\{0\}$ are bounded. 
By compactness and up to an extraction we may assume that 
	$$\lim_{n\to\infty}\tfrac{l_{n}}{\langle l_n\rangle}=\widetilde{c}\neq 0\quad\hbox{and}\quad \lim_{n\to\infty}\alpha_{n}=\widetilde{\alpha}.
				$$
 Passing to the limit in \eqref{alphan} for $n\to\infty$ we deduce that 
 $$\forall k \in \N,\quad \partial_\alpha^k\{{\omega}_{\textnormal{Eq}}(\widetilde\alpha)\cdot \widetilde{c}\}=0\quad{\rm with}\quad \widetilde{c}\neq 0.$$  We conclude that the real analytic  function $\alpha  \to {\omega}_{\textnormal{Eq}}(\alpha)\cdot \widetilde{c}\,$ is identically zero. This is in contradiction with Lemma \ref{non-degeneracy}.
 
 \smallskip

${\bf{(ii)}}$ We shall first check the result for  the case $l=0, j\in\mathbb{N}^*$. Obviously one has from \eqref{inf-V0}, 				\begin{align*}
						\quad\inf_{\alpha\in[0,\overline\alpha]}\max_{k\leqslant  q_{0}}\big|\partial_{\alpha}^{k}\big( jV_{0,\alpha}(\alpha)\big)\big|&\geqslant \inf_{\alpha\in[0,\overline\alpha]}\big| V_{0,\alpha}(\alpha)\big|\geqslant\rho_{0}\langle l\rangle,
						\end{align*}
for some $\rho_0>0.$ Next, we shall consider the case  $j\in\mathbb{N}$, $l\in\mathbb{Z}^d\setminus\{0\}$.
By the triangle inequality combined with the the boundedness of  ${\omega}_{\textnormal{Eq}}$ and the bound  \eqref{inf-V0} we find
$$|{\omega}_{\textnormal{Eq}}(\alpha)\cdot l+jV_{0,\alpha}|\geqslant|j||V_{0,\alpha}|-|{\omega}_{\textnormal{Eq}}(\alpha)\cdot l|\geqslant c|j|-C|l|\geqslant |l|$$
provided that  $|j|\geqslant C_{0}|l|$ for some $C_{0}>0.$ Thus, we shall restrict the proof to  indices  $j$ and $l$ with
\begin{equation}\label{parameter condition 10}
 |j|\le C_{0}|l|, \quad j\in\mathbb{N}, \quad l\in\mathbb{Z}^d\setminus\{0\}.
\end{equation}
Arguing by contradiction as in the previous case, we may assume the existence of sequences $l_{n}\in \Z^d\backslash\{0\}$, $j_n \in  {\mathbb{N}}$ satisfying \eqref{parameter condition 10} and $\alpha_{n}\in[0,\overline\alpha]$ such that 
$$\max_{k\leqslant n}\left|\partial_{\alpha}^{k}\left({\omega}_{\textnormal{Eq}}(\alpha_{n})\cdot\tfrac{l_{n}}{\langle l_{n}\rangle}+\tfrac{{j_{n}}V_{0,\alpha_{n}}}{\langle l_{n}\rangle}\right)\right|<\tfrac{1}{1+n}$$ 
and therefore
\begin{equation}\label{Rossemann 01}
\forall k\in\mathbb{N},\quad\forall n\geqslant k,\quad\left|\partial_{\alpha}^{k}\left({\omega}_{\textnormal{Eq}}(\alpha_{n})\cdot\tfrac{l_{n}}{\langle l_{n}\rangle}+\tfrac{{j_{n}}V_{0,\alpha_{n}}}{\langle l_{n}\rangle}\right)\right|<\tfrac{1}{1+n}.
\end{equation}
Since sequences $\{\alpha_n\}$, $\{b_n\}\triangleq \big\{\frac{j_n}{\langle l_n\rangle}\big\}$  and $\{c_n\}\triangleq \big\{\frac{l_n}{\langle l_n\rangle}\big\}$ are bounded, then up to an extraction we can assume
$$
\lim_{n\to\infty}\alpha_{n}=\widetilde{\alpha},\quad \lim_{n\to\infty}b_{n}=\widetilde{b}\neq 0\quad\hbox{and}\quad \lim_{n\to\infty}c_{n}=\widetilde{c}\neq 0.
				$$
Hence, letting  $n\rightarrow+\infty$ in \eqref{Rossemann 01} and  using $\alpha\mapsto V_{0,\alpha}$ is smooth   we obtain

$$\forall k\in\mathbb{N},\quad\partial_{\alpha}^{k}\left({\omega}_{\textnormal{Eq}}({\alpha})\cdot\widetilde{ c}+{\widetilde{b}}\,V_{0,{\alpha}}\right)_{|\alpha=\widetilde{\alpha}}=0.$$
Thus, the real analytic  function $\alpha\mapsto {\omega}_{\textnormal{Eq}}({\alpha})\cdot\widetilde{ c}+{\widetilde{b}}\,V_{0,{\alpha}}$ with $(\widetilde{ c},{\widetilde{b}})\neq (0,0)$ is identically zero
and this contradicts Lemma  \ref{non-degeneracy}.

 \smallskip

${\bf{(iii)}}$ Consider $(l,j)\in\mathbb{Z}^{d }\times (\mathbb{N}\setminus\mathbb{S}_0)$. Then applying  the  triangle inequality and Lemma \ref{lem-asym}-${\rm{(iv)}}$, we get
$$|{\omega}_{\textnormal{Eq}}(\alpha)\cdot l+\Omega_{j}(\alpha)|\geqslant|\Omega_{j}(\alpha)|-|{\omega}_{\textnormal{Eq}}(\alpha)\cdot l|\geqslant \Omega j-C|l|\geqslant \langle l\rangle$$
provided that  $|j|\geqslant C_{0}|l|$ for some $C_{0}>0.$ So we shall restrict the proof to integers  $j$ with 
\begin{equation}\label{parameter condition 1}
0\leqslant j< C_{0} \langle l\rangle,\quad j\in\mathbb{N}\setminus\mathbb{S}_0\quad\hbox{and}\quad l\in\mathbb{Z}^d\backslash\{0\}.
\end{equation}
Arguing  by contradiction, one can check  that for all $n\in\mathbb{N}$, there exists $l_{n}\in \Z^d\backslash\{0\}, j_n \notin  {\mathbb{S}}_0$ and $\alpha_{n}\in[0,\overline\alpha]$ such that 
$$\max_{k\leqslant n}\left|\partial_{\alpha}^{k}\left({\omega}_{\textnormal{Eq}}(\alpha)\cdot\tfrac{l_{n}}{\langle l_{n}\rangle}+\tfrac{\Omega_{j_{n}}(\alpha)}{\langle l_{n}\rangle}\right)_{|_{\alpha=\alpha_n}}\right|<\frac{1}{1+n}$$ 
and therefore
\begin{equation}\label{Rossemann 1}
\forall k\in\mathbb{N},\quad\forall n\geqslant k,\quad\left|\partial_{\alpha}^{k}\left({\omega}_{\textnormal{Eq}}(\alpha)\cdot\tfrac{l_{n}}{\langle l_{n}\rangle}+\tfrac{\Omega_{j_{n}}(\alpha)}{\langle l_{n}\rangle}\right)_{|\alpha={\alpha}_n}\right|<\frac{1}{1+n}\cdot
\end{equation}
Since sequences $\{\alpha_n\}$  and $\{c_n\}\triangleq \big\{\frac{l_n}{\langle l_n\rangle}\big\}$ are bounded, then up to an extraction we can assume
$$
\lim_{n\to\infty}\alpha_{n}=\widetilde{\alpha}\quad\hbox{and}\quad \lim_{n\to\infty}c_{n}=\widetilde{c}\neq 0.
				$$
Now we shall distinguish two cases.

\smallskip

$\bullet$ {\it Case $1$}: $(l_{n})_{n}$ is bounded. Then from \eqref{parameter condition 1} and up to an extraction the sequences $(l_n)$ and $(j_n)$ are stationary. So we can assume that for any $n\in\N$, we have $l_n=\widetilde{l}$ and $j_n=\widetilde{j}$, with $\widetilde{l}\in  \Z^d\backslash\{0\}$ and $\widetilde{j}\in \mathbb{N}.$ 
Hence, taking the limit as $n\rightarrow+\infty$ in \eqref{Rossemann 1} yields
$$\forall k\in\mathbb{N},\quad\partial_{\alpha}^{k}\left({\omega}_{\textnormal{Eq}}({\alpha})\cdot\widetilde{l}+\Omega_{\widetilde{j}}({\alpha})\right)_{|_{\alpha=\widetilde{\alpha}}}=0.$$
Thus, the real analytic function $\alpha\mapsto\omega_{\textnormal{Eq}}(\alpha)\cdot\widetilde{l}+\Omega_{\widetilde{j}}(\alpha)$ with $(\widetilde{l},1)\neq (0,0)$ is identically zero which contradicts
 Lemma \ref{non-degeneracy}.

\smallskip

$\bullet$ {\it Case $2$}:  $( l_{n})$ is unbounded. Up to a subsequence, we can assume that $ \displaystyle\lim_{n\to\infty} | l_{n}|=\infty$ and 
$\displaystyle \lim_{n\to\infty}\frac{ l_{n}}{\langle  l_{n}\rangle} =\widetilde{c}\in \RR^d\backslash\{0\}.$
We shall distinguish two subclasses depending on whether the sequence $(j_n)_n$ is bounded or not. When it is bounded, then up to an extraction we may assume that this sequence of integers  is stationary.  Then, taking the limit $n\rightarrow+\infty$ in \eqref{Rossemann 1}, we find
$$\forall k\in\mathbb{N},\quad \partial_{\alpha}^{k}{\omega}_{\textnormal{Eq}}({\alpha})_{|_{\alpha=\widetilde{\alpha}}}\cdot\widetilde{c}=0.$$
As before we conclude that the real analytic function $\alpha\mapsto{\omega}_{\textnormal{Eq}}(\alpha)\cdot\widetilde{c}$, with $\widetilde{c}\neq 0$,   is identically zero , which is a contradiction with the Lemma  \ref{non-degeneracy}.\\
It remains to explore the case where $(j_{n})_{n}$ is unbounded. Then up to an extraction we can assume that   $\displaystyle \lim_{n\to\infty} j_{n}=\infty$. According to Lemma \ref{lem-asym}-{\rm (iii)} we find
\begin{equation}\label{quo1}
\tfrac{\Omega_{j_{n}}(\alpha)}{\langle  l_{n}\rangle}=\tfrac{j_{n}}{\langle  l_{n}\rangle}\left(\Omega+V_{0,\alpha}- W_{0,\alpha}\mathtt{W}(j_n,\alpha)\right).
\end{equation}
Moreover, by \eqref{parameter condition 1}, the sequence $\left(\frac{j_{n}}{\langle  l_{n}\rangle}\right)_{n}$ is bounded, so up to a subsequence, we can assume that it converges to $\widetilde{d}.$ Therefore, differentiating then taking the limit in \eqref{quo1} we obtain 
$$\lim_{n\to+\infty}\tfrac{\partial_{\alpha}^q\Omega_{j_{n}}(\alpha_{n})}{\langle  l_{n}\rangle}=\partial_{\alpha}^q\big(\widetilde{d}\big(\Omega+V_{0,{\alpha}}\big)\big)_{|_{\alpha=\widetilde{\alpha}}},
$$
where we have used in the last identity the estimate in Lemma \ref{Stirling-formula}-{\rm (i)}.
Hence, taking the limit $j\rightarrow+\infty$ in \eqref{Rossemann 1} we get :
$$\forall k\in\mathbb{N},\quad \partial_{\alpha}^{k}\left({\omega}_{\textnormal{Eq}}({\alpha})\cdot\widetilde{c}-\widetilde{d}\big(\Omega+V_{0,{\alpha}}\big)\right)_{|_{\alpha=\widetilde{\alpha}}}=0.$$
Thus, the real analytic function $\alpha\mapsto{\omega}_{\textnormal{Eq}}(\alpha)\cdot\widetilde{c}-\widetilde{d}\big(\Omega+V_{0,{\alpha}}\big)$  is identically zero. This is  a contradiction with the Lemma  \ref{non-degeneracy} because  $\widetilde{c}\neq 0$.

\smallskip

\textbf{(iv)} Consider $ l\in\mathbb{Z}^{d }, j,j^\prime\in\mathbb{N}\setminus\mathbb{S}_0$  with $(l,j)\neq(0,j^\prime).$ Then applying the  triangle inequality combined with Lemma \ref{lem-asym}-${\rm{(v)}}$, we infer that
$$|{\omega}_{\textnormal{Eq}}(\alpha)\cdot  l+\Omega_{j}(\alpha)\pm\Omega_{j'}(\alpha)|\geqslant|\Omega_{j}(\alpha)\pm\Omega_{j'}(\alpha)|-|{\omega}_{\textnormal{Eq}}(\alpha)\cdot  l|\geqslant c|j\pm j'|-C|l|\geqslant \langle l\rangle$$
provided $|j-j'|\geqslant C_{0}| l|$ for some $C_{0}>0.$ In this case the desired estimate is trivial. So we shall restrict the proof to integers  such that
\begin{equation}\label{parameter condition 2}
					|j\pm j'|< c_{0}\langle l\rangle,\quad  l\in\mathbb{Z}^{d }\backslash\{0\},\quad  j,j^\prime\in\mathbb{N}\setminus\mathbb{S}_0.
					\end{equation}
We will argue by contradiction as in the previous cases.  We assume that for all $n\in\mathbb{N}$, there exists $( l_n,j_n)\neq(0,j_n^\prime)\in \Z^{d+1}$ satisfying \eqref{parameter condition 2} and $\alpha_{n}\in[0,\overline\alpha]$ such that 
$$
\max_{k\leqslant n}\left|\partial_{\alpha}^{k}\left({\omega}_{\textnormal{Eq}}(\alpha)\cdot\tfrac{ l_{n}}{\langle  l_{n}\rangle}+\tfrac{\Omega_{j_{n}}(\alpha)\pm\Omega_{j'_{n}}(\alpha)}{\langle  l_{n}\rangle}\right)_{|_{\alpha=\alpha_n}}\right|<\frac{1}{1+n}$$ 
and therefore
\begin{equation}\label{Rossemann 2}
\forall k\in\mathbb{N},\quad\forall n\geqslant k,\quad\left|\partial_{\alpha}^{k}\left({\omega}_{\textnormal{Eq}}(\alpha)\cdot\tfrac{ l_{n}}{\langle  l_{n}\rangle}+\tfrac{\Omega_{j_{n}}(\alpha)\pm\Omega_{j'_{n}}(\alpha)}{\langle  l_{n}\rangle}\right)_{|_{\alpha=\alpha_n}}\right|<\frac{1}{1+n}\cdot
\end{equation}
Since the sequences $\left\{\frac{ l_{n}}{\langle  l_{n}\rangle}\right\}_{n}$ and $\{\alpha_{n}\}_{n}$ are bounded, then  by compactness  we can assume that ,  $\displaystyle \lim_{n\to\infty}\frac{ l_{n}}{\langle  l_{n}\rangle}=\widetilde{c}\neq 0$ and $\displaystyle \lim_{n\to\infty}\alpha_{n}=\widetilde{\alpha}.$ As before we distinguish two cases :\\
$\bullet$ {\it Case $1$}: $( l_{n})_{n}$ is bounded. We shall only focus on the most delicate case associated to the difference $\Omega_{j_n}-\Omega_{j^\prime_n}$. Up to an extraction we may assume  that this sequence of integers  is stationary, that is, $ l_{n}=\widetilde l.$ Now looking at \eqref{parameter condition 2} we have two sub-cases depending on whether or not the sequences  $(j_{n})_{n}$ and $(j'_{n})_{n}$ are bounded. 

\smallskip

$\bullet$ Sub-case \ding{172} : $(j_{n})_{n}$ and $(j'_{n})_{n}$ are bounded. Up to an extraction we can assume that they are stationary, that is, $j_n=\widetilde{j}, j_n^\prime=\widetilde{j}^\prime$. In addition from the assumption we have also $(\widetilde l, \widetilde{j})\neq(0,\widetilde{j}^\prime)$ and $\widetilde{j},\widetilde{j}^\prime\notin\mathbb{S}_0$.
Hence taking the limit $n\rightarrow+\infty$ in \eqref{Rossemann 2}, we get
$$\forall k\in\mathbb{N},\quad\partial_{\alpha}^{k}\left({\omega}_{\textnormal{Eq}}({\alpha})\cdot\widetilde{ l}+\Omega_{\widetilde{j}}({\alpha})-\Omega_{\widetilde{j}^\prime}({\alpha})\right)_{|_{\alpha=\widetilde{\alpha}}}=0.$$
Thus,  the real analytic  function $\alpha\mapsto\Omega(\alpha)\cdot\widetilde{ l}+\Omega_{\widetilde{j}}(\alpha)-\Omega_{\widetilde{j}^\prime}(\alpha)$  is identically zero.  If $\widetilde{j}=\widetilde{j^\prime}$ then this  contradicts Lemma \ref{non-degeneracy} since $\widetilde{l}\neq 0.$ However in the case  $\widetilde{j}\neq \widetilde{j^\prime}\in\mathbb{N}\setminus\mathbb{S}_0$ this still contradicts this lemma applied with the vector frequency $(\omega_{\textnormal{Eq}},\Omega_{\widetilde{j}},\Omega_{\widetilde{j'}})$ instead of $\omega_{\textnormal{Eq}}.$

\smallskip

$\bullet$ Sub-case \ding{173} :   $(j_n)_{n}$ and $(j'_{n})_{n}$ are unbounded. Then up to an extraction we can assume that  $\displaystyle \lim_{n\to\infty}j_{n}= \lim_{m\to\infty}j'_{n}=\infty$. 
Without loss of generality, we can assume that for a given $n$ we have $j_n\geqslant j_n^\prime$. Then according to Lemma \ref{lem-asym}-{\rm (iii)} we get the splitting
\begin{align}\label{split}
\tfrac{\partial_{\alpha}^{k}\left(\Omega_{j_{n}}(\alpha)-\Omega_{j'_{n}}(\alpha)\right)}{\langle  l_{n}\rangle}
&=\partial_{\alpha}^{k}\big(\Omega+V_{0,{\alpha}}\big)\tfrac{j_n-j_n^\prime}{\langle  l_n\rangle}-\partial_{\alpha}^{k}\big(W_{0,{\alpha}}\mathtt{W}(j_n,\alpha)\big)\tfrac{j_n-j_n^\prime}{\langle  l_n\rangle}\nonumber\\ &\quad-\tfrac{j'_n}{\langle  l_n\rangle}\partial_{\alpha}^{k}\big[W_{0,\alpha}\big(\mathtt{W}(j'_n,\alpha)-\mathtt{W}(j_n,\alpha)\big)\big].
\end{align}
Using  Taylor formula combined with  Lemma \ref{Stirling-formula}-{\rm (i)}, in a similar way to \eqref{estm:dif},  give for \mbox{any $\alpha\in[0,\overline\alpha],$}
  \begin{align}\label{estm:dif22}
j_n'  \big| \partial_{\alpha}^{k}\mathtt{W}(j'_n,\alpha)- \partial_{\alpha}^{k}\mathtt{W}(j_n,\alpha)\big| &\leqslant C\Big|j_n' \int_{j'_n}^{j_n}\frac{dx}{x^{2-\alpha-\epsilon}}\Big|\nonumber\\ &\leqslant C {|j_n-j_n'|}{j_n^{\overline\alpha+\epsilon-1}}\cdot
  \end{align}
Using once again \eqref{parameter condition 2},   up to an extraction, we assume that $\displaystyle\lim_{n\to\infty}\frac{j'_{n}-j_{n}}{\langle l_{n}\rangle}=\widetilde{d}$. Therefore, combining \eqref{split}, \eqref{estm:dif22} and Lemma \ref{Stirling-formula}-{\rm (i)}, we conclude that
$$\lim_{n\to\infty}\partial_{\alpha}^{k}\left(\frac{\Omega_{j_{n}}(\alpha)-\Omega_{j'_{n}}(\alpha)}{\langle  l_{n}\rangle}\right)_{|_{\alpha=\alpha_n}}=\widetilde{d}\,\partial_{\alpha}^{k}\big(\Omega+V_{0,{\alpha}}\big)_{|_{\alpha=\widetilde\alpha}}.$$
So taking the limit $n\rightarrow+\infty$ in \eqref{Rossemann 2}, yields
$$\forall k\in\mathbb{N},\quad\partial_{\alpha}^{k}\left({\omega}_{\textnormal{Eq}}({\alpha})\cdot\widetilde{c}+\widetilde{d}\,\big(\Omega+V_{0,{\alpha}}\big)\right)_{|_{\alpha=\widetilde{\alpha}}}=0.$$
Thus, the real analytic function $\alpha\mapsto{\omega}_{\textnormal{Eq}}({\alpha})\cdot\widetilde{c}+\widetilde{d}\,\big(\Omega+V_{0,\alpha}\big)$ with $(\widetilde{c},\widetilde{d})\neq (0,0)$ is identically zero 
 which is in a contradiction with Lemma \ref{non-degeneracy}.

\smallskip

$\bullet$ {\it Case} $2$: $( l_{n})_{n}$ is unbounded. Up to an extraction  we can assume that $\displaystyle \lim_{n\to\infty}|l_{n}|=\infty.$\\
				We shall distinguish three sub-cases.\\
				$\bullet$ Sub-case \ding{172}. The sequences  $(j_{n})_{n}$ and $(j'_{n})_{n}$ are bounded. Up to an extraction  they will  converge and then taking the limit in \eqref{Rossemann 2} yields,
				$$
				\forall k\in\mathbb{N},\quad\partial_{\alpha}^{k}{\omega}_{\textnormal{Eq}}(\bar{\alpha})\cdot\widetilde{c}=0.
				$$
				which leads to a contradiction as before. \\
				$\bullet$ Sub-case \ding{173}. The sequences  $(j_{n})_{n}$ and $(j'_{n})_{n}$ are both unbounded. This is similar to  the sub-case \ding{173} of the case 1.\\
				$\bullet$ Sub-case \ding{174}. The sequence $(j_{n})_{n}$ is unbounded and $(j'_{n})_{n}$ is bounded (the symmetric case is similar).  Without loss of generality we can assume that $\displaystyle \lim_{n\to\infty}j_n=\infty$ and $  j_{n}^\prime=\widetilde{j}.$ By \eqref{parameter condition 2} and  up to an extraction one gets  $\displaystyle \lim_{n\to\infty}\frac{j_{n}\pm j'_{n}}{| l_{n}|}=\widetilde{d}.$ Using \eqref{split}
				combined with  \eqref{estm:dif22} and Lemma \ref{Stirling-formula}-{\rm (i)} in order to get for any $k\in\mathbb{N},$
				\begin{align*}
					\lim_{n\to\infty}\langle  l_{n}\rangle^{-1}
					\partial_\alpha^k\Big(\Omega_{j_n}(\alpha)\pm\Omega_{j_{n}^\prime}(\alpha)-(j_n\pm j^\prime_{n})\big(\Omega+V_{0,\alpha}\big)\Big)_{|\alpha=\alpha_n}&=\\
					\lim_{n\to\infty}
					\partial_\alpha^k\left(\tfrac{(j_n\pm j^\prime_n)}{\langle  l_{n}\rangle}W_{0,{\alpha}}\mathtt{W}(j_n,\alpha)\pm \tfrac{j^\prime_{n}}{\langle  l_{n}\rangle}W_{0,{\alpha}}\Big(\mathtt{W}(j_n,\alpha)-\mathtt{W}(j'_n,\alpha)\right)_{|\alpha=\alpha_n}&=0.
				\end{align*}
				Hence, taking the limit  in \eqref{Rossemann 2} implies 
				$$\forall k\in\mathbb{N},\quad \partial_{\alpha}^{k}\left({\omega}_{\textnormal{Eq}}({\alpha})\cdot\widetilde{c}+\widetilde{d}\big(\Omega+V_{0,\alpha}\big)\right)_{\alpha=\widetilde\alpha}=0.$$
				Thus, the real analytic function $\alpha\mapsto{\omega}_{\textnormal{Eq}}(\alpha)\cdot\widetilde{c}+\widetilde{d}(\Omega+V_{0,\alpha})$ is identically zero with $(\widetilde{c},\widetilde{d})\neq0$ which contradicts Lemma \ref{non-degeneracy}. This completes the proof of the lemma.

\end{proof}

\subsection{Linear quasi-periodic solutions}
Note that all the solutions \eqref{rev-sol} for \eqref{linearized-op} are periodic, quasi-periodic or almost periodic in
time, with linear frequencies of oscillations $\Omega_{j}(\alpha)$, depending
on the irrationality properties of the frequencies
$\Omega_{j}(\alpha)$ and how many normal mode amplitudes $h_j(0)$ are
not zero. We shall prove that if $h_j(0) = 0$ for any index $j$ except at a finite set $\mathbb{S}$ (tangential sites) then the linear solutions \eqref{rev-sol} are
quasi-periodic in time.

\begin{lemma}\label{lemma-Line-mes}
Let $\Omega>0, \overline\alpha\in(0,1)$,  $\mathbb{S}\subset{\mathbb{N}\backslash\{0\}}$, with $\#\mathbb{S}=d\geqslant1.$
There exists  a Cantor-like set $\mathtt{G}_{0}\subset(0,\overline\alpha)$ with  full Lebesgue  measure 
such that for any $\alpha\in\mathtt{G}_{0}$, the linearized vortex patch equation \eqref{linearized-op} admits  time quasi-periodic solutions 
 with a non-resonant  frequency vector ${\omega}_{\textnormal{Eq}}(\alpha)\triangleq (\Omega_{j}(\alpha))_{j\in\mathbb{S}}\in\mathbb{R}^{d}$ satisfying \eqref{non-res}  and taking  the form
\begin{equation}\label{q-r-l}
h(t,\theta)=\sum_{j\in\mathbb{S}}{h}_{j}\cos\big(j\theta-\Omega_{j}(\alpha)t\big)
\end{equation}
where 
 $\Omega_{j}(\alpha)$ are the equilibrium frequencies defined in \eqref{omega}.
\end{lemma}
\begin{proof}
We consider the sets
\begin{equation}\label{set-U0se3}
\mathtt{C}_{0}\triangleq\left\lbrace\alpha\in(0,\overline\alpha);\; {\omega}_{\textnormal{Eq}}(\alpha) \cdot l\neq 0, \; \forall l\in \mathbb{Z}^d\setminus\{0\}\right\rbrace,
\end{equation}
\begin{equation}\label{set-U1se3}
 (0,\overline\alpha)\backslash\mathtt{C}_{0}\subset\bigcup_{l\in\mathbb{Z}^{d}\backslash\{0\}}\mathtt{G}_{l}\quad\mbox{ with }\quad\mathtt{G}_{l}\triangleq\left\lbrace\alpha\in(0,\overline\alpha);\; |{\omega}_{\textnormal{Eq}}(\alpha) \cdot l|\leqslant\frac{\kappa}{\langle l\rangle^{\tau_{1}}}\right\rbrace.
\end{equation}
Applying Lemma \ref{Piralt}  and Proposition \ref{lemma transversality}, we find 
$$
|\mathtt{G}_{l} |\lesssim\kappa^{\frac{1}{q_0}}\langle l\rangle^{-1-\frac{\tau_{1}+1}{q_0}}.
$$
It follows from \eqref{set-U1se3} that
\begin{align*}
\big|(0,\overline\alpha)\backslash\mathtt{C}_{0}\big|&\leqslant \sum_{l\in\mathbb{Z}^{d}\backslash\{0\}}\big|\mathtt{G}_{l}\big|\\ &\lesssim \sum_{l\in\mathbb{Z}^{d}}\kappa^{\frac{1}{q_0}}\langle l\rangle^{-1-\frac{\tau_{1}+1}{q_0}}.
\end{align*}
Then assuming $\tau_1>(d-1)q-1$
we deduce that
\begin{equation}\label{set-U4sec3}
\big|(0,\overline\alpha)\backslash\mathtt{C}_{0}\big|\lesssim  \kappa^{\frac{1}{q_0}}
\end{equation}
and thus
\begin{equation}\label{set-U5sec3}
\forall \kappa\in(0,1),\quad  \overline\alpha\geqslant \big|\mathtt{C}_{0}\big|\geqslant   \overline\alpha- C\kappa^{\frac{1}{q_0}}.
\end{equation}
Taking the  limit when $\overline \gamma$ goes to zero we  conclude the proof of Lemma \ref{lemma-Line-mes}.
\end{proof}

\section{Functional tools}\label{Sec-F-Spaces1}
In this section we set the functional framework, recall  basic  definitions and gather some technical  results that will be used along the paper.

\smallskip 

{\bf Notations.}
We denote $\mathbb{N}\triangleq \{0,1,\ldots\}$ and $\mathbb{N}^*\triangleq \{1,2,\dots\}$.
Along this paper we shall make use of the following parameters
			\begin{equation}\label{initial parameter condition}
				\kappa\in(0,1),\quad(q,d)\in(\mathbb{N}^{*})^{2},\quad S\geqslant s\geqslant s_{0}>{\tfrac{d+1}{2}}+q,
			\end{equation}
			where $S$ is a fixed large number.

\subsection{Function spaces}
We shall  introduce various function   spaces  that will be used frequently throughout this paper. 
The first one is the classical Sobolev space $H^{s}(\mathbb{T}^{d}\times\T;\C)$ which  is the set of all complex periodic  functions $h:\mathbb{T}^{d}\times\T\to\C$ such that 
$$
 h=\sum_{(l,j)\in \Z^{d+1}}h_{l,j}{\bf{e}}_{l,j},\quad h_{l,j}\in\C,\quad \quad \|h\|_{H^s}^2=\sum_{(l,j)\in \Z^{d+1}}\langle l,j\rangle^{2s}|h_{l,j}|^2,
$$
where
$$
{\bf{e}}_{l,j}(\varphi,\theta)=e^{\ii (l\cdot\varphi+j\cdot\theta)}\quad\hbox{and}\quad \langle l,j\rangle\triangleq\max(1,|l|,|j|).
$$
 For some reasons connected with the reversibility of the system we need to distinguish the following sub-spaces,  
\begin{equation}\label{FS-even}
 H^s_{\textnormal{even}}=\Big\{h\in H^{s}(\mathbb{T}^{d+1};\RR);\;\, h(-\varphi,-\theta)=h(\varphi,\theta),\;\, \forall (\varphi,\theta)\in\T^{d+1}\Big\}
\end{equation}
and
$$
H^s_{\textnormal{odd}}=\Big\{h\in H^{s}(\mathbb{T}^{d+1};\RR);\;\,  h(-\varphi,-\theta)=-h(\varphi,\theta), \;\, \forall (\varphi,\theta)\in\T^{d+1}\Big\}
$$
We  denote $\displaystyle{H^\infty\triangleq \cap_{s\in\RR} H^s}$ and similarly we define the subspaces $H^\infty_{\textnormal{even}}$ and $H^\infty_{\textnormal{even}}.$

\smallskip

To implement Nash-Moser scheme for the nonlinear equation we need to measure the dependence of the solutions with respect to the exterior parameters $\lambda=(\omega,\alpha)\in\mathbb{R}^{d+1}$ and for this aim we need  to define the weighted Sobolev spaces.  
\begin{definition}\label{Def-WS}
Let $\kappa\in(0,1), s\in\RR, q\in\mathbb{N},d\in\mathbb{N}^\star.$ Let $\mathcal{O}$ be an open bounded subset of $\mathbb{R}^{d+1}.$\\
We define the following spaces
$$
W^{q,\infty}_\gamma(\mathcal{O},H^{s}(\T^{d+1}))=\Big\{ h:\mathcal{O}\rightarrow H^{s};\,\|h\|_{s}^{q,\kappa}<+\infty\Big\}
$$
with
\begin{equation}\label{Def-weit-norm}
\|h\|_{s}^{q,\kappa}\triangleq \sup_{\lambda\in{\mathcal{O}},\\ 
\atop |\beta|\leqslant q}\kappa^{|\beta|}\|\partial_{\lambda}^{\beta}h(\lambda,\cdot)\|_{H^{s-|\beta|}(\T^{d+1})}
\end{equation}
and
$$W^{q,\infty}_\gamma(\mathcal{O},\mathbb{C})=\Big\{ h:\mathcal{O}\rightarrow\mathbb{C};\,\|h\|^{q,\kappa}<+\infty\Big\}$$
with
\begin{equation}\label{Def-abu1}
\|h\|^{q,\kappa}\triangleq \sup_{\lambda\in{\mathcal{O}},\\ 
\atop |\beta|\leqslant q}\kappa^{|\beta|}\sup_{\lambda\in{\mathcal{O}}}|\partial_{\lambda}^{\beta}h(\lambda)|.
\end{equation}
\end{definition}
\subsection{Classical operations}
The main aim of this section is to recall some classical results related  the law products and composition laws. 
In what follows  $\mathcal{O}$ denotes an open bounded  subset of $\mathbb{R}^{d+1}.$ 
The firs result is standard and whose proof can be found in   \cite{BertiMontalto} and the references therein.
\begin{lemma}\label{Law-prodX1}
Let $\kappa\in(0,1), q\in\mathbb{N}$ and $ d\in\mathbb{N}^\star$. Then the following assertions hold true.
\begin{enumerate}
 \item Let $s \geqslant0, s_0>\frac{d+1}{2}$ and $f,g\in H^{s}(\T^{d+1})\cap H^{s_0}(\T^{d+1}),$  then 
$$
\|fg\|_{H^{s}}\lesssim \|f\|_{H^{s}}\|g\|_{H^{s_0}}+\|g\|_{H^{{s}}}\|f\|_{H^{s_0}}.
$$

\item Let $s\geqslant q,   s_{0}>\frac{d+1}{2}+q$ and  $f,g\in W^{q,\infty}_{\kappa}(\mathcal{O},H^{s}(\T^{d+1})),$ then $fg\in W^{q,\infty}_{\kappa}(\mathcal{O},H^{s}(\T^{d+1}))$ with 
$$\|fg\|_{s}^{q,\kappa}\lesssim\|f\|_{s_{0}}^{q,\kappa}\|g\|_{s}^{q,\kappa}+\|f\|_{s}^{\gamma,q}\|g\|_{s_{0}}^{q,\kappa}.$$
\item  Let $f,g\in W^{q,\infty}_{\kappa}(\mathcal{O},\mathbb{C}),$ 
then $fg\in W^{q,\infty}_{\kappa}(\mathcal{O},\mathbb{C})$ with
$$\|fg\|^{q,\kappa}\lesssim\|f\|^{q,\kappa}\|g\|^{q,\kappa}.$$
\item Let $s\geqslant q$ and $(f,g)\in W^{q,\infty}_{\kappa}(\mathcal{O},\mathbb{C})\times W^{q,\infty}_{\kappa}(\mathcal{O},H^{s}(\T^{d+1})),$ 
then $fg\in W^{q,\infty}_{\kappa}(\mathcal{O},H^{s}(\T^{d+1}))$ with 
$$\|fg\|_{s}^{q,\kappa}\lesssim\|f\|^{q,\kappa}\|g\|_{s}^{q,\kappa}.$$
\end{enumerate}
\end{lemma}
Next we shall introduce a cut-off frequency operator which has the advantage to smooth out functions. For $N\in\mathbb{N}^{*},$ we define   the orthogonal projections  on the space $H^s(\T^{d+1})$ through
$$\Pi_{N} h=\sum_{\underset{\langle l,j\rangle\leqslant N}{(l,j)\in\mathbb{Z}^{d+1}}} h_{l,j}\mathbf{e}_{l,n}\quad \mbox{ and }\quad \Pi^{\perp}_{N}=\textnormal{Id}-\Pi_{N}.$$
The next result is elementary and can be easily checked. 
\begin{lemma}\label{orthog-Lem1}
Let $N\in\mathbb{N}^{*}, s\in\RR$ and $\mu\in\mathbb{R}_{+}.$ Then 
\begin{equation*}
\|\Pi_{N}h\|_{H^{s+\mu}}\leqslant N^{\mu}\|h\|_{H^{s}},\quad \|\Pi_{N}^{\perp}h\|_{H^{s}}\leqslant N^{-\mu}\|h\|_{H^{s+\mu}}.
\end{equation*}
and
\begin{equation*}
\|\Pi_{N}h\|_{s+\mu}^{q,\kappa}\leqslant N^{\mu}\|h\|_{s}^{q,\kappa},\quad \|\Pi_{N}^{\perp}h\|_{s}^{q,\kappa}\leqslant N^{-\mu}\|h\|_{s+\mu}^{q,\kappa}.
\end{equation*}
\end{lemma}
Now we shall recall the  classical interpolation inequalities, see for instance \cite{BertiMontalto}.
\begin{lemma}\label{interpolation-In}
{Let $\kappa\in(0,1), q\in\mathbb{N}, d\in\N^\star$ and $s_{1}\leqslant s\leqslant s_{2}\in\mathbb{R}$ such that $s=\theta s_{1}+(1-\theta)s_{2}$ with $\theta\in[0,1]$. Then, for any   $h\in W^{q,\infty}_{\kappa}(\mathcal{O},H^{s_{2}})$,
$$\|h\|_{s}^{q,\kappa}\lesssim \left(\|h\|_{q,s_{1}}^{q,\kappa}\right)^{\theta}\left(\|h\|_{q,s_{2}}^{q,\kappa}\right)^{1-\theta}.$$}
\end{lemma}
The proof of the next result dealing with some composition laws can be found for instance in \cite{BertiMontalto}.
\begin{lemma}\label{Compos-lemm}
Let $s\geqslant s_0>\frac{d+1}{2}+q $ and $\kappa\in(0,1).$ Let $f:\mathcal{O}\times\RR\to \RR$ be $\mathscr{C}^\infty$ and  $h\in W^{q,\infty}_{\kappa}(\mathcal{O},H^{s}(\T^{d+1}))$ such that $\|h\|_{s_{0}}^{q,\kappa}\leqslant C_0$ for and arbitrary $C_0>0$ and define the point-wise composition
$$
\forall\, (\lambda,x)\in\mathcal{O}\times\T^{d+1},\quad F(h)(\lambda,x)\triangleq F(\lambda, h(\lambda,x)).
$$
Then $F(h)\in  W^{q,\infty}_{\kappa}(\mathcal{O},H^{s}(\T^{d+1}))$ with  
$$
\|F(h)-F(\lambda,0)\|_{s}^{q,\kappa}\leqslant C(s,F, C_0) \|h\|_{s}^{q,\kappa}.
$$
\end{lemma}
We shall also make use of the following classical result that can proved directly from Taylor formula combined with  the chain rule (Fa\'a di Bruno's formula) and interpolation inequalities in $W^{k,\infty}$.
\begin{lemma}\label{Compos-lemm-VM}
Let $q\in\N^\star$ and  $F:\RR\to\RR$ be a $C^\infty$ function with bounded derivatives and $h:\mathcal{O}\to \RR$ be  a function in $W^{q,\infty}(\mathcal{O})$. Then the composition $F(h) \in W^{q,\infty}(\mathcal{O})$ with
$$
\|F(h)-F(0)\|_{W^{q,\infty}}\leqslant C(q,F) \|u\|_{W^{q,\infty}}\big(1+\|h\|_{L^\infty}^{q-1}\big).
$$
\end{lemma}

The next result will be useful later in studying the asymptotic structure  of the linearized operator.
 \begin{lemma}\label{Lemm-RegZ}
Let $\epsilon>0, d\in\N^\star$ and $q\in\N$, then the following assertions  hold true.
\begin{enumerate}
\item 
Let $s\geqslant 0,$ there exists $C>0$ such that if $f\in H^{s+\frac52+\epsilon}(\T^2)$ and  satisfying in addition the symmetry property  $f(\theta,\eta)=f(\eta,\theta)$, then there exists $g:\T^2\to\RR$ such that
$$
f(\theta,\eta)=\tfrac12\Big(f(\theta,\theta)+f(\eta,\eta)\Big)+\sin^2\left(\tfrac{\eta-\theta}{2}\right)g(\theta,\eta)
$$
with 
$$
\sup_{\eta\in\T}\|g(\cdot,\cdot+\eta\|_{H^{s}(\T)}\leqslant C\|f\|_{H^{s+\frac{5}{2}+\epsilon}(\T^2)}.
$$
More generally, we have the estimate: for any $p_1,p_2\in\N$
\begin{align*}
\sup_{\eta\in\T}\|(\partial_\theta^{p_1}\partial_\eta^{p_2}g)(\cdot,\cdot+\eta)\|_{H^s(\T)}\lesssim \|f\|_{H^{s+\epsilon+\frac52+p_1+p_2}(\T^2)}.
\end{align*}

\item
Let $s\geqslant q,$ there exists $C>0$ such that if  $f\in W^{q,\infty}_{\kappa}\big(\mathcal{O};H^{s+\frac52+\epsilon}(\T^{d}\times\T^2)\big)$ being a symmetric function: $f(\lambda,\varphi,\theta,\eta)=f(\lambda,\varphi,\eta,\theta)$ then there exists $g$ such that
$$
f(\lambda,\varphi,\theta,\eta)=\tfrac12\Big(f(\lambda,\varphi,\theta,\theta)+f(\lambda,\varphi,\eta,\eta)\Big)+\sin^2\left(\tfrac{\eta-\theta}{2}\right)g(\lambda,\varphi,\theta,\eta)
$$
with 
\begin{align*}
\sup_{\eta\in\T}\big\|g\big(\ast,\cdot,\centerdot,\centerdot+\eta\big)\big\|_{s}^{q,\kappa}&\leqslant  C\|\Delta_{\theta,\eta}f\|_{s+\frac12+\epsilon}^{q,\kappa}\\
&\leqslant  C\|f\|_{s+\frac52+\epsilon}^{q,\kappa}
\end{align*}
and where  $\Delta_{\theta,\eta}=\partial_{\theta}^2+\partial_\eta^2.$ More generally, we have the estimate: for any $p_1,p_2\in\N$
\begin{align*}
\sup_{\eta\in\T}\big\|(\partial_\theta^{p_1}\partial_\eta^{p_2}g)\big(\ast,\cdot,\centerdot,\centerdot+\eta\big)\big\|_{s}^{q,\kappa}&\leqslant  C\|f\|_{s+\epsilon+\frac52+p_1+p_2}^{q,\kappa}.
\end{align*}
Here the symbols $\ast,\cdot,\centerdot$ denote $\lambda,\varphi,\theta,$ respectively.
\end{enumerate}
\end{lemma}
\begin{proof}
${\bf{(i)}}$ Let us split $f$ through its Fourier expansion,
$$
f(\theta,\eta)=\sum_{j,j'\in\Z}a_{j,k}e^{\ii (j\theta+j'\eta)}
$$
Since $f$ is symmetric then one finds  from some algebraic manipulations, 
\begin{align*}
f(\theta,\eta)=&\frac12\sum_{j' ,j\in\Z}a_{j,j' }\big[e^{\ii (j \theta+j' \eta)}+e^{\ii (j \eta+j' \theta)}\big]\\
=&\frac12\big(f(\theta,\theta)+f(\eta,\eta)\big)+f_1(\theta,\eta)
\end{align*}
where
\begin{equation}\label{RestZ00}
f_1(\theta,\eta)=\frac12\sum_{j' ,j\in\Z}a_{j,j' }\mu_{j,j' }(\theta,\eta)e^{\ii (j \theta+j' \eta)}
\end{equation}
and with
$$
\mu_{j,j' }(\theta,\eta)=\big(1-e^{\ii j(\eta-\theta)}\big)\big(1-e^{-\ii  j' (\eta-\theta)}\big).
$$
Using elementary  trigonometric identities combined using  Chebyshev polynomials of second order $U_n$  that
\begin{align*}
\mu_{j,j' }(\theta,\eta)=&-4 \sin\left(j' \tfrac{\theta-\eta}{2}\right)\sin\left( j\tfrac{\theta-\eta}{2}\right)e^{\frac{\ii}{2}(j' -j)(\theta-\eta)}\\
=&-4 \sin^2\left(\tfrac{\theta-\eta}{2}\right)U_{j' -1}\left(\cos\left( \tfrac{\theta-\eta}{2}\right)\right)U_{j-1}\left(\cos\left( \tfrac{\theta-\eta}{2}\right)\right)e^{\frac{\ii}{2}(j' -j)(\theta-\eta)}\\
=&-4 \sin^2\left(\tfrac{\theta-\eta}{2}\right)\mu_{j,j' }^1(\theta,\eta)
\end{align*}
with
\begin{align}\label{TchebZ1}
\mu_{j,j' }^1(\theta,\eta)\triangleq U_{j' -1}\left(\cos\left( \tfrac{\theta-\eta}{2}\right)\right)U_{j-1}\left(\cos\left( \tfrac{\theta-\eta}{2}\right)\right)e^{\frac12\ii(j' -j)(\theta-\eta)}.
\end{align}
Here, we have changed slightly the convention for $U_j$, it is defined as the polynomial such that
$$
\sin(j\theta)=\sin(\theta)U_j(\cos(\theta)), \forall \theta\in \RR.
$$
With this convention we have  $U_0\equiv0$ and it is still defined for negative integer, with the identity: $U_{-j}=-U_j$. Plugging this into \eqref{RestZ00} allows to get
$$
f_1(\theta,\eta)=\sin^2\left(\tfrac{\theta-\eta}{2}\right)g(\theta,\eta)
$$
with
\begin{align}\label{Form-RP1}
g(\theta,\eta)=-2\sum_{j' ,j\in\Z}a_{j,j' }\mu_{j,j' }^1(\theta,\eta)e^{\ii (j\theta+j' \eta)}.
\end{align}
It follows that
$$
g(\theta,\theta+\eta)=-2\sum_{j' ,j\in\Z}a_{j,j' }\,\mu_{j,j' }^1(\theta,\theta+\eta)\,e^{\ii j' \cdot\eta}e^{\ii (j+j' )\theta}.
$$
One can check easily that
\begin{align}\label{mujj}
\mu_{j,j'}^1(\theta,\theta+\eta)&=U_{j'-1}\big(\cos\left( \tfrac{\eta}{2}\right)\big)U_{j-1}\big(\cos\left( \tfrac{\eta}{2}\right)\big)e^{\frac12\ii(j-j')\eta}\\
&=\mu_{j,j'}^1(0,\eta).\nonumber
\end{align}
Thus we find after a change of indices
\begin{align*}
g(\theta,\theta+\eta)&=-2\sum_{j',j\in\Z^d}a_{j,j'}\,\mu_{j,j'}^1(0,\eta)\,e^{\ii j'\cdot\eta}e^{\ii (j+j')\theta}\\
&=-2\sum_{j',j\in\Z}a_{j-j',j'}\,\mu_{j-j',j'}^1(0,\eta)\,e^{\ii j'\eta}e^{\ii j \theta}.
\end{align*}
Therefore
$$
\|g(\cdot,\cdot+\eta)\|_{H^s}^2=4\sum_{j\in\Z}\langle j\rangle^{2s}\Big|\sum_{j'\in \Z} a_{j-j',j'}\,\mu_{j-j',j'}^1(0,\eta)e^{\ii j'\eta}\Big|^2.
$$
Using Cauchy-Schwarz inequality we obtain
\begin{align*}
\Big|\sum_{j'\in \Z} a_{j-j',j'}\,\mu_{j-j',j'}^1(0,\eta)e^{\ii j'\cdot\eta}\Big|^2\leqslant&\left(\sum_{j'\in\Z}\frac{1}{\langle j'-j\rangle^{1+2\epsilon}}\right)\sum_{j'\in \Z} \langle j'-j\rangle^{1+2\epsilon} |a_{j-j',j'}|^2\,|\mu_{j-j',j'}^1(0,\eta)|^2\\
\lesssim&\sum_{j'\in \Z} \langle j'-j\rangle^{1+2\epsilon} |a_{j-j',j'}|^2\,|\mu_{j-j',j'}^1(0,\eta)|^2.
\end{align*}
Now, from \eqref{mujj} we get, for any $k\in\N$ there exists $C>0$ such that 
\begin{align*}
\forall\, j\in\Z,\quad \sup_{\theta\in[0,2\pi]}\big|\partial_\theta^k\big(U_j(\cos \theta)\big)\big|\leqslant C |j|^{k+1}.
\end{align*}
Combining this with \eqref{TchebZ1} together with Leibniz formula yield
\begin{align}\label{techeb-be}
\sup_{\theta,\eta\in\T}\big(|\partial_\theta^{k}\mu_{j,j'}^1(\theta,\eta)|+|\partial_\eta^{k}\mu_{j,j'}^1(\theta,\eta)|\big)\lesssim \left(|j|+|j'|\right)^{{k}+2}.
\end{align}
Applying this for ${k}=0$ implies
\begin{align*}
\sum_{j'\in \Z} \langle j'-j\rangle^{1+2\epsilon} |a_{j-j',j'}|^2|\mu_{j-j',j'}^1(0,\eta)|^2&\lesssim\sum_{j'\in \Z} \langle j'-j\rangle^{1+2\epsilon}\left(|j-j'|+| j'|\right)^4 |a_{j-j',j'}|^2.
\end{align*}
It follows that for $s\geqslant0$
\begin{align*}
\|g(\cdot,\cdot+\eta)\|_{H^s}^2\lesssim& \sum_{j,j'\in\Z}\langle j'-j\rangle^{1+2\epsilon}\left(| j-j'|+| j'|\right)^4 \langle j\rangle^{2s} |a_{j-j',j'}|^2\\
\lesssim& \sum_{j,j'\in\Z}\langle j\rangle^{1+2\epsilon}\left(|j|+|j'|\right)^4 \langle j'+j\rangle^{2s} |a_{j,j'}|^2\\
\lesssim& \sum_{j,j'\in\Z}\langle j',j\rangle^{5+2\epsilon+2s} |a_{j,j'}|^2.
\end{align*}
Thus we deduce that
$$
\sup_{\eta\in\T}\|g(\cdot,\cdot+\eta)\|_{H^s}\lesssim \|f\|_{H^{s+\frac52+\epsilon}}.
$$
Implementing the same approach as before  and using \eqref{techeb-be} we get
\begin{align*}
\sup_{\eta\in\T}\left(\|\partial_\theta^{k}g(\cdot,\cdot+\eta)\|_{H^s}+\|(\partial_\eta^{k}g)(\cdot,\cdot+\eta)\|_{H^s}\right)\lesssim \|f\|_{H^{s+\frac52+{k}+\epsilon}}.
\end{align*}
The proof of the following estimates can be done in a similar way,
\begin{align*}
\sup_{\eta\in\T}\|\partial_\theta^{p_1}\partial_\eta^{p_2}g(\cdot,\cdot+\eta)\|_{H^s}\lesssim \|f\|_{H^{s+\frac52+p_1+p_2+\epsilon}}.
\end{align*} 
${\bf{(ii)}}$ Splitting as before $f$ into its Fourier expansion 
$$
f(\lambda, \varphi,\theta,\eta)=\sum_{j,j'\in\Z\\
\atop l\in \Z^d}a_{j,j',l}(\lambda)e^{\ii (j\theta+j'\eta+l\cdot\varphi)}.
$$
Then proceeding as in the first point ${\bf{(i)}}$ we may obtain  $$
f(\lambda,\varphi,\theta,\eta)=\tfrac12\big(f(\lambda,\varphi,\theta,\theta)+f(\lambda,\varphi,\eta,\eta)\big)+\sin^2\left(\tfrac{\eta-\theta}{2}\right)g(\lambda,\varphi,\theta,\eta)
$$
and $g$ takes a similar  form to  \eqref{Form-RP1}
\begin{align}\label{Form-RP2}
g(\lambda,\varphi,\theta,\eta)=-2\sum_{j',j\in\Z\atop l\in\Z^d}a_{j,j'}(\lambda)\mu_{j,j'}^1(\theta,\eta)e^{\ii (j\theta+j'\eta+l\cdot\varphi)}.
\end{align}
The coefficients $\mu_{j,j'}^1$ are given by \eqref{TchebZ1}. Thus we get
\begin{align}\label{Form-RP3}
g(\lambda,\varphi,\theta,\theta+\eta)=-2\sum_{j',j\in\Z\atop l\in\Z^d}a_{j,j',l}(\lambda)\mu_{j,j'}^1(0,\eta)e^{\ii j'\cdot\eta}e^{\ii ((j+j')\theta+l\cdot\varphi)}.
\end{align}
By changing the indices
\begin{align*}
g(\lambda,\varphi,\theta,\theta+\eta)&=-2\sum_{j',j\in\Z\atop l\in\Z^d}a_{j,j',l}(\lambda)\mu_{j,j'}^1(0,\eta)e^{\ii j'\cdot\eta}e^{\ii ((j+j')\theta+l\cdot\varphi)}\\
&=-2\sum_{j',j\in\Z\atop l\in\Z^d}a_{j-j',j',l}(\lambda)\mu_{j-j',j'}^1(0,\eta)e^{\ii j'\eta}e^{\ii (j \theta+l\cdot\varphi)}.
\end{align*}
Then for ${k}\in\N^{d+1}, |{k}|\leqslant q$
\begin{align*}
\big\|\partial_\lambda^{k}g\big(\lambda,\cdot,\centerdot,\centerdot+\eta\big)\big\|_{H^{s-|{k}|}}^2&= \sum_{j\in\Z,l\in\Z^d}\langle j,l\rangle^{2(s-|{k}|)}\Big|\sum_{j'\in \Z} \partial_\lambda^{\mathtt{j}}a_{j-j',j',l}(\lambda)\mu_{j-j',j'}^1(0,\eta)e^{\ii j'\eta}\Big|^2.
\end{align*}
Using Cauchy-Schwarz inequality we obtain for any $\epsilon>0$
\begin{align*}
\Big|\sum_{j'\in \Z} \partial_\lambda^{k}a_{j-j',j',l}(\lambda)\mu_{j-j',j'}^1(0,\eta)e^{\ii j'\cdot\eta}\Big|^2\lesssim&\sum_{j'\in \Z} \langle j'-j\rangle^{1+2\epsilon} \big|\partial_\lambda^{k}a_{j-j',j',l}(\lambda)\mu_{j-j',j'}^1(0,\eta)\big|^2.
\end{align*}
Combining this with \eqref{techeb-be} yields
\begin{align*}
\sum_{j'\in \Z} \langle j'-j\rangle^{1+2\epsilon} \big| \partial_\lambda^{k}a_{j-j',j',l}(\lambda)\big|^2| \mu_{j-j',j'}^1(0,\eta)|^2&\lesssim\sum_{j'\in \Z} \langle j'-j\rangle^{1+2\epsilon}\left(|j-j'|+|j'|\right)^4 \big| \partial_\lambda^{k}a_{j-j',j',l}(\lambda)\big|^2.
\end{align*}
It follows that for $s\geqslant q\geqslant |{k}|,$
\begin{align}\label{Tel-11}
\nonumber\big\|\partial_\lambda^{k}g\big(\lambda,\cdot,\centerdot,\centerdot+\eta\big)\big\|_{H^{s-|{k}|}}^2\lesssim& \sum_{j,j'\in\Z\atop l\in\Z^d}\langle j'-j\rangle^{1+2\epsilon}\left(| j-j'|+| j'|\right)^4 \langle j,l\rangle^{2(s-|{k}|)}  \big| \partial_\lambda^{k}a_{j-j',j',l}(\lambda)\big|^2\\
\lesssim& \sum_{j,j'\in\Z\atop l\in\Z^d}\langle j\rangle^{1+2\epsilon}\left(|j|+|j'|\right)^4 \langle j'+j,l\rangle^{2(s-|{k}|)} |\partial_\lambda^{k}a_{j,j',l}(\lambda)|^2\\
\nonumber\lesssim& \sum_{j,j'\in\Z\atop l\in\Z^d}\langle j',j,l\rangle^{2(s-|{k}|+\frac52+\epsilon)}  |\partial_\lambda^{k}a_{j,j',l}(\lambda)|^2.
\end{align}
It follows that for any $\lambda\in\mathcal{O}, \eta\in\T$ and ${k}\in\N^{d+1}, |{k}|\leqslant q$
$$
\big\|\partial_\lambda^{k}g\big(\lambda,\cdot,\centerdot,\centerdot+\eta\big)\big\|_{H^{s-|{k}|}(\T^d\times\T)}\lesssim \big\|\partial_\lambda^{k}f\big\|_{H^{s-|{k}|+\frac52+\epsilon}(\T^d\times\T^2)}.
$$
Consequently
$$
\sup_{\eta\in\T}\big\|g\big(\ast,\cdot,\centerdot,\centerdot+\eta\big)\big\|_{s}^{q,\kappa}\leqslant  C\|f\|_{s+\frac52+\epsilon}^{q,\kappa}.
$$
We can also get from \eqref{Tel-11} that 
\begin{align}\label{Est-Pr01}
\sup_{\eta\in\T}\big\|g\big(\ast,\cdot,\centerdot,\centerdot+\eta\big)\big\|_{s}^{q,\kappa}\leqslant  C\|\Delta_{\theta,\eta}f\|_{s+\frac12+\epsilon}^{q,\kappa}.
\end{align}
Now let us move to the  estimate of $(\partial_\theta^pg).$ %
Applying  Leibniz formula to \eqref{Form-RP2} yields, for $p\in\N$,
\begin{equation*}
\partial_\theta^{p}g(\lambda,\varphi,\theta,\eta)=-2\sum_{0\leqslant m\leqslant p\\
\atop 
 l,j\in\Z,\\
 l\in\Z^d}\left(_{{m}}^{{p}}\right) a_{j,j'}(\lambda)(\ii j)^{p-m}\partial_\theta^{m}\mu_{j,j'}^1(\theta,\eta)e^{\ii (j\theta+j'\eta+l\cdot\varphi)}.
\end{equation*}
It follows that
\begin{align*}
(\partial_\theta^pg)(\lambda,\varphi,\theta,\theta+\eta)&=-2\sum_{m=0}^p\sum_{j',j\in\Z, l\in\Z^d}\left(_m^p\right)a_{j,j'}(\lambda)(\ii j)^{p-m}(\partial_\theta^m\mu_{j,j'}^1)(0,\eta)e^{\ii j'\cdot\eta}e^{\ii ((j+j')\theta+l\cdot\varphi)}\\
&=-2\sum_{m=0}^p\sum_{j',j\in\Z; l\in\Z^d}\left(_m^p\right)a_{j-j',j'}(\lambda)(\ii j)^{p-m}(\partial_\theta^m\mu_{j-j',j'}^1)(0,\eta)e^{\ii j'\cdot\eta}e^{\ii (j\theta+l\cdot\varphi)}.
\end{align*}
Then from \eqref{techeb-be} and proceeding as before we get
\begin{align}\label{Tel-L1}
\nonumber\big\|\partial_\lambda^{k} (\partial_\theta^pg)\big(\lambda,\cdot,\centerdot,\centerdot+\eta\big)\big\|_{H^{s-|{k}|}}\lesssim& \sum_{j,j'\in\Z\atop l\in\Z^d}\langle j'-j\rangle^{1+2\epsilon}\left(| j-j'|+| j'|\right)^{4+2p} \langle j,l\rangle^{2(s-|{k}|)} \big| \partial_\lambda^{k}a_{j-j',j',l}(\lambda)\big|^2\\
\lesssim& \sum_{j,j'\in\Z\atop l\in\Z^d}\langle j',j,l\rangle^{2(s-|{k}|+\frac52+\epsilon+p)}  |\partial_\lambda^{k}a_{j,j',l}(\lambda)|^2.
\end{align}
Hence,  for any $\lambda\in\mathcal{O}, \eta\in\T$ and $p\in\N$,
$$
\big\|\partial_\lambda^{k} (\partial_\theta^pg)\big(\lambda,\cdot,\centerdot,\centerdot+\eta\big)\big\|_{H^{s-|{k}|}(\T^d\times\T)}\lesssim \big\|\partial_\lambda^{k}f\big\|_{H^{s_1-|{k}|+\frac52+p+\epsilon}(\T^d\times\T^2)}.
$$
Consequently
$$
\sup_{\eta\in\T}\big\| (\partial_\theta^pg)\big(\ast,\cdot,\centerdot,\centerdot+\eta\big)\big\|_{s}^{q,\kappa}\leqslant  C\|f\|_{s+\frac52+p+\epsilon}^{q,\kappa}.
$$
Following the same lines  by applying  in particular Leibniz formula to \eqref{Form-RP2} yields, for $p_1,p_2\in\N$,
\begin{equation*}
\partial_\theta^{p_1}\partial_\eta^{p_2}g(\lambda,\varphi,\theta,\eta)=2\sum_{0\leqslant m_1,m_2\leqslant p\\
\atop 
 j',j\in\Z\\
 l\in\Z^d}\left(_{{m_1}}^{{p_1}}\right) \left(_{{m_2}}^{{p_2}}\right)a_{j,j'}(\lambda)(\ii n)^{p_1-m_1}(\ii j')^{p_2-m_2}\partial_\theta^{m_1}\partial_\eta^{m_2}\mu_{j,j'}^1(\theta,\eta)e^{\ii (j\theta+j'\eta+l\cdot\varphi)}.
\end{equation*}
This gives the desired estimate.
\end{proof}
Next we shall prove the following lemma which turns out to be  very useful.
\begin{lemma} \label{lem-Reg1}
 Let $f:\mathcal{O}\times\mathbb{T}^{d}\times\mathbb{T}\rightarrow\mathbb{C}$ be a smooth function and  define  ${g}:\mathcal{O}\times\mathbb{T}^{d}\times\mathbb{T}^2\rightarrow\mathbb{C}$ by 
 $${g}(\lambda,\varphi,\theta,\eta)=\left\lbrace\begin{array}{ll}
\frac{f(\lambda,\varphi,\eta+\theta)-f(\lambda,\varphi,\theta)}{\tan\left(\frac{\eta}{2}\right)} & \mbox{if }\eta\in\RR\backslash2\pi\Z\\
2\partial_{\theta}f(\lambda,\varphi,\theta) & \mbox{if }\eta\in 2\pi\Z.
\end{array}\right.$$
Then the followings assertions hold true.
\begin{enumerate}
\item For any $s\in\RR$, we have
$$
\sup_{\eta\in\T}\|g(\ast,\cdot,\centerdot, \eta)\|_{s}^{\gamma,q}\lesssim \|f\|_{s+1}^{q,\kappa}.
$$
\item For any $s\geqslant0$
\begin{align*}
\|{g}(\ast,\cdot,\centerdot)\|_{s}^{q,\kappa}&\lesssim\|\partial_{\theta}f\|_{s}^{q,\kappa}\\
&\lesssim\|f\|_{s+1}^{q,\kappa},
\end{align*}
where the symbols $\ast,\cdot,\centerdot$ denote  $\lambda, \varphi$  and  $(\theta,\eta)$. 
\end{enumerate}
\end{lemma}
\begin{proof}
{\bf{(i)}} 
We  start with  expanding $f$ into its Fourier series,
\begin{align}\label{F-Exp}
f(\varphi,\theta)=\sum_{(l,j)\in\mathbb{Z}^{d+1 }}f_{l,j}e^{\ii (l\cdot\varphi+j\theta)}.
\end{align}
Therefore we get from the Bessel's identity combined with some trigonometric identities  
$$\begin{array}{rcl}
\|g(\cdot,\cdot,\eta)\|_{H_{\varphi,\theta}^{s}}^{2} & = & \displaystyle\Big\|\sum_{(l,j)\in\mathbb{Z}^{d+1}}f_{l,j}\frac{e^{\ii j\eta}-1}{\tan\left(\frac{\eta}{2}\right)}e^{\ii (l\cdot\varphi+j\theta)}\Big\|_{H^{s}_{\varphi,\theta}}^{2}\\
& \leqslant & \displaystyle\Big\|\sum_{(l,j)\in\mathbb{Z}^{d+1 }}f_{l,j}e^{\ii \frac{j\eta}{2}}\tfrac{\sin\left(\frac{j\eta}{2}\right)}{\tan\left(\frac{\eta}{2}\right)}e^{\ii (l\cdot\varphi+j\theta)}\Big\|_{H^{s}_{\varphi,\theta}}^{2}\\
& \lesssim & \displaystyle\sum_{(l,j)\in\mathbb{Z}^{d+1 }}\langle l,j\rangle^{2s}|f_{l,j}|^{2}\tfrac{\sin^{2}\left(\frac{j\eta}{2}\right)}{\sin^{2}\left(\frac{\eta}{2}\right)}.
\end{array}$$
Using the bound $\displaystyle \tfrac{\sin^{2}\left(\frac{j\eta}{2}\right)}{\sin^{2}\left(\frac{\eta}{2}\right)}\lesssim |j| $ following from basic properties of F\'ejer kernel, we infer
$$\begin{array}{rcl}
\|g(\cdot,\cdot,\eta)\|_{H_{\varphi,\theta}^{s}}^{2} 
& {\lesssim} & \displaystyle\sum_{(l,j)\in\mathbb{Z}^{d}\times\mathbb{Z}}\langle l,j\rangle^{2s}|j|^{2}|f_{l,j}|^{2}\\
& \lesssim & \displaystyle\|\partial_{\theta}f\|_{H^{s}}^{2}.
\end{array}$$
This ensures  the desired result for $q=0.$

\smallskip

{\bf{(ii)}} 
 We first split $f$ into its Fourier expansion,
$$
f(\lambda,\varphi,\theta)=\sum_{(l,j)\in \Z^{d+1}}c_{l,j}(\lambda)e^{\ii ( l\cdot \varphi+j\theta)}.
$$
Consequently
$$
g(\lambda,\varphi,\theta,\eta)=\sum_{(l,j)\in \Z^{d+1}}\mu_j(\eta)\,c_{j',j}(\lambda)e^{\ii (l\cdot \varphi+j\theta)},\quad \mu_j(\eta)\triangleq \frac{e^{\ii  j\eta}-1}{\tan(\frac\eta2)}\cdot
$$
The latter function is  $\mathscr{C}^\infty$ and $2\pi-$periodic and therefore its Fourier expansion takes the form
$$
\mu_j(\eta)=\sum_{j'\in\Z}A_{j,j'}e^{\ii j'\eta}.
$$
It follows that, for any ${k}\in\N^{d+1}$ with  $|{k}|\leqslant q$,
$$
\partial_\lambda^{k}g(\lambda,\varphi,\theta,\eta)=\sum_{(l,j,j')\in \Z^{d+2}}\,(\partial_\lambda^{k}c_{l,j}(\lambda))A_{j,j'}e^{\ii (l\cdot \varphi+j\theta+j'\eta)}.
$$
Thus we find
\begin{align}\label{Sam-X1}
\nonumber \|\partial_\lambda^{k}g(\lambda,\cdot,\centerdot)\|_{H^{s-|{k}|}}^2&=\sum_{(l,j,j')\in \Z^{d+2}}\langle l,j',j\rangle^{2(s-|{k}|)}\,\big|\partial_\lambda^{k}c_{l,j}(\lambda)\big|^2|A_{j,j'}|^2\\
&\lesssim \sum_{(l,j)\in \Z^{d+1}}\nu_{l,j}\big|\partial_\lambda^{k}c_{l,j}(\lambda)\big|^2
\end{align}
with
$$
\nu_{l,j}\triangleq \sum_{j'\in \Z}\langle l,j',j\rangle^{2(s-|{k}|)}\, |A_{j,j'}|^2.
$$
Since $s\geqslant q\geqslant|{k}|,$ then we can easily check  from the inequalities $(|a|+|b|)^{s}\lesssim |a|^{s}+|b|^{s}$ that
\begin{align*}
\nu_{l,j}&\lesssim\sum_{j'\in \Z}\langle l,j\rangle^{2(s-|{k}|)}\, |A_{j,j'}|^2+\sum_{j'\in \Z}\langle j'\rangle^{2(s-|{k}|)}\, |A_{j,j'}|^2.
\end{align*}
This can be written in the form
\begin{align*}
\nu_{l,j}&\lesssim\langle l,j\rangle^{2(s-|{k}|)}\,\|\mu_j\|_{L^2}^2+\, \|\mu_j\|_{H^{s-|{k}|}}^2.
\end{align*}
Next we claim that for any $s\in\RR_+$
\begin{align}\label{F-S1}
\forall \, j\in\Z,\quad \|\mu_j\|_{H^s}\lesssim |j|^{\frac12+s}.
\end{align}
Let us give the main ideas  of the proof. The case $s=0$ can be done in the following way
\begin{align*}
\|\mu_j\|_{L^2}^2\lesssim&\int_{0}^\pi\frac{\sin^2(j\eta/2)}{\sin^2(\eta/2)}d\eta\\
\lesssim&\int_{0}^{\frac\pi2}\frac{\sin^2(j\eta)}{\sin^2(\eta)}d\eta.
\end{align*}
Then using the inequality $\forall \eta\in[0,\frac\pi2],\,\sin(\eta)\geqslant \frac{2\eta}{\pi}$, we find through change of variables 
\begin{align*}
\|\mu_j\|_{L^2}^2&\lesssim\int_{0}^{\frac\pi2}\frac{\sin^2(j\eta)}{\eta^2}d\eta\\
&\lesssim |j|\int_{0}^{|j|\frac\pi2}\frac{\sin^2(\eta)}{\eta^2}d\eta\\
&\lesssim |j|.
\end{align*}
The second step  is  to establish the estimate for $s\in\N^*$ by direct differentiation and using  similar arguments as for $s=0$. The  extension of the estimate to  $s\in\RR_+$ can be obtained by interpolation inequalities.\\
Now using \eqref{F-S1} we obtain the estimate
\begin{align*}
\nu_{l,j}&\lesssim\langle l,j\rangle^{2(s-|{k}|)}|j|+\langle j\rangle^{2(s-|{k}|)}|j|\\
&\lesssim\langle l,j\rangle^{2(s-|{k}|)}|j|.
\end{align*}
Inserting this estimate into \eqref{Sam-X1} yields
\begin{align}\label{Sam-X2}
\nonumber\|\partial_\lambda^{k}g(\lambda,\cdot,\centerdot)\|_{H^{s-|{k}|}}^2
&\lesssim \sum_{(l,j)\in \Z^{d+1}}\langle l,j\rangle^{2(s-|{k}|)}|j|\big|\partial_\lambda^{k}c_{l,j}(\lambda)\big|^2\\
&\lesssim \big\|\partial_\lambda^{k} \partial_\theta f(\lambda,\cdot,\centerdot) \big\|_{H^{s-|{k}|}}.
\end{align}
From \eqref{Sam-X2} and the Definition \ref{Def-WS} we infer
\begin{align*}
\|g\|_{s}^{q,\kappa}\lesssim& \|\partial_\theta f\|_{s}^{q,\kappa}\\
\lesssim&\| f\|_{s+1}^{q,\kappa}.
\end{align*}
This achieves the proof.
\end{proof}

\subsection{Modified periodic fractional Laplacian}
The main goal of this section is to explore some basic properties of a modified fractional Laplacian defined as follows. 
Let  $h:\mathbb{T}^{d+1}\to \RR$ and consider $\alpha\in(0,1),$ we define the modified partial fractional Laplacian by  
\begin{equation}\label{fract1}
|\textnormal D|^{-\alpha}h(\varphi,\theta)=\frac{1}{2\pi}\bigintsss_{0}^{2\pi} \frac{h(\varphi,\theta-\eta)}{|\sin\big(\frac{\eta}{2}\big)|^{1-\alpha}} d\eta.
\end{equation}
As we shall see for instance in \eqref{Frac-Mult}, this operator is different from the standard fractional Laplacian denoted by $\Lambda^{-\alpha}=(-\Delta)^{-\frac{\alpha}{2}}$ but asymptotically in Fourier side they have the same behavior.\\
To state some results connected with the action of the modified fractional Laplacian we need to  introduce Sobolev  anisotropic spaces. First, we define the anisotropic Fourier multiplier  $\Lambda^{s_1,s_2}$  as  follows: For $s_1,s_2\in\RR$
\begin{equation}\label{Fmult}
\Lambda^{s_1,s_2}{\bf e}_{l,j}=\langle l,j\rangle^{s_1}\langle j\rangle^{s_2}{\bf e}_{l,j}\quad\hbox{with}\quad \quad \langle j\rangle=\max(1,|j|).
\end{equation} 
  The anisotropic Sobolev space $H^{s_1,s_2}(\mathbb{T}^{d+1};\C)$ is  the set of   functions $h:\mathbb{T}^{d+1}\to\mathbb{C}$ such that
$$
 h=\sum_{(l,j)\in \Z^{d+1}}f_{l,j}{\bf{e}}_{l,j},\quad h_{l,j}\in\C\quad\hbox{and}\quad \|h\|_{H^{s_1,s_2}}^2=\sum_{(l,j)\in \Z^{d+1}}\langle l,j\rangle^{2s_1}\langle j\rangle^{2s_2}|h_{l,j}|^2<\infty.
$$
The main goal here is to collect some  elementary results related to the action of   the modified  fractional Laplacian to anisotropic spaces.
\begin{lemma} \label{Laplac-frac0} Let    $s_1,s_2 \in\RR$ and $\alpha\in (0,1)$, then the following properties hold true.
\begin{enumerate}
\item For any  $p\in\N$, we have $ \partial_\theta^p: H^{s_1,s_2}\to H^{s_1,s_2-p}$ is continuous.
\item We have  $ |\textnormal D|^{-\alpha}: H^{s_1,s_2}\to H^{s_1, s_2+\alpha}$ is continuous

\item For any  $ k\in\N$, then 
$$ \partial_\alpha^k|\textnormal D|^{-\alpha} : H^{s_1,s_2}\to H^{s_1, s_2+\alpha-\epsilon}$$ is continuous, with 
$$
\partial_\alpha^k|\textnormal D|^{-\alpha}h(\varphi,\theta)=\frac{1}{2\pi}\bigintsss_{0}^{2\pi} \frac{\ln^{k}\big[\sin\big({\eta}/{2}\big)\big]}{|\sin\big(\frac{\eta}{2}\big)|^{1-\alpha}} \,h(\varphi,\theta-\eta)d\eta.
$$

\end{enumerate}\end{lemma}
\begin{proof}
{\bf{(i)}} This can be easily  obtained  from
$$
h=\sum_{(l,j)\in \Z^{d+1}}h_{n,k}{\bf{e}}_{l,j}\Longrightarrow \partial_\theta^kh=\sum_{(l,j)\in \Z^{d+1}}(\ii j)^kh_{l,j}{\bf{e}}_{l,j}.
$$
Thus
\begin{eqnarray*}
\|\partial_\theta^p h\|_{H^{s_1,s_2-p}}^2&=&\sum_{(l,j)\in \Z^{+1}}\langle l,j\rangle^{2s_1}\langle j\rangle^{2(s_2-p)}|j|^{2l}|h_{l,j}|^2\\
&\le&\|h\|_{H^{s_1,s_2}}^2.
\end{eqnarray*}
{\bf{(ii)}} Using the definition and a change of variables  we find
\begin{eqnarray*}
|\textnormal D|^{-\alpha}{\bf e}_{l,j}(\varphi,\theta)
&=&\frac{1}{2\pi}\bigintsss_{0}^{2\pi} \frac{e^{-\ii  j\,\eta}}{|\sin\big(\frac{\eta}{2}\big)|^{1-\alpha}} d\eta\,\, {\bf e}_{l,j}(\varphi,\theta)\\
&=&\frac{1}{\pi}\bigintsss_{0}^{\pi} \frac{e^{2\ii j\eta}}{|\sin{\eta}|^{1-\alpha}} d\eta\,\, {\bf e}_{l,j}(\varphi,\theta).
\end{eqnarray*}
Applying  \eqref{int} and using \eqref{Pocc16}  give
\begin{align}\label{Fo-c}
|\textnormal D|^{-\alpha}{\bf e}_{l,j}(\varphi,\theta)&=\frac{2^{1-\alpha}\Gamma(\alpha)}{\Gamma(\frac{1-\alpha}{2})\Gamma(\frac{1+\alpha}{2})}\frac{\Gamma\big(|j|+\frac{1-\alpha}{2}\big)}{\Gamma\big(|j|+\frac{1+\alpha}{2}\big)}{\bf e}_{l,j}(\varphi,\theta)\\ &
=\frac{2^{1-\alpha}\Gamma(\alpha)}{\Gamma(\frac{1-\alpha}{2})\Gamma(\frac{1+\alpha}{2})}\mathtt{W}(|j|,1-\alpha){\bf e}_{l,j}(\varphi,\theta)
\end{align}
where $\mathtt{W}$ was  defined in Lemma \ref{lem-asym}-(iii).
According to Lemma  \ref{Stirling-formula} the asymptotic of Wallis quotient  for large $j$ is given by
$$
\mathtt{W}(|j|,1-\alpha)=\frac{1}{|j|^\alpha}+O\bigg(\frac{1}{|j|^{2+\alpha}}\bigg).
$$
It follows that 
\begin{equation}\label{Frac-Mult}
|\textnormal D|^{-\alpha}{\bf e}_{l,j}=\frac{2^{1-\alpha}\Gamma(\alpha)}{\Gamma(\frac{1+\alpha}{2})\Gamma(\frac{1-\alpha}{2})}\bigg[\frac{1}{|j|^\alpha}+O\bigg(\frac{1}{|j|^{2+\alpha}}\bigg)\bigg]{\bf e}_{l,j}.
\end{equation}
In particular, the operator $|\textnormal D|^{-\alpha}$ acts like the standard fractional Laplacian $(-\Delta)^{-\frac\alpha2}$ and one gets
\begin{eqnarray*}
\|\textnormal D|^{-\alpha}h\|_{H^{s_1,s_2+\alpha}}^2&\le&C\sum_{(l,j)\in \Z^d\times \Z}\langle (l,j)\rangle^{2s_1}\langle j\rangle^{2(s_2+\alpha)}\langle j\rangle ^{-2\alpha}|h_{j,k}|^2\\
&\le&C\|h\|_{H^{s_1,s_2}}^2.
\end{eqnarray*}

\smallskip

{\bf{(iii)}} We can write by definition and change of variables 
\begin{align*}
\partial_\alpha^k|\textnormal D|^{-\alpha}{\bf e}_{l,j}(\varphi,\theta)
= a_j\, {\bf e}_{l,j}(\varphi,\theta)\quad\hbox{with}\quad 
a_j=\frac{1}{\pi}\bigintsss_{0}^{\pi} \frac{e^{2\ii j\eta}}{|\sin{\eta}|^{1-\alpha}}\ln^k\big(\sin{\eta} \big) d\eta.
\end{align*}
 Let $\delta\in(0,\frac12)$ and write
\begin{align*}
a_j=\frac{1}{2\ii j\pi}\bigintsss_{\sin\eta\geqslant\delta} \frac{\big(e^{2\ii j\eta}\big)^\prime}{|\sin{\eta}|^{1-\alpha}}\ln^k\big(\sin{\eta} \big) d\eta+\frac{1}{\pi}\bigintsss_{0\leqslant \sin\eta\leqslant\delta} \frac{e^{2\ii j\eta}}{|\sin{\eta}|^{1-\alpha}}\ln^k\big(\sin{\eta} \big) d\eta.
\end{align*}
Using integration by parts implies
\begin{align*}
|a_j|\lesssim |j|^{-1}\delta^{\alpha-1}{\ln^k(1/\delta)}+|j|^{-1}\bigintsss_{\sin\eta\geqslant\delta} \frac{\big|\ln^{k}\big(\sin{\eta} \big)\big|}{|\sin{\eta}|^{2-\alpha}} d\eta+\frac{1}{\pi}\bigintsss_{0\leqslant \sin\eta\leqslant\delta} \frac{\big|\ln^k\big(\sin{\eta} \big)\big|}{|\sin{\eta}|^{1-\alpha}}d\eta.
\end{align*}
Therefore using the fact that for any $\epsilon>0,$ $\displaystyle{\sup_{x\in(0,1)}|x|^{\epsilon}|\ln x|<\infty}$ we get
\begin{align*}
|a_j|\lesssim |j|^{-1}\delta^{\alpha-1-\epsilon}+\delta^{\alpha-\epsilon}.
\end{align*}
By taking $\delta=|j|^{-1}$ we deduce that
$$
|a_j|\lesssim |j|^{-\alpha+\epsilon}.
$$
This gives the desired result, that is,
\begin{eqnarray*}
\big\|\partial_\alpha^k|\textnormal D|^{-\alpha}h\big\|_{H^{s_1,s_2+\alpha-\epsilon}}^2
&\le&C\|h\|_{H^{s_1,s_2}}^2.
\end{eqnarray*}
The  proof of the lemma is now achieved.
\end{proof}

\section{Hamiltonian toolbox}

We shall reformulate   the equation \eqref{HS}  as  a perturbation of the Hamiltonian linear equation given by \eqref{linearized-op}. More precisely, we write \eqref{HS} in the form 
\begin{equation}\label{NLeq}
\partial_{t}r=\partial_{\theta}\, \mathrm{L}(\alpha)(r)+X_{P}(r)
\end{equation}
where $\mathrm{L}(\alpha)$  is defined in \eqref{defL} and  $X_{P}$ is the  Hamiltonian vector field generated by  the higher order terms $P\triangleq  H_{\geqslant 3}(  r )$, namely  
\begin{equation}\label{XP}
X_{P}(r)=(\Omega+V_{0,\alpha})\,\partial_\theta r-\partial_{\theta}\, \mathbb{K}_{0,\alpha} r- F_\alpha[r],
\end{equation}
where the non-local contribution $ F_\alpha[r]$ is introduced in  \eqref{eq},  the constant $V_{0,\alpha}$ is defined in \eqref{defTheta} and the operator $\mathbb{K}_{0,\alpha}$ is given by \eqref{k0}. \\
Since we are going to construct small quasi-periodic solutions with small amplitude then it is convenient to rescale the unknown as follows $ r \mapsto \varepsilon r $, leading to
\begin{align}\label{formaHep}
{\mathcal H}_{\varepsilon} (r) 
& \triangleq \varepsilon^{-2} H(\varepsilon r ) =
H_{\rm L} (r) +  \varepsilon P_\varepsilon  (r) 
\qquad {\rm where} \qquad P_\varepsilon \triangleq \varepsilon^{-3} P( \varepsilon r )
\end{align}
and  $ H_{\rm L}(r) $ is the quadratic  Hamiltonian in \eqref{defQH}. 
Therefore the equation \eqref{equation} writes
 \begin{align}\label{NLeq-eps}
\partial_tr&=\partial_\theta \, \mathrm{L}(\alpha)\,r+\varepsilon X_{P_\varepsilon}(r), \qquad {\rm with} \qquad X_{P_\varepsilon}(r)\triangleq\varepsilon^{-2} X_{P}(\varepsilon r).
\end{align}

\subsection{Action-angle coordinates}
We consider finitely many "tangential" sites
\begin{equation}\label{tangent-set}
{\mathbb S} \triangleq \{ n_1, \ldots, n_d \} ,\quad
{1 \leqslant n_1} < n_2 < \ldots < n_d
\end{equation}
and we introduce the set
\begin{equation}\label{tangent-set2}
  \quad \overline{\mathbb S}\triangleq \big\{ \pm n,\,n\in {\mathbb S}\big\}\quad\hbox{and}\quad {\mathbb S_0}\triangleq\overline{\mathbb S}\cup\{0\}.
\end{equation}
Denote
 the unperturbed tangential and normal frequency vectors by 
\begin{equation}\label{tan-nor}
{\omega}_{\textnormal{Eq}} \triangleq (\Omega_j (\alpha))_{j\in \mathbb{S}},\qquad {\Omega}(\alpha) \triangleq (\Omega_j (\alpha))_{j\in \mathbb{Z}\setminus {\mathbb{S}}_0}
\end{equation}
where $\Omega_j (\alpha) $ is given by \eqref{omega}.
We decompose the phase space  
$$
L^2_0 (\mathbb{T})\triangleq  \bigg\{ r(\theta) = 
\sum_{ j\in \mathbb{Z}\setminus\{0\}} r_j e^{\ii j\theta};\;\,   \overline{r_j}=r_{-j}, \;\, \| r \|_s^2 = \sum_{j \in \mathbb{Z}\setminus\{0\}}|r_j|^2  < + \infty \bigg\}
$$ as the direct sum 
\begin{equation}\label{decoacca}
\begin{aligned}
L^2_0 (\mathbb{T}) = \mathbb{H}_{\overline{\mathbb{S}}}\oplus \mathbb{H}^{\bot}_{\mathbb{S}_0}\quad\hbox{with}\quad
 \mathbb{H}_{\overline{\mathbb{S}}}&\triangleq \bigg\{v=\sum_{ j\in \overline{ \mathbb{S}}} r_j e^{\ii j\theta} ;\;   \overline{r_j}=r_{-j}  \bigg\}, \\
 \mathbb{H}^{\bot}_{\mathbb{S}_0}&\triangleq \bigg\{ z = 
\sum_{ j\in \mathbb{Z}\setminus\mathbb{S}_0} z_j e^{\ii j\theta}  \in L^2_0 (\T)   \bigg\} \, .
\end{aligned}
\end{equation}
The orthogonal projectors $\Pi_{\overline{\mathbb{S}}}$ and $ \Pi^\bot_{\mathbb S_0}$ are defined as follows, for $r\in L_0^2(\mathbb{T})$,
\begin{equation}\label{projectors-tan-normal}
 \Pi_{\overline{\mathbb S}}r\triangleq\sum_{ j\in \overline{\mathbb S}} r_j e^{\ii j\theta},\qquad \Pi^\bot_{\mathbb S_0}r\triangleq\sum_{ j\in \mathbb{Z}\setminus\mathbb{S}_0} r_j e^{\ii j\theta}.
\end{equation}
 On the finite dimensional space $ \mathbb{H}_{\overline{\mathbb{S}}}$ we shall use a local parametrization through  the action-angle variables.
Fix some positive  amplitudes $(\mathtt{a}_{j})_{j\in\mathbb{S}}\in(\mathbb{R}_{+}^{*})^{d}$ such that $\mathtt{a}_{-j}=\mathtt{a}_{j}$ and define the finite dimensional torus,
$$
\mathbf{T}^{d}\triangleq\Big\{ \theta\in\mathbb{R}\mapsto\sum_{j\in{\overline{\mathbb{S}}}} \mathtt{c}_j\mathtt{a}_j e^{\ii\,j\theta};\; \mathtt{c}_j\in\mathbb{C},\; |\mathtt{c}_j|=1,\; c_j=\overline{c_j}  \Big\}
$$ Now we shall  introduce action-angle variables $(\vartheta,I)=\big((\vartheta_j)_{j\in\overline{\mathbb{S}}},(I_j)_{j\in\overline{\mathbb{S}}}\big)$ allowing to describe the phase space around the linear  torus $\mathbf{T}^{d}$, as follows\begin{equation}\label{ham-syst-k}
\begin{cases}
r_j  =
\sqrt{\mathtt{a}_{j}^2+\frac{|j|}{2\pi}I_j}\,  e^{\ii \vartheta_j}, \quad\textnormal{if}& \quad  j\in \overline{\mathbb{S}},
\\
r_j = z_j, \quad\textnormal{if}& \quad  j\in \mathbb{Z}\setminus\mathbb{S}_0
\end{cases} 
\end{equation}
where
\begin{equation}
I_{-j}=I_j, \quad \vartheta_{-j}=-\vartheta_j, \quad \vartheta_j,I_j\in \mathbb{R}, \quad j\in \overline{\mathbb{S}}.
\end{equation}
Therefore we construct a local parameterization of the phase space around the linear torus $\mathbf{T}^{d}$ as  \begin{equation}\label{aacoordinates}
\begin{aligned}
 r(\theta) &= \mathbf{A}( \vartheta,I,z) \triangleq  
v (\vartheta,I)+ z   \quad
 {\rm where} \\
v (\vartheta, I) &\triangleq \sum_{j \in \overline{\mathbb{S}}}    
 \sqrt{\mathtt{a}_{j}^2+\tfrac{|j|}{2\pi}I_j}\,  e^{\ii (\vartheta_j+j\theta)}\quad\hbox{and} \quad
 z&=\sum_{ j\in \mathbb{Z}\setminus\mathbb{S}_0} r_j e^{\ii j\theta}.
\end{aligned}
\end{equation}
Remark that in the new  coordinates system, $(\vartheta,I)=(-{\omega}_{\textnormal{Eq}} t,0)$ corresponds to the  solution to the linear system \eqref{linearized-op} given by \eqref{q-r-l}  which  belongs to the torus  $\mathbf{T}^{d}.$ The symplectic  $2$-form in \eqref{sy2form-fourier2} becomes
\begin{equation}\label{sympl_form}
{\mathcal W} =  
\sum_{n \in {\mathbb S}} d\vartheta_n  \wedge  d I_n   +\frac{1}{2\ii}\sum_{j \in \Z \setminus \mathbb{S}_0}\frac{1}{j}dr_j\wedge dr_{-j}=
\Big(\sum_{n \in {\mathbb S}}d\vartheta_n  \wedge  d I_n   \Big)  \oplus {\mathcal W}_{| \mathbb{H}^{\bot}_{\mathbb{S}_0}} 
\end{equation}
where ${\mathcal W}_{| \mathbb{H}^{\bot}_{\mathbb{S}_0}}$ denotes the restriction of $\mathcal{W}$ to $ \mathbb{H}^{\bot}_{\mathbb{S}_0}$. 
We observe  that ${\mathcal W}  $ is an exact $ 2 $-form as 
$$
\mathcal{W} = d \Lambda 
$$
where $ \Lambda $ is the Liouville $ 1$-form
\begin{equation}\label{Lambda 1 form}
\Lambda_{( \vartheta, I, z)}[
 \widehat \vartheta, \widehat I, \widehat z] \triangleq-\sum_{n \in {\mathbb S}} I_n   \widehat\vartheta_n 
+ \frac12 ( \partial_\theta^{-1}  z, \widehat z )_{L^2 (\mathbb{T})} \, .    
\end{equation}
The Poisson bracket are given by
\begin{equation}\label{poisson-bracket}
\{F,G\}=\mathcal{W}(X_F,X_G)=\langle\nabla F, \mathbf{J}\nabla G\rangle
\end{equation}
where $\langle \cdot,\cdot\rangle$ is  the inner product, defined by
$$
\langle ( \vartheta_1,  I_1,  z_1),(\vartheta_2,I_2, z_2)\rangle=\vartheta_1\cdot\vartheta_2+  I_1\cdot I_2+(z_1,  z_2 )_{L^2 (\mathbb{T})}.
$$
The Poisson structure $\mathbf{J}$ corresponding to $\mathcal{W}$, defined by the identity \eqref{poisson-bracket} , is the unbounded operator 
$$
\mathbf{J}: (\vartheta,  I,  z)\mapsto (  I,- \vartheta, \partial_\theta  z).
$$ 

\smallskip

Now we shall  study the Hamiltonian system  generated by the Hamiltonian $ {\mathcal H}_\varepsilon   $ in \eqref{formaHep}, 
in  the action-angle and normal  coordinates 
$ 
(\vartheta, I, z) \in  \mathbb{T}^\nu \times \mathbb{R}^\nu \times  \mathbb{H}^{\bot}_{\mathbb{S}_0} \, . 
$ 
We consider the Hamiltonian $ H_{\varepsilon} (\vartheta, I, z )$ defined by 
\begin{equation}\label{Hepsilon}
H_{\varepsilon} =\mathcal{H}_{\varepsilon} \circ   \mathbf{A}
\end{equation}
where  
$  \mathbf{A} $ is the map defined in \eqref{aacoordinates}.  Since $ \mathrm{L}(\alpha) $  in \eqref{defL} preserves the subspace $ \mathbb{H}^{\bot}_{\mathbb{S}_0}$ then the quadratic Hamiltonian $ H_{\rm L}$ in \eqref{defQH} (see \eqref{defQH1})
in the variables $ (\vartheta, I, z) $ reads, up to a constant,
\begin{align}\label{QHAM}
H_N \circ  \mathbf{A} =  -\sum_{n\in\mathbb{S}} \, \Omega_n(\alpha)I_n+ \frac12  (  \mathrm{L}(\alpha)\, z, z )_{L^2} =   -{\omega}_{\textnormal{Eq}}\cdot I
+ \frac12  (  \mathrm{L}(\alpha)\, z, z )_{L^2}  
\end{align}
where $ {\omega}_{\textnormal{Eq}} \in \mathbb{R}^d $ is the unperturbed 
tangential frequency vector.
By \eqref{formaHep} and \eqref{QHAM}, 
the Hamiltonian $H_{\varepsilon} $ in \eqref{Hepsilon} reads
\begin{equation}\label{cNP}
\begin{aligned}
&  H_{\varepsilon} = 
 {\mathcal N} + \varepsilon \mathcal{ P}_{\varepsilon}  \qquad {\rm with}  \\
&
 {\mathcal N} \triangleq   -{\omega}_{\textnormal{Eq}}\cdot I  + \frac12  (  \mathrm{L}(\alpha)\, z, z )_{L^2},  
 \quad 
 \quad \mathcal{ P}_{\varepsilon} \triangleq   P_\varepsilon \circ  \mathbf{A}.  
\end{aligned}
\end{equation}
We look for an embedded invariant torus
\begin{equation}\label{rev-torus}
i :\mathbb{T}^d \rightarrow
\mathbb{R}^d \times \mathbb{R}^d \times  \mathbb{H}^{\bot}_{\mathbb{S}_0} 
\,, \quad \varphi \mapsto i(\varphi)\triangleq (  \vartheta(\varphi), I(\varphi), z(\varphi)),  
\end{equation}

of the Hamiltonian vector field 
$$ 
X_{H_{\eps}} \triangleq 
(\partial_I H_{\eps} , -\partial_\vartheta H_{\eps} , \Pi_{\mathbb{S}_0}^\bot
 \partial_\theta \nabla_{z} H_{\eps} ) 
 $$ 
filled by quasi-periodic solutions with Diophantine frequency 
vector $\omega$. 
Note that for the  value   $\varepsilon=0$, the Hamiltonian system 
$$\omega\cdot\partial_\varphi i (\varphi) = X_{H_\eps} ( i (\varphi))$$   possesses, for any value of the parameter $\alpha\in [0,1)$ , the invariant torus  $$i_{\textnormal{flat}}(\varphi)=(\varphi,0,0),
$$
provided that $\omega=-{\omega}_{\textnormal{Eq}}.$
Now we consider the family  of Hamiltonians,
\begin{equation}\label{H alpha}
\begin{aligned}
H_\eps^\tc \triangleq {\mathcal N}_\tc +\varepsilon  {\mathcal P}_{\eps} \, , \quad  {\mathcal N}_\tc \triangleq  
\tc \cdot I 
+ \tfrac12 ( \mathrm{L}(\alpha)\, z, z)_{L^2}, 
\end{aligned}
\end{equation}
which depend on the constant vector $ \tc \in \mathbb{R}^d $. For the value $\tc=-{\omega}_{\textnormal{Eq}}$ we have $H_\eps^\tc= H_\eps$.
The parameter $\tc$ is introduced in order to control the average in the $I$-component of the linearized equations.
We look for zeros of the nonlinear operator
\begin{align}
 {\mathcal F} (i, \tc ) 
& \triangleq   {\mathcal F} (i, \tc, \omega,\alpha,  \eps )  \triangleq \omega\cdot\partial_\varphi i (\varphi) - X_{H_\eps^\tc} ( i (\varphi))
=  \omega\cdot\partial_\varphi i (\varphi) -  (X_{{\mathcal N}_\tc}  +\varepsilon   X_{\mathcal{P}_{\eps}})  (i(\varphi) ) \nonumber \\
&  \triangleq  \left(
\begin{array}{c}
\omega\cdot\partial_\varphi \vartheta (\varphi) - \tc -  \varepsilon  \partial_I \mathcal{P}_{\eps} ( i(\varphi)   )  \\
\omega\cdot\partial_\varphi I (\varphi) + \varepsilon  \partial_\vartheta \mathcal{P}_{\eps}( i(\varphi)  )  \\
\omega\cdot\partial_\varphi z (\varphi) 
-  \partial_\theta {\rm L}(\alpha) z(\varphi) - \varepsilon  \partial_\theta \nabla_z \mathcal{P}_{\eps} ( i(\varphi) )   
\end{array}
\right)  \label{operatorF} 
\end{align}
where $  \vartheta (\varphi) - \varphi $ is a  $ (2 \pi)^d $-periodic function.
 We denote by 
\begin{equation}\label{periodic-comp}
\mathfrak{I}(\varphi)\triangleq i(\varphi)-(\varphi,0,0)=
( \vartheta (\varphi)-\varphi, I(\varphi), z(\varphi))
\end{equation}
the periodic component of the torus $\varphi\mapsto i(\varphi)$.  Note that the involution $ {\mathfrak S} $,  described in Proposition \ref{prop-hamilt}-${\rm (iii)}$, becomes 
\begin{equation}\label{rev_aa}
	\widetilde{\mathfrak S}:  (   \vartheta,I, z)\mapsto ( -\vartheta,I,  {\mathfrak S} z ).
\end{equation}
Moreover,  we  can easily check that the Hamiltonian vector field $X_{H_\eps^\tc}$ is reversible with respect  $\widetilde {\mathfrak S} $.
Thus, it is natural to  look for reversible solutions of $ {\mathcal F}(i, \tc) = 0 $,  namely satisfying, 
\begin{equation}\label{parity solution}
\vartheta(-\varphi) = - \vartheta (\varphi) \, , \,  \
I(-\varphi) = I(\varphi) \, , \, \ 
z (- \varphi ) = ( {\mathfrak S} z)(\varphi) \, . 
\end{equation}
The norm of the periodic component of the embedded torus 
\eqref{periodic-comp}
is 
\begin{equation}\label{def:norma-cp}
\|  {\mathfrak I}  \|_s^{q,\overline \gamma} 
\triangleq \| \vartheta-\textnormal{Id} \|_{s}^{q,\overline \gamma} +  \| I  \|_{s}^{q,\overline \gamma} 
+  \| z \|_s^{q,\overline \gamma}\, . 
\end{equation}
We emphasize that we have used for the tangential variables, depending only on the  variable $\varphi$,  the norm \eqref{Def-weit-norm} associated to functions with two variables $(\varphi,\theta)$. In this case the norm $H^s(\T^{d+1})$ coincides with the norm $H^s(\T^d).$  
Concerning the  index regularity  $q$ it is an arbitrary  integer that will be fixed only at the end of the paper in Section \ref{Section 6.2} devoted to the final Cantor set measure. There we shall make the choice 
\begin{equation}
q\triangleq q_0+1.
\end{equation}
with $q_0$  being  the  non-degeneracy index  provided by Proposition \ref{lemma transversality}, which
only depends on the linear unperturbed frequencies.

\subsection{Hamiltonian regularity}

We  shall analyze in this section some  regularity aspects of the Hamiltonian vector field $X_{P}$  described by   \eqref{XP} as well as  the rescaled one associated to $\mathcal{P}_{\varepsilon}$  given by \eqref{cNP}. The arguments are classical and the computations turn out to be long and straightforward. For this reason we will only sketch the proofs only  for some essential points.
The first main result  reads as follows.
\begin{lemma}\label{lemma estimates vector field XP}
Let $(\kappa,q,s_{0},s)$ satisfying \eqref{initial parameter condition}.
{  {There exists $\varepsilon_{0}\in(0,1]$ such that if
\begin{equation}\label{small-C15}
\| r\|_{s_0+2}^{q,\kappa}\leqslant\varepsilon_{0},
\end{equation}
then} the vector field $X_{P}$, given by  \eqref{XP}, satisfies the following estimates :
\begin{enumerate}
\item $\| X_{P}(r)\|_{s}^{q,\kappa}\lesssim \| r\|_{s+1}^{q,\kappa}.$
\item $\| d_{r}X_{P}(r)[h]\|_{s}^{q,\kappa}\lesssim\|h\|_{s+1}^{q,\kappa}+\| r\|_{s+2}^{q,\kappa}\|h\|_{s_{0}+1}^{q,\kappa}.$
\item 
$\| d_r^{2}X_{P}(r)[h,h]\|_{s}^{q,\kappa}\lesssim \|h\|_{s+2}^{q,\kappa}\|h\|_{s_{0}+1}^{q,\kappa}+\| r\|_{s+2}^{q,\kappa}\big(\|h\|_{s_{0}+1}^{q,\kappa}\big)^2$.
\end{enumerate}}
\end{lemma}

\begin{proof}
{\rm (i)} Recall from \eqref{XP} that
 \begin{equation*}
X_{P}(r)=(\Omega+V_{0,\alpha})\,\partial_\theta r-\partial_{\theta}\, \mathbb{K}_{0,\alpha} r- F_\alpha[r].
\end{equation*}
where we have by   \eqref{eq}, 
\begin{equation}\label{Falpha}
F_\alpha[r](\varphi,\theta)= \displaystyle \frac{C_\alpha}{2\pi} \bigintssss_{0}^{2\pi}\frac{\partial_{\theta\eta}^{2}\left(R(\varphi,\theta)R(\varphi,\eta)\sin(\eta-\theta)\right)}{A_r^\frac\alpha2(\varphi,\theta,\eta)}d\eta. 
\end{equation}
and from  \eqref{defV0} and \eqref{fract1}, 
\begin{align}\label{K0alpha}
\mathbb{K}_{0,\alpha} &=2^{-\alpha}C_\alpha |\textnormal D|^{\alpha-1}.
\end{align}
We shall first estimate the term $\partial_\theta \mathbb{K}_{0,\alpha} r$.  Using Lemma \ref{Laplac-frac0}-{\rm (iii)}  combined with Leibniz formula we obtain, for all $\epsilon\in (0,1)$ and 
  $\lambda=(\omega,\alpha) \in \mathcal{O}$ ,
\begin{align*}
\|\partial_\lambda^k \partial_\theta [ |\textnormal D|^{\alpha-1}r(\lambda, \centerdot)]\|_{H^{s-|k|}}
&\lesssim_k  \sum_{j\leqslant k}\|\partial_\lambda^{j}  |\textnormal D|^{\alpha-1}\partial_\lambda^{k-j} (\partial_\theta r)(\lambda, \centerdot)\|_{H^{s-|k|}}
\\
&\lesssim_k  \sum_{j\leqslant k}\| \partial_\lambda^{k-j}   r(\lambda, \centerdot)\|_{H^{s+\alpha+\epsilon-|k|}}
.
\end{align*}
Combining the last estimate together  with \eqref{K0alpha} and Definition \ref{Def-WS} we find
\begin{equation}\label{estimate kernel equilibrium}
\|\partial_{\theta}\mathbb{K}_{0,\alpha} r\|_{s}^{q,\kappa}\lesssim \| r\|_{s+\overline\alpha+\epsilon}^{q,\kappa}.
\end{equation}
Next we shall move to the estimate of $F_\alpha[r]$.
From \eqref{A} we deduce the identity
\begin{align*}
A_r(\varphi,\theta,\eta)&=\sin^2\left(\tfrac{\eta-\theta}{2}\right) \left[\left(\tfrac{R(\varphi,\eta)-R(\varphi,\theta)}{\sin\big(\frac{\eta-\theta}{2}\big)}\right)^2+4R(\varphi,\eta)R(\varphi,\theta)\right]\\
&=\sin^2\left(\tfrac{\eta-\theta}{2}\right) \left[\left(\tfrac{R(\varphi,\eta)-R(\varphi,\theta)}{\tan\big(\frac{\eta-\theta}{2}\big)}\right)^2+\big(R(\varphi,\eta)+R(\varphi,\theta)\big)^2\right]. 
\end{align*} 
It follows that
\begin{align}\label{Ar-decomp}
C_\alpha A_r^{-\frac\alpha2}(\varphi,\theta,\eta)&=\left|\sin\left(\tfrac{\eta-\theta}{2}\right)\right|^{-\alpha} {\mathscr  A}_{r,\alpha}(\varphi,\theta,\eta)
\end{align} 
with
\begin{align}\label{Ar-22}
{\mathscr  A}_{r,\alpha}(\varphi,\theta,\eta)
 &=C_\alpha\left[\left(\tfrac{R(\varphi,\eta)-R(\varphi,\theta)}{\tan\big(\frac{\eta-\theta}{2}\big)}\right)^2+\big(R(\varphi,\eta)+R(\varphi,\theta)\big)^2\right]^{-\frac\alpha2}.
\end{align}
Applying Lemmata \ref{Compos-lemm}-\ref{Law-prodX1}-\ref{lem-Reg1} combined with   the smallness condition \eqref{small-C1},
we obtain
\begin{align}\label{ScrA15}
\sup_{\eta\in\T}\|{\mathscr  A}_{r,\alpha}(\cdot,\centerdot,\eta+\centerdot)-2^{-\alpha}C_\alpha\|_{s}^{q,\kappa}
&\leqslant C\|r\|_{s+1}^{q,\kappa}.
\end{align}
According to \eqref{Ar-decomp} combined with a change of variables in the expression of $F_\alpha[r]$, one gets
\begin{equation}\label{F-alpha-r}
F_\alpha[r](\varphi,\theta)=\displaystyle \frac{1}{2\pi} \bigintssss_{0}^{2\pi}\frac{\mathscr{B}_{r,\alpha}^0(\varphi,\theta,\eta)}{|\sin\big(\frac{\eta}{2}\big)|^\alpha}d\eta 
\end{equation}
with 
\begin{align*}
\mathscr{B}_{r,\alpha}^0(\varphi,\theta,\eta)&\triangleq g_0(\varphi,\theta,\eta){\mathscr  A}_{r,\alpha}(\varphi,\theta,\eta+\theta),\\
g_0(\varphi,\theta,\eta)&\triangleq \big[\partial_{\theta}R(\varphi,\theta)(\partial_{\eta}R)(\varphi,\eta+\theta)+R(\varphi,\theta)R(\varphi,\eta+\theta)\big]\sin(\eta)\\ &\,\,+\big[\partial_{\theta}R(\varphi,\theta)R(\varphi,\eta+\theta)-(\partial_{\eta}R)(\varphi,\eta+\theta)R(\varphi,\theta)\big]\cos(\eta).
\end{align*}
We can easily check by symmetry that 
\begin{align}\label{RHH-1}
\nonumber F_\alpha[0](\varphi,\theta)&= \frac{2^{-\alpha} C_\alpha}{2\pi} \bigintssss_{0}^{2\pi}\frac{\sin \eta}{|\sin\big(\frac{\eta}{2}\big)|^\alpha}d\eta \\
&=0.
\end{align}
By the law products in Lemma  \ref{Law-prodX1}, the composition law in Lemma \ref{Compos-lemm} and  \eqref{small-C15}  we infer
\begin{align}\label{estimate partialthetaeta}
\sup_{\eta\in\T}\| g_0(\cdot,\centerdot,\eta)-{\sin(\eta)}\|_{s}^{q,\kappa}&\lesssim\| r\|_{s+1}^{q,\kappa}.
\end{align}
Hence using once again the  law products of Lemma \ref{Law-prodX1}, \eqref{small-C15}, \eqref{ScrA15} and \eqref{estimate partialthetaeta} we find
\begin{align}\label{estimate B0}
\sup_{\eta\in\T}\| \mathscr{B}_{r,\alpha}^0(\cdot,\centerdot,\eta)-\mathscr{B}_{0,\alpha}^0(\cdot,\centerdot,\eta)\|_{s}^{q,\kappa}&=\sup_{\eta\in\T}\| \mathscr{B}_{r,\alpha}^0(\cdot,\centerdot,\eta)-2^{-\alpha}C_\alpha{\sin(\eta)}\|_{s}^{q,\kappa}\notag\\
 &\lesssim\| r\|_{s+1}^{q,\kappa}.
\end{align}
Finally, from    \eqref{F-alpha-r}, \eqref{RHH-1}, \eqref{estimate B0}  and  straightforward computations lead to
\begin{align}
\|F_{\alpha}[r]\|_{s}^{q,\kappa}&\lesssim\bigintsss_{0}^{2\pi} \frac{\left(1+|\ln(\sin(\eta/2))|^q\right)}{|\sin(\frac\eta2)|^{\frac12}}  \sup_{\eta\in\T}\|\mathscr{B}_{r,\alpha}^0(\cdot,\centerdot,\eta)-\mathscr{B}_{0,\alpha}^0\|_{s}^{q,\kappa}d\eta\notag\\ &\lesssim\|r\|_{s+1}^{q,\kappa}.\label{estimateFalpha}
\end{align} 
Putting together this with \eqref{estimate kernel equilibrium} and \eqref{XP} achieves the proof of the first point. 

\smallskip

\textbf{(ii)} 
 From Proposition \ref{lin-eq-r}  we have
$$d_{r}F_{\alpha}[r]h=\partial_{\theta}\left({V}_{r,\alpha}h-\mathbb{K}_{r,\alpha}h\right).$$
Thus, by \eqref{XP} we get
\begin{equation}\label{expression of dXP(r)}
d_{r}X_{P}(r)h=-\partial_{\theta}\left((V_{r,\alpha}-V_{0,\alpha})h\right)+\partial_{\theta}\left(\mathbb{K}_{r,\alpha}h-\mathbb{K}_{0,\alpha}h\right).
\end{equation}
Combining \eqref{Ar-decomp} with \eqref{defV} and using a change of variables
  we find
  \begin{equation}\label{V-eq80}
V_{r,\alpha}(\varphi,\theta)=\displaystyle \frac{1}{2\pi} \bigintssss_{0}^{2\pi}\frac{\mathscr{B}_{r,\alpha}^1(\varphi,\theta,\eta)}{|\sin\big(\frac{\eta}{2}\big)|^\alpha}d\eta 
\end{equation}
with
\begin{align*}
\mathscr{B}_{r,\alpha}^1(\varphi,\theta,\eta)&\triangleq g_1(\varphi,\theta,\eta){\mathscr  A}_{r,\alpha}(\varphi,\theta,\eta+\theta),\\
g_1(\varphi,\theta,\eta)&\triangleq \frac{(\partial_{\eta}R)(\varphi,\eta+\theta)\sin(\eta)+ R(\varphi,\eta+\theta)\cos(\eta)}{R(\varphi,\theta)}.
\end{align*}
Proceeding in a similar way to \eqref{estimate partialthetaeta} and \eqref{estimate B0} we obtain
\begin{align*}
\sup_{\eta\in\T}\| \mathscr{B}_{r,\alpha}^1(\cdot,\centerdot,\eta)-\mathscr{B}_{0,\alpha}^1(\cdot,\centerdot,\eta)\|_{s}^{q,\kappa} &\lesssim\| r\|_{s+1}^{q,\kappa}.
\end{align*}
Then, straightforward computations lead to
\begin{align}\label{estimate V}
\|V_{r,\alpha}-V_{0,\alpha}\|_{s}^{q,\kappa}\lesssim&\bigintsss_{0}^{2\pi} \frac{\left(1+|\ln(\sin(\eta/2))|^q\right)}{|\sin(\frac\eta2)|^{\frac12}} \displaystyle\sup_{\eta\in\T} \|\mathscr{B}_{r,\alpha}^1(\cdot,\centerdot,\eta)-\mathscr{B}_{0,\alpha}^1(\cdot,\centerdot,\eta)\|_{s}^{q,\kappa}d\eta\notag\\
&\lesssim \displaystyle\sup_{\eta\in\T} \|\mathscr{B}_{r,\alpha}^1(\cdot,\centerdot,\eta)-\mathscr{B}_{0,\alpha}^1(\cdot,\centerdot,\eta)\|_{s}^{q,\kappa}\notag d\eta\notag\\
&\quad\lesssim \|r\|_{s+1}^{q,\kappa}.
\end{align}
Using the  product laws  in Lemma \ref{Law-prodX1} together with \eqref{estimate V} and  \eqref{small-C15},  we obtain
\begin{align*}
\|\partial_{\theta}\left((V_{r,\alpha}-V_{0,\alpha})h\right)\|_s^{q,\kappa}&\lesssim \|V_{r,\alpha}-V_{0,\alpha}\|_{s+1}^{q,\kappa}\|h\|_{s_{0}}^{q,\kappa}+\|V_{r,\alpha}-V_{0,\alpha}\|_{s_0}^{q,\kappa}\|h\|_{s+1}^{q,\kappa}\\
&\lesssim \|r\|_{s+2}^{q,\kappa}\|h\|_{s_{0}}^{q,\kappa}+\|h\|_{s+1}^{q,\kappa}.
\end{align*}
Coming back to the second term of \eqref{expression of dXP(r)}. From \eqref{Ar-decomp} and \eqref{defK} we have the following identity
\begin{align}\label{K0-22}
\mathbb{K}_{r,\alpha} h (\varphi,\theta)&= \frac{ 1}{2\pi }\bigintsss_{0}^{2\pi} \frac{h(\varphi,\eta+\theta){\mathscr  A}_{r,\alpha}(\varphi,\theta,\eta+\theta)}{|\sin(\frac{\eta}{2})|^{\alpha}} d\eta.
\end{align}
In a similar way to \eqref{estimate V} we obtain
\begin{align*}
\|\partial_{\theta}\left(\mathbb{K}_{r,\alpha}h-\mathbb{K}_{0,\alpha}h\right)\|_{s}^{q,\kappa}&\lesssim
\|\mathbb{K}_{r,\alpha}h-\mathbb{K}_{0,\alpha}h\|_{s+1}^{q,\kappa}\\ &\lesssim \|r\|_{s+2}^{q,\kappa}\|h\|_{s_0}^{q,\kappa}+\|h\|_{s+1}^{q,\kappa}.
\end{align*}
This concludes the proof of \textbf{(ii)}.

\smallskip

\textbf{(iii)} Differentiating in $r$ the identity  \eqref{expression of dXP(r)} yields,
\begin{equation}\label{Second-diff}
d_{r}^{2}X_{P}(r)[{h},{h}]=\partial_{\theta}\left(d_{r}\big(\mathbb{K}_{r,\alpha}h\big)[h]\right)-\partial_{\theta}\left(\big(d_{r}{V}_{r,\alpha}[h]\big)h\right).
\end{equation}
For the first term of the right-hand side we write, by \eqref{K0-22} and  \eqref{Ar-22}, 
\begin{align*}
d_{r}\big(\mathbb{K}_{r,\alpha}h\big)[h] (\varphi,\theta)&= \frac{ 1}{2\pi }\bigintsss_{0}^{2\pi} \frac{h(\varphi,\eta+\theta)d_{r}{\mathscr  A}_{r,\alpha}[h](\varphi,\theta,\eta+\theta)}{|\sin(\frac{\eta}{2})|^{\alpha}} d\eta.
\end{align*}
with
\begin{align*}
d_r{\mathscr  A}_{r,\alpha}[h](\varphi,\theta,\eta)
 &=2C_\alpha\Bigg[\tfrac{\frac{h(\varphi,\eta)}{R(\varphi,\eta)}-\frac{h(\varphi,\theta)}{R(\varphi,\theta)}}{\tan\big(\frac{\eta-\theta}{2}\big)}\tfrac{R(\varphi,\eta)-R(\varphi,\theta)}{\tan\big(\frac{\eta-\theta}{2}\big)} +\left(\tfrac{h(\varphi,\eta)}{R(\varphi,\eta)}+\tfrac{h(\varphi,\theta)}{R(\varphi,\theta)}\right)\big(R(\varphi,\eta)+R(\varphi,\theta)\big)\Bigg]\\ &\quad \times \left[\left(\tfrac{R(\varphi,\eta)-R(\varphi,\theta)}{\tan\big(\frac{\eta-\theta}{2}\big)}\right)^2+\big(R(\varphi,\eta)+R(\varphi,\theta)\big)^2\right]^{-\frac\alpha2-1}.
\end{align*}
Applying Lemmata \ref{Compos-lemm}-\ref{Law-prodX1}-\ref{lem-Reg1} combined with   the smallness condition \eqref{small-C15},
we obtain
\begin{align*}
\sup_{\eta\in\T}\|d_r{\mathscr  A}_{r,\alpha}[h](\cdot,\centerdot,\eta+\centerdot)\|_{s}^{q,\kappa}
&\lesssim {\|h\|_{s+1}^{q,\kappa}+\|h\|_{s_0+1}^{q,\kappa}\|r\|_{s+1}^{q,\kappa}}.
\end{align*}
Then, direct computations yield
\begin{align}
\|d_{r}\big(\partial_{\theta}\mathbb{K}_{r,\alpha}h\big)[h]\|_{s}^{q,\kappa}&\lesssim
\|d_{r}\big(\mathbb{K}_{r,\alpha}h\big)[h]\|_{s+1}^{q,\kappa}\notag\\ &\lesssim \|h\|_{s+2}^{q,\kappa}\|h\|_{s_{0}+1}^{q,\kappa}+\| r\|_{s+2}^{q,\kappa}\big(\|h\|_{s_{0}+1}^{q,\kappa}\big)^2.\label{drkr}
\end{align}
Finally, differentiating \eqref{V-eq80} with respect to $r$ in the direction $h$ gives
  \begin{align*}
d_rV_{r,\alpha}[h](\varphi,\theta)&=\displaystyle \frac{1}{2\pi} \bigintssss_{0}^{2\pi}\frac{d_r g_1[h](\varphi,\theta,\eta){\mathscr  A}_{r,\alpha}(\varphi,\theta,\eta+\theta)}{|\sin\big(\frac{\eta}{2}\big)|^\alpha}d\eta\\
&+\displaystyle \frac{1}{2\pi} \bigintssss_{0}^{2\pi}\frac{g_1(\varphi,\theta,\eta)d_r{\mathscr  A}_{r,\alpha}[h](\varphi,\theta,\eta+\theta)}{|\sin\big(\frac{\eta}{2}\big)|^\alpha}d\eta.
\end{align*}
Arguing as for \eqref{estimateFalpha} and \eqref{drkr} we get
\begin{align*}
\|\partial_{\theta}\left(d_rV_{r,\alpha}[h]h\right)\|_s^{q,\kappa}&\lesssim  \|h\|_{s+2}^{q,\kappa}\|h\|_{s_{0}+1}^{q,\kappa}+\| r\|_{q,s+2}^{q,\kappa}\big(\|h\|_{s_{0}+1}^{\kappa}\big)^2.
\end{align*}
Combining the last estimate with \eqref{drkr} and \eqref{Second-diff} we conclude the proof of {\rm (iii)}.
\end{proof}

Finally we state tame estimates for the composition operator induced by the Hamiltonian vector
field 
$$
X_{\mathcal{P}_{\varepsilon}}=(\partial_{I}\mathcal{P}_{\varepsilon},-\partial_{\vartheta}\mathcal{P}_{\varepsilon},\Pi_{\mathbb{S}_0}^{\perp}\partial_{\theta}\nabla_{z}\mathcal{P}_{\varepsilon}),$$
where $\mathcal{P}_{\varepsilon}$ is defined in \eqref{cNP}. %
\begin{lemma}\label{tame estimates for the vector field XmathcalPvarepsilon}
{Let $(\kappa,q,s_{0},s)$ satisfying \eqref{initial parameter condition}.
 There exists $\varepsilon_0\in(0,1)$ such that if 
$$
\|\mathfrak{I}\|_{{s_0+2}}^{q,\kappa}\leqslant \varepsilon_0,
$$ 
then  the  perturbed Hamiltonian vector field $X_{\mathcal{P}_{\varepsilon}}$ satisfies the following tame estimates,
\begin{enumerate}
\item $\| X_{\mathcal{P}_{\varepsilon}}(i)\|_{s}^{q,\kappa}\lesssim 1+\|\mathfrak{I}\|_{s+2}^{q,\kappa}.$
\item $\big\| d_{i}X_{\mathcal{P}_{\varepsilon}}(i)[\,\widehat{i}\,]\big\|_{s}^{q,\kappa}\lesssim \|\,\widehat{i}\,\|_{s+2}^{q,\kappa}+\|\mathfrak{I}\|_{s+2}^{q,\kappa}\|\,\widehat{i}\,\|_{s_{0}+2}^{q,\kappa}.$
\item $\big\| d_{i}^{2}X_{\mathcal{P}_{\varepsilon}}(i)[\,\widehat{i},\widehat{i}\,]\big\|_{s}^{q,\kappa}\lesssim \|\,\widehat{i}\,\|_{s+2}^{q,\kappa}\|\,\widehat{i}\,\|_{s_{0}+2}^{q,\kappa}+\|\mathfrak{I}\|_{s+2}^{q,\kappa}\left(\|\,\widehat{i}\,\|_{s_{0}+2}^{q,\kappa}\right)^{2}.$
\end{enumerate}}
\end{lemma}
\begin{proof}
The proof is a straightforward  and slightly long. It is based on Taylor expansion in $r$ around the origin  of the Hamiltonian vector field $X_P$ detailed in \eqref{XP}. Then the estimates for  the rescaled vector field in \eqref{NLeq-eps} follow from    Lemma \ref{lemma estimates vector field XP},  the law products and the  following fact 
\begin{align}\label{v-theta1}
 \forall a,b\in\mathbb{N}^{d},\quad |a|+|b|\leqslant 3,\quad \|\partial_{\vartheta}^{a}\partial_{I}^{b}v(\vartheta,I)\|_{s}^{q,\kappa}&\lesssim  1+\|\mathfrak{I}\|_{s}^{q,\kappa},
\end{align}
for all $\vartheta,I\in H^s(\mathbb{T}^d,\mathbb{R}^d)$, $\|\vartheta \|_{s_0},\| I \|_{s_0}\leqslant 1$. For more details we can refer to  Lemma 5.1 of \cite{BertiMontalto} where we explore a similar situation.
\end{proof}

\section{Berti-Bolle approach for the approximate inverse}\label{sec:Approximate-inverse}
In order to apply a  modified Nash-Moser scheme,
we need to construct an  approximate
right inverse of the linearized operator  at an arbitrary torus $i_0$ close to the flat one  and an arbitrary vector-valued function   $\tc_0:\mathcal{O}\to \mathbb{R}^d$,
\begin{equation}\label{Lineqrized-op-F}
d_{(i, \tc)} {\mathcal F}(i_0, \tc_0 )[\widehat \imath \,, \widehat \tc ] =
\omega \cdot \partial_\varphi \widehat \imath - d_i X_{H_\eps^{\tc_0}} ( i_0 (\varphi) ) [\widehat \imath ] - (\widehat \tc,0, 0 ),
\end{equation}
where $ {\mathcal F}  $ is the nonlinear operator defined  in \eqref{operatorF}. 

We follow the strategy presented by Berti and Bolle in \cite{BB13},  which reduces the search of an approximate  inverse of \eqref{Lineqrized-op-F} to 
 the task of inverting the operator in the normal direction. Their main idea consists in 
  linearizing  the functional ${\mathcal F}$ around an isotropic torus  sufficiently close to the torus $i_0$ and    make use  of a convenient  symplectic transformation  in such a way the linearized equations become "approximately" decoupled through a trianglular system in the action-angle components and the normal ones.
  
  \smallskip

The main novelty related to this part is  to bypass the use of isotropic torus and    work directly with the original one  $i_0$. Although the transformation $G_0$ introduced in \eqref{trasform-sympl}  is not symplectic and the nonlinear Hamiltonian structure is no longer  preserved, the  conjugation of the linearized operator via the linear change of variables $DG_0(G_0^{-1}
(i_0))$
leads a triangular system with small errors of size $Z\triangleq{\mathcal F}(i_0, \tc_0)$, see Proposition \ref{Prop-Conjugat} and Lemma \ref{L-Invert0}. Yet we  emphasize that the Hamiltonian structure of the original system is of  paramount importance, namely  in the  construction of the transformation $G_0$  leading to the final triangular system.
Therefore, along this section we shall focus on this linear change of variable. The estimates that we shall  perform  are
very similar to those in \cite{Baldi-berti,BB13, BertiMontalto}, with minor differences due to  the accumulation of different
extra errors induced by the isotropic torus as for example in the various Cantor sets emerging during the reducibility scheme. Thus, we often  refer to these papers for some technical points in order to avoid here too much details.

\smallskip

 We shall first fix  
some notations and definitions that will be used in this section.
Let $\mathtt{H}$ be an Hilbert space equipped with the inner  product $\langle \cdot,\cdot \rangle_{\mathtt{H}}$. Given a linear  operator $A\in  {\mathcal L}({\mathbb R}^d, \mathtt{H})$
 we define the transposed operator $A^{\top}:\mathtt{H} \to\mathbb{R}^d$  by the duality relation
\begin{align}\label{Def-dua}
 \big\langle A^{\top}u,v\big\rangle_{\mathbb{R}^d} =\big\langle u,  Av\big\rangle_{\mathtt{H}}, \quad \quad \forall \, u\in {\mathtt{H}}\, , \,  v\in\mathbb{R}^d.
\end{align}
In accordance with the notation introduced in \eqref{periodic-comp} we denote by
 \begin{equation}
i_0(\varphi)\triangleq ( \vartheta_0 (\varphi), I_0(\varphi), z_0(\varphi)), \qquad\mathfrak{I}_0(\varphi)\triangleq i_0(\varphi)-(\varphi,0,0).
\end{equation}
We shall first  assume the following hypothesis:  
the map $\lambda\mapsto \mathfrak{I}_0 (\varphi; \lambda) $ 
is $q$-times differentiable with respect to the parameters $\lambda=(\omega, \alpha) \in \mathcal{O}\subset \RR^d \times [\underline \alpha, \overline \alpha] $
and 
 there exists $\varepsilon_{0}\in(0,1)$ such that if
\begin{equation}\label{ansatz 0}
 \|  \mathfrak{I}_0\|_{s_0+2}^{q,\kappa}+ \| \tc_0 - \omega\|^{q,\kappa}\leqslant\varepsilon_{0}
\end{equation}
where  $\varepsilon_{0}$ is assumed  to be small enough. We point out that this assumption is enough to provide tame estimates on the tangential part. Later on, a stronger assumption will be required  for  the normal  contribution, see \eqref{small-C3}. 
\smallskip

We introduce  the  diffeomorpshim 
$ G_0 : (\phi, y, w) \to (\vartheta, I, z)$ of the phase space $\mathbb{T}^d \times \mathbb{R}^d \times  \mathbb{H}^{\bot}_{\mathbb{S}_0}$ defined by
\begin{equation}\label{trasform-sympl}
\begin{pmatrix}
\vartheta \\
I \\
z
\end{pmatrix} \triangleq G_0 \begin{pmatrix}
\phi \\
y \\
w
\end{pmatrix} \triangleq 
\begin{pmatrix}
\!\!\!\!\!\!\!\!\!\!\!\!\!\!\!\!\!\!\!\!\!\!\!\!\!\!\!\!\!\!\!\!\!
\!\!\!\!\!\!\!\!\!\!\!\!\!\!\!\!\!\!\!\!
 \vartheta_0(\phi) \\
I_0 (\phi) + L_1(\phi) y + L_2(\phi) w \\
\!\!\!\!\!\!\!\!\!\!\!\!\!\!\!\!\!\!\!\!\!\!\!\!\!\!\!\!\!
\!\!\!\!\!\!\!\!\!\!\!\!
z_0(\phi) + w
\end{pmatrix} 
\end{equation}
where
\begin{equation}\label{l1l2}
\begin{aligned}
L_1(\phi)&\triangleq  [\partial_\varphi \vartheta_0(\varphi)]^{-\top},
\\
 L_2(\phi) &\triangleq - [(\partial_\vartheta \widetilde{z}_0)(\vartheta_0(\phi))]^\top \partial_\theta ^{-1}\quad \textnormal{with} \quad  \widetilde{z}_0 (\vartheta) \triangleq z_0 (\vartheta_0^{-1} (\vartheta)).
\end{aligned}
\end{equation}
In the new coordinates, $i_0$ becomes the trivial
embedded torus $(\phi,y,w)=(\varphi,0,0)$, namely
$$
G_0(\varphi,0,0)=i_0(\varphi).
$$
In what follows we shall denote by $${\mathtt u}_0(\varphi)\triangleq G_0^{-1}
(i_0)(\varphi)=(\varphi,0,0)$$ the trivial
torus and  ${\mathtt u}=(\phi, y,w)$ the  coordinates induced by $G_0$ in \eqref{trasform-sympl}. We shall also denote by 
$$
 \widetilde{G}_0 ({\mathtt u}, \tc) \triangleq  \big( G_0 ({\mathtt u}), \tc \big) 
 $$ 
the diffeomorphism with the identity on the $ \tc $-component.

We quantify how an embedded torus $i_0(\mathbb{T})$ is approximately invariant for the Hamiltonian vector field $X_{H_\varepsilon^\tc}$ in terms of the "error function"     
\begin{equation} \label{def Z}
Z(\varphi) \triangleq  (Z_1, Z_2, Z_3) (\varphi) \triangleq {\mathcal F}(i_0, \tc_0) (\varphi) =
\omega \cdot \partial_\varphi i_0(\varphi) - X_{H^{\tc}_\varepsilon}(i_0(\varphi), \tc_0) \, .
\end{equation}

\subsection{Linear change of variables}
The main task in this subsection is to  conjugate the linearized operator $d_{i,\tc} {\mathcal F} (i_0,{\mathtt c}_0)$ in \eqref{Lineqrized-op-F}, via the linear change of variables
\begin{equation}\label{DG0}
D G_0({\mathtt u}_0(\varphi))
\begin{pmatrix}
\widehat \phi \, \\
\widehat y \\
\widehat w
\end{pmatrix} 
=
\begin{pmatrix}
\partial_\varphi \vartheta_0(\varphi) & 0 & 0 \\
\partial_\varphi I_0(\varphi) &L_1(\varphi) & 
L_2(\varphi) \\
\partial_\varphi z_0(\varphi) & 0 & I
\end{pmatrix}
\begin{pmatrix}
\widehat \phi \, \\
\widehat y \\
\widehat w
\end{pmatrix} ,
\end{equation}
to a triangular system up to a remainder of order $O(Z)$. As it was shown in \cite{BB13} this transformation is symplectic up to errors of type $Z$. Our main results reads as follows.
\begin{proposition}\label{Prop-Conjugat}
Under the linear change of variables $D G_0({\mathtt u}_0)$ the linearized operator $d_{i,\tc} {\mathcal F} (i_0,{\mathtt c}_0)$ is transformed into
\begin{align}\label{bbL}
& [D G_0({\mathtt u}_0)]^{-1}  d_{(i,\tc)} {\mathcal F} (i_0,\tc_0) D\widetilde G_0({\mathtt u}_0)
 [\widehat \phi, \widehat y, \widehat w, \widehat \tc ]
 = \mathbb{D} [\widehat \phi, \widehat y, \widehat w, \widehat \tc ]+\mathbb{E} [\widehat \phi, \widehat y, \widehat w] 
\end{align}
where
\begin{enumerate} 
\item the operator $ \mathbb{D}$ has the triangular form
\begin{align*}
 \mathbb{D} [\widehat \phi, \widehat y, \widehat w, \widehat \tc ]\triangleq\left(
\begin{array}{c}
\omega\cdot\partial_\varphi  \widehat\phi-\big[K_{20}(\varphi) \widehat y+K_{11}^\top(\varphi) \widehat w+L_1^\top (\varphi)\widehat \tc\big]
\\
\omega\cdot\partial_\varphi  \widehat y+\mathcal{B}(\varphi) \widehat \tc \\
 \omega\cdot\partial_\varphi \widehat w-\partial_\theta\big[K_{11}(\varphi) \widehat y+K_{02}(\varphi)\widehat w -L_2^{\top}(\varphi) \widehat\tc   \big] 
\end{array}
\right),
\end{align*}
$\mathcal{B}(\varphi)$ and   $K_{20}(\varphi) $ are   $d \times d$ real matrices, 
\begin{align*}
\mathcal{B}(\varphi)& \triangleq[\partial_\varphi \vartheta_0(\varphi)]^\top \partial_\varphi I_0(\varphi) L_1^\top(\varphi) +[\partial_\varphi z_0(\varphi)]^{\top} L_2^\top(\varphi) ,
\\
K_{20}(\varphi)&\triangleq \varepsilon L_1^\top(\varphi) ( \partial_{II} \mathcal{P}_\varepsilon)(i_0(\varphi))  L_1(\varphi) \, ,
\end{align*}
$K_{02}(\varphi)$ is a linear self-adjoint operator of $  \mathbb{H}^{\bot}_{\mathbb{S}_0}$, given by 
\begin{align*}
   K_{02}(\varphi)& \triangleq ( \partial_{z}\nabla_z H_\varepsilon^\tc) (i_0(\varphi))  +\varepsilon L_2^\top(\varphi) ( \partial_{II} \mathcal{P}_\varepsilon)(i_0(\varphi)) L_2(\varphi) 
\\ &\quad+\varepsilon L_2^\top(\varphi)( \partial_{zI} \mathcal{P}_\varepsilon) (i_0(\varphi))     + \varepsilon  (\partial_I\nabla_z \mathcal{P}_\varepsilon) (i_0(\varphi))  L_2(\varphi),
\end{align*}
and $K_{11}(\varphi)  \in {\mathcal L}({\mathbb R}^d,   \mathbb{H}^{\bot}_{\mathbb{S}_0})$, 
\begin{align*}
K_{11}(\varphi)&\triangleq \varepsilon  L_2^\top(\varphi)  ( \partial_{II} \mathcal{P}_\varepsilon)(i_0(\varphi))L_1(\varphi) +\varepsilon ( \partial_I\nabla_z \mathcal{P}_\varepsilon) (i_0(\varphi)) L_1(\varphi) ,
\end{align*}
\item the remainder $\mathbb{E} $ is given by
\begin{align*}
\mathbb{E} [\widehat \phi, \widehat y, \widehat w] &\triangleq
 [D G_0({\mathtt u}_0)]^{-1}    \partial_\varphi Z(\varphi) \widehat\phi  \\  &\quad + \left(
\begin{array}{c}
0 
\\
 \mathcal{A}(\varphi)\big[K_{20}(\varphi) \widehat y+K_{11}^\top(\varphi) \widehat w\big]-R_{10}(\varphi) \widehat y -R_{01}(\varphi) \widehat w 
 \\
0   
\end{array}
\right)
\end{align*}
where $\mathcal{A}(\varphi)$ and   $R_{10}(\varphi) $ are   $d \times d$ real matrices, 
\begin{align*}
\mathcal{A}(\varphi)& \triangleq[\partial_\varphi \vartheta_0(\varphi)]^\top \partial_\varphi I_0(\varphi)-[\partial_\varphi I_0(\varphi)]^\top \partial_\varphi \vartheta_0(\varphi)  +[\partial_\varphi z_0(\varphi)]^{\top} \partial_\theta ^{-1} \partial_\varphi z_0(\varphi),
\\
R_{10}(\varphi)&\triangleq  [\partial_\varphi Z_1(\varphi)]^{\top}  L_1(\varphi), 
\end{align*}
and $R_{01}(\varphi)\in {\mathcal L}( \mathbb{H}^{\bot}_{\mathbb{S}_0},{\mathbb R}^d)$, 
\begin{align*}
R_{01}(\varphi)&\triangleq -[\partial_\varphi Z_1(\varphi)]^{\top} L_2(\varphi)+ [\partial_\varphi  Z_3(\varphi)]^{\top} \partial_\theta^{-1}
 \, .
\end{align*}
\end{enumerate}
\end{proposition}
\begin{remark}
A priori, it is not clear from the expression of $\mathcal{A}$ that the error term $\mathbb{E}$ is of order $O(Z)$.
This is the purpose of Lemma \ref{lem:est-akj} where we quantify the size of $\mathcal{A}$ in terms of
the error function $Z$ defined in \eqref{def Z}.
\end{remark}
\begin{proof}
Under the  map $G_0$, the nonlinear operator $\mathcal{F}$  in \eqref{operatorF}  is transformed into
\begin{align}
 {\mathcal F} (G_0({\mathtt u}(\varphi)),\tc) 
= \omega\cdot\partial_\varphi \big( G_0({\mathtt u}(\varphi))\big)- X_{H_\varepsilon^\tc} \big(G_0({\mathtt u}(\varphi))\big).  
 \label{operatorF3} 
\end{align}
Differentiating \eqref{operatorF3}  at $({\mathtt u}_0,\tc_0)$  in the direction  $(\widehat {\mathtt u}, \widehat \tc)$ gives
\begin{align}\label{df-comp}
&d_{({\mathtt u},\tc)} ({\mathcal F} \circ G_0)({\mathtt u}_0,\tc_0)
[(\widehat {\mathtt u}, \widehat \tc )](\varphi)
=  \omega\cdot\partial_\varphi \big( DG_0({\mathtt u}_0)\widehat {\mathtt u} \big)- \partial_\phi\big[ X_{H_\varepsilon^\tc} \big(G_0({\mathtt u}(\varphi))\big)\big]_{{\mathtt u}={\mathtt u}_0} \widehat\phi
\\ &\qquad\qquad\qquad\quad -  \partial_y\big[ X_{H_\varepsilon^\tc} \big(G_0({\mathtt u}(\varphi))\big)\big]_{{\mathtt u}={\mathtt u}_0} \widehat y- \partial_w\big[ X_{H_\varepsilon^\tc} \big(G_0({\mathtt u}(\varphi))\big) \big]_{{\mathtt u}={\mathtt u}_0} \widehat w
- \left(
\begin{array}{c}
\widehat\tc  \\
0
 \\
 0  
\end{array}
\right).\nonumber
\end{align}
From the expression of $DG_0({\mathtt u}_0)$ in \eqref{DG0},  we obtain 
 \begin{align}\label{dgu0}
 \omega\cdot\partial_\varphi \big( DG_0({\mathtt u}_0)[\widehat {\mathtt u} ](\varphi)\big)
  &=DG_0({\mathtt u}_0) \, \omega\cdot\partial_\varphi \widehat {\mathtt u} + \partial_\varphi\big(\omega\cdot\partial_\varphi  i_0 \big)\widehat\phi \\ & \quad+\left(
\begin{array}{c}
0
\\
(\omega\cdot\partial_\varphi  L_1(\varphi))\widehat y+\big(\omega\cdot\partial_\varphi  L_2(\varphi)\big) \widehat w
 \\
0
\end{array}
\right). \notag
\end{align} 
In view of  \eqref{l1l2} and \eqref{def Z} we have
\begin{equation}\label{l1}
\begin{aligned}
\omega\cdot\partial_\varphi  L_1(\varphi)&=-[\partial_\varphi \vartheta_0(\varphi)]^{-\top}(\omega\cdot\partial_\varphi  [\partial_\varphi \vartheta_0(\varphi)]^{\top}) [\partial_\varphi \vartheta_0(\varphi)]^{-\top}\\
&=-[\partial_\varphi \vartheta_0(\varphi)]^{-\top} \Big(\big[\partial_\varphi Z_1(\varphi)\big]^{\top} + \big[ \partial_\varphi\big( ( \partial_{I} H_\varepsilon^\tc) (i_0(\varphi))\big) \big]^{\top}\Big) [\partial_\varphi \vartheta_0(\varphi)]^{-\top}.
\end{aligned}
\end{equation}
By \eqref{l1l2} we can easily check that
\begin{equation}\label{z00}
\partial_\varphi z_0(\varphi)=(\partial_\vartheta \widetilde{z}_0)(\vartheta_0(\varphi))\partial_\phi \vartheta_0(\varphi)
\end{equation}
and thus, we may write the operator $L_2(\varphi)$ in term of the matrix $L_1(\varphi)$,
\begin{equation}\label{l2}
L_2(\varphi)=-[\partial_\varphi \vartheta_0(\varphi)]^{-\top}[\partial_\varphi z_0(\varphi)]^{\top}\partial_\theta^{-1}=-L_1(\varphi)[\partial_\varphi z_0(\varphi)]^{\top}\partial_\theta^{-1}.
\end{equation}
Then, by \eqref{l1}, \eqref{l2} we have
\begin{align}
\omega\cdot\partial_\varphi  L_2(\varphi)
 &=[\partial_\varphi \vartheta_0(\varphi)]^{-\top}(\omega\cdot\partial_\varphi  [\partial_\varphi \vartheta_0(\varphi)]^{\top}) [\partial_\varphi \vartheta_0(\varphi)]^{-\top}[\partial_\varphi z_0(\varphi)]^{\top}\partial_\theta^{-1}\notag
 \\
&  -[\partial_\varphi \vartheta_0(\varphi)]^{-\top}  [\partial_\varphi    (\omega\cdot\partial_\varphi z_0)(\varphi)]^{\top} \partial_\theta^{-1}\notag
 \end{align}
and from \eqref{def Z}  we get
 \begin{align}\label{l22}
\omega\cdot\partial_\varphi  L_2(\varphi)
 &=-[\partial_\varphi \vartheta_0(\varphi)]^{-\top}\Big(\big[\partial_\varphi Z_1(\varphi)\big]^{\top} + \big[ \partial_\varphi\big( ( \partial_{I} H_\varepsilon^\tc) (i_0(\varphi))\big) \big]^{\top}\Big) L_2(\varphi)  \notag
\\ &\quad  -[\partial_\varphi \vartheta_0(\varphi)]^{-\top}   \Big( \big[\partial_\varphi  Z_3(\varphi)\big]^{\top}+\big[ \partial_\varphi\big( ( \nabla_{z} H_\varepsilon^\tc) (i_0(\varphi))\big) \big]^{\top}\Big) \partial_\theta^{-1}.
 \end{align}
Putting together \eqref{dgu0}, \eqref{l1} and  \eqref{l22} we conclude that
 \begin{align}\label{omdphdg}
 \omega\cdot\partial_\varphi \big( DG_0({\mathtt u}_0)[\widehat {\mathtt u}](\varphi) \big)&=DG_0({\mathtt u}_0) \, \omega\cdot\partial_\varphi \widehat {\mathtt u} + \partial_\varphi\big(\omega\cdot\partial_\varphi  i_0 \big)\widehat\phi \nonumber
\\ & -  \left(
\begin{array}{c}
0
\\
~ [\partial_\varphi \vartheta_0(\varphi)]^{-\top}\big[\mathcal{C}_I(\varphi) L_1(\varphi)  +R_{10}(\varphi)  \big] \widehat y
 \\
0
\end{array}
\right)\nonumber\\ & -  \left(
\begin{array}{c}
0
\\
~  [\partial_\varphi \vartheta_0(\varphi)]^{-\top}\big[\mathcal{C}_I(\varphi) L_2(\varphi)+ \mathcal{C}_z(\varphi)   \partial_\theta^{-1}  + R_{01}(\varphi) \big] \widehat w
 \\
0
\end{array}
\right),
\end{align} 
where $R_{10}(\varphi)$ and  $R_{01}(\varphi)$ are given by {\rm (ii)} 
 and
\begin{align}
\mathcal{C}_I(\varphi) &\triangleq \big[ \partial_\varphi\big( ( \partial_{I} H_\varepsilon^\tc) (i_0(\varphi))\big) \big]^{\top}\notag\\ 
&=[\partial_\varphi I_0(\varphi)]^\top( \partial_{II} H_\varepsilon^\tc)(i_0(\varphi))+ [\partial_\varphi \vartheta_0(\varphi)]^\top ( \partial_{\vartheta I} H_\varepsilon^\tc)(i_0(\varphi)) \notag\\ &\quad+ [\partial_\varphi z_0(\varphi)]^\top ( \partial_{I} \nabla_z H_\varepsilon^\tc)(i_0(\varphi)),\label{CI}
 \\ 
 \mathcal{C}_z(\varphi) & \triangleq \big[ \partial_\varphi\big( ( \nabla_{z} H_\varepsilon^\tc) (i_0(\varphi))\big) \big]^{\top}\notag\\
 &=[\partial_\varphi I_0(\varphi)]^\top( \partial_{zI} H_\varepsilon^\tc)(i_0(\varphi))+ [\partial_\varphi \vartheta_0(\varphi)]^\top ( \partial_{z\vartheta } H_\varepsilon^\tc)(i_0(\varphi))  \notag\\ &\quad+ [\partial_\varphi z_0(\varphi)]^\top ( \partial_{z} \nabla_z H_\varepsilon^\tc)(i_0(\varphi)).\label{Cz}
\end{align}

On the other hand, in view \eqref{operatorF} and  \eqref{trasform-sympl} we obtain 
\begin{align}
 \partial_\phi\big[ X_{H_\varepsilon^\tc} \big(G_0({\mathtt u}(\varphi))\big)\big]_{{\mathtt u}={\mathtt u}_0} \widehat\phi &
= \partial_\varphi\big[ X_{H_\varepsilon^\tc} (i_0(\varphi)))\big] \widehat\phi \label{dx1},
\\
  \partial_y\big[ X_{H_\varepsilon^\tc} \big(G_0({\mathtt u}(\varphi))\big)\big]_{{\mathtt u}={\mathtt u}_0} \widehat y &
= \left(
\begin{array}{c}
( \partial_{II} H_\varepsilon^\tc)(i_0(\varphi))  L_1(\varphi) \widehat y 
\\
-( \partial_{I\vartheta} H_\varepsilon^\tc) (i_0(\varphi))  L_1(\varphi) \widehat y 
 \\
   \partial_\theta \big[ ( \partial_I\nabla_z H_\varepsilon^\tc) (i_0(\varphi))  L_1(\varphi) \widehat y \big]     
\end{array}
\right), \label{dx2}
\\
 \partial_w\big[ X_{H_\varepsilon^\tc} \big(G_0({\mathtt u}(\varphi))\big) \big]_{{\mathtt u}={\mathtt u}_0} \widehat w &
=\left(
\begin{array}{c}
( \partial_{II} H_\varepsilon^\tc)(i_0(\varphi))L_2(\varphi) \widehat w+( \partial_{zI} H_\varepsilon^\tc) (i_0(\varphi))  \widehat w 
\\
-( \partial_{I\vartheta} H_\varepsilon^\tc) (i_0(\varphi)) L_2(\varphi) \widehat w- ( \partial_{z\vartheta} H_\varepsilon^\tc) (i_0(\varphi)) \widehat w
 \\
   \partial_\theta \big[( \partial_I\nabla_z H_\varepsilon^\tc) (i_0(\varphi)) L_2(\varphi) \widehat w  + ( \partial_{z}\nabla_z H_\varepsilon^\tc) (i_0(\varphi)) \widehat w   \big]  
\end{array}
\right). \label{dx3}
\end{align}
Plugging \eqref{omdphdg}, \eqref{dx1}, \eqref{dx2} and \eqref{dx3} into \eqref{df-comp} we find
 \begin{align}\label{dfcirg}
&d_{({\mathtt u},\tc)} ({\mathcal F} \circ G_0)({\mathtt u}_0,\tc_0)
[(\widehat {\mathtt u}, \widehat \tc )]= DG_0({\mathtt u}_0) \, \omega\cdot\partial_\varphi \widehat {\mathtt u} + \partial_\varphi\big[{\mathcal F}(i_0(\varphi)) \big] \widehat\phi  \nonumber
  \\ 
  &   + \left(
\begin{array}{c}
-( \partial_{II} H_\varepsilon^\tc)(i_0(\varphi))  L_1(\varphi) \widehat y 
\\
( \partial_{I\vartheta} H_\varepsilon^\tc) (i_0(\varphi))  L_1(\varphi)\widehat y  -[\partial_\varphi \vartheta_0(\varphi)]^{-\top} [\mathcal{C}_I(\varphi) L_1(\varphi) +R_{10}(\varphi) ]\widehat y 
 \\
-  \partial_\theta ( \partial_I\nabla_z H_\varepsilon^\tc) (i_0(\varphi))  L_1(\varphi) \widehat y      
\end{array}
\right) 
 \nonumber
\\ 
&+\left(
\begin{array}{c}
-( \partial_{II} H_\varepsilon^\tc)(i_0(\varphi))L_2(\varphi) \widehat w-( \partial_{zI} H_\varepsilon^\tc) (i_0(\varphi))  \widehat w 
\\
\big[( \partial_{I\vartheta} H_\varepsilon^\tc) (i_0(\varphi)) L_2(\varphi) + ( \partial_{z\vartheta} H_\varepsilon^\tc) (i_0(\varphi))\big] \widehat w
 \\
 -  \partial_\theta \big[( \partial_I\nabla_z H_\varepsilon^\tc) (i_0(\varphi)) L_2(\varphi) \widehat w  + \partial_{\theta} ( \partial_{z}\nabla_z H_\varepsilon^\tc) (i_0(\varphi)) \widehat w \big]     
\end{array}
\right)\nonumber
\\ &-  \left(
\begin{array}{c}
0
\\
~ [\partial_\varphi \vartheta_0(\varphi)]^{-\top}\big[\mathcal{C}_I(\varphi) L_2(\varphi)\widehat w+ \mathcal{C}_z(\varphi)   \partial_\theta^{-1} \widehat w +R_{01}(\varphi) \widehat w\big] 
 \\
0
\end{array}
\right)- \left(
\begin{array}{c}
\widehat\tc  \\
0
 \\
 0  
\end{array}
\right).
\end{align}
According to  \eqref{DG0} and using the identities \eqref{z00} and \eqref{l2},  the inverse of the linear operator $D G_0({\mathtt u}_0)$  is given by
\begin{align}\label{dg-1}
& [D G_0({\mathtt u}_0)]^{-1} 
= 
 \begin{pmatrix}
[\partial_\varphi \vartheta_0(\varphi)]^{-1} & 0 & 0 \\
-\mathcal{B}(\varphi) & [\partial_\varphi \vartheta_0(\varphi)]^\top & 
[\partial_\varphi z_0(\varphi)]^\top \partial_\theta^{-1}   \\
(\partial_\vartheta \widetilde{z}_0)(\vartheta_0(\varphi))    & 0 & I
\end{pmatrix}
\end{align}
where $\mathcal{B}(\varphi)$  is given by {\rm(i)}.
Applying $[D G_0({\mathtt u}_0)]^{-1}$ to \eqref{dfcirg} and using the identities \eqref{CI}, \eqref{Cz} and the fact that
\begin{align}
\mathcal{B}(\varphi)
&=\mathcal{A}(\varphi) [\partial_\varphi \vartheta_0(\varphi)]^{-1}+[\partial_\varphi I_0(\varphi)]^\top,\label{idAB}
\end{align}
where $\mathcal{A}(\varphi)$  is defined in {\rm (ii)},
 we obtain
\begin{align*}
& [D G_0({\mathtt u}_0)]^{-1}  d_{({\mathtt u},\tc)} ({\mathcal F} \circ G_0)({\mathtt u}_0,\tc_0)
[\widehat {\mathtt u}, \widehat \tc ]
=  \omega\cdot\partial_\varphi \widehat {\mathtt u}+  [D G_0({\mathtt u}_0)]^{-1}   \partial_\varphi\big[{\mathcal F}(i_0(\varphi)) \big]  \widehat\phi 
\\
 &  + \left(
\begin{array}{c}
-K_{20}(\varphi)\widehat y
\\
\mathcal{A}(\varphi) K_{20}(\varphi)  \widehat y -R_{10}(\varphi) \widehat y  \\
-\partial_\theta K_{11}(\varphi) \widehat y       
\end{array}
\right)
+ \left(
\begin{array}{c}
-K_{11}^\top(\varphi) \widehat w 
\\
\mathcal{A}(\varphi) K_{11}^\top(\varphi) \widehat w -R_{01}(\varphi) \widehat w  \\
-\partial_\theta K_{02}(\varphi)\widehat w   
\end{array}
\right)
+ \left(
\begin{array}{c}
- L_1^\top(\phi) \widehat \tc
\\
\mathcal{B}(\varphi) \widehat\tc  \\
\partial_\theta L_2^\top(\varphi) \widehat\tc      
\end{array}
\right),
\end{align*}
with
\begin{align*}
K_{20}(\varphi)&\triangleq L_1^\top(\varphi) ( \partial_{II} H_\varepsilon^\tc)(i_0(\varphi))  L_1(\varphi) \, ,
\\
K_{11}(\varphi)&\triangleq L_2^\top(\varphi)  ( \partial_{II} H_\varepsilon^\tc)(i_0(\varphi))L_1(\varphi) + ( \partial_I\nabla_z H_\varepsilon^\tc) (i_0(\varphi)) L_1(\varphi)\, , 
   \\ 
   K_{02}(\varphi)& \triangleq ( \partial_{z}\nabla_z H_\varepsilon^\tc) (i_0(\varphi))  +L_2^\top(\varphi) ( \partial_{II} H_\varepsilon^\tc)(i_0(\varphi)) L_2(\varphi)
+L_2^\top(\varphi)( \partial_{zI} H_\varepsilon^\tc) (i_0(\varphi))     \\ &+   (\partial_I\nabla_z H_\varepsilon^\tc) (i_0(\varphi))  L_2(\varphi). \nonumber 
\end{align*}
Finally, by \eqref{H alpha} we conclude the desired identity and this ends the proof of  Proposition~\ref{Prop-Conjugat}.
\end{proof}
The next aim  is to  estimate  the induced composition operators $[DG_0({\mathtt u}_0)]^{\pm 1}$, given by \eqref{DG0}, \eqref{dg-1} and the coefficients   $R_{10}$,  $ R_{01}$, $K_{20}$, $ K_{11}$, $ K_{11}^\top$ and ${\mathcal A}$ defined in Proposition~\ref{Prop-Conjugat}. More precisely, arguing as in Lemmata  5.6-5.7  in \cite{BertiMontalto}, by using    law product of Lemma \ref{Law-prodX1}-(ii)  and the smallness condition  \eqref{ansatz 0}, we obtain the following result.
\begin{lemma} \label{lem:est-G} 
Let $(q,d,\kappa,s,s_0)$ as in \eqref{initial parameter condition}, then the  following assertions hold true.
\begin{enumerate}
\item The operator  $DG_0({\mathtt u}_0)$ and $[DG_0({\mathtt u}_0)]^{- 1} $ satisfy for all $\widehat{{\mathtt u}}=(\widehat{\phi},\widehat{y},\widehat{w})$,
 \begin{align*}
&\| [DG_0({\mathtt u}_0)]^{\pm 1}[\widehat{{\mathtt u}}]\|_{s}^{q,\kappa}
\lesssim\|\widehat{{\mathtt u}}\|_{s}^{q,\kappa}+\|\mathfrak{I}_{0}\|_{s+1}^{q,\kappa}\|\widehat{{\mathtt u}}\|_{s_0}^{q,\kappa}\, . 
\end{align*}
\item  The operators $R_{10}$ and  $ R_{01}$ satisfy the estimates
\begin{align*}
& \| R_{10} \widehat y \|_s^{q, \overline \gamma} \lesssim_s \| Z \|_{s+ 1 }^{q, \kappa} 
+ \| Z \|_{s_0 + 1 }^{q, \kappa} \big( \|\mathfrak{I}_{0}\|_{s+1}^{q,\kappa}\| \widehat y \|_{s_0 }^{q,\kappa}+\| \widehat y \|_{s }^{q,\kappa}\big),
\\
& \| R_{01} \widehat w  \|_s^{q, \overline \gamma} \lesssim_s  \| Z \|_{s+ 1 }^{q, \kappa}  \| \widehat w \|_{s_0 }^{q, \kappa}
+ \| Z \|_{s_0 + 1 }^{q, \kappa}\big( \|\mathfrak{I}_{0}\|_{s+1}^{q,\kappa}\| \widehat w \|_{s_0 }^{q,\kappa}+\| \widehat w \|_{s }^{q,\kappa}\big)
\, . 
\end{align*}

\item The matrix $K_{20}$ and the operators $ K_{11}, K_{11}^\top$ satisfy the estimates
\begin{align*}
& \| K_{20}  \|_s^{q, \overline \gamma} \lesssim_s \varepsilon \big( 1 + \| \mathfrak{I}_{0}\|_{s +2
}^{q, \overline \gamma} \big) \, ,  \\ 
& \| K_{11} \widehat y \|_s^{q, \overline \gamma} 
\lesssim_s \varepsilon \big(\| \widehat y \|_{s+2}^{q, \overline \gamma}
+ \| \mathfrak{I}_{0} \|_{s +2}^{q, \overline \gamma}  
\| \widehat y \|_{s_0+2}^{q, \overline \gamma} \big) \, , \\ 
&  \| K_{11}^\top \widehat w \|_s^{q, \overline \gamma}
\lesssim_s \varepsilon \big(\| \widehat w \|_{s + 2}^{q, \overline \gamma}
+  \| \mathfrak{I}_{0} \|_{s + 2}^{q, \overline \gamma}
\| \widehat w \|_{s_0 + 2}^{q, \overline \gamma} \big)\, . 
\end{align*}
\item The matrices $\mathcal{A}$ and $\mathcal{B}$ satisfy the estimates
\begin{align*}
& \| \mathcal{A}  \|_s^{q, \overline \gamma} +\| \mathcal{B}  \|_s^{q, \overline \gamma}\lesssim_s  \| \mathfrak{I}_{0} \|_{s+ 1 }^{q, \kappa} .
\end{align*}

\end{enumerate}
\end{lemma}
\subsection{Defect of the symplectic structure}
In this subsection we shall prove that the matrix  $\mathcal{A}$, defined in  Proposition \ref{Prop-Conjugat}-{\rm (ii)},  is zero at an exact solution on some Cantor like set, up to an exponentially small remainder. 
Recall that the coefficients of the matrix $\mathcal{A}$ are given by 
\begin{align}\label{def:akj}
\mathcal{A}_{kj}(\varphi)
&= \partial_{\varphi_k} I_0(\varphi)\cdot\partial_{\varphi_j} \vartheta_0(\varphi)
-\partial_{\varphi_k} \vartheta_0(\varphi)\cdot\partial_{\varphi_j} I_0(\varphi)  
+\big(\partial_\theta ^{-1} \partial_{\varphi_k}  z_0(\varphi) ,\partial_{\varphi_j}  z_0(\varphi)\big)_{L^2}.
\end{align}
We point out that those coefficients coincide with the "lack of isotropy coefficients" introduced in  Lemma 5 of \cite{ BB13} where the following identitity is proved 
\begin{align}\label{Akjz}
\omega\cdot\partial_\varphi \mathcal{A}_{jk}(\varphi)=&  {\mathcal W}\big( \partial_\varphi Z(\varphi) \underline{e}_k ,  \partial_\varphi i_0(\varphi)  \underline{e}_j \big) 
+ 
{\mathcal W} \big(\partial_\varphi i_0(\varphi) \underline{e}_k , \partial_\varphi Z(\varphi) \underline{e}_j \big) 
\\
&= \partial_{\varphi_k} Z_2(\varphi)\cdot\partial_{\varphi_j} \vartheta_0(\varphi)
-\partial_{\varphi_k} Z_1(\varphi) \cdot\partial_{\varphi_j} I_0(\varphi)
+\big(\partial_\theta ^{-1} \partial_{\varphi_k} Z_3(\varphi),\partial_{\varphi_j}  z_0(\varphi)\big)_{L^2}\nonumber
\\
&\quad+ \partial_{\varphi_k} I_0(\varphi)\cdot\partial_{\varphi_j} Z_2(\varphi)
-\partial_{\varphi_k} \vartheta_0(\varphi)\cdot\partial_{\varphi_j} Z_1(\varphi) 
+\big(\partial_\theta ^{-1} \partial_{\varphi_k}  z_0(\varphi) ,\partial_{\varphi_j} Z_3(\varphi)\big)_{L^2}.\nonumber
\end{align}
Here $(\underline{e}_1,\ldots,\underline{e}_d)$ denotes the canonical basis of $\mathbb{R}^d$.  
In order  to solve \eqref{Akjz} we  need to discuss some elementary results on the invertibility of Fourier multiplier operator in the presence of a small divisor problem.
Given  $\kappa \in(0,1]$ and $\tau_1>0$,    we introduce  the Cantor set  
  \begin{equation}\label{DC tau gamma}
 \mathtt {DC}(\kappa, \tau_1) \triangleq \Big\{ \omega \in \RR^d \, : \, 
|\omega \cdot l | \geqslant \tfrac{\kappa}{\langle l \rangle^{\tau_1}}, \quad \forall\,  l\in\mathbb{Z}^{d}\backslash\{0\}\Big\}.
\end{equation}
If  $\lambda\in\mathtt {DC}(\kappa, \tau_1)$ then, for any smooth function $h:\mathcal{O}\times\T^{d}\to \RR$ with zero average, the equation $ \omega\cdot\partial_\varphi u=h$
has a periodic solution $h:\mathbb{T}^{d}\to \RR$ given by
 $$
 u(\lambda,\varphi)=-\ii\sum_{l\in\Z^{d}\backslash\{0\}}\frac{h_{l}(\lambda)}{\omega\cdot l}{\bf{e}}_{l}(\varphi)\quad {\rm where}\quad h=\sum_{l\in\Z^{d}\backslash\{0\}}h_{l}(\lambda) {\bf{e}}_{l}, \quad {\bf e}_l(\varphi)\triangleq e^{\ii l\cdot \varphi}.
 $$
For all $\omega\in\mathcal{O}$  we define the smooth extension of $u$ by
\begin{align}\label{Extend001}
(\omega\cdot\partial_\varphi)_{\textnormal{ext}}^{-1}h\triangleq-\ii\sum_{l\in\Z^{d}\backslash\{0\}}\frac{\chi\big((\omega\cdot l)\kappa^{-1}\langle l\rangle^{\tau_1} \big)h_{l}(\lambda)}{\omega\cdot l}{\bf{e}}_{l},
\end{align}
where $\chi\in\mathscr{C}^\infty(\RR,\RR)$ is an even positive cut-off function  such that 
\begin{equation}\label{chi-def-1} 
\chi(\xi)=\left\{ \begin{array}{ll}
0\quad \hbox{if}\quad |\xi|\leqslant\frac13,&\\
 1\quad \hbox{if}\quad |\xi|\geqslant\frac12.
  \end{array}\right.
\end{equation}
Notice that this operator is well-defned in the whole set of parameters $\mathcal{O}$ and coincides with
the formal inverse of $(\omega\cdot\partial_\varphi)^{-1}$ when the frequency $\omega$ belongs to $ \mathtt {DC}(\kappa, \tau_1)$. On the other hand, we note that at each iterative step of the Nash-Moser iteration -and correspondingly for the reduction of the linearized operator in Section \ref{Reducibility of the linearized operator}- we only require that the frequency vector $\omega\in \mathbb{R}^d$ satisfies finitely
many non-resonance diophantine conditions. More precisely we assume at
the $n$-th step that $\omega$ belongs to
the truncated Cantor set
  \begin{equation}\label{DC tau gamma N}
\mathtt {DC}_{N_n} (\kappa, \tau_1) \triangleq \Big\{ \omega \in \RR^d;\;\, 
|\omega \cdot l | \geqslant \tfrac{\kappa}{\langle l \rangle^{\tau_1}}, \;\,   \forall  |l |\leqslant N_n \Big\} \, 
\end{equation}
where for any $n\in\N\cup\{-1\}$ the sequence $(N_n)$ is defined by
\begin{equation}\label{definition of Nm}
N_{-1}\triangleq 1,\quad \forall n\in\mathbb{N},\quad N_{n}\triangleq N_{0}^{\left(\frac{3}{2}\right)^{n}}\quad \hbox{with} \quad N_0\geqslant2.
\end{equation}
We shall make use of the following  result. 
\begin{lemma}\label{L-Invert0}
Let $\kappa\in(0,1], q\in\N$. Then the following holds true.
\begin{enumerate}
\item For any $ s\geqslant q$ we have
$$
\big\|(\omega\cdot\partial_\varphi)_{\textnormal{ext}}^{-1}h\big\|_{s}^{q,\kappa}\leqslant C\kappa^{-1}  \|h\|_{s+\tau_1(q+1)}^{q,\kappa}.
$$
\item For any $n\in\NN$ and  any $\omega\in \mathtt {DC}_{N_n} (\kappa, \tau_1)$ we have
$$
(\omega\cdot\partial_\varphi)(\omega\cdot\partial_\varphi)_{\textnormal{ext}}^{-1}\Pi_{N_n}=\Pi_{N_n},
$$
where $\Pi_{N_n}$ is the orthogonal projection defined by
$$
\Pi_{N_n}\sum_{l\in\Z^{d}}h_{l} {\bf{e}}_{l}=\sum_{l\in\Z^{d}\setminus\{0\}\atop |l|\leqslant N_n} h_{l} {\bf{e}}_{l} .
$$
\end{enumerate}
\end{lemma}

\begin{proof}
{\bf(i)}The proof of the first point can be done using Fa\`a di Bruno's formula in a similar way to   \cite[Lemma 2.5]{Baldi-berti}. 

\smallskip

{\bf(ii)} By construction, one has  for $\omega\in \mathtt {DC}_{N_n} (\kappa, \tau_1) $  and $|l|\leqslant N_n$,
$$
\chi\big((\omega\cdot l)\kappa^{-1}\langle l\rangle^{\tau_1} \big)=1,
$$
Thus, according to the explicit extension \eqref{Extend001},
\begin{align}\label{id:omg-ext}
(\omega\cdot\partial_\varphi)_{\textnormal{ext}}^{-1}\Pi_{N_n}h&=-\ii\sum_{l\in\Z^{d}\backslash\{0\}\atop |l|\leqslant N_n}\frac{\chi\big((\omega\cdot l)\kappa^{-1}\langle l\rangle^{\tau_1} \big)h_{l}(\lambda)}{\omega\cdot l}{\bf{e}}_{l}\nonumber\\
&=-\ii\sum_{l\in\Z^{d}\backslash\{0\}\atop |l|\leqslant N_n}\frac{h_{l}(\lambda)}{\omega\cdot l}{\bf{e}}_{l}.
\end{align}
Therefore, we obtain 
\begin{align*}
(\omega\cdot\partial_\varphi)(\omega\cdot\partial_\varphi)_{\textnormal{ext}}^{-1}\Pi_{N_n}h&=\sum_{l\in\Z^{d}\backslash\{0\}\atop |l|\leqslant N_n}{h_{l}(\lambda)}{\bf{e}}_{l}\\
&=\Pi_{N_n} h.
\end{align*}
This concludes the proof of the lemma.
\end{proof}

Now, we are   in a position to solve \eqref{Akjz}. The following lemma  is proved in \cite[Lemma 5.3]{BertiMontalto},
but for the sake of completeness we will give some key ingredients of the proof.
\begin{lemma}\label{lem:est-akj}
The coefficients $\mathcal{A}_{jk}$, defined in \eqref{def:akj}, admit the decomposition
\begin{equation}\label{akj-decomp}
\mathcal{A}_{kj}=\mathcal{A}_{kj}^{(n)}+\mathcal{A}_{kj}^{(n),\perp}\quad\hbox{with} \qquad \mathcal{A}_{kj}^{(n)}\triangleq \Pi_{N_n}\mathcal{A}_{kj}\qquad \hbox{and}\qquad \mathcal{A}_{kj}^{(n),\perp}\triangleq \Pi_{N_n}^\perp\mathcal{A}_{kj},
\end{equation}
where
\begin{enumerate}
\item The function  $\mathcal{A}_{kj}^{(n),\perp}$ satisfies  for any $s\in \mathbb{R}$,
  \begin{equation*}
  \| \mathcal{A}_{k j}^{(n),\perp} \|_{s}^{q,\kappa} \lesssim_{s, b} N_n^{-b}  \|  {\mathfrak I}_0 \|_{s+1+b}^{q,\kappa}, \quad \forall b\geqslant 0.
\end{equation*}
\item There exist functions $\mathcal{A}_{kj}^{(n),\textnormal{ext}}$ defined for any $\lambda=(\omega,\alpha)\in \mathcal{O}$, $q$-times differentiable with respect to $\lambda$ and satisfying, 
 for any $s\geqslant s_0$, the estimate
\begin{equation*}
\|  \mathcal{A}_{k j}^{(n),\textnormal{ext}} \|_s^{q,\kappa} \lesssim_s \kappa^{-1}
\big(\| Z \|_{s+\tau_1(q + 1)+1 }^{q,\kappa} + \| Z \|_{s_0+1}^{q,\kappa} \|  {\mathfrak I}_0 \|_{s+\tau_1(q + 1) +1}^{q,\kappa} \big)\,.
\end{equation*}
Moreover,  $\mathcal{A}_{k j}^{(n),\textnormal{ext}}$ coincides with $ \mathcal{A}_{k j}^{(n)}$ on the Cantor set $\mathtt {DC}_{N_n} (\kappa, \tau_1) $.
\end{enumerate}
\end{lemma}
\begin{proof}
The point  ${\bf (i)}$ follows immediately from \eqref{def:akj}, \eqref{akj-decomp} and  Lemma \ref{orthog-Lem1}.

\smallskip

${\bf (ii)}$ 
Applying the projector to the identity \eqref{Akjz} we obtain
\begin{align*}
\omega\cdot\partial_\varphi \mathcal{A}_{jk}^{(n)}(\varphi)&=  \Pi_{N_n}\big[{\mathcal W}\big( \partial_\varphi Z(\varphi) \underline{e}_k ,  \partial_\varphi i_0(\varphi)  \underline{e}_j \big)+ 
{\mathcal W} \big(\partial_\varphi i_0(\varphi) \underline{e}_k , \partial_\varphi Z(\varphi) \underline{e}_j \big) \big].
\end{align*}
Then, by \eqref{ansatz 0}, Lemma \ref{orthog-Lem1} and Lemma \ref{Law-prodX1} we get
\begin{equation*}
\big\|  \Pi_{N_n}\big[{\mathcal W}\big( \partial_\varphi Z(\varphi) \underline{e}_k ,  \partial_\varphi i_0(\varphi)  \underline{e}_j \big)+ 
{\mathcal W} \big(\partial_\varphi i_0(\varphi) \underline{e}_k , \partial_\varphi Z(\varphi) \underline{e}_j \big) \big] \big\|_s^{q,\kappa} \lesssim_s 
\big(\| Z \|_{s+1 }^{q,\kappa} + \| Z \|_{s_0+1}^{q,\kappa} \|  {\mathfrak I}_0 \|_{s+1}^{q,\kappa} \big)\,.
\end{equation*}
We define the the function $ \mathcal{A}_{k j}^{(n),\textnormal{ext}}$ as
$$
 \mathcal{A}_{k j}^{(n),\textnormal{ext}}(\varphi)\triangleq (\omega\cdot\partial_\varphi)_{\textnormal{ext}}^{-1}\Pi_{N_n}\big[{\mathcal W}\big( \partial_\varphi Z(\varphi) \underline{e}_k ,  \partial_\varphi i_0(\varphi)  \underline{e}_j \big)+ 
{\mathcal W} \big(\partial_\varphi i_0(\varphi) \underline{e}_k , \partial_\varphi Z(\varphi) \underline{e}_j \big) \big].
$$
Applying   Lemma \ref{L-Invert0}  concludes the proof of the Lemma.

\end{proof}

\subsection{Construction of an approximate inverse} Since the error term $\mathbb{E}$  is zero at an exact solution, up to an exponentially small remainder (see Proposition \ref{Prop-Conjugat}-{\rm (ii)} and Lemma \ref{lem:est-akj}) on the Cantor set $\mathtt {DC}_{N_n} (\kappa, \tau_1) $, then in order to find an approximate inverse of the linear operator in \eqref{bbL}
(and so of $ d_{(i,\tc)} {\mathcal F} (i_0,\tc_0) $) it is sufficient to almost invert the operator $\mathbb{D}$. 
The linear system $\mathbb{D}[\widehat u]=(g_1,g_2,g_3)$ may be solved in a triangular way, first inverting the action-component equation, which is decoupled from the other equations, 
 $$
 {\omega\cdot\partial_\varphi  \widehat y+\mathcal{B}(\varphi) \widehat \tc=g_2.}
$$
As we require only the finitely many non-resonance conditions \eqref{DC tau gamma N}, for any $\omega\in\mathbb{R}^d$ we also  decompose $\omega\cdot \partial_\varphi$ as
\begin{equation}\label{omeg-phi-decomp}
\begin{aligned}
&\omega\cdot \partial_\varphi=\mathcal{D}_{(n)} +\mathcal{D}_{(n)}^{\perp},\\
\mathcal{D}_{(n)}\triangleq \omega\cdot \partial_\varphi\, &\Pi_{N_n}+ \Pi_{N_n,\mathtt{g}}^\perp \qquad 
\mathcal{D}_{(n)}^{\perp}\triangleq  \omega\cdot \partial_\varphi\, \Pi_{N_n}^\perp- \Pi_{N_n, \mathtt{g}}^\perp 
\end{aligned}
\end{equation}
 where
$$
\Pi_{N_n, \mathtt{g}}^\perp\sum_{l\in\Z^{d}\setminus\{0\}}h_{l} {\bf{e}}_{l}\triangleq\sum_{l\in\Z^{d}\setminus\{0\}\atop |l|> N_n} g(l) h_{l} {\bf{e}}_{l}.
$$
and the function $\mathtt{g}:\mathbb{Z}^d\setminus \{0\}\to \{-1,1\}$ is defined, for all $l=(l_1,\cdots,l_d)\in \mathbb{Z}^d\setminus \{0\}$,   as the sign of the first non-zero component in the vector $l$. Thus, it satisfies 
$$\mathtt{g}(-l)=-\mathtt{g}(l)\quad \forall l\in \mathbb{Z}^d\setminus \{0\} .$$
Notice that for the high frequency we use the projector  $\Pi_{N_n, \mathtt{g}}^\perp$ instead of $\Pi_{N_n}^\perp$ in order to preserve the reversibility. We shall need the following Lemma.  

\begin{lemma}\label{lem:est-Dn}
The following holds true.
\begin{enumerate}
\item The operator $\mathcal{D}_{(n)}^{\perp}$ satisfies
$$
\forall b\geqslant 0,\, s\in \mathbb{R},\quad \|\mathcal{D}_{(n)}^{\perp}h\|_{s}^{q,\kappa }\lesssim N_n^{-b}\|h\|_{s+b+1}^{q,\kappa }. $$
\item There exists a family of linear operators $\big([\mathcal{D}_{(n)}]_{\textnormal{ext}}^{-1} \big)_n$ satisfying, for any $g\in H_0^s(\mathbb{T}^d)$,  
$$
\forall  s\geqslant s_0,\quad \sup_{n\in\mathbb{N}} \|[\mathcal{D}_{(n)}]_{\textnormal{ext}}^{-1} g\|_{s}^{q,\kappa }\lesssim \kappa^{-1}\|g\|_{s+q(\tau_1+1)}^{q,\kappa }.
$$
Moreover, for all $\omega\in  \mathtt {DC}_{N_n} (\kappa, \tau_1) $ one has the identity
$$
 \mathcal{D}_{(n)}[\mathcal{D}_{(n)}]_{\textnormal{ext}}^{-1}=\textnormal{Id}.
$$
\end{enumerate}
\end{lemma}

\begin{proof}
The point  ${\bf (i)}$ follows immediately from Lemma \ref{orthog-Lem1}.

\smallskip

${\bf (ii)}$  We define the operator $[\mathcal{D}_{(n)}]_{\textnormal{ext}}^{-1}$ as
\begin{equation}\label{dn-1}
[\mathcal{D}_{(n)}]_{\textnormal{ext}}^{-1}\triangleq (\omega\cdot\partial_\varphi)_{\textnormal{ext}}^{-1}\Pi_{N_n}+ \Pi_{N_n,\frac{1}{\mathtt{g}}}^\perp.
\end{equation}
From \eqref{id:omg-ext}, \eqref{omeg-phi-decomp} and \eqref{dn-1}  we get,
for all $\omega\in  \mathtt {DC}_{N_n} (\kappa, \tau_1) $,
\begin{align*}
\mathcal{D}_{(n)}[\mathcal{D}_{(n)}]_{\textnormal{ext}}^{-1}  &=\omega\cdot \partial_\varphi \,\Pi_{N_n}\big[(\omega\cdot\partial_\varphi)_{\textnormal{ext}}^{-1}\Pi_{N_n}+ \Pi_{N_n,\frac{1}{\mathtt{g}}}^\perp\big]+ \Pi_{N_n,\mathtt{g}}^\perp\big[(\omega\cdot\partial_\varphi)_{\textnormal{ext}}^{-1}\Pi_{N_n}+ \Pi_{N_n,\frac{1}{\mathtt{g}}}^\perp\big]
\\ 
&=\omega\cdot \partial_\varphi \,(\omega\cdot\partial_\varphi)_{\textnormal{ext}}^{-1}\Pi_{N_n}+ \Pi_{N_n}^\perp.
\end{align*}
 Applying  Lemma \ref{L-Invert0}-{\rm (ii)} we conclude that
$$
 \mathcal{D}_{(n)}[\mathcal{D}_{(n)}]_{\textnormal{ext}}^{-1}=\Pi_{N_n}+ \Pi_{N_n}^\perp=\textnormal{Id}.
$$
The estimate on $[\mathcal{D}_{(n)}]_{\textnormal{ext}}^{-1} $ follows from \eqref{dn-1}, Lemma \ref{orthog-Lem1} and  Lemma \ref{L-Invert0}.
\end{proof}

\smallskip

The second step in the resolution of the Linear equation  $\mathbb{D}[\widehat u]=g$, where $\mathbb{D}$  is given by Proposition \ref{Prop-Conjugat}-{\rm (i)}, is to solve the last normal-component equation
$$
 \omega\cdot\partial_\varphi \widehat w- 
\partial_\theta  K_{02}(\varphi)\widehat w =g_3+\partial_\theta\big[K_{11}(\varphi) \widehat y -L_2^{\top}(\varphi) \widehat\tc   \big].
 $$
For this aim we need to find an approximate  right  inverse of the linearized operator in the normal direction 
\begin{equation}\label{Lomega def}
\widehat{\mathcal{L}}_{\omega} \triangleq \Pi_{\mathbb{S}_0}^\bot \big(\omega\cdot \partial_\varphi   - 
\partial_\theta  K_{02}(\varphi) \big)\Pi_{{\mathbb{S}}_0}^\bot
\end{equation}
 when the set of parameters is restricted to a Cantor like set. Here the projector $\Pi_{\mathbb{S}_0}^\bot$ is the one defined in \eqref{projectors-tan-normal}.  For the sake of clarity we shall
give a brief statement about the invertibility in the normal direction; for a precise statement with a detailed description of Cantor like sets  see  Thoerem~\ref{inversion of the linearized operator in the normal directions}.

\begin{proposition}\label{thm:inversion of the linearized operator in the normal directions}
Let $(\kappa,q,d,\tau_{1},\tau_{2},s_{0})$ satisfy  \eqref{Conv-T2}, $(\mu_2,s_h)$  satisfy \eqref{cond-diman1} and assume the smallness condition \eqref{small-C3}. Then there exist $\sigma_5=\sigma_5(\tau_1,\tau_2,d,q)>0$, a family of linear  operator $(\mathtt{T}_{\omega,n})_n$ satisfying
\begin{equation}\label{estimate mathcalTomega0}
\forall \, s\in\,[ s_0, S] ,\quad\sup_{n\in\mathbb{N}}\|\mathtt{T}_{\omega,n}h\|_{s}^{q,\kappa }\lesssim\kappa^{-1}\left(\|h\|_{s+{\sigma_5}}^{q,\kappa }+\| \mathfrak{I}_{0}\|_{s+{\sigma_5}}^{q,\kappa }\|h\|_{s_{0}+{\sigma_5}}^{q,\kappa }\right)
\end{equation}
and   a  Cantor set $\mathtt{G}_n=\mathtt{G}_n(\kappa,\tau_{1},\tau_{2},i_{0})\subset \mathtt {DC}_{N_n} (\kappa, \tau_1) \times (\underline \alpha, \overline \alpha)$ in which 
we have the decomposition
$$
\widehat{\mathcal{L}}_{\omega}=\widehat{\mathtt{L}}_{\omega,n}+\widehat{\mathtt{R}}_n, 
$$
with
\begin{equation}\label{id-inv}
 \widehat{\mathtt{L}}_{\omega,n}\mathtt{T}_{\omega,n}=\textnormal{Id} 
\end{equation}
where the operators $\widehat{\mathtt{L}}_{\omega,n}$ and $\widehat{\mathtt{R}}_n$ are defined in the whole set $\mathcal{O}$ with the estimates
\begin{align*}
\forall\, s\in [s_0,S],\quad &\|\widehat{\mathtt{L}}_{\omega,n} h\|_{s}^{q,\kappa}\lesssim \|h\|_{s+1}^{q,\kappa}+{{\varepsilon\kappa^{-3}}\| \mathfrak{I}_{0}\|_{q,s+\sigma_5}^{q,\kappa }\|h\|_{s_{0}+1}^{q,\kappa }},\\
\forall\, b\in [0,S],\quad & \|\widehat{\mathtt{R}}_nh\|_{s_0}^{q,\kappa}
\lesssim N_n^{-b}\kappa^{-1}\Big( \|h\|_{s_0+b+\sigma_5}^{q,\kappa}+{\varepsilon\kappa^{-3}}\| \mathfrak{I}_{0}\|_{s_0+b+\sigma_5}^{q,\kappa}\|h\|_{s_0+{\sigma_5}}^{q,\kappa} \Big)\\ &\qquad\qquad\qquad\qquad\qquad\qquad\qquad+ \varepsilon\kappa^{-4}N_{0}^{{\mu}_{2}}{N_{n}^{-\mu_{2}}} \|h\|_{s_0+\sigma_5}^{q,\kappa}.
\end{align*}
\end{proposition}

\smallskip
According to Proposition~\ref{Prop-Conjugat}, the identities \eqref{akj-decomp}-\eqref{omeg-phi-decomp} and  Proposition~\ref{thm:inversion of the linearized operator in the normal directions} we may decompose the operator $[D G_0({\mathtt u}_0)]^{-1}  d_{(i,\tc)} {\mathcal F} (i_0,\tc_0) D\widetilde G_0({\mathtt u}_0)$  as follows
\begin{equation}\label{decomposition}
[D G_0({\mathtt u}_0)]^{-1}  d_{(i,\tc)} {\mathcal F} (i_0,\tc_0) D\widetilde G_0({\mathtt u}_0)={\mathbb D}_n+{\mathbb E}_n+ {\mathscr P}_n+{\mathscr Q}_n
\end{equation}
with
\begin{align}
&{\mathbb D}_n [\widehat \phi, \widehat y, \widehat w, \widehat \tc ]  \triangleq
 \left(
\begin{array}{c}
\mathcal{D}_{(n)}  \widehat\phi-K_{20}\widehat y-K_{11}^\top \widehat w -L_1^{\top}\widehat \tc
\\
\mathcal{D}_{(n)} \widehat y+ \mathcal{B} \widehat \tc \\
\widehat{\mathtt{L}}_{\omega,n} \widehat w -\partial_\theta\big[K_{11} \widehat y -L_2^{\top} \widehat\tc   \big] 
\end{array}
\right), \label{def:Dn}
\\ &{\mathbb E}_{n} [\widehat \phi, \widehat y, \widehat w]  \triangleq
 [D G_0({\mathtt u}_0)]^{-1}    [\partial_\varphi Z] \widehat\phi+\left(
\begin{array}{c}
0 
\\
 \mathcal{A}^{(n)}\big[K_{20} \widehat y+K_{11}^\top \widehat w\big]-R_{10}\widehat y-R_{01}\widehat w
 \\
0   
\end{array}
\right) \label{def:Rn0}
\end{align}
and
 \begin{align}
&{\mathscr P}_n [\widehat \phi, \widehat y, \widehat w ]  \triangleq
 \left(
\begin{array}{c}
\mathcal{D}_{(n)}^{\perp}  \widehat\phi
\\
\mathcal{D}_{(n)}^{\perp} \widehat y + \mathcal{A}^{(n),\perp}\big[K_{20} \widehat y+K_{11}^\top \widehat w\big]
 \\
0
\end{array}
\right),\quad {\mathscr Q}_n  [\widehat \phi, \widehat y, \widehat w ]  \triangleq
 \left(
\begin{array}{c}
0
\\
0
 \\
\widehat{\mathtt{R}}_n[ \widehat  w]
\end{array}
\right),
 \label{def:Pn}
\end{align}
where $\mathcal{A}^{(n)}$ and $\mathcal{A}^{(n),\perp}$ are the matrices with coefficients 
$\mathcal{A}_{kj}^{(n)}$ and $\mathcal{A}_{kj}^{(n),\perp}$ respectively, see \eqref{akj-decomp}.
Because the coefficients $\mathcal{A}_{jk}^{(n)}$ do not vanish at  exact solutions on whole set of parameters $\mathcal{O}$ and the $\mathcal{A}_{kj}^{(n),\textnormal{ext}}$ do (see Lemma \ref{lem:est-akj}), we shall replace the operator  $\mathbb{E}_n$ by the extension 
\begin{align}
{\mathbb E}_{n}^{{\rm ext}} [\widehat \phi, \widehat y, \widehat w ] & \triangleq
 [D G_0({\mathtt u}_0)]^{-1}    [\partial_\varphi Z] \widehat\phi+\left(
\begin{array}{c}
0 
\\
 \mathcal{A}^{(n),\textnormal{ext}}\big[K_{20}(\varphi) \widehat y+K_{11}^\top \widehat w\big]-R_{10}\widehat y-R_{01}\widehat w
 \\
0   
\end{array}
\right), \label{def:Rn}
\end{align}
with $\mathcal{A}^{(n),\textnormal{ext}}$ is the matrix with coefficients 
$\mathcal{A}_{kj}^{(n),\textnormal{ext}}$. Thus, we define the linear operator $\mathbb{L}_{\textnormal{ext}}$ as \begin{equation}\label{def:lext}
\mathbb{L}_{\textnormal{ext}}= \mathbb{D}_n+{\mathbb E}_{n}^{{\rm ext}}+{\mathscr P}_n+{\mathscr Q}_n.
\end{equation}
The operator $\mathbb{L}_{\textnormal{ext}}$ is  defined on the whole set $\mathcal{O}$ and, by construction,  coincides with the linear operator in \eqref{decomposition} on the Cantor set $\mathtt{G}_n$,
\begin{equation}\label{lext-f}
\forall \lambda \in \mathtt{G}_n,\quad\mathbb{L}_{\textnormal{ext}}=[D G_0({\mathtt u}_0)]^{-1}  d_{({\mathtt u},\tc)} {\mathcal F} (i_0,\tc_0) D \tilde G_0({\mathtt u}_0).
\end{equation}

Now, we shall estimate the error terms ${\mathbb E}_{n}^{{\rm ext}}$, ${\mathscr P}_n$, ${\mathscr Q}_n$ and find  an exact inverse of the principal term  $\mathbb{D}_n$. This will be the subject of the following proposition.

\begin{proposition}\label{prop:decomp-lin}
Let $(\kappa,q,d,\tau_{1},s_{0})$ satisfy  \eqref{Conv-T2}, $(\mu_2,s_h)$  satisfy \eqref{cond-diman1} and assume the smallness condition \eqref{small-C3}. Then the following assertions hold true.
\begin{enumerate}

\item The operator ${\mathbb E}_{n}^{{\rm ext}}$ satisfies the estimate 
\begin{align*} 
\| {\mathbb E}_{n}^{{\rm ext}} [\widehat{\mathtt u}] \|_{s_0}^{q, \kappa}
&\lesssim_s \| Z \|_{s_0+ 1 }^{q, \kappa}  \| \widehat{\mathtt u} \|_{s_0 + 1 }^{q, \kappa}.
\end{align*}

\item For all $ b\geqslant 0$ the operator $ {\mathscr P}_n $ satisfies the estimate
\begin{align*} 
\|  {\mathscr P}_n [\widehat{\mathtt u}] \|_{s_0}^{q, \kappa} &\lesssim N_n^{-b}  \big( \| \widehat{\mathtt u} \|_{s_0 + 2 + b }^{q,\kappa}+
\| \mathfrak{I}_{0} \|_{s_0
+  2  +b    }^{q,\kappa} \big \| \widehat{\mathtt u} \|_{s_0 +2}^{q,\kappa} \big).
\end{align*}
\item  For all $ b\in [0,S]$,
the operator $ {\mathscr Q}_n $ satisfies the estimate
\begin{align*} \|{\mathscr Q}_n [\widehat{\mathtt u}]\|_{s_0}^{q,\kappa}
\lesssim N_n^{-b}\kappa^{-1}\Big( \| \widehat{\mathtt u}\|_{s_0+b+\sigma_5}^{q,\kappa}+{\varepsilon\kappa^{-3}}\| \mathfrak{I}_{0}\|_{s_0+b+\sigma_5}^{q,\kappa}\| \widehat{\mathtt u}\|_{s_0+{\sigma_5}}^{q,\kappa} \Big)+ \varepsilon\kappa^{-4}N_{0}^{{\mu}_{2}}{N_{n}^{-\mu_{2}}} \| \widehat{\mathtt u}\|_{s_0+\sigma_5}^{q,\kappa}.
\end{align*}

\item There exist $\sigma_6=\sigma_6(\tau_1,\tau_2,d,q)>0$ and  a family of operators $\big([{\mathbb D}_n]_{\textnormal{ext}}^{-1}\big)_n$  such that
 for all $ g \triangleq (g_1, g_2, g_3) $ 
satisfying  the symmetry  
\begin{equation*}
g_1(\varphi) = g_1(- \varphi)\,,\quad g_2(\varphi) = - g_2(- \varphi)\,,\quad g_3(\varphi) = - ({\mathcal S} g_3)(\varphi)\,,
\end{equation*}
 the function 
$ [{\mathbb D}_n]_{\textnormal{ext}}^{-1} g $ 
  satisfies the estimate, for any $s_0 \leqslant s \leqslant S$,
\begin{equation*} 
\| [{\mathbb D}_n]_{\textnormal{ext}}^{-1}g \|_s^{q, \kappa}
\lesssim_S \kappa^{-1} \big( \| g \|_{s + \sigma_6}^{q, \kappa}
+  \| {\mathfrak I}_0  \|_{s + \sigma_6}^{q, \kappa}
 \| g \|_{s_0 + \sigma_6}^{q, \kappa}  \big)
\end{equation*}
and for all $\lambda \in \mathtt{G}_n$ one has
$$
 {\mathbb D}_n [{\mathbb D}_n]_{\textnormal{ext}}^{-1} =\textnormal{Id}.
$$
\end{enumerate}
\end{proposition}
\begin{proof}
{\bf (i)} From \eqref{def:Rn}, Lemma \ref{lem:est-G}, Lemma \ref{Law-prodX1}-(ii),  Lemma \ref{lem:est-akj}-{\rm (ii)}, we get the estimate on  ${\mathbb E}_{n}^{{\rm ext}}$.

\smallskip

{\bf (ii)} The estimate on   ${\mathscr P}_n$ can be easily obtained from \eqref{def:Pn}, Lemma \ref{lem:est-Dn}-{\rm (i)},   Lemma \ref{Law-prodX1}-(ii),  Lemma \ref{lem:est-akj}-{\rm (i)}, Lemma \ref{lem:est-G}-{\rm (ii)}.

\smallskip

{\bf (iii)} Follows immediately from \eqref{def:Pn} and Theorem \ref{thm:inversion of the linearized operator in the normal directions}.

\smallskip

{\bf (iv)} We shall look for an exact inverse of ${\mathbb D}_n$ 
by solving the system 
\begin{equation}\label{operatore inverso approssimato proiettato}
{\mathbb D}_n [\widehat \phi, \widehat y, \widehat w, \widehat \tc]  
= \begin{pmatrix}
g_1  \\
g_2  \\
g_3 
\end{pmatrix}
\end{equation}
where $(g_1, g_2, g_3)$ satisfy the reversibility property 
\begin{equation}\label{parita g1 g2 g3}
g_1(\varphi) = g_1(- \varphi)\,,\quad g_2(\varphi) = - g_2(- \varphi)\,,\quad g_3(\varphi) = - ({\mathcal S} g_3)(- \varphi)\,.
\end{equation}
 We first consider the second equation in \eqref{operatore inverso approssimato proiettato} which takes the form, according to \eqref{def:Dn},
$$ \mathcal{D}_{(n)} \, \widehat y  =
  g_2  -\mathcal{B}(\varphi) \widehat \tc. $$ 
Since $g_2$ is odd,  
 the $\varphi$-average of the right hand side of this equation is zero then by Lemma \ref{lem:est-Dn}-${\rm (iv)}$  its solution on the Cantor set $\mathtt{ DC}_{N_n}$ is  
\begin{equation}\label{soleta}
\widehat y \triangleq [\mathcal{D}_{(n)}]_{\textnormal{ext}}^{-1} \big(
  g_2  -\mathcal{B}(\varphi) \widehat \tc\big)\,.  
\end{equation}
Then we consider the third equation
$\widehat{\mathtt{L}}_{\omega,n} \widehat w = g_3 + \partial_\theta K_{11}(\varphi) \widehat y + \partial_\theta L_2^\top(\varphi)  \widehat \tc $,
which, by Theorem \ref{thm:inversion of the linearized operator in the normal directions}, has a solution 
\begin{equation}\label{normalw}
\widehat w \triangleq \mathtt{T}_{\omega,n} \big( g_3 + \partial_\theta K_{11}(\varphi) \widehat y + \partial_\theta L_2^\top(\varphi)  \widehat \tc\big) \, .  
\end{equation}
Finally, we solve the first equation in \eqref{operatore inverso approssimato proiettato}, 
which, substituting \eqref{soleta}, \eqref{normalw}, becomes
\begin{equation}\label{equazione psi hat}
\mathcal{D}_{(n)} \widehat \phi  = 
g_1 +  M_1(\varphi)\widehat \tc + M_2(\varphi) g_2 + M_3(\varphi) g_3\,,
\end{equation}
where
\begin{align}\label{M1}
 M_1(\varphi) &\triangleq L_1^\top (\varphi)- M_2(\varphi) \mathcal{B}(\varphi)  + M_3(\varphi)  \partial_\theta L_2^\top(\varphi)  \,, \\
 \label{cal M2}
M_2(\varphi) &\triangleq  K_{20}(\varphi)  [\mathcal{D}_{(n)}]_{\textnormal{ext}}^{-1}  + K_{11}^\top(\varphi)\mathtt{T}_{\omega,n} \partial_\theta K_{11}(\varphi)[\mathcal{D}_{(n)}]_{\textnormal{ext}}^{-1}  \, , \\
M_3(\varphi) &\triangleq K_{11}^\top (\varphi) \mathtt{T}_{\omega,n} \, .    \label{cal M3}
\end{align}
In order to solve equation \eqref{equazione psi hat} we have 
to choose $ \widehat \tc $ such that the right hand side  has zero average.
In fact, by Lemma \ref{lem:est-G}, \eqref{ansatz 0}, \eqref{estimate mathcalTomega0}, Lemma \ref{lem:est-Dn}-{\rm (ii)} we get
\begin{equation*} 
\| M_2(\varphi)g \|_{s}^{q, \kappa}+\| M_3(\varphi)g \|_{s}^{q, \kappa}
\lesssim \varepsilon\Big(  \| g \|_{s + \sigma}^{q, \kappa}
+  \| {\mathfrak I}_0  \|_{s + \sigma}^{q, \kappa}
 \| g \|_{s_0 + \sigma}^{q, \kappa}\Big) .
\end{equation*}
Then from Lemma \ref{lem:est-G}, \eqref{small-C3}, the $\phi$-averaged matrix is
$ \langle M_1 \rangle = {\rm Id} + O( \eps \kappa^{-1 }) $.  Therefore, for $ \eps \kappa^{- 1} $ small enough,  
$ \langle M_1 \rangle$ is invertible and $\langle M_1 \rangle^{-1} = {\rm Id} 
+ O(\eps \kappa^{- 1})=O(1)$. 
 Thus we define 
\begin{equation}\label{sol alpha}
\widehat \tc  \triangleq - \langle M_1 \rangle^{-1} 
( \langle g_1 \rangle + \langle M_2 g_2 \rangle + \langle M_3 g_3 \rangle ) \, .
\end{equation}
Therefore we deduce that
\begin{equation}\label{estimate-c}
\|\widehat \tc\|^{q, \kappa}\lesssim  \| g \|_{s_0 + \sigma }^{q, \kappa}
\end{equation}
With this choice of $ \widehat \tc $,
equation \eqref{equazione psi hat} has the solution
\begin{equation}\label{sol psi}
\widehat \phi \triangleq
[\mathcal{D}_{(n)}]_{\textnormal{ext}}^{-1} \big( g_1 + M_1(\varphi)[\widehat \tc] + M_2(\varphi) g_2 + M_3(\varphi) g_3 \big) \, . 
\end{equation}
Coming back to \eqref{soleta} and using Lemma \ref{lem:est-G}-{\rm (iv)},  Lemma \ref{lem:est-Dn}-{\rm (ii)}, \eqref{estimate-c}    we obtain
\begin{equation}\label{estimate-y}
\| \widehat y \|_{s}^{q, \kappa}
\lesssim \kappa^{-1} \big( \| g \|_{s + \sigma}^{q, \kappa}
+  \| {\mathfrak I}_0  \|_{s + \sigma}^{q, \kappa}
 \| g \|_{s_0 + \sigma}^{q, \kappa}  \big).
\end{equation}
The bound in {\rm (iv)}  for $\widehat w$ can be easily obtained from  \eqref{normalw}, \eqref{estimate-c}, \eqref{estimate mathcalTomega0}, \eqref{estimate-y}, Lemma \ref{lem:est-G}-{\rm (iv)} and the smallness condition. Consequently, we find that  $\widehat \phi$ satisfies the estimate in {\rm (iv)} using \eqref{sol psi}, Lemma \ref{lem:est-G},  Lemma \ref{lem:est-Dn}-{\rm (ii)}  and \eqref{estimate mathcalTomega0}. Finally, the identity  
$$
\forall \lambda\in  \mathtt{G}_n,\quad   {\mathbb D}_n [{\mathbb D}_n]_{\textnormal{ext}}^{-1} =\textnormal{Id}
$$
follows from Lemma \ref{lem:est-Dn} and \eqref{id-inv}.
\end{proof}

\smallskip

Finally, we shall prove that the operator 
\begin{equation}\label{definizione T} 
{\rm T}_0 \triangleq {\rm T}_0(i_0) \triangleq D { \widetilde G}_0({\mathtt u}_0)\, [{\mathbb D}_n]_{\textnormal{ext}}^{-1}\,[D G_0({\mathtt u}_0)]^{-1} 
\end{equation}
is an approximate right  inverse for $d_{i,\tc} {\mathcal F}(i_0 )$. 

\begin{theorem}  \label{thm:stima inverso approssimato}
{\bf (Approximate inverse)}

Let $(\gamma,q,d,\tau_{1},s_{0})$ satisfy  \eqref{Conv-T2}, $(\mu_2,s_h)$  satisfy \eqref{cond-diman1} and assume  \eqref{small-C3}. 
Then there exists $ { \overline\sigma}= { \overline\sigma}(\tau_1,\tau_2,d,q)>0$ such that
for smooth $ g = (g_1, g_2, g_3) $, satisfying \eqref{parita g1 g2 g3},  
the operator $ {\rm T}_0  $ defined in \eqref{definizione T} satisfies 
\begin{equation}\label{estimate on T0}
\forall s\in [s_0,S],\quad \| {\rm T}_0 g\|_{s}^{q,\kappa}\lesssim\kappa^{-1}\left(\|g\|_{s+{\overline\sigma}}^{q,\kappa}+\|\mathfrak{I}_{0}\|_{s+{\overline\sigma}}^{q,\kappa}\|g\|_{s_{0}+{\overline\sigma}}^{q,\kappa}\right).
\end{equation}
 Moreover  ${\rm T}_0$ is an almost-approximate  
 right 
 inverse of $d_{i, \tc} 
 \mathcal{ F}(i_0, \tc_0)$ on the Cantor set $ \mathtt{G}_n $. More precisely,   for all $ \lambda\in \mathtt{G}_n $ one has
\begin{equation}\label{splitting per approximate inverse}
d_{i , \tc} \mathcal{ F} (i_0,\tc_0)  {\rm T}_0
- {\rm Id} = \mathcal{E}^{(n)}_1+\mathcal{E}^{(n)}_2+\mathcal{E}^{(n)}_3,
\end{equation}
where the operators $\mathcal{E}^{(n)}_1$, $\mathcal{E}^{(n)}_2$ and $\mathcal{E}^{(n)}_3$ are defined in the whole set $\mathcal{O}$ with the estimates
\begin{align*}
\|{\mathcal{E}_1^{(n)}} h \|_{s_0}^{q,\kappa} & \lesssim \overline \gamma^{-1 } \| \mathcal{ F}(i_0, \tc_0) \|_{s_0 +\overline\sigma}^{q,\kappa} \| h \|_{s_0 + \overline\sigma}^{q,\kappa},
 \\
\| \mathcal{E}_2^{(n)} h \|_{s_0}^{q,\kappa}& \lesssim_{ b} 
\kappa^{- 1} N_n^{- b } \big( \| h \|_{s_0 + \overline\sigma + b }^{q,\kappa}+
\| \mathfrak{I}_{0} \|_{s_0
+  \overline\sigma  +b    }^{q,\kappa} \big \| h \|_{s_0 +\overline\sigma}^{q,\kappa} \big)\,,\,\,
\quad\forall   b\geqslant 0,
\\
\| \mathcal{E}_3^{(n)} h \|_{s_0}^{q,\kappa}& \lesssim N_n^{-b}\kappa^{-2}\Big( \|h\|_{s_0+b+\overline\sigma}^{q,\kappa}+{\varepsilon\kappa^{-3}}\| \mathfrak{I}_{0}\|_{s_0+b+\overline\sigma}^{q,\kappa}\|h\|_{s_0+{\overline\sigma}}^{q,\kappa} \Big)+ \varepsilon\kappa^{-5}N_{0}^{{\mu}_{2}}{N_{n}^{-\mu_{2}}} \|h\|_{s_0+\overline\sigma}^{q,\kappa},\; \forall b\in [0,S].  
\end{align*}
\end{theorem}

\begin{proof}
The bound \eqref{estimate on T0} follows from \eqref{definizione T}, Proposition \ref{prop:decomp-lin}-{\rm (iv)} and Lemma \ref{lem:est-G}-{\rm (i)}. Then, according to  \eqref{def:lext} and 
\eqref{lext-f}, on the Cantor set $ \mathtt{G}_n $ we have the decomposition
\begin{equation*}
\begin{aligned}
 d_{i,\tc}\mathcal{F}(i_{0},\tc_{0})&=DG_{0}({\mathtt u}_0) \,\mathbb{L}_{\textnormal{ext}}\,  [D\widetilde{G}_{0}({\mathtt u}_0)]^{-1}
 \\
 &=DG_{0}({\mathtt u}_0) \, {\mathbb{D}}_n\, [D\widetilde{G}_{0}({\mathtt u}_0)]^{-1}+ DG_{0}({\mathtt u}_0) \, {\mathbb E}_{n}^{{\rm ext}} \, [D\widetilde{G}_{0}({\mathtt u}_0)]^{-1}\\ &+DG_{0}({\mathtt u}_0)\,   {\mathscr P}_n\, [D\widetilde{G}_{0}({\mathtt u}_0)]^{-1}+DG_{0}({\mathtt u}_0) \, {\mathscr Q}_n\, [D\widetilde{G}_{0}({\mathtt u}_0)]^{-1}.
\end{aligned}
\end{equation*}
Applying ${\rm T}_0$, defined in \eqref{definizione T},  to the right of the last identity   we get for all $\lambda=(\omega,\alpha)\in  \mathtt{G}_n$
$$d_{i,\tc}\mathcal{F}(i_{0},\tc_{0}) {\rm T}_{0}-\textnormal{Id}=\mathcal{E}_1^{(n)}+\mathcal{E}_2^{(n)}+\mathcal{E}_3^{(n)}
$$
with
\begin{align*}
&\mathcal{E}_1^{(n)}\triangleq DG_{0}({\mathtt u}_0) \,  {\mathbb E}_{n}^{{\rm ext}} \, [D\widetilde{G}_{0}({\mathtt u}_0)]^{-1}{\rm T}_{0},
\\
&\mathcal{E}_2^{(n)}\triangleq DG_{0}({\mathtt u}_0)\,  {\mathscr P}_n\, [D\widetilde{G}_{0}({\mathtt u}_0)]^{-1}{\rm T}_{0},
\\
&\mathcal{E}_3^{(n)}\triangleq DG_{0}({\mathtt u}_0)\,  {\mathscr Q}_n\, [D\widetilde{G}_{0}({\mathtt u}_0)]^{-1}{\rm T}_{0}.
\end{align*}
The estimates on $\mathcal{E}^{(n)}_1$, $\mathcal{E}^{(n)}_2$ and $\mathcal{E}^{(n)}_3$  follow from \eqref{estimate on T0}, Proposition \ref{prop:decomp-lin} and Lemma \ref{lem:est-G}-{\rm (i)}. 

\end{proof}

\section{Toroidal pseudo-differential operators}
In this  section we shall collect some tools from function spaces and toroidal  pseudo-differential operators. In particular, we shall introduce suitable functions spaces, discuss some operator topologies  related to pseudo-differential operators and establish  some useful commutator estimates connected with modified fractional Laplacian. 
\subsection{Operators and reversibility }
This section is devoted to some discussions on parameter-dependent operator norms. Consider a smooth family of bounded operators acting on Sobolev spaces $H^s(\T^{d+1})$, that is a smooth map 
$$\mathcal{A}: \lambda\in \mathcal{O}\to \mathcal{A}(\lambda)\in  \mathcal{L}(H^s(\T^{d+1}); H^s(\T^{d+1})) 
$$ of linear continuous operators $\mathcal{A}(\lambda): H^s(\T^{d+1})\to H^s(\T^{d+1}) $. Then  $\mathcal{A}(\lambda)$ can be represented through  the infinite  matrix $\left(\mathcal{A}_{l_{0},j_{0}}^{l,j}(\lambda)\right)_{\underset{(j,j_{0})\in\mathbb{Z}^{2}}{(l,l_{0})\in(\mathbb{Z}^{d})^{2}}}$ with $$\mathcal{A}(\lambda)\mathbf{e}_{l_{0},j_{0}}=\sum_{(l,j)\in\mathbb{Z}^{d}\times\mathbb{Z}}\mathcal{A}_{l_{0},j_{0}}^{l,j}(\lambda)\mathbf{e}_{l,j}$$
and
\begin{eqnarray*}\mathcal{A}_{l_{0},j_{0}}^{l,j}(\lambda)&\triangleq&\langle \mathcal{A}(\lambda)\mathbf{e}_{l_{0},j_{0}},\mathbf{e}_{l,j}\rangle_{L^{2}(\mathbb{T}^{d+1})}\\
&=&\tfrac{1}{(2\pi)^{d+1}}\int_{\mathbb{T}^{d+1}}\big(\mathcal{A}(\lambda)\mathbf{e}_{l_{0},j_{0}}\big)(\varphi,\theta)\,\mathbf{e}_{-l,-j}(\varphi,\theta)d\varphi d\theta.
\end{eqnarray*}
We say that the  operator $\mathcal{A}(\lambda)$ is T\"oplitz in time (actually in the variable $\varphi$) if its Fourier coefficients satisfy,
$$
\mathcal{A}_{l_{0},j_{0}}^{l,j}(\lambda)=\mathcal{A}_{0,j_0}^{l-l_0,j}(\lambda).
$$
Now we say that an element $\mathcal{A}$ belongs to the  function space $W^{q,\infty}_{\kappa}(\mathcal{O}; \mathcal{L}(H^s; H^s)) $  if and only if
\begin{align}\label{strong-Top1}
\|\mathcal{A}\|_{s}^{q,\kappa}=\max_{\underset{|\alpha|\leqslant q}{\alpha\in\mathbb{N}^{d}}}\gamma^{|\alpha|}\sup_{\lambda\in{\mathcal{O}}}\|\partial_{\lambda}^{\alpha}\mathcal{A}(\lambda)\|_{\mathcal{L}(H^{s-|\alpha|})}<\infty.
\end{align}
The next result deals with the algebra structure of the preceding space.
\begin{lemma}
Let $\mathcal{A}_1, \mathcal{A}_2\in W^{q,\infty}_{\kappa}(\mathcal{O}; \mathcal{L}(H^s; H^s)) $, then  $\mathcal{A}_1\mathcal{A}_2\in W^{q,\infty}_{\kappa}(\mathcal{O}; \mathcal{L}(H^s; H^s))$ with
$$
\|\mathcal{A}_1 \mathcal{A}_2\|_{s}^{q,\kappa}\leqslant C\|\mathcal{A}_1\|_{s}^{q,\kappa}\|\mathcal{A}_2\|_{s}^{q,\kappa}.
$$

\end{lemma}
We remark that this class of operators is not useful in our case  because it is too weak to get tame estimates. Therefore it will be replaced by a stronger one that will be detailed in Section \ref{Se-Toroidal pseudo-differential operators}.
Notice also  that  along this paper the operators and the test functions may depend on the same  parameter $\lambda$ and thus the action of a family $\mathcal{A}=(\mathcal{A}(\lambda))$ on an element $h\in W^{q,\infty}_{\kappa}(\mathcal{O},H^{s}(\T^{d}))$  is by convention defined through
$$
(\mathcal{A}h)(\lambda,\varphi,\theta)\triangleq \mathcal{A}(\lambda) h(\lambda,\varphi,\theta),
$$
meaning that we observe both objects at the same point $\lambda$. This is justified by the fact that  with respect to the  equations that we consider here, $\lambda$ is an external parameter in the model.

In what follows we shall   collect some definitions and properties about reversible operators.
We recall the following definition, see for instance \cite{Baldi-berti,BFM21,BFM,BertiMontalto}.
We recall  the involution $\mathcal{S}$ introduced in Proposition \ref{prop-hamilt},
\begin{equation}\label{definition of the involution mathcal S}
\forall\, (\varphi,\theta)\in\T^{d+1},\quad (\mathcal{S}h)(\varphi,\theta)= h(-\varphi,-\theta).
\end{equation}
\begin{definition}\label{Def-Rev}
Consider a $\lambda$-dependent family of operators  $\mathcal{A}(\lambda)$.
We say that $\mathcal{A}$ is 
\begin{enumerate}
\item real if for all $h\in L^{2}(\mathbb{T}^{d+1},\mathbb{C}),$ we have 
$$\overline{h}=h\quad\Longrightarrow\quad\overline{\mathcal{A} h}=\mathcal{A}h.$$
\item reversible if
$$\mathcal{A}(\lambda)\circ\mathcal{S}=-\mathcal{S}\circ \mathcal{A}(\lambda).$$
\item reversibility preserving if
$$\mathcal{A}(\lambda)\circ\mathcal{S}=\mathcal{S}\circ \mathcal{A}(\lambda).
$$
\end{enumerate}
\end{definition}
We now give  a  useful  characterization   based on Fourier coefficients and the proof is straightforward.
\begin{proposition}\label{characterization of real operator by its Fourier coefficients}
Let $\mathcal{A}$ be an operator, then we have the following results. 
\begin{enumerate}
\item $\mathcal{A}$ is  real if and only if
$$\forall\, l,l_{0},j_{0},\,j \in\, \mathbb{Z},\quad \mathcal{A}_{-l_{0},-j_{0}}^{-l,-j}=\overline{\mathcal{A}_{l_{0},j_{0}}^{l,j}}.$$
{\item  $\mathcal{A}$ is  reversible if and only if
$$\forall\,l,l_{0},j_{0},\,j\, \in\, \mathbb{Z},\quad \mathcal{A}_{-l_{0},-j_{0}}^{-l,-j}=-\mathcal{A}_{l_{0},j_{0}}^{l,j}.$$
\item $\mathcal{A}$ is reversibility-preserving if and only if
$$\forall\,l,l_{0},j_{0},\,j\, \in\, \mathbb{Z},\quad \mathcal{A}_{-l_{0},-j_{0}}^{-l,-j}=\mathcal{A}_{l_{0},j_{0}}^{l,j}.$$}
\end{enumerate}
\end{proposition}

\subsection{Symbol class topology}\label{Se-Toroidal pseudo-differential operators}
In this section we shall collect some classical results on pseudo-differential operators on the periodic setting. Let   $\mathcal{A}:\mathscr{C}^\infty(\T)\to \mathscr{C}^\infty(\T)$ be a linear operator, we define its symbol $\sigma_{\mathcal{A}}$ by
$$
\forall\,(\theta,\xi)\in\mathbb{T }\times\Z,\quad \sigma_{\mathcal{A}}(\theta,\xi)\triangleq {\bf e}_{-\xi}(\theta) (\mathcal{A}{\bf e}_{\xi})(\theta),\quad {\bf{e}}_{\xi}(\theta)={ e}^{\ii \theta\,\xi}.
$$
Then for $\displaystyle h(\theta)=\sum_{\xi\in\Z}h_\xi e^{\ii  \theta\,\xi}$ we get
$$
\mathcal{A}h(\theta)=\sum_{\xi\in\Z}  \sigma_{\mathcal{A}}(\theta,\xi)h_\xi e^{\ii  \theta\,\xi}.
$$
One also gets  the kernel representation
$$
\mathcal{A}h(\theta)=\int_{\T} K(\theta,\eta)h(\eta)d\eta,\quad\text{with}\quad
K(\theta,\eta)={\frac{1}{2\pi}}\sum_{\xi\in\Z}  \sigma_{\mathcal{A}}(\theta,\xi)e^{\ii (\theta-\eta)\xi}.
$$
By the Fourier inversion formula we may recover the symbol from the kernel through the formula 
$$
\sigma_{\mathcal{A}}(\theta,\xi)=\int_\T K(\theta,\theta+\eta)e^{\ii \eta\xi}d\eta.
$$
Next, we shall introduce the difference operators in the periodic setting. Consider a    discrete function $h:\Z\to\C$ and define  the forward difference operator  $\Delta_\xi$    by
\begin{equation}\label{Difference-op78}
\Delta_\xi h(\xi)=h(\xi+1)-h(\xi)\quad\hbox{and}\quad \Delta_\xi^\ell=\underbrace{{\Delta_\xi\circ\ldots\circ\Delta_\xi}}_{\ell-\text{ times}}. 
\end{equation}
The backward difference $\overline\Delta_\xi$ is defined by 
$$
\overline\Delta_\xi h(\xi)\triangleq h(\xi)-h(\xi-1),\quad \overline\Delta_\xi^\ell=\underbrace{{\overline\Delta_\xi\circ\ldots\circ\overline\Delta_\xi}}_{\ell-\text{ times}}. 
$$
We will recall some useful identities of the difference operator, see for instance \cite[Chapter 3]{Ruzhansky}. First, the discrete Leibniz formula reads as follows,
\begin{align}\label{Leibn-disc}
\forall\, \gamma\in\N,\quad\Delta_\xi^\gamma(fg)(\xi)=\sum_{0\leqslant \beta\leqslant\gamma}\left(_\beta^\gamma\right)\Delta_\xi^\beta f(\xi)\,\Delta_\xi^{\gamma-\beta}g(\xi+\beta).
\end{align} 
The next one is the summation by parts
\begin{align}\label{Sum-parts}
\sum_{\xi\in\Z}f(\xi)\Delta_\xi^\gamma(g)(\xi)=(-1)^{|\gamma|}\sum_{\xi\in\Z}\left(\overline\Delta_\xi^\gamma f(\xi)\right)\, g (\xi).
\end{align}
The next property is easy to check by induction
\begin{align}\label{It-action}
\Delta_\xi^\gamma e^{\ii \xi \eta}=\big(e^{\ii  \eta}-1\big)^\gamma e^{\ii \xi \eta}.
\end{align}
Throughout this paper  we shall make use of   multi-parameter pseudo-differential  operators $\varphi\in \T^d\mapsto \mathcal{A}(\varphi)$. This yields to  one-parameter symbol $\sigma_{\mathcal{A}}(\varphi)$ defined in the following way: for a smooth periodic function  $\displaystyle h(\varphi,\theta)=\sum_{\xi\in\Z}h_\xi(\varphi) e^{\ii \theta\,\xi}$
\begin{align}\label{A-phi}
\mathcal{A}(\varphi)h(\varphi,\theta)=\sum_{\xi\in\Z}  \sigma_{\mathcal{A}}(\varphi,\theta,\xi)h_\xi(\varphi) e^{\ii \theta\,\xi}.
\end{align}
The symbol can be recovered from the operator  as follows
$$
\forall\,(\theta,\xi)\in\mathbb{T }\times\Z,\quad \sigma_{\mathcal{A}}(\varphi,\theta,\xi)\triangleq {\bf e}_{-\xi}(\theta) (\mathcal{A}(\varphi){\bf e}_{\xi})(\theta),\quad {\bf e}_{\xi}(\theta)=e^{\ii \theta\,\xi}
$$
and the kernel representation takes  the form
$$
\mathcal{A}(\varphi)h(\varphi,\theta)=\int_{\T} K(\varphi,\theta,\eta)h(\varphi,\eta)d\eta,\quad\text{with}\quad
K(\varphi,\theta,\eta)={\frac{1}{2\pi}}\sum_{\xi\in\Z} \sigma_{\mathcal{A}}(\varphi,\theta,\xi) e^{\ii(\theta-\eta)\xi} .
$$
By the Fourier inversion formula we get 
\begin{equation}\label{symb-kern}
\sigma_{\mathcal{A}}(\varphi,\theta,\xi)=\int_\T K(\varphi,\theta,\theta+\eta)e^{\ii \eta\xi}d\eta.
\end{equation}
Let $s,m\in\mathbb{R},\gamma\in\mathbb{N}$ and define  the norm over the symbol class
\begin{align}\label{Def-Norm-M1}
\interleave \mathcal{A}\interleave_{m,s,\gamma}&\triangleq\sup_{\xi\in\Z\atop
0\leqslant\ell\leqslant\gamma}\langle \xi\rangle^{-m+\ell}\big\|\Delta_\xi^\ell\sigma_{\mathcal{A}}(\cdot,\centerdot,\xi)\|_{H^{s}(\T^{d+1})}.
\end{align}
The case $m=\gamma=0$ plays a central  role along the paper, needed especially in the remainder KAM reduction, and for the sake of simple notation we simply write
$$
\interleave \mathcal{A}\interleave_{s}\triangleq \interleave \mathcal{A}\interleave_{0,s,0}.
$$
We shall give an equivalent form to this norm using Fourier coefficients which turns out to be very useful later at some points. Define
$$
\mathcal{A}_{j'}^{j}(l)=\langle \mathcal{A}{\bf{e}}_{0,j'},{\bf{e}}_{l,j}\rangle_{L^2(\T^{d+1})} \quad\hbox{with}\quad {\bf{e}}_{l,j}(\varphi,\theta)=e^{\ii (l\cdot\varphi+j\theta)}.
$$
Then it is easy to check that for any $j\in\mathbb{Z}$ (the operator is T\"oplitz in $\varphi$)
\begin{align*}
\sigma_{\mathcal{A}}(\varphi,\theta,j)&={\bf e}_{0,-j}(\varphi,\theta) \big(\mathcal{A}{\bf e}_{0,j}\big)(\varphi,\theta)\\
&=\sum_{(l,j')\in\mathbb{Z}^{d+1}}\mathcal{A}_{j}^{j'}(l){\bf e}_{l,j'-j}(\varphi,\theta).
\end{align*}
Consequently
\begin{align*}
\|\sigma_{\mathcal{A}}(\cdot,\centerdot,j)\|_{H^s}^2
&=\sum_{(l,j')\in\mathbb{Z}^{d+1}}\big|\mathcal{A}_{j}^{j+j'}(l)\big|^2\langle l,j'\rangle^{2s},
\end{align*}
which implies that
\begin{align}\label{equivalent-norm}
\nonumber \interleave \mathcal{A}\interleave_{s}^2&=\sup_{j\in\Z}\|\sigma_{\mathcal{A}}(\cdot,\centerdot,j)\|_{H^s}^2\\
&=\sup_{j\in\Z}\sum_{(l,j')\in\mathbb{Z}^{d+1}}\langle l,j'\rangle^{2s}\big|\mathcal{A}_{j}^{j+j'}(l)\big|^2.
\end{align}
We shall also make use of weighted pseudo-differential operators $(\lambda,\varphi)\in \mathcal{O}\times\T^d\mapsto \mathcal{A}(\lambda, \varphi)$ where $ \mathcal{O}$ is a bounded open set of $\RR^{d+1}$. First the multi-parameter symbol of $\mathcal{A}$  is defined as follows,
\begin{align}\label{A-phi-lambda}
 h(\lambda,\varphi,\theta)=\sum_{\xi\in\Z}h_\xi(\lambda,\varphi) e^{\ii  \theta\,\xi}\Longrightarrow\,\, \mathcal{A}(\lambda,\varphi)h(\lambda,\varphi,\theta)=\sum_{\xi\in\Z}  \sigma_{\mathcal{A}}(\lambda,\varphi,\theta,\xi)h_\xi(\lambda,\varphi) e^{\ii  \theta\,\xi}.
\end{align}
In this case we have the kernel representation,
\begin{align}\label{kernel-phi-lambda}
\mathcal{A}(\lambda,\varphi)h(\lambda,\varphi,\theta)=\int_{\T} K(\lambda,\varphi,\theta,\eta)h(\lambda,\varphi,\eta)d\eta
\end{align}
with
$$
K(\lambda,\varphi,\theta,\eta)=\frac{1}{2\pi}\sum_{\xi\in\Z} e^{\ii (\theta-\eta)\xi} \sigma_{\mathcal{A}}(\lambda,\varphi,\theta,\xi).
$$
Take  $m, s\in\RR, \gamma\in\mathbb{N}$ and $\kappa\in(0,1)$ and define the weighted  norm 
\begin{align}\label{Def-pseud-w}
\interleave \mathcal{A}\interleave_{m,s,\gamma}^{q,\kappa}\triangleq \max_{\underset{|\beta|\leqslant q}{\beta\in\mathbb{N}^{d+1 }}}\kappa^{|\beta|}\sup_{\lambda\in{\mathcal{O}}}\interleave \partial_\lambda^{\beta}\mathcal{A}\interleave_{m,s-|\beta|,\gamma}.
\end{align}
We shall adopt for the particular case $m=\gamma=0$  the following notation
$$
\interleave \mathcal{A}\interleave_{s}^{q,\kappa}\triangleq \interleave \mathcal{A}\interleave_{0,s,0}^{q,\kappa}.
$$
Using \eqref{equivalent-norm} we get the following characterization,
\begin{align}\label{Top-NormX}
\left(\interleave \mathcal{A}\interleave_{s}^{q,\kappa}\right)^2&=\max_{|\beta|\leqslant q}\sup_{\lambda\in\mathcal{O}}\sup_{j\in\Z}\sum_{(l,j')\in\mathbb{Z}^{d+1}}\kappa^{2|\beta|}\langle l,j'\rangle^{2(s-|\beta|)}\big|\partial_\lambda^\beta\mathcal{A}_{j}^{j+j'}(\lambda,l)\big|^2.
\end{align}
The next lemma is very useful.
\begin{lemma}\label{lemma-Sym-R}
{Let $q,\gamma\in\mathbb{N},\,m,s\in\RR$, then the following assertions hold true.
\begin{enumerate}
\item Let $T_M$ be a multiplication operator by a real-valued function $M$, that is, $T_M h=Mh$. Then \begin{enumerate}
\item If $M(\lambda,-\varphi,-\theta)=M(\lambda,\varphi,\theta)$, then $T_M$ is  real and reversibility preserving T\"oplitz in time and space operator.
\item If $M(\lambda,-\varphi,-\theta)=-M(\lambda,\varphi,\theta)$, then $T_M$ is  real and reversible T\"oplitz in time and space operator. Moreover, 
$$\interleave T_M\interleave_{0,s,\gamma}^{q,\kappa}=\| M\|_{s}^{q,\kappa}.$$
\end{enumerate}

\item Let $\mathcal{A}$ be an integral operator with a real-valued kernel $K$ as in \eqref{kernel-phi-lambda}.
\begin{enumerate}
\item If $K(\lambda,-\varphi,-\theta,-\eta)=K(\lambda,\varphi,\theta,\eta)$, then $\mathcal{A}$ is   real and reversibility preserving T\"oplitz in time operator.
\item If $K(\lambda,-\varphi,-\theta,-\eta)=-K(\lambda,\varphi,\theta,\eta)$, then $\mathcal{A}$ is   real and reversible T\"oplitz in time operator. In addition, 
$$\interleave \mathcal{A}\interleave_{m,s,\gamma}^{q,\kappa}\lesssim  \sum_{0\leqslant \ell\leqslant\gamma}\bigintssss_{\T}\big|\sin(\eta/2)\big|^\ell\|\partial_\eta^\ell(-\Delta_\eta)^{-\frac{m}{2}} K(\cdot,\centerdot,\centerdot+\eta)\|_{s}^{q,\kappa} d\eta.
$$

\end{enumerate}
\item Let $\mathcal{A}$ as in \eqref{kernel-phi-lambda} and $m<-\frac12$ then 
\begin{align*}
\int_{\T}\big(\|K(\cdot,\centerdot,\centerdot+\eta)\|_{s}^{q,\kappa}\big)^2d\eta
&\lesssim\big(\interleave \mathcal{A}\interleave_{m,s,0}^{q,\kappa}\big)^2.
\end{align*}

\end{enumerate}}

\end{lemma}
\begin{proof}
{\bf{(i)-(ii)}} The proofs are starightfoward.\\
{\bf{(iii)}}
Recall that the link between the kernel and the symbol is given by
$$
K(\varphi,\theta,\theta+\eta)={\frac{1}{2\pi}}\sum_{\xi\in\Z}\sigma_{\mathcal{A}}(\varphi,\theta,\xi)e^{-\ii  \,\eta\xi}.
$$
Thus from Bessel identity we get for any $(\varphi,\theta)\in\T^{d+1}$
\begin{align*}
\int_{\T}\big|K(\varphi,\theta,\theta+\eta)\big|^2d\eta=&\sum_{\xi\in\Z}\big|\sigma_{\mathcal{A}}(\varphi,\theta,\xi)\big|^2.
\end{align*}
Integrating in $(\varphi,\theta)$ we get
\begin{align*}
\int_{\T}\|K(\cdot,\centerdot,\centerdot+\eta)\|_{L^2_{\varphi,\theta}}^2d\eta=&\sum_{\xi\in\Z}\|\sigma_{\mathcal{A}}(\cdot,\centerdot,\xi)\|_{L^2_{\varphi,\theta}}^2.
\end{align*}
Then differentiating in $\varphi, \theta$ and using the same argument as before we get for $m<-\frac12$,
\begin{align*}
\int_{\T}\|K(\cdot,\centerdot,\centerdot+\eta)\|_{H^s_{\varphi,\theta}}^2d\eta &\lesssim\sum_{\xi\in\Z}\|\sigma_{\mathcal{A}}(\cdot,\centerdot,\xi)\|_{H^s_{\varphi,\theta}}^2\\
&\lesssim\interleave {\mathcal{A}}\interleave_{m,s,0}^2\sum_{\xi\in\Z} \langle \xi\rangle^{2m}\\
& \lesssim\interleave {\mathcal{A}}\interleave_{m,s,0}^2.
\end{align*}
Similarly we obtain for any $q\in\NN$,
\begin{align*}
\int_{\T}\big(\|K(\cdot,\centerdot,\centerdot+\eta)\|_{s}^{q,\kappa}\big)^2d\eta &\lesssim\sum_{\xi\in\Z}\big(\| \sigma_{\mathcal{A}}(\cdot,\centerdot,\xi)\|_{s}^{q,\kappa}\big)^2\\
&\lesssim\big(\interleave {\mathcal{A}}\interleave_{m,s,0}^{q,\kappa}\big)^2.
\end{align*}
This ends the proof.
\end{proof}
\subsection{Continuity and law products}
This section is devoted to some classical results. First, we define the frequency cut-off projectors $(P_N)_ {N\in\mathbb{N}^{*}}$ as follows, 
\begin{equation}\label{definition of projections for operators}
\left(P_{N} \mathcal{A}(\lambda)\right)\mathbf{e}_{l_{0},j_{0}}=\sum_{\underset{|l-l_{0}|,|j-j_{0}|\leqslant N}{(l,j)\in\mathbb{Z}^{d }\times\mathbb{Z}}} \mathcal{A}_{l_{0},j_{0}}^{l,j}(\lambda)\mathbf{e}_{l,j}\quad\mbox{and}\quad P_{N}^{\perp} \mathcal{A}= \mathcal{A}-P_{N} \mathcal{A},
\end{equation}
where $\mathcal{A}$ is a pseudo-differential operator as in Section \ref{Se-Toroidal pseudo-differential operators}.
In the next lemma we shall gather classical results some of the  proofs can be found for instance in \cite{BertiMontalto,Ruzhansky}.
	
\begin{lemma}\label{Lem-Rgv1}
Let $\mathcal{A}$  be a pseudo-differential operator as in \eqref{A-phi}, then the following assertions hold.
\begin{enumerate}
\item Let $ {s_0>{d+1}}$, then
$$
\|\mathcal{A}h\|_{L^2(\T^{d+1})}\lesssim   \interleave \mathcal{A}\interleave_{s_0}\|h\|_{L^2(\T^{d+1})}.
$$
\item If $s\geqslant 0, s_0>\frac{d+1}{2}$ then 
$$
\|\mathcal{A}h\|_{H^s}\lesssim  \interleave \mathcal{A}\interleave_{s_0}\|h\|_{H^s}+ \interleave \mathcal{A}\interleave_{s}\|h\|_{H^{s_0}}.
$$
\item Let  $\mathcal{A}$ as in \eqref{A-phi-lambda} and $s\geqslant q, s_0>\frac{d+1}{2}+q$. Then 
$$
\|\mathcal{A}h\|_{s}^{q,\overline \gamma}\lesssim \interleave \mathcal{A}\interleave_{s_0}^{q,\kappa} \|h\|_{s}^{q,\kappa}+\interleave\mathcal{A}\interleave_{s}^{q,\kappa} \|h\|_{s_0}^{q,\kappa}.
$$
 \item Let $N\in\mathbb{N}^{*}, s\in\RR, t\geqslant 0,$ then
$$\interleave P_{N}\mathcal{A}\interleave _{s+t}^{q,\kappa}\leqslant N^{t}\interleave \mathcal{A}\interleave_{s}^{q,\kappa}\quad\mbox{and}\quad\interleave P_{N}^\perp\mathcal{A}\interleave _{s}^{q,\kappa}\leqslant N^{-t}\interleave \mathcal{A}\interleave_{s+t}^{q,\kappa}.$$
\item Interpolation inequality : Let $q<s_{1}\leqslant s_{3}\leqslant s_{2}, t\in[0,1]$ with $s_{3}=ts_{1}+(1-t)s_{2}.$  Then 
$$\interleave \mathcal{A}\interleave _{s_3}^{q,\kappa}\leqslant C(q,d )\left(\interleave \mathcal{A}\interleave _{s_1}^{q,\kappa}\right)^{t}\left(\interleave \mathcal{A}\interleave _{s_2}^{q,\kappa}\right)^{1-t}.$$
\item {Adjoint estimate}: the $L^2$ adjoint $ \mathcal{A}^*$ satisfies
$$
\interleave \mathcal{A}^*\interleave _{m,s,0}^{q,\kappa}\lesssim \interleave \mathcal{A}\interleave _{m,s+s_0+|m|,0}^{q,\kappa}.
$$
\end{enumerate}
\end{lemma}
Now we recall the following result, see for instance \cite{BertiMontalto}.
\begin{lemma}\label{comm-pseudo1}
 Let  $\mathcal{A}, {\mathcal{B}}$ as in \eqref{A-phi-lambda} and $(\kappa,d,q,s,s_0)$  as in \eqref{initial parameter condition} and consider  $m_1,m_2\in\RR, \gamma\in\N.$ Then we have
\begin{align*}
\interleave\mathcal{A}{\mathcal{B}}\interleave_{m_1+m_2,s,\gamma}^{q,\kappa}&\lesssim \interleave \mathcal{A}\interleave_{m_1,s,\gamma}^{q,\kappa}\interleave {\mathcal{B}}\interleave_{m_2,s_0+\gamma+|m_1|,\gamma}^{q,\kappa}+\interleave \mathcal{A}\interleave_{m_1,s_0,\gamma}^{q,\kappa} \interleave {\mathcal{B}}\interleave_{m_2,s+\gamma+|m_1|,\gamma}^{q,\kappa}.\end{align*}
In particular we get
$$
\interleave\mathcal{A}{\mathcal{B}}\interleave_{s}^{q,\kappa}\lesssim \interleave \mathcal{A}\interleave_{s_0}^{q,\kappa} \interleave{\mathcal{B}}\interleave_{s}^{q,\kappa}+\mathcal{A}\interleave_{s}^{q,\kappa} \interleave{\mathcal{B}}\interleave_{s_0}^{q,\kappa}.
$$
\end{lemma}
Now we intend to establish  refined estimates needed later to derive tame estimates in Section \ref{Flows and Egorov theorem type}.
\begin{lemma}\label{comm-pseudo}
 Let  $\mathcal{A}, {\mathcal{B}}$ as in \eqref{A-phi-lambda} and $(\kappa,d,q,s,s_0)$  as in \eqref{initial parameter condition}. Then the following assertions hold true.
\begin{enumerate}
\item Let $m_1,m_2\in\RR, $ then for any $\gamma\in\N$ and $\epsilon>0$
\begin{align*}
&\interleave\mathcal{A}{\mathcal{B}}\interleave_{m_1+m_2,s,\gamma}^{q,\kappa}\lesssim \\ &\sum_{0\leqslant\beta\leqslant\gamma}\Big(\interleave \mathcal{A}\interleave_{m_1,s,\beta}^{q,\kappa}\interleave{\mathcal{B}}\interleave_{m_2,s_0+\frac12+m_1^++\epsilon,\gamma-\beta}^{q,\kappa}+\interleave \mathcal{A}\interleave_{m_1,s_0,\beta}^{q,\kappa}\interleave{\mathcal{B}}\interleave_{m_2,s+\frac12+m_1^++\epsilon,\gamma-\beta}^{q,\kappa}\Big)\\
&+\sum_{ 0\leqslant\beta\leqslant\gamma}\left(\interleave \mathcal{A}\interleave_{m_1,s,0}^{q,\kappa}\interleave {\mathcal{B}}\interleave_{m_2,s_0+\frac12+\beta-m_1^-,\gamma-\beta}^{q,\kappa} +\interleave \mathcal{A}\interleave_{m_1,s_0,0}^{q,\kappa}\interleave {\mathcal{B}}\interleave_{m_2,s+\frac12+\beta-m_1^-,\gamma-\beta}^{q,\kappa} \right).\end{align*}
\item If  $m_1,m_2\leqslant1$ and $\,\gamma\in\N,$ then for any $\epsilon>0$ 
\begin{align*}
&\interleave[\mathcal{A},{\mathcal{B}}]\interleave_{m_1+m_2-1,s,\gamma}^{q,\kappa}\lesssim\\ &
\sum_{0\leqslant\beta\leqslant\gamma}\interleave \mathcal{A}\interleave_{m_1,s,1+\beta}^{q,\kappa}\interleave{\mathcal{B}}\interleave_{m_2,s_0+\frac32+\epsilon,\gamma-\beta}^{q,\kappa}+\interleave \mathcal{A}\interleave_{m_1,s_0,1+\beta}^{q,\kappa}\interleave{\mathcal{B}}\interleave_{m_2,s+\frac32+\epsilon,\gamma-\beta}^{q,\kappa}\\
&+\sum_{ 0\leqslant\beta\leqslant\gamma}\interleave \mathcal{A}\interleave_{m_1,s,0}^{q,\kappa}\interleave {\mathcal{B}}\interleave_{m_2,s_0+\frac32+\beta-m_1^-,\gamma-\beta}^{q,\kappa} +\interleave \mathcal{A}\interleave_{m_1,s_0,0}^{q,\kappa}\interleave {\mathcal{B}}\interleave_{m_2,s+\frac32+\beta-m_1^-,\gamma-\beta}^{q,\kappa} \\
&+
\sum_{0\leqslant\beta\leqslant\gamma}\interleave {\mathcal{B}}\interleave_{m_2,s,1+\beta}^{q,\kappa}\interleave\mathcal{A}\interleave_{m_1,s_0+\frac32+\epsilon,\gamma-\beta}^{q,\kappa}+\interleave {\mathcal{B}}\interleave_{m_2,s_0,1+\beta}^{q,\kappa}\interleave\mathcal{A}\interleave_{m_1,s+\frac32+\epsilon,\gamma-\beta}^{q,\kappa}\\
&+\sum_{ 0\leqslant\beta\leqslant\gamma}\interleave {\mathcal{B}}\interleave_{m_2,s,0}^{q,\kappa}\interleave \mathcal{A}\interleave_{m_1,s_0+\frac32+\beta-m_2^-,\gamma-\beta}^{q,\kappa} +\interleave {\mathcal{B}}\interleave_{m_2,s_0,0}^{q,\kappa}\interleave \mathcal{A}\interleave_{m_1,s+\frac32+\beta-m_2^-,\gamma-\beta}^{q,\kappa}.
\end{align*}
\item Let   $m_1+m_2\leqslant 1$ and $\gamma\in\N,$ then 
\begin{align*}
&\interleave[\mathcal{A},{\mathcal{B}}]\interleave_{0,s,\gamma}^{q,\kappa}\lesssim \\ &\sum_{0\leqslant\beta\leqslant\gamma}\interleave \mathcal{A}\interleave_{m_1,s,1+\gamma-\beta}^{q,\kappa}\interleave{\mathcal{B}}\interleave_{m_2,s_0+{\mu},\beta}^{q,\kappa}+\interleave \mathcal{A}\interleave_{m_1,s_0,1+\gamma-\beta}^{q,\kappa}\interleave{\mathcal{B}}\interleave_{m_2,s+{\mu},\beta}^{q,\kappa}
\\
&+ \sum_{0\leqslant\beta\leqslant\gamma}\interleave \mathcal{A}\interleave_{m_1,s,0}^{q,\kappa} \interleave {\mathcal{B}}\interleave_{m_2,s_0+{\mu}_\beta,\gamma-\beta}^{q,\kappa}+\interleave \mathcal{A}\interleave_{m_1,s_0,0}^{q,\kappa} \interleave {\mathcal{B}}\interleave_{m_2,s+{\mu}_\beta,\gamma-\beta}^{q,\kappa}\\
&+ \sum_{0\leqslant\beta\leqslant\gamma}\interleave\mathcal{A}\interleave_{m_1,s_0+{\mu},\beta}^{q,\kappa}\interleave {\mathcal{B}}\interleave_{m_2,s,1+\gamma-\beta}^{q,\kappa}+\interleave\mathcal{A}\interleave_{m_1,s+{\mu},\beta}^{q,\kappa}\interleave {\mathcal{B}}\interleave_{m_2,s_0,1+\gamma-\beta}^{q,\kappa}
\\
&+\sum_{0\leqslant\beta\leqslant\gamma}\interleave \mathcal{A}\interleave_{m_1,s_0+\overline{\mu}_\beta,\gamma-\beta}^{q,\kappa}\interleave {\mathcal{B}}\interleave_{m_2,s,0}^{q,\kappa} +\interleave \mathcal{A}\interleave_{m_1,s+\overline{\mu}_\beta,\gamma-\beta}^{q,\kappa}\interleave {\mathcal{B}}\interleave_{m_2,s_0,0}^{q,\kappa},
\end{align*}
provided that $\mu, \mu_\beta$ and $\overline{\mu}_\beta$ satisfy 
$$
{\mu}>\frac12+\max\big(0,m_1,m_2,m_1+m_2\big),\quad  {\mu}_\beta\geqslant\frac12+m_1^++m_2+\beta,\quad \overline{\mu}_\beta\geqslant \frac12+m_1+m_2^++\beta,
 $$ 
 with the notation $m^+=\max(m,0)$ and $m^-=\min(m,0).$
\end{enumerate}

\end{lemma}
\begin{proof}
{$($\bf{i}$)$} We shall implement the proof for the case $q=0$ and the general case $q\geqslant1$ can be done similarly using Leibniz rule. Indeed, one writes
$$
\kappa^{|\beta|}\partial_\lambda^\beta\mathcal{A}{\mathcal{B}}=\sum_{0\leqslant\beta^\prime\leqslant\beta}\left(^\beta_{\beta^\prime}\right)\kappa^{|\beta^\prime|}(\partial_\lambda^{\beta^\prime}\mathcal{A}) \kappa^{|\beta-\beta^\prime|}\partial_\lambda^{\beta-\beta^\prime}{\mathcal{B}}
$$
and then we may apply the case $q=0$ in order to get the desired result. This remark concerns also the points {$($\bf{ii}$)$} and {$($\bf{iii}$)$}. Now let us prove the estimate for $q=0$. First notice that the symbol of the operator $\mathcal{A}{\mathcal{B}}$ is given by
\begin{align*}
\sigma_{\mathcal{A}{\mathcal{B}}}(\varphi,\theta,\xi)&=\sum_{\eta\in \Z}e^{\ii  \eta\,\theta}\sigma_{\mathcal{A}}(\varphi,\theta,\xi+\eta)\widehat{\sigma}_{{\mathcal{B}}}(\varphi,\eta,\xi)
\end{align*}
where $\widehat{\sigma}_{{\mathcal{B}}}(\varphi,\eta,\xi)$ is the Fourier coefficients of the partial periodic function $\theta\in\T\mapsto \sigma_A (\varphi,\theta,\xi).$
Then one gets the decomposition,
\begin{align*}
\sigma_{\mathcal{A}{\mathcal{B}}}(\varphi,\theta,\xi)
&=\sigma_{\mathcal{A}{\mathcal{B}}}^{\textnormal{H}}(\varphi,\theta,\xi)+\sigma_{\mathcal{A}{\mathcal{B}}}^{\textnormal{L}}(\varphi,\theta,\xi)
\end{align*}
with 
\begin{align*}
\sigma_{\mathcal{A}{\mathcal{B}}}^{\textnormal{H}}(\varphi,\theta,\xi)&=\sum_{ |\eta|\leqslant \frac{|\xi|}{2}\atop
|\eta|\geqslant  2|\xi|}\frac{\sigma_{\mathcal{A}}(\varphi,\theta,\xi+\eta)}{|\eta|^{\mu}} e^{\ii  \eta\,\theta}\widehat{\Lambda_\theta^{\mu}\sigma}_{{\mathcal{B}}}(\varphi,\eta,\xi)
\end{align*}
and
\begin{align*}
\sigma_{\mathcal{A}{\mathcal{B}}}^{\textnormal{L}}(\varphi,\theta,\xi)&=\sum_{\frac{|\xi|}{2}< |\eta|< 2|\xi|}\frac{\sigma_{\mathcal{A}}(\varphi,\theta,\xi+\eta)}{|\eta|^{\mu_\beta}} e^{\ii  \eta\,\theta}\widehat{\Lambda_\theta^{\mu_\beta}\sigma}_{{\mathcal{B}}}(\varphi,\eta,\xi).
\end{align*}
Here,  $\Lambda_\theta=\sqrt{-\Delta_\theta}$ and $\mu, \mu_\beta$ are free real numbers. Let us start with studying the high frequency part. 
Then using Leibniz formula \eqref{Leibn-disc} yields
\begin{align*}
\Delta_\xi^\gamma\sigma_{\mathcal{A}{\mathcal{B}}}^{\textnormal{H}}(\varphi,\theta,\xi)&=\sum_{\beta\leqslant\gamma}\sum_{ |\eta|\leqslant \frac{|\xi|}{2}\atop
|\eta|\geqslant  2|\xi|}\left(_\beta^\gamma\right)\frac{\Delta_\xi^\beta\sigma_{\mathcal{A}}(\varphi,\theta,\xi+\eta)}{|\eta|^{\mu}}e^{\ii  \eta\,\theta}\Delta_\xi^{\gamma-\beta}\widehat{\Lambda_\theta^{\mu}\sigma}_{{\mathcal{B}}}(\varphi,\eta,\xi+\beta).
\end{align*}
Applying  Sobolev  law products combined with \eqref{Rdd11}, Cauchy-Schwarz inequality and the choice 
$$\mu=\frac12+m_1^++\epsilon, \epsilon>0
$$ we get
\begin{align}\label{Yo0T1}
\nonumber \left\|\Delta_\xi^\gamma\sigma_{\mathcal{A}{\mathcal{B}}}^{\textnormal{H}}(\cdot,\centerdot,\xi)\right\|_{H^{s}(\T^{d+1})}&\lesssim {C}_0\sum_{ |\eta|\leqslant \frac{|\xi|}{2}\atop
|\eta|\geqslant  2|\xi|}c_\eta(\xi) |\eta|^{-\mu}\langle \xi+\eta\rangle^{m_1-\beta} \langle \xi\rangle^{m_2-\gamma+\beta}\\
&\lesssim {C}_0 \langle \xi\rangle^{m_1+m_2-\gamma},
\end{align}
with $\sum_{\eta} c_\eta^2(\xi)=1$ and
$$
{C}_0=\sum_{0\leqslant\beta\leqslant\gamma}\interleave \mathcal{A}\interleave_{m_1,s,\beta}\interleave{\mathcal{B}}\interleave_{m_2,s_0+\frac12+m_1^++\epsilon,\gamma-\beta}+\interleave \mathcal{A}\interleave_{m_1,s_0,\beta}\interleave{\mathcal{B}}\interleave_{m_2,s+\frac12+m_1^++\epsilon,\gamma-\beta}.
$$
Now let us move to the low frequency part. First we have 
\begin{align*}
\Delta_\xi^\gamma\sigma_{\mathcal{A}{\mathcal{B}}}^{\textnormal{L}}(\varphi,\theta,\xi)&=\sum_{\beta\leqslant\gamma}\sum_{\frac{|\xi|}{2}< |\eta|< 2|\xi|}\left(_\beta^\gamma\right)\frac{\Delta_\xi^\beta\sigma_{\mathcal{A}}(\varphi,\theta,\xi+\eta)}{|\eta|^{\mu_\beta}}e^{\ii  \eta\,\theta}\Delta_\xi^{\gamma-\beta}\widehat{\Lambda_\theta^{\mu_\beta}\sigma}_{{\mathcal{B}}}(\varphi,\eta,\xi+\beta).
\end{align*}
Then by definition we get  for $\frac{|\xi|}{2}\leqslant |\eta|\leqslant 2|\xi|$
\begin{align*}
\left\|\Delta_\xi^{\beta}\left(\sigma_{\mathcal{A}}(\cdot,\centerdot,\xi+\eta)\right)\right\|_{H^s}&\lesssim \langle \xi\rangle^{m_1^+}\interleave \mathcal{A}\interleave_{m_1,s,0}.
\end{align*}
Combined with the inequality \eqref{Rdd11}, proved later,  and taking  
\begin{align*}\mu_\beta=\tfrac12+\beta-m_1^-
\end{align*}
 we find by virtue of Cauchy-Schwarz inequality
\begin{align}\label{YoU1}
\nonumber \|\Delta_\xi^\gamma\mathcal{R}_1^{\textnormal{L}}(\cdot,\centerdot,\xi)\|_{H^s}&\lesssim\sum_{\frac{|\xi|}{2}< |\eta|< 2|\xi|\atop 0\leqslant\beta\leqslant\gamma} c_\eta(\xi) |\eta|^{-\mu_\beta}\langle \xi\rangle^{m_1^++m_2-\gamma+\beta}\,{C}_{1,\beta}\\
\nonumber &\lesssim C_1 \langle \xi\rangle^{m_1+m_2-\gamma}\sum_{\frac{|\xi|}{2}< |\eta|< 2|\xi|} c_\eta(\xi) |\eta|^{-\frac12}\\
&\lesssim C_1 \langle \xi\rangle^{m_1+m_2-\gamma}
\end{align}
with
\begin{align*}
C_1&=\sum_{0\leqslant \beta\leqslant\gamma}{C}_{1,\beta}\\
C_{1,\beta}\triangleq \interleave \mathcal{A}\interleave_{m_1,s,0}\interleave {\mathcal{B}}\interleave_{m_2,s_0+\frac12+\beta-m_1^-,\gamma-\beta} &+\interleave \mathcal{A}\interleave_{m_1,s_0,0}\interleave {\mathcal{B}}\interleave_{m_2,s+\frac12+\beta-m_1^-,\gamma-\beta}.
\end{align*}
This achieves the proof of the desired result.

\smallskip

{$($\bf{ii}$)$} We shall start with using the following decomposition 
\begin{align}\label{abou-s1}
\sigma_{\mathcal{A}{\mathcal{B}}}(\varphi,\theta,\xi)&=\sum_{\eta\in \Z}e^{\ii  \eta\,\theta}\sigma_{\mathcal{A}}(\varphi,\theta,\xi+\eta)\widehat{\sigma}_{{\mathcal{B}}}(\varphi,\eta,\xi)\\
\nonumber&=\sum_{\eta\in \Z}e^{\ii  \eta\,\theta}\Big(\sigma_{\mathcal{A}}(\varphi,\theta,\xi+\eta)-\sigma_{\mathcal{A}}(\varphi,\theta,\xi)\Big)\widehat{\sigma}_{{\mathcal{B}}}(\varphi,\eta,\xi)\\
\nonumber&\quad+\sigma_{\mathcal{A}}(\varphi,\theta,\xi)\sigma_{{\mathcal{B}}}(\varphi,\theta,\xi)\\
\nonumber &\triangleq\mathcal{R}_1(\varphi,\theta,\xi)+\sigma_{\mathcal{A}}(\varphi,\theta,\xi)\sigma_{{\mathcal{B}}}(\varphi,\theta,\xi).
\end{align}
 Let $\mu$ be an arbitrary  parameter then we have the splitting,
\begin{align*}
\mathcal{R}_1(\varphi,\theta,\xi)&=\mathcal{R}_1^{\textnormal{H}}(\varphi,\theta,\xi)+\mathcal{R}_1^{\textnormal{L}}(\varphi,\theta,\xi)
\end{align*}
with 
\begin{align*}
\mathcal{R}_1^{\textnormal{H}}(\varphi,\theta,\xi)&=\sum_{ |\eta|\leqslant \frac{|\xi|}{2}\atop
|\eta|\geqslant  2|\xi|}\frac{\sigma_{\mathcal{A}}(\varphi,\theta,\xi+\eta)-\sigma_{\mathcal{A}}(\varphi,\theta,\xi)}{|\eta|^{\mu}} e^{\ii  \eta\,\theta}\widehat{\Lambda_\theta^{\mu}\sigma}_{{\mathcal{B}}}(\varphi,\eta,\xi)
\end{align*}
and
\begin{align*}
\mathcal{R}_1^{\textnormal{L}}(\varphi,\theta,\xi)&=\sum_{\frac{|\xi|}{2}< |\eta|< 2|\xi|}\left({\sigma_{\mathcal{A}}(\varphi,\theta,\xi+\eta)-\sigma_{\mathcal{A}}(\varphi,\theta,\xi)}\right) e^{\ii  \eta\,\theta}\widehat{\sigma}_{{\mathcal{B}}}(\varphi,\eta,\xi).
\end{align*}
Let us start with estimating the first term $\mathcal{R}_1^{\textnormal{H}}$. Concerning  the region $|\eta|\leqslant \frac{|\xi|}{2}$ we may use discrete Taylor formulae, see \cite[Chapter 3]{Ruzhansky}, leading to the following inequality 
\begin{equation}\label{WX1}
\left\|\Delta_\xi^{\beta}\sigma_{\mathcal{A}}(\cdot,\centerdot,\xi+\eta)-\Delta_\xi^{\beta}\sigma_{\mathcal{A}}(\cdot,\centerdot,\xi)\right\|_{H^s(\T^{d+1})}\lesssim |\eta|\sup_{|\nu|\leqslant|\eta|}\left\|\Delta_\xi^{1+\beta}\sigma_{\mathcal{A}}(\cdot,\centerdot,\xi+\nu)\right\|_{H^s(\T^{d+1})}.
 \end{equation}
Therefore, we get in view of \eqref{Def-Norm-M1},
\begin{align}\label{Loz1}
\nonumber \left\|\Delta_\xi^{\beta}\big[\sigma_{\mathcal{A}}(\cdot,\centerdot,\xi+\eta)-\sigma_{\mathcal{A}}(\cdot,\centerdot,\xi)\big]\right\|_{H^s(\T^{d+1})}
\nonumber &\lesssim |\eta| \sup_{|\nu|\leqslant|\eta|}\langle \xi+\nu\rangle^{m_1-1-\beta} \interleave \mathcal{A}\interleave_{m_1,s,1+\beta}\\
&\lesssim |\eta|\langle \xi\rangle^{m_1-1-\beta}\interleave \mathcal{A}\interleave_{m_1,s,1+\beta}.
\end{align}
However in the region $|\eta|\geqslant 2{|\xi|}$ one deduces from the triangle inequality combined with \eqref{Def-Norm-M1}
\begin{align}\label{Loz2}
\nonumber \left\|\Delta_\xi^{\beta}\big[\sigma_{\mathcal{A}}(\cdot,\centerdot,\xi+\eta)-\sigma_{\mathcal{A}}(\cdot,\centerdot,\xi)\big]\right\|_{H^s(\T^{d+1})}&\leqslant \left\|\Delta_\xi^{\beta}\sigma_{\mathcal{A}}(\cdot,\centerdot,\xi+\eta)\right\|_{H^s(\T^{d+1})}+\left\|\Delta_\xi^{\beta}\sigma_{\mathcal{A}}(\cdot,\centerdot,\xi)\right\|_{H^s(\T^{d+1})}\\
\nonumber &\lesssim \interleave \mathcal{A}\interleave_{m_1,s,\beta}\left( \langle \xi+\eta\rangle^{m_1-\beta}+\langle \xi\rangle^{m_1-\beta}\right)
\\
&\lesssim \interleave \mathcal{A}\interleave_{m_1,s,1+\beta}\left( \langle \eta\rangle^{m_1-\beta}+\langle \xi\rangle^{m_1-\beta}\right).
\end{align}
Next, we shall  use the following  identity where we make appeal to integration by parts and change of variables
\begin{align*}
\Delta_\xi^{\beta}\, e^{\ii  \eta\,\theta }\widehat{\Lambda_\theta^{\mu_1}\sigma}_{{\mathcal{B}}}(\varphi,\eta,\xi)&=\frac{1}{2\pi}\int_{\T}\big(\Lambda_{\theta^\prime}^{\mu_1}\Delta_\xi^\beta\sigma_{{\mathcal{B}}}\big)(\varphi,\theta^\prime,\xi)e^{\ii \eta(\theta-\theta^\prime)}d\theta^\prime \\
&=\frac{1}{2\pi}\int_{\T}\big(\Lambda_{\theta}^{\mu_1}\Delta_\xi^\beta\sigma_{{\mathcal{B}}}\big)(\varphi,\theta-\theta^\prime,\xi)e^{\ii \eta\theta^\prime}d\theta^\prime.
\end{align*}
Using  Bessel identity in the variable $\theta$ one deduces
\begin{align*}
\sum_{\eta\in\Z}\left\|\Delta_\xi^\beta e^{\ii  \eta\,\centerdot }\widehat{\Lambda_\theta^{\mu_1}\sigma}_{{\mathcal{B}}}(\varphi,\eta,\xi)\right\|_{H^{s}_\theta}^2&={2\pi}\left\|\Lambda_{\theta}^{\mu_1}\Delta_\xi^\beta\sigma_{{\mathcal{B}}}(\varphi,\centerdot,\xi)\right\|_{H^s_\theta}^2
\end{align*}
and integrating in $\varphi$ yields
\begin{align*} 
\sum_{\eta\in\Z}\left\|\Delta_\xi^\beta e^{\ii  \eta\,\centerdot }\widehat{\Lambda_\theta^{\mu_1}\sigma}_{{\mathcal{B}}}(\cdot,\eta,\xi)\right\|_{L^2_\varphi H^{s}_\theta}^2&={2\pi}\left\|\Lambda_{\theta}^{\mu_1}\Delta_\xi^\beta\sigma_{{\mathcal{B}}}(\varphi,\cdot,\xi)\right\|_{L^2_\varphi H^s_\theta}^2.
\end{align*}
Similarly we get 
\begin{align*} 
\sum_{\eta\in\Z}\left\|\Delta_\xi^\beta e^{\ii  \eta\,\centerdot }\widehat{\Lambda_\theta^{\mu_1}\sigma}_{{\mathcal{B}}}(\cdot,\eta,\xi)\right\|_{H^s_\varphi L^2_\theta}^2&={2\pi}\left\|\Lambda_{\theta}^{\mu_1}\Delta_\xi^\beta\sigma_{{\mathcal{B}}}(\varphi,\cdot,\xi)\right\|_{H^s_\varphi L^2_\theta}^2.
\end{align*}
Consequently
\begin{align*} 
\sum_{\eta\in\Z}\left\|\Delta_\xi^\beta e^{\ii  \eta\,\centerdot }\widehat{\Lambda_\theta^{\mu_1}\sigma}_{{\mathcal{B}}}(\cdot,\eta,\xi)\right\|_{H^s_{\varphi,\theta} }^2&\lesssim \left\|\Lambda_{\theta}^{\mu_1}\Delta_\xi^\beta\sigma_{{\mathcal{B}}}(\varphi,\cdot,\xi)\right\|_{H^s_{\varphi,\theta}}^2\\
&\lesssim \langle \xi\rangle^{2m_2-2\beta} \interleave {\mathcal{B}}\interleave_{m_2,s+\mu_1,\beta}^2.
\end{align*}
Therefore we find 
\begin{align*}
\left\|\Delta_\xi^\beta e^{\ii  \eta\,\centerdot }\widehat{\Lambda_\theta^{\mu}\sigma}_{{\mathcal{B}}}(\cdot,\eta,\xi)\right\|_{H^{s}}&\lesssim c_\eta(\xi) \langle \xi\rangle^{m_2-\beta} \interleave {\mathcal{B}}\interleave_{m_2,s+\mu,\beta}
\end{align*}
with $\displaystyle{\sum_\eta c^2_\eta(\xi)=1},$ for any $\xi\in\Z$.
As a consequence,
\begin{align}\label{Rdd11}
\| \Delta_\xi^{\gamma-\beta}e^{\ii  \eta\,\centerdot }\widehat{\Lambda_\theta^{\mu}\sigma}_{{\mathcal{B}}}(\cdot,\centerdot,\xi+\beta)\|_{H^s}&\lesssim c_\eta(\xi) \langle \xi\rangle ^{m_2-\gamma+\beta}\,\interleave {\mathcal{B}}\interleave_{m_2,s+\mu,\gamma-\beta}.
\end{align}
Hence we get 
in view of  Cauchy-Schwarz inequality  and by fixing $\mu=\frac32+\epsilon, \epsilon>0$ and for $m_1\leqslant1$
\begin{align}\label{Yo01}
\nonumber \left\|\Delta_\xi^\gamma \mathcal{R}_1^{\textnormal{H}}(\cdot,\centerdot,\xi)\right\|_{H^{s}}&\lesssim \mathcal{C}_0\sum_{ |\eta|\leqslant \frac{|\xi|}{2}} c_\eta(\xi)|\eta|^{-\frac12-\epsilon} \langle \xi\rangle^{m_1+m_2-1-\gamma}\\
\nonumber&\quad+\mathcal{C}_0\sum_{ 
|\eta|\geqslant  2|\xi|\atop 0\leqslant\beta\leqslant\gamma}  c_\eta(\xi)|\eta|^{-\frac32-\epsilon} \langle \xi\rangle^{m_2-\gamma+\beta}\left( \langle\eta\rangle^{m_1-\beta}+\langle \xi\rangle^{m_1-\beta}\right)\\
&\lesssim\mathcal{C}_0 \langle \xi\rangle^{m_1+m_2-1-\gamma},
\end{align}
with
$$
\mathcal{C}_0=\sum_{0\leqslant\beta\leqslant\gamma}\interleave \mathcal{A}\interleave_{m_1,s,1+\beta}\interleave{\mathcal{B}}\interleave_{m_2,s_0+\frac32+\epsilon,\gamma-\beta}+\interleave \mathcal{A}\interleave_{m_1,s_0,1+\beta}\interleave{\mathcal{B}}\interleave_{m_2,s+\frac32+\epsilon,\gamma-\beta}.
$$
It remains to investigate the lower interaction term described through 
\begin{align*}
\Delta_\xi^\gamma\mathcal{R}_1^{\textnormal{L}}(\varphi,\theta,\xi)&=\sum_{\frac{|\xi|}{2}< |\eta|< 2|\xi|,\atop 0\leqslant  \beta\leqslant\gamma}\big(_\beta^\alpha\big)\left(\Delta_\xi^{\beta}{\sigma_{\mathcal{A}}(\varphi,\theta,\xi+\eta)-\Delta_\xi^{\beta}\sigma_{\mathcal{A}}(\varphi,\theta,\xi)}\right) e^{\ii  \eta\,\theta}\Delta_\xi^{\gamma-\beta}\widehat{\sigma}_{{\mathcal{B}}}(\varphi,\eta,\xi+\beta).
\end{align*}
It is plain that for any family of real numbers $\big\{\mu_\beta,0\leqslant\beta\leqslant\gamma\big\}$, we may write
\begin{align*}
\Delta_\xi^\gamma\mathcal{R}_1^{\textnormal{L}}(\varphi,\theta,\xi)&=\sum_{\frac{|\xi|}{2}< |\eta|< 2|\xi|,\atop 0\leqslant  \beta\leqslant\gamma}\big(_\beta^\alpha\big)\frac{\Delta_\xi^\beta\sigma_{\mathcal{A}}(\varphi,\theta,\xi+\eta)-\Delta_\xi^{\beta}\sigma_{\mathcal{A}(\varphi,\theta,\xi)}}{|\eta|^{\mu_\beta}} e^{\ii  \eta\,\theta}\Delta_\xi^{\gamma-\beta}\widehat{\Lambda_\theta^{\mu_\beta}\sigma}_{{\mathcal{B}}}(\varphi,\eta,\xi+\beta).
\end{align*}
Using the  definition  \eqref{Def-Norm-M1} we infer  for $\frac{|\xi|}{2}\leqslant  |\eta|\leqslant 2|\xi|$
\begin{align*}
\left\|\Delta_\xi^{\beta}\left(\sigma_{\mathcal{A}}(\cdot,\centerdot,\xi+\eta)-\sigma_{\mathcal{A}}(\cdot,\centerdot,\xi)\right)\right\|_{H^s}&\lesssim \langle \xi\rangle^{m_1^+}\interleave \mathcal{A}\interleave_{m_1,s,0}.
\end{align*}
Combining this estimate  with \eqref{Rdd11}  and fixing the choice 
\begin{align*}\mu_\beta=\frac32+\beta-m_1^-\quad\hbox{with}\quad m_1^-\triangleq \min(m_1,0),
\end{align*}
 we find 
\begin{align}\label{Yo1}
\nonumber \|\Delta_\xi^\gamma\mathcal{R}_1^{\textnormal{L}}(\cdot,\centerdot,\xi)\|_{H^s}&\lesssim\sum_{\frac{|\xi|}{2}< |\eta|< 2|\xi|\atop 0\le\beta\le\gamma}  c_\eta(\xi)|\eta|^{-\mu_\beta}\langle \xi\rangle^{m_1^++m_2-\gamma+\beta}\,\mathcal{C}_{1,\beta}\\
\nonumber &\lesssim \langle \xi\rangle^{m_1+m_2-\gamma-1}\sum_{\frac{|\xi|}{2}< |\eta|< 2|\xi|\atop 0\leqslant\beta\leqslant\gamma}  c_\eta(\xi)|\eta|^{-\frac12}\,\mathcal{C}_{1,\beta}\\
& \lesssim\mathcal{C}_1 \langle \xi\rangle^{m_1+m_2-\gamma-1}
\end{align}
with
$$
\mathcal{C}_{1,\beta}\triangleq \interleave \mathcal{A}\interleave_{m_1,s,0}\interleave {\mathcal{B}}\interleave_{m_2,s_0+\frac32+\beta-m_1^-,\gamma-\beta} +\interleave \mathcal{A}\interleave_{m_1,s_0,0}\interleave {\mathcal{B}}\interleave_{m_2,s+\frac32+\beta-m_1^-,\gamma-\beta} 
$$
and
\begin{align*}
\mathcal{C}_1&\triangleq \sum_{ 0\leqslant\beta\leqslant\alpha}\mathcal{C}_{1,\beta}.
\end{align*}
Notice that we have used in the last inequality of \eqref{Yo1} Cauchy-Schwarz inequality as follows
\begin{align*}
\sum_{\frac{|\xi|}{2}< |\eta|< 2|\xi|}  c_\eta(\xi)|\eta|^{-\frac12}&\leqslant\Big(\sum_{\frac{|\xi|}{2}< |\eta|< 2|\xi|}  c_\eta^2(\xi)\Big)^{\frac12}\Big(\sum_{\frac{|\xi|}{2}< |\eta|< 2|\xi|}  |\eta|^{-1}\Big)^{\frac12}\\
&\leqslant \Big(\sum_{\frac{|\xi|}{2}< |\eta|< 2|\xi|}  |\eta|^{-1}\Big)^{\frac12}\\
&\leqslant C
\end{align*}
with $C$ being  independent of $\eta.$
Combining \eqref{Yo01} and \eqref{Yo1} we get for $m_1\leqslant1$ 
 \begin{align*}
\|\Delta_\xi^\gamma\mathcal{R}_1^{\textnormal{L}}(\cdot,\centerdot,\xi)\|_{H^s}
&\lesssim \big(\mathcal{C}_0+\mathcal{C}_1\big)\, \langle \xi\rangle ^{m_1+m_2-\gamma-1}.\end{align*}
Consequently, we get
\begin{align*}
\interleave\mathcal{R}_1^{\textnormal{L}}\interleave_{m_1+m_2-1,s,\gamma}\lesssim \mathcal{C}_0+\mathcal{C}_1.
\end{align*}
In a similar way to \eqref{abou-s1} one gets, by exchanging $\mathcal{A}$ and ${\mathcal{B}}$, the decomposition
\begin{align*}
\sigma_{{\mathcal{B}}\mathcal{A}}(\varphi,\theta,\xi)
&\triangleq \mathcal{R}_2(\varphi,\theta,\xi)+\sigma_{\mathcal{A}}(\varphi,\theta,\xi)\sigma_{{\mathcal{B}}}(\varphi,\theta,\xi)
\end{align*}
with the estimate, under the assumption $m_2\leqslant 1$,
\begin{align*}
\interleave \mathcal{R}_2\interleave_{m_1+m_2-1,s,\gamma}&\lesssim \mathcal{C}^\prime_0+\mathcal{C}^\prime_1\end{align*}
where
$$
\mathcal{C}^\prime_0\triangleq \sum_{0\leqslant\beta\leqslant\gamma}\interleave {\mathcal{B}}\interleave_{m_2,s,1+\beta}\interleave\mathcal{A}\interleave_{m_1,s_0+\frac32+\epsilon,\gamma-\beta}+\interleave {\mathcal{B}}\interleave_{m_2,s_0,1+\beta}\interleave\mathcal{A}\interleave_{m_1,s+\frac32+\epsilon,\gamma-\beta}.
$$
and
\begin{align*}
\mathcal{C}^\prime_1&\triangleq \sum_{ 0\leqslant\beta\leqslant\gamma}\left(\interleave {\mathcal{B}}\interleave_{m_2,s,0}\interleave \mathcal{A}\interleave_{m_1,s_0+\frac32+\beta-m_2^-,\gamma-\beta} +\interleave {\mathcal{B}}\interleave_{m_2,s_0,0}\interleave \mathcal{A}\interleave_{m_1,s+\frac32+\beta-m_2^-,\gamma-\beta} \right).
\end{align*}
This achieves the desired estimate for  the commutator.

\smallskip

{$($\bf{iii}$)$} We proceed in a similar way to the preceding point $($ii$)$. 
Putting together \eqref{Rdd11}, \eqref{Loz1} and \eqref{Loz2} yields ( similarly to \eqref{Yo01}) under the assumptions 
\begin{align}\label{Assum-mu}
\mu>\frac12+\max\big(0,m_1,m_2,m_1+m_2\big),\,\quad  m_1+m_2\leq 1
\end{align}
 and in view of  Cauchy-Schwarz inequality to 
\begin{align*}
\left\|\Delta_\xi^\gamma \mathcal{R}_1^{\textnormal{H}}(\cdot,\centerdot,\xi)\right\|_{H^{s}}&\lesssim \mathcal{C}_2\sum_{ |\eta|\leqslant \frac{|\xi|}{2}} c_\eta(\xi)|\eta|^{1-\mu} \langle \xi\rangle^{m_1+m_2-1-\gamma}\\
&\quad+\mathcal{C}_2\sum_{ 
|\eta|\geqslant  2|\xi|\atop 0\leqslant\beta\leqslant\gamma} c_\eta(\xi)|\eta|^{-\mu} \langle \xi\rangle^{m_2-\gamma+\beta}\left( \langle\eta\rangle^{m_1-\beta}+\langle \xi\rangle^{m_1-\beta}\right)\\
&\lesssim \mathcal{C}_2 \langle \xi\rangle^{-\gamma},
\end{align*}
with
$$
\mathcal{C}_2=\sum_{0\leqslant\beta\leqslant\gamma}\interleave \mathcal{A}\interleave_{m_1,s,1+\beta}\interleave{\mathcal{B}}\interleave_{m_2,s_0+\mu,\gamma-\beta}+\interleave \mathcal{A}\interleave_{m_1,s_0,1+\beta}\interleave{\mathcal{B}}\interleave_{m_2,s+\mu,\gamma-\beta}.
$$
As to the lower interaction term, we proceed as in  \eqref{Yo1}. One gets
in view of  \eqref{Rdd11},  Cauchy-Schwarz inequality and with the choice
\begin{align}\label{Assum-mu1}\mu_\beta= m_1^++m_2+\frac12+\beta
\end{align}
 that
\begin{align*}
\|\Delta_\xi^\gamma\mathcal{R}_1^{\textnormal{L}}(\cdot,\centerdot,\xi)\|_{H^s}&\lesssim\sum_{\frac{|\xi|}{2}< |\eta|< 2|\xi|\atop 0\leqslant\beta\leqslant\gamma}  c_\eta(\xi)|\eta|^{-\mu_\beta}\langle \xi\rangle^{m_1^++m_2-\gamma+\beta}\,\mathcal{C}_{3,\beta}\\
&\lesssim \langle \xi\rangle^{-\gamma}\sum_{\frac{|\xi|}{2}< |\eta|< 2|\xi|\atop 0\leqslant\beta\leqslant\gamma}  c_\eta(\xi)|\eta|^{-\mu_\beta+m_1+m_2+\beta}\,\mathcal{C}_{3,\beta}\\
&\lesssim\mathcal{C}_3 \langle \xi\rangle^{-\gamma}
\end{align*}
with
\begin{align*}
\mathcal{C}_3&\triangleq \sum_{ 0\leqslant\beta\leqslant\gamma}\mathcal{C}_{3,\beta},\\
\mathcal{C}_{3,\beta}&\triangleq \interleave \mathcal{A}\interleave_{m_1,s,0}\interleave {\mathcal{B}}\interleave_{m_2,s_0+\mu_\beta,\gamma-\beta} +\interleave \mathcal{A}\interleave_{m_1,s_0,0}\interleave {\mathcal{B}}\interleave_{m_2,s+\mu_\beta,\gamma-\beta}.
\end{align*}
Putting together  the preceding estimates we get under the assumptions \eqref{Assum-mu} and \eqref{Assum-mu1}
\begin{align*}
\interleave \mathcal{R}_1\interleave_{0,s,\gamma}&\lesssim\mathcal{C}_2+\mathcal{C}_3 \\
&\lesssim\sum_{0\leqslant\beta\leqslant\gamma}\interleave \mathcal{A}\interleave_{m_1,s,1+\gamma-\beta}\interleave{\mathcal{B}}\interleave_{m_2,s_0+\mu,\beta}+\interleave \mathcal{A}\interleave_{m_1,s_0,1+\gamma-\beta}\interleave{\mathcal{B}}\interleave_{m_2,s+\mu,\beta}
\\
&\quad +\sum_{ 0\leqslant\beta\leqslant\gamma}\left(\interleave \mathcal{A}\interleave_{m_1,s,0}\interleave {\mathcal{B}}\interleave_{m_2,s_0+\mu_\beta,\gamma-\beta} +\interleave \mathcal{A}\interleave_{m_1,s_0,0}\interleave {\mathcal{B}}\interleave_{m_2,s+\mu_\beta,\gamma-\beta} \right).
  \end{align*}
Similarly, we get
\begin{align*}
\sigma_{{\mathcal{B}}\mathcal{A}}(\varphi,\theta,\xi)
&=\sum_{\eta\in \Z}e^{\ii  \eta\,\theta}\Big(\sigma_{{\mathcal{B}}}(\varphi,\theta,\xi+\eta)-\sigma_{{\mathcal{B}}}(\varphi,\theta,\xi)\Big)\widehat{\sigma}_{\mathcal{A}}(\varphi,\eta,\xi)\\
&\quad+\sigma_{{\mathcal{B}}}(\varphi,\theta,\xi)\sigma_{\mathcal{A}}(\varphi,\theta,\xi)\\
&\triangleq \mathcal{R}_2(\varphi,\theta,\xi)+\sigma_{\mathcal{A}}(\varphi,\theta,\xi)\sigma_{{\mathcal{B}}}(\varphi,\theta,\xi)
\end{align*}
and $\mathcal{R}_2$ satisfies under the assumptions \eqref{Assum-mu} and \eqref{Assum-mu1} the estimate
\begin{align*}
\interleave \mathcal{R}_2\interleave_{0,s,\gamma}&\lesssim\sum_{0\leqslant\beta\leqslant \gamma}\interleave {\mathcal{B}}\interleave_{m_1,s,1+\gamma-\beta}\interleave\mathcal{A}\interleave_{m_2,s_0+\mu,\beta}+\interleave {\mathcal{B}}\interleave_{m_1,s_0,1+\alpha-\beta}\interleave\mathcal{A}\interleave_{m_2,s+\mu,\beta}
\\
&\quad+\sum_{ 0\leqslant \beta\leqslant\gamma}\left(\interleave {\mathcal{B}}\interleave_{m_1,s,0}\interleave \mathcal{A}\interleave_{m_2,s_0+\overline{\mu}_\beta,\gamma-\beta} +\interleave {\mathcal{B}}\interleave_{m_1,s_0,0}\interleave \mathcal{A}\interleave_{m_2,s+\overline{\mu}_\beta,\gamma-\beta} \right)
 \end{align*}
 with
 $$
 \overline{\mu}_\beta= m_1+m_2^++\frac12+\beta.
 $$
Since $\sigma_{[\mathcal{A},{\mathcal{B}}]}=\sigma_{\mathcal{A}{\mathcal{B}}}-\sigma_{{\mathcal{B}}\mathcal{A}}=\mathcal{R}_1-\mathcal{R}_2$ then 
$$
\interleave\sigma_{[\mathcal{A},{\mathcal{B}}]}\interleave_{0,s,\gamma}\le \interleave\mathcal{R}_1\interleave_{0,s,\gamma}+\interleave\mathcal{R}_2\interleave_{0,s,\gamma}
$$
and the result  follows from the preceding estimates.
 \end{proof}

\subsection{Commutators}
Our purpose in this section is to establish some commutator estimates between particular pseudo-differential operators of the following type: For a smooth function $h$ we  define the operator
$$
\mathbb{A}_\rho\triangleq \partial_\theta\big(\rho   |\textnormal D|^{\alpha-1}+|\textnormal D|^{\alpha-1}\rho\big),
$$
where  the action of the fractional Laplacian $|\textnormal{D}|^{\alpha-1}$ is given by
$$
|\textnormal{D}|^{\alpha-1}\rho(\lambda,\varphi,\theta)\triangleq \frac{1}{2\pi}\bigintsss_{0}^{2\pi} \frac{\rho(\lambda,\varphi,\theta-\eta)}{|\sin\big(\frac{\eta}{2}\big)|^{\alpha}} d\eta, \quad\lambda=(\omega,\alpha)\in\mathcal{O}.
$$
In  this section $\mathcal{O}$ denotes an open set of $\RR^{d+1}$ taking the form
$$
\mathcal{O}=\mathscr{U}\times(0,\overline\alpha)\subset\RR^{d+1}, \overline\alpha\in(0,\tfrac12)
$$
with $\mathscr{U}$ is an open bounded set of $\RR^{d}.$
We shall first detail the action of fractional Laplacian $|\textnormal{D}|^{\alpha-1}$ on the spaces $W^{q,\infty}_{\kappa}(\mathcal{O},H^{s}(\T^{d+1}))$ introduced in Definition \ref{Def-WS}. 
\begin{lemma} \label{Laplac-frac} Let $s \in\RR, q,\gamma,\ell\in\NN, \kappa\in(0,1),$ then the  following assertions hold true.
\begin{enumerate}
\item For any $ \epsilon>0$, there exists $C>0$ such that
$$
\interleave \partial_\theta^\ell|\textnormal D|^{\alpha-1}\interleave_{\overline\alpha+\ell+\epsilon-1,s,\gamma}^{q,\kappa}\leqslant C.
$$
\item For any $ \epsilon>0$, there exists $C>0$ 
\begin{align*}
\interleave\partial_\theta  |\textnormal D|^{\alpha-1} |\textnormal D|^{\alpha-1}\interleave_{2\overline\alpha+\epsilon-1,s,\gamma}^{q,\kappa}\leqslant C.
\end{align*}
\end{enumerate}\end{lemma}
\begin{proof}
${\bf{(i)}}$ According to \eqref{Fo-c} and the formula $x\Gamma(x)=\Gamma(1+x)$ one finds 
\begin{align}\label{Fo-c1}
\nonumber |\textnormal D|^{\alpha-1}{\bf e}_{n}(\theta)=&\,\frac{2^{\alpha}\Gamma(1-\alpha)}{\Gamma(1-\frac\alpha2)\Gamma(\frac{\alpha}{2})}\frac{\Gamma(|n|+\frac{\alpha}{2})}{\Gamma(|n|+1-\frac{\alpha}{2})}{\bf e}_{n}(\theta)\\
\nonumber=&\frac{2^{\alpha-1}\alpha\Gamma(1-\alpha)}{\Gamma(1-\frac\alpha2)\Gamma(1+\frac{\alpha}{2})}\frac{\Gamma(|n|+\frac{\alpha}{2})}{\Gamma(|n|+1-\frac{\alpha}{2})}{\bf e}_{n}(\theta)\\
=&\mu(\alpha)\mathtt{W}(|n|,\alpha){\bf e}_{n}(\theta),
\end{align}
where $\mathtt{W}(n,\alpha)$ was  defined in Lemma \ref{lem-asym}-(iii).
Then the symbol of $\mathcal{A}\triangleq \partial_\theta^\ell|\textnormal D|^{\alpha-1}$ is given by
$$
\sigma_{\mathcal{A}}(\lambda,\varphi,\theta,n)=(\ii n)^\ell \mu(\alpha)\mathtt{W}(|n|,\alpha).
$$
 From classical properties of Gamma function we deduce that  $\mu$ is smooth in $\alpha$ and
$$
\sup_{0\leqslant k\leqslant j\\\atop
\alpha\in(0,\overline\alpha)}|\mu^{(k)}(\alpha)|<\infty.
$$
According to Lemma  \ref{Stirling-formula}-{\rm (ii)} one has
\begin{align*}
\forall\, n,\gamma\in\NN,\quad\max_{k\in\llbracket 0,{q}\rrbracket\atop \alpha\in (0,\overline\alpha) }  |\partial_{\alpha}^k\Delta_n^\gamma\mathtt{W}(n,\alpha)|&\leqslant C(\gamma,q,\epsilon)\,\langle n\rangle^{\overline\alpha+\epsilon-1-{\gamma}}.
\end{align*}
Therefore,  according to  Leibniz rule combined with the preceding estimates  we obtain 
\begin{align}\label{Der-E1}
\forall\, n,\gamma\in\NN,\quad \quad \max_{k\in\llbracket 0,{q}\rrbracket\atop \alpha\in (0,\overline\alpha) }\langle n\rangle^{1-\ell-\overline\alpha-\epsilon+\gamma}\big|\partial_\alpha^k\Delta_n^\gamma\sigma_{\mathcal{A}}(\lambda,\varphi,\theta,n)\big|\leqslant C(\gamma,q,\ell,\epsilon).
\end{align}
This implies in view of \eqref{Def-pseud-w}
$$
\interleave  \partial_\theta^\ell|\textnormal D|^{\alpha-1}\interleave_{\overline\alpha+\ell+\epsilon-1,s,\gamma}^{q,\kappa}\leqslant C(\gamma,q,\ell,\epsilon),$$
which ensures the desired result. We point out that  for $q=0$ we have no the loss $\epsilon$. More precisely,
\begin{align}\label{EPPLL1}
\forall \alpha\in(0,\overline\alpha),\quad \interleave  \partial_\theta^\ell|\textnormal D|^{\alpha-1}\interleave_{\alpha+\ell-1,s,\gamma}\leqslant C(\gamma,\ell,\epsilon).
\end{align}

\smallskip

${\bf{(ii)}}$ Iterating  \eqref{Fo-c1} we obtain
\begin{align*}
 \partial_\theta  |\textnormal D|^{\alpha-1} |\textnormal D|^{\alpha-1}{\bf e}_{n}
&=\ii n \mu^2(\alpha)\, \mathtt{W}^2(|n|,\alpha)\,{\bf e}_{n}.
\end{align*}
Similarly to \eqref{Der-E1} we deduce that for any $\epsilon>0,$ 
\begin{align*}
\forall\, n,\gamma\in\NN,\quad \quad \max_{k\in\llbracket 0,{q}\rrbracket\atop \alpha\in (0,\overline\alpha) }\langle n\rangle^{1-2\overline\alpha-\epsilon+\gamma}\big|\partial_\alpha^k\Delta_n^\gamma\big(n \mu^2(\alpha)\, \mathtt{W}^2(|n|,\alpha)\big|\leqslant C(\gamma,q,\epsilon).
\end{align*}
It follows from \eqref{Def-pseud-w} that
$$
\interleave \partial_\theta  |\textnormal D|^{\alpha-1}|\textnormal D|^{\alpha-1}\interleave_{2\overline\alpha+\epsilon-1,s,\gamma}^{q,\kappa}\leqslant C(\gamma,q,\epsilon),$$
which  achieves the proof of Lemma \ref{Laplac-frac}.
\end{proof}
As an application, we shall prove the following result.
\begin{lemma}\label{Lem-Commutator}
Let $ \epsilon>0$ be small enough, $\overline\alpha\in(0,\frac12)$ and $(\kappa,d,q,s,s_0)$  as in \eqref{initial parameter condition}. Then the  following assertions hold true. 
\begin{enumerate}
\item For any $\alpha\in(0,\overline\alpha)$, there exists $C>0$ such that
$$
\interleave \mathbb{A}_\rho \interleave_{\alpha,s,\gamma}\leqslant C \|\rho\|_{s+\gamma+1}.
$$

\item There exists $C>0$ such that
$$
\interleave\mathbb{A}_\rho \interleave_{\overline\alpha+\epsilon,s,\gamma}^{q,\kappa}\leqslant C \|\rho\|_{s+\gamma+1}^{q,\kappa}.
$$
\item There exists $C>0$ such that
\begin{eqnarray*}
\nonumber&&\interleave [\mathbb{A}_{\rho_1},\mathbb{A}_{\rho_2}]\interleave_{2\overline\alpha-1+\epsilon,s,\gamma}^{q,\kappa}\leqslant C\big(\|\rho_1\|_{s_0+3}^{q,\kappa}+\|\rho_2\|_{s_0+3}^{q,\kappa}\big)\big(\|\rho_1\|_{s+\gamma+3}^{q,\kappa}+\|\rho_2\|_{s+\gamma+3}^{q,\kappa}\big).
\end{eqnarray*}
\item  Let   $\mathcal{A}$  as in \eqref{A-phi-lambda} and $m\in[-1,0]$. Then for any $s\geqslant   s_0$ there exits $C>0$ 
\begin{align*}
\interleave[\mathcal{A},\mathbb{A}_{\rho}]\interleave_{m,s,\gamma}^{q,\kappa}&\lesssim
\sum_{0\leqslant\beta\leqslant\gamma}\interleave \mathcal{A}\interleave_{m,s_0+2,1+\beta}^{q,\kappa}\|\rho\|_{s+4+\gamma-\beta}^{q,\kappa}+\interleave \mathcal{A}\interleave_{m,s+2,1+\beta}^{q,\kappa}\rho\|_{s_0+4+\gamma-\beta}^{q,\kappa}\\
&+\sum_{ 0\leqslant\beta\leqslant\gamma}\interleave \mathcal{A}\interleave_{m,s_0+2+\beta,\gamma-\beta}^{q,\kappa}\|\rho\|_{s+1}^{q,\kappa} +\interleave \mathcal{A}\interleave_{m,s+2+\beta,\gamma-\beta}^{q,\kappa}\|\rho\|_{s_0+1}^{q,\kappa}.
\end{align*}
\item  Let   $\overline{s}_0>\frac{d+5}{2}, s\geqslant0,$ $\rho \in H^{s+\overline{s}_0}$ and $h\in H^{s}\cap H^{\overline{s}_0}$. Then
$$
\left\|\partial_\theta\left[\rho , \Lambda^{s}\right] h\right\|_{L^2(\T^{d+1})}\lesssim \|h\|_{H^{s}(\T^{d+1})}\|\rho\|_{H^{\overline{s}_0}(\T^{d+1})}+ \|h\|_{H^{\overline{s}_0}(\T^{d+1})}\|\rho\|_{H^{s+\overline{s}_0}(\T^{d+1})}.
$$

\end{enumerate}
\end{lemma}
\begin{proof}
{\bf{(i)}} We write 
\begin{align*}
\mathbb{A}_\rho&=\rho  \partial_\theta|\textnormal D|^{\alpha-1}+(\partial_\theta \rho)   |\textnormal D|^{\alpha-1}+(\partial_\theta|\textnormal D|^{\alpha-1})\rho\\
&\triangleq \mathbb{A}_{\rho,1}+\mathbb{A}_{\rho,2}+\mathbb{A}_{\rho,3}.
\end{align*}
Denote by $T_\rho$ the multiplicative operator defined by $T_\rho h=\rho h$ then $\sigma_{T_\rho}=\rho$ and therefore 
$$
\interleave  T_\rho \interleave_{0,s,\gamma}^{q,\kappa}\lesssim \|\rho\|_{s}^{q,\kappa}.
$$
Applying  Lemma \ref{comm-pseudo}-(i) combined with \eqref{EPPLL1} we infer for any $s\geqslant s_0$
\begin{align*}
\interleave \mathbb{A}_{\rho,1}\interleave_{\alpha,s,\gamma}&\lesssim\|\rho\|_{s}.
\end{align*}
Similarly we get
\begin{align*}
\interleave \mathbb{A}_{\rho,2}\interleave_{\alpha,s,\gamma}&\lesssim \interleave T_{\partial_\theta \rho} \interleave_{0,s,\gamma}\\
&\lesssim \|\rho\|_{s+1}.
\end{align*}
For the last term we use once again Lemma \ref{comm-pseudo}-(i) combined with \eqref{EPPLL1} and the fact that $\alpha\in(0,\frac12),$ in order to get for any $s\geqslant s_0$
\begin{align*}
\interleave\mathbb{A}_{\rho,3}\interleave_{\alpha,s,\gamma}&\lesssim\|\rho\|_{s+\frac12+\alpha}+\|\rho\|_{s+\frac12+\gamma}\\
&\lesssim\|\rho\|_{s+\gamma+1}.
\end{align*}
This achieves the proof of the first point.

\smallskip

{\bf{(ii)}} With the notations of the first point one has according to  Lemma \ref{comm-pseudo}-(i) combined with Lemma \ref{Laplac-frac} we infer for any $s\geqslant s_0$
\begin{align*}
\forall\, s\geqslant s_0,\quad\interleave\mathbb{A}_{\rho,1}\interleave_{\overline\alpha+\epsilon,s,\gamma}^{q,\kappa}&\lesssim\|\rho\|_{s}^{q,\kappa}.
\end{align*} 
Similar arguments yield
\begin{align*}
\interleave \mathbb{A}_{\rho,2}\interleave_{\overline\alpha+\epsilon,s,\gamma}^{q,\kappa}
&\lesssim\interleave T_{\partial_\theta \rho} \interleave_{0,s,\gamma}^{q,\kappa}\\
&\lesssim \|\rho\|_{s+1}^{q,\kappa}.
\end{align*}
For the last term we implement  once again Lemma \ref{comm-pseudo}-(i) combined with Lemma \ref{Laplac-frac} and the fact that $\overline\alpha\in(0,\frac12),$ in order to get for any $s\geqslant s_0$
\begin{align*}
\interleave \mathbb{A}_{\rho,3}\interleave_{\overline\alpha+\epsilon,s,\gamma}^{q,\kappa}&\lesssim\|\rho\|_{s+\frac12+\overline\alpha+\epsilon}^{q,\kappa}+\|\rho\|_{s+\frac12+\gamma}^{q,\kappa}\\
&\lesssim\|\rho\|_{s+\gamma+1}^{q,\kappa},
\end{align*}
where $\epsilon$ is taken small enough.
 Putting together the preceding estimates gives
\begin{equation}\label{Hm-id1}
\interleave \mathbb{A}_\rho \interleave_{\overline\alpha+\epsilon,s,\gamma}^{q,\kappa}\leqslant C \|\rho\|_{s+1+\gamma}^{q,\kappa},
\end{equation}
which completes the proof of the desired result. 

\smallskip

${(\bf{iii})}$ Applying Lemma \ref{comm-pseudo}-(ii) with $m_1=m_2=\overline\alpha+\frac{\epsilon}{2}$ and using  the estimate \eqref{Hm-id1} combined with Sobolev embeddings
\begin{align*}
\interleave[\mathbb{A}_{\rho_1},\mathbb{A}_{\rho_2}]\interleave_{2\overline\alpha-1+\epsilon,s,\gamma}^{q,\kappa}&\lesssim
\sum_{0\leqslant\beta\leqslant\gamma} \|\rho_1\|_{s+2+\beta}^{q,\kappa}\ \|\rho_2\|_{s_0+3+\gamma-\beta}^{q,\kappa}+ \|\rho_1\|_{s_0+2+\beta}^{q,\kappa}\ \|\rho_2\|_{s+3+\gamma-\beta}^{q,\kappa}\\
&+\sum_{ 0\leqslant\beta\leqslant\gamma}\|\rho_1\|_{s+1}^{q,\kappa}\ \|h_2\|_{s_0+3+\gamma}^{q,\kappa}+\|h_1\|_{s_0+1}^{q,\kappa}\ \|\rho_2\|_{s+3+\gamma}^{q,\kappa} \\
&+
\sum_{0\leqslant\beta\leqslant\gamma} \|\rho_2\|_{s+2+\beta}^{q,\kappa}\ \|\rho_1\|_{s_0+3+\gamma-\beta}^{q,\kappa}+ \|\rho_2\|_{s_0+2+\beta}^{q,\kappa}\ \|\rho_1\|_{s+3+\gamma-\beta}^{q,\kappa}\\
&+\sum_{ 0\leqslant\beta\leqslant\gamma}\|\rho_2\|_{s+1}^{q,\kappa}\ \|\rho_1\|_{s_0+3+\gamma}^{q,\kappa}+\|\rho_2\|_{s_0+1}^{q,\kappa}\ \|\rho_1\|_{s+3+\gamma}^{q,\kappa}.
\end{align*}
Hence we find from Sobolev embeddings
\begin{align}\label{Uw-X1}
\interleave[\mathbb{A}_{\rho_1},\mathbb{A}_{\rho_2}]\interleave_{2\overline\alpha-1+\epsilon,s,\gamma}^{q,\kappa}&\lesssim
\sum_{0\leqslant\beta\leqslant\gamma}\Big( \|\rho_1\|_{s+3+\beta}^{q,\kappa}\ \|\rho_2\|_{s_0+3+\gamma-\beta}^{q,\kappa}+ \|\rho_1\|_{s_0+3+\beta}^{q,\kappa}\ \|\rho_2\|_{s+3+\gamma-\beta}^{q,\kappa}\Big).
\end{align}
Applying interpolation inequality we find 
$$
 \|\rho_1\|_{s+3+\beta}^{q,\kappa}\ \|\rho_2\|_{s_0+3+\gamma-\beta}^{q,\kappa}\lesssim \big(\|\rho_1\|_{s_0+3}^{q,\kappa}+ \|\rho_2\|_{s_0+3}^{q,\kappa}\big)\big(\|\rho_1\|_{s+3+\gamma}^{q,\kappa}+ \|\rho_2\|_{s+3+\gamma}^{q,\kappa}\big),
$$
which implies by virtue of \eqref{Uw-X1}
\begin{align*}
\interleave[\mathbb{A}_{\rho_1},\mathbb{A}_{\rho_2}]\interleave_{2\overline\alpha-1+\epsilon,s,\gamma}^{q,\kappa}&\lesssim\big(\|\rho_1\|_{s_0+3}^{q,\kappa}+ \|\rho_2\|_{s_0+3}^{q,\kappa}\big)\big(\|\rho_1\|_{s+3+\gamma}^{q,\kappa}+ \|\rho_2\|_{s+3+\gamma}^{q,\kappa}\big).
\end{align*}
This ends the proof of the third point.
\smallskip

${(\bf{iv})}$ First we write according to the definition and since $\overline\alpha+\epsilon-1\leqslant 0$
$$
\interleave[\mathcal{A} , \mathbb{A}_{h}] \interleave_{m,s,\gamma}^{q,\kappa} \lesssim \interleave[\mathcal{A} , \mathcal{U}_{h}] \interleave_{m+\overline\alpha+\epsilon-1,s,\gamma}^{q,\kappa}.
$$
Then applying Lemma \ref{comm-pseudo}-(ii) combined with \eqref{Hm-id1} and Sobolev embeddings
\begin{align*}
\interleave[\mathcal{A},\mathbb{A}_{\rho}]\interleave_{m+\overline\alpha+\epsilon-1,s,\gamma}^{q,\kappa}&\lesssim
\sum_{0\leqslant\beta\leqslant\gamma}\interleave \mathcal{A}\interleave_{m,s,1+\beta}^{q,\kappa}\|\rho\|_{s_0+3+\gamma-\beta}^{q,\kappa}+\interleave \mathcal{A}\interleave_{m,s_0,1+\beta}^{q,\kappa}\|\rho\|_{s+3+\gamma-\beta}^{q,\kappa}\\
&+\interleave \mathcal{A}\interleave_{m,s,0}^{q,\kappa}\|\rho\|_{s_0+4+\gamma}^{q,\kappa} +\interleave \mathcal{A}\interleave_{m,s_0,0}^{q,\kappa}\|\rho\|_{s+4+\gamma}^{q,\kappa} \\
&+
\sum_{0\leqslant\beta\leqslant\gamma} \interleave\mathcal{A}\interleave_{m,s_0+2,\beta}^{q,\kappa}\|\rho\|_{s+2+\gamma-\beta}^{q,\kappa}+\interleave\mathcal{A}\interleave_{m,s+2,\beta}^{q,\kappa}\|\rho\|_{s_0+2+\gamma-\beta}^{q,\kappa}\\
&+\sum_{ 0\leqslant\beta\leqslant\gamma}\interleave \mathcal{A}\interleave_{m,s_0+2+\beta,\gamma-\beta}^{q,\kappa}\|\rho\|_{s+1}^{q,\kappa} +\interleave \mathcal{A}\interleave_{m,s+2+\beta,\gamma-\beta}^{q,\kappa}\|\rho\|_{s_0+1}^{q,\kappa}.
\end{align*}
Applying Sobolev embeddings we infer
\begin{align*}
\interleave[\mathcal{A},\mathbb{A}_{\rho}]\interleave_{m,s,\gamma}^{q,\kappa}&\lesssim
\sum_{0\leqslant\beta\leqslant\gamma}\interleave \mathcal{A}\interleave_{m,s_0+2,1+\beta}^{q,\kappa}\|\rho\|_{s+4+\gamma-\beta}^{q,\kappa}+\interleave \mathcal{A}\interleave_{m,s+2,1+\beta}^{q,\kappa}\|\rho\|_{s_0+4+\gamma-\beta}^{q,\kappa}\\
&+\sum_{ 0\leqslant\beta\leqslant\gamma}\interleave \mathcal{A}\interleave_{m,s_0+2+\beta,\gamma-\beta}^{q,\kappa}\|\rho\|_{s+1}^{q,\kappa} +\interleave \mathcal{A}\interleave_{m,s+2+\beta,\gamma-\beta}^{q,\kappa}\|\rho\|_{s_0+1}^{q,\kappa}.
\end{align*}
This proves the desired result.

\smallskip

{\bf{(v)}}  Denote $f\triangleq \partial_\theta\left[\rho, \Lambda^{s}\right] h,$ then  its Fourier  expansion writes
$$
f=\sum_{(l,j)\in\Z^{d+1}}c_{l,j}(f){\bf{e}}_{l,j},\quad {\bf{e}}_{l,j}(\varphi,\theta)=e^{\ii(l\cdot \varphi+j\theta)}
$$
with
$$
c_{l,j}(f)=\ii j \sum_{(l^\prime,j^\prime)\in\Z^{d+1}}\mu_{l^\prime, j^\prime}^{l,j}\,c_{l^\prime,j^\prime}(h) c_{l-l^\prime,j-j^\prime}(\rho)
$$
and
$$
\mu_{l^\prime, j^\prime}^{l,j}=\langle l^\prime,j^\prime\rangle^{s}-{\langle l,j\rangle^{s}}.
$$
Using  the classical  inequality
\begin{align}\label{PS0}
\forall a,b>0,\quad  \big| a^{s}-b^{s}\big|\lesssim_s |a-b|\big(a^{s-1}+|a-b|^{s-1}\big)
\end{align}
with $a=\langle l^\prime,j^\prime\rangle$ and $ b=\langle l,j\rangle$ yields from the  triangle inequality
\begin{align*}
\big|\langle l^\prime,j^\prime\rangle^{s_1}-\langle l,j\rangle^{s_1} \big|&\lesssim \langle l-l^\prime,j-j^\prime\rangle \Big(\langle l^\prime,j^\prime\rangle^{s-1}+\langle l-l^\prime,j-j^\prime\rangle^{s-1}\Big)\\
&\lesssim \langle l-l^\prime,j-j^\prime\rangle \langle l^\prime,j^\prime\rangle^{s-1}+\langle l-l^\prime,j-j^\prime\rangle^{s}.
\end{align*}
It follows that
\begin{align*}
 |\mu_{l^\prime, j^\prime}^{l,j}|&\lesssim \langle l-l^\prime,j-j^\prime\rangle \langle l^\prime,j^\prime\rangle^{s-1}+\langle l-l^\prime,j-j^\prime\rangle^{s}.
\end{align*}
Combining this inequality with $|j|\leqslant |j-j^\prime|+|j^\prime|$ and $|j^\prime|\leqslant  \langle l^\prime,j^\prime\rangle$  allows to get for $s\geqslant0$
\begin{align}\label{PS1}
\nonumber |j| |\mu_{l^\prime, j^\prime}^{l,j}|\lesssim&\, \langle l-l^\prime,j-j^\prime\rangle^2 \langle l^\prime,j^\prime\rangle^{s-1}+\langle l-l^\prime,j-j^\prime\rangle^{s+1} \\
&+\langle l-l^\prime,j-j^\prime\rangle \langle l^\prime,j^\prime\rangle^s+\langle l-l^\prime,j-j^\prime\rangle^{s} \langle l^\prime,j^\prime\rangle.
\end{align}
Using the convolution laws allows to get
\begin{align*}
\|f\|_{L^2}\lesssim& \|h\|_{H^{s-1}}\sum_{l,j}\langle l,j\rangle^2|c_{l,j}(\rho)|+ \|h\|_{L^2}\sum_{l,j}\langle l,j\rangle^{s+1}|c_{l,j}(\rho)|\\
&+\|h\|_{H^s}\sum_{l,j}\langle l,j\rangle|c_{l,j}(\rho)+\|h\|_{H^1}\sum_{l,j}\langle l,j\rangle^{s}|c_{l,j}(\rho)|.
\end{align*}
Then using Cauchy-Schwarz inequality leads for any $\epsilon>0$
\begin{align*}
\|f\|_{L^2}\lesssim& \|h\|_{H^{s-1}}\|\rho\|_{H^{\frac{d+5}{2}+\epsilon}(\T^{d+1})}+ \|h\|_{L^2}\|\rho\|_{H^{s+\frac{d+3}{2}+\epsilon}}\\
&+ \|h\|_{H^{s}}\|\rho\|_{H^{\frac{d+3}{2}+\epsilon}(\T^{d+1})}+ \|h\|_{H^1}\|\rho\|_{H^{s+\frac{d+1}{2}+\epsilon}}.
\end{align*}
To get the desired inequality, it is enough to use Sobolev embeddings and this achieves the proof of Lemma \ref{Lem-Commutator}.
\end{proof}

\section{New approach related to  Egorov theorem}\label{Flows and Egorov theorem type}

This section is devoted to the construction of  hyperbolic flows in infinite dimensional spaces associated to nonlocal pseudo-differential operators and explore the conjugation of pseudo-differential operators with this flow. More precisely, we are concerned with solving the evolution equation
\begin{equation} \label{eqn:omega0}
\left\{ \begin{array}{ll}
  \partial_t \Phi =\partial_\theta\big(\rho \,|\textnormal  D|^{\alpha-1}+|\textnormal  D|^{\alpha-1}\rho)\Phi(t)\triangleq \partial_\theta\mathcal{T}\Phi(t)\triangleq \mathbb{A}\Phi(t)  ,&\\ 
   \Phi(0)=\textnormal{Id},
  \end{array}\right.
\end{equation}
where $|\textnormal  D|^{\alpha-1}$ is the modified fractional Laplacian introduced in \eqref{fract1} and  $\rho:\mathcal{O}\times\mathbb{T}^{d+1}\to \RR$ is a smooth function that enjoys the following symmetry,
$$
\forall (\lambda,\varphi,\theta)\in\mathcal{O}\times\T^{d+1},\quad \rho(\lambda,-\varphi,-\theta)=-\rho(\lambda,\varphi,\theta).
$$
We shall discuss throughout Section \ref{Sect-Nonlocal pseudo-differential} some basic properties of the flow such as its continuity over Sobolev spaces or its Lipschitz dependence with respect to $\rho$.  We point out  that this flow will be used later  in Proposition \ref{prop-constant-coe} in order to reduce the nonlocal part of the linearized operator to a Fourier multiplier, and this requirement will definitively   fix the function $\rho.$ 
  Afterwards, we shall investigate in Section \ref{Egorov's theorem type}  the conjugation of  pseudo-differential operators by these hyperbolic flows. In particular, if  we take an operator $\mathscr{A}$ of order $m$ satisfying \eqref{Def-pseud-w} and  we conjugate it with the flow $\Phi(t)$, that is to consider $\mathscr{A}_t\triangleq \Phi^{-1}(t)\mathscr{A}\Phi(t)$, then from Egorov theorem  $\mathscr{A}_t$ remains a pseudo-differential operator with order $m$. The delicate point    is to estimate $\mathscr{A}_t$ with respect to the topology given by \eqref{Def-pseud-w} with suitable tame estimates. The main result in this direction  will be explained in Theorem \ref{Prop-EgorV}.  This theorem, which is interesting in itself,   is of great importance  and will be used   in Section \ref{section-KAM-Red} during the KAM scheme implemented to reduce the remainder to a diagonal part.  There,  we should check that the new remainders obtained through the different transformations  (transport and the nonlocal reductions) give rise to good operators of order zero and enjoying suitable tame estimates in the topology \eqref{Top-NormX}. Then the required estimates will stem from Theorem \ref{Prop-EgorV}. 
  
  \subsection{Nonlocal pseudo-differential hyperbolic equation}\label{Sect-Nonlocal pseudo-differential}
  In this section we shall be concerned with the continuity of the hyperbolic flows generated by the nonlocal pseudo-differential equation \eqref{eqn:omega0}. At this level, we take for granted the existence of such flows and restrict the discussion  to the a priori  estimates. Notice that the flow construction, which is  omitted here, can be done in a classical way by mollifying the equation and getting  the suitable a priori estimates combined with  compactness arguments,  see for instance \cite{BertiMontalto}. 
\begin{proposition}\label{flowmap00}
Let $d\in\N^\star, q\in\N$ and $ s_0>\frac{d+5}{2}$, then the followings assertions hold true.

\begin{enumerate}
\item  Assume that  $\,\,\|\rho\|_{H^{2 s_0}}\lesssim1,\,\,$ then   for any $|t|\leqslant 1$ and $s\geqslant 0$
$$
\| \Phi(t)h\|_{H^{s}}+  \| \Phi^{\star}(t)h\|_{H^{s}}\lesssim \| h\|_{H^{s}}+\|\rho\|_{H^{s+ s_0}}\| h\|_{H^{ s_0}}.
$$
Moreover the flow map is reversiblity preserving in the sense of Definition $\ref{Def-Rev}.$

\item If $\,\,\|\rho\|_{2s_0+q+1}^{q,\kappa}\lesssim 1,\,\,$ then for any $|t|\leqslant 1$ and $s\geqslant q$
$$
\| \Phi(t)h\|_{s}^{q,\kappa}+\| \Phi^{\star}(t)h\|_{s}^{q,\kappa}\lesssim \| h\|_{s}^{q,\kappa}+\|\rho\|_{s+s_0}^{q,\kappa}\| h\|_{s_0+q}^{q,\kappa}.
$$
Moreover, we have
$$
{\| (\Phi(t)-\textnormal{Id})h\|_{s}^{q,\kappa}+\| (\Phi^{\star}(t)-\textnormal{Id})h\|_{s}^{q,\kappa}\lesssim \| h\|_{s+1}^{q,\kappa}\|\rho\|_{s_0}^{q,\kappa}+\|\rho\|_{s+s_0}^{q,\kappa}\| h\|_{s_0+1+q}^{q,\kappa}
.}
$$
\item Assume that
$\displaystyle\max_{i=1,2}\|\rho_i\|_{2s_0+1+q}^{q,\kappa}\lesssim1,
$ then for any  $s\geqslant \max(s_0,q)$  and $|t|\leqslant 1$
\begin{align*}
\| \Delta_{12}\Phi(t)h\|_{{s}}^{q,\kappa}+\| \Delta_{12}\Phi^{\star}(t)h\|_{{s}}^{q,\kappa}&\lesssim \|\Delta_{12}\rho\|_{{s}+1}^{q,\kappa} \| h\|_{s+1+q}^{q,\kappa}\Big(1+\|\rho_1\|_{s+{s}_0+1}^{q,\kappa}+\|\rho_2\|_{s+{s}_0+1}^{q,\kappa}\Big).\end{align*}
\end{enumerate}
\end{proposition}
\begin{proof}
${\bf{(i)}}$
As we have mentioned before, we shall skip the details about the  existence of the flow map and just   focus on  the a priori estimates.
Let $h$ be a smooth function and denote by $h_t\triangleq\Phi(t)h$  the unique  solution of the pseudo-differential equation
\begin{equation} \label{eqn:omegaa}
\left\{ \begin{array}{ll}
  \partial_t h_t  =  \partial_\theta\mathcal{T}(\varphi,\theta)h_t  ,&\\ 
   h(0)=h.
  \end{array}\right.
\end{equation}
where $\mathcal{T}$ is the self-adjoint operator $\mathcal{T}\triangleq\rho \,|\textnormal  D|^{\alpha-1}+|\textnormal  D|^{\alpha-1}\rho.$
We shall first  start with the case $s=0$. Using  integration by parts and the fact that $\mathcal{T}$ is self-adjoint we get
\begin{eqnarray*}
\frac12 \frac{d}{dt}\|h_t\|_{L^2}^2&=&\int_{\T^{d+1}} h_t(\varphi,\theta)\partial_\theta\mathcal{T}h_t(\varphi,\theta) d\varphi d\theta\\
&&=-\int_{\T^{d+1}} h_t(\varphi,\theta)\mathcal{T}\partial_\theta h_t(\varphi,\theta) d\varphi d\theta\\
&&\quad \triangleq I_t.
\end{eqnarray*}
Thus by summation,
\begin{eqnarray*}
I_t&=&\frac12\int_{\T^{d+1}} h_t\big[\partial_\theta,\mathcal{T} \big]h_t d\varphi d\theta\\
&=&\frac12\int_{\T^{d+1}} h_t(\partial_\theta\rho)|\textnormal  D|^{\alpha-1}h_t d\varphi d\theta+\frac12\int_{\T^{d+1}} h_t|\textnormal  D|^{\alpha-1}\big((\partial_\theta\rho)h_t\big) d\varphi d\theta.
\end{eqnarray*}
Using Lemma \ref{Laplac-frac0} we get for $\alpha\in [0,1)$,
$$
\||\textnormal  D|^{\alpha-1}h_t \|_{L^2}\lesssim \|h_t\|_{L^2}\quad\hbox{and}\quad \||\textnormal  D|^{\alpha-1}(\partial_\theta\rho) h_t \|_{L^2}\lesssim \|\partial_\theta\rho\|_{L^\infty} \|h_t\|_{L^2}.
$$
Consequently 
\begin{eqnarray*}
|I_t|&\lesssim&\|\partial_\theta\rho\|_{L^\infty}\|h_t\|_{L^2}^2.
\end{eqnarray*}
Therefore we get after simplification 
\begin{align}\label{Gron01}
\frac{d}{dt}\|h_t\|_{L^2}\lesssim \|\partial_\theta\rho\|_{L^\infty}\|h_t\|_{L^2}.
\end{align}

Hence we infer from  Gronwal inequality 
\begin{align}\label{GronX1}
\sup_{t\in[-1,1]}\|h_t\|_{L^2}\leqslant \|h\|_{L^2}e^{C\|\partial_\theta\rho\|_{L^\infty}}.
\end{align}
To get the desired estimate it is enough to make appeal to Sobolev embeddings.
Let us now move to the estimate of $\|h_t\|_{H^{s}}$ for $s\geqslant0, $. Applying the Fourier multiplier $\Lambda^{s}$ to \eqref{eqn:omegaa} yields
\begin{equation*}
\left\{ \begin{array}{ll}
  \partial_t \Lambda^{s}h_t  = \partial_\theta\mathcal{T} \Lambda^{s}h_t+\big[\Lambda^{s},\partial_\theta\mathcal{T}\big]h_t  ,&\\ 
   h(0)=h\in H^{s}(\T^{d+1}).
  \end{array}\right.
\end{equation*}
Proceeding as before for $s=0$ we get similarly to \eqref{Gron01}
\begin{align}\label{Com-J0}
\frac{d}{dt}\| \Lambda^{s}h_t\|_{L^2}\lesssim \|\partial_\theta \rho\|_{L^\infty}\| \Lambda^{s}h_t\|_{L^2}+\big\|\big[\Lambda^{s},\partial_\theta\mathcal{T}\big]h_t\big\|_{L^2}.
\end{align}
It remains to estimate the commutator which can be written in the form
$$
\big[\Lambda^{s},\partial_\theta\mathcal{T}\big]=\partial_\theta\big[\Lambda^{s},\rho\big]|\textnormal{D}|^{\alpha-1}+ |\textnormal{D}|^{\alpha-1}\partial_\theta\big[\Lambda^{s},\rho\big].
$$
Combining   Lemma \ref{Lem-Commutator}-(v) with   Lemma \ref{Laplac-frac0} we infer for any $s_0>\frac{d+5}{2}$
\begin{align*}
\big\|\partial_\theta\big[\Lambda^{s},\rho\big]|\textnormal{D}|^{\alpha-1}h_t\big\|_{L^2}\lesssim&\||\textnormal{D}|^{\alpha-1}h\|_{H^{s}}\|\rho\|_{H^{s_0}}+ \||\textnormal{D}|^{\alpha-1}h\|_{H^{s_0}}\|\rho\|_{H^{s+s_0}}\\
\lesssim&\|h\|_{H^{s}}\|\rho\|_{H^{s_0}}+ \|h\|_{H^{s_0}}\|\rho\|_{H^{s+s_0}}
\end{align*}
and in a similar way 
\begin{align*}
\big\||\textnormal{D}|^{\alpha-1}\partial_\theta\big[\Lambda^{s},\rho\big]h_t\big\|_{L^2}\lesssim&\|h\|_{H^{s}}\|\rho\|_{H^{s_0}}+ \|h\|_{H^{s_0}}\|\rho\|_{H^{s+s_0}}.
\end{align*}
Putting together the preceding estimates we obtain
\begin{align}\label{Com-J1}
\big\|\big[\Lambda^{s},\partial_\theta\mathcal{T}\big]h_t\big\|_{L^2}\lesssim&\|h_t\|_{H^{s}}\|\rho\|_{H^{s_0}}+ \|h_t\|_{H^{s_0}}\|\rho\|_{H^{s+s_0}}.
\end{align}
Inserting \eqref{Com-J1} into \eqref{Com-J0} we obtain from Sobolev embeddings
\begin{align}\label{Com-J2}
\nonumber \frac{d}{dt}\| h_t\|_{H^{s}}\lesssim& \|\partial_\theta \rho\|_{L^\infty}\| h_t\|_{H^{s}}+\|h_t\|_{H^{s}}\|\rho\|_{H^{s_0}}+ \|h_t\|_{H^{s_0}}\|\rho\|_{H^{s+s_0}}\\
\lesssim&\|h_t\|_{H^{s}}\|\rho\|_{H^{s_0}}+ \|h_t\|_{H^{s_0}}\|\rho\|_{H^{s+s_0}}.
\end{align}
By taking $s=s_0, $ and applying Gronwall Lemma we find for $t\in\RR$
\begin{align}\label{Com-J3}
 \| h_t\|_{H^{s_0}}\leqslant & e^{C|t|\|\rho\|_{H^{2s_0}}}\| h\|_{H^{s_0}}.
\end{align}
Applying once again Gronwall lemma to  \eqref{Com-J2}  and using \eqref{Com-J3} yields 
\begin{align*}
\nonumber \| h_t\|_{H^{s}}\leqslant & e^{C|t|\|\rho\|_{H^{s_0}}}\left(\| h\|_{H^{s}}+\|\rho\|_{H^{s+s_0}}\left|\int_0^t\|h_{t^\prime}\|_{H^{s_0}}dt^\prime\right| \right)\\
\leqslant&e^{C|t|\|\rho\|_{H^{2s_0}}}\Big(\| h\|_{H^{s}}+|t|\|\rho\|_{H^{s+s_0}}\| h\|_{H^{s_0}} \Big).
\end{align*}
This can be written in the form 
\begin{align}\label{Com-J5}
 \| \Phi(t)h\|_{H^{s}}\leqslant &e^{C|t|\|\rho\|_{H^{2s_0}}}\Big(\| h\|_{H^{s}}+|t|\|\rho\|_{H^{s+s_0}}\| h\|_{H^{s_0}} \Big).
\end{align}
This gives the desired  statement when $\|\rho\|_{H^{2s_0}}\lesssim 1$ and $|t|\leqslant 1.$  Next we shall estimate the   adjoint operator $\Phi^{\star}(t)$. We can check that $\Phi^{\star}(t)$ satisfies the equation 
\begin{equation*} 
\left\{ \begin{array}{ll}
  \partial_t \Phi^{\star}(t) =\Phi^{\star}(t)\mathbb{A}^*\,&\\ 
   \Phi(0)=\textnormal{Id}.
  \end{array}\right.
\end{equation*}
Since $\Phi(t)$ commutes with its generator $\mathbb{A}$, that is, $\Phi(t)\mathbb{A}=\mathbb{A}\Phi(t)$ and then taking the adjoint yields  $\Phi^{\star}(t)\mathbb{A}^\star=\mathbb{A}^\star\Phi^{\star}(t)$. Therefore from the identity $\mathbb{A}^\star=-\mathcal{T}\partial_\theta$ we find
\begin{equation} \label{eqn-adjoint}
\left\{ \begin{array}{ll}
  \partial_t \Phi^{\star}(t) =-\big(\rho \,|\textnormal  D|^{\alpha-1}+|\textnormal  D|^{\alpha-1}\rho)\partial_\theta \Phi^{\star}(t) ,&\\ 
   \Phi(0)=\textnormal{Id}.
  \end{array}\right.
\end{equation}
It follows that the equation of $\Phi^{\star}(t)$ is similar to $\Phi(t)$ and its generator enjoys the same properties as the generator $\mathbb{A}$.  Therefore, following  the same steps one gets  at the end  the same estimate \eqref{Com-J5} for the adjoint operator $\Phi^{\star}(t)$.
The reversibility preserving can be  easily checked from the symmetry of $\rho$ and the uniqueness of the Cauchy problem.
\\
${\bf{(ii)}}$
First, observe that the estimate  \eqref{Com-J5} gives the desired result  with $q=0$.
Now, we need to establish  a similar estimate for $\partial_\lambda^\beta \Phi(t)h$, with $\lambda=(\omega,\alpha)$ and $|\beta|\leqslant q.$ We point out that  the delicate manipulation  concerns the estimates when we want to differentiate with respect to $\alpha$ because the order of the generator $\mathcal{T}$ depends on $\alpha$ and not on $\omega$.   Therefore,  we shall restrict the discussion to this case and estimate $\partial_\alpha^j \Phi(t)h$ for $1\leqslant j\leqslant q.$ To start, we differentiate \eqref{eqn:omegaa}  according to Leibniz formula 
\begin{equation} \label{eqnM0}
  \partial_t \partial_\alpha^jh_t  = \partial_\theta \mathcal{T}(\varphi,\theta)\partial_\alpha^jh_t+\sum_{k=0}^{j-1}\left(_k^j\right)\partial_\theta(\partial_\alpha^{j-k}\mathcal{T}) \partial_\alpha^k h_t .
\end{equation}
Then similarly to \eqref{Com-J0} and \eqref{Com-J1} we deduce for $s\geqslant q$ 
\begin{align}\label{Com-J00}
\nonumber \frac{d}{dt}\|  \partial_\alpha^jh_t\|_{{H^{s-j}}}&\lesssim \| \partial_\alpha^jh_t\|_{H^{s-j}}\|\rho\|_{H^{s_0}(\T^{d+1})}+ \| \partial_\alpha^jh_t\|_{H^{s_0}}\|\rho\|_{H^{s-j+s_0}}\\
&+
\sum_{k=0}^{j-1}\big\|\partial_\theta(\partial_\alpha^{j-k}\mathcal{T}) \partial_\alpha^k h_t\big\|_{H^{s-j}}.
\end{align}
Using once again  Leibniz formula we may write
\begin{align}\label{Com-JL00}
( \partial_\alpha^{j-k}\mathcal{T})=&
\sum_{m=0}^{j-k}\left(_m^{j-k}\right)\Big((\partial_\alpha^m\rho)(\partial_\alpha^{j-k-m}|\textnormal{D}|^{\alpha-1})+(\partial_\alpha^{j-k-m}|\textnormal{D}|^{\alpha-1})(\partial_\alpha^m\rho)\Big).
\end{align}
Using Lemma \ref{Law-prodX1} and Lemma \ref{Laplac-frac0} we get
\begin{align}\label{Com-JL000}
\nonumber\big\|\partial_\theta((\partial_\alpha^m\rho)(\partial_\alpha^{j-k-m}|\textnormal{D}|^{\alpha-1})g)\big\|_{H^{s-j}}&\lesssim \big\|(\partial_\alpha^m\rho)(\partial_\alpha^{j-k-m}|\textnormal{D}|^{\alpha-1})g\big\|_{H^{s-j+1}}\\
&\lesssim \big\|\partial_\alpha^m\rho\|_{H^{s_0}}\|g\big\|_{H^{s-j+1}}+\big\|\partial_\alpha^m\rho\|_{H^{s-j+1}}\|g\big\|_{H^{s_0}}.
 \end{align}
Similarly we obtain
\begin{align}\label{Com-JL01}
\nonumber \big\|\partial_\theta(\partial_\alpha^{j-k-m}|\textnormal{D}|^{\alpha-1})(\partial_\alpha^m\rho)g\big\|_{H^{s-j}}&\lesssim \big\|\partial_\alpha^m\rho\|_{H^{s_0}}\|g\big\|_{H^{s-j+1}}\\
&+\big\|\partial_\alpha^m\rho\|_{H^{s-j+1}}\|g\big\|_{H^{s_0}}.
\end{align}
Therefore putting together \eqref{Com-JL00},\eqref{Com-JL000} and  \eqref{Com-JL01} allows to get,
\begin{align}\label{Com-P01}
\nonumber \big\|\partial_\theta(\partial_\alpha^{j-k}\mathcal{T}) \partial_\alpha^k h_t\big\|_{H^{s-j}}&\lesssim \sum_{m=0}^{j-k} \big\|\partial_\alpha^m\rho\|_{H^{s_0}}\| \partial_\alpha^k h_t\big\|_{H^{s-j+1}}\\
&+\sum_{m=0}^{j-k} \big\|\partial_\alpha^m\rho\|_{H^{s-j+1}}\| \partial_\alpha^k h_t\big\|_{H^{s_0}}.
\end{align}
Plugging \eqref{Com-P01} into \eqref{Com-J00} yields
\begin{align}\label{Com-JX0}
 \nonumber \frac{d}{dt}\|  \partial_\alpha^jh_t\|_{{H^{s-j}}}&\lesssim \| \partial_\alpha^jh_t\|_{H^{s-j}}\|\rho\|_{H^{s_0}(\T^{d+1})}+ \| \partial_\alpha^jh_t\|_{H^{s_0}}\|\rho\|_{H^{s-j+s_0}}\\
 &+
\sum_{k=0}^{j-1}\sum_{m=0}^{j-k} \Big(\big\|\partial_\alpha^m\rho\|_{H^{s_0}}\| \partial_\alpha^k h_t\big\|_{H^{s-j+1}}+ \big\|\partial_\alpha^m\rho\|_{H^{s-j+1}}\| \partial_\alpha^k h_t\big\|_{H^{s_0}}\Big).
\end{align}
Using Gronwall inequality we get for $t\geqslant0$
\begin{align*}
 \nonumber\|  \partial_\alpha^jh_t\|_{{H^{s-j}}}&\lesssim e^{Ct\|\rho\|_{H^{s_0}}}\Bigg[\|  \partial_\alpha^jh\|_{{H^{s-j}}}+\|\rho\|_{H^{s-j+s_0}}\int_0^t e^{-C\tau\|\rho\|_{H^{s_0}}}\| \partial_\alpha^jh_\tau\|_{H^{s_0}}d\tau\\
\nonumber &+
\sum_{k=0}^{j-1}\sum_{m=0}^{j-k} \big\|\partial_\alpha^m\rho\|_{H^{s_0}}\int_0^t e^{-C\tau\|\rho\|_{H^{s_0}}}\| \partial_\alpha^k h_\tau\big\|_{H^{s-j+1}}d\tau\\
&+\sum_{k=0}^{j-1}\sum_{m=0}^{j-k} \big\|\partial_\alpha^m\rho\|_{H^{s-j+1}}\int_0^t e^{-C\tau\|\rho\|_{H^{s_0}}}\| \partial_\alpha^k h_\tau\big\|_{H^{s_0}}d\tau\Bigg].
\end{align*}
Multiplying by $\kappa^j$ with $\kappa\in(0,1)$   we infer
\begin{align*}
 \nonumber\kappa^j\|  \partial_\alpha^jh_t\|_{{H^{s-j}}}&\lesssim e^{Ct\|\rho\|_{H^{s_0}}}\Bigg[\kappa^j\|  \partial_\alpha^jh\|_{{H^{s-j}}}+\|\rho\|_{H^{s-j+s_0}}\kappa^j\int_0^te^{-C\tau\|\rho\|_{H^{s_0}}}\| \partial_\alpha^jh_\tau\|_{H^{s_0}} d\tau\\
\nonumber &+
\sum_{k=0}^{j-1}\sum_{m=0}^{j-k} \kappa^m\big\|\partial_\alpha^m\rho\|_{H^{s_0}}\kappa^k\int_0^t e^{-C\tau\|\rho\|_{H^{s_0}}} \| \partial_\alpha^k h_\tau\big\|_{H^{s-j+1}}d\tau\\
&+\sum_{k=0}^{j-1}\sum_{m=0}^{j-k} \kappa^m\big\|\partial_\alpha^m\rho\|_{H^{s-j+1}}\kappa^k\int_0^t e^{-C\tau\|\rho\|_{H^{s_0}}}\| \partial_\alpha^k h_\tau\big\|_{H^{s_0}}d\tau\Bigg].
\end{align*} 
Thus we find from Sobolev embeddings
\begin{align*}
 \nonumber\kappa^j\|  \partial_\alpha^jh_t\|_{{H^{s-j}}}&\lesssim e^{Ct\|\rho\|_{H^{s_0}}}\Bigg[\kappa^j\|  \partial_\alpha^jh\|_{{H^{s-j}}}+\|\rho\|_{H^{s+s_0}}\kappa^j\int_0^te^{-C\tau\|\rho\|_{H^{s_0}}}\| \partial_\alpha^jh_\tau\|_{H^{s_0}} d\tau\\
\nonumber &+
\sum_{k=0}^{j-1}\sum_{m=0}^{j-k} \kappa^m\big\|\partial_\alpha^m\rho\|_{H^{s_0+j-m}}\kappa^k\int_0^t e^{-C\tau\|\rho\|_{H^{s_0}}} \| \partial_\alpha^k h_\tau\big\|_{H^{s-k}}d\tau\\
&+\sum_{k=0}^{j-1}\sum_{m=0}^{j-k} \kappa^m\big\|\partial_\alpha^m\rho\|_{H^{s+1-m}}\kappa^k\int_0^t e^{-C\tau\|\rho\|_{H^{s_0}}}\| \partial_\alpha^k h_\tau\big\|_{H^{s_0+j-1-k}}d\tau\Bigg].
\end{align*} 
Consequently, according to the Definition \ref{Def-WS} we get for any $s\geqslant q$,
\begin{align}\label{Com-JX01}
 \nonumber\| h_t\|_{s}^{q,\kappa}&\lesssim e^{Ct\|\rho\|_{s_0}^{0,\kappa}}\Bigg[\| h\|_{s}^{q,\kappa}+\|\rho\|_{s+s_0}^{0,\kappa}\int_0^t e^{-C\tau\|\rho\|_{s_0}^{0,\kappa}}\| h_\tau\|_{s_0}^{q,\kappa}d\tau\\
\nonumber  &+
\|\rho\|_{s_0+q}^{q,\kappa}\int_0^t e^{-C\tau\|\rho\|_{s_0}^{0,\kappa}} \|  h_\tau\big\|_{s}^{q-1,\kappa}d\tau\\
&+\|\rho\|_{s+1}^{q,\kappa}\int_0^t e^{-C\tau\|\rho\|_{s_0}^{0,\kappa}} \|  h_\tau\big\|_{s_0+q-1}^{q-1,\kappa}d\tau\Bigg].
\end{align}
Assume that
\begin{align}\label{Hou-1}
 \|\rho\|_{2s_0+q+1}^{q,\kappa}\lesssim1,
\end{align}
then we get from the foregoing estimate combined with Sobolev embeddings
\begin{align*}
 \forall\,|t|\lesssim 1,\quad\| h_t\|_{s}^{q,\kappa}&\lesssim \| h\|_{s}^{q,\kappa}+\|\rho\|_{s+s_0}^{q,\kappa}\int_0^t\| h_\tau\|_{s_0+q}^{q,\kappa}d\tau+
\|\rho\|_{s_0+q}^{q,\kappa}\int_0^t \|  h_\tau\big\|_{s}^{q,\kappa}d\tau.
\end{align*}
It follows from Gronwall inequality and \eqref{Hou-1} that 
\begin{align}\label{It-P00}
 \forall\,|t|\lesssim 1,\,\forall s\geqslant q,\quad\| h_t\|_{s}^{q,\kappa}&\lesssim \| h\|_{s}^{q,\kappa}+\|\rho\|_{s+s_0}^{q,\kappa}\int_0^t\| h_\tau\|_{s_0+q}^{q,\kappa}d\tau.
 \end{align}
 By taking $s=s_0+q$ and using once again Gronwall inequality combined with  \eqref{Hou-1} we infer
 \begin{align*}
 \forall\,|t|\lesssim 1,\quad\| h_t\|_{s_0+q}^{q,\kappa}&\lesssim \| h\|_{s_0+q}^{q,\kappa}.
 \end{align*}
 Plugging this inequality into \eqref{It-P00} yields
 \begin{align*}
 \forall\,|t|\lesssim 1,\,\forall s\geqslant q,\quad\| h_t\|_{s}^{q,\kappa}&\lesssim \| h\|_{s}^{q,\kappa}+\|\rho\|_{s+s_0}^{q,\kappa}\| h\|_{s_0+q}^{q,\kappa}.
 \end{align*}
Using the flow mapping $\Phi$, the preceding estimate can be written in the form
 \begin{align}\label{ham-pod1}
 \forall\,|t|\lesssim 1,\,\forall s\geqslant q,\quad\|\Phi(t) h\|_{s}^{q,\kappa}&\lesssim \| h\|_{s}^{q,\kappa}+\|\rho\|_{s+s_0}^{q,\kappa}\| h\|_{s_0+q}^{q,\kappa}.
 \end{align} 
 The estimate of $\Phi^{\star}$ which satisfies the equation \eqref{eqn-adjoint} can be done in a similar way  to $ \Phi(t)$. The estimates of $g(t)\triangleq (\Phi(t)-\hbox{Id})h$ can be done using the following steps. First we write the equation
\begin{equation*} 
\left\{ \begin{array}{ll}
  \partial_t g  =  \partial_\theta\mathcal{T}g+\partial_\theta\mathcal{T}h,\\ 
   g(0)=0.
  \end{array}\right.
\end{equation*}
Then using Duhamel formula we get
$$
g(t)=\int_0^t\Phi(t-\tau)\partial_\theta\mathcal{T}hd\tau.
$$
Applying Proposition \ref{flowmap00}-$($ii$)$ yields   under the assumption \eqref{Hou-1}  
that
\begin{align*}
\forall\, |t|\leqslant 1,\quad \| g(t)\|_{s}^{q,\kappa}&\lesssim \| \partial_\theta\mathcal{T}h\|_{s}^{q,\kappa} + \|\rho\|_{s+s_0}^{q,\kappa} \|\partial_\theta\mathcal{T}h\|_{s_0+q}^{q,\kappa}\\
&\lesssim \|\mathcal{T}h\|_{s+1}^{q,\kappa} + \|\rho\|_{s+s_0+1}^{q,\kappa} \|\mathcal{T}h\|_{s_0+q+1}^{q,\kappa}
\end{align*}
Then applying the first result in Lemma \ref{Law-prodX1}-(ii) and Lemma \ref{Laplac-frac0},  one may get under the smallness condition \eqref{Hou-1} \begin{align*}
\|g(t)\|_{s}^{q,\kappa}&\lesssim \|\rho\|_{s_0}^{q,\kappa} \| h\|_{s+1}^{q,\kappa}+ \|\rho\|_{s+s_0}^{q,\kappa}\| h\|_{s_0+q+1}^{q,\kappa}.
\end{align*}
 The estimate $ (\Phi^\star(t)-\hbox{Id})h$ can be done in a similar way and this ends the proof of the second point.
 
 \smallskip

{\bf{(iii)}} First notice that the operator $\Phi_{12}\triangleq \Delta_{12}\Phi$ satisfies the equation
\begin{equation} \label{Taou-taou1}
\left\{ \begin{array}{ll}
  \partial_t \Phi_{12} =\mathbb{A}_1\Phi_{12}(t)+\mathbb{A}_{12}\Phi_{2}(t)  ,&\\ 
   \Phi(0)=0,
  \end{array}\right.
\end{equation}
Then using Duhamel formula we get
$$
\Phi_{12}(t)h=\int_0^t\Phi_1(t-\tau)\mathbb{A}_{12}\Phi_{2}(\tau)hd\tau.
$$
Consequently, Proposition \ref{flowmap00}-$($ii$)$ yields   under \eqref{Hou-1} 
\begin{align*}
\forall\, |t|\leqslant 1,\quad \| \Phi_{12}h\|_{s}^{q,\kappa}&\lesssim  \int_0^t\| \mathbb{A}_{12}\Phi_{2}(\tau)h\|_{s}^{q,\kappa}d\tau + \|\rho_1\|_{s+s_0}^{q,\kappa} \int_0^t\| \mathbb{A}_{12}\Phi_{2}(\tau)h\|_{s_0+q}^{q,\kappa} d\tau.
\end{align*}
On the other hand using the law products of Lemma \ref{Law-prodX1} one gets
\begin{align*}
\| \mathbb{A}_{12}\Phi_{2}(\tau)h\|_{s}^{q,\kappa}&\lesssim \|\Delta_{12}\rho\|_{s+1}^{q,\kappa} \|\Phi_{2}(\tau)h\|_{s_0}^{q,\kappa}+ \|\Delta_{12}\rho\|_{s_0}^{q,\kappa} \|\Phi_{2}(\tau)h\|_{s+1}^{q,\kappa}.
\end{align*}
Using once again  Proposition \ref{flowmap00}-$($ii$)$  combined with \eqref{Laplac-frac0} we find
 for any $\tau\in[-1,1]$
\begin{align*}
\| \mathbb{A}_{12}\Phi_{2}(\tau)h\|_{s}^{q,\kappa}&\lesssim \|\Delta_{12}\rho\|_{s+1}^{q,\kappa}\| h\|_{s_0+q}^{q,\kappa}
+ \|\Delta_{12}\rho\|_{s_0}^{q,\kappa} \Big[\| h\|_{s+1}^{q,\kappa}+\|\rho_2\|_{s+1+s_0}^{q,\kappa}\| h\|_{s_0+q}^{q,\kappa}\Big].
\end{align*}
Combining the preceding few estimates and using \eqref{Hou-1}  we find
\begin{align*}
\forall\, |t|\leqslant 1,\quad \| \Phi_{12}(t)h\|_{{s}_0}^{q,\kappa}&\lesssim \|\Delta_{12}\rho\|_{{s}_0+q+1}^{q,\kappa} \| h\|_{{s}_0+q+1}^{q,\kappa}.
\end{align*}
In addition, we get for any $s\geqslant s_0$
\begin{align*}
\| \Phi_{12}(t)h\|_{{s}}^{q,\kappa}&\lesssim \|\Delta_{12}\rho\|_{{s}+1}^{q,\kappa} \| h\|_{s+1+q}^{q,\kappa}\Big(1+\|\rho_1\|_{s+{s}_0+1}^{q,\kappa}+\|\rho_2\|_{s+{s}_0+1}^{q,\kappa}\Big).
\end{align*}
The same estimate  holds  for $\Phi_{12}^*(t)$ by using similar arguments combined with \eqref{eqn-adjoint}.
The proof of the Proposition \ref{flowmap00} is now achieved.
\end{proof}

 \subsection{T\"oplitz  matrices and binomial convolution}\label{Sect-Top-M}
 We shall discuss some basic properties on T\"oplitz matrix operators that will  play a crucial role later when we will investigate the proof of an Egorov type theorem stated in Theorem \ref{Prop-EgorV}.  \\
Given  $n\in\NN$ and $ k\in\{0,1,..,n+1\}$, we denote  by $\DDD_k(n)$  the set of  lower triangluar matrix operators $\mathcal{M}$ such that
\begin{equation}\label{M-spe}\mathcal{M}=\begin{pmatrix}
    m_{0,0} &0&..&..&..&0  \\
m_{1,1}&m_{0,1}&0&..&0&0 \\
   m_{2,2}& m_{1,2}&m_{0,2}&.. &..&0\\
    ..&..&..&..&0&0\\
    m_{n-1,n-1}& m_{n-2,n-1}&..&..& m_{0,n-1}&0\\
   m_{n,n}& m_{n-1,n}&..&..& m_{1,n}&m_{0,n}
  \end{pmatrix},
\end{equation}
supplemented with the following constraints
$$
m_{i,j}=\left(_i^j\right)\mathcal{\mu}_\mathcal{M}(i), \forall\, i\leqslant j\in\{0,..,n\}\quad\hbox{and}\quad m_{i,j}=0,\,\forall\, 0\leqslant  i\leqslant k-1, 
$$
where $\left(_i^j\right)$ are the usual binomial coefficients and $\mathcal{\mu}_\mathcal{M}(i)$ is a scalar operator that depends only on  the index $i$ and not on $k$. The index $k$ gives the number of vanishing sub-diagonals ( including the main diagonal) in the square matrix operator. We call $\mu_\mathcal{M}$ the law of $\mathcal{M}.$  Typically,  the main examples that we shall encounter during the proof of Theorem \ref{Prop-EgorV}   are pseudo-differential operators of type $\mathcal{\mu}_\mathcal{M}(i)=\rho^{(i)} |\textnormal{D}|^{\alpha}$, where the index $i$ measures the spatial regularity of some suitable functions. Notice that by construction one has the inclusions
$$
\DDD_0(n)\supset \DDD_1(n)\supset..\supset\DDD_{n+1}(n)=\{0\}.
$$
For further purposes, we need to define  for $q\in\NN$ the class  $\DDD_k(n,q)$ of  bloc matrix operators taking  the form \eqref{M-spe} and such that each entry $\mathcal{\mu}_\mathcal{M}(i)$ is a matrix operator in $ \DDD_{0}(q).$ It is obvious from the definitions that   $\DDD_k(n,0)=\DDD_k(n)$. We point out that the introduction of such kind of operators through the index $q$ appears natural during  the proof of Theorem  \ref{Prop-EgorV}  when we will track  the  regularity estimates with respect to the external parameters. As to  the index $n$ it is associated to the time/spatial regularity.\\
At this stage, we do not care about the  topological structure of  the operators lying in the class $\DDD_k(n,q)$,  we are simply concerned with some of their algebraic structures.  One of the main results  discussed below shows that  commutators between matrix operators in the class $\DDD_k(n,q)$ still in the same class and whose  entries  are commutators too. This property fails in general even with standard real matrices, unless we impose additional structures as for example  with T\"oplitz matrices. Notice that getting commutators in the entries is crucial to reduce the order of the operators during  the proof of Theorem  \ref{Prop-EgorV}.
\begin{lemma}\label{Lem-top}
The following assertions hold true.
\begin{enumerate}
\item Let $n,q\in\NN$ and $0\leqslant k\leqslant n+1$, then the set $\DDD_k(n,q)$ is stable under multiplication, that is, for any $\mathcal{M},\mathcal{N}\in\DDD_k(n,q)$ we have $\mathcal{M}\mathcal{N}\in \DDD_k(n,q)$ and its law $\mathcal{\mu}_\mathcal{M N}$ is given by the binomial convolution,
\begin{align*}
\mathcal{\mu}_\mathcal{M N}(i)&=\sum_{\ell=0}^i\left(_\ell^i\right)\mathcal{\mu}_{\mathcal{M}}(i-\ell)\mathcal{\mu}_{\mathcal{N}}(\ell)\\
&\triangleq \left(\mathcal{\mu}_{\mathcal{M}}\varoast \mathcal{\mu}_{\mathcal{N}}\right)(i).
\end{align*}
In addition,  if $\mathcal{M}\in \DDD_k(n,q),\mathcal{N}\in\DDD_m(n,q)$, then $\mathcal{M}\mathcal{N}\in \DDD_{\min(k+m,n+1)}(n,q).$ 
\item For any $\mathcal{M},\mathcal{N}\in\DDD_k(n,q)$ the law of the commutator  $[\mathcal{M},\mathcal{N}]$ is given by
 \begin{align*}
\mathcal{\mu}_{[\mathcal{M},\mathcal{ N}]}(i)&=\sum_{\ell=0}^i\left(_\ell^i\right)\big[\mathcal{\mu}_{\mathcal{M}}(i-\ell),\mathcal{\mu}_{\mathcal{N}}(\ell)\big]\\
&\triangleq (\mathcal{\mu}_{\mathcal{M}}\boxast \mathcal{\mu}_{\mathcal{N}})(i).
\end{align*}
Furthermore, the law of the iterated commutator $\textnormal{Ad}^n_{\mathcal{M}} \mathcal{N}$ is given by
 \begin{align*}
\mathcal{\mu}_{\textnormal{Ad}^n_{\mathcal{M}} \mathcal{N}}&=\underbrace{\mathcal{\mu}_{\mathcal{M}}\boxast \mathcal{\mu}_{\mathcal{M}}\boxast...\boxast\mathcal{\mu}_{\mathcal{M}}\boxast \mathcal{\mu}_{\mathcal{N}}}_{n\,\textnormal{ times}}.
\end{align*}
\end{enumerate}
\end{lemma}
\begin{remark}\label{Rmq-1}
As a byproduct of this lemma, if $\mathscr{M},\mathscr{N}\in \DDD_0(n)$ defined in \eqref{M-spe} and  their entries  $m_{i,j}, n_{k,\ell}$ commute together, that is, $[m_{i,j}, n_{k,\ell}]=0$, then $\big[\mathscr{M},\mathscr{N} \big]=0.$ In addition, the stability point  stated in $(i)$ asserts that the product of two elements  contains more vanishing sub-diagonals. In particular, if $\mathscr{M},\mathscr{N}\in \DDD_1(n)$  and $1\leqslant j\leqslant n+1$ then 
$$
\mathscr{M}\mathscr{N}\in \DDD_2(n)\quad\hbox{and}\quad \mathscr{M}^j\in \DDD_j(n).
$$
From the last point we infer that the matrix operator $\mathscr{M}$ is nilpotent and $\mathscr{M}^{n+1}=0.$
\end{remark}
\begin{proof}
$($i$)$
Let  $\mathcal{M}=(\overline{m}_{i,j})_{0\leqslant i,j\leqslant n}\in \DDD_k(n,q)$ and $\mathcal{N}=(\overline{n}_{i,j})_{0\leqslant i,j\leqslant n}\in \DDD_k(n,q)$ and set $\mathcal{P}=\mathcal{M}\mathcal{N}$, then with these notations 
$$
\forall\, 0\leqslant j\leqslant i,\quad \overline{m}_{i,j}={m}_{i-j,i},\quad \overline{n}_{i,j}={n}_{i-j,i}.
$$ 
Thus we may check that
$$\mathcal{P}=\begin{pmatrix}
    p_{0,0} &0&..&..&..&0  \\
p_{1,1}&p_{0,1}&0&..&0&0 \\
   p_{2,2}& p_{1,2}&p_{0,2}&.. &..&0\\
    ..&..&..&..&0&0\\
    p_{n-1,n-1}& p_{n-2,n-1}&..&..& p_{0,n-1}&0\\
   p_{n,n}& p_{n-1,n}&..&..& p_{1,n}&p_{0,n}
  \end{pmatrix}=(\overline{p}_{i,j})_{0\leqslant i,j\leqslant n}
$$
with $\overline{p}_{i,j}=0$ for $i<j$, and for $0\leqslant j\leqslant i$ we have
\begin{align*}
\overline{p}_{i,j}&=\sum_{\ell=j}^{i}\overline{m}_{i,\ell}\overline{n}_{\ell,j}\\
&=\sum_{\ell=j}^{i}{m}_{i-\ell,i}{n}_{\ell-j,\ell}.
\end{align*}
Consequently we find through elementary combinatorics formula 
\begin{align*}
\overline{p}_{i,j}&=\sum_{\ell=j}^{i} \left(_{i-\ell}^i\right) \left(_{\ell-j}^\ell\right)\mathcal{\mu}_\mathcal{M}(i-\ell)\mathcal{\mu}_\mathcal{N}(\ell-j)\\
&=\left(_{j}^i\right) \left(\mathcal{\mu}_{\mathcal{M}}\varoast \mathcal{\mu}_{\mathcal{N}}\right)(i-j)
\end{align*}
and therefore for $0\leqslant i \leqslant j$ 
\begin{align*}
{p}_{i,j}&=\overline{p}_{j,j-i}\\
&=\left(_{i}^j\right) \left(\mathcal{\mu}_{\mathcal{M}}\varoast \mathcal{\mu}_{\mathcal{N}}\right)(i).
\end{align*}
Thus the law of  $\mathcal{M}\mathcal{N}$  is  given by the binomial convolution. 
Now since $\mathcal{M},\mathcal{N}\in\DDD_k(n,q)$ then by construction  $\mathcal{\mu}_\mathcal{M}(i),\mathcal{\mu}_\mathcal{N}(i)\in\DDD_0(q)$ and 
$$
\mathcal{\mu}_\mathcal{M}(i)=\mathcal{\mu}_\mathcal{N}(i)=0,\quad \forall\, 0\leqslant i\leqslant k-1.
$$ 
This implies that $\mu_{\mathcal{M}\mathcal{N}}(i)\in\DDD_0(q)$ with  $\mu_{\mathcal{M}\mathcal{N}}(i)=0,$ for any  $ 0\leqslant i\leqslant k-1$. Thus the bloc matrix operator $\mathcal{M}\mathcal{N}\in\DDD_k(n,q). $ Similarly, if $\mathcal{M}\in\DDD_k(n,q)$ and $ \mathcal{N}\in\DDD_m(n,q)$ then from straightforward matrix computations we find
$$
\mathcal{\mu}_{\mathcal{M N}}(i)=0,\quad \forall\, 0\leqslant i\leqslant \min(k+m-1,n).
$$
Hence $\mathcal{M N}\in\DDD_{\min(k+m,n+1)}(n,q).$

$($ii$)$ Let  $\mathcal{M},\mathcal{N}\in\DDD_k(n,q)$ then from the preceding point $($i$)$
 the commutator $[\mathcal{M},\mathcal{N}]$ remains in $\DDD_k(n,q)$ and its law is described by
\begin{eqnarray*}
\mathcal{\mu}_{[\mathcal{M},\mathcal{N}]}(i)&=&\sum_{\ell=0}^{i}\left(_\ell^i\right)[\mathcal{\mu}_{\mathcal{M}}(i-\ell),\mathcal{\mu}_{\mathcal{N}}(\ell)]\\
&\triangleq &\mathcal{\mu}_{\mathcal{M}}\boxast \mathcal{\mu}_{\mathcal{N}}(i)
\end{eqnarray*}
where $[\mathcal{\mu}_{\mathcal{M}}(i-\ell),\mathcal{\mu}_{\mathcal{N}}(\ell)]$ is the  commutator between the  matrix  operators $\mathcal{\mu}_{\mathcal{M}}(i-\ell)$ and $ \mathcal{\mu}_{\mathcal{N}}(\ell)$ that belong to $\DDD_0(q)$ and thus 
$$
[\mathcal{\mu}_{\mathcal{M}}(i-\ell),\mathcal{\mu}_{\mathcal{N}}(\ell)]\in \DDD_0(q).
$$
Notice that the entries of the  matrix operator $[\mathcal{\mu}_{\mathcal{M}}(i-\ell),\mathcal{\mu}_{\mathcal{N}}(\ell)]$ are linear combination of scalar commutators. This ends the proof of the desired result.
\end{proof}
 \subsection{Refined Egorov's theorem}\label{Egorov's theorem type}
 The main concern  of this section is to deal with the conjugation of some class of pseudo-differential operators with the hyperbolic flow given in \eqref{eqn:omega0} and  constructed in Proposition \ref{flowmap00}. We shall establish in Theorem \ref{Prop-EgorV} below a refined  version of the classical Egorov theorem  in a suitable  topology of pseudo-differential operators of order $m\in[-1,0].$  Notice that the tame estimates, which are crucial in the KAM reduction and along Nash-Moser scheme,  will be delicately  established. This part is  the hardest  technical point raised in this paper and  will be performed through  a new approach based on  both kernel and symbol  representations of operators  combined with T\"oplitz structure on some matrices involved in the dynamics of  the kernel derivatives.    \\
 Let us start with fixing the problem that we wish to solve and directly connected with  KAM reduction of the remainder  that will be carefully investigated  in Section \ref{section-KAM-Red}.
Consider a toroidal pseudo-differential  operator  $\mathcal{R}_0: \mathscr{C}^\infty(\T^{d+1})\mapsto \mathscr{C}(\T^{d+1})$ with the following kernel representation 
\begin{align}\label{Kernel-op}
\mathcal{R}_0h(\varphi,\theta)=\int_\T \mathcal{K}_0(\varphi,\theta,\eta) h(\varphi,\eta) d\eta.
\end{align}
This structure is very  special and more  specifically, the variable $\varphi$ acts as a local parameter in this description. It allows to cover   a large class of operators that  turns out to be sufficient for the applications in this paper.  The symbol $a_0$ associated to $\mathcal{R}_0$ is related to the kernel $\mathcal{K}_0$ according to the formula \eqref{symb-kern} 
$$
\sigma_{\mathcal{R}_0}(\varphi,\theta,\xi)=\int_{\T}\mathcal{K}_0(\varphi,\theta,\theta+\eta) e^{\ii  \eta\xi} d\eta.
$$
Let  $\Phi(t)$ be the flow defined as the solution  to the equation \eqref{eqn:omega0} and whose generator is given by the pseudo-differential operator 
\begin{align}\label{AA-S1}
\mathbb{A}&=\partial_\theta\big(\rho|\textnormal{D}|^{\alpha-1}+|\textnormal{D}|^{\alpha-1}\rho\big)\\
\nonumber&\triangleq\partial_\theta \mathcal{T}_\theta.
\end{align}
The  function $\rho$ is smooth enough and may depend on $(\varphi,\theta)\in\T^{d+1}$ or on the variables $(\lambda,\varphi,\theta)\in\mathcal{O}\times\T^{d+1}$. We remind that  the modified fractional Laplacian  $|\textnormal{D}|^{\alpha-1}$ was   introduced \mbox{in \eqref{fract1}} and from which we infer
$$
\mathbb{A}h(\varphi,\theta)=\partial_\theta\int_\T K(\varphi,\theta,\eta)h(\varphi,\eta)d\eta,\quad K(\varphi,\theta,\eta)=\frac{1}{2\pi}\frac{\rho(\varphi,\theta)+\rho(\varphi,\eta)}{|\sin(\frac{\theta-\eta}{2})|^{\alpha}}\cdot
$$
Let us precise an elementary observation that will be frequently used later.  As one can see  the kernel $K$ is singular at the diagonal $\theta=\eta$ but  the modified one given by  $(\varphi,\theta,\eta)\mapsto K(\varphi,\theta,\theta+\eta)$ is smooth on $\theta$ and  singular only at $\eta=0$ which is the integration variable. By this way, the singularity in $\eta$ which  is not so violent can be merely  treated by integration. 
Here we shall fix the rectangular open box
$$
\mathcal{O}=\mathscr{U}\times (0,\overline\alpha)\quad\hbox{with}\quad \overline\alpha\in\left(0,\frac12\right)
$$
and $\mathscr{U}$ is an open bounded  set of $\RR^d$. The main task  here  is to analyze the conjugation of $\mathcal{R}_0$ by the flow $\Phi(t)$, that is to consider $\Phi(t)\mathcal{R}_0\Phi(-t),$ and show that this operator remains  in the same class of pseudo-differential operators  with   the suitable tame estimates. We refer to  \mbox{Section \ref{Se-Toroidal pseudo-differential operators}} for the definition of  the class of symbols that will used below.  Our main statement reads as follows.  
\begin{theorem}\label{Prop-EgorV}
Let $m\in[-1,0], q\in\N$ and $ \mathcal{R}_0$ be the  operator in \eqref{Kernel-op}. Then the following assertions hold true.
\begin{enumerate}
\item  {Let $s\geqslant0,\, s_0>\frac{d+1}{2}$ and assume  the existence of an increasing  $\log$-convex function $F_0$ such that  
$$\forall\, s^\prime\geqslant 0,\,\forall \,\gamma\in\N\quad\hbox{with}\quad   s^\prime+\gamma\leqslant s_0+27+s\Longrightarrow \interleave \mathcal{R}_0\interleave_{m,s^\prime,\gamma}\leqslant   F_0(s^\prime+\gamma).
$$ 
}
Then for any ${t\in[0,1]}$ 
$$
\interleave \Phi(t)\mathcal{R}_0\Phi(-t)\interleave_{m,s,0}\leqslant {C e^{C\mu^{s+1}(0)}}\mu(s),
$$
with
$$
\mu(s)=F_0(s+27+s_0)+\|\rho\|_{H^{s+27+s_0}}
$$
and $C$ being a constant depending  only on $s$. We also get
$$
\interleave \Phi(t)\mathcal{R}_0\Phi(-t)\interleave_{m,s,0}\leqslant C e^{C\mu^{s+1}s(0)}F_0(s+27+s_0)\big(1+\mu(s)\big).
$$
\item  Let $s\geqslant0,\, s_0>\frac{d+1}{2}+q$ and assume  the existence of an increasing  $\log$-convex function $\overline{F}_0$  such that   {
$$\forall\, s'\geqslant 0,\,\forall\,\gamma\in\N\quad\hbox{with}\quad   s'+\gamma\leqslant s_0+27+s\Longrightarrow  \interleave \mathcal{R}_0\interleave_{m,s',\gamma}^{q,\kappa}\leqslant \overline{F}_0(s'+\gamma).
$$ 
}
Then for any  ${t\in[0,1]}$
$$
\interleave \Phi(t)\mathcal{R}_0\Phi(-t)\interleave_{m,s,0}^{q,\kappa}\leqslant C e^{C\overline\mu^{s+1}(0)}\overline\mu(s)
$$
with
$$
\overline\mu(s)=\overline{F}_0(s+27+s_0)+\|\rho\|_{{s+27+s_0}}^{q,\kappa}.
$$
We also get
$$
\interleave \Phi(t)\mathcal{R}_0\Phi(-t)\interleave_{m,s,0}^{q,\kappa}\leqslant C e^{C\overline\mu^{s+1}(0)}\overline{F}_0(s+27+s_0)\big(1+\overline\mu(s)\big).
$$
\item Denote $\Delta_{12}\Phi(t)\mathcal{R}_0\Phi(-t)\triangleq \Phi_1(t)\mathcal{R}_0\Phi_1(-t)-\Phi_2(t)\mathcal{R}_0\Phi_2(-t)$ where $\Phi_i$ is the flow associated to $\rho_i$.  Let $s_0>\frac{d+1}{2}+q$ and assume  the existence of an increasing  $\log$-convex function $\overline{F}_0$ such for such that   {
$$\forall\, s'\geqslant 0,\,\forall\,\gamma\in\N\quad\hbox{with}\quad   s'+\gamma\leqslant s_0+10\Longrightarrow  \interleave \mathcal{R}_0\interleave_{m,s',\gamma}^{q,\kappa}\leqslant \overline{F}_0(s'+\gamma).
$$ 
}
Then for any $t\in[0,1]$
$$
\interleave \Delta_{12}\Phi(t)\mathcal{R}_0\Phi(-t)\interleave_{0,0,0}^{q,\kappa}\leqslant C e^{C\big(\|\rho_2\|_{H^{s_0+10}}+\|\rho_1\|_{H^{s_0+10}}+\overline{F}_0(s_0+10)\big)}\|\rho_1-\rho_2\|_{{s_0+5}}^{q,\kappa}.
$$
\end{enumerate}
\end{theorem}
Before moving to the proof, some remarks are in order.
\begin{remark}
\begin{enumerate}
\item 
The statement of the theorem is restricted for $m\in[-1,0]$ which quite  sufficient in the application during the remainder reduction that will be explored in \mbox{Section \ref{section-KAM-Red}.} The extension to the values of  $m$ in the complement set is likely possible  provided we impose the suitable adaptations. 
\item It is worthy to point out that the constraints on the initial symbol covers a range of values on the parameters $s,\gamma$ but the conclusion deals only with the estimate with the norm $\|\cdot\|_{m,s,0}$ where $\gamma=0$. This fact is due to a diffusion mechanism in the regularity persistence of the symbol estimates: this stems in part  from the commutator estimates  detailed in \mbox{Lemma \ref{comm-pseudo}} where we win in the order of the operator but we lose in space-phase variables.  
\item In the statement $(iii)$ the estimates on the difference are measured in the topology $\interleave \cdot\interleave_{0,0,0}^{q,\kappa}$ and not with the expected one  $\interleave \cdot\interleave_{m,s,0}^{q,\kappa}$.  This is due to delicate technical points that could be solved with the same approach  but with more refined analysis. Notice that, later during the application in the study of  the spectrum dependence of the linearized operator, this  weak estimate will provide us by  interpolation arguments  with only a  H\"older continuity  dependence with respect to the torus $i$ instead of Lipschitz dependence.  This slight loss of regularity has no important  effects in the  arguments used to establish  the stability   of Cantor sets during the last phase in Nash-Moser scheme related the measure of the final Cantor set  associated to the emergence of invariant nonlinear torus. 
\end{enumerate}
\end{remark}
\begin{proof}
Let us briefly explore the main ideas of the proof.  We shall write down  the kernel equation of the the operator $\mathcal{R}(t)\triangleq  \Phi(t)\mathcal{R}_0\Phi(-t)$. Then the symbol estimates are transformed into looking for  suitable kernel estimates in weighted Sobolev spaces. This will be implemented by successive differentiations of  the time evolution of the kernel equation combined with energy estimates. This procedure leads to a collection of coupled  transport equations with non dissipative positive order operators acting as source terms that cannot be directly treated by energy estimates. The key observation is to understand  the  global structure of the full set of equations that can be recast into a vectorial transport equations with an additional  positive order matrix operator. The particularity of this matrix is to be nilpotent and of T\"opliz structure as in the definition \eqref{M-spe}. Then we perform energy estimates through a new mixed norm and the use of refined commutator estimates leading to several technical difficulties that will be made clear throughout the proof.

 {\bf{(i)}}  It seems to be convenient to   start with proving the result for $m=0$ and move later to the case $m\in[-1,0)$ by making the suitable adaptations.\\
  $\blacktriangleright$ {\it{ Case}} $m=0.$
Denote by $\mathcal{R}(t)\triangleq \Phi(t)\mathcal{R}_0\Phi(-t)$ then we can check directly from the flow \mbox{equation \eqref{eqn:omega0}}  that  this operator  satisfies the Heisenberg equation
$$
\partial_t \mathcal{R}(t)=[\mathbb{A},\mathcal{R}(t)].
$$
Now, set  $\mathcal{R}_1(t)\triangleq \mathcal{R}(t)- \mathcal{R}_0$, then we deduce from the preceding equation that
\begin{equation}\label{Heisenberg-11}
\left\{ \begin{array}{ll}
  \partial_t \mathcal{R}_1(t)=[\mathbb{A},\mathcal{R}_1(t)]+[\mathbb{A},\mathcal{R}_0]&\\
  \mathcal{R}_1(0)=0.
  \end{array}\right.
\end{equation}
To show the desired estimate for $\mathcal{R}(t)$ it is enough  to prove, according to Sobolev embeddings,  that the operator $\mathcal{R}_1(t)$  satisfies for any $t\in[0,1]$ and $s\in\NN$
\begin{align}\label{T-j-L}
\interleave \mathcal{R}_1(t)\interleave_{0,s,0}\lesssim C{\big(1+F_0^2(s_0+12)}\big) {e^{C\|\rho\|_{H^{s_0+12}}^{s}}}\Big(F_0(s+12)+\|\rho\|_{H^{s+12}}\Big).
\end{align}
Since the variable $\varphi$ acts as a parameter on the generator  $\mathbb{A}$, then it will act in a similar way on  the associated flow $\Phi(t).$  This implies that one may view the operator $\mathcal{R}_1(t)$ as a collection of operators $\mathcal{R}_1(t,\varphi)$ parametrized by $\varphi$ and acting pointwisely on test functions  as $\mathcal{R}_1(t,\varphi)h(\varphi,\cdot)$. As a consequence  one gets    the following kernel representation of $\mathcal{R}_1(t)$,
$$
\mathcal{R}_1(t)h(\varphi,\theta)=\int_\T \mathcal{K}_t(\varphi,\theta,\eta) h(\varphi,\eta) d\eta
$$
and therefore the   dynamics is encoded  in the new kernel $\mathcal{K}_t$ that we need to estimate carefully. 
Coming back to  the norm definition \eqref{Def-Norm-M1}, we may write
\begin{align*}
\interleave\mathcal{R}_1(t)\interleave_{0,s,0}&=\sup_{\xi\in\Z}\|\sigma_t(\cdot,\centerdot,\xi)\|_{H^s(\T^{d+1})}
\end{align*}
with $\sigma_t$ being the symbol of $\mathcal{R}_1(t)$ that can be recovered from the kernel according to the \mbox{identity \eqref{symb-kern}}
\begin{align}\label{Symb-Kern}
\nonumber \sigma_t(\varphi,\theta,\xi)&=e^{-\ii  \xi\theta}\mathcal{R}_1(t)e^{\ii \xi\theta}\\
&=\int_\T\mathcal{K}_t(\varphi,\theta,\theta+\eta)e^{\ii \eta\xi}d\eta.
\end{align}
Then applying  Cauchy-Schwarz inequality and using a change of variables give with $s=n\in\N$
\begin{align*}
\interleave\mathcal{R}_1(t)\interleave_{0,s,0}&\leqslant \int_\T\|\mathcal{K}_t(\cdot,\centerdot,\centerdot+\eta)\|_{H^{s}_{\varphi,\theta}}d\eta\\
&\leqslant \int_\T\|\mathcal{K}_t(\cdot,\centerdot,\centerdot+\eta)\|_{H^{s}_{\theta}L^2_\varphi}d\eta+\int_\T\|\mathcal{K}_t(\cdot,\centerdot,\centerdot+\eta)\|_{L^2_\theta H^{s}_{\varphi}}d\eta\\
&\lesssim\sum_{i=0}^n\big(\|\partial_\chi^i\mathcal{K}_t\|_{L^2(\T^{d+2})}+\|\partial_\varphi^i\mathcal{K}_t\|_{L^2(\T^{d+2})}\big)
\end{align*}
with 
$$
\partial_\chi=\partial_\theta+\partial_\eta.
$$
The estimates of the two terms of the right-hand side are quite similar and we shall start the first one dealing the mixed derivatives.

\smallskip

$\bullet$ {\it  Estimate of $\|\partial_\chi^i\mathcal{K}_t\|_{L^2(\T^{d+2})}$}. The main goal  is to prove the following estimate: $\forall\, i\in\{0,1,..,n\},$
\begin{align}\label{Est-kernel-i}
 \forall t\in[0,1],\quad \| \partial_\chi^i\mathcal{K}_t\|_{L^2(\T^{d+2})}
   &\leqslant C \big(1+G^{i+1}(0)\big)\, e^{C\|\rho\|_{H^{s_0+3}}^n} G(1+i)
\end{align}
with
\begin{align*}
G(i)\triangleq 
 \big(\|\rho\|_{H^{s_0}}+F_0(s_0)\big)\big(\|\rho\|_{H^{s_0+11+i}}+F_0(s_0+11+i)\big)+\|\rho\|_{H^{s_0+7+i}}.
\end{align*}
Recall that $F_0$ is a log-convex function that controls some suitable symbol norm as  assumed in Theorem \ref{Prop-EgorV}-(i).
To proceed, we shall start with writing the evolution equation governing the kernel $\mathcal{K}_t$.  From the integral representation one gets by definition of the commutator
\begin{align}\label{Iden-H-C}
\nonumber[\mathbb{A},\mathcal{R}_1(t)]h(\varphi,\theta)&=\bigintsss_{\T}\big(\mathbb{A}_\theta\mathcal{K}_t(\varphi,\theta,\eta)\big) h(\varphi,\eta)d\eta
-\int_{\T}\mathcal{K}_t(\varphi,\theta,\eta) \big(\mathbb{A}_\eta h\big)(\varphi,\eta) d\eta\\
&=\int_{\T}\big(\mathbb{A}_\theta\mathcal{K}_t-\mathbb{A}_\eta^\star\mathcal{K}_t\big)(\varphi,\theta,\eta) h(\varphi,\eta)d\eta.
\end{align} 
The notation $\mathbb{A}_\theta$ stands for  the action of the operator $\mathbb{A}$ in the variable $\theta$, that is,
$$
\mathbb{A}_\theta h(\varphi,\theta)=\partial_\theta\big(\rho(\varphi,\theta)|\textnormal{D}_\theta|^{\alpha-1}h+|\textnormal{D}_\theta|^{\alpha-1}\rho(\varphi,\cdot)h\big)(\varphi,\theta)
$$
 and similarly we define $\mathbb{A}_\eta$. We also denote by   $\mathbb{A}_\eta^\star$ the partial adjoint of $\mathbb{A}_\eta$ with respect to Hilbert structure of $L^2_\eta(\T).$  
 Then combining  Heisenberg equation \eqref{Heisenberg-11} with \eqref{Iden-H-C} we get the evolution nonlocal PDE for the kernel
\begin{equation*}
\left\{ \begin{array}{ll}
  \partial_t\mathcal{K}_t(\varphi,\theta,\eta)=\mathbb{A}_\theta \mathcal{K}_t(\varphi,\theta,\eta)-\mathbb{A}_\eta^\star\mathcal{K}_t(\varphi,\theta,\eta)+K_0(\varphi,\theta,\eta)
&\\
  \mathcal{K}_0=0,
  \end{array}\right.
\end{equation*}
where $K_0$ is the kernel associated to $ [\mathbb{A},\mathcal{R}_0]$ in the sense
$$
[\mathbb{A},\mathcal{R}_0]h(\varphi,\theta)=\int_{\T}K_0(\varphi,\theta,\eta)h(\varphi,\eta) d\eta,
$$
and one can check the following relation 
$$
K_0(\varphi,\theta,\eta)=\mathbb{A}_\theta \mathcal{K}_0(\varphi,\theta,\eta)-\mathbb{A}_\eta^\star\mathcal{K}_0(\varphi,\theta,\eta).
$$
By virtue of \eqref{AA-S1}, it is clear that  the operator  $\mathcal{T}_\eta$ is self-adjoint and then we can easily check that 
\begin{equation}\label{adj-st}\mathbb{A}_\eta^\star=-\mathbb{A}_\eta+\big[\partial_\eta,\mathcal{T}_\eta\big].
\end{equation} 
Consequently,  the kernel equation becomes 
\begin{equation}\label{kernel-dyna}
\partial_t\mathcal{K}_t(\varphi,\theta,\eta)-\big(\mathbb{A}_\theta+\mathbb{A}_\eta\big) \mathcal{K}_t(\varphi,\theta,\eta)+\big[\partial_\eta,\mathcal{T}_\eta \big]\mathcal{K}_t(\varphi,\theta,\eta)=K_0(\varphi,\theta,\eta).
\end{equation}
 Direct computations yield  the identity
\begin{equation}\label{diffM1}
\mathcal{S}_0\triangleq \big[\partial_\eta,\mathcal{T}_\eta \big]=(\partial_\eta\rho)|\textnormal{D}_\eta|^{\alpha-1}+|\textnormal{D}_\eta|^{\alpha-1}(\partial_\eta\rho).
\end{equation}
Now, we shall prove the estimate \eqref{Est-kernel-i} for $i=0$. For this purpose, we take  the $L^2(\T^{d+2})$ inner product of the equation \eqref{kernel-dyna} with $\mathcal{K}_t$ 
\begin{align*}
\frac{d}{dt}\int_{\T^{d+2}}|\mathcal{K}_t(\varphi,\theta,\eta)|^2d\varphi d\theta d\eta&=\int_{\T^{d+2}}\mathcal{K}_t(\varphi,\theta,\eta)\left(\big[\mathcal{T}_\eta,\partial_\eta \big]-\big[\mathcal{T}_\theta,\partial_\theta \big]\right)\mathcal{K}_t(\varphi,\theta,\eta)d\varphi d\theta d\eta\\
&+2\int_{\T^{d+2}}\mathcal{K}_t(\varphi,\theta,\eta)K_0(\varphi,\theta,\eta)d\varphi d\theta d\eta,
\end{align*}
where we have used by virtue of \eqref{adj-st} the identity
\begin{align*}
\langle \mathbb{A}_\theta \mathcal{K}_t,\mathcal{K}_t\rangle_{L^2(\T^{d+2})}&=\tfrac12 \langle (\mathbb{A}_\theta+\mathbb{A}_\theta^\star) \mathcal{K}_t,\mathcal{K}_t\rangle_{L^2(\T^{d+2})}\\
&=\tfrac12 \langle [\partial_\theta,\mathcal{T}_\theta ] \mathcal{K}_t,\mathcal{K}_t\rangle_{L^2(\T^{d+2})}.
\end{align*}
It follows that
\begin{align*}
\frac{d}{dt}\int_{\T^{d+2}}|\mathcal{K}_t(\varphi,\theta,\eta)|^2d\varphi d\theta d\eta
&\lesssim \|\partial_\theta\rho\|_{L^\infty} \|\mathcal{K}_t\|_{L^2(\T^{d+2})}^2+ \|\mathcal{K}_t\|_{L^2(\T^{d+2})} \|{K}_0\|_{L^2(\T^{d+2})}.
\end{align*}
Applying  Gronwall inequality yields
\begin{align}\label{Hyp-est0}
\|\mathcal{K}_t\|_{L^2(\T^{d+2})}\lesssim \|{K}_0\|_{L^2(\T^{d+2})}t\,e^{Ct \|\rho\|_{W^{1,\infty}}}.
\end{align}
Notice that his energy estimate can be extended to cover a more  general case with time dependent source terms. Indeed, let  $f$ and $g$ satisfy
\begin{align}\label{Tr01}
\nonumber \mathscr{L}(f)&\triangleq \partial_tf+\widehat{\mathscr{L}}(f)\\
&\triangleq \partial_tf-\widehat{\mathbb{A}} f+\mathcal{S}_0 f=g,\quad  \widehat{\mathbb{A}}\triangleq \mathbb{A}_\theta+\mathbb{A}_\eta,
\end{align}
then for any $t\geqslant0$
\begin{align}\label{Hyp-est}
\|f(t)\|_{L^2(\T^{d+2})}\leqslant e^{Ct\|\rho\|_{W^{1,\infty}}}\Big(\|f(0)\|_{L^2(\T^{d+2})} +\int_0^{t}e^{-C\tau\|\rho\|_{W^{1,\infty}}}\|g(\tau)\|_{L^2(\T^{d+2})}d\tau\Big).
\end{align}
The next step aims at extending   the estimate \eqref{Hyp-est0} to  higher order regularity. First, we remark that the equation \eqref{kernel-dyna} can be recast in the form
\begin{align}\label{kernel-dyna1}
\mathscr{L}(\mathcal{K}_t)=K_0.
\end{align}
Then applying $\partial_\chi^i$ to \eqref{kernel-dyna1} yields 
\begin{align}\label{Ker-N-F1}
\mathscr{L}\partial_{\chi}^i\mathcal{K}_t=\big[\partial_{\chi}^i,\widehat{\mathbb{A}}\big]\mathcal{K}_t-\big[\partial_{\chi}^i,\mathcal{S}_0\big] \mathcal{K}_t+\partial_{\chi}^i{K}_0.
\end{align}
By virtue of the abstract Leibniz rule we obtain
\begin{align}\label{Leib-R-P}
\big[\partial_{\chi}^i,\widehat{\mathbb{A}}\big]=\sum_{k=0}^{i-1}\left(_k^i\right)\textnormal{Ad}_{\partial_{\chi}}^{i-k}(\widehat{\mathbb{A}})\, \partial_{\chi}^{k},
\end{align}
where the adjoint mapping $\textnormal{Ad}_A$ and its iterated ones $\textnormal{Ad}^i_A$ are defined by
$$
\textnormal{Ad}_AB=[A,B]\quad\hbox{and}\quad \textnormal{Ad}_A^i=\underbrace{\textnormal{Ad}_A\circ..\circ\textnormal{Ad}_A}_{i\,\, times}.
$$
Using the special structure of the operator $\widehat{\mathbb{A}}$ detailed thtough \eqref{AA-S1} and \eqref{Tr01} we  can check by induction the following  identity, using in particular the general  relation $\textnormal{Ad}_A^i=[A,\textnormal{Ad}_A^{i-1}],$
$$
\textnormal{Ad}_{\partial_{\chi}}^{k}(\widehat{\mathbb{A}})=\widehat{\mathbb{A}}_{(k)},$$
with
\begin{align}\label{Ima1}
\widehat{\mathbb{A}}_{(k)}=\mathbb{A}_{\theta,(k)}+\mathbb{A}_{\eta,(k)},\quad \mathbb{A}_{\theta,(k)}\triangleq\partial_\theta\big((\partial_\theta^k\rho)|\textnormal{D}_\theta|^{\alpha-1}+|\textnormal{D}_\theta|^{\alpha-1}(\partial_\theta^k\rho)\big).
\end{align}
Notice that we have also used some   commutation relations of the type $[\partial_\theta,{\mathbb{A}}_{\eta,(k)}]=[\partial_\eta,{\mathbb{A}}_{\theta,(k)}]=0.$
Therefore we get by virtue of \eqref{Leib-R-P}
\begin{align}\label{HamA1}
\big[\partial_{\chi}^i,\widehat{\mathbb{A}}\big]=\sum_{k=0}^{i-1}\left(_k^i\right)\widehat{\mathbb{A}}_{(i-k)}\, \partial_{\chi}^{k}.
\end{align}
Similarly  we get from \eqref{diffM1}
\begin{align}\label{HamA2}
\nonumber\big[\partial_{\chi}^i,\mathcal{S}_0\big]=&\sum_{k=0}^{i-1}\left(_k^i\right)\textnormal{Ad}_{\partial_{\chi}}^{i-k}(\mathcal{S}_0)\, \partial_{\chi}^{k}\\
=&\sum_{k=0}^{i-1}\left(_k^i\right)\mathbb{S}_{(i-k)}\, \partial_{\chi}^{k}
\end{align}
with
\begin{align}\label{Ima2}
\mathbb{S}_{(k)}\triangleq(\partial_\eta^{(k+1)}\rho)|\textnormal{D}_\eta|^{\alpha-1}+|\textnormal{D}_\eta|^{\alpha-1}(\partial_\eta^{(k+1)}\rho).
\end{align}
Plugging the identities \eqref{HamA1} and \eqref{HamA2} into \eqref{Ker-N-F1} gives for any $i\in\{0,1,..,n\}$
\begin{align}\label{Ham1}
\nonumber \mathscr{L}\partial_{\chi}^{i}\mathcal{K}_t&=\sum_{k=0}^{{i}-1}\left(_{i-k}^{i}\right)\big(\widehat{\mathbb{A}}_{({i}-k)}-\mathbb{S}_{(i-k)}\big)\, \partial_{\chi}^{k} \mathcal{K}_t+\partial_{\chi}^{j} {K}_0\\
&=\sum_{k=0}^{{i}-1}m_{i-k,{i}}\, \partial_{\chi}^{k} \mathcal{K}_t+\partial_{\chi}^{i} {K}_0
\end{align}
with
\begin{align}\label{Bnk}
m_{k,i}&\triangleq\left(_k^i\right)\big(\widehat{\mathbb{A}}_{(k)}-\mathbb{S}_{(k)}\big).
\end{align}
We observe that for $i=0$ the first term in the right hand side disappears. Now, we  intend to write the system of equations \eqref{Ham1} in the matrix form. For this goal, let us introduce the lower triangular  matrix operator
\begin{equation}\label{Matrix-y}
 \mathcal{M}_n=\begin{pmatrix}
    0 &0&..&..&..&0  \\
 m_{1,1}&0&..&..&0&0 \\
   m_{2,2}& m_{1,2}&0&.. &..&0\\
    ..&..&..&..&0&0\\
 ..& ..&..&m_{1,n-1}& 0&0\\
 m_{n,n}& m_{n-1,n}&..&..& m_{1,n}&0
  \end{pmatrix}
\end{equation}
and the vectors
\begin{equation}\label{Matrix-yz}
\mathcal{X}_n(t)=\begin{pmatrix}
     \mathcal{K}_t &  \\
    \partial_\chi \mathcal{K}_t& \\
    ..&\\
    ..&\\
    \partial_\chi^n  \mathcal{K}_t
  \end{pmatrix},\quad\mathcal{Y}_n=\begin{pmatrix}
{K}_0 &  \\
    \partial_\chi{K}_0& \\
    ..&\\
    ..&\\
    \partial_\chi^n  {K}_0
  \end{pmatrix}.
  \end{equation}
Then \eqref{Ham1} can be recast in the matrix form as follows
\begin{align}\label{Hamb1}
\mathcal{L}_n\mathcal{X}_n(t)&=\mathcal{M}_n\mathcal{X}_n(t)+\mathcal{Y}_n,\quad \mathcal{L}_n\triangleq \mathscr{L}\textnormal{I}_{n+1}.
\end{align}
It is worthy  to mention that  the time-independent matrix operator  $\mathcal{M}_n$  belongs to the class $\DDD_1$  introduced in  belongs is  \eqref{M-spe}. This algebraic structure will be of  important  use along the proof. In particular, one deduces  that this matrix is nilpotent with  
\begin{align}\label{Nilpot}
\mathcal{M}_n^{n+1}=0.
\end{align}
This property is very important because the entries of $\mathcal{M}_n$ are operators with strictly positive orders that cannot be directly treated by energy estimates but the nilpotent structure gives somehow a cancellation and one may expect to close energy estimates by making appeal to a suitable norms combination. In order to take advantage of the nilpotent structure, we shall apply successively the matrix 
$\mathcal{M}_n $ to the equation \eqref{Hamb1} and write down the suitable evolution equations for  the iterated vectors $\mathcal{M}_n^k\mathcal{X}_n(t)$. To perform this, we first   apply the  abstract Leibniz rule used in \eqref{Leib-R-P} leading to 
\begin{align*}
\big[\mathcal{M}_n^k,\mathcal{L}_n\big]&=\sum_{j=0}^{k-1}\left(_j^k\right)\textnormal{Ad}_{\mathcal{M}_n}^{k-j}(\mathcal{L}_n)\, \mathcal{M}_n^{j}\\
=&\sum_{j=0}^{k-1}\left(_j^k\right)\textnormal{Ad}_{\mathcal{M}_n}^{k-j}(\widehat{\mathcal{L}}_n)\, \mathcal{M}_n^{j},\quad  \widehat{\mathcal{L}}_n\triangleq \widehat{\mathscr{L}}\,\,\textnormal{I}_{n+1}.
\end{align*}
We have used in the last line  the identity  $\textnormal{Ad}_{\mathcal{M}_n}^{k-j}(\widehat{\mathcal{L}}_n)=\textnormal{Ad}_{\mathcal{M}_n}^{k-j}(\mathcal{L}_n)$ which is a consequence of the definition  \eqref{Tr01} and the fact that the matrix $\mathcal{M}_n$ is time-independent. 
 Applying $\mathcal{M}_n^k$ to the equation   \eqref{Hamb1} and using the preceding identity  allow to get   for any $k\in\{1,.,n\}$
\begin{align}\label{Hamb2}
\nonumber \mathcal{L}_n\mathcal{M}_n^k\mathcal{X}_n(t)&=-\big[\mathcal{M}_n^k,\mathcal{L}_n\big]\mathcal{X}_n(t)+\mathcal{M}_n^{k+1}\mathcal{X}_n(t)+\mathcal{M}_n^{k}\mathcal{Y}_n\\
&=-\sum_{j=0}^{k-1}\left(_j^k\right)\textnormal{Ad}_{\mathcal{M}_n}^{k-j}(\widehat{\mathcal{L}}_n)\, \mathcal{M}_n^{j}\mathcal{X}_n(t)+\mathcal{M}_n^{k+1}\mathcal{X}_n(t)+\mathcal{M}_n^{k}\mathcal{Y}_n.
\end{align}
Define  $\mathcal{Z}_k(t)=\mathcal{M}_n^k\mathcal{X}_n(t)$, and denote by   $Y(t,i)$   the i-th component of the vector  $Y(t)$.
Then applying \eqref{Hyp-est}  with \eqref{Hamb2} we deduce  for any $k\in\{1,..,n\}$
\begin{align}\label{Hamb3}
 \nonumber\|\mathcal{Z}_{k}(t,i)\|_{L^2(\T^{d+2})}&\lesssim e^{Ct\|\rho\|_{W^{1,\infty}}}\Bigg(\int_0^t e^{-C\tau\|\rho\|_{W^{1,\infty}}}\sum_{j=0}^{k-1}\|\big(\textnormal{Ad}_{\mathcal{M}_n}^{k-j}(\widehat{\mathcal{L}}_n)\, \mathcal{M}_n^{j}\mathcal{X}_n\big)(\tau,i)\|_{L^2(\T^{d+2})}d\tau\\
 &+\int_0^t e^{-C\tau\|\rho\|_{W^{1,\infty}}}\|\mathcal{Z}_{k+1}(t,i)\|_{L^2(\T^{d+2})}d\tau+t \|\mathcal{M}_n^{k}\mathcal{Y}_n(i)\|_{L^2(\T^{d+2})} \Bigg)
\end{align}
and for $k=0$ we infer
\begin{align}\label{Es-z0-I}
\|\mathcal{Z}_{0}(t,i)\|_{L^2(\T^{d+2})}\lesssim e^{Ct\|\rho\|_{W^{1,\infty}}}\left( \int_0^t e^{-C\tau\|\rho\|_{W^{1,\infty}}}\|\mathcal{Z}_{1}(t,i)\|_{L^2(\T^{d+2})}d\tau+ t \|\mathcal{Y}_n(i)\|_{L^2(\T^{d+2})}\right).
\end{align}
At this level we shall make appeal to  Lemma \ref{lem-com-it2} allowing  to get  for any $i\in\{2,..,n+1\}$
$$
\left\|\left(\textnormal{Ad}_{\mathcal{M}_n}^{k-j}(\mathcal{L}_n)\, \mathcal{M}_n^{j}\mathcal{X}_n\right)(\tau,i)\right\|_{L^2(\T^{d+2})}\lesssim \|\rho\|_{H^{s_0+3}(\T^{d+1})}^{k-j}\sum_{\ell=1}^{i-1}\|\rho\|_{H^{s_0+3+i-\ell}(\T^{d+1})}\,\,\left\|\mathcal{Z}_j(\tau,\ell)\right\|_{L^2(\T^{d+2})}.
$$
Notice that for $i=1$ the first component $\left(\textnormal{Ad}_{\mathcal{M}_n}^{k-j}(\mathcal{L}_n)\, \mathcal{M}_n^{j}\mathcal{X}_n\right)(\tau,1)$ is vanishing because 
$$\textnormal{Ad}_{\mathcal{M}_n}^{k-j}(\mathcal{L}_n)\, \mathcal{M}_n^{j}\in \DDD_1(n)
$$
which follows from Lemma \ref{Lem-top}.
Inserting the foregoing  estimate into \eqref{Hamb3} yields 
\begin{align}\label{Hamb30}
 \nonumber&\|\mathcal{Z}_{k}(t,i)\|_{L^2(\T^{d+2})} \lesssim  e^{Ct\|\rho\|_{W^{1,\infty}}}\Bigg(\int_0^t e^{-C\tau\|\rho\|_{W^{1,\infty}}}\|\mathcal{Z}_{k+1}(t,i)\|_{L^2(\T^{d+2})}d\tau+t \|\mathcal{M}_n^{k}\mathcal{Y}_n(i)\|_{L^2(\T^{d+2})}\\
 &\qquad +\bigintsss_0^t e^{-C\tau\|\rho\|_{W^{1,\infty}}}\sum_{j=0}^{k-1}\|\rho\|_{H^{s_0+3}(\T^{d+1})}^{k-j}\sum_{\ell=1}^{i-1}\|\rho\|_{H^{s_0+3+i-\ell}(\T^{d+1})}\,\,\left\|\mathcal{Z}_j(\tau,\ell)\right\|_{L^2(\T^{d+2})}d\tau \Bigg).
\end{align}
Let us introduce the quantity 
$$
\zeta(t,i)\triangleq e^{-Ct \|\rho\|_{W^{1,\infty}}}\sum_{k=0}^{n} \|\mathcal{Z}_{k}(t,i)\|_{L^2}.
$$
Then summing  up \eqref{Es-z0-I} and \eqref{Hamb30}  over $k=1,..,n$ and using the relation  $ \mathcal{M}_n^{n+1}=0$, giving $\mathcal{Z}_{n+1}(t)=0$, we find 
\begin{align*}
 \nonumber\zeta(t,i) \lesssim  &\bigintsss_0^t e^{-C\tau\|\rho\|_{W^{1,\infty}}} \sum_{k=1}^{n}\sum_{j=0}^{k-1}\|\rho\|_{H^{s_0+3}}^{k-j}\sum_{\ell=1}^{i-1}\|\rho\|_{H^{s_0+3+i-\ell}}\,\,\left\|\mathcal{Z}_j(\tau,\ell)\right\|_{L^2(\T^{d+2})}d\tau\\
 &+\int_0^t \zeta(\tau,i)\|_{L^2}d\tau+t\sum_{k=0}^n \|\mathcal{M}_n^{k}\mathcal{Y}_n(i)\|_{L^2}.
\end{align*}
Consequently, we get from Fubini's principle that for $i\in\{1,..,n+1\}$
\begin{align*}
 \nonumber\zeta(t,i)
  \lesssim  &\Big(\|\rho\|_{H^{s_0+3}}+\|\rho\|_{H^{s_0+3}}^n\Big) \sum_{\ell=1}^{i-1}\|\rho\|_{H^{s_0+3+i-\ell}}\,\,\int_0^t\zeta(\tau,\ell)d\tau\\
 &+\int_0^t \zeta(\tau,i)\|_{L^2}d\tau+t\sum_{k=0}^n \|\mathcal{M}_n^{k}\mathcal{Y}_n(i)\|_{L^2}.
\end{align*}
Notice that for $i=1$ the first sum part in the right-hand side disappears.
Then Gronwall inequality gives 
\begin{align}\label{Es-z14}
\zeta(t,i)
  \leqslant &C_0 \sum_{\ell=1}^{i-1}\|\rho\|_{H^{s_0+3+i-\ell}}\,\,\int_0^t e^{C(t-\tau)}\zeta(\tau,\ell)d\tau+ C  e^{Ct}\sum_{k=0}^n \|\mathcal{M}_n^{k}\mathcal{Y}_n(i)\|_{L^2},
\end{align}
with
$$
C_0\triangleq C(1+\|\rho\|_{H^{s_0+3}}^{n}).
$$
For $i=1$ we get
\begin{align}\label{YY11}
\zeta(t,1)
 \leqslant  C  e^{Ct}\sum_{k=0}^n \|\mathcal{M}_n^{k}\mathcal{Y}_n(1)\|_{L^2}.
\end{align}
From Lemma \ref{V-13} we have 
\begin{align}\label{GKL01}
\|\mathcal{M}_n^k\mathcal{Y}_n(i)\|_{L^2(\T^{d+2})}\leqslant C\big(1+G^k(0)\big)G(i),
\end{align}
with
\begin{align}\label{GKL1}
G(i)\triangleq F_1\big(s_0+5+i\big)+\|\rho\|_{H^{s_0+6+|m|+i}},
\end{align}
where $F_1$ should be  an increasing  $\log-$convex function such that 
\begin{align}\label{WH1}
\forall s\geqslant s_0,\gamma\in\N\quad\hbox{with}\quad s+\gamma\leqslant s_0+n+1\Longrightarrow\interleave [\mathbb{A},\mathcal{R}_0]\interleave_{m,s,\gamma}\le F_1(s+\gamma),
\end{align}
for some fixed $m<-\frac12$. The next  goal is devoted to an explicit construction of $F_1$. According to Lemma \ref{comm-pseudo}-$($ii$)$ applied with $q=0$  one gets
\begin{align}\label{es-zk6}
\nonumber \interleave[\mathbb{A},\mathcal{R}_0]\interleave_{\alpha-1,s,\gamma}\lesssim&
\sum_{0\leqslant\beta\leqslant\gamma}\interleave\mathbb{A}\interleave_{\alpha,s,1+\beta}\interleave\mathcal{R}_0\interleave_{0,s_0+2,\gamma-\beta}+\interleave\mathbb{A}\interleave_{\alpha,s_0,1+\beta}\interleave\mathcal{R}_0\interleave_{0,s+2,\gamma-\beta}\\
\nonumber&+\sum_{ 0\le\beta\le\gamma}\interleave \mathbb{A}\interleave_{\alpha,s,0}\interleave \mathcal{R}_0\interleave_{0,s_0+2+\beta,\gamma-\beta} +\interleave \mathbb{A}\interleave_{0,s_0,0}\interleave\mathcal{R}_0\interleave_{0,s+2+\beta,\gamma-\beta} \\
\nonumber&\quad+
\sum_{0\leqslant\beta\leqslant\gamma}\interleave \mathcal{R}_0\interleave_{0,s,1+\beta}\interleave\mathbb{A}\interleave_{\alpha,s_0+2,\gamma-\beta}+\interleave\mathcal{R}_0\interleave_{0,s_0,1+\beta}\interleave\mathcal{A}\interleave_{\alpha,s+2,\gamma-\beta}\\
&\qquad+\sum_{ 0\leqslant\beta\leqslant\gamma}\interleave \mathcal{R}_0\interleave_{0,s,0}\interleave \mathbb{A}\interleave_{\alpha,s_0+2+\beta,\gamma-\beta} +\interleave \mathcal{R}_0\interleave_{0,s_0,0}\interleave\mathbb{A}\interleave_{\alpha,s+2+\beta,\gamma-\beta}.
\end{align}
Using \eqref{Hm-id1} one gets for   $\alpha\in(0,1)$
\begin{align}\label{WH-P1}
\interleave \mathbb{A}\interleave_{\alpha,s,\gamma}\lesssim \|\rho\|_{H^{s+\gamma+1}}.
\end{align}
Since $\mathcal{R}_0$ satisfies
\begin{align*}
\forall s\geqslant s_0,\gamma\in\N\quad\hbox{with}\quad s+\gamma\leqslant s_0+n+3\Longrightarrow\interleave \mathcal{R}_0\interleave_{m,s,\gamma}\le F_0(s+\gamma),
\end{align*}
then  plugging this into \eqref{es-zk6} allows to get for any $ s\geqslant s_0,\gamma\in\N$ with $ s+\gamma\leqslant s_0+n+1$
\begin{align*}
\nonumber\interleave[\mathbb{A},\mathcal{R}_0]\interleave_{\alpha-1,s,\gamma}\lesssim& \sum_{ 0\leqslant\beta\leqslant\gamma}
\|\rho\|_{H^{s+\beta+2}}F_0(s_0+2+\gamma-\beta)+\|\rho\|_{H^{s_0+\beta+2}}F_0(s+2+\gamma-\beta)\\
\nonumber&+\sum_{ 0\leqslant\beta\leqslant\gamma}\|\rho\|_{H^{s+1}}F_0(s_0+2+\gamma)+\|\rho\|_{H^{s_0+1}}F_0(s+2+\gamma) \\
\nonumber&\quad+
\sum_{0\leqslant\beta\leqslant\gamma}\|\rho\|_{H^{s_0+3+\gamma-\beta}}F_0(s+1+\beta)+\|\rho\|_{H^{s+3+\gamma-\beta}}F_0(s_0+1+\beta)\\
&\qquad+\sum_{ 0\leqslant\beta\leqslant\gamma}\|\rho\|_{H^{s_0+3+\gamma}}F_0(s)+\|\rho\|_{H^{s+3+\gamma}}F_0(s_0).
\end{align*}
Applying  Sobolev embeddings and the monotonicity of $F_0$ we obtain
\begin{align}\label{Rec-XK1}
\interleave[\mathbb{A},\mathcal{R}_0]\interleave_{\alpha-1,s,\gamma}&\lesssim \sum_{ 0\leqslant\beta\leqslant\gamma}
\|\rho\|_{H^{s+\beta+3}}F_0(s_0+2+\gamma-\beta)+\|\rho\|_{H^{s_0+\beta+3}}F_0(s+2+\gamma-\beta).
\end{align}
By virtue of  the $\log-$convexity of Sobolev norms and the function $F_0$ we infer 
\begin{align}\label{day-98}
\nonumber\|\rho\|_{H^{s+\beta+3}}F_0(s_0+2+\gamma-\beta)\leqslant& \|\rho\|_{H^{s_0}}^{\delta}\|\rho\|_{H^{s+\gamma+5}}^{1-\delta}F_0^{1-\delta}(s_0)F_0^{\delta}(s+5+\gamma)\\
\nonumber&\leqslant  
{F_0(s_0)\|\rho\|_{H^{s+\gamma+5}}+F_0(s+5+\gamma)\|\rho\|_{H^{s_0}}}
\\
&\quad \leqslant \Big(\|\rho\|_{H^{s_0}}+F_0(s_0)\Big)\Big(\|\rho\|_{H^{s+\gamma+5}}+F_0(s+5 +\gamma)\Big),
\end{align}
with $\delta\in(0,1)$.  Therefore we obtain
\begin{align}\label{D-K-G1}
\forall s\geqslant s_0,\gamma\in\N\quad\hbox{with}\quad s+\gamma\leqslant s_0+n+1\Longrightarrow\interleave[\mathbb{A},\mathcal{R}_0]\interleave_{\alpha-1,s,\gamma}&\leqslant F_1(s+\gamma),
\end{align}
with
\begin{align}\label{D-K-12}
F_1(s)\triangleq C
 \big(\|\rho\|_{H^{s_0}}+F_0(s_0)\big)\big(\|\rho\|_{H^{s+5}}+F_0(s+5)\big).
\end{align} 
Coming back to \eqref{WH1}, one may fix from  \eqref{D-K-12} the function  $F_1$ with the constraint $m=\alpha-1<-\frac12$ which is satisfied since $\alpha\in[0,\frac12).$ Notice that $F_1$ is $\log-$convex since the sum of $\log-$convex functions is $\log-$convex too. Inserting \eqref{D-K-12}  into \eqref{GKL1} allows to make the choice
\begin{align}\label{GKL2}
G(i)=C
 \big(\|\rho\|_{H^{s_0}}+F_0(s_0)\big)\big(\|\rho\|_{H^{s_0+11+i}}+F_0(s_0+11+i)\big)+\|\rho\|_{H^{s_0+7+i}}.
\end{align}
Applying Lemma \ref{Lem-top}-$($i$)$ we get $\mathcal{M}_n^k\in\DDD_k(n) $ and therefore
$$
\forall\, 1\leqslant i\leqslant k,\quad (\mathcal{M}_n^k\mathcal{Y}_n)(i)=0.
$$
Combining  this with  \eqref{GKL01} we find 
\begin{align}\label{Tyu1}
\nonumber\sum_{k=0}^{n}\|\mathcal{M}_n^k\mathcal{Y}_n(i)\|_{L^2(\T^{d+2})}&=\sum_{k=0}^{i-1}\|\mathcal{M}_n^k\mathcal{Y}_n(i)\|_{L^2(\T^{d+2})}\\
&\leqslant C 
\big(1+G^i(0)\big)G(i).
\end{align}
Now plugging \eqref{Tyu1} into \eqref{YY11} allows to get
\begin{align}\label{Ind-step1}
\zeta(t,1)
 &\leqslant C  e^{Ct}
\big(1+G(0)\big)G(1).
\end{align}
We intend to prove by induction  that for any $i\in\{1,..,n+1\},$
\begin{align}\label{Est-In13}
\zeta(t,i)
 \leqslant 2 C \big(1+G^i(0)\big)G(i) e^{\lambda\,C_0t}\quad\hbox{with}\quad C_0=C(1+\|\rho\|_{H^{s_0+3}}^{n}),
\end{align}
where $\lambda \geqslant 1$ to be fixed later.
According to \eqref{Ind-step1} this is true for $i=1$. Now assume that \eqref {Est-In13} is true  from $1$ up to $i-1$ and let us check that it remains  true at the order $i$. According to   \eqref{Es-z14} and \eqref{Tyu1} we may write
\begin{align*}
\zeta(t,i)
   \nonumber&\leqslant  C\big(1+G^i(0)\big)G(i) e^{Ct}\\
 \nonumber &+2C(\lambda-1)^{-1}\, e^{\lambda\, C_0t} \sum_{\ell=1}^{i-1}\big(1+G^{\ell}(0)\big)\|\rho\|_{H^{s_0+3+i-\ell}}G(\ell).
\end{align*}
Using interpolation inequality, since $G$ defined in \eqref{GKL2} is $\log-$convex, we get for $\ell\in\{0,..,i\}$
\begin{align*}
\|\rho\|_{H^{s_0+3+i-\ell}}G(\ell)&\leqslant G(i-\ell)G(\ell)\\
&\leqslant G(0) G(i).
\end{align*}
Therefore we obtain for any $\lambda\geqslant 1$
\begin{align*}
 \nonumber\zeta(t,i)
  &\leqslant   C \big(1+G^i(0)\big) G(i)\,e^{\lambda C_0t}+2C(\lambda-1)^{-1}\, i e^{\lambda C_0t}\big(1+G^{i}(0)\big) G(i).
\end{align*}
By fixing  $\lambda=1+2\,n,$ we find the estimate \eqref{Est-In13}. Hence we obtain  for any $i\in\{0,1,..n\}$
\begin{align}\label{Y-K1}
\nonumber \|\partial_\chi^i \mathcal{K}_t\|_{L^2(\T^{d+2})}=&\|\mathcal{Z}_{0}(t,i+1)\|_{L^2}\\
  \nonumber &\leqslant e^{Ct\|\rho\|_{W^{1,\infty}}} \zeta(t,i+1)\\
  &\quad \leqslant  C_1 \big(1+G^{i+1}(0)\big) e^{C_1(1+\|\rho\|_{H^{s_0+3}}^n)t} G(1+i).
\end{align}
for some constant $C_1$ that depends only on $n$. This achieves the proof of \eqref{Est-kernel-i} for $t\in[0,1]$.\\
Notice that making slight modifications of the preceding computations using in particular the second estimate in \eqref{day-98} and the second estimate of Lemma \ref{V-13}-(ii) yield instead of \eqref{Y-K1}
\begin{align}\label{Y-KL1}
 \|\partial_\chi^i \mathcal{K}_t\|_{L^2(\T^{d+2})}&\leqslant  C_1 \big(1+G^{i}(0)\big) e^{C_1(1+\|\rho\|_{H^{s_0+3}}^n)t}G(1+i)F_0(s_0+11+i).
\end{align}

$\bullet$ {\it  Estimate of $\|\partial_\varphi^i\mathcal{K}_t\|_{L^2(\T^{d+2})}$}. It can be done in in a similar way to \eqref{Est-kernel-i} and to avoid redundancy   we shall  merely explain the main steps  to estimate $\partial_\varphi^i\mathcal{K}_t$.  The goal is to  check  that for any  $\,i\in\{0,1,..,n\},$
\begin{align}\label{Est-kernel-phi}
 \forall t\in[0,1],\,\| \partial_\varphi^i\mathcal{K}_t\|_{L^2(\T^{d+2})}
  &\leqslant  C_1 \big(1+G^{i+1}(0)\big) e^{C_1\|\rho\|_{H^{s_0+3}}^n} G(1+i).
\end{align}
Proceeding as for \eqref{Ham1} we may write
\begin{align*}
 \mathscr{L}\partial_{\varphi}^{i}\mathcal{K}_t
&=\sum_{k=0}^{{i}-1}\widehat{m}_{i-k,{i}}\, \partial_{\varphi}^{k} \mathcal{K}_t+\partial_{\varphi}^{i} {K}_0
\end{align*}
with
\begin{align}
\widehat{m}_{k,i}&\triangleq \left(_k^i\right)\big(\widehat{\mathbb{B}}_{k}-\mathbb{T}_k\big)\\
\nonumber \widehat{\mathbb{B}}_{k}=\mathbb{B}_{\theta,(k)}+\mathbb{B}_{\eta,(k)},\quad& \mathbb{B}_{\theta,(k)}\triangleq \partial_\theta\big((\partial_\varphi^k\rho)|\textnormal{D}_\theta|^{\alpha-1}+|\textnormal{D}_\theta|^{\alpha-1}(\partial_\varphi^k\rho)\big)
\end{align}
and
\begin{align*}
\mathbb{T}_{k}\triangleq (\partial_\eta\partial_\varphi^{k}\rho)|\textnormal{D}_\eta|^{\alpha-1}+|\textnormal{D}_\eta|^{\alpha-1}(\partial_\eta\partial_\varphi^{k}\rho).
\end{align*}
Next, we  introduce the matrix operator
\begin{equation*}
 \widehat{\mathcal{M}}_n=\begin{pmatrix}
    0 &0&..&..&..&0  \\
 \widehat{m}_{1,1}&0&..&..&0&0 \\
   \widehat{m}_{2,2}& \widehat{m}_{1,2}&0&.. &..&0\\
    ..&..&..&..&0&0\\
\widehat{m}_{n,n}& \widehat{m}_{n-1,n}&..&..& \widehat{m}_{1,n}&0
  \end{pmatrix}
\end{equation*}
and
\begin{equation*}
\widehat{\mathcal{X}}_n(t)=\begin{pmatrix}
     \mathcal{K}_t &  \\
    \partial_\varphi \mathcal{K}_t& \\
    ..&\\
    ..&\\
    \partial_\varphi^n  \mathcal{K}_t
  \end{pmatrix},\quad\widehat{\mathcal{Y}}_n=\begin{pmatrix}
{{K}}_0 &  \\
    \partial_\varphi{K}_0& \\
    ..&\\
    ..&\\
    \partial_\varphi^n  {K}_0
  \end{pmatrix}.
  \end{equation*}
Then \eqref{Ham1} can be written in the matrix form as follows
\begin{align*}
\mathcal{L}_n\widehat{\mathcal{X}}_n(t)&=\widehat{\mathcal{M}}_n\,\widehat{\mathcal{X}}_n(t)+\widehat{\mathcal{Y}}_n,\quad \mathcal{L}_n\triangleq \mathscr{L}\textnormal{I}_{n+1}.
\end{align*}
Then we get the same structure as for $\mathcal{X}_n$ and therefore we may implement  exactly the same approach to deduce the estimate \eqref{Est-kernel-phi}.
This achieves the estimate  of the Theorem \ref{Prop-EgorV}-$($i$)$ in the case $m=0$.

$\blacktriangleright$ {\bf{Case} $m\in[-1,0).$ }
We shall first remove from $ \mathcal{R}(t)$ the first terms given by Taylor formula, that is, to consider the operator
$$\mathcal{R}_2(t)\triangleq \mathcal{R}(t)- \sum_{k=0}^2\textnormal{Ad}_{\mathbb{A}}^k\mathcal{R}_0 \frac{t^k}{k!},
$$ then it is straightforward to check 
\begin{equation*}
\left\{ \begin{array}{ll}
  \partial_t \mathcal{R}_2(t)=[\mathbb{A},\mathcal{R}_2(t)]+\textnormal{Ad}_{\mathbb{A}}^3\mathcal{R}_0 \frac{t^3}{3!}&\\
  \mathcal{R}_1(0)=0.
  \end{array}\right.
\end{equation*}
Then to  show the desired estimate it is enough  to check it for   $\mathcal{R}_2(t)$ and $\displaystyle{\sum_{k=0}^2\textnormal{Ad}_{\mathbb{A}}^k\mathcal{R}_0 \frac{t^k}{k!}}$.

$\bullet$ Estimate of $\mathcal{R}_2(t)$. Denote by $\mathcal{K}_{1,t} $ the kernel associated to $\mathcal{R}_2(t)$ and $K_1(t)$ the kernel associated to $\textnormal{Ad}_{\mathbb{A}}^3\mathcal{R}_0 \frac{t^3}{3!} $. Then using \eqref{Symb-Kern} combined with integration by parts and fixing $s=n\in\N$
\begin{align*}
\interleave\mathcal{R}_2(t)\interleave_{-1,s,0}&\leqslant \int_\T\|\partial_\eta \mathcal{K}_{1,t}(\cdot,\centerdot,\centerdot+\eta)\|_{H^{s}_{\varphi,\theta}}d\eta\\
&\lesssim\sum_{i=0}^n\big(\|\partial_\chi^i\partial_\eta \mathcal{K}_{1,t}\|_{L^2(\T^{d+2})}+\|\partial_\varphi^i\partial_\eta \mathcal{K}_{1,t}\|_{L^2(\T^{d+2})}\big).
\end{align*}
Set
$$\mathcal{K}_{2,t}\triangleq \partial_\eta \mathcal{K}_{1,t}\quad\hbox{and}\quad K_2\triangleq \partial_\eta K_1(t).
$$
Then differentiating \eqref{kernel-dyna} with respect to $\eta$ and using the identity
$$
\partial_\eta\left(-\mathbb{A}_\eta +\big[\partial_\eta,\mathcal{T}_\eta \big]\right)=-\mathbb{A}_\eta\partial_\eta
$$
 yield
\begin{equation*}
\partial_t\mathcal{K}_{2,t}(\varphi,\theta,\eta)-\mathbb{A}_\theta \mathcal{K}_{2,t}(\varphi,\theta,\eta)-\mathbb{A}_\eta \mathcal{K}_{2,t}(\varphi,\theta,\eta)=K_2(\varphi,\theta,\eta).
\end{equation*}
which is similar to the equation \eqref{kernel-dyna}. Notice that the operator associated to the kernel $K_2$ is given by $- \frac{t^3}{3!}\textnormal{Ad}_{\mathbb{A}}^3\mathcal{R}_0\partial_\theta$.   Then from \eqref{Y-K1}, \eqref{GKL1} and \eqref{WH1} we deduce that 
\begin{align}\label{Y-K2}
\forall t\in[0,1],\quad\|\partial_\chi^i \mathcal{K}_{2,t}\|_{L^2(\T^{d+2})}
  &\leqslant C_1 \big(1+G_1^{i+1}(0)\big) e^{C_1\|\rho\|_{H^{s_0+3}}^nt} G_1(1+i),
\end{align}
with
\begin{align}\label{GKL10}
G_1(i)\triangleq F\big(s_0+5+i\big)+\|\rho\|_{H^{s_0+6+|m|+i}},
\end{align}
where $F$ should be  an increasing  $\log-$convex function such that 
\begin{align}\label{WH10}
\forall s\geqslant s_0,\gamma\in\N\quad\hbox{with}\quad s+\gamma\leqslant s_0+n+1\Longrightarrow \interleave\textnormal{Ad}_{\mathbb{A}}^3\mathcal{R}_0\partial_\theta\interleave_{m,s,\gamma}\le F(s+\gamma),
\end{align}
for some fixed $m<-\frac12$. The next  goal is devoted to an explicit construction of $F$. According to Lemma \ref{comm-pseudo}-$($i$)$ one has for any $s\geqslant s_0,\gamma\in\N$
\begin{align*}
 \interleave\textnormal{Ad}_{\mathbb{A}}^3\mathcal{R}_0\partial_\theta\interleave_{3\alpha-2,s,\gamma}\lesssim  \interleave\textnormal{Ad}_{\mathbb{A}}^3\mathcal{R}_0\interleave_{3\alpha-3,s,\gamma}.\end{align*}
 Using  Lemma \ref{comm-pseudo}-$($ii$)$ combined with \eqref{WH-P1}, \eqref{D-K-G1}, \eqref{D-K-12} and $\alpha\in(0,\frac12)$
\begin{align*}
\interleave\textnormal{Ad}_{\mathbb{A}}^2\mathcal{R}_0\interleave_{2\alpha-2,s,\gamma}&\lesssim
\sum_{0\leqslant\beta\leqslant\gamma}\|\rho\|_{H^{s+\beta+2}}F_1(s_0+2+\gamma-\beta)+\|\rho\|_{H^{s_0+\beta+2}}F_1(s+2+\gamma-\beta)\\
&+\|\rho\|_{H^{s+1}}F_1(s_0+2+\gamma)+\|\rho\|_{H^{s_0+1}}F_1(s+2+\gamma) \\
&\quad+
\sum_{0\leqslant\beta\leqslant\gamma}\|\rho\|_{H^{s_0+3+\gamma-\beta}}F_1(s+\beta+1)+\|\rho\|_{H^{s+3+\gamma-\beta}}F_1(s_0+\beta+1)\\
&\qquad+\|\rho\|_{H^{s_0+4+\gamma}}F_1(s)+\|\rho\|_{H^{s+4+\gamma}}F_1(s_0).
\end{align*}
Thus by Sobolev embeddings we infer
\begin{align*}
\interleave\textnormal{Ad}_{\mathbb{A}}^2\mathcal{R}_0\interleave_{2\alpha-2,s,\gamma}&\lesssim
\sum_{0\leqslant\beta\leqslant\gamma}\|\rho\|_{H^{s+\beta+4}}F_1(s_0+2+\gamma-\beta)+\|\rho\|_{H^{s_0+\beta+4}}F_1(s+2+\gamma-\beta).
\end{align*}
Then using the $\log-$convexity of Sobolev norms and $F_1$ we find
for $s\geqslant s_0$ and $0\leqslant\beta\leqslant\gamma,$
\begin{align*}
\|\rho\|_{H^{s+\beta+4}}F_0(s_0+2+\gamma-\beta)&\le \|\rho\|_{H^{s_0}}^{\delta}\|\rho\|_{H^{s+\gamma+6}}^{1-\delta}F_0^{1-\delta}(s_0)F_0^{\delta}(s+6+\gamma)\\
&\le \Big(\|\rho\|_{H^{s_0}}+F_1(s_0)\Big)\Big(\|\rho\|_{H^{s+\gamma+6}}+F_1(s+6+\gamma)\Big)
\end{align*}
with $\delta\in(0,1)$. By setting
\begin{align}\label{Ta-M01}
F_2(s)\triangleq C\Big(\|\rho\|_{H^{s_0}}+F_1(s_0)\Big)\Big(\|\rho\|_{H^{s+6}}+F_1(s+6)\Big)
\end{align}
we get that $F_2$ is $\log-$ convex and 
\begin{align}\label{Ta-M0D1}
\interleave\textnormal{Ad}_{\mathbb{A}}^2\mathcal{R}_0\interleave_{2\alpha-2,s,\gamma}&\leqslant F_2(s+\gamma).
\end{align}
Implementing the same approach we also get
\begin{align}\label{UFF1}
\interleave\textnormal{Ad}_{\mathbb{A}}^3\mathcal{R}_0\interleave_{3\alpha-3,s,\gamma}&\lesssim
 F_3(s+\gamma)
\end{align}
with
\begin{align}\label{Ta-M1}
F_3(s)\triangleq C\Big(\|\rho\|_{H^{s_0}}+F_2(s_0)\Big)\Big(\|\rho\|_{H^{s+7}}+F_2(s+7)\Big).
\end{align}
Putting together  \eqref{D-K-12}, \eqref{Ta-M01} and \eqref{Ta-M1} implies
$$
F_3(s)\leqslant C\Big(1+\|\rho\|_{H^{s_0+17}}^3+F_2^3(s_0+17)\Big)\Big(\|\rho\|_{H^{s+17}}+F_0(s+17)\Big).
$$
Define
\begin{align*}
F(s)\triangleq C\Big(1+\|\rho\|_{H^{s_0+17}}^3+F_0^3(s_0+17)\Big)\Big(\|\rho\|_{H^{s+17}}+F_0(s+17)\Big)
\end{align*}
we get that $F$ is $\log-$convex and for any $\alpha\in(0,\frac12)$
\begin{align*}
\interleave\textnormal{Ad}_{\mathbb{A}}^3\mathcal{R}_0\interleave_{3\alpha-3,s,\gamma}&\leqslant
 F(s+\gamma).
\end{align*}
This gives \eqref{WH10} with  $m=3\alpha-3$ and the assumption $m<-\frac12$ is satisfied when $\alpha<\frac12.$
From \eqref{GKL10} we infer that
\begin{align*}
G_1(i)&\lesssim \|\rho\|_{H^{s_0+9+i}}+F(s_0+9+i)\\
&\lesssim \|\rho\|_{H^{s_0+9+i}}+\Big(1+\|\rho\|_{H^{s_0+17}}^3+F_0^3(s_0+17)\Big)\Big(\|\rho\|_{H^{s_0+26+i}}+F_0(s_0+26+i)\Big).
\end{align*}
Inserting this inequality  into \eqref{Y-K2} yields in view of Sobolev embeddings and for $t\in[0,1],$
 \begin{align*}
\|\partial_\chi^{i}\mathcal{K}_{2,t}\|_{L^2}&\leqslant C\big(1+F_0^{4(s+1)}{(s_0+27)}\big) {e^{C\|\rho\|_{H^{s_0+27}}^{s}}}\Big(F_0(s+27+s_0)+\|\rho\|_{H^{s+27+s_0}}\Big)\\
&\leqslant Ce^{C\mu^s(0)}\,\mu(s),
\end{align*}
with
$$
\mu(s)\triangleq F_0(s+27+s_0)+\|\rho\|_{H^{s+27+s_0}}.
$$

$\bullet$ {\it Estimate of }$\displaystyle{\sum_{k=0}^2\textnormal{Ad}_{\mathbb{A}}^k\mathcal{R}_0 \frac{t^k}{k!}}$. First for $k=0$ we know by assumption that 
\begin{align*}
\interleave\textnormal{Ad}_{\mathbb{A}}^0\mathcal{R}_0\interleave_{m,s,\gamma}=\interleave\mathcal{R}_0\interleave_{m,s,\gamma}\leqslant F_0(s+\gamma).
\end{align*}
Concerning $k=2$ we use \eqref{Ta-M0D1} which gives, since $2\alpha-2\leqslant -1\leqslant m$,  
\begin{align*}
\nonumber \interleave\textnormal{Ad}_{\mathbb{A}}^2\mathcal{R}_0\interleave_{m,s,\gamma}&\leqslant\interleave\textnormal{Ad}_{\mathbb{A}}^2\mathcal{R}_0\interleave_{2\alpha-2,s,\gamma}\\
&\leqslant F_2(s+\gamma).
\end{align*}
As to case $k=1$ we proceed similarly to \eqref{es-zk6}. Indeed, according to Lemma \ref{comm-pseudo}-$($ii$)$ applied with $q=0$  one gets for $m\in[-1,0]$
\begin{align*}
\nonumber \interleave[\mathbb{A},\mathcal{R}_0]\interleave_{\alpha+m-1,s,\gamma}&\lesssim
\sum_{0\leqslant\beta\leqslant\gamma}\interleave\mathbb{A}\interleave_{\alpha,s,1+\beta}\interleave\mathcal{R}_0\interleave_{m,s_0+2,\gamma-\beta}+\interleave\mathbb{A}\interleave_{\alpha,s_0,1+\beta}\interleave\mathcal{R}_0\interleave_{m,s+2,\gamma-\beta}\\
\nonumber&\quad+\sum_{ 0\le\beta\le\gamma}\interleave \mathbb{A}\interleave_{\alpha,s,0}\interleave \mathcal{R}_0\interleave_{m,s_0+2+\beta,\gamma-\beta} +\interleave \mathbb{A}\interleave_{0,s_0,0}\interleave\mathcal{R}_0\interleave_{m,s+2+\beta,\gamma-\beta} \\
\nonumber&\qquad+
\sum_{0\leqslant\beta\leqslant\gamma}\interleave \mathcal{R}_0\interleave_{m,s,1+\beta}\interleave\mathbb{A}\interleave_{\alpha,s_0+2,\gamma-\beta}+\interleave\mathcal{R}_0\interleave_{m,s_0,1+\beta}\interleave\mathcal{A}\interleave_{\alpha,s+2,\gamma-\beta}\\
&\qquad\quad\,\,+\sum_{ 0\leqslant\beta\leqslant\gamma}\interleave \mathcal{R}_0\interleave_{0,s,0}\interleave \mathbb{A}\interleave_{\alpha,s_0+3+\beta,\gamma-\beta} +\interleave \mathcal{R}_0\interleave_{0,s_0,0}\interleave\mathbb{A}\interleave_{\alpha,s+3+\beta,\gamma-\beta}.
\end{align*}
Therefore implementing the same arguments as for getting \eqref{D-K-G1} we obtain
\begin{align*}
\interleave[\mathbb{A},\mathcal{R}_0]\interleave_{m,s,\gamma}&\leqslant\interleave[\mathbb{A},\mathcal{R}_0]\interleave_{\alpha+m-1,s,\gamma}\\
 &\leqslant F_1(s+\gamma+1)
\end{align*}
where $F_1$ is defined in \eqref{D-K-12}.

\smallskip

${\bf{(ii)}}$ The case $q=0$ has been done in the first point $(i)$. We observe from \eqref{Symb-Kern} that for any $0\leqslant j\leqslant q$ and for $s=n\in\NN$
\begin{align*}
\interleave\partial_\lambda^j\mathcal{R}_1(t)\interleave_{0,s,0}&\leqslant \int_\T\|\partial_\lambda^j\mathcal{K}_t(\cdot,\centerdot,\centerdot+\eta)\|_{H^{s}_{\varphi,\theta}}d\eta\\
&\lesssim\sum_{i=0}^n\big(\|\partial_\lambda^j\partial_\chi^i\mathcal{K}_t\|_{L^2(\T^{d+2})}+\||\partial_\lambda^j\partial_\varphi^i\mathcal{K}_t\|_{L^2(\T^{d+2})}\big).
\end{align*}
The estimates of $\partial_\lambda^j\partial_\chi^i\mathcal{K}_t$ and $\partial_\lambda^j\partial_\varphi^i\mathcal{K}_t$ can be done in a similar way. Therefore  we shall restrict the discussion to the estimate of $\partial_\lambda^j\partial_\chi^i\mathcal{K}_t$ .

{$\bullet$ \bf Estimate of  $\|\partial_\lambda^j\partial_\chi^i\mathcal{K}_t\|_{{L^2(\T^{d+2})}},\, 0\leqslant j\leqslant q, 0\leqslant i\leqslant n$.} 
By applying $\partial_\lambda^j$ to the equation  \eqref{kernel-dyna1} we get after straightforward algebraic computations 
\begin{align}\label{Ham1X} 
 \kappa^j\mathscr{L}\partial_{\lambda}^{j}\mathcal{K}_t
&=\sum_{k=0}^{{j}-1}\overline{m}_{j-k,{j}}\,  \kappa^k\partial_{\lambda}^{k} \mathcal{K}_t+ \kappa^j\partial_{\lambda}^{j} {K}_0,
\end{align}
where
\begin{align}\label{LuX1}
\overline{m}_{k,j}&\triangleq \left(_k^j\right)\big(\overline{\mathbb{B}}_{k}-\overline{\mathbb{T}}_k\big)
\end{align}
with
\begin{align*}
\overline{\mathbb{B}}_{k}=\overline{\mathbb{B}}_{\theta,(k)}+\overline{\mathbb{B}}_{\eta,(k)},\quad \overline{\mathbb{B}}_{\theta,(k)}=\kappa^k\big[\partial_\lambda^k,\mathbb{A}_\theta \big]
\end{align*}
and
\begin{align*}
\overline{\mathbb{T}}_{k}\triangleq \kappa^k\big[\partial_\lambda^k,\mathcal{S}_0    \big].
\end{align*}
Recall that the operator  $\mathcal{S}_0$ is defined in \eqref{diffM1}.
Set
\begin{align}\label{LUX2}
\nonumber{X}_q&=\begin{pmatrix}
     \mathcal{K}_t &  \\
    \gamma\partial_\lambda \mathcal{K}_t& \\
    ..&\\
    ..&\\
    \gamma^q\partial_\lambda^q  \mathcal{K}_t
  \end{pmatrix},\quad{Y}_q=\begin{pmatrix}
{{K}}_0 &  \\
    \gamma\partial_\lambda{K}_0& \\
    ..&\\
    ..&\\
    \gamma^q\partial_\lambda^q  {K}_0
  \end{pmatrix}\\
  \quad\hbox{and}\quad 
 \overline{{M}}_q&=\begin{pmatrix}
    0 &0&..&..&..&0  \\
 \overline{m}_{1,1}&0&..&..&..&0 \\
  \overline{m}_{2,2}&  \overline{m}_{1,2}&0&.. &..&0\\
    ..&..&..&..&0&0\\
      ..&..&..&\overline{m}_{1,q-1}&0&0\\
 \overline{m}_{q,q}&  \overline{m}_{q-1,q}&..&..& \overline{m}_{1,q}&0
  \end{pmatrix}.
\end{align}
Then \eqref{Ham1X} can be written in the matrix form
\begin{align}\label{LUXPP2}
L_q {X}_q(t)&= \overline{{M}}_q\,{X}_q+{{Y}}_q,\quad{L}_q\triangleq \mathscr{L}\textnormal{I}_{q+1}.
\end{align}
Define
\begin{equation}\label{ourta1}
\overline{\mathcal{X}}_n(t)=\begin{pmatrix}
    {X}_q &  \\
    \partial_\chi{X}_q& \\
    ..&\\
    ..&\\
    \partial_\chi^n {X}_q
  \end{pmatrix}\quad\hbox{and}\quad\overline{\mathcal{Y}}_n=\begin{pmatrix}
{{Y}}_q &  \\
    \partial_\chi {{Y}}_q& \\
    ..&\\
    ..&\\
    \partial_\chi^n  {{Y}}_q
  \end{pmatrix}.
  \end{equation}
Therefore successive differentiation of \eqref{LUXPP2} yields to  the matrix form
\begin{align*}
\overline{\mathcal{L}}_n\overline{\mathcal{X}}_n(t)&=\overline{\mathcal{M}}_n\,\overline{\mathcal{X}}_n(t)+\overline{\mathcal{Y}}_n,\quad \overline{\mathcal{L}}_n\triangleq {L}_q\mathbb{I}_{n+1},
\end{align*}
with the bloc matrix operator
\begin{equation}\label{blocmatrix}
\overline{\mathcal{M}}_n=\begin{pmatrix}
     \overline{{M}}_q &0&..&..&..&0  \\
 \overline{M}_{1,1}& \overline{{M}}_q&0&..&..&0 \\
  \overline{M}_{2,2}&  \overline{M}_{1,2}& \overline{{M}}_q&0 &..&0\\
  .. &..&..&..&..&0  \\
    ..&..&..&\overline{M}_{1,n-1}& \overline{{M}}_q&0\\
 \overline{M}_{n,n}&  \overline{M}_{n-1,n}&..&..& \overline{M}_{1,n}& \overline{{M}}_q
  \end{pmatrix}\triangleq \overline{M}_q\mathbb{I}_{n+1}+N_n\in \DDD_0(n,q).  \end{equation}
  and
  \begin{align}\label{blocM}
\overline{M}_{k,i}&\triangleq \left(_k^i\right)\big(\textnormal{Ad}_{\partial_\chi}^k {\overline{M}}_q-\textnormal{Ad}_{\partial_\chi}^k\overline{L}_q\big), \,\overline{L}_q=\widehat{\mathscr{L}}\,\,\textnormal{I}_{q+1},
\end{align}
where $\widehat{\mathscr{L}}$ was defined in \eqref{Tr01}.
Since $
\overline{M}_q^{q+1}=0$ and $ N_n^{n+1}=0,
$ then
\begin{equation}\label{Nilp-ot}
\overline{\mathcal{M}}_n^{q+n+1}=0.
\end{equation}
Indeed, we observe that  $\overline{\mathcal{M}}_n^{q+n+1}$ can be expanded as the sum of elements in the form
$$
\prod_{i=1}^m \overline{M}_q^{\alpha_i}\mathbb{I}_{n+1}N_n^{\beta_i}\quad\hbox{with}\quad \sum_{i=1}^m(\alpha_i+\beta_i)=n+q+1.
$$
According to  Lemma \ref{Lem-top}, since $N_n\in\DDD_{1}(n,q)$  then  $N_n^{\beta_i}\in\DDD_{\beta_i}(n,q)$ and $\overline{M}_q^{\alpha_i}\mathbb{I}_{n+1}N_n^{\beta_i}\in\DDD_{\beta_i}(n,q)$. The same lemma gives
$$
P_m\triangleq \prod_{i=1}^m \overline{M}_q^{\alpha_i}\mathbb{I}_{n+1}N_n^{\beta_i}\in\DDD_{\min(\sum_{i=1}^m\beta_i,n+1)}(n,q).
$$
Consequently, if $\displaystyle{\sum_{i=1}^m\beta_i\geqslant n+1}$ then we get $
P_m\in \DDD_{n+1}(n,q)=\{0\}.$
 In the case  $\displaystyle{\sum_{i=1}^m\beta_i\leqslant n}$  we have necessary $\displaystyle{\sum_{i=1}^m\alpha_i\geqslant q+1}.$ Moreover the bloc entries of  $P_m$ is the sum of terms in the form
 $$
 \prod_{i=1}^m \overline{M}_q^{\alpha_i}C_{n,i},\quad\hbox{with}\quad C_{n,i}\in\DDD_0(q).
 $$
 Since $ \overline{M}_q\in\DDD_1(q)$ then applying  once again  Lemma \ref{Lem-top} we find successively  that $ \overline{M}_q^{\alpha_i}\in\DDD_{\min(\alpha_i,q+1)}(q),$ $  \overline{M}_q^{\alpha_i}C_{n,i}\in\DDD_{\min(\alpha_i,q+1)}(q)$ and 
 $$
 \prod_{i=1}^m \overline{M}_q^{\alpha_i}C_{n,i}\in\DDD_{\min(\sum_{i=1}^m\alpha_i,q+1)}(q)=\DDD_{q+1)}(q)=\{0\}.
 $$
 Thus we get that $P_m=0$ and this achieves the proof of \eqref{Nilp-ot}. We point out that similar arguments  give the following result
 \begin{equation}\label{Nilp-ot-gen}
\forall\, k\in\{0,..,n+1\},\quad \overline{\mathcal{M}}_n^{q+k}\in\DDD_k(n,q).
\end{equation}
Now, in a similar way to  \eqref{Hamb2} we deduce  for any $k\in\{1,.,n+q\}$
\begin{align}\label{Hamb2X}
\nonumber \overline{\mathcal{L}}_n \overline{\mathcal{M}}_n^k \overline{\mathcal{X}}_n(t)&=-\big[ \overline{\mathcal{M}}_n^k, \overline{\mathcal{L}}_n\big]\mathcal{X}_n(t)+ \overline{\mathcal{M}}_n^{k+1} \overline{\mathcal{X}}_n(t)+ \overline{\mathcal{M}}_n^{k} \overline{\mathcal{Y}}_n\\
&=-\sum_{j=0}^{k-1}\left(_j^k\right)\textnormal{Ad}_{ \overline{\mathcal{M}}_n}^{k-j}( \overline{\mathcal{L}}_n)\,  \overline{\mathcal{M}}_n^{j} \overline{\mathcal{X}}_n(t)+ \overline{\mathcal{M}}_n^{k+1} \overline{\mathcal{X}}_n(t)+ \overline{\mathcal{M}}_n^{k} \overline{\mathcal{Y}}_n.
\end{align}
Set $ \overline{\mathcal{Z}}_k(t)= \overline{\mathcal{M}}_n^k \overline{\mathcal{X}}_n(t)$ and denote by $ \overline{\mathcal{Z}}_k(t,i)$ the $i^{\textnormal{th}}$ bloc  vector of size $q+1$  .
Then applying the estimate \eqref{Hyp-est}  we get for any $k\in\{1,..,n+q\}$
\begin{align}\label{Hamb3XX}
 \nonumber\| \overline{\mathcal{Z}}_{k}(t,i)\|_{L^2}&\lesssim e^{Ct\|\rho\|_{W^{1,\infty}}}\Bigg(\bigintsss_0^t e^{-C\tau\|\rho\|_{W^{1,\infty}}}\mathlarger{\sum_{j=0}^{k-1}}\|\big(\textnormal{Ad}_{ \overline{\mathcal{M}}_n}^{k-j}(\overline{\mathcal{L}}_n)\,  \overline{\mathcal{M}}_n^{j}\overline{\mathcal{X}}_n\big)(\tau,i)\|_{L^2}d\tau\\
 &+\bigintsss_0^t e^{-C\tau\|\rho\|_{W^{1,\infty}}}\| \overline{\mathcal{Z}}_{k+1}(t,i)\|_{L^2}d\tau+t \| \overline{\mathcal{M}}_n^{k} \overline{\mathcal{Y}}_n(i)\|_{L^2} \Bigg)
\end{align}
and for $k=0$ it becomes
\begin{align*}
\|\overline{\mathcal{Z}}_{0}(t,i)\|_{L^2}\lesssim e^{Ct\|\rho\|_{W^{1,\infty}}}\left( \bigintsss_0^t e^{-C\tau\|\rho\|_{W^{1,\infty}}}\|\overline{\mathcal{Z}}_{1}(t,i)\|_{L^2}d\tau+ t \|\overline{\mathcal{Y}}_n(i)\|_{L^2}\right).
\end{align*}
By applying Lemma \ref{lem-com-it2} we deduce for $i\in\{1,..,n+1\}$
$$
\left\|\left(\textnormal{Ad}_{\overline{\mathcal{M}}_n}^{k-j}(\overline{\mathcal{L}}_n)\, \overline{\mathcal{M}}_n^{j}\overline{\mathcal{X}}_n\right)(\tau,i)\right\|_{L^2(\T^{d+2})}\lesssim  \big(\|\rho\|_{{s_0+3}}^{q,\kappa}\big)^{k-j}\mathlarger{\sum_{\ell=1}^{i}}  \|\rho\|_{{s_0+3+i-\ell}}^{q,\kappa}\,\,\,\left\|\overline{\mathcal{Z}}_j(\tau,\ell)\right\|_{L^2(\T^{d+2})}.
$$
Then setting 
$$
\overline{\zeta}(t,i)=e^{-Ct \|\rho\|_{W^{1,\infty}}}\sum_{k=0}^{n+q} \|\overline{\mathcal{Z}}_{k}(t,i)\|_{L^2}
$$
and  summing  up \eqref{Hamb3XX} over $k$, using in a crucial way  the relation  \eqref{Nilp-ot}, which gives   $\overline{\mathcal{Z}}_{n+q+1}(t)=0$, we find 
\begin{align*}
 \nonumber\overline{\zeta}(t,i) \lesssim  &\bigintsss_0^t \mathlarger{\sum_{k=1}^{n+q}\sum_{j=0}^{k-1}}\big(\|\rho\|_{{s_0+3}}^{q,\kappa}\big)^{k-j}\sum_{\ell=1}^{i}  \|\rho\|_{{s_0+3+i-\ell}}^{q,\kappa}\,\,\,\left\|\overline{\mathcal{Z}}_j(\tau,\ell)\right\|_{L^2(\T^{d+2})}d\tau\\
 &+\int_0^t \overline{\zeta}(\tau,i)\|_{L^2}d\tau+t\sum_{k=0}^{n+q} \|\overline{\mathcal{M}}_n^{k}\overline{\mathcal{Y}}_n(i)\|_{L^2}.
\end{align*}
From this we deduce that  for $i\in\{1,..,n+1\}$
\begin{align*}
 \nonumber\overline{\zeta}(t,i)
  \lesssim  &\big(\|\rho\|_{{s_0+3}}^{q,\kappa}+\big(\|\rho\|_{{s_0+3}}^{q,\kappa)}\big)^{n+q}\big) \sum_{\ell=1}^{i}\|\rho\|_{{s_0+3+i-\ell}}^{q,\kappa}\,\,\int_0^t\overline{}\overline{\zeta}(\tau,\ell)d\tau\\
 &+\int_0^t \overline{\zeta}(\tau,i)\|_{L^2}d\tau+t\sum_{k=0}^{n+q} \|\overline{\mathcal{M}}_n^{k}\overline{\mathcal{Y}}_n(i)\|_{L^2}.
\end{align*}
Then Gronwall inequality gives 
\begin{align}\label{Es-z140}
\overline{\zeta}(t,i)
  \leqslant & C_0\sum_{\ell=1}^{i-1}\|\rho\|_{H^{s_0+3+i-\ell}}\,\,\int_0^t e^{C(t-\tau)}\overline{\zeta}(\tau,\ell)d\tau+ C te^{Ct}\sum_{k=0}^{n+q} \|\overline{\mathcal{M}}_n^{k}\overline{\mathcal{Y}}_n(i)\|_{L^2},
\end{align}
with
$$
C_0\triangleq C\big(1+\big(\|\rho\|_{{s_0+3}}^{q,\kappa}\big)^{n+q}\big).
$$
For $i=1$ the preceding inequality becomes
\begin{align}\label{YY111}
\zeta(t,1)
 \leqslant  C t e^{C_0t}\sum_{k=0}^{n+q} \|\overline{\mathcal{M}}_n^{k}\overline{\mathcal{Y}}_n(1)\|_{L^2}.
\end{align}
From Lemma \ref{V-13}-$($ii$)$ we have 
\begin{align}\label{GKL01M}
\|\overline{\mathcal{M}}_n^{k}\overline{\mathcal{Y}}_n(i)\|_{L^2(\T^{d+2})}\leqslant C\big(1+\overline{G}^k(0)\big)\overline{G}(i),
\end{align}
with
\begin{align*}
\overline{G}(i)\triangleq\overline{F}_1\big(s_0+5+i\big)+\|\rho\|_{{s_0+6+|m|+i}}^{q,\kappa},
\end{align*}
where $\overline{F}_1$ should be  an increasing  $\log-$convex function such that 
\begin{align*}
\forall s\geqslant s_0,\gamma\in\N\quad\hbox{with}\quad s+\gamma\leqslant s_0+n+1\Longrightarrow\interleave [\mathbb{A},\mathcal{R}_0]\interleave_{m,s,\gamma}^{q,\kappa}\le \overline{F}_1(s+\gamma),
\end{align*}
for some fixed $m<-\frac12$.  The construction of $F_1$ can be done in a similar way to that of  \eqref{WH1} and one gets similarly to \eqref{D-K-G1} and \eqref{D-K-12}
\begin{align*}
\overline{F}_1(s)=C
 \big(\|\rho\|_{{s_0}}^{q,\kappa}+\overline{F}_0(s_0)\big)\big(\|\rho\|_{{s+5}}^{q,\kappa}+\overline{F}_0(s+5)\big).
\end{align*} 
Then we conclude using an induction principle as  for \eqref{Est-In13}.

\smallskip

{\bf{(iii)}} With the notation  $\mathcal{R}^i(t)\triangleq  \Phi_i(t)\mathcal{R}_0\Phi_i(-t)-\mathcal{R}_0$ and $\mathcal{R}(t)\triangleq\mathcal{R}^1(t)-\mathcal{R}^2(t)$ then we write from Heisenberg equations that 
\begin{equation}\label{Heisenberg-11-M}
\left\{ \begin{array}{ll}
  \partial_t \mathcal{R}^i(t)=[\mathbb{A}^i,\mathcal{R}^i(t)]+[\mathbb{A}^i,\mathcal{R}_0]&\\
  \mathcal{R}(0)=0,
  \end{array}\right.
\end{equation}
where $\mathbb{A}^i$ is the operator \eqref{AA-S1} associated to $\rho_i$.
Therefore the kernel $\mathcal{K}_t$ of $\mathcal{R}(t)$ satisfies in a similar way  to \eqref{Tr01} and \eqref{kernel-dyna1}
\begin{align}\label{Tr0Z1}
\partial_t\mathcal{K}_t+\widehat{\mathscr{L}}_1(\mathcal{K}_t)&=(\widehat{\mathscr{L}}_2-\widehat{\mathscr{L}}_1)\mathcal{K}_t^2+K_{1,2},
\end{align}
where  $\widehat{\mathscr{L}}_i=-\widehat{\mathbb{A}} f+\mathcal{S}_0 $ being  the operator associated to $\rho_i$, $\mathcal{K}_t^i$ the kernel associated to $\mathcal{R}^i(t)$ and   $K_{1,2}$ the kernel associated to $[\mathbb{A}^2-\mathbb{A}^1,\mathcal{R}_0].$  Notice that we have used that the  kernel $\mathcal{K}_t^i$ satisfies the equation
\begin{align}\label{Tr0Z2}
\partial_t\mathcal{K}_t^i+\widehat{\mathscr{L}}_i(\mathcal{K}_t^i)&=K_i
\end{align}
where $K_i$ is the kernel of $[\mathbb{A}^i,\mathcal{R}_0]$.
Applying the energy estimate \eqref{Hyp-est} with \eqref{Tr0Z1} allows to get
\begin{align}\label{Hyp-est-MP}
\|\mathcal{K}_t\|_{L^2(\T^{d+2})}\leqslant e^{Ct\|\rho_1\|_{W^{1,\infty}}}\Big(\int_0^{t}\big\|(\widehat{\mathscr{L}}_2-\widehat{\mathscr{L}}_1)\mathcal{K}_\tau^2\big\|_{L^2(\T^{d+2})}d\tau+t\,\|K_{1,2}\|_{L^2(\T^{d+2})}\Big).
\end{align}
Then from \eqref{Tr0Z2} we deduce that
\begin{align}\label{Tr0ZZZ2}
\big(\partial_t+\widehat{\mathscr{L}}_2\big)(\widehat{\mathscr{L}}_2-\widehat{\mathscr{L}}_1)\mathcal{K}_t^2&=(\widehat{\mathscr{L}}_2-\widehat{\mathscr{L}}_1)K_2+\big[\widehat{\mathscr{L}}_2,(\widehat{\mathscr{L}}_2-\widehat{\mathscr{L}}_1) \big]\mathcal{K}_t^2.
\end{align}
Applying once again the energy estimate \eqref{Hyp-est}
\begin{align}\label{Hyp-est-MPP}
\big\|(\widehat{\mathscr{L}}_2-\widehat{\mathscr{L}}_1)\mathcal{K}_t^2\big\|_{L^2(\T^{d+2})}&\lesssim t e^{Ct\|\rho_2\|_{W^{1,\infty}}}\big\|(\widehat{\mathscr{L}}_2-\widehat{\mathscr{L}}_1)K_2\big\|_{L^2(\T^{d+2})}\\
\nonumber &+e^{Ct\|\rho_2\|_{W^{1,\infty}}}\int_0^{t}\big\|\big[\widehat{\mathscr{L}}_2,(\widehat{\mathscr{L}}_2-\widehat{\mathscr{L}}_1) \big]\mathcal{K}_\tau^2\big\|_{L^2(\T^{d+2})}d\tau.
\end{align}

Applying \eqref{Uw-X1}  with $q=\gamma=0$ and $s=s_0$ we find 
\begin{align*}
\interleave[\widehat{\mathscr{L}}_2,(\widehat{\mathscr{L}}_2-\widehat{\mathscr{L}}_1) ]\interleave_{0,s_0,0}&\lesssim
\|\rho_2\|_{H^{s_0+3}}\|\rho_1-\rho_2\|_{H^{s_0+3}}.
\end{align*}
By virtue of  Lemma \eqref{Lem-Rgv1}-(i) we deduce that
\begin{align}\label{Co-dida1}
\big\|\big[\widehat{\mathscr{L}}_2,(\widehat{\mathscr{L}}_2-\widehat{\mathscr{L}}_1) \big]\mathcal{K}_t^2\big\|_{L^2(\T^{d+2})}\lesssim
\|\rho_2\|_{H^{s_0+3}} \|\rho_1-\rho_2\|_{H^{s_0+3}}\big\|\mathcal{K}_t^2\big\|_{L^2(\T^{d+2})}.
\end{align}
Applying \eqref{Hyp-est} we find 
\begin{align}\label{Hyp-estM0}
\|\mathcal{K}_t^2\|_{L^2}\lesssim \|{K}_2\|_{L^2}t\,e^{Ct \|\rho_2\|_{W^{1,\infty}}}.
\end{align}
Recall that $K_2$ is the kernel of the commutator $[\mathbb{A}^2,\mathscr{R}_0]$ where $\mathbb{A}^2$ is the nonlocal operator constructed in \eqref{AA-S1}  from $\rho_2$. Therefore we deduce from \eqref{D-K-G1} and \eqref{D-K-12}
\begin{align}\label{DTR1}
\interleave[\mathbb{A},\mathcal{R}_0]\interleave_{\overline\alpha-1,s_0,0}&\leqslant  C\big(\|\rho_2\|_{H^{s_0+5}}+\overline{F}_0(s_0+5)\big)^2,
\end{align}
 Since $\overline\alpha-1<-\frac12$,
we obtain from Lemma \ref{lemma-Sym-R}-(iii) applied with $q=0$
\begin{align*}
\int_{\T} \|K_2(\cdot,\centerdot,\centerdot+\eta)\|_{s_0}^2d\eta
&\lesssim \interleave[\mathbb{A},\mathcal{R}_0]\interleave_{\overline\alpha-1,s,0}^2\\
&\lesssim \big(\|\rho_2\|_{H^{s_0+5}}+\overline{F}_0(s_0+5)\big)^4.
\end{align*}
Plugging this estimate into \eqref{Hyp-estM0} yields
\begin{align}\label{Hyp-estMM1}
\|\mathcal{K}_t^2\|_{L^2}\lesssim t\,e^{Ct \|\rho\|_{W^{1,\infty}}}\big(\|\rho_2\|_{H^{s_0+5}}+\overline{F}_0(s_0+5)\big)^2.
\end{align}
Therefore we deduce from \eqref{Co-dida1} and Sobolev embeddings that
\begin{align}\label{Kern-zT1}
\big\|\big[\widehat{\mathscr{L}}_2,(\widehat{\mathscr{L}}_2-\widehat{\mathscr{L}}_1) \big]\mathcal{K}_t^2\big\|_{L^2(\T^{d+2})}&\lesssim
\,e^{Ct \|\rho_2\|_{H^{s_0+3}}}\big(\|\rho_2\|_{H^{s_0+5}}+\overline{F}_0(s_0+5)\big)^2 \|\rho_1-\rho_2\|_{H^{s_0+3}}.
\end{align}
Let us now move to the estimate of $(\widehat{\mathscr{L}}_2-\widehat{\mathscr{L}}_1)K_2$. Then applying  Lemma \ref{lemm-iter1}-(iii) applied with $q=0$ and using \eqref{D-K-G1} and \eqref{D-K-12} yield
 \begin{align}\label{TJ-1}
\int_{\T}\big\|(\widehat{\mathscr{L}}_2-\widehat{\mathscr{L}}_1)K_2(\cdot,\centerdot,\centerdot+\eta)\big\|_{H^{s_0}}^2d\eta&\lesssim \overline{F}_0\big(s_0+10\big) \ \|\rho_1-\rho_2\|_{H^{s_0+5}}.
\end{align}
Putting together \eqref{TJ-1}, \eqref{Co-dida1} and \eqref{Hyp-est-MPP}, and using Sobolev embeddings give for $t\in[0,1]$
\begin{align}\label{Hyp-est-MPP2}
\big\|(\widehat{\mathscr{L}}_2-\widehat{\mathscr{L}}_1)\mathcal{K}_t^2\big\|_{L^2(\T^{d+2})}&\leqslant C e^{C\big(\|\rho_2\|_{H^{s_0+10}}+\overline{F}_0(s_0+10)\big)} \|\rho_1-\rho_2\|_{H^{s_0+5}}.
\end{align}
It remains to estimate $\|K_{1,2}\|_{L^2(\T^{d+2})}$. Then using \eqref{Rec-XK1} combined with Sobolev embeddings
\begin{align*}
\interleave[\mathbb{A}^1-\mathbb{A}^2,\mathcal{R}_0]\interleave_{\overline\alpha-1,s_0,0}&\lesssim 
\|\rho_1-\rho_2\|_{H^{s_0+3}}\overline{F}_0(s_0+3)).
\end{align*}
By virtue of Lemma \ref{lemma-Sym-R}-(iii) and Sobolev embeddings we deduce that
\begin{align}\label{Rec-XKK1}
\|K_{1,2}\|_{L^2(\T^{d+2})}\lesssim \|\rho_1-\rho_2\|_{H^{s_0+3}}\overline{F}_0(s_0+3)). 
\end{align}
Inserting \eqref{Rec-XKK1} and \eqref{Hyp-est-MPP2} into \eqref{Hyp-est-MP} implies for any $t\in[0,1]$
\begin{align}\label{Hyp-est-MY}
\|\mathcal{K}_t\|_{L^2(\T^{d+2})} \leqslant C e^{C\big(\|\rho_2\|_{H^{s_0+10}}+\|\rho_1\|_{H^{s_0+10}}+\overline{F}_0(s_0+10)\big)} \|\rho_1-\rho_2\|_{H^{s_0+5}}.
\end{align}
This achieves the proof for $q=0$. However for the genral case we may proceed as in \eqref{Ham1X}. Since the analysis is quite similar we shall only restrict the discussion to the main lines. First we start from \eqref{Tr0Z1} and differentiate it  successively with respect to $\lambda$  giving rise to an analogous system to \eqref{LUXPP2}. Actually, we obtain
\begin{align}\label{LUXPPR2}
L_q {X}_q(t)&= \overline{{M}}_q\,{X}_q+{{Y}}_q,\quad{L}_q\triangleq \mathscr{L}\textnormal{I}_{q+1}.
\end{align}
where 
\begin{equation*}
{X}_q=\begin{pmatrix}
     \mathcal{K}_t &  \\
    \gamma\partial_\lambda \mathcal{K}_t& \\
    ..&\\
    ..&\\
    \gamma^q\partial_\lambda^q  \mathcal{K}_t
  \end{pmatrix},\quad{Y}_q=\begin{pmatrix}
{{K}}_0 &  \\
    \gamma\partial_\lambda{K}_0& \\
    ..&\\
    ..&\\
    \gamma^q\partial_\lambda^q  {K}_0
  \end{pmatrix}\quad\hbox{and}\quad 
 \overline{{M}}_q=\begin{pmatrix}
    0 &0&..&..&..&0  \\
 \overline{m}_{1,1}&0&..&..&..&0 \\
  \overline{m}_{2,2}&  \overline{m}_{1,2}&0&.. &..&0\\
    ..&..&..&..&0&0\\
      ..&..&..&\overline{m}_{1,q-1}&0&0\\
 \overline{m}_{q,q}&  \overline{m}_{q-1,q}&..&..& \overline{m}_{1,q}&0
  \end{pmatrix}.
\end{equation*}
with $K_0=(\widehat{\mathscr{L}}_2-\widehat{\mathscr{L}}_1)\mathcal{K}_t^2+K_{1,2}$ and the entries $\overline{m}_{ij}$ are defined ina similar way to \eqref{LuX1}. Then we proceed in a straightforward way following the proof of   the first point $($i$)$. This achieved the proof.

 \end{proof}
 
\subsubsection{Tame estimates of the iterated commutators}
In this section we shall discuss an important result used during  the proof of Theorem \ref{Prop-EgorV} related to some tame estimates governing iterated  commutators of  pseudo-differential operators of special structure.  First we shall recall some operators seen in the preceding sections. The operator $\widehat{\mathcal{L}}_n$ is described through \eqref{Hamb1} and \eqref{Tr01} and  the matrix $\mathcal{M}_n$ is defined in    \eqref{Matrix-y} and takes the form
$$
\mathcal{M}_n=\begin{pmatrix}
    0 &0&..&..&..&0  \\
m_{1,1}&0&..&..&0&0 \\
   m_{2,2}&m_{1,2}&0&.. &..&0\\
    ..&..&..&..&0&0\\
    ..& ..&..&m_{1,n-1}& 0&0\\
 m_{n,n}& m_{n-1,n}&..&..& m_{1,n}&0
  \end{pmatrix}\in\DDD_1(n).
$$
Using \eqref{Bnk}, \eqref{Ima1} and \eqref{Ima2} one finds that
\begin{align}\label{BnkU}
{m}_{k,i}=\left(_k^i\right)\big(\mu(\varphi,\theta,k)+\nu(\varphi,\eta,k)\big)
\end{align}
with
\begin{align}\label{BnkV}
\mu(\varphi,\theta,k)&\triangleq\partial_\theta\big((\partial_\theta^k\rho)|\textnormal{D}_\theta|^{\alpha-1}+|\textnormal{D}_\theta|^{\alpha-1}(\partial_\theta^k\rho)\big)\\
\nonumber \nu(\varphi,\eta,k)&\triangleq
(\partial_\eta^k\rho)\partial_\eta|\textnormal{D}_\eta|^{\alpha-1}+|\textnormal{D}_\eta|^{\alpha-1}(\partial_\eta^k\rho)\partial_\eta.
\end{align}
Define
\begin{equation}\label{Ck-m}
\mathcal{C}_{k}=\textnormal{Ad}_{\mathcal{M}_n}^{k}(\widehat{\mathcal{L}}_n).
\end{equation}
Then we have the recursive relation 
$$
\mathcal{C}_{k+1}=\big[\mathcal{M}_n, \mathcal{C}_{k}\big], 
$$
with the  initial relation 
\begin{align*}
\ \mathcal{C}_{0}&=\widehat{\mathscr{L}}\,\textnormal{I}_{n+1}\\
&=-\big(\mu(\varphi,\theta,0)+\nu(\varphi,\eta,0)\big) \textnormal{I}_{n+1}.
\end{align*}
The last identity comes from the definition of  $\widehat{\mathscr{L}}$ seen in \eqref{Tr01}. Direct computations show that 
\begin{align}
\nonumber\mathcal{C}_{1}&= \begin{pmatrix}
    0 &0&..&..&..&0  \\
 \big[m_{1,1},{\widehat{\mathscr{L}}}\,\,\big]&0&..&..&0&0 \\
  \big[m_{2,2},{\widehat{\mathscr{L}}}\,\,\big]& \big[m_{1,2},{\widehat{\mathscr{L}}}\,\,\big]&0&.. &..&0\\
    ..&..&..&..&0&0\\
   \big[m_{n,n},{\widehat{\mathscr{L}}}\,\,\big]& \big[m_{n-1,n},{\widehat{\mathscr{L}}}\,\,\big]&..&..& \big[m_{1,n},\widehat{\mathscr{L}}\,\,\big]&0
  \end{pmatrix}\in\DDD_1(n).
\end{align}
According to Lemma \ref{Lem-top} one has  
\begin{equation}\label{Ma-S}
\mathcal{C}_k=\begin{pmatrix}
    0 &0&..&..&..&0  \\
c^k_{1,1}&0&..&..&0&0 \\
   c^k_{2,2}& c^k_{1,2}&0&.. &..&0\\
    ..&..&..&..&0&0\\
   c^k_{n,n}& c^k_{n-1,n}&..&..& c^k_{1,n}&0
  \end{pmatrix}\in\DDD_k(n)
\end{equation}
with
\begin{align}\label{Nq-1}
c^k_{i,j}=\left(_i^j\right)\big(\mathcal{\mu}_k(\varphi,\theta,i)+\mathcal{\nu}_k(\varphi,\eta,i)\big), \,\,\forall\, i\leqslant j\in\{1,..,n\}
\end{align}
and the law $\big[\mathcal{M}_n, \mathcal{C}_{k}\big]$ is given by 
\begin{align*}
\mathcal{\mu}_{\mathcal{M}_n}\boxast \mathcal{\mu}_{\mathcal{C}_{k}}={\mu}\boxast{\mu}_k+{\nu}\boxast {\nu}_k.
\end{align*}
Thus we find the recursive relations
\begin{align*}
\mathcal{\mu}_{k+1}={\mu}\boxast{\mu}_k&,\quad \mu_0(\varphi,\theta, i)=- \mu(\varphi,\theta, 0)\\
\mathcal{\nu}_{k+1}={\nu}\boxast {\nu}_k&,\quad \nu_0(\varphi,\eta, i)=- \nu(\varphi,\eta, 0).
\end{align*}
The estimates of $\mathcal{\mu}_{k+1}$ and $\mathcal{\nu}_{k+1}$ are similar and  we shall restrict the discussion only to  the first one. By definition of the convolution law one may write
$$
\mathcal{\mu}_{k+1}(\varphi,\theta,i)=\mathlarger{\sum_{\ell=0}^{i}}\left(_\ell^i\right)\big[\mathcal{\mu}(\varphi,\theta,i-\ell),\mathcal{\mu}_{k}(\varphi,\theta,\ell)\big].
$$
Since $\mathcal{M}_n\in\DDD_1(n)$ and $\mathcal{C}_k\in \DDD_k(n)$ then  
\begin{align}\label{Rec-eq}
\forall i\in\{k+1,..,n\},\quad \mathcal{\mu}_{k+1}(\varphi,\theta,i)&=\sum_{\ell=k}^{i-1}\left(_\ell^i\right)\big[\mathcal{\mu}(\varphi,\theta,i-\ell),\mathcal{\mu}_{k}(\varphi,\theta,\ell)\big],\\
\nonumber \forall i\in\{0,..,k\},\quad \mathcal{\mu}_{k+1}(\varphi,\theta,i) &=0.
\end{align}
Our first main result reads as follows.
\begin{lemma}\label{lem-com-it1}
Let $\gamma\in\N, s\geqslant s_0>\frac{d+1}{2},$ then for any $n\in\N^\star,\,k\in\{1,..,n\}, i\in\{k,..,n\}$
$$
 \interleave \mu_k(i)\interleave_{0,s,\gamma}+\interleave \nu_k(i)\interleave_{0,s,\gamma}\lesssim \|\rho\|_{H^{s+3+i+\gamma}}\|\rho\|_{H^{s_0+3}}^{k}.
$$
In addition, we have for any $s_0>d+1$
$$
\|\mu_k(i) h\|_{L^2(\T^{d+2})}+\|\nu_k(i) h\|_{L^2(\T^{d+2})}\lesssim \|\rho\|_{H^{s_0+3+i}(\T^{d+1})}\,\|\rho\|_{H^{s_0+3}(\T^{d+1})}^{k}\,\| h\|_{L^2(\T^{d+2})}.
$$
\end{lemma}
\begin{proof}
To alleviate the notation we denote $\mu_k(\varphi,\theta,i)=\mu_k(i)$. The goal is to prove the estimate using the induction  principle in $k$. Let us first check the estimate for $k=1$. In this case,
\begin{align*}\mathcal{\mu}_{1}(\varphi,\theta,i)&=\sum_{\ell=k}^{i-1}\left(_\ell^i\right)\big[\mathcal{\mu}(\varphi,\theta,i-\ell),\mathcal{\mu}_{0}(\varphi,\theta,\ell)\big]\\
&=-\sum_{\ell=0}^{i-1}\left(_\ell^i\right)\big[\mathcal{\mu}(\varphi,\theta,i-\ell),\mathcal{\mu}(\varphi,\theta,0)\big].
\end{align*}
Applying Lemma \ref{comm-pseudo}-$($iii$)$ with $\mathcal{A}=\mathcal{\mu}(i-\ell),{\mathcal{B}}=\mu(0), q=0$ and
$$
m_1=m_2=\alpha\in\left(0,\frac12\right),\quad {\mu}=\frac32,\mu_\beta=\overline{\mu}_\beta=\frac32+\beta
$$
yields
\begin{align*}
\interleave[\mathcal{\mu}(i-\ell),\mu(0)]\interleave_{0,s,\gamma}&\le \sum_{0\leqslant\beta\leqslant\gamma}\interleave \mathcal{A}\interleave_{\alpha,s,1+\gamma-\beta}\interleave{\mathcal{B}}\interleave_{\alpha,s_0+\mu,\beta}+\interleave \mathcal{A}\interleave_{\alpha,s_0,1+\gamma-\beta}\interleave{\mathcal{B}}\interleave_{\alpha,s+\mu,\beta}
\\
&\quad+ \sum_{0\leqslant\beta\leqslant\gamma}\interleave \mathcal{A}\interleave_{\alpha,s,0} \interleave {\mathcal{B}}\interleave_{\alpha,s_0+\mu_\beta,\gamma-\beta}+\interleave \mathcal{A}\interleave_{\alpha,s_0,0} \interleave {\mathcal{B}}\interleave_{\alpha,s+\mu_\beta,\gamma-\beta}\\
&\qquad+ \sum_{0\leqslant\beta\leqslant\gamma}\interleave\mathcal{A}\interleave_{\alpha,s_0+\mu,\beta}\interleave {\mathcal{B}}\interleave_{\alpha,s,1+\gamma-\beta}+\interleave\mathcal{A}\interleave_{\alpha,s+\mu,\beta}\interleave {\mathcal{B}}\interleave_{\alpha,s_0,1+\gamma-\beta}
\\
&\quad\qquad+\sum_{0\leqslant\beta\leqslant\gamma}\interleave \mathcal{A}\interleave_{\alpha,s_0+\mu_\beta,\gamma-\beta}\interleave {\mathcal{B}}\interleave_{\alpha,s,0} +\interleave \mathcal{A}\interleave_{\alpha,s+\mu_\beta,\gamma-\beta}\interleave {\mathcal{B}}\interleave_{\alpha,s_0,0} .
\end{align*}
Using \eqref{BnkV} and \eqref{WH-P1}, we obtain
\begin{align}\label{ES-M1}
\interleave\mathcal{\mu}(\ell)\interleave_{\alpha,s,\beta}&\lesssim \|\rho\|_{H^{s+\beta+\ell+1}}.
\end{align}
Therefore, combining the preceding two  estimates yields
\begin{align}\label{Init-01}
\nonumber \interleave[\mathcal{\mu}(i-\ell),\mu(0)]\interleave_{0,s,\gamma}&\le \sum_{0\leqslant\beta\leqslant\gamma}\|\rho\|_{H^{s+2+\gamma-\beta+i-\ell}}\|\rho\|_{H^{s_0+\frac52+\beta}}+\|\rho\|_{H^{s_0+2+\gamma-\beta+i-\ell}}\|\rho\|_{H^{s+\frac52+\beta}}
\\
\nonumber&\quad+\|\rho\|_{H^{s+1+i-\ell}}\|\rho\|_{H^{s_0+\frac52+\gamma}}+\|\rho\|_{H^{s_0+1+i-\ell}}\|\rho\|_{H^{s+\frac52+\gamma}}\\
\nonumber&\qquad+ \sum_{0\leqslant\beta\leqslant\gamma}\|\rho\|_{H^{s_0+\frac52+\beta+i-\ell}}\|\rho\|_{H^{s+2+\gamma-\beta}}+\|\rho\|_{H^{s+\frac52+\beta+i-\ell}}\|\rho\|_{H^{s_0+2+\gamma-\beta}}
\\
&\quad\qquad+\|\rho\|_{H^{s_0+\frac52+\gamma+i-\ell}}\|\rho\|_{H^{s+1}}+\|\rho\|_{H^{s+\frac52+\gamma+i-\ell}}\|\rho\|_{H^{s_0+1}} .
\end{align}
Using Sobolev embeddings and interpolation inequality we get for $s\geqslant s_0$
\begin{align*}
\|\rho\|_{H^{s+2+\gamma-\beta+i-\ell}}\|\rho\|_{H^{s_0+\frac52+\beta}}&\lesssim \|\rho\|_{H^{s+2+\gamma-\beta+i}}\|\rho\|_{H^{s_0+\frac52+\beta}}\\
&\lesssim \|\rho\|_{H^{s+3+i+\gamma}}\|\rho\|_{H^{s_0+\frac32}}
\end{align*}
and
\begin{align*}
\|\rho\|_{H^{s_0+2+\gamma-\beta+i-\ell}}\|\rho\|_{H^{s+\frac52+\beta}}&\lesssim \|\rho\|_{H^{s_0+2+\gamma-\beta+i}}\|\rho\|_{H^{s+\frac52+\beta}}\\
&\lesssim \|\rho\|_{H^{s+3+i+\gamma}}\|\rho\|_{H^{s_0+\frac32}}.
\end{align*}
Similarly we get
\begin{align*}
\|\rho\|_{H^{s+1+i-\ell}}\|\rho\|_{H^{s_0+\frac52+\gamma}}&\lesssim \|\rho\|_{H^{s+1+i}}\|\rho\|_{H^{s_0+\frac52+\gamma}}\\
&\lesssim \|\rho\|_{H^{s+3+i+\gamma}}\|\rho\|_{H^{s_0+\frac12}}
\end{align*}
and
\begin{align*}
\|\rho\|_{H^{s_0+1+i-\ell}}\|\rho\|_{H^{s+\frac52+\gamma}}&\lesssim \|\rho\|_{H^{s_0+1+i}}\|\rho\|_{H^{s+\frac52+\gamma}}\\
&\lesssim \|\rho\|_{H^{s+3+i+\gamma}}\|\rho\|_{H^{s_0+\frac12}}.
\end{align*}
On the other hand
\begin{align*}
\|\rho\|_{H^{s_0+\frac52+\beta+i-\ell}}\|\rho\|_{H^{s+2+\gamma-\beta}}&\lesssim \|\rho\|_{H^{s_0+\frac52+\beta+i}}\|\rho\|_{H^{s+2+\gamma-\beta}}\\
&\lesssim \|\rho\|_{H^{s+3+i+\gamma}}\|\rho\|_{H^{s_0+\frac32}}
\end{align*}
and
\begin{align*}
\|\rho\|_{H^{s+\frac52+\beta+i-\ell}}\|\rho\|_{H^{s_0+2+\gamma-\beta}}&\lesssim \|\rho\|_{H^{s+\frac52+\beta+i}}\|\rho\|_{H^{s_0+2+\gamma-\beta}}\\
&\lesssim \|\rho\|_{H^{s+3+i+\gamma}}\|\rho\|_{H^{s_0+\frac32}}.
\end{align*} 
As to the last terms in \eqref{Init-01} we implement similar arguments as before yielding 
\begin{align*}
\|\rho\|_{H^{s_0+\frac52+\gamma+i-\ell}}\|\rho\|_{H^{s+1}}&\lesssim\|\rho\|_{H^{s_0+\frac52+\gamma+i}}\|\rho\|_{H^{s+1}}\\
&\lesssim \|\rho\|_{H^{s+3+i+\gamma}}\|\rho\|_{H^{s_0+\frac12}}
\end{align*}
and
\begin{align*}
\|\rho\|_{H^{s+\frac52+\gamma+i-\ell}}\|\rho\|_{H^{s_0+1}}&\lesssim\|\rho\|_{H^{s+\frac52+\gamma+i}}\|\rho\|_{H^{s_0+1}}\\
&\lesssim \|\rho\|_{H^{s+3+i+\gamma}}\|\rho\|_{H^{s_0+\frac12}}.
\end{align*}
Putting together the preceding estimates and using Sobolev embeddings yield 
\begin{align*}
\interleave[\mathcal{\mu}(i-\ell),\mu(0)]\interleave_{0,s,\gamma}&\lesssim \|\rho\|_{H^{s+3+i+\gamma}}\|\rho\|_{H^{s_0+\frac32}}\\
&\lesssim \|\rho\|_{H^{s+3+i+\gamma}}\|\rho\|_{H^{s_0+3}}.
\end{align*}
This gives the result of Lemma \ref {lem-com-it1} for $k=1$.
 Now let us assume that this result is true at the order $k$ and  prove  it at the order $k+1$. According to the relation \eqref{Rec-eq} on gets
 \begin{align}\label{Big-k0}
\interleave\mathcal{\mu}_{k+1}(i)\interleave_{0,s,\gamma}\lesssim \sum_{\ell=k}^{i-1}\interleave[\mathcal{\mu}(i-\ell),\mathcal{\mu}_{k}(\ell)]\interleave_{0,s,\gamma}.
\end{align}
 Applying Lemma \ref{comm-pseudo}-$(\textnormal{iii})$ with $m_1=\alpha\in(0,\frac12), m_2=0, q=0$ and fixing   the choice 
$$
{\mu}=1, \quad \mu_\beta=\overline{\mu}_\beta=1+\beta
$$
 we find 
 \begin{align}\label{Big-k}
&\interleave[\mathcal{\mu}(i-\ell),\mathcal{\mu}_{k}(\ell)]\interleave_{0,s,\gamma}
\lesssim \sum_{0\leqslant\beta\leqslant\gamma}\interleave\mathcal{\mu}(i-\ell)\interleave_{\alpha,s,1+\gamma-\beta}\interleave\mathcal{\mu}_{k}(\ell)\interleave_{0,s_0+1,\beta}\\
\nonumber&\quad+\sum_{0\leqslant\beta\leqslant\gamma}\interleave\mathcal{\mu}(i-\ell)\interleave_{\alpha,s_0,1+\gamma-\beta}\interleave\mathcal{\mu}_{k}(\ell)\interleave_{0,s+1,\beta}
+ \interleave \mathcal{\mu}(i-\ell)\interleave_{\alpha,s,0} \interleave\mathcal{\mu}_{k}(\ell)\interleave_{0,s_0+1+\beta,\gamma-\beta}\\
\nonumber&\qquad+ \sum_{0\leqslant\beta\leqslant\gamma}\interleave \mathcal{\mu}(i-\ell)\interleave_{\alpha,s_0,0} \interleave \mathcal{\mu}_{k}(\ell)\interleave_{0,s+1+\beta,\gamma-\beta}+\interleave\mathcal{\mu}(i-\ell)\interleave_{\alpha,s_0+1,\beta}\interleave \mathcal{\mu}_{k}(\ell)\interleave_{0,s,1+\gamma-\beta}\\
\nonumber&\qquad\quad+\sum_{0\leqslant\beta\leqslant\gamma}\interleave\mathcal{\mu}(i-\ell)\interleave_{\alpha,s+1,\beta}\interleave \mathcal{\mu}_{k}(\ell)\interleave_{0,s_0,1+\gamma-\beta}
+\interleave\mathcal{\mu}(i-\ell)\interleave_{\alpha,s_0+1+\beta,\gamma-\beta}\interleave \mathcal{\mu}_{k}(\ell)\interleave_{0,s,0} \\
\nonumber&\qquad \qquad \quad+\sum_{0\leqslant\beta\leqslant\gamma} \interleave\mathcal{\mu}(i-\ell)\interleave_{\alpha,s+1+\beta,\gamma-\beta}\interleave \mathcal{\mu}_{k}(\ell)\interleave_{0,s_0,0}\triangleq  \sum_{0\leqslant\beta\leqslant\gamma}\sum_{j=1}^8\mathcal{I}_j^\beta.
\end{align}
Now let us estimate the first term $\mathcal{I}_1^\beta$. One has from the induction assumption and \eqref{ES-M1},
\begin{align*}
\mathcal{I}_1^\beta\lesssim&  \|\rho\|_{H^{s+2+\gamma-\beta+i-\ell}}\interleave\mathcal{\mu}_{k}(\ell)\interleave_{0,s_0+1,\beta}\\
\lesssim&\|\rho\|_{H^{s+2+\gamma-\beta+i-\ell}} \|\rho\|_{H^{s_0+4+\beta+\ell}} \|\rho\|_{H^{s_0+3}}^{k}.
\end{align*}
Using interpolation inequality we find, since $\ell\leqslant i-1,$ that for any $s\geqslant s_0$
$$
\|\rho\|_{H^{s+2+\gamma-\beta+i-\ell}}\|\rho\|_{H^{s_0+4+\beta+\ell}} \lesssim\|\rho\|_{H^{s_0+3}} \|\rho\|_{H^{s+3+\gamma+i}}.
$$
Consequently, 
\begin{align}\label{J-1}
\mathcal{I}_1^\beta\lesssim&\|\rho\|_{H^{s+3+\gamma+i}}\|\rho\|_{H^{s_0+3}}^{k+1}.
\end{align}
Concerning the second term $\mathcal{I}_2^\beta$ we write
\begin{align*}
\mathcal{I}_2^\beta\lesssim&  \|\rho\|_{H^{s_0+2+\gamma-\beta+i-\ell}}\interleave\mathcal{\mu}_{k}(\ell)\interleave_{0,s+1,\beta}\\
\lesssim&\|\rho\|_{H^{s_0+2+\gamma-\beta+i-\ell}} \|\rho\|_{H^{s+4+\ell+\beta}}\|\rho\|_{H^{s_0+3}}^{k}
\end{align*}
Interpolation inequalities allow to get since $\ell\leqslant i-1$
\begin{align*}
\|\rho\|_{s_0+2+\gamma-\beta+i-\ell} \|\rho\|_{H^{s+4+\ell+\beta}} \lesssim&\|\rho\|_{s_0+3} \|\rho\|_{H^{s+3+i+\gamma}}
\end{align*}
yields
\begin{align}\label{J-2}
 \mathcal{I}_2^\beta
\lesssim& \|\rho\|_{H^{s+3+i+\gamma}}\|\rho\|_{H^{s_0+3}}^{k+1}.
\end{align}
For the third term, we write
\begin{align*}
\mathcal{I}_3^\beta=& \interleave \mathcal{\mu}(i-\ell)\interleave_{\alpha,s,0} \interleave\mathcal{\mu}_{k}(\ell)\interleave_{0,s_0+1+\beta,\gamma-\beta}\\
\lesssim&  \|\rho\|_{H^{s+1+i-\ell}} \|\rho\|_{H^{s_0+4+\gamma+\ell}} \|\rho\|_{H^{s_0+3}}^{k}.
\end{align*}
Using the interpolation inequality, since $\ell\leqslant i-1$, we find  that for any $s\geqslant s_0$
$$
 \|\rho\|_{H^{s+1+i-\ell}} \|\rho\|_{H^{s_0+4+\gamma+\ell}}\lesssim  \|\rho\|_{H^{s_0+2}} \|\rho\|_{H^{s+3+\gamma+i}}.
$$
Hence, Sobolev embeddings imply
\begin{align}\label{J-3}
 \mathcal{I}_3^\beta\lesssim& \|\rho\|_{H^{s+3+i+\gamma}}\|\rho\|_{H^{s_0+3}}^{k+1}.
\end{align}
Let us move the term $ \mathcal{I}_4^\beta$. Then one has 
\begin{align*}
\mathcal{I}_4^\beta=&
\interleave \mathcal{\mu}(i-\ell)\interleave_{\alpha,s_0,0} \interleave \mathcal{\mu}_{k}(\ell)\interleave_{0,s+1+\beta,\gamma-\beta}\\
\lesssim& \|\rho\|_{H^{s_0+1+i-\ell}} \|\rho\|_{H^{s+4+\ell+\gamma}}\|\rho\|_{H^{s_0+3}}^{k}.
\end{align*}
As before, since $\ell\leqslant i-1$ then using the following  interpolation inequality 
$$
 \|\rho\|_{H^{s_0+1+i-\ell}} \|\rho\|_{H^{s+4+\ell+\gamma}}\lesssim  \|\rho\|_{H^{s_0+2}} \|\rho\|_{H^{s+3+i+\gamma}}.
 $$
Therefore we obtain 
\begin{align}\label{J-4}
 \mathcal{I}_4^\beta\lesssim& \|\rho\|_{H^{s+3+i+\gamma}}\|\rho\|_{H^{s_0+3}}^{k+1}.
\end{align}
As to the term $ \mathcal{I}_5^\beta$, one has 
\begin{align*}
\mathcal{I}_5^\beta=&
\interleave\mathcal{\mu}(i-\ell)\interleave_{\alpha,s_0+1,\beta}\interleave \mathcal{\mu}_{k}(\ell)\interleave_{0,s,1+\gamma-\beta}\\
\lesssim&\|\rho\|_{H^{s_0+2+i-\ell+\beta}}\|\rho\|_{H^{s+4+\ell+\gamma-\beta}}\|\rho\|_{H^{s_0+3}}^{k}.
\end{align*}
Now we use the interpolation inequality, since $\ell\leqslant i-1,$
$$
\|\rho\|_{H^{s_0+2+i-\ell+\beta}}\|\rho\|_{H^{s+4+\ell+\gamma-\beta}}\lesssim \|\rho\|_{H^{s_0+3}}\|\rho\|_{H^{s+3+i+\gamma}}.
$$
Hence
\begin{align}\label{J-5}
 \mathcal{I}_5^\beta\lesssim& \|\rho\|_{H^{s+3+i+\gamma}}\|\rho\|_{H^{s_0+3}}^{k}.
\end{align}
For the term  $ \mathcal{I}_6^\beta$, one has 
\begin{align*}\mathcal{I}_6^\beta&=
\interleave\mathcal{\mu}(i-\ell)\interleave_{\alpha,s+1,\beta}\interleave \mathcal{\mu}_{k}(\ell)\interleave_{0,s_0,1+\gamma-\beta}\\
&\lesssim  \|\rho\|_{H^{s+2+i-\ell+\beta}}\|\rho\|_{H^{s_0+4+\ell+\gamma-\beta}}\|\rho\|_{H^{s_0+3}}^{k}.
\end{align*}
 From interpolation inequality we infer
$$
\|\rho\|_{H^{s+2+i-\ell+\beta}}\|\rho\|_{H^{s_0+4+\ell+\gamma-\beta}}\lesssim \|\rho\|_{H^{s+3+i+\gamma}}\|\rho\|_{H^{s_0+3}}
$$
and therefore
\begin{align}\label{J-6}
 \mathcal{I}_6^\beta\lesssim& \|\rho\|_{H^{s+3+i+\gamma}}\|\rho\|_{H^{s_0+3}}^{k+1}.
\end{align}
Concerning the term  $ \mathcal{I}_7^\beta$, one has by interpolation
\begin{align*}\mathcal{I}_7^\beta&=
\interleave\mathcal{\mu}(i-\ell)\interleave_{\alpha,s_0+1+\beta,\gamma-\beta}\interleave \mathcal{\mu}_{k}(\ell)\interleave_{0,s,0}\\
&\lesssim \|\rho\|_{H^{s_0+2+\gamma+i-\ell}}\|\rho\|_{H^{s+3+\ell}}\|\rho\|_{H^{s_0+3}}^{k}\\
&\lesssim \|\rho\|_{H^{s_0+2}}\|\rho\|_{H^{s+3+\gamma+i}}\|\rho\|_{H^{s_0+3}}^{k}.
\end{align*}
Thus
\begin{align}\label{J-7}
 \mathcal{I}_7^\beta\lesssim& \|\rho\|_{H^{s+3+i+\gamma}}\|\rho\|_{H^{s_0+3}}^{k+1}.
\end{align}
For the last term, we write 
\begin{align*}\mathcal{I}_8^\beta&=\interleave\mathcal{\mu}(i-\ell)\interleave_{\alpha,s+1+\beta,\gamma-\beta}\interleave \mathcal{\mu}_{k}(\ell)\interleave_{0,s_0,0}\\
&\lesssim  \|\rho\|_{H^{s+2+\gamma+i-\ell}}\|\rho\|_{H^{s_0+3+\ell}}\|\rho\|_{H^{s_0+3}}^{k}.
\end{align*}
By interpolation we get
$$
\|\rho\|_{H^{s+2+\gamma+i-\ell}}\|\rho\|_{H^{s_0+3+\ell}}\lesssim \|\rho\|_{H^{s+3+i+\gamma}}\|\rho\|_{H^{s_0+2}}.
$$
It follows that
\begin{align}\label{J-8}
 \mathcal{I}_8^\beta\lesssim& \|\rho\|_{H^{s+3+i+\gamma}}\|\rho\|_{H^{s_0+3}}^{k+1}.
\end{align}
Inserting the estimate  \eqref{J-1}-\eqref{J-2}-..-\eqref{J-8} into \eqref{Big-k} we deduce for any $0\leqslant \ell\leqslant i-1$
$$
\interleave[\mathcal{\mu}(i-\ell),\mathcal{\mu}_{k}(\ell)]\interleave_{0,s,\gamma}\lesssim \|\rho\|_{H^{s+3+i+\gamma}}\|\rho\|_{H^{s_0+3}}^{k+1}.
$$
Plugging this into \eqref{Big-k0} gives
$$
\interleave\mathcal{\mu}_{k+1}(i)\interleave_{0,s,\gamma}\lesssim \|\rho\|_{H^{s+3+i+\gamma}}\|\rho\|_{H^{s_0+3}}^{k+1},
$$
which achieves the induction argument. To prove the second part of Lemma \ref{lem-com-it1} it suffices to combine the first part with Lemma \ref{Lem-Rgv1}-$($i$)$.

\end{proof}
Now we shall establish the following result based on Lemma \ref{lem-com-it1} and used before in the proof of Theorem \ref{Prop-EgorV}.
Notice that we shall use below the same notations introduced in \eqref{Matrix-y}, \eqref{Matrix-yz}, \eqref{Hamb1}, \eqref{ourta1} and  \eqref{blocmatrix}.
\begin{lemma}\label{lem-com-it2}
Let $0\leqslant j< k\leqslant  n$ and $s_0>\frac{d}{2}+1$.  Then  we have  the following assertions.
\begin{enumerate}
\item For any $i\in\{2,..,n+1\}$
$$
\left\|\left(\textnormal{Ad}_{\mathcal{M}_n}^{k-j}(\mathcal{L}_n)\, \mathcal{M}_n^{j}\mathcal{X}_n\right)(i)\right\|_{L^2(\T^{d+2})}\lesssim \|\rho\|_{H^{s_0+3}(\T^{d+1})}^{k-j}\sum_{\ell=1}^{i-1}\|\rho\|_{H^{s_0+3+i-\ell}(\T^{d+1})}\,\,\left\|\left(\mathcal{M}_n^{j}\mathcal{X}_n\right)(\ell)\right\|_{L^2(\T^{d+2})}.
$$
 The notation $Y(i)$ means the the i-th component of the vector $Y$. Notice that for $i=1$ the first component of the left hand side member  is vanishing.
 \item  For any $i\in\{1,..,n+1\}$
$$
\left\|\left(\textnormal{Ad}_{\overline{\mathcal{M}}_n}^{k-j}(\overline{\mathcal{L}}_n)\, \overline{\mathcal{M}}_n^{j}\overline{\mathcal{X}}_n\right)(i)\right\|_{L^2(\T^{d+2})}\lesssim \big(\|\rho\|_{{s_0+3}}^{q,\kappa}\big)^{k-j}\sum_{\ell=1}^{i}  \|\rho\|_{{s_0+3+i-\ell}}^{q,\kappa}\,\,\,\left\|\left(\overline{\mathcal{M}}_n^{j}\overline{\mathcal{X}}_n\right)(\ell)\right\|_{L^2(\T^{d+2})}.
$$
Here, the notation $Y(i)$ means the the i-th bloc vector of size $q+1$.
 \end{enumerate}
\end{lemma}
\begin{proof}
$(\bf{i})$ 
Using \eqref{Ck-m} and \eqref{Ma-S}, we can write 
$$\mathcal{C}_{k-j}\triangleq \textnormal{Ad}_{\mathcal{M}_n}^{k-j}(\mathcal{L}_n)=(a_{i,\ell}^{k-j})_{1\leqslant i,\ell\leqslant n+1}\in\DDD_{k-j}(n),
$$ with
\begin{align*}
a_{i,\ell}^{k-j}&=c_{i-\ell,i-1}^{k-j}, \quad\hbox{if}\quad i\geqslant \ell+k-j\\
a_{i,\ell}^{k-j}&=0, \quad\hbox{otherwise}.
\end{align*}
Denote $\mathcal{Z}_j=\mathcal{M}_n^j \mathcal{X}_n$, then we obtain for $0\leqslant j<k$
\begin{align*}
 \left(\mathcal{C}_{k-j}\mathcal{M}_n^{j}\mathcal{X}_n\right)(i)&=\sum_{\ell=1}^{i+j-k}a_{i,\ell}^{k-j}\mathcal{Z}_j(\ell)\\
 &=\sum_{\ell=1}^{i-1}a_{i,\ell}^{k-j}\mathcal{Z}_j(\ell).
\end{align*}
According to \eqref{Nq-1} and Lemma \ref{lem-com-it1} we deduce 
$$
\|a_{i,\ell}^{k-j}\mathcal{Z}_j(\ell)\|_{L^2(\T^{d+2})}\lesssim \|\rho\|_{H^{s_0+3+i-\ell}(\T^{d+1})}\,\|\rho\|_{H^{s_0+3}(\T^{d+1})}^{k-j}\,\| \mathcal{Z}_j(\ell)\|_{L^2(\T^{d+2})}.
$$
Therefore we get 
\begin{align*}
\left\| \left(\mathcal{C}_{k-j}\mathcal{M}_n^{j}\mathcal{X}_n\right)(i)\right\|_{L^2(\T^{d+2})}&\lesssim \,\|\rho\|_{H^{s_0+3}(\T^{d+1})}^{k-j}\sum_{\ell=1}^{i-1} \|\rho\|_{H^{s_0+3+i-\ell}(\T^{d+1})}\,\| \mathcal{Z}_j(\ell)\|_{L^2(\T^{d+2})},
\end{align*}
which is the desired result.

\smallskip

$(\bf{ii})$ Using \eqref{blocmatrix} and Lemma \ref{Lem-top} we can write 
$$\overline{\mathcal{C}}_{k-j}\triangleq \textnormal{Ad}_{\overline{\mathcal{M}}_n}^{k-j}(\overline{\mathcal{L}}_n)=(\overline{a}_{i,\ell}^{k-j})_{1\leqslant i,\ell\leqslant n+1}\in\DDD_{0}(n,q)
$$ with $\overline{a}_{i,\ell}^{k-j}$  a T\"oplitz matrix  belonging to $\DDD_0(q).$
Denote $\overline{\mathcal{Z}}_j=\overline{\mathcal{M}}_n^j \overline{\mathcal{X}}_n$. Hence we obtain for $0\leqslant j<k$
\begin{align*}
 \left(\overline{\mathcal{C}}_{k-j}\overline{\mathcal{M}}_n^{j}\overline{\mathcal{X}}_n\right)(i) &=\sum_{\ell=1}^{i}\overline{a}_{i,\ell}^{k-j}\overline{\mathcal{Z}}_j(\ell).
\end{align*}
The main goal is to prove the following estimate, 
\begin{equation}\label{WL90}
\|\overline{a}_{i,\ell}^{k-j}\overline{\mathcal{Z}}_j(\ell)\|_{L^2(\T^{d+2})}\lesssim  \|\rho\|_{{s_0+3+i-\ell}}^{q,\kappa}\big(\|\rho\|_{{s_0+3}}^{q,\kappa}\big)^{k-j}\,\| \overline{\mathcal{Z}}_j(\ell)\|_{L^2(\T^{d+2})}.
\end{equation}
Assume for a while this estimate, then 
\begin{align*}
\left\| \left(\overline{\mathcal{C}}_{k-j}\overline{\mathcal{M}}_n^{j}\overline{\mathcal{X}}_n\right)(i)\right\|_{L^2(\T^{d+2})}&\lesssim \,\big(\|\rho\|_{{s_0+3}}^{q,\kappa}\big)^{k-j}\sum_{\ell=1}^{i}  \|\rho\|_{{s+3+i-\ell}}^{q,\kappa}\,\| \overline{\mathcal{Z}}_j(\ell)\|_{L^2(\T^{d+2})}.
\end{align*}
Remind that we have the recursive relation 
$$
\overline{\mathcal{C}}_{k+1}=\big[\overline{\mathcal{M}}_n, \overline{\mathcal{C}}_{k}\big], 
$$
with the  initial condition 
\begin{align*}
\ \overline{\mathcal{C}}_{0}={L}_q \mathbb{I}_{n+1}=-\big(\mu(\varphi,\theta,0)+\nu(\varphi,\eta,0)\big) \textnormal{I}_{(n+1)(q+1)}.
\end{align*}
From \eqref{blocmatrix} and \eqref{blocM} one has that the law of $\overline{\mathcal{M}}_n$ is given by 
\begin{equation}\label{LL1}
{\mu}_{\overline{\mathcal{M}}_n}(i)=\left\{\begin{array}{ll}
          	\textnormal{Ad}_{\partial_\chi}^i {\overline{M}}_q-\textnormal{Ad}_{\partial_\chi}^i\overline{L}_q,&1\leqslant i\leqslant  n\\
	{\overline{M}}_q,& i=0
       \end{array}\right.
       \end{equation}
       and one can check from straightforward computations based on \eqref{LuX1} and \eqref{LUX2}  that we have the splitting 
$$
{\mu}_{\overline{\mathcal{M}}_n}(i)\triangleq \overline{\mu}(i)+ \overline{\nu}(i)
$$
where for $1\leqslant i\leqslant n$
\begin{equation}\label{LUX20}
 \overline{\mu}(i)=\begin{pmatrix}
    \overline{\mu}_{0,0}(i) &0&0&..&..&0  \\
 \overline{\mu}_{1,1}(i)&\overline{\mu}_{0,1}(i)&0&..&..&0 \\
  \overline{\mu}_{2,2}(i)&  \overline{\mu}_{1,2}(i)&\overline{\mu}_{0,2}(i)&0 &..&0\\
    ..&..&..&..&0&0\\
      ..&..&..&\overline{\mu}_{1,q-1}(i)&\overline{\mu}_{0,q-1}(i)&0\\
 \overline{\mu}_{q,q}(i)&  \overline{\mu}_{q-1,q}(i)&..&..& \overline{\mu}_{1,q}(i)&\overline{\mu}_{0,q}(i)
  \end{pmatrix}, \quad
\end{equation}
such that  for $1\leqslant k\leqslant q$
\begin{align*}
 \overline{\mu}_{k,j}(i)&= \left(_k^j\right)\textnormal{Ad}_{\partial_\chi}^i\overline{\mathbb{B}}_{\theta,(k)}\\
 &=\left(_k^j\right)\kappa^k\big[\partial_\lambda^k,\mathbb{A}_{\theta,(i)}\big]
\end{align*}
and for $k=0$, using \eqref{Tr01},
\begin{align*}
\overline{\mu}_{0,j}(i)&=-\textnormal{Ad}_{\partial_\chi}^i\widehat{\mathscr{L}}\\
&=\mathbb{A}_{\theta,(i)}.
\end{align*}
In addition for $i=0$ we have
\begin{equation}\label{LUX200}
 \overline{\mu}(0)=\begin{pmatrix}
    0&0&0&..&..&0  \\
 \overline{\mu}_{1,1}(0)&0&0&..&..&0 \\
  \overline{\mu}_{2,2}(0)&  \overline{\mu}_{1,2}(0)&0&0 &..&0\\
    ..&..&..&..&0&0\\
      ..&..&..&\overline{\mu}_{1,q-1}(0)&0&0\\
 \overline{\mu}_{q,q}(0)&  \overline{\mu}_{q-1,q}(0)&..&..& \overline{\mu}_{1,q}(0)&0
  \end{pmatrix}, \quad  \overline{\mu}_{k,j}(0)=\left(_k^j\right)\kappa^k\big[\partial_\lambda^k,\mathbb{A}_{\theta,(0)}\big].
\end{equation}
In a similar way we get 
\begin{equation}\label{LUX21}
 \overline{\nu}(i)=\begin{pmatrix}
    \overline{\nu}_{0,0}(i) &0&..&..&..&0  \\
 \overline{\nu}_{1,1}(i)&\overline{\nu}_{0,1}(i)&0&..&..&0 \\
  \overline{\nu}_{2,2}(i)&  \overline{\nu}_{1,2}&\overline{\nu}_{1,2}(i)&0 &..&0\\
    ..&..&..&..&0&0\\
      ..&..&..&\overline{\nu}_{1,q-1}(i)&\overline{\nu}_{0,q-1}(i)&0\\
 \overline{\mu}_{q,q}(i)&  \overline{\nu}_{q-1,q}(i)&..&..& \overline{\nu}_{1,q}(i)&\overline{\nu}_{0,q}(i)
  \end{pmatrix}, \quad
\end{equation}
where for $1\leqslant  i\leqslant n$
\begin{equation*}  \overline{\nu}_{k,j}(i)=
\left\{ \begin{array}{ll}
  \left(_k^j\right)\textnormal{Ad}_{\partial_\chi}^i\big(\overline{\mathbb{B}}_{\eta,(k)}-\overline{\mathbb{T}}_{k}\big),& 1\leqslant k\leqslant q\\ 
 \overline{\mathbb{B}}_{\eta,(0)}-\overline{\mathbb{T}}_{0},& k=0
  \end{array}\right.
\end{equation*}
and for $i=0$ we make the adaptation as before.
 By applying Lemma \ref{Lem-top} one finds 
\begin{equation}\label{blocmatrixM}
\overline{\mathcal{C}}_k=\begin{pmatrix}
      \overline{C}^k_{0,0} &0&..&..&..&0  \\
 \overline{C}^k_{1,1}&   \overline{C}^k_{0,1}&0&..&..&0 \\
  \overline{C}^k_{2,2}&  \overline{C}_{1,2}& \overline{C}^k_{0,2}&0 &..&0\\
  .. &..&..&..&..&0  \\
    ..&..&..&\overline{C}^k_{1,n-1}&\overline{C}^k_{0,n-1}&0\\
 \overline{C}^k_{n,n}&  \overline{C}^k_{n-1,n}&..&..& \overline{C}^k_{1,n}& \overline{C}^k_{0,n}
  \end{pmatrix}.  \end{equation}  
with
\begin{align}\label{Nq-M1}
 \overline{C}^k_{i,j}=\left(_i^j\right)\big(\overline{\mathcal{\mu}}_k(\varphi,\theta,i)+\overline{\mathcal{\nu}}_k(\varphi,\eta,i)\big), \,\,\forall\, i\leqslant j\in\{1,..,n\}
\end{align}
and the law $\big[\overline{\mathcal{M}}_n,\overline{ \mathcal{C}}_{k}\big]$ is given by 
\begin{align}\label{Moupa-01}
(\overline{\mu}+\overline{\nu})\boxast(\overline{\mu}_k+ \overline{\nu}_k)&=
\overline{\mu}\boxast\overline{\mu}_k+\overline{\nu}\boxast \overline{\nu}_k
\end{align}
where in the last identity we use the commutation relations
\begin{align}\label{Moupa-1}
(\overline{\mu}\boxast\overline{\mu}_k)(i)=0=(\overline{\nu}\boxast\overline{\nu}_k)(i).
\end{align}
Indeed, one has 
\begin{align*}
(\overline{\mu}\boxast\overline{\mu}_k)(i)&=\sum_{\ell=0}^i\left(_\ell^i\right)\big[\overline{\mu}(i-\ell),\overline{\mu}_k(\ell)\big].
\end{align*}
Since $\overline{\mu}(i-\ell),\overline{\mu}_k(\ell)\in\DDD_0(q)$ and the entries of the matrix operators $\overline{\mu}(i-\ell)$ and $\overline{\mu}_k(\ell)$ commute mutually   then from Remark \ref{Rmq-1}  we deduce  that  $\big[\overline{\mu}(i-\ell),\overline{\nu}_k(\ell)\big]=0$. This proves \eqref{Moupa-1}.
Thus we obtain from \eqref{Moupa-01} the recursive relations
\begin{align*}
\overline{\mathcal{\mu}}_{k+1}=\overline{\mu}\boxast\overline{\mu}_k&,\quad \overline{\mu}_0( i)=- {\mu}(\varphi,\theta, 0)\textnormal{I}_{q+1}\\
\overline{\mathcal{\nu}}_{k+1}=\overline{\nu}\boxast \overline{\nu}_k&,\quad \overline{\nu}_0( i)=- \nu(\varphi,\theta, 0)\textnormal{I}_{q+1}.
\end{align*}
The estimates of $\overline{\mathcal{\mu}}_{k+1}$ and $\overline{\mathcal{\nu}}_{k+1}$ are similar and  we shall restrict the discussion only to  the first one. Hence we get  for $i\in\{0,..,n\}$
\begin{align}\label{IterKL1}
\nonumber \overline{\mathcal{\mu}}_{k+1}(i)&=\big[\overline{\mathcal{\mu}}(0),\overline{\mathcal{\mu}}_{k}(i)\big]+\sum_{\ell=0}^{i-1}\left(_\ell^i\right)\big[\overline{\mathcal{\mu}}(i-\ell),\overline{\mathcal{\mu}}_{k}(\ell)\big]\\
&=\textnormal{Ad}_{\overline{\mathcal{\mu}}(0)}\,\overline{\mathcal{\mu}}_{k}(i)+\sum_{\ell=0}^{i-1}\left(_\ell^i\right)\big[\overline{\mathcal{\mu}}(i-\ell),\overline{\mathcal{\mu}}_{k}(\ell)\big].
\end{align}
Iterating this yields 
\begin{align*}
\overline{\mathcal{\mu}}_{k+q+1}(i)&=\textnormal{Ad}^{q+1}_{\overline{\mathcal{\mu}}(0)}\,\overline{\mathcal{\mu}}_{k}(i)+\sum_{\ell=0}^{i-1}\left(_\ell^i\right)\sum_{j=0}^{q}\textnormal{Ad}^{q-j}_{\overline{\mathcal{\mu}}(0)}\big[\overline{\mathcal{\mu}}(i-\ell),\overline{\mathcal{\mu}}_{k+j}(\ell)\big]
\end{align*}
On the other hand, from \eqref{LUX200} one gets $\overline{\mu}(0)\in\DDD_1(q)$. Therefore applying Lemma \ref{Lem-top} we deduce that
$$
\textnormal{Ad}^{m}_{\overline{\mathcal{\mu}}(0)}\,\overline{\mathcal{\mu}}_{k}(i)\in\DDD_{\min\{m,q+1\}}(q).
$$
In particular we find that $\textnormal{Ad}^{q+1}_{\overline{\mathcal{\mu}}(0)}\,\overline{\mathcal{\mu}}_{k}(i)=0$ and consequently
\begin{align}\label{Houma1V}
\overline{\mathcal{\mu}}_{k+q+1}(i)&=\sum_{\ell=0}^{i-1}\left(_\ell^i\right)\sum_{j=0}^{q}\textnormal{Ad}^{q-j}_{\overline{\mathcal{\mu}}(0)}\big[\overline{\mathcal{\mu}}(i-\ell),\overline{\mathcal{\mu}}_{k+j}(\ell)\big].
\end{align}
We point out that from \eqref{IterKL1}, the matrix operator $\overline{\mu}_k(i)$ is a Topelitz matrix in $\DDD_0(q)$ and therefore the entries of the commutators matrices   $\big[\overline{\mathcal{\mu}}(i-\ell),\overline{\mathcal{\mu}}_{k}(\ell)\big]$  are linear combination of scalar commutators and then we can  use the estimates of Lemma \ref{comm-pseudo}. Thus straightforward computations   with slight adaptations  of the proof of Lemma \ref{lem-com-it1}  allow to check  by  induction, using in particular  the identity \eqref{Houma1V},  that for any  $k\in\N, j\in\{0,1,..,q\},\gamma\in\N,\, s\geqslant s_0$
$$
\interleave\overline{\mathcal{\mu}}_{k+j}(i)\interleave_{0,s,\gamma}\lesssim \|\rho\|_{{s_0+3+i+j+\gamma}}^{q,\kappa}\big(\|\rho\|_{{s_0+3}}^{q,\kappa}\big)^{k+j},
$$
which gives the estimate \eqref{WL90} in view of Lemma \ref{Lem-Rgv1}-$($i$)$.
\end{proof}
 \subsubsection{Iterated kernel estimates}\label{Sect-It-K-Es}
 This section is devoted to some estimates on the  iterated kernels used before during the proof of Theorem \ref{Prop-EgorV}.
 Let us first introduce  the problem that we want to solve. Given a sequence of smooth functions $\rho_n:\mathcal{O}\times\T^{d+1}\to\RR$ and define the  pseudo-differential operator of order $\alpha\in(0,\frac12)$
 $$
 \mathbb{A}_{\theta,n}=\partial_\theta\big(\rho_n(\varphi,\theta)|\textnormal{D}|^{\alpha-1}+|\textnormal{D}|^{\alpha-1}\rho_n(\varphi,\theta)\big).
 $$
 Let $(\varphi,\theta,\eta)\mapsto K_0(\varphi,\theta,\eta)$ be a  function with some prescribed regularity, essentially smooth in the variable $\varphi$ but slightly singular in the diagonal $\theta=\eta$,  and consider the iterative kernels 
 $$
 K_{n+1}=\left(\mathbb{A}_{\theta,n}+\mathbb{A}_{\eta,n}-\mathbb{S}_{n}\right) K_{n}, n\in\N,
 $$
 with 
 $$
\mathbb{S}_{n}\triangleq (\partial_\eta\rho_n)|\textnormal{D}_\eta|^{\alpha-1}+|\textnormal{D}_\eta|^{\alpha-1}(\partial_\eta\rho_n).
$$
 We want to find suitable  recursive estimates for the sequence $(K_n)_n$ that allow to deduce later  tame estimates. At this stage, it is worthy to point out that, because  the initial kernel $K_0$ may be singular at the diagonal  $\theta=\eta$, we cannot deal separately with the operators  $\mathbb{A}_{\theta,n}$ and $\mathbb{A}_{\eta,n}$ and one should keep together their collective behavior in order to  cancel the diagonal singularity. This can be easily understood with  the following simple toy model. Take $K_0(\varphi,\theta,\eta)=\Psi(\theta-\eta)$ and $\Psi$ is singular and let  $\mathbb{A}_{\theta}=\partial_\theta$. Then both new kernels $\mathbb{A}_{\theta}K_0 $ and $\mathbb{A}_{\eta}K_0$ are more singular than $K_0$ and their iterated kernels $\mathbb{A}_{\theta}^nK_0 $ and $\mathbb{A}_{\eta}^nK_0 $ become more and more singular. However, it is clear that $(\mathbb{A}_{\theta}+\mathbb{A}_{\eta})K_0=0$ which means that the singularity is canceled due to the mutual actions of the involved operators. This example is very particular because we are dealing with differential operators which act point-wisely, which is not the case with our nonlocal pseudo-differential operators.  Then to capture this cancellation, we find convenient to move to the symbol representation, find  suitable recursive estimates and later come back to the kernel. According to the formula \eqref{symb-kern}, the symbol associated to $K_n$ and denoted here by $a_n$ can be recovered     through the formula
 \begin{align}\label{Symb-Rec1}
 a_n(\varphi,\theta,\xi)=\int_\T K_n(\varphi,\theta,\theta+\eta)e^{\ii \, \eta\xi}d\eta,
\end{align}
 with the inverse Fourier formula
\begin{align}\label{KKernel-S34}
 K_n(\varphi,\theta,\eta)=\frac{1}{2\pi}\sum_{\xi\in\Z}  a_n(\varphi,\theta,\xi)e^{\ii \, (\theta-\eta)\xi}.
 \end{align}
 Now we shall use the following splitting 
\begin{align}\label{Formul1X}
 K_{n+1}=\left(\mathcal{A}_{\theta,n}+\mathcal{A}_{\eta,n}\right) K_{n}+\mathcal{R}_n K_{n},\quad  n\in\N
\end{align}
 with
 $$
 \mathcal{A}_{\theta,n}=2\rho_n(\varphi,\theta)\partial_\theta|\textnormal{D}_\theta|^{\alpha-1}
 $$
 and
\begin{align}\label{Remaind1}
\nonumber \mathcal{R}_n=&2(\partial_\theta\rho_n)|\textnormal{D}_\theta|^{\alpha-1}+(\partial_\eta\rho_n)|\textnormal{D}_\eta|^{\alpha-1}\\
\nonumber &\quad+\partial_\theta\big[|\textnormal{D}_\theta|^{\alpha-1},\rho_n\big]+\partial_\eta\big[|\textnormal{D}_\eta|^{\alpha-1},\rho_n\big]-|\textnormal{D}_\eta|^{\alpha-1}(\partial_\eta\rho_n)\\
 &\qquad\triangleq \mathcal{R}_{n,1}+\mathcal{R}_{n,2}+\mathcal{R}_{n,3}+\mathcal{R}_{n,4}+\mathcal{R}_{n,5}.
\end{align}
For the sake of simplicity we  shall make in the statement below the following  convention: if $\sigma_\mathscr{A}$ is the symbol of an operator $\mathscr{A}$ as in \eqref{kernel-phi-lambda} then we denote by $\interleave \sigma_ \mathcal{A}\interleave_{m,s,\gamma}^{q,\kappa}$ the same norm $\interleave \mathcal{A}\interleave_{m,s,\gamma}^{q,\kappa}$ defined in \eqref{Def-pseud-w}.
Our  first main result reads as follows.
\begin{lemma}\label{lemm-iter1}
 Let $q\in\N$, $ m\in\RR, \gamma\in\N$ and $  s\geqslant  s_0>\frac{d+1}{2}{+q},$ then the following assertions  hold true.
\begin{enumerate}
\item For any $n\in\N$ we have 
\begin{align*}
\nonumber \interleave a_{n+1}\interleave_{m,s,\gamma}^{q,\kappa}\lesssim&\|\rho_n\|_{{s_0+3}}^{q,\kappa}\interleave a_{n}\interleave_{m,s+1,\gamma}^{q,\kappa}+\|\rho_n\|_{{s+3}}^{q,\kappa}\interleave a_{n}\interleave_{m,s_0+1,\gamma}^{q,\kappa}\\
\nonumber&+\|\rho_n\|_{{s_0+3}}^{q,\kappa}\interleave a_n\interleave_{m,s,\gamma+1}^{q,\kappa}+\|\rho_n\|_{{s+3}}^{q,\kappa}\interleave a_n\interleave_{m,s_0,\gamma+1}^{q,\kappa} \\
\nonumber&\quad+\sum_{0\leqslant\beta\leqslant\gamma}\|\rho_n\|_{{s_0+3+|m|+\gamma-\beta}}^{q,\kappa} \interleave a_n\interleave_{m,s,\beta}^{q,\kappa}+\|\rho_n\|_{{s_0+3}}^{q,\kappa}\interleave a_n\interleave_{m,s+\gamma-\beta,\beta}^{q,\kappa}\\
\nonumber&\qquad+
\sum_{0\leqslant\beta\leqslant\gamma}\|\rho_n\|_{{s+3+|m|+\gamma-\beta}}^{q,\kappa} \interleave a_n\interleave_{m,s_0,\beta}^{q,\kappa}+\|\rho_n\|_{{s+3}}^{q,\kappa}\interleave a_n\interleave_{m,s_0+\gamma-\beta,\beta}^{q,\kappa}.
\end{align*}
\item Let $m\in\RR$, $s\geqslant s_0, \gamma,n\in\N$ and assume the existence of   an   increasing $\log$-convex function $F$ such that for   
\begin{align}\label{Assu1}
\forall\, s'\geqslant s_0,\,\forall \,\gamma'\in\N\quad\hbox{with}\quad   s'+\gamma'\leqslant s+\gamma+n\Longrightarrow\interleave a_{0}\interleave_{m,s^\prime,\gamma^\prime}^{q,\kappa}\leqslant F(s^\prime+\gamma^\prime).
\end{align}
Then,
\begin{align*}
 \interleave a_{n}\interleave_{m,s,\gamma}^{q,\kappa}&\lesssim F\big(s_0+3\big)\sum_{i=0}^{n-1} \|\rho_i\|_{{s+n+4+|m|+\gamma}}^{q,\kappa}\prod_{j=0\atop  j\neq i}^{n-1} \|\rho_j\|_{{s_0+3+|m|}}^{q,\kappa}\\
 &\quad+F\big(s+\gamma+n+4\big) \prod_{j=0}^{n-1} \|\rho_j\|_{{s_0+3+|m|}}^{q,\kappa},
\end{align*}
with the convention $\displaystyle{\sum_{i=0}^{-1}..=0}$\quad and \quad $\displaystyle{\prod_{j=0}^{-1}=1}$.

\item Let $m<-\frac12$, $s\geqslant s_0$ and $n\in\N$. Assume that \eqref{Assu1} occurs, then    
\begin{align*}
\left(\int_{\T}\left(\|K_n(\cdot,\centerdot,\centerdot+\eta)\|_{s}^{q,\kappa}\right)^2d\eta\right)^{\frac12}&\lesssim F\big(s_0+3\big)\sum_{i=0}^{n-1} \|\rho_i\|_{{s+n+4+|m|}}^{q,\kappa}\prod_{j=0\atop  j\neq i}^{n-1} \|\rho_j\|_{{s_0+3+|m|}}^{q,\kappa}\\
 &\quad+F\big(s+n+4\big) \prod_{j=0}^{n-1} \|\rho_j\|_{{s_0+3+|m|}}^{q,\kappa}.
\end{align*}

\end{enumerate}
\end{lemma}
\begin{proof}
$\bf{(i)}$ 
First, by  the integral representation \eqref{fract1}  one gets by straightforward computations 
$$
\partial_\theta|\textnormal{D}_\theta|^{\alpha-1}h(\varphi,\theta)=\textnormal{p.v.}\int_{\T}\mathscr{K}_0(\theta^\prime)h(\varphi,\theta-\theta^\prime)d\theta^\prime,\quad \mathscr{K}_0(\theta)=\frac{-\alpha\,\mathrm{cotan}(\frac{\theta}{2})}{4\pi|\sin(\frac\theta2)|^\alpha}\cdot
$$
Then using in particular \eqref{KKernel-S34} we find
\begin{align*}
\mathcal{A}_{\theta,n} K_n(\varphi,\theta,\eta)&=2\rho_n(\varphi,\theta)\int_{\T}  \mathscr{K}_0(\theta^\prime)K_n(\varphi,\theta-\theta',\eta)d\theta'
\\
&=\frac{1}{\pi}\rho_n(\varphi,\theta)\sum_{\xi^\prime\in\Z}\int_{\T}  \mathscr{K}_0(\theta^\prime)a_n(\varphi,\theta-\theta^\prime,\xi^\prime)e^{\ii \, (\theta-\theta^\prime-\eta)\xi^\prime}d\theta^\prime.
\end{align*}
Similarly, we get 
\begin{align}\label{Ahl-Bar1}
\nonumber\mathcal{A}_{\eta,n} K_n(\varphi,\theta,\eta)&=\frac1\pi\rho_n(\varphi,\eta)\sum_{\xi^\prime\in\Z}\int_{\T}  \mathscr{K}_0(\theta^\prime)a_n(\varphi,\theta,\xi^\prime)e^{\ii \, (\theta+\theta^\prime-\eta)\xi^\prime}d\theta^\prime
\\&=-\frac1\pi\rho_n(\varphi,\eta)\sum_{\xi^\prime\in\Z}\int_{\T} a_n(\varphi,\theta,\xi^\prime) \mathscr{K}_0(\theta^\prime)e^{\ii \, (\theta-\theta^\prime-\eta)\xi^\prime}d\theta^\prime,
\end{align}
where we have used in the last line the fact that $\mathscr{K}_0$ is odd.
Therefore the  symbol associated to $ \mathcal{A}_{\theta,n} K_n$ and  given by \eqref{Symb-Rec1} takes the form 
\begin{align*}
&\frac1\pi\sum_{\xi^\prime \in\Z}\rho_n(\varphi,\theta)\iint_{\T^2}  \mathscr{K}_0(\theta^\prime)a_n(\varphi,\theta-\theta^\prime,\xi^\prime)e^{\ii \,(\xi-\xi^\prime)\eta}e^{-\ii  \,\theta^\prime\xi^\prime}d\theta^\prime d\eta 
\end{align*}
and the one  of $\mathcal{A}_{\eta,n} K_n$ is given by
$$
-\frac1\pi\sum_{\xi^\prime \in\Z}a_n(\varphi,\theta,\xi^\prime)\iint_{\T^2}  \mathscr{K}_0(\theta^\prime)\rho_n(\varphi,\theta+\eta)e^{\ii \,(\xi-\xi^\prime)\eta}e^{-\ii  \,\theta^\prime\xi^\prime}d\theta^\prime d\eta. 
$$
Hence the symbol associated to  $ \mathcal{A}_{\theta,n} K_n+\mathcal{A}_{\eta,n} K_n$ takes the form
\begin{align*}
&\frac1\pi\sum_{\xi^\prime \in\Z}\rho_n(\varphi,\theta)\iint_{\T^2}  \mathscr{K}_0(\theta^\prime)\big[a_n(\varphi,\theta-\theta^\prime,\xi^\prime)-a_n(\varphi,\theta,\xi^\prime)\big]e^{\ii \,(\xi-\xi^\prime)\eta}e^{-\ii  \,\theta^\prime\xi^\prime}d\theta^\prime d\eta \\
&\quad -\frac1\pi\sum_{\xi^\prime \in\Z}\lambda(\xi^\prime)a_n(\varphi,\theta,\xi^\prime)\int_{\T}  \big[\rho_n(\varphi,\theta+\eta)-\rho_n(\varphi,\theta)\big]e^{\ii \,(\xi-\xi^\prime)\eta} d\eta\\
&\qquad \triangleq \sigma_{n,1}(\varphi,\theta,\xi)+\sigma_{n,2}(\varphi,\theta,\xi) \end{align*}
with
\begin{align}\label{Eq-bv1}
\lambda(\xi)\triangleq \int_{\T}  \mathscr{K}_0(\theta^\prime)e^{-\ii  \,\theta^\prime\xi}d\theta^\prime.
\end{align}
One may easily write by integration in $\eta$
\begin{align*}
\sigma_{n,1}(\varphi,\theta,\xi)=&2\rho_n(\varphi,\theta)\int_{\T}  \mathscr{K}_0(\theta^\prime)\big[a_n(\varphi,\theta-\theta^\prime,\xi)-a_n(\varphi,\theta,\xi)\big]e^{-\ii  \,\theta^\prime\xi}d\theta^\prime 
 \end{align*}
 and
\begin{align*}
\sigma_{n,2}(\varphi,\theta,\xi)=&\frac1\pi\bigintsss_{\T} \frac{\rho_n(\varphi,\theta+\eta)-\rho_n(\varphi,\theta)}{e^{\ii \eta}-1}\sum_{\xi^\prime \in\Z}\lambda(\xi+\xi^\prime)a_n(\varphi,\theta,\xi+\xi^\prime) \big(1-e^{\ii \eta}\big)e^{-\ii  \,\xi^\prime\eta} d\eta. \end{align*}
Using the identities \eqref{It-action} and \eqref{Sum-parts}
$$
\big(1-e^{\ii \eta}\big)e^{-\ii  \,\xi^\prime\eta}=\overline{\Delta}_{\xi^\prime}e^{-\ii  \,\xi^\prime\eta}\quad\hbox{and}\quad \sum_{\xi^\prime\in\Z}\overline{\Delta}_{\xi^\prime}f(\xi^\prime)g(\xi^\prime)=-\sum_{\xi^\prime\in\Z}f(\xi^\prime){\Delta}_{\xi^\prime}g(\xi^\prime)
$$
we deduce that 
\begin{align}\label{sigma-n2-Ta}
\sigma_{n,2}(\varphi,\theta,\xi)=&\frac1\pi\bigintsss_{\T} \frac{\rho_n(\varphi,\theta+\eta)-\rho_n(\varphi,\theta)}{1-e^{\ii \eta}}\sum_{\xi^\prime \in\Z}e^{\ii \,(\xi-\xi^\prime)\eta}{\Delta}_{\xi^\prime}\big[\lambda(\xi^\prime)a_n(\varphi,\theta,\xi^\prime)\big]  d\eta.
 \end{align}
According to the relation \eqref{Formul1X} we get that the symbol $a_n$  satisfies the recursive relation
\begin{equation}\label{sigma10}
a_{n+1}(\varphi,\theta,\xi)=\sigma_{n,1}(\varphi,\theta,\xi)+\sigma_{n,2}(\varphi,\theta,\xi)+\mu_n
\end{equation}
where $\mu_n$ is the symbol associated to $\mathcal{R}_nK_n$.
Let us start with estimating  $\sigma_{n,1}$. For $\gamma\in\N$, we use the Leibniz rule \eqref{Leibn-disc}
\begin{align*}
\Delta_\xi^\gamma\sigma_{n,1}(\varphi,\theta,\xi)=2\rho_n(\varphi,\theta)\sum_{0\leqslant\beta\leqslant\gamma}\left(_\beta^\gamma \right)&\int_{\T}  \mathscr{K}_0(\theta^\prime)\Delta_\xi^\beta\big[a_n(\varphi,\theta-\theta^\prime,\xi)-a_n(\varphi,\theta,\xi)\big]\\
&\qquad \Delta_\xi^{\gamma-\beta}[e^{-\ii  \,\theta^\prime\xi}](\xi+\beta)d\theta^\prime.
 \end{align*}
Applying \eqref{It-action} gives for large $\xi$ (such that $\xi+\beta\neq 0$)
\begin{align}\label{Formg1}
 \Delta_\xi^{\gamma-\beta}[e^{-\ii  \,\theta^\prime\xi}](\xi+\beta)&=\left(e^{-\ii  \theta^\prime}-1\right)^{\gamma-\beta}e^{-\ii  \,\theta^\prime(\xi+\beta)}\\
\nonumber &= i^{\gamma-\beta}\left(e^{-\ii  \theta^\prime}-1\right)^{\gamma-\beta}(\xi+\beta)^{\beta-\gamma}\partial_{\theta^\prime}^{\gamma-\beta}e^{-\ii  \,\theta^\prime(\xi+\beta)}.
 \end{align}
Hence,  integration by parts allows to get for large $\xi$
  \begin{align*}
\langle\xi\rangle^{-m+\gamma}\Delta_\xi^\gamma\sigma_{n,1}(\varphi,\theta,\xi)=&2\rho_n(\varphi,\theta)\sum_{0\leqslant\beta\leqslant\gamma\atop0\leqslant\varrho\leqslant\gamma-\beta}(-i)^{\alpha-\beta}\left(_\beta^\gamma\right)\left(_\varrho^{\gamma-\beta}\right)\int_{\T}  \partial_{\theta^\prime}^{\gamma-\beta-\varrho}\Big[(1-e^{-\ii  \theta^\prime})^{\gamma-\beta}\mathscr{K}_0(\theta^\prime)\Big]\\
&\frac{\langle\xi\rangle^{-m+\gamma}}{(\xi+\beta)^{\gamma-\beta}}\partial_{\theta^\prime}^{\varrho}\Big[{\Delta_\xi^\beta  a_n(\varphi,\theta-\theta^\prime,\xi)-\Delta_\xi^\beta a_n(\varphi,\theta,\xi)}\Big]e^{-\ii  \,\theta^\prime(\xi+\beta)}d\theta^\prime \\
&\qquad \triangleq 2 \rho_n(\varphi,\theta)\sum_{0\leqslant\beta\leqslant\gamma\atop0\leqslant\varrho\leqslant\gamma-\beta}\left(_\beta^\gamma\right)\left(_\varrho^{\gamma-\beta}\right)\mathcal{T}_n^{\beta,\varrho}(\varphi,\theta,\xi).
 \end{align*}
From straightforward computations we get for any $\epsilon>0$
 $$
 \left|\partial_{\theta^\prime}^{\gamma-\beta-\varrho}\Big[(1-e^{-\ii  \theta^\prime})^{\gamma-\beta}\partial_\alpha^j\mathscr{K}_0(\theta^\prime)\Big]\right|\lesssim |\sin(\theta^\prime/2)|^{\varrho-1-\alpha-\epsilon}.
 $$
 Hence we get from Sobolev embeddings  that for $\gamma-\beta\geqslant\varrho\geq 1$ 
 \begin{align*}
\sup_{\alpha\in(0,\overline\alpha),\atop 0\leqslant j\leqslant q}\kappa^j\left\|\partial_\alpha^j\mathcal{T}_n^{\beta,\varrho}(\cdot,\centerdot,\xi)\right\|_{H^{s-j}}\lesssim&\sup_{\alpha\in(0,\overline\alpha),\atop 0\leqslant j\leqslant q}\int_{\T}  |\sin(\theta^\prime/2)|^{-\alpha-\epsilon}\langle\xi\rangle^{-m+\beta}\kappa^j\left\|\Delta_\xi^\beta \partial_{\alpha}^j a_n(\cdot,\centerdot,\xi)\right\|_{H^{s+\varrho-j}}d\theta^\prime \\
&\lesssim \interleave a_n\interleave_{m,s+\varrho,\beta}^{q,\kappa}\\
&\quad\lesssim \interleave a_n\interleave_{m,s+\gamma-\beta,\beta}^{q,\kappa}.
 \end{align*}
  For $\varrho=0$ we proceed as follows using in particular Lemma \ref{lem-Reg1},
 \begin{align*}
\sup_{\alpha\in(0,\overline\alpha),\atop 0\leqslant j\leqslant q}\kappa^j\left\|\partial_\alpha^q\mathcal{T}_n^{\beta,0}(\cdot,\centerdot,\xi)\right\|_{H^{s-j}}\lesssim&\sup_{\alpha\in(0,\overline\alpha),\atop 0\leqslant j\leqslant q}\bigintsss_{\T}  |\sin(\theta^\prime/2)|^{-\alpha-\epsilon}\langle\xi\rangle^{-m+\beta}\\
&\times \kappa^j\left\|\frac{\Delta_\xi^\beta  \partial_\alpha^ja_n(\cdot,\centerdot-\theta^\prime,\xi)-\Delta_\xi^\beta \partial_\alpha^ja_n(\cdot,\centerdot,\xi)}{1-e^{-\ii  \theta^\prime}}\right\|_{H^{s-j}}d\theta^\prime \\
&\quad\lesssim \interleave a_n\interleave_{m,s+1,\beta}^{q,\kappa}\\
&\qquad \lesssim \interleave a_n\interleave_{m,s+1,\gamma}^{q,\kappa}.
 \end{align*}
Notice that the treatment  of the quantities related to the differentiation with respect to the parameter $\omega$ instead of $\alpha$ can be done in a similar way and we get the same kind of estimates. Therefore, performing the law products in Lemma \ref{Law-prodX1}  with the preceding estimates  yields
 \begin{align}\label{sigma11}
\nonumber\interleave\sigma_{n,1}\interleave_{m,s,\gamma}^{q,\kappa}&\lesssim\|\rho_n\|_{{s_0}}^{q,\kappa}\interleave a_{n}\interleave_{m,s+1,\gamma}^{q,\kappa}+\|\rho_n\|_{{s}}^{q,\kappa}\interleave a_{n}\interleave_{m,s_0+1,\gamma}^{q,\kappa}\\
&+\sum_{0\leqslant\beta\leqslant\gamma}\big(\|\rho_n\|_{{s_0}}^{q,\kappa}\interleave a_n\interleave_{m,s+\gamma-\beta,\beta}^{q,\kappa}+\|\rho_n\|_{{s}}^{q,\kappa}\interleave a_n\interleave_{m,s_0+\gamma-\beta,\beta}^{q,\kappa}\big).
 \end{align}
 Next, we shall  move to the estimate of the term  $\sigma_{n,2}$ described in \eqref{sigma-n2-Ta}. Define
 \begin{align}\label{ddeff}
\zeta_n(\varphi,\theta,\xi)=&\frac{1}{\pi}\bigintsss_{\T} \frac{\rho_n(\varphi,\theta+\eta)-\rho_n(\varphi,\theta)}{1-e^{\ii \eta}} e^{\ii \,\xi\eta} d\eta,
 \end{align}
then it is clear that
\begin{align*}
\sigma_{n,2}(\varphi,\theta,\xi)=&\sum_{\xi^\prime \in\Z}\zeta_n(\varphi,\theta,\xi-\xi^\prime){\Delta}_{\xi^\prime}\big[\lambda(\xi^\prime)a_n(\varphi,\theta,\xi^\prime)\big].
 \end{align*}
From the convolution structure, we infer
\begin{align*}
\Delta_\xi^\gamma\sigma_{n,2}(\varphi,\theta,\xi)=&\sum_{\xi^\prime \in\Z}\zeta_n(\varphi,\theta,\xi-\xi^\prime){\Delta}_{\xi^\prime}^{1+\gamma}\big[\lambda(\xi^\prime)a_n(\varphi,\theta,\xi^\prime)\big].
 \end{align*}
On the other hand, using   \eqref{Eq-bv1}, \eqref{fract1} and \eqref{Fo-c}, we may write
\begin{align*}
\lambda(\xi)&=\big(\partial_\theta|\textnormal D|^{\alpha-1}e^{-\ii  \,\theta\xi}\big)_{|\theta=0}\\
&=-\ii\frac{2^{\alpha}\Gamma(1-\alpha)}{\Gamma(\frac{\alpha}{2})\Gamma(1-\frac{\alpha}{2})}\frac{\xi\, \Gamma\big(|\xi|+\frac{\alpha}{2}\big)}{\Gamma\big(|\xi|+1-\frac{\alpha}{2}\big)}\cdot
\end{align*} 
From the  expansion of Gamma function, similarly  to \eqref{Frac-Mult} and based on Lemma \ref{Stirling-formula}-{\rm (ii)}, one finds by straightforward computations that for any $j,\beta\in\NN, \epsilon>0$
\begin{align}\label{Eq-bv2}
\sup_{\alpha\in(0,\overline\alpha)}\big|\Delta_\xi^\beta\partial_\alpha^j \lambda(\xi)\big|\lesssim \langle \xi\rangle^{\overline\alpha+\epsilon-\beta}.
\end{align}
Consequently, applying Leibniz formula \eqref{Leibn-disc} leads to
$$
\left\|{\Delta}_{\xi^\prime}^{1+\gamma}\big[\lambda(\xi^\prime)a_n(\cdot,\centerdot,\xi^\prime)\big]\right\|_{s}^{q,\kappa}\lesssim \langle\xi^\prime\rangle^{m+\overline\alpha+\epsilon-1-\gamma}\interleave a_n\interleave_{m,s,1+\gamma}^{q,\kappa}
$$
and 
$$
\left\|{\Delta}_{\xi^\prime}^{1+\gamma}\big[\lambda(\xi^\prime)a_n(\cdot,\centerdot,\xi^\prime)\big]\right\|_{s}^{q,\kappa}\lesssim \langle\xi^\prime\rangle^{m+\overline\alpha+\epsilon}\interleave a_n\interleave_{m,s,0}^{q,\kappa}.
$$
Integration by parts  in \eqref{ddeff}, using Lemma \ref{lem-Reg1}, implies for any $N\in\N$
$$
 \langle\xi\rangle^N\left\|\zeta_n(\cdot,\centerdot,\xi)\right\|_{s}^{q,\kappa}\lesssim\|\rho_n\|_{{s+1+N}}^{q,\kappa}.
$$
Putting together the preceding estimates with  the law products of Lemma \ref{Law-prodX1} we find  for any $N_1,N_2\in\N$
\begin{align}\label{sum-biz-125}
\nonumber\left\|\Delta_\xi^\gamma\sigma_{n,2}(\cdot,\centerdot,\xi)\right\|_{s}^{q,\kappa}\lesssim&\|\rho_n\|_{{s_0+1+N_1}}^{q,\kappa}\interleave a_n\interleave_{m,s,0}^{q,\kappa}\sum_{|\xi^\prime|\leqslant\frac12|\xi| \atop|\xi^\prime|\geqslant2|\xi|}
\langle\xi-\xi^\prime\rangle^{-N_1} \langle\xi^\prime\rangle^{m+\overline\alpha+\epsilon}\\
\nonumber+&\|\rho_n\|_{{s_0+1+N_2}}^{q,\kappa}\interleave a_n\interleave_{m,s,1+\gamma}^{q,\kappa}\sum_{\frac12|\xi|\leqslant|\xi^\prime|\leqslant 2|\xi|}\langle\xi-\xi^\prime\rangle^{-N_2} \langle\xi^\prime\rangle^{m+\overline\alpha+\epsilon-1-\gamma}\\
\nonumber+&\|\rho_n\|_{{s+1+N_1}}^{q,\kappa}\interleave a_n\interleave_{m,s_0,0}^{q,\kappa}\sum_{|\xi^\prime|\leqslant\frac12|\xi| \atop|\xi^\prime|\geqslant2|\xi|}
\langle\xi-\xi^\prime\rangle^{-N_1} \langle\xi^\prime\rangle^{m+\overline\alpha+\epsilon}\\
+&\|\rho_n\|_{{s+1+N_2}}^{q,\kappa}\interleave a_n\interleave_{m,s_0,1+\gamma}^{q,\kappa}\sum_{\frac12|\xi|\leqslant|\xi^\prime|\leqslant 2|\xi|}\langle\xi-\xi^\prime\rangle^{-N_2} \langle\xi^\prime\rangle^{m+\overline\alpha+\epsilon-1-\gamma}.
 \end{align}
Since $\overline\alpha\in[0,1)$ then one can check that with the choice $N_1=\gamma+2+|m|$ and $N_2=2$,  then 
$$
\sup_{\xi\in\Z}\langle\xi\rangle^{-m+\gamma}\sum_{|\xi^\prime|\leqslant\frac12|\xi| \atop|\xi^\prime|\geq2|\xi|}
\langle\xi-\xi^\prime\rangle^{-N_1} \langle\xi^\prime\rangle^{m+\overline\alpha+\epsilon}<\infty
$$
and
$$
\sup_{\xi\in\Z}\langle\xi\rangle^{-m+\gamma}\sum_{\frac12|\xi|\le|\xi^\prime|\le 2|\xi|}\langle\xi-\xi^\prime\rangle^{-N_2} \langle\xi^\prime\rangle^{m+\overline\alpha+\epsilon-1-\gamma}<\infty.
$$
Hence
\begin{align}\label{sigma12}
\nonumber\interleave\sigma_{n,2}\interleave_{m,s,\gamma}^{q,\kappa}\lesssim&\|\rho_n\|_{{s_0+3+\gamma+|m|}}^{q,\kappa}\interleave a_n\interleave_{m,s,0}^{q,\kappa}+\|\rho_n\|_{{s_0+3}}^{q,\kappa}\interleave a_n\interleave_{m,s,\gamma+1}^{q,\kappa}\\
+&\|\rho_n\|_{{s+3+\gamma+|m|}}^{q,\kappa}\interleave a_n\interleave_{m,s_0,0}+\|\rho_n\|_{{s+3}}^{q,\kappa}\interleave a_n\interleave_{m,s_0,\gamma+1}^{q,\kappa}. \end{align}
Next, we shall  move to the estimate of the symbol  associated to  $\mathcal{R}_{n,1}K_n$ introduced in \eqref{Remaind1}. This symbol will be denoted  by $\mu_{n,1}$ and can be recovered from the formula \eqref{Symb-Rec1}.  By direct computations based on \eqref{fract1} and \eqref{KKernel-S34} we can check that
\begin{align*}
\mu_{n,1}(\varphi,\theta,\xi)&=2\partial_\theta \rho_n(\varphi,\theta)\int_{\T}  \frac{a_n(\varphi,\theta-\theta^\prime,\xi)}{|\sin(\theta^\prime/2)|^\alpha}e^{-\ii  \,\theta^\prime\xi}d\theta^\prime\\ 
&\triangleq 2\partial_\theta \rho_n(\varphi,\theta)\mathcal{T}_n(\varphi,\theta,\xi).
 \end{align*}
Then  proceeding  as for the term $\sigma_{n,1}$ one gets
\begin{align*}
\langle\xi\rangle^{-m+\gamma}\left\|\Delta_\xi^\gamma\mathcal{T}_n(\cdot,\centerdot,\xi)\right\|_{s}^{q,\kappa}&\lesssim \sum_{0\leqslant\beta\leqslant\gamma}\interleave a_n\interleave_{m,s+\gamma-\beta,\beta}^{q,\kappa} \end{align*}
 and by the law products we infer 
 \begin{align}\label{sigma13}
\interleave\mu_{n,1}\interleave_{m,s,\gamma}^{q,\kappa}\lesssim&\sum_{0\leqslant\beta\leqslant\gamma}\big(\|\rho_n\|_{{s_0+1}}^{q,\kappa}\interleave a_n\interleave_{m,s+\gamma-\beta,\beta}^{q,\kappa}+\|\rho_n\|_{{s+1}}^{q,\kappa}\interleave a_n\interleave_{m,s_0+\gamma-\beta,\beta}^{q,\kappa}\big).
 \end{align}
The symbol $\mu_{n,2}$ associated to  $\mathcal{R}_{n,2}K_n$ defined in \eqref{Remaind1} can be recovered from $a_n$ in a similar way to \eqref{Ahl-Bar1} and one gets
\begin{align}\label{conv-2021}
\nonumber \mu_{n,2}(\varphi,\theta,\xi)=&\frac{1}{2\pi}\bigintsss_{\T} {\partial_\eta\rho_n(\varphi,\theta+\eta)}\sum_{\xi^\prime \in\Z}\lambda_1(\xi^\prime)a_n(\varphi,\theta,\xi^\prime) e^{\ii \,(\xi-\xi^\prime)\eta} d\eta\\
&=\sum_{\xi^\prime \in\Z}\lambda_1(\xi^\prime)a_n(\varphi,\theta,\xi^\prime) \zeta_{n,1}(\varphi,\theta,\xi-\xi^\prime),
\end{align}
with 
\begin{align*}
\zeta_{n,1}(\varphi,\theta,\xi)=&\frac{1}{2\pi}\bigintsss_{\T} \partial_\eta\rho_n(\varphi,\theta+\eta) e^{\ii \,\xi\eta} d\eta \end{align*}
and
\begin{align*}
\lambda_1(\xi)&=\frac{1}{2\pi}\bigintsss_{\T}\frac{e^{\ii \theta^\prime \xi}}{|\sin(\theta^\prime/2)|^\alpha}d\theta^\prime.
\end{align*}
By applying   \eqref{fract1} and \eqref{Fo-c}, we find that
\begin{align*}
\lambda_1(\xi)&=\big(|\textnormal D|^{\alpha-1}e^{-\ii  \,\theta\xi}\big)_{|\theta=0}\\
&=\frac{2^{\alpha}\Gamma(1-\alpha)}{\Gamma(\frac{\alpha}{2})\Gamma(1-\frac{\alpha}{2})}\frac{ \Gamma\big(|\xi|+\frac{\alpha}{2}\big)}{\Gamma\big(|\xi|+1-\frac{\alpha}{2}\big)}\cdot
\end{align*} 
{From Lemma \ref{Stirling-formula} we get  that for any $j,\beta\in\NN, \epsilon>0$}
\begin{align}\label{Eq-bvm2}
\sup_{\alpha\in(0,\overline\alpha)}\big|\Delta_\xi^\beta\partial_\alpha^j \lambda_1(\xi)\big|\lesssim \langle \xi\rangle^{\overline\alpha-1+\epsilon-\beta}.
\end{align}
Now, according to the  convolution structure in the sum \eqref{conv-2021} we infer 
\begin{align*}
\Delta_\xi^\gamma\mu_{n,2}(\varphi,\theta,\xi)&=\sum_{\xi^\prime \in\Z} \zeta_{n,1}(\varphi,\theta,\xi-\xi^\prime)\,\Delta_{\xi^\prime}^\gamma\big[\lambda_1(\xi^\prime)a_n(\varphi,\theta,\xi^\prime)\big].
\end{align*}
Using the law products  with \eqref{Eq-bvm2} and performing similar arguments as in \eqref{sum-biz-125} we get for any $N_1,N_2\in\N$
\begin{align*}
\left\|\Delta_\xi^\gamma\mu_{n,2}(\cdot,\centerdot,\xi)\right\|_{s}^{q,\kappa}\lesssim&\|\rho_n\|_{{s_0+N_1+1}}^{q,\kappa}\interleave a_n\interleave_{m,s,0}^{q,\kappa}\sum_{|\xi^\prime|\leqslant\frac12|\xi| \atop|\xi^\prime|\geqslant2|\xi|}
\langle\xi-\xi^\prime\rangle^{-N_1} \langle\xi^\prime\rangle^{m+\overline\alpha+\epsilon-1}\\
&+\|\rho_n\|_{{s_0+N_2+1}}^{\kappa}\interleave a_n\interleave_{m,s,\gamma}^{q,\kappa}\sum_{\frac12|\xi|\leqslant|\xi^\prime|\leqslant 2|\xi|}\langle\xi-\xi^\prime\rangle^{-N_2} \langle\xi^\prime\rangle^{m+\overline\alpha+\epsilon-1-\gamma}\\
&\quad+\|\rho_n\|_{{s+N_1+1}}^{q,\kappa}\interleave a_n\interleave_{m,s_0,0}^{q,\kappa}\sum_{|\xi^\prime|\leqslant\frac12|\xi| \atop|\xi^\prime|\geqslant2|\xi|}
\langle\xi-\xi^\prime\rangle^{-N_1} \langle\xi^\prime\rangle^{m+\overline\alpha+\epsilon-1}\\
&\qquad+\|\rho_n\|_{{s+N_2+1}}^{q,\kappa}\interleave a_n\interleave_{m,s_0,\gamma}^{q,\kappa}\sum_{\frac12|\xi|\leqslant|\xi^\prime|\leqslant 2|\xi|}\langle\xi-\xi^\prime\rangle^{-N_2} \langle\xi^\prime\rangle^{m+\overline\alpha+\epsilon-1-\gamma}.
 \end{align*}
Fixing  $N_1=\gamma+1+|m|$ and $N_2=2$  then 
\begin{align}\label{sigma14}
\nonumber \interleave\mu_{n,2}\interleave_{m,s,\gamma}^{q,\kappa}\lesssim&\|\rho_n\|_{{s_0+2+\gamma+|m|}}^{q,\kappa}\interleave a_n\interleave_{m,s,0}^{q\kappa}+\|\rho_n\|_{{s_0+3}}^{q,\kappa}\interleave a_n\interleave_{m,s,\gamma}^{q,\kappa}\\
+&\|\rho_n\|_{{s+2+\gamma+|m|}}^{q,\kappa}\interleave a_n\interleave_{m,s_0,0}^{q,\kappa}+\|\rho_n\|_{{s+3}}^{q,\kappa}\interleave a_n\interleave_{m,s_0,\gamma}^{q,\kappa}.
 \end{align}
As to the term $\mathcal{R}_{n,3}K_n$ in \eqref{Remaind1} we first write in view of \eqref{fract1}  
\begin{align*}
\partial_\theta\big[|\textnormal{D}_\theta|^{\alpha-1},\rho_n\big]K_n(\varphi,\theta,\eta)&=\frac{1}{2\pi}\partial_\theta\bigintsss_{\T}\frac{\rho_n(\varphi,\theta-\theta^\prime)-\rho_n(\varphi,\theta)}{|\sin(\theta^\prime/2)|^\alpha}K_n(\varphi,\theta-\theta^\prime,\eta)d\theta^\prime\\
&=\bigintsss_{\T}\Psi(\varphi,\theta,\theta^\prime)K_n(\varphi,\theta-\theta^\prime,\eta)d\theta^\prime,
\end{align*}
with
\begin{align}\label{Psi-sym1}
 \Psi(\varphi,\theta,\theta^\prime)=&-\frac{\alpha}{4\pi}\frac{\rho_n(\varphi,\theta-\theta^\prime)-\rho_n(\varphi,\theta)}{\sin(\theta^\prime/2)}\frac{\cos(\theta^\prime/2)}{|\sin(\theta^\prime/2)|^\alpha}-\frac{\partial_\theta\rho_n(\varphi,\theta)}{2\pi|\sin(\theta^\prime/2)|^\alpha}\cdot
\end{align}
Therefore
\begin{align*}
\partial_\theta\big[|\textnormal{D}_\theta|^{\alpha-1},\rho_n\big]K_n(\varphi,\theta,\eta)&=\frac{1}{2\pi}\sum_{\xi\in\Z}e^{\ii (\theta-\eta)\xi}\bigintsss_{\T}\Psi(\varphi,\theta,\theta^\prime)a_n(\varphi,\theta-\theta^\prime,\xi)e^{-\ii  \theta^\prime\xi}d\theta^\prime\\
&=\sum_{\xi\in\Z}\mu_{n,3}(\varphi,\theta,\xi)e^{\ii (\theta-\eta)\xi}
\end{align*}
where the symbol $\mu_{n,3}$ is given by
$$
\mu_{n,3}(\varphi,\theta,\xi)=\frac{1}{2\pi}\bigintsss_{\T}\Psi(\varphi,\theta,\theta^\prime)a_n(\varphi,\theta-\theta^\prime,\xi)e^{-\ii  \theta^\prime\xi}d\theta^\prime
$$
We shall proceed as for $\sigma_{n,1}$ and one first writes by \eqref{Leibn-disc}
\begin{align*}
\Delta_\xi^\gamma\mu_{n,3}(\varphi,\theta,\xi)=&\frac{1}{2\pi}\sum_{0\le\beta\leqslant\gamma}\left(_\beta^\gamma\right)\int_{\T} \Psi(\varphi,\theta,\theta^\prime)\Delta_\xi^\beta\big[a_n(\varphi,\theta-\theta^\prime,\xi)\big]\Delta_\xi^{\gamma-\beta}[e^{-\ii  \,\theta^\prime\xi}](\xi+\beta)d\theta^\prime.
 \end{align*}
 Using \eqref{Formg1} then integration by parts implies
 \begin{align*}
\langle\xi\rangle^{-m}\xi^\gamma\Delta_\xi^\gamma\mu_{n,3}(\varphi,\theta,\xi)=&\sum_{0\leqslant\beta\leqslant\gamma\atop0\leqslant\eta\leqslant\gamma-\beta}\left(_\beta^\gamma\right)\left(_\eta^{\gamma-\beta}\right)i^{\gamma-\beta}\frac{\xi^{\gamma-\beta}}{(\xi+\beta)^{\gamma-\beta}}\int_{\T}  \partial_{\theta^\prime}^{\gamma-\beta-\eta}\Big[(1-e^{-\ii  \theta^\prime})^{\gamma-\beta} \Psi(\varphi,\theta,\theta^\prime)\Big]\\
&\qquad\qquad\langle\xi\rangle^{-m}\xi^\beta\partial_{\theta^\prime}^{\eta}\Big[{\Delta_\xi^\beta  a_n(\varphi,\theta-\theta^\prime,\xi)}\Big]e^{-\ii  \,\theta^\prime(\xi+\beta)}d\theta^\prime. 
 \end{align*}
From \eqref{Psi-sym1} we get through straightforward calculus 
 $$
 \left\|\partial_{\theta^\prime}^{\gamma-\beta-\eta}\Big[(1-e^{-\ii  \theta^\prime})^{\gamma-\beta}\Psi(\cdot,\centerdot,\theta^\prime)\Big]\right\|_{s}^{q,\kappa}\lesssim \frac{\|\rho_n\|_{{s+1+\gamma-\beta-\eta}}^{q,\kappa}}{|\sin(\theta^\prime/2)|^{\overline \alpha+\epsilon}}\cdot
 $$
 Hence we deduce from the law products stated in Lemma \ref{Law-prodX1}
 \begin{align}\label{sigma15}
\nonumber  \left\|\mu_{n,3}\right\|_{m,s,\gamma}^{q,\kappa}\lesssim&\sum_{0\leqslant\beta\leqslant\gamma\atop0\leqslant\eta\leqslant\gamma-\beta}\|\rho_n\|_{{s_0+1+\gamma-\beta-\eta}}^{q,\kappa}\interleave a_n\interleave_{m,s+\eta,\beta}^{q,\kappa}\\
+&\sum_{0\leqslant\beta\leqslant\gamma\atop0\leqslant\eta\leqslant\gamma-\beta}\|\rho_n\|_{{s+1+\gamma-\beta-\eta}}^{q,\kappa}\interleave a_n\interleave_{m,s_0+\eta,\beta}^{q,\kappa}.
 \end{align}
Concerning  the term $\mathcal{R}_{n,4}K_n$ in \eqref{Remaind1} we first write 
\begin{align}\label{sigma105}
\nonumber \partial_\eta\big[|\textnormal{D}_\eta|^{\alpha-1},\rho_n\big]K_n(\varphi,\theta,\eta)&=\frac{1}{2\pi}\partial_\eta\bigintsss_{\T}\frac{\rho_n(\varphi,\eta-\eta^\prime)-\rho_n(\varphi,\eta)}{|\sin(\eta^\prime/2)|^\alpha}K_n(\varphi,\theta,\eta-\eta^\prime)d\eta^\prime\\
&=\bigintsss_{\T}\Psi(\varphi,\eta,\eta^\prime)K_n(\varphi,\theta,\eta-\eta^\prime)d\eta^\prime
\end{align}
where $\Psi$ is introduced in \eqref{Psi-sym1}.
Therefore the associated symbol is given by
\begin{align*}
\mu_{n,4}(\varphi,\theta,\xi)&=\frac{1}{2\pi}\sum_{\xi^\prime\in\Z}a_n(\varphi,\theta,\xi^\prime)\bigintsss_{\T}\left(\int_\T\Psi(\varphi,\theta+\eta,\eta^\prime)e^{\ii \eta^\prime\xi^\prime}d\eta^\prime\right)e^{\ii (\xi-\xi^\prime)\eta} d\eta.
\end{align*}
From the summation by parts \eqref{Sum-parts} and Leibniz rule \eqref{Leibn-disc} we infer
\begin{align*}
\Delta_\xi^\gamma\mu_{n,4}(\varphi,\theta,\xi)&=\frac{1}{2\pi}\sum_{\xi'\in\Z}\bigintsss_{\T}\Delta_{\xi^\prime}^\gamma\left[a_n(\varphi,\theta,\xi^\prime)\left(\int_\T\Psi(\varphi,\theta+\eta,\eta^\prime)e^{\ii \eta^\prime\xi^\prime}d\eta^\prime\right)\right]e^{\ii (\xi-\xi^\prime)\eta} d\eta\\
&=\frac{1}{2\pi}\sum_{\xi'\in\Z\atop 0\leqslant\beta\leqslant\gamma}\left(_\beta^\gamma\right)\Delta_{\xi^\prime}^\beta a_n(\varphi,\theta,\xi^\prime)\bigintssss_{\T}\left(\Delta_{\xi^\prime}^{\gamma-\beta} \mathcal{P}(\varphi,\theta,\eta,\xi^\prime+\beta)\right)e^{\ii (\xi-\xi^\prime)\eta} d\eta,
\end{align*}
with
$$
\mathcal{P}(\varphi,\theta,\eta,\xi)\triangleq \int_\T\Psi(\varphi,\theta+\eta,\eta^\prime)e^{\ii \eta^\prime\xi^\prime}d\eta^\prime.
$$
By virtue of \eqref{Formg1} combined with integration by parts, we find for any $N\in\N$ and $0\leqslant \varrho\leqslant\gamma-\beta$
\begin{align*}
\left\|\int_{\T}\left(\Delta_{\xi^\prime}^{\gamma-\beta} \mathcal{P}(\varphi,\theta,\eta,\xi^\prime+\beta)\right)e^{-\ii  (\xi-\xi^\prime)\eta} d\eta\right\|_{s}^{q,\kappa}&\lesssim\int_\T\big\|\partial_{\eta^\prime}^{\varrho}\partial_\eta^N\big((1-e^{-\ii  \eta^\prime})^{\varrho}\Psi(\cdot,\centerdot+\eta,\eta^\prime)\big)\big\|_{s}^{q,\kappa}d\eta^\prime d\eta\\
&\times \langle\xi-\xi^\prime\rangle^{-N}\langle\xi^\prime\rangle^{-\varrho}.
\end{align*}
Using \eqref{Psi-sym1} one obtains
\begin{align*}
\int_\T\big\|\partial_{\eta^\prime}^{\varrho}\partial_\eta^N\big((1-e^{-\ii  \eta^\prime})^{\gamma-\beta}\Psi(\cdot,\centerdot+\eta,\eta^\prime)\big)\big\|_{s}^{q,\kappa}d\eta^\prime d\eta
&\lesssim \|\rho_n\|_{{s+1+N+\varrho}}^{q,\kappa}\int_\T\frac{d\eta^\prime}{|\sin(\eta^\prime/2)|^{\overline\alpha+\epsilon}}\\
&\lesssim \|\rho_n\|_{{s+1+N+\varrho}}^{q,\kappa}\end{align*}
Hence using the law products we get for any $N_1,N_2\in\N$
\begin{align*}
\left\|\Delta_\xi^\gamma\mu_{n,4}(\cdot,\centerdot,\xi)\right\|_{s}^{q,\kappa}\lesssim&\|\rho_n\|_{{s_0+1+N_1}}^{q,\kappa}\interleave a_n\interleave_{m,s,0}^{q,\kappa}\sum_{|\xi^\prime|\leqslant\frac12|\xi| \atop|\xi^\prime|\geqslant2|\xi|}
\langle\xi-\xi^\prime\rangle^{-N_1} \langle\xi^\prime\rangle^m\\
&+\sum_{0\leqslant\beta\leqslant \gamma}\interleave a_n\interleave_{m,s,\beta}^{q,\kappa}\|\rho_n\|_{{s_0+1+N_2+\gamma-\beta}}^{q,\kappa}\sum_{\frac12|\xi|\leqslant|\xi^\prime|\leqslant 2|\xi|}\langle\xi-\xi^\prime\rangle^{-N_2} \langle\xi^\prime\rangle^{m-\gamma}\\
&\quad+\|\rho_n\|_{{s+1+N_1}}^{q,\kappa}\interleave a_n\interleave_{m,s_0,0}^{q,\kappa}\sum_{|\xi^\prime|\leqslant\frac12|\xi| \atop|\xi^\prime|\geqslant 2|\xi|}
\langle\xi-\xi^\prime\rangle^{-N_1}\langle\xi^\prime\rangle ^m\\
&\qquad+\sum_{0\leqslant\beta\leqslant \gamma}\interleave a_n\interleave_{m,s_0,\beta}^{q,\kappa}\|\rho_n\|_{{s+1+N_2+\gamma-\beta}}^{q,\kappa}\sum_{\frac12|\xi|\leqslant|\xi^\prime|\leqslant 2|\xi|}\langle\xi-\xi^\prime\rangle^{-N_2} \langle\xi^\prime\rangle^{m-\gamma}.
 \end{align*}
Then making the choice  $N_1=\gamma+2+|m|$ and $N_2=2$  yields
\begin{align}\label{sigma16}
\nonumber \interleave\mu_{n,4}\interleave_{m,s,\gamma}^{q,\kappa}\lesssim&\|\rho_n\|_{{s_0+3+\gamma+|m|}}^{q,\kappa}\interleave a_n\interleave_{m,s,0}^{q,\kappa}+\|\rho_n\|_{{s+3+\gamma+|m|}}^{q,\kappa}\interleave a_n\interleave_{m,s_0,0}^{q,\kappa}\\
 +& \sum_{0\leqslant\beta\leqslant \gamma}\interleave a_n\interleave_{m,s,\beta}\|\rho_n\|_{{s_0+3+\gamma-\beta}}^{q,\kappa}+ \sum_{0\leqslant\beta\leqslant \gamma}\interleave a_n\interleave_{m,s_0,\beta}^{q,\kappa}\|\rho_n\|_{{s+3+\gamma-\beta}}^{q,\kappa}.
\end{align}
It remains to analyze the last term  $ \mu_{n,5}$ in \eqref{Remaind1}. 
Similarly to \eqref{sigma105} we get 
\begin{align}\label{sigma106}
\nonumber |\textnormal{D}_\eta|^{\alpha-1}(\partial_\eta\rho)\,K_n(\varphi,\theta,\eta)&=\frac{1}{2\pi}\bigintsss_{\T}\frac{\partial_\eta\rho_n(\varphi,\eta-\eta^\prime)}{|\sin(\eta^\prime/2)|^\alpha}K_n(\varphi,\theta,\eta-\eta^\prime)d\eta^\prime\\
&\triangleq \bigintsss_{\T}\Psi_1(\varphi,\eta,\eta^\prime)K_n(\varphi,\theta,\eta-\eta^\prime)d\eta^\prime.
\end{align}
Therefore the associated symbol takes the form,
\begin{align*}
\mu_{n,5}(\varphi,\theta,\xi)&=\frac{1}{2\pi}\sum_{\xi^\prime\in\Z}a_n(\varphi,\theta,\xi^\prime)\bigintsss_{\T}\left(\int_\T\Psi_1(\varphi,\theta+\eta,\eta^\prime)e^{\ii \eta^\prime\xi^\prime}d\eta^\prime\right)e^{-\ii  (\xi-\xi^\prime)\eta} d\eta.
\end{align*}
By reproducing exactly the same lines of the proof of \eqref{sigma16} we find
\begin{align}\label{sigma17}
\nonumber \interleave\mu_{n,5}\interleave_{m,s,\gamma}^{q,\kappa}\lesssim&\|\rho_n\|_{{s_0+3+\gamma+|m|}}^{q,\kappa}\nonumber\interleave a_n\interleave_{m,s,0}^{q,\kappa}\|\rho_n\|_{{s+3+\gamma+|m|}}^{q,\kappa}\interleave a_n\interleave_{m,s_0,0}^{q,\kappa}\\
\nonumber &\quad+ \sum_{0\leqslant\beta\leqslant \gamma}\interleave a_n\interleave_{m,s,\beta}^{q,\kappa}\|\rho_n\|_{{s_0+3+\gamma-\beta}}^{q,\kappa}\\
&\qquad+ \sum_{0\leqslant\beta\leqslant \gamma}\interleave a_n\interleave_{m,s_0,\beta}^{q,\kappa}\|\rho_n\|_{{s+3+\gamma-\beta}}^{q,\kappa}.
\end{align}
Putting together \eqref{sigma10}, \eqref{sigma11}, \eqref{sigma12}, \eqref{sigma13}, \eqref{sigma14}, \eqref{sigma15}, \eqref{sigma16} and \eqref{sigma17} leads to
\begin{align}\label{Tfg1}
\nonumber \interleave a_{n+1}\interleave_{m,s,\gamma}^{q,\kappa}\lesssim&\|\rho_n\|_{{s_0+3+\gamma+|m|}}^{q,\kappa}\interleave a_n\interleave_{m,s,0}^{q,\kappa}\|\rho_n\|_{{s+3+\gamma+|m|}}^{q,\kappa}\interleave a_n\interleave_{m,s_0,0}^{q,\kappa}\\
\nonumber&+\|\rho_n\|_{{s_0}}^{q,\kappa}\interleave a_{n}\interleave_{m,s+1,\gamma}^{q,\kappa}+\|\rho_n\|^{q,\kappa}_{s}\interleave a_{n}\interleave_{m,s_0+1,\gamma}^{q,\kappa}\\
\nonumber&\quad+\|\rho_n\|_{{s_0+3}}^{q,\kappa}\interleave a_n\interleave_{m,s,\gamma+1}^{q,\kappa}+\|\rho_n\|_{q,{s+3}}^{q,\kappa}\interleave a_n\interleave_{m,s_0,\gamma+1}^{q,\kappa} \\
\nonumber&\qquad+\sum_{0\leqslant\beta\leqslant\gamma\atop0\leqslant\eta\leqslant\gamma-\beta}\|\rho_n\|_{{s_0+3+\gamma-\beta-\eta}}^{q,\kappa}\interleave a_n\interleave_{m,s+\eta,\beta}^{q,\kappa}\\
&\qquad\quad+\sum_{0\leqslant\beta\leqslant\gamma\atop0\leqslant\eta\leqslant\gamma-\beta}\|\rho_n\|_{{s+3+\gamma-\beta-\eta}}^{q,\kappa}\interleave a_n\interleave_{m,s_0+\eta,\beta}^{q,\kappa}. 
\end{align}
Now by interpolation we get for any $0\leqslant \eta\leqslant\gamma-\beta$
\begin{align*}
\interleave a_n\interleave_{m,s+\eta,\beta}^{q,\kappa}\le \left(\interleave a_n\interleave_{m,s,\beta}^{q,\kappa}\right)^{1-\frac{\eta}{\gamma-\beta}}\left(\interleave a_n\interleave_{m,s+\gamma-\beta,\beta}^{q,\kappa}\right)^{\frac{\eta}{\gamma-\beta}}.
\end{align*}
and
\begin{align*}
\|\rho_n\|_{{s_0+3+\gamma-\eta-\beta}}^{q,\kappa}\le\left(\|\rho_n\|_{{s_0+3+\gamma-\beta}}^{q,\kappa}\right)^{1-\frac{\eta}{\gamma-\beta}}\left(\|\rho_n\|_{{s_0+3}}^{q,\kappa}\right)^{\frac{\eta}{\gamma-\beta}}.
\end{align*}
Consequently
\begin{align*}
\interleave a_n\interleave_{m,s+\eta,\beta}^{q,\kappa}\|\rho_n\|_{{s_0+3+\gamma-\eta-\beta}}^{q,\kappa}&\le\left(\|\rho_n\|_{{s_0+3+\gamma-\beta}}^{q,\kappa} \interleave a_n\interleave_{m,s,\beta}^{q,\kappa}\right)^{1-\frac{\eta}{\gamma-\beta}}\big(\|\rho_n\|_{{s_0+3}}^{q,\kappa}\interleave a_n\interleave_{m,s+\gamma-\beta,\beta}^{q,\kappa}\big)^{\frac{\eta}{\gamma-\beta}}\\
&\lesssim\|\rho_n\|_{{s_0+3+\gamma-\beta}}^{q,\kappa} \interleave a_n\interleave_{m,s,\beta}^{q,\kappa}+\|\rho_n\|_{{s_0+3}}^{q,\kappa}\interleave a_n\interleave_{m,s+\gamma-\beta,\beta}^{q,\kappa}.
\end{align*}
Similar arguments yield
\begin{align*}
\|\rho_n\|_{{s+3+\gamma-\beta-\eta}}^{q,\kappa} \interleave a_n\interleave_{m,s_0+\eta,\beta}^{q,\kappa}\lesssim\|\rho_n\|_{{s+3+\gamma-\beta}}^{q,\kappa} \interleave a_n\interleave_{m,s_0,\beta}^{q,\kappa}+\|\rho_n\|_{{s+3}}^{q,\kappa}\interleave a_n\interleave_{m,s_0+\gamma-\beta,\beta}^{q,\kappa}.
\end{align*}
Therefore we deduce  from  \eqref{Tfg1} and  Sobolev embeddings
\begin{align}\label{Tfg2}
\nonumber \interleave a_{n+1}\interleave_{m,s,\gamma}^{q,\kappa}\lesssim&\|\rho_n\|_{{s_0+3}}^{q,\kappa}\interleave a_{n}\interleave_{m,s+1,\gamma}^{q,\kappa}+\|\rho_n\|_{{s+3}}^{q,\kappa}\interleave a_{n}\interleave_{m,s_0+1,\gamma}^{q,\kappa}\\
\nonumber&+\|\rho_n\|_{{s_0+3}}^{q,\kappa}\interleave a_n\interleave_{m,s,\gamma+1}^{q,\kappa}+\|\rho_n\|_{{s+3}}^{q,\kappa}\interleave a_n\interleave_{m,s_0,\gamma+1}^{q,\kappa} \\
\nonumber&\quad+\sum_{0\leqslant\beta\leqslant\gamma}\|\rho_n\|_{{s_0+3+|m|+\gamma-\beta}}^{q,\kappa} \interleave a_n\interleave_{m,s,\beta}^{q,\kappa}+\|\rho_n\|_{{s_0+3}}^{q,\kappa}\interleave a_n\interleave_{m,s+\gamma-\beta,\beta}^{q,\kappa}\\
\nonumber&\qquad+
\sum_{0\leqslant\beta\leqslant\gamma}\|\rho_n\|_{{s+3+|m|+\gamma-\beta}}^{\kappa} \interleave a_n\interleave_{m,s_0,\beta}^{q,\kappa}+\|\rho_n\|_{{s+3}}^{q,\kappa}\interleave a_n\interleave_{m,s_0+\gamma-\beta,\beta}^{q,\kappa}\\
&\quad\qquad\triangleq \sum_{\ell=1}^6\mathcal{T}_\ell.
\end{align}
This achieves the proof of the first point.

\smallskip

$\bf{(ii)}$ Let  $s\geqslant s_0, \gamma, n\in\N$  be fixed numbers and assume \eqref{Assu1}. We will check that for any integer $k\in[0,n] $ we have  the following property: $\forall s^\prime\geqslant s_0, \gamma^\prime\in\N$ with $ s^\prime+\gamma^\prime\leqslant s+\gamma+n-k,\,$ we have
\begin{align}\label{Es-Ind-P89}
 \nonumber \interleave a_{k}\interleave_{m,s',\gamma'}^{q,\kappa}&\lesssim F\big(s_0+3\big)\sum_{i=0}^{k-1} \|\rho_i\|_{{s'+\gamma'+k+4+|m|}}^{q,\kappa}\prod_{j=0\atop  j\neq i}^{k-1} \|\rho_j\|_{{s_0+3+|m|}}^{q,\kappa}\\
 &\quad+F\big(s'+\gamma'+k+4\big) \prod_{j=0}^{k-1} \|\rho_j\|_{{s_0+3+|m|}}^{q,\kappa},
\end{align}
where we take for $k=0$ the convention $\displaystyle{\sum_{i=0}^{-1}..=0}$\quad and \quad $\displaystyle{\prod_{j=0}^{-1}=1}$ in such a way that the inequality remains true for $k=0$ according to the assumption \eqref{Assu1}. We shall proceed with finite  induction principle over $k$.   Assume that the estimate \eqref{Es-Ind-P89} is true at the order $k\leqslant n-1$ and let us show that it remains true at the order $k+1.$ Set
\begin{align*}U_k(s,\gamma)\triangleq &
F(s_0+3)\sum_{i=0}^{k-1} \|\rho_i\|_{{s+\gamma+k+4+|m|}}^{q,\kappa}\prod_{j=0\atop  j\neq i}^{k-1} \|\rho_j\|_{{s_0+3+|m|}}^{q,\kappa}\\
 &+F(s+\gamma+k+4) \prod_{j=0}^{k-1} \|\rho_j\|_{{s_0+3+|m|}}^{q,\kappa}.
 \end{align*}
We intend to prove that for any $ s^\prime\geqslant s_0, \gamma^\prime\in\N$ with $ s^\prime+\gamma^\prime\leqslant s+\gamma+n-k-1,\,$ we have
 $$
  \interleave a_{k+1}\interleave_{m,s',\gamma'}^{q,\kappa}\lesssim  U_{k+1}(s',\gamma').
 $$
Let us estimate each term $\mathcal{T}_\ell$ in \eqref{Tfg2} with exchanging  $n$ by $k$ and $(s,\gamma)$ by $(s',\gamma')$. For the first one we simply write in view of \eqref{Es-Ind-P89}, since $ s^\prime+1+\gamma^\prime\leqslant s+\gamma+n-k,\,$
\begin{align*}
\mathcal{T}_1=\|\rho_k\|_{{s_0+3}}^{q,\kappa}\interleave a_{k}\interleave_{m,s'+1,\gamma'}^{q,\kappa}&\le F(s_0+3)\|\rho_k\|_{{s_0+3}}^{q,\kappa}\sum_{i=0}^{k-1} \|\rho_i\|_{{s'+k+5+\gamma'+|m|}}^{q,\kappa}\prod_{j=0\atop  j\neq i}^{k-1} \|\rho_j\|_{{s_0+3+|m|}}^{q,\kappa}\\
 &\quad+F(s'+\gamma'+k+5) \|\rho_k\|_{{s_0+3}}^{q,\kappa}\prod_{j=0}^{k-1} \|\rho_j\|_{{s_0+3+|m|}}^{q,\kappa}\\
 & \qquad \le F(s_0+3)\sum_{i=0}^{k} \|\rho_i\|_{{s'+(k+1)+4+\gamma'+|m|}}^{q,\kappa}\prod_{j=0\atop  j\neq i}^{k} \|\rho_j\|_{{s_0+3+|m|}}^{q,\kappa}\\
 &\quad\qquad+F(s'+\gamma'+k+5) \prod_{j=0}^{k} \|\rho_j\|_{{s_0+3+|m|}}^{q,\kappa}.
\end{align*}
Hence, we obtain
$$
\mathcal{T}_1\leqslant U_{k+1}(s',\gamma').
$$
For the second term we write according to \eqref{Es-Ind-P89}, which can be applied since 
$$s_0+1+\gamma'\leqslant s'+1+\gamma' \leqslant s+\gamma+n-k,
$$
 combined with Sobolev embeddings
\begin{align*}
\mathcal{T}_2=\|\rho_k\|_{{s'+3}}^{q,\kappa'}\interleave a_{k}\interleave_{m,s_0+1,\gamma'}^{q,\kappa}&\le F(s_0+3)\|\rho_k\|_{{s'+3+|m|}}^{q,\kappa}\sum_{i=0}^{k-1} \|\rho_i\|_{{s_0+k+5+\gamma'+|m|}}^{q,\kappa}\prod_{j=0\atop  j\neq i}^{k-1} \|\rho_j\|_{{s_0+3+|m|}}^{q,\kappa}\\
 &\quad+F(s_0+\gamma'+k+5) \|\rho_k\|_{{s'+3+|m|}}^{q,\kappa}\prod_{j=0}^{k-1} \|\rho_j\|_{{s_0+3+|m|}}^{q,\kappa}.
\end{align*}
Applying once again Sobolev embeddings and interpolation inequality we find
\begin{align*}
 \|\rho_k\|_{{s'+3+|m|}}^{q,\kappa}\|\rho_i\|_{{s_0+k+5+\gamma'+|m|}}^{q,\kappa}
&\leqslant \eta \|\rho_k\|_{{s_0+3+|m|}}^{q,\kappa}\|\rho_i\|_{{s'+k+5+\gamma'+|m|}}^{q,\kappa}\\
&+(1-\eta)  \|\rho_i\|_{{s_0+3+|m|}}^{q,\kappa}\|\rho_k\|_{{s'+k+5+\gamma'+|m|}}^{q,\kappa}, \\
\eta&=\frac{s'-s_0}{s'-s_0+\gamma'+k+2}\cdot
\end{align*}
Consequently, we deduce that
\begin{align}\label{Covex1}
\nonumber& F(s_0+3)\|\rho_k\|_{{s'+3+|m|}}^{q,\kappa}\sum_{i=0}^{k-1} \|\rho_i\|_{{s_0+k+5+\gamma'+|m|}}^{q,\kappa}\prod_{j=0\atop  j\neq i}^{k-1} \|\rho_j\|_{{s_0+3+|m|}}^{q,\kappa}\\
\nonumber& \quad\leqslant F(s_0+3)\sum_{i=0}^{k} \|\rho_i\|_{{s'+k+5+\gamma'+|m|}}^{q,\kappa}\prod_{j=0\atop  j\neq i}^{k} \|\rho_j\|_{{s_0+3+|m|}}^{q,\kappa}\\
 &\qquad \leqslant U_{k+1}(s',\gamma').
\end{align}
Now from the log-convexity assumption on $F$ we obtain
\begin{align*}
 F(s_0+\gamma'+k+5) \leqslant\,& F^{\eta}(s_0+3)F^{1-\eta}(s'+\gamma'+k+5).
  \end{align*}
Combined this with the interpolation inequality
\begin{align*}
\|\rho_k\|_{{s'+3+|m|}}^{q,\kappa}&\le\big( \|\rho_k\|_{{s_0+3+|m|}}^{q,\kappa}\big)^{1-\eta}\big(\|\rho_k\|_{{s'+k+5+\gamma'+|m|}}^{q,\kappa}\big)^\eta&
\end{align*}
we find
\begin{align}\label{Convv}
\nonumber F(s_0+\gamma'+k+5) \|\rho_k\|_{{s'+3+|m|}}^{q,\kappa}\leqslant\,& \eta\,F(s_0+3) \|\rho_k\|_{{s'+\gamma'+|m|+k+5}}^{q,\kappa}\\
 &+(1-\eta)F(s'+\gamma'+k+5) \|\rho_k\|_{{s_0+3+|m|}}^{q,\kappa}.
\end{align}
It follows that
\begin{align*}
F(s_0+\gamma'+k+5)& \|\rho_k\|_{{s'+3+|m|}}^{q,\kappa}\prod_{j=0}^{k-1} \|\rho_j\|_{{s_0+3+|m|}}^{q,\kappa}\leqslant \eta\,F(s_0+3) \|\rho_k\|_{{s'+\gamma'+|m|+k+5}}^{q,\kappa} \\
&\times\prod_{j=0}^{k-1} \|\rho_j\|_{{s_0+3+|m|}}^{q,\kappa}+(1-\eta)F(s'+\gamma'+k+5) \prod_{j=0}^{k} \|\rho_j\|_{{s_0+3+|m|}}^{q,\kappa}\\
&\quad\leqslant \eta\, F(s_0+3)\sum_{i=0}^{k} \|\rho_i\|_{{s'+\gamma'+k+5+|m|}}^{q,\kappa}\prod_{j=0\atop  j\neq i}^{k} \|\rho_j\|_{{s_0+3+|m|}}^{q,\kappa}\\
&\qquad+(1-\eta)F(s'+\gamma'+k+5) \prod_{j=0}^{k} \|\rho_j\|_{{s_0+3+|m|}}^{q,\kappa}.
\end{align*}
Thus we find
\begin{align}\label{Covex2}
F(s_0+\gamma'+k+5) \|\rho_k\|_{{s'+3+|m|}}^{q,\kappa}\prod_{j=0}^{k-1} \|\rho_j\|_{{s_0+3+|m|}}^{q,\kappa}&\le U_{k+1}(s',\gamma')
\end{align}
Therefore we obtain
$$
\mathcal{T}_2\leqslant 2\,U_{k+1}(s',\gamma').
$$
Let us now move to $\mathcal{T}_3$, then the induction assumption \eqref{Es-Ind-P89} gives
\begin{align*}
 \mathcal{T}_3\leqslant\,& F(s_0+3)\sum_{i=0}^{k-1} \|\rho_i\|_{{s'+\gamma'+k+5+|m|}}^{q,\kappa}\prod_{j=0\atop  j\neq i}^{k} \|\rho_j\|_{{s_0+3+|m|}}^{q,\kappa}\\
 &\quad+F(s'+\gamma'+k+5)  \prod_{j=0}^{k} \|\rho_j\|_{{s_0+3+|m|}}^{q,\kappa}\\
 &\qquad\leqslant U_{k+1}(s',\gamma').
\end{align*}
As to $\mathcal{T}_4$, we may write using the induction assumption \eqref{Es-Ind-P89} 
\begin{align*}
 \mathcal{T}_4\leqslant\,& F(s_0+3)\|\rho_k\|_{{s'+3+|m|}}^{q,\kappa}\sum_{i=0}^{k-1} \|\rho_i\|_{{s_0+k+5+\gamma'+|m|}}^{q,\kappa}\prod_{j=0\atop  j\neq i}^{k-1} \|\rho_j\|_{{s_0+3+|m|}}^{q,\kappa}\\
 &+F(s_0+\gamma'+k+5) \|\rho_k\|_{{s'+3+|m|}}^{q,\kappa}\prod_{j=0}^{k-1} \|\rho_j\|_{{s_0+3+|m|}}^{q,\kappa}.\end{align*}
Applying \eqref{Covex1} and \eqref{Covex2} we infer 
\begin{align*}
 \mathcal{T}_4\leqslant& 2 U_{k+1}(s',\gamma').
 \end{align*}
Concerning the term $\mathcal{T}_5$ in \eqref{Tfg2} we first split it as follows
$$
\mathcal{T}_5=\sum_{0\leqslant\beta\leqslant\gamma'}\mathcal{T}_5^{1,\beta}+\mathcal{T}_5^{2,\beta}
$$
with
\begin{align*}
\mathcal{T}_5^{1,\beta}&\triangleq \|\rho_k\|_{{s_0+3+|m|+\gamma'-\beta}}^{q,\kappa} \interleave a_k\interleave_{m,s',\beta}^{q,\kappa}\\
\mathcal{T}_5^{2,\beta}&\triangleq \|\rho_k\|_{{s_0+3}}^{q,\kappa}\interleave a_k\interleave_{m,s'+\gamma'-\beta,\beta}^{q,\kappa}.
\end{align*}
Let us start with estimating $\mathcal{T}_5^{1,\beta}$. Applying Sobolev embeddings, the monotonicity of $F$ and \eqref{Es-Ind-P89}
\begin{align*}
 \mathcal{T}_5^{1,\beta}\leqslant& F(s_0+3)\|\rho_k\|_{{s_0+3+|m|+\gamma'-\beta}}^{q,\kappa}\sum_{i=0}^{k-1} \|\rho_i\|_{{s'+k+5+|m|+\beta}}^{q,\kappa}\prod_{j=0\atop  j\neq i}^{k-1} \|\rho_j\|_{{s_0+3+|m|}}^{q,\kappa}\\
 &+F(s'+\beta+k+5) \|\rho_k\|_{{s_0+3+|m|+\gamma'-\beta}}^{q,\kappa}\prod_{j=0}^{k-1} \|\rho_j\|_{{s_0+3+|m|}}^{q,\kappa}.
\end{align*}
Using the convexity inequality we infer
\begin{align*}
\|\rho_k\|_{{s_0+3+|m|+\gamma'-\beta}}^{q,\kappa}\|\rho_i\|_{{s'+k+5+|m|+\beta}}^{q,\kappa}&\lesssim \|\rho_k\|_{{s_0+3+|m|}}^{q,\kappa}\|\rho_i\|_{{s'+k+5+|m|+\gamma'}}^{q,\kappa}\\
&+\|\rho_i\|_{{s_0+3+|m|}}^{q,\kappa}\|\rho_k\|_{{s'+k+5+|m|+\gamma'}}^{q,\kappa}.
\end{align*}
Similar argument using the $\log$-convexity of $F$ yields
\begin{align*}
F(s'+\beta+k+5) \|\rho_k\|_{H^{s_0+3+|m|+\gamma'-\beta}}
&\lesssim \|\rho_k\|_{{s_0+3+|m|}}^{q,\kappa}\,F(s'+\gamma'+k+5) \\
&+\|\rho_k\|_{{s'+\gamma'+k+5+|m|}}^{q,\kappa}\,F(s_0+3).
\end{align*}
Consequently we get
\begin{align*}
 \mathcal{T}_5^{1,\beta}\lesssim& \,U_{k+1}(s',\gamma').
 \end{align*}
 As to the term $ \mathcal{T}_5^{2,\beta}$ we may write in view of \eqref{Es-Ind-P89}
 \begin{align*}
 \mathcal{T}_5^{2,\beta}\leqslant& \,F(s_0+3)\|\rho_k\|_{{s_0+3+|m|}}^{q,\kappa}\sum_{i=0}^{k-1} \|\rho_i\|_{{s'+\gamma'+k+5+|m|}}^{q,\kappa}\prod_{j=0\atop  j\neq i}^{k-1} \|\rho_j\|_{{s_0+3+|m|}}^{q,\kappa}\\
 &\quad+F(s'+\gamma'+k+5) \|\rho_k\|_{{s_0+3+|m|}}^{q,\kappa}\prod_{j=0}^{k-1} \|\rho_j\|_{{s_0+3+|m|}}^{q,\kappa}\\
&\qquad  \leqslant \,\,U_{k+1}(s',\gamma').\end{align*}
 Putting together the preceding estimates yields to
 $$
 \mathcal{T}_5\lesssim \,U_{k+1}(s',\gamma').
 $$
 It remains to analyze the last term $\mathcal{T}_6$ in \eqref{Tfg2} which can be decomposed as follows 
 $$
\mathcal{T}_6=\sum_{0\leqslant\beta\leqslant\gamma'}\mathcal{T}_6^{1,\beta}+\mathcal{T}_6^{2,\beta}
$$
with
\begin{align*}
\mathcal{T}_6^{1,\beta}&\triangleq \|\rho_k\|_{{s'+3+|m|+\gamma'-\beta}}^{q,\kappa}\interleave a_k\interleave_{m,s_0,\beta}^{q,\kappa}\\
\mathcal{T}_6^{2,\beta}&\triangleq \|\rho_k\|_{{s'+3}}^{q,\kappa}\interleave a_k\interleave_{m,s_0+\gamma'-\beta,\beta}^{q,\kappa}.
\end{align*}
To estimate the first one we make appeal to  \eqref{Es-Ind-P89}
 \begin{align*}
 \mathcal{T}_6^{1,\beta}\leqslant& F(s_0+3)\|\rho_k\|_{{s'+3+|m|+\gamma'-\beta}}^{q,\kappa}\sum_{i=0}^{k-1} \|\rho_i\|_{{s_0+k+5+|m|+\beta}}^{q,\kappa}\prod_{j=0\atop  j\neq i}^{k-1} \|\rho_j\|_{{s_0+3+|m|}}^{q,\kappa}\\
 &+F(s_0+k+5+\beta) \|\rho_k\|_{{s'+3+|m|+\gamma'-\beta}}^{q,\kappa}\prod_{j=0}^{k-1} \|\rho_j\|_{{s_0+3+|m|}}^{q,\kappa}.
 \end{align*}
 Then using  the convex inequalities 
 \begin{align*}
\|\rho_k\|_{{s'+3+|m|+\gamma'-\beta}}^{q,\kappa}\|\rho_i\|_{{s_0+k+5+|m|+\beta}}^{q,\kappa}&\lesssim \|\rho_k\|_{{s_0+3+|m|}}^{q,\kappa}\|\rho_i\|_{{s'+k+5+\gamma'+|m|}}^{q,\kappa}\\
&+ \|\rho_k\|_{{s'+k+5+\gamma'+|m|}}^{q,\kappa}\,\|\rho_i\|_{{s_0+3+|m|}}^{q,\kappa}
\end{align*}
 and
 \begin{align*}
\|\rho_k\|_{{s+3+|m|+\gamma'-\beta}}^{q,\kappa}F(s_0+k+5+\beta) &\lesssim  \|\rho_k\|_{{s_0+3+|m|}}^{q,\kappa}F(s+k+5+\gamma')\\
&+\|\rho_k\|_{{s'+k+5+\gamma'+|m|}}^{q,\kappa}F(s_0+3),
\end{align*}
one deduces that
\begin{align*}
 \mathcal{T}_6^{1,\beta}\lesssim &\,
 \, U_{k+1}(s',\gamma').
 \end{align*}
 The estimate of  the last $\mathcal{T}_6^{2,\beta}$ can be done in the same spirit. 
Finally, putting together the preceding estimates  we get: $\forall s^\prime\geqslant s_0, \gamma^\prime\in\N$ with $ s^\prime+\gamma^\prime\leqslant s+\gamma+n-k-1,\,$
\begin{align*}
\nonumber \interleave a_{k+1}\interleave_{m,s',\gamma'}^{q,\kappa}\lesssim & \, U_{k+1}(s',\gamma')
\end{align*}
and this achieves the induction in \eqref{Es-Ind-P89}.

$\bf{(iii)}$ 
Since $m<-\frac12$ then we can apply Lemma \ref{lemma-Sym-R}-(iii) allowing to get
\begin{align*}
\int_{\T}\big(\|K_n(\cdot,\centerdot,\centerdot+\eta)\|_{s}^{q,\kappa}\big)^2d\eta
&\lesssim\big(\interleave a_{n}\interleave_{m,s,0}^{q,\kappa}\big)^2.
\end{align*}
Then it suffices to apply    Lemma \ref{lemm-iter1}-(ii)  to get the desired result.

\end{proof}
\subsubsection{Iterated T\"oplitz matrix operators}
In this section we shall be concerned with the establishment of tame  estimates for   the quantities $
\|(\mathcal{M}_n^k\mathcal{Y}_n)(i)\|_{L^2}$ and $ \|(\overline{\mathcal{M}}_n^k\overline{\mathcal{Y}}_n)(i)\|_{L^2}
 $ encountered  before during the proof of Theorem \ref{Prop-EgorV}. Before stating ours result let us 
 recall that the T\"oplitz matrix operator  $\mathcal{M}_n$ was defined in \eqref{Matrix-y} whose entries are given in \eqref{Bnk}. More precisely, $$
 \mathcal{M}_n=\begin{pmatrix}
    0 &0&..&..&..&0  \\
m_{1,1}&0&..&..&0&0 \\
   m_{2,2}& m_{1,2}&0&.. &..&0\\
    ..&..&..&..&0&0\\
   m_{n,n}& m_{n-1,n}&..&..& m_{n,n}&0
  \end{pmatrix},\quad \mathcal{Y}_n=\begin{pmatrix}
  {K}_0 &  \\
    \partial_\chi {K}_0& \\
    ..&\\
    ..&\\
    \partial_\chi^n  {K}_0
  \end{pmatrix},\quad
$$
with
\begin{align*}
m_{l,j}&=\left(_k^n\right) \big(\widehat{\mathbb{A}}_{(k)}-\mathbb{S}_{(k)}\big)\\
&\triangleq \left(_k^n\right){\mathcal{D}_k},
\end{align*}
where
$$
\widehat{\mathbb{A}}_{(k)}=\mathbb{A}_{\theta,(k)}+\mathbb{A}_{\eta,(k)},\quad \mathbb{A}_{\theta,(k)}=\partial_\theta\big((\partial_\theta^{k}\rho)|\textnormal{D}_\theta|^{\alpha-1}+|\textnormal{D}_\theta|^{\alpha-1}(\partial_\theta^{k}\rho)\big)
$$
and
$$
\mathbb{S}_{(k)}=(\partial_\eta^{k+1}\rho)|\textnormal{D}_\eta|^{\alpha-1}+|\textnormal{D}_\eta|^{\alpha-1}(\partial_\eta^{k+1}\rho).
$$
Denote by $ a_{0}$  the symbol associated to ${K}_0$ defined through the relation,
\begin{align*}
{K}_0(\varphi,\theta,\eta)&=\sum_{\xi\in\Z}a_{0}(\varphi,\theta,\xi)\,e^{\ii (\theta-\eta)\xi}.
\end{align*}
Concerning the matrices $\overline{\mathcal{M}}_n$, they are  defined in \eqref{blocmatrix} and  the vectors $\overline{\mathcal{Y}}_n$ are given by \eqref{ourta1} and \eqref{LUX2}.
In what follows we shall state the  main result of this section where the T\"oplitz structure of the aforementioned matrices  plays a crucial role in getting tame estimates. To be more precise, we intend to  establish the following lemma. 
\begin{lemma}\label{V-13}
Let  $ \,n\in\N$  and $ m<-\frac12,$ then the following assertions hold true.
\begin{enumerate}
\item Let $s_0>\frac{d+1}{2}$ and assume the existence of   an   increasing $\log$-convex function $F_0$ such that
\begin{align*}\forall s\geqslant s_0,\gamma\in\N\quad\hbox{with}\quad s+\gamma\leqslant s_0+n+1\Longrightarrow\interleave a_{0}\interleave_{m,s,\gamma}\leqslant F_0(s+\gamma).
\end{align*}
Then we have for any $ k\in\{0,..,n\}$ and $i\in\{1,..,n+1\}$
$$
\big\|\mathcal{M}_n^k\mathcal{Y}_n(i)\big\|_{L^2(\T^{d+2})}\leqslant C\Big(F_0\big(s_0+5\big)+\|\rho\|_{H^{s_0+6+|m|}}\Big)^k\Big(F_0\big(s_0+5+i\big)+\|\rho\|_{H^{s_0+6+|m|+i}}\Big).
$$
We also get for $k\in\{1,..,n\},$
$$
\big\|\mathcal{M}_n^k\mathcal{Y}_n(i)\big\|_{L^2(\T^{d+2})}\leqslant C\|\rho\|_{H^{s_0+6+|m|}}^{k-1}\|\rho\|_{H^{s_0+6+|m|+i}}F_0\big(s_0+5+i\big).
$$
and for $k=0,\, i\in\{1,..,n+1\}$
$$
\big\|\mathcal{Y}_n(i)\big\|_{L^2(\T^{d+2})}\leqslant CF_0\big(s_0+5+i\big).
$$
\item  Let $q\in\N, s_0>\frac{d+1}{2}+q$ and assume the existence of an increasing  $\log$-convex function $F_0$ such that
$$
\forall s\geqslant s_0,\gamma\in\N\quad\hbox{with}\quad s+\gamma\leqslant s_0+n+1\Longrightarrow\interleave a_{0}\interleave_{m,s,\gamma}^{q,\kappa}\le F_0(s+\gamma).
$$ 
Then we have for any $ k\in\{0,..,n+q\}, i\in\{1,..,n+1\}$
$$
\big\|\overline{\mathcal{M}}_n^k\overline{\mathcal{Y}}_n(i)\big\|_{L^2(\T^{d+2})}\leqslant C\Big(F_0\big(s_0+5\big)+\|\rho\|_{{s_0+6+|m|}}^{q,\kappa}\Big)^k\Big(F_0\big(s_0+5+i\big)+\|\rho\|_{{s_0+6+|m|+i}}^{q,\kappa}\Big).
$$
We also get
$$
\big\|\overline{\mathcal{M}}_n^k\overline{\mathcal{Y}}_n(i)\big\|_{L^2(\T^{d+2})}\leqslant  C\big(\|\rho\|_{{s_0+6+|m|}}^{q,\kappa}\big)^{k-1}
\|\rho\|_{s_0+6+|m|+i}^{q,\kappa}F_0\big(s_0+5+i\big)
$$
and for $k=0,$
$$
\big\|\overline{\mathcal{Y}}_n(i)\big\|_{L^2(\T^{d+2})}\leqslant CF_0\big(s_0+5+i\big).
$$
Here the bloc matrix operator $\overline{\mathcal{M}}_n$ was defined in \eqref{blocmatrix} and $\overline{\mathcal{Y}}_n$ in \eqref{ourta1} and \eqref{LUX2}.
\end{enumerate}
\end{lemma}
\begin{proof}
$(${\bf{i}}$)$
First remark that $\mathcal{M}_n\in \DDD_1$ and according to Lemma \ref{Lem-top} we have  $\mathcal{M}_n^k\in \DDD_k$ for any $k\in\{1,..,n\}$ and  $\mathcal{M}_n^{n+1}=0.$ In particular, we may write  $\mathcal{M}_n^k=(a_{i,j}^{k})_{1\leqslant i,j\leqslant n+1}\ $ with
$$
a_{i,j}^k=0,\quad \hbox{if}\quad  i-j< k\quad \hbox{or}\quad  i\leqslant  k.
$$
In addition, according to Lemma \ref{Lem-top} each coefficient $a_{i,j}^k$ is a linear combination of the iterative operator
\begin{align}\label{Tgv1}
\prod_{\ell=1}^k \mathcal{D}_{\alpha_\ell}\quad\hbox{such that}\quad \sum_{\ell=1}^k\alpha_\ell=i-j,\quad \alpha_\ell\geqslant1.
\end{align}
Denote ${Z}=\big(Z_1,...,Z_{n+1}\big)\triangleq \mathcal{M}_n^k\mathcal{Y}_n$, then
$$
{Z}_i=\sum_{j=1}^{i-k}a_{i,j}^k\partial_\chi^j  K_0.
$$
Consequently ${Z}_i$  is a linear combination of 
\begin{align}\label{decom-Lin12}
\prod_{\ell=1}^k \mathcal{D}_{\alpha_\ell}\partial_\chi^j  {K}_0\quad\hbox{such that}\quad \sum_{\ell=1}^k\alpha_\ell=i-j,\quad 1\leqslant j\leqslant i-k, \quad \alpha_\ell\geqslant1.
\end{align}
Using Lemma \ref{lemm-iter1}-$($iii$)$ with $q=0$ and  $m<-\frac12$, we get for using  Sobolev embeddings and change of variables 
\begin{align}\label{TatuXDD}
\nonumber \left\|\prod_{\ell=1}^k \mathcal{D}_{\alpha_\ell}\partial_\chi^j  {K}_0\right\|_{L^2_{\varphi,\theta,\eta}}&\lesssim \left(\bigintss_{\T}\left\|\left(\prod_{\ell=1}^k \mathcal{D}_{\alpha_\ell}\partial_\chi^j {K}_0\right)(\cdot,\centerdot,\centerdot+\eta)\right\|_{H^{s_0}_{\varphi,\theta}}^2d\eta\right)^{\frac12}
\\
\nonumber &\quad \lesssim F\big(s_0+3\big)\sum_{\ell=1}^{k} \|\rho\|_{H^{s_0+k+4+|m|+\alpha_\ell}}\prod_{p=1\atop  p\neq \ell}^{k} \|\rho\|_{H^{s_0+3+|m|+\alpha_p}}\\
 &\qquad +F\big(x{s_0}+k+4\big) \prod_{\ell=1}^{k} \|\rho\|_{H^{s_0+3+|m|+\alpha_\ell}},
\end{align}
where $F$ is defined as an increasing log-convex function such that if  $ a_{j,0}$ is  the symbol associated to the kernel  $\partial_\chi^j  {K}_0$ then
$$
\forall s\geqslant s_0,\gamma\in\N\quad\hbox{with}\quad s+\gamma\leqslant s_0+k\Longrightarrow \interleave a_{j,0}\interleave_{m,s,\gamma}\leqslant F(s+\gamma).
$$
Recall that  $ a_{0}$ is the symbol associated to ${K}_0$ and one gets easily the relation 
\begin{align*}
\partial_\chi^j  {K}_0(\varphi,\theta,\eta)&=\partial_\chi^j\sum_{\xi\in\Z}a_{0}(\varphi,\theta,\xi)\,e^{\ii (\theta-\eta)}\\
&=\sum_{\xi\in\Z}\partial_\theta^j a_{0}(\varphi,\theta,\xi)\,e^{\ii (\theta-\eta)},
\end{align*}
which implies that 
$$a_{j,0}(\varphi,\theta,\xi)=\partial_\theta^j a_{0}(\varphi,\theta,\xi).
$$
 Thus
 $$
 \interleave a_{j,0}\interleave_{m,s,\gamma}\lesssim \interleave a_{0}\interleave_{m,s+j,\gamma}.
$$
Consequently, if $F_0$ is  an increasing log-convex function such that
$$
\forall s\geqslant s_0,\gamma\in\N\quad\hbox{with}\quad s+\gamma\leqslant s_0+k+j\Longrightarrow\interleave a_{0}\interleave_{m,s,\gamma}\le F_0(s+\gamma)
$$
then we get 
\begin{align}\label{est-k=0-L}
\forall s\geqslant s_0,\gamma\in\N\quad\hbox{with}\quad s+\gamma\leqslant s_0+k\Longrightarrow\interleave a_{j,0}\interleave_{m,s,\gamma}\leqslant F_0(s+\gamma+j)
\end{align}
and the choice $F(x)\triangleq F_0(x+j)$ gives an increasing and log-convex function. Notice that from \eqref{decom-Lin12} one has $k+j\leqslant i\leqslant n+1$ and therefore it suffices to impose to $a_0$ the assumption
$$
\forall s\geqslant s_0,\gamma\in\N\quad\hbox{with}\quad s+\gamma\leqslant s_0+n+1\Longrightarrow\interleave a_{0}\interleave_{m,s,\gamma}\le F_0(s+\gamma)
$$
Consequently, we get from the estimate \eqref{TatuXDD}
\begin{align}\label{TatuX1}
\nonumber \left\|\prod_{\ell=1}^k \mathcal{D}_{\alpha_\ell}\partial_\chi^j {K}_0\right\|_{L^2_{\varphi,\theta,\eta}}^2&\lesssim F_0\big(s_0+3+j\big)\sum_{\ell=1}^{k} \|\rho\|_{H^{s_0+k+4+|m|+\alpha_\ell}}\prod_{p=1\atop  p\neq \ell}^{k} \|\rho\|_{H^{s_0+3+|m|+\alpha_p}}\\
 &+F_0\big(s_0+k+4+j\big) \prod_{\ell=1}^{k} \|\rho\|_{H^{s_0+3+|m|+\alpha_\ell}},
\end{align}
Now we shall deal with the first term of the right hand side. Fix $\ell\in\{1,..,k\},$ since  $\alpha_p\geqslant 1,$  then  interpolation inequality yields
$$
\|\rho\|_{H^{s_0+3+|m|+\alpha_p}}\leqslant \|\rho\|_{H^{s_0+4+|m|}}^{1-\mu_p}\|\rho\|_{H^{s_0+6+|m|+i}}^{\mu_p},\quad \mu_p=\frac{\alpha_p-1}{i+2}\cdot
$$
Then 
\begin{align*}
\prod_{p=1\atop  p\neq \ell}^{k} \|\rho\|_{H^{s_0+3+|m|+\alpha_p}}&\leqslant \|\rho\|_{H^{s_0+4+|m|}}^{{\Theta_{1,\ell}}}\|\rho\|_{H^{s_0+6+|m|+i}}^{\Theta_{2,\ell}},\\
\Theta_{1,\ell}=\sum_{p=1\atop p\neq \ell}^k(1-\mu_p)&,\quad \Theta_{2,\ell}=\sum_{p=1\atop p\neq \ell}^k\mu_p.
\end{align*}
 Similarly
$$
\|\rho\|_{H^{s_0+4+k+|m|+\alpha_\ell}}\leqslant \|\rho\|_{H^{s_0+4+|m|}}^{1-\mu_\ell}\|\rho\|_{H^{s_0+6+|m|+i}}^{\mu_\ell},\quad \mu_\ell=\frac{\alpha_\ell+k}{i+2}\cdot
$$
Hence
\begin{align}\label{sec-point1}
\nonumber\|\rho\|_{H^{s_0+4+k+|m|+\alpha_\ell}}\prod_{p=1\atop  p\neq \ell}^{k} \|\rho\|_{H^{s_0+3+|m|+\alpha_p}}&\leqslant \|\rho\|_{H^{s_0+4+|m|}}^{{\Theta_{1,\ell}+1-\mu_\ell}}\,\,\|\rho\|_{H^{s_0+6+|m|+i}}^{\Theta_{2,\ell}+\mu_\ell}\\
&\leqslant \|\rho\|_{H^{s_0+6+|m|}}^{{\Theta_{1,\ell}+1-\mu_\ell}}\,\,\|\rho\|_{H^{s_0+6+|m|+i}}^{\Theta_{2,\ell}+\mu_\ell}.
\end{align}
Now we may check in view of \eqref{Tgv1}
\begin{align*}
\Theta_{2,\ell}+\mu_\ell&=\frac{1}{i+2}\sum_{p=1\atop p\neq \ell}^k(\alpha_p-1)+\frac{\alpha_\ell+k}{i+2}\\
&=\frac{1}{i+2}\sum_{p=1}^k\alpha_p+\frac{1}{i+2}\\
&=\frac{i-j+1}{i+2}\cdot
\end{align*}
Using the convexity inequality and the monotonicity of $F_0$
\begin{align*}
F_0\big(s_0+3+j\big)&\leqslant F_0^{1-\frac{j}{i+2}}\big(s_0+3\big)F_0^{\frac{j}{i+2}}\big(s_0+5+i\big)\\
&\leqslant F_0^{1-\frac{j}{i+2}}\big(s_0+5\big)F_0^{\frac{j}{i+2}}\big(s_0+5+i\big).
\end{align*}
Therefore we find
\begin{align*}
  F_0\big(s_0+3+j\big)\sum_{\ell=1}^{k} \|\rho\|_{H^{s_0+k+4+|m|+\alpha_\ell}}\prod_{p=1\atop  p\neq \ell}^{k} \|\rho\|_{H^{s_0+3+|m|+\alpha_p}}&\lesssim  G_0^{k+\frac{1}{i+2}}(0)G_0^{\frac{i+1}{i+2}}(i)\\
 &\lesssim G_0^{k}(0)\,G_0(i),
\end{align*}
with 
$$
G_0(i)=F_0\big(s_0+5+i\big)+\|\rho\|_{H^{s_0+6+|m|+i}}
$$
where we have used  in the last inequality the fact $G_0(0)\leqslant G_0(i).$ The last term in \eqref{TatuX1} can be  treated in a similar way as before. Indeed, 
\begin{align}\label{STT-12}
\nonumber \prod_{\ell=1}^{k} \|\rho\|_{H^{s_0+3+|m|+\alpha_\ell}}&\leqslant \|\rho\|_{H^{s_0+4+|m|}}^{k-\frac{i-j-k}{i+2}}\,\,\|\rho\|_{H^{s_0+6+|m|+i}}^{\frac{i-j-k}{i+2}}\\
&\leqslant \big(G_0(0)\big)^{k-\frac{i-j-k}{i+2}}\,\,\big(G_0(i)\big)^{{\frac{i-j-k}{i+2}}}.
\end{align}
and
\begin{align*}
F_0\big(s_0+4+j+k\big)&\leqslant F_0^{\frac{1+i-k-j}{i+2}}\big(s_0+3\big)F_0^{\frac{1+k+j}{i+2}}\big(s_0+5+i\big)\\
&\leqslant \big(G_0(0)\big)^{\frac{1+i-k-j}{i+2}}\big(G_0(i)\big)^{\frac{1+k+j}{i+2}}.
\end{align*}
Consequently we find
\begin{align*}
  F_0\big(s_0+k+4+j\big) \prod_{\ell=1}^{k} \|\rho\|_{H^{s_0+3+|m|+\alpha_\ell}} &\lesssim G_0^{k+\frac{1}{i+2}}(0)G_0^{\frac{i+1}{i+2}}(i)\\
 &\lesssim G_0^{k}(0)\,G_0(i).
\end{align*}
This completes the proof of the first point of $(i)$. As to the second estimate we simply use for $k\geqslant1$ the inequality  \eqref{sec-point1} combined with Sobolev embeddings and the fact $\Theta_{2,\ell}+\mu_\ell<1$  in order to get into 
\begin{align}\label{sec-point2}
\|\rho\|_{H^{s_0+4+k+|m|+\alpha_\ell}}\prod_{p=1\atop  p\neq \ell}^{k} \|\rho\|_{H^{s_0+3+|m|+\alpha_p}}&\leqslant \|\rho\|_{H^{s_0+6+|m|}}^{k-1}\,\,\|\rho\|_{H^{s_0+6+|m|+i}}.
\end{align}\\
Concerning the estimate \eqref{STT-12} it is replaced in a similar way by
  \begin{align}\label{STT-13}
 \prod_{\ell=1}^{k} \|\rho\|_{H^{s_0+3+|m|+\alpha_\ell}}&\leqslant \|\rho\|_{H^{s_0+4+|m|}}^{k-1}\,\,\|\rho\|_{H^{s_0+6+|m|+i}}.
\end{align}
Combining both estimates yields to the desired estimate for $k\geqslant1$. Concerning the case $k=0$, it can derived from \eqref{est-k=0-L} and Lemma \ref{lemm-iter1}-$($iii$)$.

\smallskip

$(${\bf ii}$)$ 
Recall from \eqref{Nilp-ot-gen}  that  the matrix $\overline{\mathcal{M}}_n^{q+k}\in{\DDD}_k(n,q)$ and the set ${\DDD}_k(n,q)$ was defined before in Section \ref{Sect-Top-M}. In particular  $\overline{\mathcal{M}}_n^{k+q}=(M_{i,j}^{k})_{1\leqslant i,j\leqslant n+1}\ $ with
$$
M_{i,j}^{k}=0,\quad \hbox{if}\quad  i-j< k\quad \hbox{or}\quad  i\leqslant  k.
$$
In addition, according to Lemma \ref{Lem-top} and \eqref{LL1} each matrix entry  $M_{i,j}^k$ is a linear combination of the iterative matrix operator
\begin{align*}
\prod_{\ell=1}^{k+q} {\mu}_{\overline{\mathcal{M}}_n}({\alpha_\ell})\quad\hbox{such that}\quad \sum_{\ell=1}^{k+q}\alpha_\ell=i-j,\quad \alpha_\ell\geqslant{0}.
\end{align*}
From \eqref{LL1} one has ${\mu}_{\overline{\mathcal{M}}_n}(0)=\overline{M}_q\in\DDD_1$ which is nilpotent with ${\mu}_{\overline{\mathcal{M}}_n}^{1+q}=0$ and its main diagonal is vanishing. In addition ${\mu}_{\overline{\mathcal{M}}_n}({\alpha_\ell})\in\DDD_0$ and therefore we deduce from Lemma \ref{Lem-top}-$($i$)$
$$
\#\Big\{\ell, \alpha_\ell=0\Big\}\geqslant q+1\Longrightarrow \prod_{\ell=1}^{k+q} {\mu}_{\overline{\mathcal{M}}_n}({\alpha_\ell})=0
$$ 
Consequently $M_{i,j}^k$ is a linear combination of the iterative matrix operator
\begin{align}\label{Tgv11}
\prod_{\ell=1}^{k+q} {\mu}_{\overline{\mathcal{M}}_n}({\alpha_\ell})\quad\hbox{such that}\quad \sum_{\ell=1}^{k+q}\alpha_\ell=i-j\quad\hbox{and}\quad \#\Big\{\ell, \alpha_\ell\geqslant 1\Big\}\geqslant k
\end{align}
Set $\overline{Z}=(\overline{Z}_1,..,\overline{Z}_2)\triangleq \overline{\mathcal{M}}_n^{k+q}\overline{\mathcal{Y}}_n,$ then
$$
\overline{Z}_i=\sum_{j=1}^{i-k}M_{i,j}^k\partial_\chi^j  Y_q.
$$
Consequently we get in view of \eqref{Tgv11} that  $\overline{Z}_i$  is a linear combination of 
\begin{align}\label{Tgv111}
\prod_{\ell=1}^{k+q} {\mu}_{\overline{\mathcal{M}}_n}({\alpha_\ell})\partial_\chi^j  Y_q\quad\hbox{such that}\quad \sum_{\ell=1}^{k+q}\alpha_\ell=i-j,\quad \#\Big\{\ell, \alpha_\ell\geqslant 1\Big\}\geqslant k
\end{align}
Applying  Lemma \ref{lemm-iter1}-$($iii$)$ with slight adaptation  we get for $m<-\frac12$ and using  Sobolev embeddings and change of variables 
\begin{align}\label{TatuX}
\nonumber \left\|\prod_{\ell=1}^{k+q} {\mu}_{\overline{\mathcal{M}}_n}({\alpha_\ell})\partial_\chi^j  Y_q\right\|_{L^2_{\varphi,\theta,\eta}}^2&\lesssim \bigintss_{\T}\left\|\left(\prod_{\ell=1}^{k+q} {\mu}_{\overline{\mathcal{M}}_n}({\alpha_\ell})\partial_\chi^j  Y_q\right)(\cdot,\centerdot,\centerdot+\eta)\right\|_{H^{s_0}_{\varphi,\theta}}^2d\eta
\\
\nonumber &\lesssim F\big(s_0+3|\big)\sum_{\ell=1}^{k} \|\rho\|_{q,{s_0+k+4+|m|+\alpha_\ell}}^{q,\kappa}\prod_{p=1\atop  p\neq \ell}^{k} \|\rho\|_{q,{s_0+3+|m|+\alpha_p}}^{q,\kappa}\\
 &+F\big({s_0}+k+4\big) \prod_{\ell=1}^{k} \|\rho\|_{q,{s_0+3+|m|+\alpha_\ell}}^{q,\kappa},
\end{align}
where $F$ is defined as an increasing log-convex function such that if  $ \overline{a}_{j,0}$ is the matrix symbol of $\partial_\chi^j  Y_q,$ where $ Y_q$ is defined in \eqref{LUX2}, then
$$
\forall s\geqslant s_0,\gamma\in\N\quad\hbox{with}\quad s+\gamma\leqslant s_0+k\Longrightarrow\interleave \overline{a}_{j,0}\interleave_{m,s,\gamma}\leqslant F(s+\gamma).
$$
Let  $ \overline{a}_{0}$  be the symbol associated to the vector  $  Y_q$,  then we have the relation 
$$
\overline{a}_{j,0}=\partial_\theta^j\overline{a}_{0}
$$
Thus
 \begin{align*}
 \interleave \overline{a}_{j,0}\interleave_{m,s,\gamma}&\lesssim \interleave \overline{a}_{0}\interleave_{m,s+j,\gamma}\\
 &\lesssim \interleave {a}_{0}\interleave_{m,s+j,\gamma}^{q,\kappa}
\end{align*}
Consequently, if $F_0$ is  an increasing log-convex function such that
$$
\forall s\geqslant s_0,\gamma\in\N\quad\hbox{with}\quad s+\gamma\leqslant s_0+k+j\Longrightarrow\interleave {a}_{0}\interleave_{m,s,\gamma}^{q,\kappa}\le F_0(s+\gamma)
$$
then we get 
$$
\forall s\geqslant s_0,\gamma\in\N\quad\hbox{with}\quad s+\gamma\leqslant s_0+k\Longrightarrow\interleave \overline{a}_{j,0}\interleave_{m,s,\gamma}\leqslant F_0(s+\gamma+j)
$$
and the choice $F(x)=F_0(x+j)$ fits with the required assumptions. Remark  that  one gets in view  of  \eqref{Tgv11} that $k+j\leqslant i\leqslant n+1$ and therefore it suffices to impose to $a_0$ the assumption
$$
\forall s\geqslant s_0,\gamma\in\N\quad\hbox{with}\quad s+\gamma\leqslant s_0+n+1\Longrightarrow\interleave {a}_{0}\interleave_{m,s,\gamma}^{q,\kappa}\le F_0(s+\gamma)
$$
At this stage,    we  proceed exactly as in the preceding point $($i$)$ using interpolation inequalities in \eqref{TatuX} in order to get the desired result.
The estimate of the second estimate of the point (ii) can be done in the same spirit as in (i) using in particular analogous estimates to
\eqref{sec-point2} and \eqref{STT-13}. The proof is now complete.
\end{proof}

\section{Reducibility of the linearized operator}\label{Reducibility of the linearized operator}
In this section we shall investigate  the construction of a right approximate inverse for the linearized operator  $\widehat{\mathcal{L}}_{\omega}$ at a state $i_0$ close to the flat torus. This is the most delicate part in the construction of quasi-periodic solutions and needed during the implementation of Nash-Moser scheme . This operator was defined in \eqref{Lomega def} and it takes the form
\begin{equation}\label{Norm-proj}
\widehat{\mathcal{L}}_{\omega}=\widehat{\mathcal{L}}_{\omega}(i_{0})=\Pi_{\mathbb{S}_0}^\bot \big(\omega\cdot \partial_\varphi   - 
\partial_\theta  K_{02}(\varphi) \big)\Pi_{{\mathbb{S}}_0}^\bot
\end{equation}
 Notice that we are dealing with  an unbounded operator with variable coefficients and when it is evaluated at the flat torus, we simply find a diagonal operator that can be formally  inverted provided that we avoid resonances in a suitable way by imposing suitable Diophantine conditions. The intuitive strategy that we shall follow and which is a common fact on the most studies around this topic  is to conjugate    the linearized operator $\widehat{\mathcal{L}}_{\omega}$ into a    constant coefficients operator. This will be implemented  in three long steps. First, we  focus on the transport part  and perform   a quasi-periodic change of coordinates allowing to get  after a suitable conjugation   a new transport part  with constant coefficients provided that the parameters $\alpha,\omega$ are restricted  in a Cantor set. The construction of this transformation is based on  KAM reducibility procedure  in the spirit of the recent works \cite{BertiMontalto, BCP}. Then the final  outcome of the first step  is a new operator whose transport part is diagonal with unbounded  nonlocal  perturbation of order $\alpha$ supplemented with  a remainder of order zero. As to the the second step, it  consists in reducing to a constant coefficient the fractional  part of order $\alpha$ using a nonlocal hyperbolic flow constructed in Proposition \ref{flowmap00}, and  of course without altering  the transport part. The final product of the second step  is linear operator whose positive part is a Fourier multiplier with a small non-diagonal  bounded  perturbation. The main concern of the last step is to reduce to a constant operator the remainder by   applying the  KAM scheme as  in  the paper  \cite{BertiMontalto}. We emphasize that during this step we need to strengthen the topology of linear continuous operators and work with a suitable topology on pseudo-differential operators described in \eqref{Def-pseud-w}. One of the delicate point that one should face  at the end of step 2 is to check that the new remainder  stays in this  strong  topology. This issue has been  solved in Theorem \ref{Prop-EgorV} following a new approach based on the kernel dynamics and  refined structures on T\"oplitz matrix operators. 
We stress  that along the different reduction processes  we need to  impose some  non-resonance  assumptions  that will be finally recast in terms of a Cantor like set. \\
Next, we shall introduce the  parameters $\kappa, q, \tau_1, \tau_2$ and $\varrho$ connected with the geometric structure of the  Cantor sets that will emerge during the reduction process. We shall also fix some  regularity levels that will be made clear during the proofs.  \begin{align}\label{Conv-T2}
\kappa\in(0,1),&\quad(q,d)\in(\mathbb{N}^{*})^{2},\quad S\geqslant s\geqslant s_{0}>\frac{d+1}{2}+q+{3},\quad \tau_{2}>\tau_{1}>d.
\end{align}
where $S$ is a fixed large number and $\sigma=\sigma(d,\tau_1,\tau_2,q)$ is a  positive  number  related  to the loss of regularities during the different steps of the reduction. It depends only on the parameters governing  the geometric structure  of the involved Cantor sets and it is independent of the regularity. Notice that this implicit number may change  during the computations from one line to another. Since the exact value of $\sigma$ does not matter, then we find it convenient  to keep  this freedom more than making precise computations on the value of $\sigma$ which is subject to continuous changes from the step $1$ to step $3$. We observe that sometimes with some statements we  can do better and  use weaker  constraints  than those of \eqref{Conv-T2}. 
We shall also need during KAM reduction and Nash-Moser scheme to make appeal to cut-offs in frequency introduced in \eqref{definition of Nm}  
Notice that throughout this section, the constants $N_{0}\geqslant 2$ and $\kappa\in(0,1)$ are independent  free parameters, but during the implementation  of  Nash-Moser scheme that will be explored later in Section \ref{N-M-section12}, they will be suitably fixed with respect to the small parameter $\varepsilon$. \\
Before starting the reduction process we shall  discuss in the next two sections some algebraic and  analytical properties of the linearized operator.

 \subsection{Structure  on the normal direction}
 In this section we shall explore some analytical properties of the the operator $\widehat{\mathcal{L}}_{\omega}$. We shall  in particular give suitable estimates of its coefficients  and view this operator  as  a small perturbation of order zero of the linearized operator computed in \mbox{Proposition \ref{lin-eq-r},} which takes the form
 \begin{align}\label{lin}
\mathcal{L}_{r,\lambda} h =\omega\cdot\partial_\varphi h +\partial_\theta\big( V_{\alpha,r}(\varphi,\theta) h  -\mathbb{K}_{r,\alpha} h\big)
\end{align}
where $\lambda=(\omega,\alpha)$ and 
\begin{align}
\mathbb{K}_{r,\alpha} h (\varphi,\theta)&\triangleq \frac{ C_\alpha }{2\pi }\bigintsss_{0}^{2\pi} \frac{h(\varphi,\eta)}{A_r^{\alpha/2}(\varphi,\theta,\eta)} d\eta,\label{K-r}
\\
V_{r,\alpha}(\varphi,\theta)&\triangleq \Omega+\frac{ C_\alpha}{2\pi R(\varphi,\theta)}\bigintsss_{0}^{2\pi} \frac{\partial_{\eta} \big[R(\varphi,\eta)\sin(\eta-\theta)\big]}{A_r(\varphi,\theta,\eta)^{\alpha/2}} d\eta, \label{V-eq}
\\
A_r(\varphi,\theta,\eta)&\triangleq \Big(R(\varphi,\eta)-R(\varphi,\theta)\Big)^2+4R(\varphi,\eta)R(\varphi,\theta)\sin^2\left(\frac{\eta-\theta}{2}\right),\label{A-rM}
\\
R(\varphi,\theta)&\triangleq \big(1+ 2r(\varphi,\theta)\big)^{\frac12}.\nonumber 
\end{align}
The first main result of this section reads as follows.
\begin{proposition}\label{lemma-GS0}
{%
The operator $\widehat{\mathcal{L}}_{\omega}$ defined in \eqref{Norm-proj} takes the form 
$$\widehat{\mathcal{L}}_{\omega}= \Pi_{\mathbb{S}_0}^\bot\left(\mathcal{L}_{\varepsilon r_0,\lambda}-\varepsilon\partial_{\theta}\mathcal{R}\right) \Pi_{\mathbb{S}_0}^\bot$$
where the operator  $\mathcal{L}_{\varepsilon r_0,\lambda}$  is defined in \eqref{lin}, with
\begin{align*}r_0(\varphi,\cdot)&\triangleq  \mathbf{A}(i_{0}(\varphi))= \mathbf{A}\big(\vartheta_{0}(\varphi),\, I_{0}(\varphi),z_0(\varphi)\big)\\
&=v\big(\vartheta_{0}(\varphi),I_0(\varphi)\big)+z_{0}(\varphi),
\end{align*}
supplemented with the reversibility assumption
\begin{equation*}
r_0(-\varphi,-\theta)=r_0(\varphi,\theta)
\end{equation*}
and  $\mathcal{R}$ is an integral operator as in \eqref{kernel-phi-lambda} with a kernel $K_0$ satisfying the property
\begin{equation*}
K_0(\lambda,-\varphi,-\theta,-\eta)=K_0(\lambda,\varphi,\theta,\eta).
\end{equation*}
Moreover,   under the assumption
					\begin{equation}\label{frakI0 bnd}
					\|\mathfrak{I}_0\|_{q,s_0}^{\gamma,\kappa}\lesssim 1,	
					\end{equation}
					 we have for all $s\geqslant s_{0}$, 
\begin{enumerate}
\item  The function $r_0$ satisfies the estimates,
\begin{equation*}
\| r_0\|_{s}^{q,\kappa}\lesssim 1+\|\mathfrak{I}_{0}\|_{s}^{q,\kappa}\quad\hbox{and}\quad \|\Delta_{12}r_0\|_{s}^{q,\kappa}\lesssim\|\Delta_{12}i\|_{s}^{q,\kappa}+{\|\Delta_{12}i\|_{s_0}^{q,\kappa}\max_{\ell=1,2}\|\mathfrak{I}_{\ell}\|_{s}^{q,\kappa}}.
\end{equation*}
\item For all $j\in \mathbb{N}$ the kernel $K_0$ satisfies the estimates,
\begin{equation*}
\sup_{\eta\in\mathbb{T}}\|(\partial_{\theta}^{j}K_0)(*,\cdot,\centerdot,\centerdot+\eta)\|_{s}^{q,\kappa}\lesssim 1+\|\mathfrak{I}_{0}\|_{s+3+j}^{q,\kappa}
\end{equation*}
and
\begin{equation*}
\sup_{\eta\in\mathbb{T}}\|\Delta_{12}(\partial_{\theta}^{j}K_0)(*,\cdot,\centerdot,\centerdot+\eta)\|_{s}^{q,\kappa}\lesssim\|\Delta_{12}i\|_{q,s+3+j}+{\|\Delta_{12}i\|_{s_0+3}^{q,\kappa}\max_{\ell=1,2}\|\mathfrak{I}_{\ell}\|_{s+3+j}^{q,\kappa},}
\end{equation*}
where $*,\cdot,\centerdot,$ denote successively the variables $\lambda,\varphi,\theta$ and $\displaystyle\mathfrak{I}_{\ell}(\varphi)=i_{\ell}(\varphi)-(\varphi,0,0).$ In addition,  for any function $f$,   $\Delta_{12} f:=f(i_1)-f(i_2)$ refers   for the difference of $f$ taken  at two different states $i_1$ and $i_2$ satisfying \eqref{frakI0 bnd}. 
\end{enumerate}} 
\end{proposition}
\begin{proof} To alleviate the notation we shall at  several stages of the proof remove the dependance  of our functions/operators with respect to $(\omega,\alpha)$ and keep it  when  it is relevant. Recall that the operator $\mathcal{L}_{\omega}$ is defined in  \eqref{Norm-proj}. 
Then we need to describe $K_{02}(\varphi)$. According to Proposition \ref{Prop-Conjugat}-{\rm (i)} and the identity \eqref{H alpha}  we may write
\begin{equation}\label{Y-2}
K_{02}(\varphi)  =  \mathrm{L}(\lambda)|_{ \mathbb{H}^{\bot}_{\mathbb{S}_0}}+\varepsilon\partial_{z}\nabla_{z}\mathcal{P}_{\varepsilon}(i_{0}(\varphi))+\varepsilon\mathcal{R}(\varphi),
\end{equation}
where 
\begin{align*}
\mathcal{R}(\varphi)\triangleq\mathcal{R}_{1}(\varphi)+\mathcal{R}_{2}(\varphi)+\mathcal{R}_{3}(\varphi),\quad \mbox{ with }\quad
\mathcal{R}_{1}(\varphi)&=L_{2}^{\top}(\varphi)\partial_{I}\nabla_I\mathcal{P}_{\varepsilon}(i_{0}(\varphi))L_{2}(\varphi),\\ \mathcal{R}_{2}(\varphi)&=L_{2}^{\top}((\varphi)\partial_{z}\nabla_{I}\mathcal{P}_{\varepsilon}(i_{0}(\varphi)),\\
\mathcal{R}_{3}(\varphi)&=\partial_{I}\nabla_{z}\mathcal{P}_{\varepsilon}(i_{0}(\varphi))L_{2}(\varphi).
\end{align*}
As we shall see  all the operators  $\mathcal{R}_{1}(\varphi),$ $\mathcal{R}_{2}(\varphi)$ and  $\mathcal{R}_{3}(\varphi)$   have  finite-dimensional rank. This property is obvious  for  the operator $L_{2}(\varphi)$ which send $ \mathbb{H}^{\bot}_{\mathbb{S}_0}$ to $\mathbb{R}^{d}$ and therefore for any $\rho\in  \mathbb{H}^{\bot}_{\mathbb{S}_0}$
$$L_{2}(\varphi)[\rho]=\sum_{k=1}^{d}\big\langle L_{2}(\varphi)[\rho],\underline{e}_{k}\big\rangle_{L^{2}(\mathbb{T})}\,\underline{e}_{k}=\sum_{k=1}^{d}\big\langle\rho,L_{2}^{\top}(\varphi)[\underline{e}_{k}]\big\rangle_{L^{2}(\mathbb{T})}\,\underline{e}_{k},$$
with $\displaystyle(\underline{e}_{k})_{k=1}^d$ being the canonical basis of $\mathbb{R}^d.$
Hence
$$\begin{array}{ll}
\displaystyle \mathcal{R}_{1}(\varphi)[\rho]=\sum_{k=1}^{d}\big\langle\rho,L_{2}^{\top}(\varphi)[\underline{e}_{k}]\big\rangle_{L^{2}(\mathbb{T})}\mathtt{A}_{1}(\varphi)[\underline{e}_{k}] & \quad\mbox{with }\quad \mathtt{A}_{1}(\varphi)=L_{2}^{\top}(\varphi)\partial_{I}\nabla_I\mathcal{P}_{\varepsilon}(i_{0}(\varphi)),\\
\displaystyle \mathcal{R}_{3}(\varphi)[\rho]=\sum_{k=1}^{d}\big\langle\rho,L_{2}^{\top}(\varphi)[\underline{e}_{k}]\big\rangle_{L^{2}(\mathbb{T})}\mathtt{A}_{3}(\varphi)[\underline{e}_{k}] &\quad \mbox{with }\quad \mathtt{A}_{3}(\varphi)=\partial_{I}\nabla_{z}\mathcal{P}_{\varepsilon}(i_{0}(\varphi)).
\end{array}$$
In the same spirit, since $\mathtt{A}_{2}(\varphi)=\partial_{z}\nabla_{I}\mathcal{P}_{\varepsilon}(i_{0}(\varphi)): \mathbb{H}^{\bot}_{\mathbb{S}_0}\rightarrow\mathbb{R}^{d},$ then we can write
$$
\mathcal{R}_{2}(\varphi)[\rho]=\sum_{k=1}^{d}\big\langle\rho,\mathtt{A}_{2}^{\top}(\varphi)[\underline{e}_{k}]\big\rangle_{L^{2}\,(\mathbb{T})}L_{2}^{\top}(\varphi)[\underline{e}_{k}].
$$
By setting
$$
g_{k,1}(\varphi,\theta)=g_{k,3}(\varphi,\theta)=\chi_{k,2}(\varphi,\theta)\triangleq L_{2}^{\top}(\varphi)[\underline{e}_{k}](\theta),\quad g_{k,2}(\varphi,\theta)\triangleq \mathtt{A}_{2}^{\top}(\varphi)[\underline{e}_{k}](\theta)
$$
and
$$
\mbox{ }\chi_{k,1}(\varphi,\theta)\triangleq \mathtt{A}_{1}(\varphi)[\underline{e}_{k}](\theta),\quad \chi_{k,3}(\varphi,\theta)\triangleq \mathtt{A}_{3}(\varphi)[\underline{e}_{k}](\theta),
$$
one can  write the operator $\mathcal{R}$ in the integral form
\begin{align*}\mathcal{R}\rho(\varphi,\theta)&=\sum_{k'=1}^{3}\sum_{k=1}^{d}\langle\rho(\varphi,\cdot),g_{k,k'}(\varphi,\cdot)\rangle_{L^{2}(\mathbb{T})}\chi_{k,k'}(\varphi,\theta)\\
&=\int_0^{2\pi}\rho(\varphi,\eta)K_0(\varphi,\theta,\eta)d\eta
\end{align*}
with 
$$K_0(\varphi,\theta,\eta)\triangleq \sum_{k'=1}^{3}\sum_{k=1}^{d}g_{k,k'}(\varphi,\eta)\chi_{k,k'}(\varphi,\theta).$$
Now we remark that  $g_{k,k'},\chi_{k,k'}\in \mathbb{H}^{\bot}_{\mathbb{S}_0}$ and satisfy the estimates 
\begin{equation}\label{estimate gkk' and chikk'}
\| g_{k,k'}\|_{s}^{q,\kappa}+\|\chi_{k,k'}\|_{s}^{q,\kappa}\lesssim 1+\|\mathfrak{I}_{0}\|_{s+3}^{q,\kappa}
\end{equation}
and 
\begin{equation}\label{estimate differential gkk' and chikk'}
\| d_{i}g_{k,k'}[\widehat{i}]\|_{s}^{q,\kappa}+\| d_{i}\chi_{k,k'}[\widehat{i}]\|_{s}^{q,\kappa}\lesssim\|\widehat{i}\|_{s+2}^{q,\kappa}+\|\mathfrak{I}_{0}\|_{s+4}^{q,\kappa}\|\widehat{i}\|_{s_{0}+2}^{q,\kappa}.
\end{equation}
By \eqref{estimate gkk' and chikk'}, combined with the law products of Lemma \ref{Law-prodX1} we deduce that 
$$\sup_{\eta\in\mathbb{T}}\| (\partial_{\theta}^{j}K_0)(\lambda,\omega,\varphi,\theta,\eta+\theta)\|_{s}^{q,\kappa}\lesssim 1+\|\mathfrak{I}_{0}\|_{s+3+j}^{q,\kappa}.$$
Applying  the mean value theorem with \eqref{estimate differential gkk' and chikk'}, yields
$$\sup_{\eta\in\mathbb{T}}\|\Delta_{12}(\partial_{\theta}^{j}K_0)(\lambda,\omega,\varphi,\theta,\eta+\theta)\|_{s}^{q,\kappa}\lesssim\|\Delta_{12}i\|_{s+3+j}^{q,\kappa}+\|\Delta_{12}i\|_{s_0+3}^{q,\kappa}\max_{\ell=1,2}\|\mathfrak{I}_{\ell}\|_{s+3+j}^{q,\kappa}.$$
This concludes the proof of {\rm (ii)}.
The symmetry property of $K_0$ is a consequence of the definition of $r$ and the reversibility condition \eqref{parity solution} imposed on the torus $i_{0}.$ Then,  using \eqref{Y-2},  \eqref{cNP} and  \eqref{formaHep} we get
$$\begin{array}{rcl}
K_{02}(\varphi) & = & \mathrm{L}(\lambda)|_{ \mathbb{H}^{\bot}_{\mathbb{S}_0}}+\varepsilon \partial_{z}\nabla_{z}\mathcal{P}_{\varepsilon}(i_{0}(\varphi))+\varepsilon\mathcal{R}(\varphi)
\\
& = & \mathrm{L}(\lambda)|_{ \mathbb{H}^{\bot}_{\mathbb{S}_0}}+\varepsilon  \Pi_{\mathbb{S}_0}^\bot\partial_{r}\nabla_{r}P_{\varepsilon}( \mathbf{A}(i_{0}(\varphi)))+\varepsilon\mathcal{R}(\varphi)
\\
& = & \Pi_{\mathbb{S}_0}^{\perp}\partial_{r}\nabla_{r}\mathcal{H}_{\varepsilon}( \mathbf{A}(i_{0}(\varphi)))+\varepsilon\mathcal{R}(\varphi)
\\
& = & \Pi_{\mathbb{S}_0}^{\perp}\partial_{r}\nabla_{r}H(\varepsilon r_0(\varphi,\cdot))+\varepsilon\mathcal{R}(\varphi)
\end{array}$$
where we have used the notation $r_0(\varphi,\cdot)\triangleq  \mathbf{A}(i_{0}(\varphi))$. According to the general form of the linearized operator stated  in Proposition \ref{lin-eq-r} one has
$$-\partial_{\theta}\partial_{r}\nabla_{r}H(\varepsilon r_0(\varphi,\cdot))=\partial_{\theta}\left(V_{\varepsilon r_0,\alpha}\cdot\right)-\partial_{\theta}\mathbb{K}_{\varepsilon r_0,\alpha}$$
which implies 
$$\begin{array}{rcl}
-\partial_\theta K_{02}(\varphi) 
= \Pi_{\mathbb{S}_0}^{\perp}\big(\partial_{\theta}\left(V_{\varepsilon r_0,\alpha}\cdot\right)-\partial_{\theta}\mathbb{K}_{\varepsilon r_0,\alpha})-\varepsilon \partial_\theta \mathcal{R}(\varphi)\big)\Pi_{\mathbb{S}_0}^{\perp}.
\end{array}$$
Plugging this identity into \eqref{Norm-proj} gives the desired result.
{Then by  \eqref{aacoordinates} and  \eqref{v-theta1} we get
\begin{align*}
\nonumber
\| r_0\|_{s}^{q,\kappa}&\lesssim \| v(\vartheta,I_0)\|_{s}^{q,\kappa}+\| z\|_{s}^{q,\kappa}\\
&\lesssim 1+\|\mathfrak{I}_{0}\|_{s}^{q,\kappa}.
\end{align*}
The bound on $\Delta_{12}r_0$ can be done in a similar way using in particular the mean value theorem. Indeed, one may write 
\begin{align*}
\nonumber
\| \Delta_{12}r_0\|_{s}^{q,\kappa}&\lesssim \| \Delta_{12}v(\vartheta,I_0)\|_{s}^{q,\kappa}+\| \Delta_{12}z\|_{s}^{q,\kappa}.
\end{align*}
Applying Taylor formula with \eqref{v-theta1} combined with the law products  allow to get 
\begin{align*}
\| \Delta_{12}v(\vartheta,I_0)\|_{s}^{q,\kappa}&\lesssim\| \Delta_{12}(\vartheta,I_0)\|_{s}^{q,\kappa}+\| \Delta_{12}(\vartheta,I_0)\|_{s_0}^{q,\kappa}\max_{\ell=1,2}\|\mathfrak{I}_{\ell}\|_{s}^{q,\kappa}\\
&\lesssim  \|\Delta_{12}i\|_{s}^{q,\kappa}+\|\Delta_{12}i\|_{s_0}^{q,\kappa}\max_{\ell=1,2}\|\mathfrak{I}_{\ell}\|_{s}^{q,\kappa},
\end{align*}
which implies in turn
\begin{align*}
\nonumber
\| \Delta_{12}r_0\|_{s}^{q,\kappa}&\lesssim\|\Delta_{12}i\|_{s}^{q,\kappa}+\|\Delta_{12}i\|_{s_0}^{q,\kappa}\max_{\ell=1,2}\|\mathfrak{I}_{\ell}\|_{s}^{q,\kappa}.
\end{align*}
}
This completes the proof of Proposition \ref{lemma-GS0}.
\end{proof}
\subsection{Asymptotic structure}
We intend here to provide  a suitable decomposition for the linearized operator described  by \eqref{lin} for a small  state $r$ and estimate its coefficients.  For the function spaces tools  used below we refer to Section \ref{Sec-F-Spaces1} and Section \ref{Se-Toroidal pseudo-differential operators}. However for the reversibility concepts we may consult Definition \ref{Def-Rev}. Our main result of this section reads as follows. 
\begin{lemma}\label{lemma-reste}
Let $q,\kappa\in \N, s_0>\frac{d+1}{2}+q$, there exists $\varepsilon_0>0$ such that  for   $r\in W^{q,\infty}_{\kappa}(\mathcal{O},H^{\infty}_{\textnormal{even}}(\T^{d+1}))$ satisfying 
\begin{equation}\label{small-C1}
\|r\|_{s_0}^{q,\kappa}\leqslant\varepsilon_0,
\end{equation}
 the  operator $\mathcal{L}_{r,\lambda}$ in \eqref{lin} is reversible and takes the form
 \begin{align*}
 \mathcal{L}_{r,\lambda}&=\omega\cdot\partial_\varphi+\partial_\theta\Big[ V_{r,\alpha}-\big(\mathscr{W}_{r,\alpha}\, |{\textnormal D}|^{\alpha-1}+|{\textnormal D}|^{\alpha-1}\mathscr{W}_{r,\alpha}\big)+\mathscr{R}_{r,\alpha}\Big]
\end{align*}
with the following properties.
\begin{enumerate}
\item  The functions  $\mathscr{W}_{r,\alpha}, V_{r,\alpha}\in W^{q,\infty}_{\kappa}(\mathcal{O},H^{\infty}_{\textnormal{even}}(\T^{d+1}))$ such that  for any $s\geqslant s_0$
$$
 \|\mathscr{W}_{r,\alpha}-C_\alpha 2^{-\alpha}\|_{s+3}^{q,\kappa}+\|\partial_\theta \mathscr{W}_{r,\alpha}\|_{s+2}^{q,\kappa}+\|V_{r,\alpha}-V_{0,\alpha}\|_{s}^{q,\kappa}\lesssim \|r\|_{s+4}^{q,\kappa}$$
where $V_{0,\alpha}$ is defined in \eqref{defTheta}. 
\item   The linear operator 
$\partial_\theta\mathscr{R}_{r,\alpha}:  W^{q,\infty}_{\kappa}(\mathcal{O},H^s_{\textnormal{even}}(\T^{d+1}))\to  W^{q,\infty}_{\kappa}(\mathcal{O}, H^s_{\textnormal{odd}}(\T^{d+1})$ is  continuous with 
$$
\interleave \partial_\theta\mathscr{R}_{r,\alpha}\interleave_{-1,s,\gamma}^{q,\kappa}\lesssim  \|r\|_{s+5+\gamma}^{q,\kappa} .
$$
\item Let $r_1,r_2$ be smooth and satisfying \eqref{small-C1}. Denote by $\Delta_{12}V_{r,\alpha}=V_{r_1,\alpha}-V_{r_2,\alpha}$ and $\Delta_{12}r=r_1-r_2$. Then  for any $s\geqslant s_0$, we have 
$$
\|\Delta_{12}V_{r,\alpha}\|_{s}^{q,\kappa}\lesssim \|\Delta_{12}r\|_{s+4}^{q,\kappa}+{ \|\Delta_{12}r\|_{s_0+4}^{q,\kappa}
\max_{\ell=1,2} \|r_\ell\|_{s+4}^{q,\kappa}.}
$$
\end{enumerate}
\end{lemma}
\begin{proof}
{\bf{$($i$)$}}  We shall first establish the operator decomposition and move later to the estimates of the coefficients. For this aim, we first  use \eqref{A-rM} in order to get
\begin{align*}
A_r(\varphi,\theta,\eta)&=\sin^2\left(\frac{\eta-\theta}{2}\right)\left[\left(\frac{R(\varphi,\eta)-R(\varphi,\theta)}{\sin\big(\frac{\eta-\theta}{2}\big)}\right)^2+4R(\varphi,\eta)R(\varphi,\theta)\right]. 
\end{align*} 
Introduce the function 
$$
f(\varphi,\theta,\eta)=C_\alpha\left[\left(\frac{R(\varphi,\eta)-R(\varphi,\theta)}{\sin\big(\frac{\eta-\theta}{2}\big)}\right)^2+4R(\varphi,\eta)R(\varphi,\theta)\right]^{-\frac\alpha2}.
$$
Then it is clear that  $f$ is symmetric, that is, $f(\varphi,\theta,\eta)=f(\varphi,\eta,\theta)$ and therefore  by  using \mbox{Lemma \ref{Lemm-RegZ}-$($ii$)$} we deduce that
\begin{equation}\label{W-eq-PP}
f(\varphi,\theta,\eta)=\mathscr{W}_{r,\alpha}(\varphi,\theta)+\mathscr{W}_{r,\alpha}(\varphi,\eta)+\sin^2\left(\frac{\eta-\theta}{2}\right)\mathscr{A}_{r,\alpha}(\varphi,\theta,\eta)
\end{equation}
with 
\begin{equation}\label{W-eq}
\mathscr{W}_{r,\alpha}(\varphi,\theta)=C_\alpha2^{-\alpha-1}\Big(R^2(\varphi,\theta)+\big(\partial_\theta R\big)^2(\varphi,\theta)\Big)^{-\frac\alpha2},
\end{equation}
where
$\mathscr{A}_{r,\alpha}$ is symmetric  and satisfies   the estimate
 \begin{align}\label{ScrA0}
\nonumber\sup_{\eta\in\T}\|\mathscr{A}_{r,\alpha}(\cdot,\centerdot, \centerdot+\eta)\|_{s}^{q,\kappa}\leqslant C\|\Delta_{\theta,\eta}f\|_{s+\frac12+\epsilon}^{q,\kappa}\\
\leqslant C\|f-2^{-\alpha}C_\alpha\|_{s+\frac52+\epsilon}^{q,\kappa}.
\end{align}
Notice that the notation $\cdot$ is for $\varphi$ and $\centerdot$ stands for $\theta$.
From the translation invariance of Lebesgue measure  we may write
$$
\|f-2^{-\alpha}C_\alpha\|_{s+\frac52+\epsilon}^{q,\kappa}=\|\tilde{f}-2^{-\alpha}C_\alpha\|_{s+\frac52+\epsilon}^{q,\kappa}
$$
with
\begin{align*}
\tilde{f}(\varphi,\theta,\eta)&=f(\varphi,\theta,\theta+\eta)\\
&=C_\alpha\left[\left(\frac{R(\varphi,\theta+\eta)-R(\varphi,\theta)}{\sin\big(\frac{\eta}{2}\big)}\right)^2+4R(\varphi,\theta+\eta)R(\varphi,\theta)\right]^{-\frac\alpha2}.
\end{align*}
Applying Lemmata \ref{Compos-lemm}-\ref{Law-prodX1}-\ref{lem-Reg1} combined with   the smallness condition \eqref{small-C1},
we obtain
\begin{align}\label{ScrA1}
\|\tilde{f}-2^{-\alpha}C_\alpha\|_{s+\frac52+\epsilon}^{q,\kappa}
&\leqslant C\|r\|_{s+\frac72+\epsilon}^{q,\kappa}.
\end{align}
Notice that in order to apply Lemma \ref{lem-Reg1}  it is enough to  write
$$
\left(\frac{R(\varphi,\theta+\eta)-R(\varphi,\theta)}{\sin\big(\frac{\eta}{2}\big)}\right)^2=\left(\frac{R(\varphi,\theta+\eta)-R(\varphi,\theta)}{\tan\big(\frac{\eta}{2}\big)}\right)^2+\left(R(\varphi,\theta+\eta)-R(\varphi,\theta)\right)^2.
$$
It follows from \eqref{ScrA0} and \eqref{ScrA1} that
\begin{align}\label{ScrA}
\sup_{\eta\in\T}\|\mathscr{A}_{r,\alpha}(\cdot,\centerdot, \centerdot+\eta)\|_{s}^{q,\kappa}&\leqslant C\|r\|_{s+\frac72+\epsilon}^{q,\kappa}.
\end{align}
 Similarly,  using Lemma \ref{Lemm-RegZ}-$($ii$)$ we get for any $p,m\in\N$
 \begin{align}\label{ScrAbis}
\sup_{\eta\in\T}\|(\partial_\theta^p\partial_\eta^m\mathscr{A}_{r,\alpha})(\cdot,\centerdot, \centerdot+\eta)\|_{s}^{q,\kappa}&\leqslant C\|r\|_{s+p+m+\frac72+\epsilon}^{q,\kappa}.
\end{align}
Coming back to \eqref{W-eq-PP}, then we have  the following decomposition
\begin{align}\label{A_r}
C_\alpha A_r^{-\frac\alpha2}(\varphi,\theta,\eta)=&\frac{\mathscr{W}_{r,\alpha}(\varphi,\theta)+\mathscr{W}_{r,\alpha}(\varphi,\eta)}{|\sin(\frac{\theta-\eta}{2})|^\alpha}+\frac{\mathscr{A}_{r,\alpha}(\varphi,\theta,\eta)}{|\sin(\frac{\eta-\theta}{2})|^{\alpha-2}}\cdot
\end{align}
Plugging this identity into  \eqref{V-eq} and using a change of variables  we find
\begin{equation}\label{V-eq0}
V_{r,\alpha}(\varphi,\theta)=\frac{ 1}{2\pi }\bigintsss_{0}^{2\pi} \frac{\mathscr{A}_{r,\alpha}^1(\varphi,\theta,\eta)}{|\sin(\frac\eta2)|^{\alpha}} d\eta
\end{equation}
with
\begin{align}\label{A_r1}
\mathscr{A}_{r,\alpha}^1(\varphi,\theta,\eta)&=\frac{\partial_{\eta} \big[R(\varphi,\eta+\theta)\sin(\eta)\big]}{R(\varphi,\theta)}\Big(\mathscr{W}_{r,\alpha}(\varphi,\theta)+\mathscr{W}_{r,\alpha}(\varphi,\eta+\theta)\Big)\\
\nonumber&+\frac{\partial_{\eta} \big[R(\varphi,\eta+\theta)\sin(\eta)\big]}{R(\varphi,\theta)}\sin^2({\eta}/{2})\mathscr{A}_{r,\alpha}(\varphi,\theta,\eta+\theta).
\end{align}
Similarly inserting \eqref{A_r} into \eqref{K-r} yields
\begin{align*}
\mathbb{K}_{r,\alpha} h (\varphi,\theta)&= \frac{ 1}{2\pi }\bigintsss_{0}^{2\pi} \frac{h(\varphi,\eta)\left(\mathscr{W}_{r,\alpha}(\varphi,\theta)+\mathscr{W}_{r,\alpha}(\varphi,\eta)\right)}{|\sin(\frac{\eta-\theta}{2})|^{\alpha}} d\eta\\
&+\frac{1 }{2\pi }\bigintsss_{0}^{2\pi} \frac{h(\varphi,\eta)\mathscr{A}_{r,\alpha}(\varphi,\theta,\eta)}{|\sin(\frac{\eta-\theta}{2})|^{\alpha-2}} d\eta.
\end{align*}
Using the definition of the modified fractional Laplacian in \eqref{fract1} we infer
\begin{align*}
\mathbb{K}_{r,\alpha} h (\varphi,\theta)&= \mathscr{W}_{r,\alpha}(\varphi,\theta)|\textnormal{D}|^{\alpha-1}h(\varphi,\theta)+|\textnormal{D}|^{\alpha-1}(\mathscr{W}_{r,\alpha} h)(\varphi,\theta)+\mathscr{R}_{r,\alpha}h(\varphi,\theta)
\end{align*}
with
\begin{align}\label{K2}
\mathscr{R}_{r,\alpha}h(\varphi,\theta)\triangleq
&\frac{ 1 }{2\pi }\bigintsss_{0}^{2\pi} \frac{\mathscr{A}_{r,\alpha}(\varphi,\theta,\eta)}{|\sin(\frac{\eta-\theta}{2})|^{\alpha-2}} h(\varphi,\eta)d\eta.
\end{align}
Direct computations yield
\begin{align}\label{linXX0}
\partial_\theta\mathscr{R}_{r,\alpha}h(\varphi,\theta)=
&\frac{ 1 }{2\pi }\bigintsss_{0}^{2\pi} \frac{ \sin(\frac{\theta-\eta}{2})}{|\sin(\frac{\eta-\theta}{2})|^{\alpha}}\mathscr{A}^2_{r,\alpha}(\varphi,\theta,\eta)\, h(\varphi,\eta) d\eta
\end{align}
with
\begin{equation}\label{Ar2}
\mathscr{A}_{r,\alpha}^2(\varphi,\theta,\eta)=\partial_\theta\mathscr{A}_{r,\alpha}(\varphi,\theta,\eta)\,\sin\left(\frac{\theta-\eta}{2}\right)+1-\frac{\alpha}{2}\mathscr{A}_{r,\alpha}(\varphi,\theta,\eta)\,\cos\left(\frac{\theta-\eta}{2}\right).
\end{equation}
Inserting this into \eqref{lin} we find
\begin{align}\label{linXX1}
 \mathcal{L}_{r,\lambda}&=\omega\cdot\partial_\varphi+\partial_\theta\Big[ V_{r,\alpha}-\big(\mathscr{W}_{r,\alpha}\, |{\textnormal D}|^{\alpha-1}+|{\textnormal D}|^{\alpha-1}\mathscr{W}_{r,\alpha}\big)+\mathscr{R}_{r,\alpha}\Big]
\end{align}
where $V_{r,\alpha}$ is given by \eqref{V-eq0} and  $\mathscr{R}_{r,\alpha}$  by \eqref{K2}. The next goal is to estimate the coefficients $\mathscr{W}_{r,\alpha}$ and $V_{r,\alpha}$. Concerning the first one, defined in  \eqref{W-eq} one has
$$
\mathscr{W}_{r,\alpha}(\varphi,\theta)=C_\alpha 2^{-\alpha-1}\left(1+2r(\varphi,\theta)+\frac{(\partial_\theta r)^2(\varphi,\theta)}{1+2r(\varphi,\theta)}\right)^{-\frac\alpha2}.
$$
Then from Lemma \ref{Compos-lemm} combined with  the smallness condition \eqref{small-C1} we obtain
\begin{equation}\label{WrC}
\|\mathscr{W}_{r,\alpha}-C_\alpha 2^{-\alpha-1}\|_{{s-1}}^{q,\kappa}\leqslant C\|r\|_{{s}}^{q,\kappa}.
\end{equation}
Similarly, computing $\partial_\theta \mathscr{W}_{r,\alpha}$ and using the  law products together  with  Lemma \ref{Compos-lemm} we find 
$$
\|\partial_\theta \mathscr{W}_{r,\alpha}\|_{{s-2}}^{q,\kappa}\leqslant C\|r\|_{{s}}^{q,\kappa}.
$$
Let us now move to the estimate of $V_{r,\alpha}$ defined by  \eqref{V-eq0} then  we find from straightforward computations
\begin{align*}
\|V_{r,\alpha}-V_{0,\alpha}\|_{s}^{q,\kappa}&\lesssim\bigintsss_{0}^{2\pi} \frac{\|\mathscr{A}_{r,\alpha}^1-\mathscr{A}_{0,\alpha}^1\|_{s}^{q,\kappa}}{|\sin(\frac\eta2)|^{\frac12}} \left(1+|\ln(\sin(\eta/2))|^q\right) d\eta\\
&\lesssim \|\mathscr{A}_{r,\alpha}^1-\mathscr{A}_{0,\alpha}^1\|_{s}^{q,\kappa}.
\end{align*}
Notice that $V_{0,\alpha}=V_0(\alpha)$ which is given by \eqref{defTheta}. 
Applying Lemma \ref{Law-prodX1} to \eqref{A_r1}   combined with  \eqref{ScrA}, \eqref{WrC} and  \eqref{small-C1} yield
\begin{align}\label{VVV-PP}
\|V_{r,\alpha}-V_0(\alpha)\|_{s}^{q,\kappa}\le  C\|r\|_{s+\frac72+\epsilon}^{q,\kappa}.
\end{align}
To get the suitable estimate it suffices to use Sobolev embeddings.\\
{\bf{$($ii$)$}}
Coming back to \eqref{linXX0} and making the  change of variables $\eta\leadsto \eta+\theta$ we find
\begin{align*}
\partial_\theta\mathscr{R}_{r,\alpha}h(\varphi,\theta)=
&\frac{ 1 }{2\pi }\bigintsss_{0}^{2\pi} \frac{ \sin(\frac{\eta}{2})}{|\sin({\eta}/{2})|^{\alpha}}\mathscr{A}^2_{r,\alpha}(\varphi,\theta,\eta+\theta)\, h(\varphi,\eta+\theta) d\eta,
\end{align*}
with
\begin{equation}\label{Ar2L}
\mathscr{A}_{r,\alpha}^2(\varphi,\theta,\eta+\theta)=-\partial_\theta\mathscr{A}_{r,\alpha}(\varphi,\theta,\eta+\theta)\,\sin\left({\eta}/{2}\right)+(1-\frac\alpha2)\mathscr{A}_{r,\alpha}(\varphi,\theta,\eta+\theta)\,\cos\left({\eta}/{2}\right).
\end{equation}
According to \eqref{symb-kern}, the  symbol of $\partial_\theta\mathscr{R}_{r,\alpha}$ can be recovered from the kernel as follows
\begin{align*}
\sigma_{\partial_\theta\mathscr{R}_{r,\alpha}}(\varphi,\theta,\xi)&=\frac{ 1 }{2\pi }\bigintsss_{0}^{2\pi} \frac{ \sin(\frac{\eta}{2})}{|\sin({\eta}/{2})|^{\alpha}}\mathscr{A}^2_{r,\alpha}(\varphi,\theta,\eta+\theta)\, e^{\ii \eta\xi}d\eta.
\end{align*}
Then integration by parts combined with \eqref{It-action} yield
\begin{align*}
\xi^{1+\gamma}\Delta_\xi^\gamma\sigma_{\partial_\theta\mathscr{R}_{r,\alpha}}(\varphi,\theta,\xi)&=\frac{ \ii^{1+\gamma} }{2\pi }\bigintsss_{0}^{2\pi} e^{\ii \eta\xi}\partial_\eta^{1+\gamma} \left[\frac{ \sin(\frac{\eta}{2})(e^{\ii \eta}-1)^\gamma}{|\sin({\eta}/{2})|^{\alpha}}\mathscr{A}^2_{r,\alpha}(\varphi,\theta,\eta+\theta)\right]\, d\eta.
\end{align*}
Therefore, Leibniz rule  gives 
\begin{align*}
\langle \xi\rangle^{1+\gamma}\big\|\Delta_\xi^\gamma\sigma_{\partial_\theta\mathscr{R}_{r,\alpha}}(\cdot,\centerdot,\xi)\big\|_{H^s}&\lesssim\sum_{0\leqslant\beta\leqslant \gamma+1}\bigintsss_{0}^{2\pi} \frac{\| \partial_\eta^{\beta}\mathscr{A}_{r,\alpha}^2(\cdot,\centerdot,\centerdot+\eta)\|_{H^s}}{|\sin({\eta}/{2})|^{\alpha}}\, d\eta.
\end{align*}
Similarly we get for any $0\leqslant |j|\leqslant q$ and $\alpha\in(0,\overline\alpha)$ and by Sobolev embeddings
\begin{align*}
\kappa^{|j|}\langle \xi\rangle^{1+\gamma}\big\|\partial_\lambda^j\Delta_\xi^\gamma\sigma_{\partial_\theta\mathscr{R}_{r,\alpha}}(\cdot,\centerdot,\xi)\big\|_{H^{s-|j|}}&\lesssim\sum_{0\leqslant\beta\leqslant \gamma+1\atop 0\leqslant j'\leqslant j}\kappa^{|j'|}\bigintsss_{0}^{2\pi} \frac{\| \partial_\lambda^{j'}\partial_\eta^{\beta}\mathscr{A}_{r,\alpha}^2(\cdot,\centerdot,\centerdot+\eta)\|_{H^{s-|j'|}}}{|\sin({\eta}/{2})|^{\frac12}}\, d\eta\\
&\lesssim \sum_{0\leqslant\beta\leqslant \gamma+1}\bigintsss_{0}^{2\pi} \frac{\| \partial_\eta^{\beta}\mathscr{A}_{r,\alpha}^2(\cdot,\centerdot,\centerdot+\eta)\|_{s}^{q,\kappa}}{|\sin({\eta}/{2})|^{\frac12}}\, d\eta.
\end{align*}
By virtue of \eqref{Ar2} and \eqref{ScrAbis} we get for any $0\leqslant\beta\leqslant 1+\gamma$
\begin{align}\label{Bas1}
\nonumber \|  \partial_\eta^{\beta} \mathscr{A}_{r,\alpha}^2(\cdot,\centerdot,\centerdot+\eta)\|_{s}^{q,\kappa}\leqslant& C\|r\|_{s+\frac72+\beta+\epsilon}^{q,\kappa}\\
\leqslant& C\|r\|_{s+5+\gamma}^{q,\kappa}.
\end{align}
Using the definition \eqref{Def-pseud-w} we obtain
$$
\interleave \partial_\theta\mathscr{R}_{r,\alpha}\interleave_{-1,s,\gamma}^{q,\kappa}\lesssim  \|r\|_{s+5+\gamma}^{q,\kappa}.$$
The symmetry of the functions  $V_{r,\alpha}$ and $\mathscr{W}_{r,\alpha}$ follows from the expressions \eqref{W-eq}, \eqref{V-eq0}. More precisely, one gets that these functions are even, that is,  $V_{r,\alpha}(-\varphi,-\theta)=V_{r,\alpha}(\varphi,\theta)$ with the same property for $\mathscr{W}_{r,\alpha}$. The reversibility of the operator $ \partial_\theta\mathscr{R}_{r,\alpha}$ in the sense of the  Definition \ref{Def-Rev} is quite similar and one gets in view the the structures \eqref{linXX0} and\eqref{Ar2},
$$
\partial_\theta\mathscr{R}_{r,\alpha}\circ\mathcal{S}=-\mathcal{S}\circ\partial_\theta\mathscr{R}_{r,\alpha}
$$
In particular we get for $ h\in H^s_{\textnormal{even}}(\T^{d+1})$ that $\partial_\theta\mathscr{R}_{r,\alpha}h\in H^s_{\textnormal{odd}}(\T^{d+1}).$\\
{\bf{(iii)}} The proof is quite similar to \eqref{VVV-PP}. We use the identities \eqref{V-eq0} and \eqref{A_r1} combined with the estimates
\begin{equation}\label{WrCHH}
\|\Delta_{12}\mathscr{W}_{r,\alpha}\|_{{s-1}}^{q,\kappa}\lesssim\|\Delta_{12}r\|_{{s}}^{q,\kappa}+\|\Delta_{12}r\|_{{s_0}}^{q,\kappa}\max_{j=1,2}\|r_j\|_{{s}}^{q,\kappa}
\end{equation}
 and
 \begin{align}\label{ScrAHH}
\sup_{\eta\in\T}\|\Delta_{12}\mathscr{A}_{r,\alpha}(\cdot,\centerdot, \centerdot+\eta)\|_{s}^{q,\kappa}&\leqslant C\|\Delta_{12}r\|_{s+\frac72+\epsilon}^{q,\kappa}+\|\Delta_{12}r\|_{{s_0+\frac72+\epsilon}}^{q,\kappa}\max_{j=1,2}\|r_j\|_{{s+\frac72+\epsilon}}^{q,\kappa}.
\end{align}
To get the estimates  \eqref{WrCHH} and \eqref{ScrAHH} we proceed in a similar way to \eqref{WrC} and \eqref{ScrA}.
 This achieves the proof of Lemma \ref{lemma-reste}.
 \end{proof}

\subsection{Reduction of the transport part }\label{sec-transport-1}
In this section, we intend to perform the reduction of the transport part of  the linearized operator $\mathcal{L}_{\varepsilon r}$ described in Proposition  \ref{lemma-GS0}. 
This topic  is now  well-developed  in KAM theory  and has been  implemented  in several papers during the last few years, especially in \cite{Baldi-Montalto21, BCP}. The formal statement  says that up to a conjugation of  the operator $\mathcal{L}_{\varepsilon r_\delta}$ by a quasi-periodic symplectic change of variables $\mathscr{B}$ described in \eqref{mathscrB}  the leading part becomes a transport operator with constant coefficients. In our context, we shall use the same techniques developed in the aforementioned papers in order to get  a slightly different  version that fits with our Cantor sets chosen to be related to the final states in the KAM reduction.
This section is organized as follows. First we discuss elementary results on the invertibility of Fourier multiplier operator  in the presence of a small divisor problem. In the second part we shall deal with the straightening of the transport equation when the coefficients are  varying slowly  around constant coefficients.

\subsubsection{Transport equation with constant coefficients}
Take  two constants $\kappa,\,\varrho\in(0,1],\,\tau_1>0$ and let   $\lambda=(\omega,\alpha)\in\mathcal{O}\mapsto c_{\lambda} \in\RR$ be   a smooth  function. We introduce  the Cantor set  
 $$
 \mathtt{C}^{\kappa,\tau_1}=\left\lbrace \lambda=(\omega,\alpha)\in \mathcal{O};\;\,\forall(l,j)\in\mathbb{Z}^{d+1}\backslash\{0\},\;\,|\omega\cdot l+j c_{\lambda}|>\tfrac{\kappa^\varrho}{\langle l\rangle^{\tau_1}}\right\rbrace.$$
 For $N\in\N^\star$ we define the truncated Cantor set
  $$\mathtt{C}^{\kappa,\tau_1}_{N}=\left\lbrace \lambda\in \mathcal{O};\;\; \forall(l,j)\in\mathbb{Z}^{d+1}\backslash\{0\} \;\;\hbox{with}\;\, |l|\leqslant N,\,|\omega\cdot l+jc_{\lambda}|>\tfrac{\kappa^\varrho}{\langle l\rangle^{\tau_1}}\right\rbrace.$$
 Given $f:\mathcal{O}\times\T^{d+1}\to \RR$ a smooth function with zero average, that can be expanded in Fourier series as follows
 $$
 h=\sum_{(l,j)\in\Z^{d+1}\backslash\{0\}}h_{l,j}(\lambda) {\bf{e}}_{l,j},\quad {\bf{e}}_{l,j}(\varphi,\theta)=e^{\ii (l\cdot\varphi+j\theta)}.
 $$
 We want  to solve the transport equation 
\begin{equation}\label{eqq-kk}
 (\omega\cdot\partial_\varphi+c_{\lambda}\partial_\theta) u=h
 \end{equation}
 in the periodic setting, that is, $u:\mathbb{T}^{d+1}\to \RR$.
 Then when $\lambda\in \mathtt{C}^{\kappa,\tau_1}$ we can solve by Fourier series and invert this operator to get
 $$
 u(\lambda,\varphi,\theta)=-\ii\sum_{(l,j)\in\Z^{d+1}\backslash\{0\}}\frac{h_{l,j}(\lambda)}{\omega\cdot l+jc_{\lambda}}{\bf{e}}_{l,j}.
 $$
 We shall give an extension of this formal inverse for the full range of $\lambda\in\mathcal{O}$. To do that, we define the smooth extension of $u$ by
\begin{align}\label{Extend00}
(\omega\cdot\partial_\varphi+c_{r,\lambda}\partial_\theta)_{\textnormal{ext}}^{-1}h\triangleq -\ii\sum_{(l,j)\in\Z^{d+1}\backslash\{0\}}\frac{\chi\big(\omega\cdot l+j c_{\lambda})\kappa^{-\varrho}\langle l\rangle^{\tau_1} \big)h_{l,j}(\lambda)}{\omega\cdot l+j c_{\lambda}}{\bf{e}}_{l,j},
\end{align}
where   the  cut-off function $\chi\in\mathscr{C}^\infty(\RR,\RR)$ is defined in \eqref{chi-def-1}.
The following result is classical and we can refer to \cite{Baldi-berti}, but  for the sake of the completeness we will give some key  ingredients of the proof.
\begin{lemma}\label{L-Invert}
Let $\kappa\in(0,1], q\in\N$ and  $\varrho\in\left(0,\tfrac{1}{q+1}\right]$. There exists $\varepsilon_0\in(0,1)$ such that if \mbox{$\|c_\lambda\|^{q,\kappa}\leqslant\varepsilon_0$,} then for any $ s\geqslant q$ we have
$$
\big\|(\omega\cdot\partial_\varphi+c_{\lambda}\partial_\theta)_{\textnormal{ext}}^{-1}h\big\|_{s}^{q,\kappa}\leqslant C\kappa^{-1} (1+\|c_\lambda\|^{q,\kappa}\big) \|h\|_{s+\tau_1(q+1)}^{q,\kappa}.
$$
In addition, for any $N\in\NN$ and for any $\lambda\in \mathtt{C}^{\kappa,\tau_1}_{N}$ we have
$$
(\omega\cdot\partial_\varphi+c_{\lambda}\partial_\theta)(\omega\cdot\partial_\varphi+c_{\lambda}\partial_\theta)_{\textnormal{ext}}^{-1}\Pi_N=\Pi_N,
$$
where $\Pi_N$ is the orthogonal projection defined by
$$
\Pi_N\sum_{(l,j)\in\Z^{d+1}}h_{l,j} {\bf{e}}_{l,j}=\sum_{(l,j)\in\Z^{d+1}\atop |l|\leqslant N}h_{l,j} {\bf{e}}_{l,j}.
$$
\end{lemma}
\begin{proof}
The proof of the first point can be done using Fa\`a di Bruno's formula in a similar way to  \cite[Lemma 2.5]{Baldi-berti}. As to the identity of the second point, it follows easily from the following observation based on the explicit extension \eqref{Extend00}
$$
(\omega\cdot\partial_\varphi+c_{\lambda}\partial_\theta)_{\textnormal{ext}}^{-1}\Pi_Nh=-\ii\sum_{(l,j)\in\Z^{d+1}\backslash\{0\}\atop |l|\leqslant N}\frac{\chi\big((\omega\cdot l+j c_{\lambda})\kappa^{-\varrho}\langle l\rangle^{\tau_1} \big)h_{l,j}(\lambda)}{\omega\cdot l+jc_{\lambda}}{\bf{e}}_{l,j}.
$$
By construction, one deduces  for $\lambda\in \mathtt{C}^{\kappa,\tau_1}_{N}$  and $|l|\leqslant N$,
$$
\chi\big((\omega\cdot l+j c_{\lambda})\kappa^{-\varrho}\langle l\rangle^{\tau_1} \big)=1,
$$
which implies that
$$
(\omega\cdot\partial_\varphi+c_{\lambda}\partial_\theta)_{\textnormal{ext}}^{-1}\Pi_Nh=-\ii\sum_{(l,j)\in\Z^{d+1}\backslash\{0\}\atop |l|\leqslant N}\frac{h_{l,j}(\lambda)}{\omega\cdot l+jc_{\lambda}}{\bf{e}}_{l,j}.
$$
Therefore, we obtain 
\begin{align*}
(\omega\cdot\partial_\varphi+c_{\lambda}\partial_\theta)(\omega\cdot\partial_\varphi+c_{\lambda}\partial_\theta)_{\textnormal{ext}}^{-1}\Pi_Nh&=\ii\sum_{(l,j)\in\Z^{d+1}\backslash\{0\}\atop |l|\leqslant N}{h_{l,j}(\lambda)}{\bf{e}}_{l,j}\\
&=\Pi_N h.
\end{align*}
This concludes the proof of the lemma.
\end{proof}
\subsubsection{Straightening of the  transport equation}
This section is  devoted to the construction of quasi-periodic change of variables  needed to conjugate the transport part of the linearized operator to a Fourier multiplier. Before stating  our   main result we need to introduce some transformations and explore  some of their basic  properties that can be found for instance in \cite{Baldi-Montalto21,BFM,FGMP19}.
Let $\beta: \mathcal{O}\times \T^{d+1}\to \mathbb{R}^d$ be a smooth function such that $\displaystyle{\sup_{\lambda\in \mathcal{O}}\|\beta(\lambda,\cdot)\|_{\textnormal{Lip}}<1}$ 
then  the mapping
$$(\varphi,\theta)\in\T^{d+1}\mapsto (\varphi, \theta+\beta(\lambda,\varphi,\theta))\in\T^{d+1}
$$ is a diffeomorphism and its  inverse takes the same  form  
 $$\displaystyle(\varphi,\theta)\in\T^{d+1}\mapsto (\varphi, \theta+\widehat{\beta}(\lambda,\omega,\varphi,\theta))\in\T^{d+1}.
 $$  The relation between $\beta$ and $\widehat\beta$ is described through,
 \begin{align}\label{IInv-77}
 y=\theta+\beta(\lambda,\varphi,\theta)\Longleftrightarrow \theta=y+\widehat\beta(\lambda,\varphi,y).
\end{align} 
 Now we define the operators
 \begin{equation}\label{mathscrB}
\mathscr{B}=(1+\partial_{\theta}\beta)\textbf{B},
\end{equation}
with $${\bf{B}}h(\lambda,\varphi,\theta)=h\big(\lambda,\varphi,\theta+\beta(\lambda,\varphi,\theta)\big).$$
Direct computations show that the inverse $\mathscr{B}^{-1}$ keeps the same form, that is, 
\begin{equation}\label{mathscrB1}
 \mathscr{B}^{-1}h(\lambda,\varphi,y)=\Big(1+\partial_y\widehat{\beta}(\lambda,\varphi,y)\Big)h\big(\lambda,\varphi,y+\widehat{\beta}(\lambda,\varphi,y)\big)
\end{equation}
and
\begin{equation}\label{mathscrBB1}
{\bf{B}}^{-1} h(\lambda,\varphi,y)=h\big(\lambda,\omega,\varphi,y+\widehat{\beta}(\lambda,\varphi,y)\big).
\end{equation}
The next algebraic properties follow from straightforward computations.
\begin{lemma}\label{algeb1}
The following assertions holds true
\begin{enumerate}
\item The action of $\mathscr{B}^{-1}$ on the derivative $\partial_\theta$  is given by
\begin{equation*}
\mathscr{B}^{-1}\partial_{\theta}=\partial_{\theta}{\bf{B}}^{-1}
\end{equation*}
\item The conjugation of a transport operator by $\mathscr{B}$  keeps the same structure
$$
 \mathscr{B}^{-1}\big(\omega\cdot\partial_\varphi+\partial_\theta(V(\varphi,\theta)\cdot)\big)\mathscr{B}=\omega\cdot\partial_\varphi+\partial_y\big(\mathscr{V}(\varphi,y)\cdot \big)
$$
with
$$
\mathscr{V}(\varphi,y)={\bf{B}}^{-1}\Big(\omega\cdot\partial_{\varphi} \beta(\varphi,\theta)+V(\varphi,\theta)\big(1+\partial_\theta \beta(\varphi,\theta)\big)\Big)
$$
\item Denote by $\mathscr{B}^\star$ the $L^2_\theta(\T)$-adjoint of $\mathscr{B}$, then
$$
\mathscr{B}^\star={\bf{B}}^{-1}\quad\hbox{and}\quad {\bf{B}}^\star=\mathscr{B}^{-1}.
$$
\end{enumerate}
\end{lemma}

The next result deals with some analytical properties  of the preceding transformations. The proof can be obtained by adapting the proof accomplished in \cite{FGMP19}.
\begin{lemma}\label{Compos1-lemm}
Let $q\in\mathbb{N}, \gamma\in(0,1),s\geqslant s_0>\frac{d+1}{2}+q+1. $ There exists $\varepsilon_0$ small enough such that if
$\|\beta \|_{{2s_0}}^{q,\kappa}\leqslant \varepsilon_0$, then we have the estimates
$$\|{\bf{B}}^{\pm 1}h\|_{s}^{q,\kappa}\leqslant  \|h\|_{s}^{q,\kappa}\big(1+C\|\beta\|_{{s_0}}^{q,\kappa}\big)+C\|\beta\|_{{s}}^{q,\kappa}\|h\|_{s_{0}}^{q,\kappa},
$$
$$\|{\mathscr{B}}^{\pm 1}h\|_{s}^{q,\kappa}\leqslant  \|h\|_{s}^{q,\kappa}\big(1+C\|\beta\|_{{s_0}}^{q,\kappa}\big)+C\|\beta\|_{{s+1}}^{q,\kappa}\|h\|_{s_{0}}^{q,\kappa}
$$
and
\begin{equation*}
\|\widehat{\beta}\|_{s}^{q,\kappa}\lesssim\|\beta\|_{s}^{q,\kappa}\quad\hbox{and}\quad \|{\bf{B}}^{\pm1}h-h\|_{s}^{q,\kappa}\lesssim\|h\|_{s_0}^{q,\kappa}\|\beta\|_{s+{1}}^{q,\kappa}+\|\beta\|_{s_0}^{q,\kappa}\|h\|_{s+1}^{q,\kappa}.
\end{equation*}
Furthermore,  let $\beta_{1},\beta_{2}\in W^{q,\infty,\gamma}(\mathcal{O},H^{\infty}(\mathbb{T}^{d+1}))$ satisfying the foregoing smallness condition and  denote 
$$\Delta_{12}\beta=\beta_{1}-\beta_{2}\quad\textnormal{ and }\quad\Delta_{12}\widehat{\beta}=\widehat{\beta}_{1}-\widehat{\beta}_{2}.
$$
Then, we have the estimate
\begin{equation}\label{link diff beta hat and diff beta}
\|\Delta_{12}\widehat{\beta}\|_{s}^{q,\kappa}\leqslant C\left(\|\Delta_{12}\beta\|_{s}^{q,\kappa}+\|\Delta_{12}\beta\|_{s_{0}}^{q,\kappa}\max_{k\in\{1,2\}}\|\beta_{i}\|_{s+1}^{q,\kappa}\right).
\end{equation}

\end{lemma}

In what follows we intend  to  state the main result of this section  concerning  the reduction of the transport part. 
\begin{proposition}\label{QP-change}
{Let   $\varrho \in\left(0,\frac{1}{q+1}\right]$ and assume \eqref{Conv-T2}. There exists ${\varepsilon}_0>0$ such that with the conditions 
\begin{equation}\label{Conv-T2N}
\begin{array}{ll}
\mu_{2}\geqslant \overline{\mu}_{2}\triangleq 4(\tau_{1}q+\tau_{1})+3, & \quad s_{l}\triangleq s_{0}+\tau_{1}q+\tau_{1}+{2},\\
\overline{s}_{h}\triangleq \frac{3}{2}\overline{\mu}_{2}+s_{l}+1, &\quad s_{h}\triangleq \frac{3}{2}\mu_{2}+s_{l}+1,\\
\qquad \qquad\qquad \qquad \sigma_{1}\triangleq\tau_{1}q+\tau_{1}+7
\end{array}
\end{equation}
 and  
\begin{align}
\label{small-C2}\|\mathfrak{I}_{0}\|_{q,s_{h}+\sigma_1}\leqslant1\quad\textnormal{and}\quad N_{0}^{\mu_{2}}\varepsilon {\kappa^{-2}}\leqslant{\varepsilon}_0,
\end{align}
 there exist
$$\lambda\mapsto c(\lambda,i_0)\in W^{q,\infty,\gamma }(\mathcal{O},\mathbb{R})\quad\mbox{ and }\quad\beta\in \bigcap_{{s\in[s_{0},S]}}W^{q,\infty,\gamma }(\mathcal{O},H_{\textnormal{{odd}}}^{s})$$
such that with $\mathscr{B}$ as in \eqref{mathscrB} one gets the following results.
\begin{enumerate}
\item The function $c(\lambda,i_0)$ satisfies the following estimate,
\begin{equation}\label{est-r1}
\| c(\lambda,i_0)-V_{0,\alpha}\|^{q,\kappa}\lesssim{\varepsilon}
\end{equation}
where $V_{0,\alpha}$ is defined in \eqref{defTheta}.
\item The transformations $\mathscr{B}^{\pm 1},{\bf{B}}^{\pm 1}, {\beta}$ and $\widehat{\beta}$ satisfy the following estimates for all $s\in[s_{0},S]$ 
\begin{equation}\label{estimate on the first reduction operator and its inverse}
\|\mathscr{B}^{\pm 1}h\|_{s}^{q,\kappa}+\|{\bf{B}}^{\pm 1}h\|_{s}^{q,\kappa}\lesssim\|h\|_{s}^{q,\kappa}+{\varepsilon\kappa ^{-1}}\| \mathfrak{I}_{0}\|_{s+\sigma_1}^{q,\kappa}\|h\|_{s_{0}}^{q,\kappa},
\end{equation}
\begin{equation*}
\big\|\left(\mathscr{B}-\textnormal{Id}\right)h\big\|_{s}^{q,\kappa}\lesssim \varepsilon\kappa^{-1}\left(\| h\|_{s+1}^{q,\kappa}+\| \mathfrak{I}_{0}\|_{s+\sigma_1}^{q,\kappa}\| h\|_{s_{0}}^{q,\kappa}\right)
\end{equation*}
and 
\begin{equation}\label{est-beta-r}
\|\widehat{\beta}\|_{s}^{q,\kappa}\lesssim\|\beta\|_{s}^{q,\kappa}\lesssim \varepsilon\kappa ^{-1}\left(1+\| \mathfrak{I}_{0}\|_{s+\sigma_1}^{q,\kappa}\right).
\end{equation}
\item On the Cantor set
$$\mathcal{O}_{\infty,n}^{\kappa,\tau_{1}}(i_{0})=\bigcap_{(l,j)\in\mathbb{Z}^{d}\times\mathbb{Z}\backslash\{(0,0)\}\atop|l|\leqslant N_{n}}\left\lbrace\lambda \in \mathcal{O};\;\, \big|\omega\cdot l+jc(\lambda,i_0)\big|> 4\kappa^{\varrho}\tfrac{\langle j\rangle}{\langle l \rangle^{\tau_{1}}}\right\rbrace$$
we have 
$$\mathscr{B}^{-1}\big(\omega\cdot\partial_{\varphi}+\partial_{\theta}\big(V_{\varepsilon r,\alpha}\cdot\big)\big)\mathscr{B}=\omega\cdot\partial_{\varphi}+c(\lambda,i_0)\partial_{\theta}+\mathtt{E}_{n}^{0}
$$
with $\mathtt{E}_{n}^{0}$ a linear operator satisfying
$$\|\mathtt{E}_{n}^{0}h\|_{s_0}^{q,\kappa}\lesssim {\varepsilon} N_{0}^{\mu_{2}}N_{n+1}^{-\mu_{2}}\|h\|_{s_{0}+2}^{q,\kappa}.$$
The function $V_{\varepsilon r,\alpha}$ was defined in Lemma $\ref{lemma-reste}.$
\item Given two tori $i_{1}$ and $i_{2}$ both satisfying \eqref{small-C2} (replacing $\mathfrak{I}_{0}$ by $\mathfrak{I}_{1}$ or $\mathfrak{I}_{2}$), then
\begin{equation}\label{difference ci}
\|\Delta_{12}c(\cdot,i)\|^{q,\kappa}\lesssim{\varepsilon}\| \Delta_{12}i\|_{\overline{s}_{h}+4}^{q,\kappa}.
\end{equation}
In addition, assume that  $0\leqslant \mathtt{a}\leqslant \frac32(\mu_2-\overline{\mu}_2)+1-\tau_1$, then 
\begin{equation}\label{difference beta}
\|\Delta_{12}\beta\|_{\overline{s}_{h}+ \mathtt{a}}^{q,\kappa}+{\|\Delta_{12}\widehat\beta\|_{\overline{s}_{h}+ \mathtt{a}}^{q,\kappa}}\lesssim\varepsilon\kappa^{-1}\|\Delta_{12}i\|_{\overline{s}_{h}+\sigma_1+ \mathtt{a}}^{q,\kappa}.
\end{equation}
\end{enumerate}
}
\end{proposition}
Before giving the proof, some remarks are in order.
\begin{remark} 
\begin{enumerate}
\item The final Cantor set $\mathcal{O}_{\infty,n}^{\kappa,\tau_{1}}(i_{0})$ is constructed over the limit coefficient $c(\lambda,{i_0})$ but it is still truncated in the time frequency, that is  $|l|\leqslant N_n$,  leading to a residual reminder with enough decay through the parameter $\mu_2$ that can be arbitrarily chosen  by selecting the regularity index  $s_h$ large enough. This induces a suitable stability property which  is crucial during the Nash-Moser scheme achieved with  the nonlinear functional.
\item The estimate \eqref{difference beta} holds if $\mu_2$ satisfies $\mathtt{a}+\tau_1-1\leqslant \frac32(\mu_2-\overline{\mu}_2)$ for some $\mathtt{a}\geqslant 0.$ This is slightly stronger than the constraint imposed in \eqref{Conv-T2N}
\item The constant $4$ that appears in the definition of the Cantor set $\mathcal{O}_{\infty,n}^{\kappa,\tau_1}(i_{0})$ is used  to ensure the inclusion of this set in all the Cantor sets built along  the KAM procedure.
\end{enumerate}
\end{remark}
 \begin{proof}
 The proof will be done in the same spirit of \cite{BFM, FGMP19} and based on  the construction of successive  iterations of linear transformations through  quasi-periodic symplectic change of coordinates. Notice that  at each step of the scheme we should extract  from  the reminder  of size $\varepsilon$ its main diagonal part leading to  a new reminder  with size $\varepsilon^2$. This can be done through solving the homological equation which requires non-resonances conditions satisfied by excision of  the external  parameters $\lambda=(\omega,\alpha)$.  Iterating this argument allows to get the desired result with a final Cantor set constructed over all the restrictions coming from the different homological equations. To be more precise, we shall describe this procedure in the KAM step and see later how to implement it.

\smallskip

 {\bf{(i)-(ii)}}
$\blacktriangleright$ \textbf{KAM step}.
Assume that we have a transport operator taking the form, 
$$
\big(\omega\cdot\partial_{\varphi}+\partial_{\theta}\big(V+f\big)\big)h= \omega\cdot\partial_{\varphi}h+\partial_{\theta}\big((V+f) h\big)
$$
 where the parameter $\lambda$  belongs to  a subset $\mathcal{O}_{-}^{\gamma}\subset\mathcal{O}$ and 
$$V=V(\lambda)\quad\textnormal{and}\quad f=f(\lambda,\varphi,\theta)$$
with $f$ being an even  function,
\begin{equation}\label{symmetry for f}
f(\lambda,-\varphi,-\theta)=f(\lambda,\varphi,\theta).
\end{equation}
Next, we introduce a symplectic quasi-periodic change of coordinates close to the identity in the form
\begin{equation}\label{change of variables}
\begin{array}{rcl}
\mathscr{B} h(\lambda,\varphi,\theta) &\triangleq& \big(1+\partial_{\theta}g(\lambda,\varphi,\theta)\big){\bf{B}} h(\lambda,\varphi,\theta)\\
& =& \big(1+\partial_{\theta}g(\lambda,\varphi,\theta)\big) h\big(\lambda,\varphi,\theta+g(\lambda,\varphi,\theta)\big)
\end{array}
\end{equation}
with $g:\mathcal{O}\times\mathbb{T}^{d+1}\rightarrow\mathbb{R}$ being a  small  function to be adjusted later  with respect to  $f.$ Then according to \mbox{Lemma \ref{algeb1},} we may write for any  $N\geqslant 2$
\begin{equation}\label{transformation KAM step transport}
\mathscr{B}^{-1}\Big(\omega\cdot\partial_{\varphi}+\partial_{\theta}\left(V+f\right)\Big)\mathscr{B}=\omega\cdot\partial_{\varphi}+\partial_{\theta}\Big({\bf B}^{-1}\left(V+\omega\cdot\partial_{\varphi}g+V\partial_{\theta}g+\Pi_{N}f+\Pi_{N}^{\perp}f+f\partial_{\theta}g\right)\cdot\Big).
\end{equation}
The main goal  is to obtain after this transformation a new transport operator in the form
\begin{equation}\label{link V,f and V+,f+}
\mathscr{B}^{-1}\Big(\omega\cdot\partial_{\varphi}+\partial_{\theta}\left(V+f\right)\Big)\mathscr{B}=\omega\cdot\partial_{\varphi}+\partial_{\theta}\big(V_{+}+f_{+}\big)
\end{equation}
where
$$V_{+}=V_{+}(\lambda)\quad\textnormal{and}\quad f_{+}=f_{+}(\lambda,\varphi,\theta)$$
with $f_{+}$ quadratically smaller than $f.$  Coming back to \eqref{transformation KAM step transport} and  in order to get rid of  the {\it 
{linear}}  terms in $f$, we shall impose  to the perturbation $g$ the following \textit{homological equation}
\begin{equation}\label{equation satisfied by g}
\omega\cdot\partial_{\varphi}g+V\partial_{\theta}g+\Pi_{N}f=\langle f\rangle_{\varphi,\theta}
\end{equation}
where the average of $f$ is defined by 
$$\langle f\rangle_{\varphi,\theta}=\frac{1}{(2\pi)^{d+1}}\int_{\mathbb{T}^{d+1}}f(\lambda,\varphi,\theta)d\varphi d\theta.$$
To solve  the {homological equation} \eqref{equation satisfied by g}, we use Fourier decomposition in order to recover  $g$ in the form
\begin{equation}\label{definition of g}
g(\lambda,\varphi,\theta)=\sum_{(l,j)\in\mathbb{Z}^{d+1}\backslash\{0\}\atop\langle l,j\rangle\leqslant N}\frac{\ii\,f_{l,j}(\lambda)}{\omega\cdot l+jV(\lambda)}e^{\ii (l\cdot\varphi+j\theta)}.
\end{equation}
Then at this level we  should deal with the  small divisors problem which is a classical issue in KAM theory. One way to fix it is to avoid  resonances through Diophantine conditions type  by restricting the exterior parameters to  the following  set
\begin{equation}\label{definition of O+}
\mathcal{O}_{+}^{\gamma}\triangleq\bigcap_{(l,j)\in\mathbb{Z}^{d+1 }\backslash\{0\}\atop\langle l,j\rangle\leqslant N}\left\lbrace\lambda=(\omega,\alpha)\in \mathcal{O}_{-}^{\gamma};\;\, \big|\omega\cdot l+jV(\lambda)\big|>\tfrac{\kappa ^{\varrho}\langle j\rangle}{\langle l\rangle^{\tau_{1}}}\right\rbrace.
\end{equation}
With this choice  we can control the size of the denominators in \eqref{definition of g} and then expect to get suitable  estimates for $g$ with some loss of regularity uniformly in $N$. Before proceeding with  this task we need  to construct  an extension of $g$ to the whole set $\mathcal{O}$ and for the sake of simplicity it will be still denoted by $g$. This can be done by extending the Fourier coefficients of $g$ in \eqref{definition of g} using the cut-off function $\chi$ defined in \eqref{chi-def-1} in the following way
\begin{equation}\label{def extension glj}
g_{l,j}(\lambda)\triangleq \ii\frac{\chi\big((\omega\cdot l+jV(\lambda))(\kappa ^{\varrho}\langle j\rangle)^{-1}\langle l\rangle^{\tau_{1}}\big)}{\omega\cdot l+jV(\lambda)}f_{l,j}(\lambda).
\end{equation}
In what follows, we shall work with this extension which is a solution to \eqref{equation satisfied by g} when the parameters are restricted to the set $\mathcal{O}_{+}^{\gamma}.$  We then define
$$V_{+}=V+\langle f\rangle_{\varphi,\theta}\quad\textnormal{and}\quad f_{+}={\bf B}^{-1}\big(\Pi_{N}^{\perp}f+f\partial_{\theta}g\big)$$
so that in the set $\mathcal{O}_{+}^{\gamma},$ the identity \eqref{link V,f and V+,f+} holds. Notice that $V_{+}$ and $f_{+}$ are well-defined in the whole set of parameters $\mathcal{O}$ and the function $g$ is smooth since it is generated by a finite number of frequencies. From the assumption \eqref{symmetry for f} we get that $g$ is odd and therefore
\begin{equation}\label{symmetry for g}
g\in\bigcap_{s\geqslant 0}W^{q,\infty}_{\kappa}(\mathcal{O},H_{\textnormal{odd}}^{s}).
\end{equation}
Next, we intend to estimate the Fourier coefficients $g_{l,j}$ defined in \eqref{def extension glj} which can be  written in the form
\begin{align}\label{g-lj-7}
g_{l,j}(\lambda)&=\ii a_{l,j}\,{\chi_1}\big(a_{l,j}A_{l,j}(\lambda)\big)\,f_{l,j}(\lambda),\quad A_{l,j}(\lambda)=\omega\cdot l+jV(\lambda),\\
\nonumber\quad a_{l,j}&=(\kappa^{\varrho}\langle j\rangle)^{-1}\langle l\rangle^{\tau_{1}}, \quad \chi_1(x)=\tfrac{\chi(x)}{x}\cdot
\end{align}
Since $\chi_1$ is $C^\infty$ with bounded derivatives and $\chi_1(0)=0$, then we may apply  Lemma \ref{Compos-lemm-VM} giving
\begin{align*}
\forall\,|\gamma|\leqslant q,\quad \|g_{l,j}\|_{W^{|\gamma|,\infty}(\mathcal{O})}\lesssim a_{l,j}^2\|A_{l,j}\|_{W^{|\gamma|,\infty}(\mathcal{O})}\Big(1+a_{l,j}^{|\gamma|-1}\|A_{l,j}\|_{L^{\infty}(\mathcal{O})}^{|\gamma|-1}\Big).
\end{align*}
It is straightforward that
\begin{align*}
\forall(l,j)\in\mathbb{Z}^{d+1},\,\forall\gamma\in\mathbb{N}^{d+1},\quad|\gamma|\leqslant q,\quad \sup_{\lambda\in\mathcal{O}}\left|\partial_{\lambda}^{\gamma}A_{l,j}(\lambda)\right|&\lesssim\langle l,j\rangle\left(1+\sup_{\lambda\in\mathcal{O}}\left|\partial_{\lambda}^{\gamma}V(\lambda)\right|\right)\\
&\lesssim\kappa^{-|\gamma|}\langle l,j\rangle\left(1+\|V\|^{q,\kappa}\right).
\end{align*}
Hence by assuming 
\begin{equation}\label{boundedness assumption on V}
\|V\|^{q,\kappa}\leqslant C
\end{equation}
we  obtain
$$
\forall(l,j)\in\mathbb{Z}^{d+1},\quad \forall\gamma\in\mathbb{N}^{d+1},\quad|\gamma|\leqslant q\quad \sup_{\lambda\in\mathcal{O}}\left|\partial_{\lambda}^{\gamma}A_{l,j}(\lambda)\right|\lesssim\kappa^{-|\gamma|}\langle l,j\rangle.
$$
Therefore we find that
\begin{align*}
 \forall |\gamma| \leqslant q,\quad \|g_{l,j}\|_{W^{|\gamma|,\infty}(\mathcal{O})}\lesssim a_{l,j}^2\kappa^{-|\gamma|}\langle l,j\rangle\Big(1+a_{l,j}^{|\gamma|-1}\langle l,j\rangle^{|\gamma|-1}\Big).
\end{align*}
Since $0\leqslant a_{l,j}\leqslant \kappa^{-\varrho}\langle l\rangle^{\tau_{1}}$ then the foregoing estimate gives
\begin{align}\label{Dent-oist}
\forall |\gamma| \leqslant q,\quad |\partial_{\lambda}^\gamma g_{l,j}|\lesssim \kappa^{-\varrho(|\gamma|+1)-|\gamma|}\langle l,j\rangle^{\tau_1(1+|\gamma|)+|\gamma|}.
\end{align}
By choosing  $\varrho$ such that
\begin{equation}\label{first choice of upsilon}
\varrho\leqslant\frac{1}{q+1}
\end{equation}
and using  Leibniz rule or the law products, we infer from the Definition \ref{Def-WS}
\begin{equation}\label{control of g by f}
\|g\|_{s}^{q,\kappa}\lesssim\kappa ^{-1}\|\Pi_{N}f\|_{s+\tau_{1} q+\tau_{1}}^{q,\kappa}.
\end{equation}
Assume 
\begin{equation}\label{smallness assumption f}
\kappa^{-1}N^{\tau_{1}q+\tau_{1} {+1}}\|f\|_{s_{0}}^{q,\kappa}\leqslant{\varepsilon}_0,
\end{equation}
then  combined with  \eqref{control of g by f} and Lemma \ref{orthog-Lem1}, we get
$$\|g\|_{s_{0}+{1}}^{q,\kappa}\leqslant C\kappa^{-1}N^{\tau_{1}q+\tau_{1}+{1}}\|f\|_{s_{0}}^{q,\kappa}\leqslant C{\varepsilon}_0.$$
Hence, taking ${\varepsilon}_0$ small enough we may guarantee  the smallness condition in Lemma \ref{Compos1-lemm} and get  that the linear operator $\mathscr{B}$ is  invertible. We now set 
$$u=\Pi_{N}^{\perp}f+f\partial_{\theta}g.$$
By Lemma \ref{Law-prodX1} and \eqref{control of g by f}, we deduce  for all $s\in[s_{0},S]$ 
$$\begin{array}{rcl}
\|u\|_{s}^{q,\kappa} 
& \leqslant & \|\Pi_{N}^{\perp}f\|_{s}^{q,\kappa}+C\Big(\|f\|_{s_0}^{q,\kappa }\|\partial_{\theta}g\|_{s}^{q,\kappa}+\|f\|_{s}^{q,\kappa}\|\partial_{\theta}g\|_{s_0}^{q,\kappa }\Big)\\
& \leqslant & \|\Pi_{N}^{\perp}f\|_{s}^{q,\kappa}+C\kappa ^{-1}N^{\tau_{1}q+\tau_{1}+1}\|f\|_{s_0}^{q,\kappa }\|f\|_{s}^{q,\kappa}.
\end{array}$$
Applying   Lemma \ref{Compos1-lemm}, Lemma \ref{Law-prodX1} and \eqref{smallness assumption f}, we obtain 
$$\begin{array}{rcl}
\|f_{+}\|_{s}^{q,\kappa} & = & \|{\bf B}^{-1}(u)\|_{s}^{q,\kappa}\\
& \leqslant & \|u\|_{s}^{q,\kappa}+C\Big(\| u\|_{s}^{q,\kappa}\|\widehat{g}\|_{s_0}^{q,\kappa }+\|\widehat{g}\|_{s}^{q,\kappa}\|u\|_{s_0}^{q,\kappa }\Big)\\
& \leqslant & \|u\|_{s}^{q,\kappa}+C\Big(\|u\|_{s}^{q,\kappa}\|g\|_{s_0}^{q,\kappa }+\|g\|_{s}^{q,\kappa}\|u\|_{s_0}^{q,\kappa }\Big)\\
& \leqslant & \|\Pi_{N}^{\perp}f\|_{s}^{q,\kappa}+C\kappa ^{-1}N^{\tau_{1} q+\tau_{1}+1}\|f\|_{s_0}^{q,\kappa }\|f\|_{s}^{q,\kappa}.
\end{array}$$
Therefore, Lemma \ref{orthog-Lem1} yields for $ s_{0} \leqslant s\leqslant \overline{s}\leqslant S$
\begin{equation}\label{estimate KAM step transport}
\|f_{+}\|_{s}^{q,\kappa}\leqslant N^{s-\overline{s}}\|f\|_{q,\overline{s}}^{\gamma ,\mathcal{O}}+C\kappa ^{-1}N^{\tau_{1} q+\tau_{1}+1}\|f\|_{s_0}^{q,\kappa }\|f\|_{s}^{q,\kappa}.
\end{equation}

\smallskip

$\blacktriangleright$ \textbf{KAM scheme}.  Assume that we have constructed $V_{m}$ and $f_{m}$ for $m\geqslant 0$,  being well-defined in the whole set of parameters $\mathcal{O}$ and verifying the assumptions \eqref{boundedness assumption on V} and \eqref{smallness assumption f}, we want to construct $V_{m+1}$ and $f_{m+1}$ still satisfying  \eqref{boundedness assumption on V} and \eqref{smallness assumption f}. For this purpose, we shall apply the KAM step with $(V,f,V_{+},f_{+},N)$ replaced by $(V_{m},f_{m},V_{m+1},f_{m+1},N_{m}).$ To be more precise, we will prove by induction the existence of a sequence $\{V_m,f_m\}_{m\in\mathbb{N}}$ such that 
\begin{equation}\label{hypothesis of induction deltam}
\delta_{m}(s_{l})\leqslant\delta_{0}(s_{h})N_{0}^{\mu_{2}}N_{m}^{-\mu_{2}},\quad\quad\delta_{m}(s_{h})\leqslant \left(2-\frac{1}{m+1}\right)\delta_{0}(s_{h})
\end{equation}
and
\begin{equation}\label{assumptions KAM iterations}
\|V_{m}\|^{q,\kappa}\leqslant C,\quad\quad N_{m}^{\tau_{1}q+\tau_{1}+1}\delta_{m}(s_{0})\leqslant{\varepsilon}_0,
\end{equation}
supplemented with the symmetry condition
\begin{equation}\label{symmetry for f_{m}}
f_{m}(\lambda,-\varphi,-\theta)=f_{m}(\lambda,\varphi,\theta),
\end{equation}
where
$$\delta_{m}(s)\triangleq \kappa ^{-1}\| f_{m}\|_{s}^{q,\kappa}.$$
We remind that the parameters $s_l$ and  $s_h$ in \eqref{hypothesis of induction deltam} are  defined in \eqref{Conv-T2N}.

\smallskip

\ding{70} \textit{Initialization.} Let us check the properties \eqref{hypothesis of induction deltam} and  \eqref{assumptions KAM iterations} with $m=0$. In this case we start with the transport operator
$$
\omega\cdot\partial_{\varphi}+\partial_{\theta}\big(V_{\varepsilon r,\alpha}\cdot)
$$
and use the following  decomposition based on Lemma \ref{lemma-reste}
\begin{align}\label{V=V0+f0}
\nonumber V_{\varepsilon r,\alpha}&=V_{0,\alpha}+(V_{\varepsilon r,\alpha}-V_{0,\alpha})\\
&\triangleq V_0(\alpha)+f_0.
\end{align}
Remind that $V_{r,\alpha}$ and $V_{0,\alpha}$ are defined in \eqref{V-eq} and \eqref{defTheta}.
Applying  Lemma  \ref{lemma-reste}-$($i$)$ and Proposition \ref{lemma-GS0} yield 
\begin{align}\label{estimate delta0 and I0}
\delta_{0}(s)=&\kappa^{-1}\|V_{\varepsilon r}-V_{0}\|_{s}^{q,\kappa}\nonumber\\
&\lesssim\varepsilon\kappa^{-1}\|r\|_{s+4}^{q,\kappa}\nonumber\\
&\quad\lesssim\varepsilon\kappa^{-1}\left(1+\|\mathfrak{I}_{0}\|_{s+{4}}^{q,\kappa}\right).
\end{align}
Then, using  the smallness condition \eqref{small-C2} we obtain
\begin{equation}\label{initial smallness condition in sh norm}
N_{0}^{\mu_{2}}\delta_{0}(s_{h})\leqslant C{\varepsilon}_0.
\end{equation}
Moreover, by \eqref{V-eq} and the symmetry of $r$ we get
\begin{equation}\label{symmetry for f0}
f_{0}(\lambda,-\varphi,-\theta)=f_{0}(\lambda,\varphi,\theta)
\end{equation}
We consider $\mathcal{O}_{0}^{\gamma}=\mathcal{O}$ and $N_{0}\geqslant 2$  and we shall check that the assumptions \eqref{boundedness assumption on V} and \eqref{smallness assumption f} are satisfied with $V_{0}$ and $f_{0}.$ First recall that $V_{0}$ is defined by \eqref{defTheta}. From the $C^\infty$-regularity of Gamma function, we easily obtain
\begin{equation}\label{boundedness of V0}
\|V_{0}\|^{q,\kappa}\leqslant C.
\end{equation}
Then the assumption \eqref{boundedness assumption on V} is satisfied for $V=V_{0}.$ Now,  applying  \eqref{initial smallness condition in sh norm} and using the assumption on $\mu_2$ we infer
\begin{align*}
\kappa^{-1}N_{0}^{\tau_{1}q+\tau_{1}+1}\|f_{0}\|_{s_{0}}^{q,\kappa}=& N_{0}^{\tau_{1}q+\tau_{1}+1}\delta_{0}(s_{0})\\
&\leqslant N_{0}^{\tau_{1}q+\tau_{1}+1-\mu_{2}}N_{0}^{\mu_{2}}\delta_{0}({s_{h}})\\
&\quad \leqslant C{\varepsilon}_0N_{0}^{-1}.
\end{align*}
Taking  $N_{0}$ large enough such that 
\begin{equation}\label{Est-N0}
CN_{0}^{-1}\leqslant 1
\end{equation}
we find
$$
\kappa^{-1}N_{0}^{\tau_{1}q+\tau_{1}+1}\|f_{0}\|_{s_{0}}^{q,\kappa}\leqslant{\varepsilon}_0.
$$
Thus, the condition \eqref{smallness assumption f} is satisfied for $f=f_{0}$ which achieves the initialization step. 

\smallskip

\ding{70} {\textit{Induction.}}  Assume  that we have constructed $V_m$ and $f_m$ with the  properties \eqref{hypothesis of induction deltam}, \eqref{assumptions KAM iterations} and \eqref{symmetry for f_{m}}, and let us construct  $V_{m+1}$ and $f_{m+1}$ and check the validity of these constraints at this order. Following the KAM step, we may consider a symplectic quasi-periodic change of variables $\mathscr{B}_{m}$ in the form 
\begin{align*}
\mathscr{B}_{m} h(\lambda,\varphi,\theta)&\triangleq\big(1+\partial_{\theta}g_{m}(\lambda,\varphi,\theta)\big){\bf B}_{m} h(\lambda,\varphi,\theta)\\
&=\big(1+\partial_{\theta}g_{m}(\lambda,\varphi,\theta)\big) h\big(\lambda,\varphi,\theta+g_{m}(\lambda,\varphi,\theta)\big)
\end{align*}
with 
\begin{align}\label{gm-Form5}
g_{m}(\lambda,\varphi,\theta)\triangleq\sum_{(l,j)\in\mathbb{Z}^{d+1}\backslash\{0\}\atop\langle l,j\rangle\leqslant N_{m}}i\frac{\chi\big((\omega\cdot l+jV_{m}(\lambda))(\kappa^{\varrho}\langle j\rangle)^{-1}\langle l\rangle^{\tau_{1}}\big)}{\omega\cdot l+jV_{m}(\lambda)}(f_{m})_{l,j}(\lambda)\,e^{\ii (l\cdot\varphi+j\theta)} 
\end{align}
where $\chi$ is the cut-off function defined in \eqref{chi-def-1} and $N_{m}$ is defined in \eqref{definition of Nm}. As it was explained during the KAM step, $g_{m}$ is well-defined on the whole set of parameters $\mathcal{O}$ and when it is  restricted to the  Cantor set 
\begin{equation}\label{definition mathcal Om+1}
\mathcal{O}_{m+1}^{\kappa}\triangleq\bigcap_{\underset{\langle l,j\rangle\leqslant N_{m}}{(l,j)\in\mathbb{Z}^{d+1}\backslash\{0\}}}\left\lbrace\lambda\in \mathcal{O}_{m}^{\kappa };\;\, \big|\omega\cdot l+jV_{m}(\lambda)\big|>\tfrac{\kappa^{\varrho}\langle j\rangle}{\langle l\rangle^{\tau_{1}}}\right\rbrace
\end{equation}
it solves  the following \textit{homological equation}
$$\omega\cdot\partial_{\varphi}g_{m}+V_{m}\partial_{\theta}g_{m}+\Pi_{N_{m}}f_{m}=\langle f_{m}\rangle_{\varphi,\theta}.
$$
Therefore  in the Cantor set $\mathcal{O}_{m+1}^{\gamma },$ we find the following reduction
\begin{align}\label{Conjug-ope-tr}
\mathscr{B}_{m}^{-1}\Big(\omega\cdot\partial_{\varphi}+\partial_{\theta}(V_{m}+f_{m})\Big)\mathscr{B}_{m}=\omega\cdot\partial_{\varphi}+\partial_{\theta}(V_{m+1}+f_{m+1})
\end{align}
where $V_{m+1}$ and $f_{m+1}$ are defined by
\begin{equation}\label{definition Vm+1 and fm+1}
\left\lbrace\begin{array}{l}
V_{m+1}=V_{m}+\langle f_{m}\rangle_{\varphi,\theta}\\
f_{m+1}={\bf B}_{m}^{-1}\Big(\Pi_{N_{m}}^{\perp}f_{m}+f_{m}\partial_{\theta}g_{m}\Big).
\end{array}\right.
\end{equation}
According to \eqref{symmetry for f_{m}}, the function  $g_m$ is odd and consequently we deduce  that $f_{m+1}$ is even which concludes the symmetry persistence during the scheme.
On the other hand, one gets similarly to \eqref{symmetry for g} that
\begin{equation}\label{symmetry for g_{m}}
g_{m}\in\bigcap_{s\geqslant 0}W^{q,\infty}_{\kappa}(\mathcal{O},H_{\textnormal{odd}}^{s})
\end{equation}
We set 
\begin{align}\label{Beta-iter-V9}
\overline{\mathscr{B}}_{-1}=\textnormal{Id}\quad\hbox{ and for m}\in\N,\quad \overline{\mathscr{B}}_{m}=\mathscr{B}_{0}\circ\mathscr{B}_{1}\circ...\circ\mathscr{B}_{m}.
\end{align} Then we may check that
\begin{align*}
\overline{\mathscr{B}}_{m}h(\lambda,\varphi,\theta)&\triangleq\big(1+\partial_{\theta}\beta_{m}(\lambda,\varphi,\theta)\big) h\big(\lambda,\varphi,\theta+\beta_{m}(\lambda,\varphi,\theta)\big)
\end{align*}
where $(\beta_{m})_{m\in\mathbb{N}}$ is defined by $\beta_{-1}=g_{-1}=0$ and  
\begin{equation}\label{definition betam}
\beta_{0}=g_{0}\quad\mbox{ and }\quad\beta_{m}(\lambda,\varphi,\theta)=\beta_{m-1}(\lambda,\varphi,\theta)+g_m\big(\lambda,\varphi,\theta+\beta_{m-1}(\lambda,\varphi,\theta)\big).
\end{equation}
By a trivial induction using \eqref{symmetry for g_{m}} and Lemma \ref{Compos-lemm} we have 
\begin{equation}\label{symmetry for h_{m}}
\beta_{m}\in\bigcap_{s\geqslant 0}W^{q,\infty}_{\kappa}(\mathcal{O},H_{\textnormal{odd}}^{s}).
\end{equation}
 From  Sobolev embeddings, \eqref{definition Vm+1 and fm+1} and the  induction hypothesis \eqref{hypothesis of induction deltam}, we have
\begin{align}\label{Cauchy Vm}
\| V_{m}-V_{m-1}\|^{q,\kappa}=&\|\langle f_{m-1}\rangle_{\varphi,\theta}\|^{q,\kappa}\nonumber\\
&\leqslant\| f_{m-1}\|_{s_0}^{q,\kappa }\nonumber=\gamma\delta_{m-1}(s_{0})\nonumber\\
&\quad \leqslant\kappa\delta_{0}(s_{h})N_{0}^{\mu_{2}}N_{m-1}^{-\mu_{2}}.
\end{align}
This implies in view of the  triangle inequality  combined with \eqref{initial smallness condition in sh norm} and with $\varepsilon_0$ small enough  
\begin{align*}
\| V_{m}\|^{q,\kappa}  \leqslant & \| V_{m-1}\|^{q,\kappa}+\kappa\delta_{0}(s_{h})N_{0}^{\mu_{2}}N_{m-1}^{-\mu_{2}}\\
& \leqslant  \displaystyle\| V_{0}\|^{q,\kappa}+\kappa\delta_{0}(s_{h})N_{0}^{\mu_{2}}\sum_{k=0}^{m-1}N_{k}^{-\mu_{2}}\\
& \quad \leqslant  \displaystyle\| V_{0}\|^{q,\kappa}+\sum_{k=0}^{+\infty}N_{k}^{-\mu_{2}}.
\end{align*}
On the other hand, by  \eqref{hypothesis of induction deltam}, \eqref{initial smallness condition in sh norm} and since  $\mu_{2}\geqslant \tau_{1}q+\tau_{1}+2$ in view of \eqref{Conv-T2}, we deduce by the assumption \eqref{Est-N0}
\begin{align}\label{useful estimate in the induction}
\delta_{m}(s_{0})N_{m}^{\tau_{1} q+\tau_{1}+1}\leqslant&\delta_{0}(s_{h})N_{0}^{\mu_{2}}N_{m}^{\tau_{1} q+\tau_{1}+1-\mu_{2}}\nonumber\\
&\leqslant C\varepsilon_0 N_{0}^{-1}\nonumber\\
&\quad \leqslant{\varepsilon_0}.
\end{align}
Therefore we obtain in view of \eqref{boundedness of V0} and \eqref{useful estimate in the induction}
\begin{equation*}
\sup_{m\in\mathbb{N}}\|V_{m}\|^{q,\kappa}\leqslant C\quad\textnormal{and}\quad \delta_{m}(s_{0})N_{m}^{\tau_{1}q+\tau_{1}+1}\leqslant{\varepsilon}_0.
\end{equation*}
At this stage, we may  apply the KAM step and the estimate \eqref{estimate KAM step transport} writes in our case
\begin{equation}\label{recurrence estimate deltam}
\delta_{m+1}(s)\leqslant N_{m}^{s-\overline{s}}\delta_{m}(\overline{s})+CN_{m}^{\tau_{1}q+\tau_{1}+1}\delta_{m}(s)\delta_{m}(s_{0}).
\end{equation}
Applying \eqref{recurrence estimate deltam} with $s=s_{l}$ and $\overline{s}=s_{h}$,  we get 
$$\delta_{m+1}(s_{l})\leqslant N_{m}^{s_{l}-s_{h}}\delta_{m}(s_{h})+CN_{m}^{\tau_{1} q+\tau_{1}+1}\delta_{m}(s_{l})\delta_{m}(s_{0}).$$
Putting together \eqref{hypothesis of induction deltam} and the fact that $s_{l}\geqslant s_{0}$ combined with Sobolev embeddings yields
\begin{align*}
\delta_{m+1}(s_{l})&\leqslant N_{m}^{s_{l}-s_{h}}\delta_{m}(s_{h})+CN_{m}^{\tau_{1} q+\tau_{1}+1}(\delta_{m}(s_{l}))^{2}\\
&\leqslant\left(2-\frac{1}{m+1}\right)N_{m}^{s_{l}-s_{h}}\delta_{0}(s_{h})+CN_{0}^{2\mu_{2}}N_{m}^{\tau_{1}q+\tau_{1}+1-2\mu_{2}}(\delta_{0}(s_{h}))^{2}\\
&\leqslant 2N_{m}^{s_{l}-s_{h}}\delta_{0}(s_{h})+CN_{0}^{2\mu_{2}}N_{m}^{\tau_{1}q+\tau_{1}+1-2\mu_{2}}(\delta_{0}(s_{h}))^{2}.
\end{align*}
If we select our parameters $s_{l},$ $s_{h}$ and $\mu_{2}$ such that
\begin{equation}\label{conv-t1}
N_{m}^{s_{l}-s_{h}}\leqslant\frac{1}{4}N_{0}^{\mu_{2}}N_{m+1}^{-\mu_{2}}\quad\textnormal{and}\quad CN_{0}^{2\mu_{2}}N_{m}^{\tau_{1}q+\tau_{1}+1-2\mu_{2}}\delta_{0}(s_{h})\leqslant\frac{1}{2}N_{0}^{\mu_{2}}N_{m+1}^{-\mu_{2}}
\end{equation}
then 
$$\delta_{m+1}(s_{l})\leqslant \delta_{0}(s_{h})N_{0}^{\mu_{2}}N_{m+1}^{-\mu_{2}}.$$
Notice that by our choice of parameters \eqref{Conv-T2}, the condition \eqref{conv-t1} is fulfilled provided that
$$4N_{0}^{-\mu_{2}}\leqslant 1\quad\textnormal{and}\quad 2C\delta_{0}(s_{h})\leqslant N_{0}^{-\mu_{2}}.$$
The first condition is automatically satisfied by taking $N_{0}$ sufficiently large and the second condition holds  by \eqref{initial smallness condition in sh norm} provided that $\overline{\varepsilon}$ is  small enough. This proves the first statement of the induction in \eqref{hypothesis of induction deltam} and we now turn to the proof of the second statement. Making use of  \eqref{recurrence estimate deltam} with $s=\overline{s}$  and using the  induction assumption \eqref{hypothesis of induction deltam}, we find
\begin{align*}
\delta_{m+1}(s_{h})&\leqslant\delta_{m}(s_{h})\left(1+CN_{m}^{\tau_{1} q+\tau_{1}+1}\delta_{m}(s_{0})\right)\\
&\leqslant\left(2-\frac{1}{m+1}\right)\delta_{0}(s_{h})\big(1+CN_{0}^{\mu_{2}}N_{m}^{\tau_{1} q+\tau_{1}+1-\mu_{2}}\delta_{0}(s_{h})\big).
\end{align*}
Therefore we deduce from \eqref{Conv-T2N} and \eqref{small-C2}
\begin{equation}\label{conv-t2}
\left(2-\frac{1}{m+1}\right)\left(1+CN_{0}^{\mu_{2}}N_{m}^{\tau_{1} q+\tau_{1}+1-\mu_{2}}\delta_{0}(s_{h})\right)\leqslant 2-\frac{1}{m+2},
\end{equation}
then we find
$$
\delta_{m+1}(s_{h})\leqslant\left(2-\frac{1}{m+2}\right)\delta_{0}(s_{h}).
$$
This estimate  achieves the induction argument of \eqref{hypothesis of induction deltam}. Now, observe that \eqref{conv-t2} is equivalent to
$$\left(2-\frac{1}{m+1}\right)CN_{0}^{\mu_{2}}N_{m}^{\tau_{1} q+\tau_{1}+1-\mu_{2}}\delta_{0}(s_{h})\leqslant\frac{1}{(m+1)(m+2)}.$$
By virtue of \eqref{Conv-T2} we have  $\mu_{2}\geqslant \tau_{1}q+\tau_{1}+2$ and thus  the preceding condition holds true if
\begin{equation}\label{conv-t3}
CN_{0}^{\mu_{2}}N_{m}^{-1}\delta_{0}(s_{h})\leqslant\frac{1}{(m+1)(m+2)}\cdot
\end{equation}
Recall that  $N_{0}\geqslant 2,$ and then in view of \eqref{definition of Nm} we may find a constant $c_{0}>0$ small enough such that
$$\forall m\in\mathbb{N},\quad c_{0}N_{m}^{-1}\leqslant\frac{1}{(m+1)(m+2)}\cdot
$$
Therefore, \eqref{conv-t3} is satisfied provided that
\begin{equation}\label{conv-t4}
CN_{0}^{\mu_{2}}\delta_{0}(s_{h})\leqslant c_{0}.
\end{equation}
By taking ${\varepsilon}$ small enough, we can ensure from \eqref{initial smallness condition in sh norm} that
\begin{align*}
CN_{0}^{\mu_{2}}\delta_{0}(s_{h})&\leqslant C{\varepsilon}_0\\
&\leqslant c_{0}.
\end{align*}
Thus, the condition \eqref{conv-t4} is satisfied. This completes the proof of \eqref{hypothesis of induction deltam}.

\smallskip

\ding{70} {\textit{Regularity persistence.}} Applying \eqref{recurrence estimate deltam} with $\overline{s}=s\in [s_{0},S]$, \eqref{hypothesis of induction deltam} and \eqref{Conv-T2}, we obtain
\begin{align*}
\delta_{m+1}(s) & \leqslant\delta_{m}(s)\left(1+CN_{m}^{\tau_{1} q+\tau_{1}+1}\delta_{m}(s_{0})\right)\\
\\
& \leqslant\delta_{m}(s)\left(1+C\delta_{0}(s_{h})N_{0}^{\mu_{2}}N_{m}^{\tau_{1} q+\tau_{1}+1-\mu_{2}}\right)\\
&\leqslant\delta_{m}(s)\left(1+CN_{m}^{-1}\right).
\end{align*}
Combining  this estimate with \eqref{estimate delta0 and I0}, we infer by a trivial induction
\begin{align}\label{uniform estimate of deltams}
\delta_{m}(s)\leqslant&\delta_{0}(s)\prod_{k=0}^{+\infty}\left(1+CN_{k}^{-1}\right)\nonumber\\
&\leqslant C\delta_{0}(s)\\
&\quad \leqslant C\varepsilon\kappa^{-1}\left(1+\|\mathfrak{I}_{0}\|_{s+4}^{q,\kappa}\right).\nonumber
\end{align}
Putting together  \eqref{control of g by f}, the interpolation inequality from  Lemma \ref{interpolation-In} and \eqref{hypothesis of induction deltam} leads to 
\begin{align*}
\|g_{m}\|_{s}^{q,\kappa}&\leqslant C\delta_{m}(s+\tau_{1}q+\tau_{1})\\
&\leqslant C\left(\delta_{m}(s_{0})\right)^{\overline{\theta}(s)}\left(\delta_{m}(s+\tau_{1}q+\tau_{1}+1)\right)^{1-\overline{\theta}(s)}\\
&\leqslant C\delta_{0}^{\overline{\theta}(s)}(s_{h})N_{0}^{\overline{\theta}(s)\mu_{2}}N_{m}^{-\overline{\theta}(s)\mu_{2}}\delta_{m}^{1-\overline{\theta}(s)}(s+\tau_{1}q+\tau_{1}+1)
\end{align*}
with $\overline{\theta}(s)\triangleq\frac{1}{s+\tau_{1}q+\tau_{1}+1-s_{0}}.$
Using \eqref{uniform estimate of deltams},  \eqref{small-C2} and \eqref{estimate delta0 and I0}
\begin{align}\label{bound gm}
\|g_{m}\|_{s}^{q,\kappa}&\leqslant C\varepsilon\kappa^{-1}\left(1+\|\mathfrak{I}_{0}\|_{s_{h}+4}^{q,\kappa}\right)\left(1+\|\mathfrak{I}_{0}\|_{s+\tau_{1}q+\tau_{1}+5}^{q,\kappa}\right))N_{0}^{\overline{\theta}(s)\mu_{2}}N_{m}^{-\overline{\theta}(s)\mu_{2}}\nonumber\\
&\leqslant C\varepsilon\kappa^{-1}\left(1+\|\mathfrak{I}_{0}\|_{s+\tau_{1}q+\tau_{1}+5}^{q,\kappa}\right)N_{0}^{\overline{\theta}(s)\mu_{2}}N_{m}^{-\overline{\theta}(s)\mu_{2}}.
\end{align}
By \eqref{definition betam} and Lemma \ref{Compos1-lemm}, we deduce for all $s\in[s_{0},S]$
\begin{equation}\label{link betam and betam-1}
\|\beta_{m}\|_{s}^{q,\kappa}\leqslant\|\beta_{m-1}\|_{s}^{q,\kappa}\left(1+C\| g_{m}\|_{s_0}^{q,\kappa }\right)+C\left(1+\| \beta_{m-1}\|_{s_{0}}^{q,\kappa}\right)\| g_{m}\|_{s}^{q,\kappa}.
\end{equation}
Applying this estimate with $s=s_{0}$ together  with Sobolev embeddings,
$$\|\beta_{m}\|_{s_0}^{q,\kappa }\leqslant\|\beta_{m-1}\|_{s_0}^{q,\kappa }\left(1+C\| g_{m}\|_{s_0}^{q,\kappa }\right)+C\| g_{m}\|_{s_0}^{q,\kappa }.$$
In order to obtain a good estimate for  $\beta_{m}$, we shall use the following result which is quite easy to prove by induction : Given three positive sequences $(a_n)_{n\in\mathbb{N}},(b_n)_{n\in\mathbb{N}}$ and $(c_n)_{n\in\mathbb{N}}$ such that
$$
\forall\, n\in\mathbb{N},\quad a_{n+1}\leqslant b_n a_n+c_n.
$$
Then
\begin{align}\label{Ind-res}
	\nonumber \forall\, n\geqslant 2,\quad a_{n}&\leqslant a_0\prod_{i=0}^{n-1}b_i+\sum_{k=0}^{n-2}c_k\prod_{i=k+1}^{n-1}b_i+c_{n-1}\\
	&\leqslant \Big(a_0+\sum_{k=0}^{n-1} c_k\Big)\prod_{i=0}^{n-1}b_i.
\end{align}
In particular if $\displaystyle \prod_{n=0}^{\infty}b_n$ and $\displaystyle \sum_{n=0}^{\infty} c_n$ converge then
\begin{align}\label{Ind-res1}
	\sup_{n\in\mathbb{N}}a_{n}&\leqslant \Big(a_0+\sum_{n=0}^{\infty} c_n\Big)\prod_{n=0}^{\infty}b_i.
\end{align}
Since \eqref{Conv-T2N} ensures that $s_{h}\geqslant s_{0}+\tau_{1}q+\tau_{1}+1$ and $\mu_{2}\geqslant \tau_{1}q+\tau_{1}+1$, then by \eqref{bound gm} and \eqref{small-C2} we get
\begin{align*}
\|g_{m}\|_{s_{0}}^{q,\kappa}&\leqslant C\varepsilon\kappa^{-1}N_{0}^{\mu_{2}}\left(1+\|\mathfrak{I}_{0}\|_{s_{0}+\tau_{1}q+\tau_{1}+5}^{q,\kappa}\right)N_{m}^{-\overline{\theta}(s_{0})\mu_{2}}\\
&\leqslant C{\varepsilon}_0N_{m}^{-1}.
\end{align*}
By \eqref{definition of Nm} and  $N_{0}\geqslant2$ we infer  $\displaystyle\sum_{m=0}^{\infty} N_{m}^{-1}<\infty$. Therefore if ${\varepsilon_0}$ is small enough in such a way  $C{\varepsilon_0}\leqslant1$, and  applying \eqref{Ind-res1} together with the fact that $\beta_{0}=g_{0}$ we find
\begin{align}\label{unif betam s0}
\sup_{m\in\mathbb{N}}\|\beta_{m}\|_{s_{0}}^{q,\kappa}&\leqslant\left(\|\beta_{0}\|_{s_{0}}^{q,\kappa}+C\sum_{k=0}^{+\infty}\|g_{k}\|_{s_{0}}^{q,\kappa}\right)\prod_{k=0}^{\infty}\left(1+C\|g_{k}\|_{s_{0}}^{q,\kappa}\right)\nonumber\\
&\leqslant \left(1+C\sum_{k=0}^{\infty}N_{k}^{-1}\right)\prod_{k=0}^{\infty}\left(1+N_{k}^{-1}\right)\leqslant C.
\end{align}
Plugging this estimate into   \eqref{link betam and betam-1} gives for all $s\in[s_{0},S]$
$$\|\beta_{m}\|_{s}^{q,\kappa}\leqslant\|\beta_{m-1}\|_{s}^{q,\kappa}\left(1+C\|g_{m}\|_{s_{0}}^{q,\kappa}\right)+C\|g_{m}\|_{s}^{q,\kappa}.$$
As above, applying \eqref{Ind-res1} and \eqref{bound gm}, we find 
\begin{align*}
\sup_{m\in\mathbb{N}}\|\beta_{m}\|_{s}^{q,\kappa}&\leqslant\left(\|\beta_{0}\|_{s}^{q,\kappa}+C\sum_{k=0}^{\infty}\|g_{k}\|_{s}^{q,\kappa}\right)\prod_{k=0}^{\infty}\left(1+C\|g_{k}\|_{s_{0}}^{q,\kappa}\right)\\
&\leqslant C\varepsilon\kappa^{-1}\left(1+\|\mathfrak{I}_{0}\|_{s+\tau_{1}q+\tau_{1}+5}^{q,\kappa}\right)\left(1+N_{0}^{\overline{\theta}(s)\mu_{2}}\sum_{k=0}^{\infty}N_{k}^{-\overline{\theta}(s)\mu_{2}}\right).
\end{align*}
Using Lemma \ref{lemma sum Nn} implies the existence of a  constant $C>0$ such that  
$$\forall s\in[s_{0},S],\quad N_{0}^{\overline{\theta}(s)\mu_{2}}\sum_{k=0}^{+\infty}N_{k}^{-\overline{\theta}(s)\mu_{2}}\leqslant C$$
and this allows to get
\begin{equation}\label{uniform estimate betam}
\forall s\in[s_{0},S],\quad\sup_{m\in\mathbb{N}}\| \beta_{m}\|_{s}^{q,\kappa}\leqslant C\varepsilon\kappa^{-1}\left(1+\|\mathfrak{I}_{0}\|_{s+\tau_{1}q+\tau_{1}+5}^{q,\kappa}\right).
\end{equation}
By Lemma \ref{Compos1-lemm}, \eqref{unif betam s0}, \eqref{control of g by f}, \eqref{hypothesis of induction deltam} and in view of  $s_{l}\geqslant s_{0}+\tau_{1}q+\tau_{1}+2$, we obtain,
\begin{align}\label{decay telescopic beta}
\|\beta_{m}-\beta_{m-1}\|_{s_{0}+2}^{q,\kappa}\leqslant &C\|g_{m}\|_{s_{0}+2}^{q,\kappa}\left(1+\|\beta_{m-1}\|_{s_{0}+2}^{q,\kappa}\right)\nonumber\\
&\leqslant C\delta_{m}(s_{l})\nonumber\\
&\quad\leqslant CN_{0}^{\mu_{2}}N_{m}^{-\mu_{2}}\delta_{0}(s_{h}).
\end{align}
Hence we get the convergence of the series
$$\sum_{m=0}^{\infty}\|\beta_{m}-\beta_{m-1}\|_{s_{0}+2}^{q,\kappa}<\infty.$$
Consequently, the sequence $(\beta_{m})_{m\in\mathbb{N}}$ converges strongly towards some function $\beta\in W^{q,\infty}_{\kappa}(\mathcal{O},H^{s_{0}+2}).$  From the uniform boundedness \eqref{uniform estimate betam} we deduce that $\beta\in W^{q,\infty}_{\kappa}(\mathcal{O},H^{s})$ and 
\begin{align}\label{estimate beta}
\forall s\in[s_{0},S],\quad\|\beta\|_{s}^{q,\kappa}&\leqslant\liminf_{m\rightarrow\infty}\|\beta_{m}\|_{s}^{q,\kappa}\nonumber\\
&\lesssim\varepsilon\kappa^{-1}\left(1+\|\mathfrak{I}_{0}\|_{s+\tau_{1}q+\tau_{1}+5}^{q,\kappa}\right).
\end{align}
Define now the quasi-periodic symplectic change of variables $\mathscr{B}$ associated with  $\beta$  by
\begin{align*}
\mathscr{B} h(\lambda,\omega,\varphi,\theta)&=\left(1+\partial_{\theta}\beta(\lambda,\omega,\varphi,\theta)\right){\bf B}h(\lambda,\omega,\varphi,\theta)\\
&=\left(1+\partial_{\theta}\beta(\lambda,\omega,\varphi,\theta)\right)h(\lambda,\omega,\varphi,\theta+\beta(\lambda,\omega,\varphi,\theta)).
\end{align*}
Notice that by \eqref{estimate beta}, hyptothesis \eqref{small-C2} and  \eqref{Conv-T2N}, we infer 
\begin{equation}\label{beta s0 bound}
\|\beta\|_{s_{0}}^{q,\kappa}\lesssim\varepsilon\kappa^{-1}\left(1+\|\mathfrak{I}_{0}\|_{s_{0}+\tau_{1}q+\tau_{1}+5}^{q,\kappa}\right)\leqslant C{\varepsilon_0}.
\end{equation}
Hence choosing ${\varepsilon_0}$ small enough, we get  by Lemma \ref{Compos1-lemm} that $\mathscr{B}$ is an invertible operator,  and combined with   \eqref{estimate beta} we obtain
\begin{equation}\label{est Bpm1}
\|\mathscr{B}^{\pm 1} h\|_{s}^{q,\kappa}\lesssim\| h\|_{s}^{q,\kappa}+\varepsilon\kappa^{-1}\|\mathfrak{I}_{0}\|_{s+\tau_{1}q+\tau_{1}+6}^{q,\kappa}\| h\|_{s_{0}}^{q,\kappa}.
\end{equation}
Notice that the symmetry of $\beta$ follows from the symmetry of the approximation $(\beta_m)$ described in   \eqref{symmetry for h_{m}}. In addition, one easily gets from \eqref{decay telescopic beta}
\begin{align}\label{difference beta and betam}
\|\beta-\beta_{m}\|_{s_{0}+2}^{q,\kappa}&\leqslant\sum_{k=m}^{\infty}\|\beta_{k+1}-\beta_{k}\|_{s_{0}+2}^{q,\kappa}\nonumber\\
&\lesssim\kappa\delta_{0}(s_{h})N_{0}^{\mu_{2}}\sum_{k=m+1}^{\infty}N_{k}^{-\mu_{2}}.
\end{align}
Applying  Lemma \ref{lemma sum Nn} yields
\begin{equation}\label{big O series Nm}
\sum_{k=m}^{+\infty}N_{k}^{-\mu_{2}}\underset{m\rightarrow\infty}{=}O\left(N_{m}^{-\mu_{2}}\right).
\end{equation}
Putting together \eqref{big O series Nm} and \eqref{difference beta and betam}, we obtain
\begin{equation}\label{rate on convergence betam to beta}
\|\beta-\beta_{m}\|_{s_{0}+2}^{q,\kappa}\lesssim\kappa\delta_{0}(s_{h})N_{0}^{\mu_{2}}N_{m+1}^{-\mu_{2}}.
\end{equation}

\smallskip

\noindent $\blacktriangleright$ \textbf{KAM conclusion}. 
 According to \eqref{Cauchy Vm} and \eqref{big O series Nm} one may write,
\begin{align*}\sum_{m=0}^{\infty}\| V_{m+1}-V_{m}\|^{q,\kappa}&\leqslant\kappa\delta_{0}(s_{h})N_{0}^{\mu_{2}}\sum_{m=0}^{\infty}N_{m}^{-\mu_{2}}\\
&\lesssim \kappa \delta_{0}(s_{h}).
\end{align*}
Thus we get   that the sequence $(V_{m})_{m\in\mathbb{N}}$ converges towards an element $\lambda\mapsto c(\lambda,i_0)\in W_{\kappa}^{q,\infty }(\mathcal{O},\RR).$ We denote $c_{i_{0}}$ its limit and one gets
\begin{align}\label{struc-ci}
c_{i_{0}}=V_0+\sum_{m\in\N} (V_{m+1}-V_m).
\end{align}
 Therefore  by applying  \eqref{estimate delta0 and I0} we find
\begin{align*}
\| c_{i_{0}}-V_{0}\|^{q,\kappa}&\leqslant\sum_{m=0}^{\infty}\| V_{m+1}-V_{m}\|^{q,\kappa}\\
&\lesssim\varepsilon\left(1+\|\mathfrak{I}_{0}\|_{s_{h}+4}^{q,\kappa}\right){\lesssim{\varepsilon}}_0.
\end{align*}
Next, we consider the truncated Cantor set
$$\mathcal{O}_{\infty,n}^{\kappa,\tau_{1}}(i_{0})=\bigcap_{(l,j)\in\mathbb{Z}^{d+1}\backslash\{0\}\atop|l|\leqslant N_{n}}\left\lbrace\lambda=(\omega,\alpha)\in\mathcal{O};\;\, \big|\omega\cdot l+jc(\lambda,{i_{0}})\big|>\tfrac{4\kappa^{\varrho}\langle j\rangle}{\langle l\rangle^{\tau_{1}}}\right\rbrace
$$
and we intend to  prove that this  set satisfies the inclusion
\begin{align}\label{kloip67}\mathcal{O}_{\infty,n}^{\kappa,\tau_{1}}(i_{0})\subset\bigcap_{m=0}^{n+1}\mathcal{O}_{m}^{\kappa}=\mathcal{O}_{n+1}^{\kappa}
\end{align}
where the intermediate Cantor sets are defined in \eqref{definition mathcal Om+1}. For this goal, we proceed by induction in $m$. First, notice that by construction we have $\mathcal{O}_{\infty,n}^{\kappa,\tau_{1}}(i_{0})\subset\mathcal{O}\triangleq\mathcal{O}_{0}^{\kappa}.$ Second,  assume that $\mathcal{O}_{\infty,n}^{\kappa,\tau_{1}}(i_{0})\subset\mathcal{O}_{m}^{\kappa}$ for some $m\leqslant n$ and let us check that
\begin{equation}\label{inclusion Oinftyn in Om+1}
\mathcal{O}_{\infty,n}^{\kappa,\tau_{1}}(i_{0})\subset\mathcal{O}_{m+1}^{\kappa}.
\end{equation}
Remark that by \eqref{Cauchy Vm} and \eqref{big O series Nm}
\begin{align}\label{estimate Vm-ci0}
|V_{m}(\lambda)-c(\lambda,i_0)|\leqslant&\| V_{m}-c(\cdot,i_0)\|^{q,\kappa}\nonumber\\
&\leqslant\sum_{l=m}^{\infty}\| V_{l+1}-V_{l}\|^{q,\kappa}\nonumber\\
&\quad \lesssim\kappa\delta_{0}(s_{h})N_{0}^{\mu_{2}}N_{m}^{-\mu_{2}}.
\end{align}
Let  $\lambda\in\mathcal{O}_{\infty,n}^{\kappa,\tau_{1}}(i_{0})$ and  $(l,j)\in\mathbb{Z}^{d+1}\backslash\{0\}$ with $0\leqslant|l|\leqslant N_{m}.$ Then $|l|\leqslant N_{n}$ and by the triangle inequality  combined with  \eqref{estimate Vm-ci0},  \eqref{initial smallness condition in sh norm} and $\varrho,\kappa\in(0,1)$  we infer
\begin{align*}
|\omega\cdot l+jV_{m}(\lambda)|  \geqslant& |\omega\cdot l+jc(\lambda,i_0)|-|j||V_{m}(\lambda)-c(\lambda,i_0)|\\
&\geqslant\,\displaystyle\frac{4\kappa^{\varrho}\langle j\rangle}{\langle l\rangle^{\tau_{1}}}-C\langle j\rangle\kappa\delta_{0}(s_{h})N_{0}^{\mu_{2}}N_{m}^{-\mu_{2}}\\
&\quad\geqslant\displaystyle\frac{4\kappa^{\varrho}\langle j\rangle}{\langle l\rangle^{\tau_{1}}}-C\langle j\rangle\kappa^{\varrho}{\varepsilon_0}\langle l\rangle^{-\mu_{2}}
\end{align*}
Then for $C{\varepsilon_0}\leqslant1$ and since $\mu_{2}\geqslant\tau_{1}$ (in view of \eqref{Conv-T2}), we get
$$|\omega\cdot l+jV_{m}(\lambda)|> \displaystyle\frac{\kappa^{\varrho}\langle j\rangle}{\langle l\rangle^{\tau_{1}}},
$$
which shows that $\lambda\in\mathcal{O}_{m+1}^{\gamma}$ and therefore the inclusion \eqref{inclusion Oinftyn in Om+1} is satisfied.

\smallskip

{\bf{(iii)}} We shall start with the following decomposition
\begin{align*}
\mathscr{B}^{-1}\Big(\omega\cdot\partial_{\varphi}+\partial_{\theta}\left(V_{0}+f_{0}\right)\Big)\mathscr{B}=&\mathscr{B}_{n}^{-1}\Big(\omega\cdot\partial_{\varphi}+\partial_{\theta}\left(V_{0}+f_{0}\right)\Big)\mathscr{B}_{n}\\
&+\left(\mathscr{B}^{-1}-\mathscr{B}_{n}^{-1}\right)\Big(\omega\cdot\partial_{\varphi}+\partial_{\theta}\left(V_{0}+f_{0}\right)\Big)\mathscr{B}\\
&\quad+{\mathscr{B}_n^{-1}}\Big(\omega\cdot\partial_{\varphi}+\partial_{\theta}\left(V_{0}+f_{0}\right)\Big)\left(\mathscr{B}-\mathscr{B}_{n}\right).
\end{align*}
By \eqref{kloip67}, \eqref{Conjug-ope-tr} and \eqref{Beta-iter-V9}we have that  on the Cantor set $\mathcal{O}_{\infty,n}^{\kappa,\tau_{1}}(i_{0})$ 
$$\mathscr{B}_{n}^{-1}\Big(\omega\cdot\partial_{\varphi}+\partial_{\theta}\left(V_{0}+f_{0}\right)\Big)\mathscr{B}_{n}=\omega\cdot\partial_{\varphi}+\partial_{\theta}\left(V_{n+1}+f_{n+1}\right).$$
Hence, in the Cantor set $\mathcal{O}_{\infty,n}^{\kappa,\tau_{1}}(i_{0}),$ the following splitting holds 
$$\mathscr{B}^{-1}\Big(\omega\cdot\partial_{\varphi}+\partial_{\theta}\left(V_{0}+f_{0}\right)\Big)\mathscr{B}=\omega\cdot\partial_{\varphi}+c(\lambda,i_0)\partial_{\theta}+\mathtt{E}_{n}^{0}(\lambda,i_{0})$$
with
\begin{align*}
\mathtt{E}_{n}^{0}(\lambda,i_{0})\triangleq&\left(V_{n+1}-c(\lambda,i_0)\right)\partial_{\theta}+\partial_{\theta}\left(f_{n+1}\cdot\right)+\left(\mathscr{B}^{-1}-\mathscr{B}_{n}^{-1}\right)\Big(\omega\cdot\partial_{\varphi}+\partial_{\theta}\left(V_{0}+f_{0}\right)\Big)\mathscr{B}\\
&+\mathscr{B}_n^{-1}\Big(\omega\cdot\partial_{\varphi}+\partial_{\theta}\left(V_{0}+f_{0}\right)\Big)\left(\mathscr{B}-\mathscr{B}_{n}\right)\\
&\quad\triangleq \mathtt{E}_{n,1}^{0}(\lambda,i_{0})+\mathtt{E}_{n,2}^{0}(\lambda,i_{0})+\mathtt{E}_{n,3}^{0}(\lambda,i_{0})+\mathtt{E}_{n,4}^{0}(\lambda,i_{0}).
\end{align*}
Applying the law products in Lemma \ref{Law-prodX1} together with  \eqref{estimate Vm-ci0}, we find 
\begin{align}\label{En11 s}
\|\mathtt{E}_{n,1}^{0}h\|_{s_0}^{q,\kappa}&\lesssim\|V_{n+1}-c(\cdot,i_{0})\|^{q,\kappa}\|h\|_{s_0+1}^{q,\kappa}\nonumber\\
&\lesssim\kappa\delta_{0}(s_{h})N_{0}^{\mu_{2}}N_{n+1}^{-\mu_{2}}\|h\|_{s_0+1}^{q,\kappa}.
\end{align}
Concerning the second term we use  Lemma \ref{Law-prodX1} together with \eqref{hypothesis of induction deltam}   and \eqref{uniform estimate of deltams}
\begin{align}\label{En12 s}
\|\mathtt{E}_{n,2}^{0}h\|_{s_0}^{q,\kappa}=&\|\partial_{\theta}\left(f_{n+1} h\right)\|_{s_0}^{q,\kappa}\nonumber\\
&\lesssim\kappa\delta_{n+1}(s_{l})\|h\|_{s_0+1}^{q,\kappa}\nonumber\\
&\lesssim\kappa\delta_{0}(s_{h})N_{0}^{\mu_{2}}N_{n+1}^{-\mu_{2}}\|h\|_{s_0+1}^{q,\kappa}.
\end{align}
We now move  to the estimate of $\mathtt{E}_{n,3}^{0}$.  Using  Lemma \ref{Law-prodX1}, we deduce that
\begin{align*}
\|\omega\cdot\partial_{\varphi}h+\partial_{\theta}\left(V_{\varepsilon r}h\right)\|_{s_0}^{q,\kappa}&\leqslant\|\omega\cdot\partial_{\varphi}h\|_{s}^{q,\kappa}+\|V_{\varepsilon r}h\|_{s_0+1}^{q,\kappa}\\
&\lesssim \| h\|_{s_0+1}^{q,\kappa}\left(1+\|V_{\varepsilon r}\|_{s_{0}+1}^{q,\kappa}\right).
\end{align*}
Putting together  \eqref{V=V0+f0}, \eqref{boundedness of V0} and \eqref{estimate delta0 and I0}, we find
\begin{align*}
\|V_{\varepsilon r}\|_{s_0+1}^{q,\kappa}&\leqslant\|V_{0}\|^{q,\kappa}+\|f_{0}\|_{s_0+1}^{q,\kappa}\\
&\leqslant C
\end{align*}
Hence by\eqref{small-C2}, we obtain
\begin{equation}\label{estimate initial transport operator}
\|\omega\cdot\partial_{\varphi}h+\partial_{\theta}\left(V_{\varepsilon r}h\right)\|_{s_0}^{q,\kappa}\lesssim\|h\|_{s_0+1}^{q,\kappa}.
\end{equation}
Then, by virtue of Taylor Formula, 
\begin{align*}
(\mathscr{B}-\mathscr{B}_{n})h(\theta)&=(1+\partial_{\theta}\beta(\theta))\left[h(\theta+\beta(\theta))-h(\theta+\beta_{n}(\theta))\right]+\partial_{\theta}\left(\beta(\theta)-\beta_{n}(\theta)\right)h(\theta+\beta_{n}(\theta))\\
&\triangleq (1+\partial_\theta\beta(\theta))(\beta(\theta)-\beta_{n}(\theta))\mathtt{I}_{n}h(\theta)+\partial_{\theta}(\beta(\theta)-\beta_{n}(\theta)){\bf B}_{n}h(\theta),
\end{align*}
with
$$\mathtt{I}_{n}h(\theta)\triangleq\int_{0}^{1}\big(\partial_{\theta}h\big)\big(\theta+\beta_{n}(\theta)+t(\beta(\theta)-\beta_{n}(\theta))\big)dt.$$
Combining Lemma \ref{Compos1-lemm} and Lemma \ref{Law-prodX1} yields 
\begin{align*}
\|\partial_{\theta}\left(\beta-\beta_{n}\right){\bf B}_nh\|_{s_0+1}^{q,\kappa}\lesssim&\|\beta-\beta_{n}\|_{s_0+2}^{q,\kappa}\|h\|_{s_{0}+1}^{q,\kappa}\left(1+\|\beta_{n}\|_{s_{0}+1}^{q,\kappa}\right)\\
&\lesssim \kappa\delta_{0}(s_{h})N_{0}^{\mu_{2}}N_{n+1}^{-\mu_{2}}\|h\|_{s_{0}}^{q,\kappa}
\end{align*}
On the other hand, performing similar arguments as before we find 
\begin{align*}
\Big\|(1+\partial_{\theta}\beta)\left(\beta-\beta_{n}\right)\mathtt{I}_{n}h\Big\|_{s_{0}+1}^{q,\kappa}&\lesssim\left(1+\|\beta\|_{s_{0}+2}^{q,\kappa}\right)\|\beta-\beta_{n}\|_{s_{0}+1}^{q,\kappa}\|h\|_{s_{0}+2}^{q,\kappa}\left(1+\|\beta_{n}\|_{s_{0}+1}^{q,\kappa}+\|\beta-\beta_{n}\|_{s_{0}+1}^{q,\kappa}\right)\\
&\lesssim\kappa\delta_{0}(s_{h})N_{0}^{\mu_{2}}N_{n+1}^{-\mu_{2}}\|h\|_{s_{0}+2}^{q,\kappa}.
\end{align*}

Putting together the preceding estimates, it follows that
\begin{align*}
\|(\mathscr{B}-\mathscr{B}_{n})h\|_{s_{0}+1}^{q,\kappa}&\lesssim\kappa\delta_{0}(s_{h})N_{0}^{\mu_{2}}N_{n+1}^{-\mu_{2}}\|h\|_{s_{0}+2}^{q,\kappa}.
\end{align*}
Thus we find by collecting   \eqref{estimate initial transport operator},  \eqref{est Bpm1}, \eqref{estimate delta0 and I0} and \eqref{small-C2}
\begin{align}\label{En14 s0}
\nonumber\|\mathtt{E}_{n,4}^{0}h\|_{s_{0}}^{q,\kappa}&\lesssim \kappa\delta_{0}(s_{h})N_{0}^{\mu_{2}}N_{n+1}^{-\mu_{2}}\|h\|_{s_{0}+2}^{q,\kappa}\\
&\lesssim \varepsilon N_{0}^{\mu_{2}}N_{n+1}^{-\mu_{2}}\|h\|_{s_{0}+2}^{q,\kappa}.
\end{align}
In a similar way, using in particular \eqref{mathscrB1} we find
\begin{align}\label{En13 s0}
\nonumber\|\mathtt{E}_{n,3}^{0}h\|_{s_{0}}^{q,\kappa}&\lesssim \kappa\delta_{0}(s_{h})N_{0}^{\mu_{2}}N_{n+1}^{-\mu_{2}}\|h\|_{s_{0}+2}^{q,\kappa}\\
&\lesssim \varepsilon N_{0}^{\mu_{2}}N_{n+1}^{-\mu_{2}}\|h\|_{s_{0}+2}^{q,\kappa}.
\end{align}
Gathering \eqref{En11 s}, \eqref{En13 s0}, \eqref{En14 s0}, yields
$$\|\mathtt{E}_{n}^{0}h\|_{s_{0}}^{q,\kappa}\lesssim\varepsilon N_{0}^{\mu_{2}}N_{n+1}^{-\mu_{2}}\|h\|_{s_{0}+2}^{q,\kappa}.$$

\smallskip

{\bf{(iv)}} We shall start with the estimate of $\Delta_{12}\beta$. For this aim we notice that, since $\beta_{-1}=0$, 
\begin{equation}\label{series Delta12beta}
\Delta_{12}\beta=\sum_{m=0}^{\infty}\Delta_{12}(\beta_{m}-\beta_{m-1}).
\end{equation}
Hence by the  triangule inequality we get 
\begin{equation}\label{estimate series Delta12 beta}
\|\Delta_{12}\beta\|_{\overline{s}_{h}+\mathtt{a}}^{q,\kappa}\leqslant\sum_{m=0}^{\infty}\|\Delta_{12}(\beta_{m}-\beta_{m-1})\|_{\overline{s}_{h}+\mathtt{a}}^{q,\kappa}.
\end{equation}
By Taylor formula and \eqref{definition betam}, we find by removing the  dependance with respect to $\lambda$ and $\varphi$ 
\begin{align*}
\Delta_{12}\beta_{m}(\theta)&=\Delta_{12}\beta_{m-1}(\theta)+{\bf B}_{m-1}^{[1]}(\Delta_{12}g_{m})(\theta)\\
&+\Delta_{12}\beta_{m-1}(\theta)\int_{0}^{1}(\partial_{\theta}g_{m}^{[2]})\big(\theta+\beta_{m-1}^{[2]}(\theta)+t\Delta_{12}\beta_{m-1}(\theta)\big)dt,
\end{align*}
which implies that
\begin{align*}
\Delta_{12}(\beta_{m}-\beta_{m-1})(\theta)&={\bf B}_{m-1}^{[1]}(\Delta_{12}g_{m})(\theta)\\
&+\Delta_{12}\beta_{m-1}(\theta)\int_{0}^{1}(\partial_{\theta}g_{m}^{[2]})\big(\theta+\beta_{m-1}^{[2]}(\theta)+t\Delta_{12}\beta_{m-1}(\theta)\big)dt.
\end{align*}
Thus, by law product in Lemma \ref{Law-prodX1}, Lemma \ref{Compos1-lemm} and Sobolev embeddings
\begin{align*}
\|\Delta_{12}(\beta_{m}-\beta_{m-1})\|_{\overline{s}_{h}+\mathtt{a}}^{q,\kappa}&\leqslant\|\Delta_{12}g_{m}\|_{\overline{s}_{h}+\mathtt{p}}^{q,\kappa}\left(1+C\|\beta_{m-1}^{[1]}\|_{s_{0}}^{q,\kappa}\right)+\|\Delta_{12}g_{m}\|_{s_{0}}^{q,\kappa}\|\beta_{m-1}^{[1]}\|_{\overline{s}_{h}+\mathtt{a}}^{q,\kappa}\\
&+C\|\Delta_{12}\beta_{m-1}\|_{s_{0}}^{q,\kappa}\|g_{m}^{[2]}\|_{\overline{s}_{h}+\mathtt{a}+1}^{q,\kappa}\left(1+\|\beta_{m-1}^{[2]}\|_{s_{0}}^{q,\kappa}+\|\Delta_{12}\beta_{m-1}\|_{s_{0}}^{q,\kappa}\right)\\
&+C\|\Delta_{12}\beta_{m-1}\|_{s_{0}}^{q,\kappa}\|g_{m}^{[2]}\|_{s_{0}+1}^{q,\kappa}\left(\|\beta_{m-1}^{[2]}\|_{\overline{s}_{h}+\mathtt{a}}^{q,\kappa}+\|\Delta_{12}\beta_{m-1}\|_{\overline{s}_{h}+\mathtt{a}}^{q,\kappa}\right)\\
&+C\|\Delta_{12}\beta_{m-1}\|_{\overline{s}_{h}+\mathtt{a}}^{q,\kappa}\|g_{m}^{[2]}\|_{s_{0}+1}^{q,\kappa}\left(1+\|\beta_{m-1}^{[2]}\|_{s_{0}}^{q,\kappa}+\|\Delta_{12}\beta_{m-1}\|_{s_{0}}^{q,\kappa}\right)
\end{align*}
and for all $s\in[s_{0},\overline{s}_{h}+\mathtt{a}]$
\begin{align*}
\|\Delta_{12}\beta_{m}\|_{s}^{q,\kappa}&\leqslant\|\Delta_{12}g_{m}\|_{s}^{q,\kappa}\left(1+C\|\beta_{m-1}^{[1]}\|_{s_{0}}^{q,\kappa}\right)+\|\Delta_{12}g_{m}\|_{s_{0}}^{q,\kappa}\|\beta_{m-1}^{[1]}\|_{s}^{q,\kappa}\\
&+C\|\Delta_{12}\beta_{m-1}\|_{s_{0}}^{q,\kappa}\|g_{m}^{[2]}\|_{s+1}^{q,\kappa}\left(1+\|\beta_{m-1}^{[2]}\|_{s_{0}}^{q,\kappa}+\|\Delta_{12}\beta_{m-1}\|_{s_{0}}^{q,\kappa}\right)\\
&+C\|\Delta_{12}\beta_{m-1}\|_{s_{0}}^{q,\kappa}\|g_{m}^{[2]}\|_{s_{0}+1}^{q,\kappa}\left(\|\beta_{m-1}^{[2]}\|_{s}^{q,\kappa}+\|\Delta_{12}\beta_{m-1}\|_{s}^{q,\kappa}\right)\\
&+\|\Delta_{12}\beta_{m-1}\|_{s}^{q,\kappa}\left(1+C\|g_{m}^{[2]}\|_{s_{0}+1}^{q,\kappa}\left(1+\|\beta_{m-1}^{[2]}\|_{s_{0}}^{q,\kappa}+\|\Delta_{12}\beta_{m-1}\|_{s_{0}}^{q,\kappa}\right)\right).
\end{align*}
Here we adopt the generic notation  $g^{[k]}$ to denote the value of $g$ at $r_k$ (or equivalently to $i_k$).
By \eqref{bound gm}, \eqref{small-C2} and \eqref{Conv-T2N}, we have since $\overline{s}_{h}+\mathtt{a}+\tau_{1}q+\tau_{1}+6\leqslant s_h+\sigma_1$
\begin{align}\label{maj gmrk}
\sup_{m\in\mathbb{N}}\max_{k\in\{1,2\}}\|g_{m}^{[k]}\|_{\overline{s}_{h}+\mathtt{a}+1}^{q,\kappa}&\leqslant C \varepsilon\kappa^{-1}\left(1+\max_{k\in\{1,2\}}\|\mathfrak{I}_{0}^{[k]}\|_{\overline{s}_{h}+\mathtt{a}+\tau_{1}q+\tau_{1}+6}^{q,\kappa}\right)\nonumber\\
& \leqslant C,
\end{align}
Actually,  we get a more precise estimate
\begin{align}\label{deacr gm}
\max_{k\in\{1,2\}}\|g_{m}^{[k]}\|_{q,\overline{s}_{h}+\mathtt{a}+1}&\leqslant C\varepsilon\kappa^{-1}\left(1+\max_{k\in\{1,2\}}\|\mathfrak{I}_{0}^{[k]}\|_{\overline{s}_{h}+\mathtt{a}+\tau_{1}q+\tau_{1}+6}^{q,\kappa}\right)N_{0}^{\overline{\theta}\left(\overline{s}_{h}+\mathtt{a}+1\right)\mu_{2}}N_{m}^{-\overline{\theta}\left(\overline{s}_{h}+\mathtt{a}+1\right)\mu_{2}}\nonumber\\
&\leqslant C\varepsilon\kappa^{-1}N_{0}^{\overline{\theta}\left(\overline{s}_{h}+\mathtt{a}+1\right)\mu_{2}}N_{m}^{-\overline{\theta}\left(\overline{s}_{h}+\mathtt{a}+1\right)\mu_{2}}.
\end{align}
The estimate  \eqref{uniform estimate betam} and \eqref{small-C2} allow to get
\begin{align}\label{maj betamrk}
\sup_{m\in\mathbb{N}}\max_{k\in\{1,2\}}\|\beta_{m}^{[k]}\|_{q,\overline{s}_{h}+\mathtt{a}}&\leqslant C.
\end{align}
Putting together  \eqref{maj gmrk}, \eqref{maj betamrk}  with  the previous two estimates imply
\begin{align}\label{itmdt dlt12 diff bt}
\|\Delta_{12}(\beta_{m}-\beta_{m-1})\|_{\overline{s}_{h}+\mathtt{a}}^{q,\kappa}&\leqslant C\left(\|\Delta_{12}g_{m}\|_{\overline{s}_{h}+\mathtt{p}}^{q,\kappa}+\|\Delta_{12}\beta_{m-1}\|_{\overline{s}_{h}+\mathtt{p}}^{q,\kappa}\|g_{m}^{[2]}\|_{\overline{s}_{h}+\mathtt{a}+1}^{q,\kappa}\right),
\end{align}
\begin{align}\label{dlt12 btm s0}
\|\Delta_{12}\beta_{m}\|_{s_{0}}^{q,\kappa}&\leqslant C\|\Delta_{12}g_{m}\|_{s_{0}}^{q,\kappa}+\|\Delta_{12}\beta_{m-1}\|_{s_{0}}^{q,\kappa}\left(1+C\|g_{m}^{[2]}\|_{s_{0}+1}^{q,\kappa}\right)
\end{align}
and
\begin{align}\label{dlt12 btm shb+1}
\|\Delta_{12}\beta_{m}\|_{\overline{s}_{h}+\mathtt{a}}^{q,\kappa}&\leqslant C\left(\|\Delta_{12}g_{m}\|_{\overline{s}_{h}+\mathtt{a}}^{q,\kappa}+\|\Delta_{12}\beta_{m-1}\|_{s_{0}}^{q,\kappa}\|g_{m}^{[2]}\|_{q,\underline{s_{h}}+\mathtt{p}+1}\right)\nonumber\\
&+\|\Delta_{12}\beta_{m-1}\|_{\overline{s}_{h}+\mathtt{a}}^{q,\kappa}\left(1+C\|g_{m}^{[2]}\|_{s_{0}+1}^{q,\kappa}\right).
\end{align}
 Applying \eqref{Ind-res1} with \eqref{dlt12 btm s0} combined with  the fact that $\beta_{0}=g_{0}$, we conclude that
$$\sup_{m\in\mathbb{N}}\|\Delta_{12}\beta_{m}\|_{s_{0}}^{q,\kappa}\leqslant\left(\|\Delta_{12}g_{0}\|_{s_{0}}^{q,\kappa}+C\sum_{m=0}^{\infty}\|\Delta_{12}g_{m}\|_{s_{0}}^{q,\kappa}\right)\prod_{m=0}^{\infty}\left(1+\|g_{m}^{[2]}\|_{s_{0}+1}^{q,\kappa}\right).$$
Therefore we obtain from  \eqref{deacr gm}, 
$$\sup_{m\in\mathbb{N}}\|\Delta_{12}\beta_{m}\|_{s_{0}}^{q,\kappa}\leqslant C\sum_{m=0}^{\infty}\|\Delta_{12}g_{m}\|_{s_{0}}^{q,\kappa}.$$
In the same way, from \eqref{dlt12 btm shb+1}, using \eqref{Ind-res1}, \eqref{deacr gm} and the previous estimate, we conclude that
$$\sup_{m\in\mathbb{N}}\|\Delta_{12}\beta_{m}\|_{\overline{s}_{h}+\mathtt{a}}^{q,\kappa}\leqslant C\sum_{m=0}^{\infty}\|\Delta_{12}g_{m}\|_{\overline{s}_{h}+\mathtt{a}}^{q,\kappa}.$$ 
Gathering the previous bounds with \eqref{itmdt dlt12 diff bt} and \eqref{deacr gm}, we finaly obtain
\begin{equation}\label{telscopic diff betam}
\|\Delta_{12}(\beta_{m}-\beta_{m-1})\|_{\overline{s}_{h}+\mathtt{a}}^{q,\kappa}\lesssim\|\Delta_{12}g_{m}\|_{\overline{s}_{h}+\mathtt{a}}^{q,\kappa}+{\varepsilon \kappa^{-1} N_0^{\overline{\theta}(\overline{s}_{h}+\mathtt{a}+1)\mu_{2}}}N_{m}^{-\overline{\theta}(\overline{s}_{h}+\mathtt{a}+1)\mu_{2}}\sum_{k=0}^{+\infty}\|\Delta_{12}g_{k}\|_{\overline{s}_{h}+\mathtt{a}}^{q,\kappa}.
\end{equation}
The next goal is   to estimate $\Delta_{12}g_{m}.$ First observe that from \eqref{gm-Form5} and  \eqref{g-lj-7} we can write
\begin{align*}
g_{m}(\lambda,\varphi,\theta)=\ii\sum_{(l,j)\in\mathbb{Z}^{d+1}\backslash\{0\}\atop\langle l,j\rangle\leqslant N_{m}}a_{l,j}\,\chi_1\big(a_{l,j}A_{l,j}(\lambda)\big)(f_{m})_{l,j}(\lambda)\,e^{\ii (l\cdot\varphi+j\theta)}.
\end{align*}
Then we find
\begin{align*}
\Delta_{12}g_{m}=&\ii\sum_{(l,j)\in\mathbb{Z}^{d+1}\backslash\{0\}\atop\langle l,j\rangle\leqslant N_{m}}a_{l,j}\,\chi_1\big(a_{l.j}A^{[1]}_{l,j}\big)(\Delta_{12}f_{m})_{l,j}\,{\bf{e}}_{l,j}\\
&+\ii\sum_{(l,j)\in\mathbb{Z}^{d+1}\backslash\{0\}\atop\langle l,j\rangle\leqslant N_{m}}a_{l,j}\,\Delta_{12}\chi_1\big(a_{l,j}A^{[1]}_{l,j}\big)(f_{m}^{[2]})_{l,j}\,{\bf{e}}_{l,j}\\
&\quad\triangleq {\bf{I}}_1+{\bf{I}}_2,
\end{align*}
where we recall that  the notation $A^{[k]}_{l,j}$ stands for the value of $A_{l,j}$  evaluated at  $r_k,$ (or the torus $i_k$) and the same thing applies  for $f_m^{[k]}$.  The estimate of the first term ${\bf{I}}_1$ is quite similar to \eqref{control of g by f} and one gets
$$
\|{\bf{I}}_1\|_{s}^{q,\kappa}\lesssim\kappa ^{-1}\|\Pi_{N_m}\Delta_{12}f_m\|_{s+\tau_{1} q+\tau_{1}}^{q,\kappa}.
$$
As to the second term ${\bf{I}}_2$ we write in view of Taylor formula
\begin{align}\label{II22}
{\bf{I}}_2=&\ii\sum_{(l,j)\in\mathbb{Z}^{d+1}\backslash\{0\}\atop\langle l,j\rangle\leqslant N_{m}}a_{l,j}^2(\Delta_{12}A_{l,j})\,\int_0^1\chi_1^\prime\left(a_{l,j}\big[(1-\tau)A^{[1]}_{l,j}+\tau A^{[2]}_{l,j}\big]\right)d\tau\,(f_{m}^{[2]})_{l,j}\,{\bf{e}}_{l,j}.
\end{align}
Then by  the same analysis developed to get \eqref{control of g by f}, using in particular,
$$
\|\Delta_{12}A_{l,j}\|^{q,\kappa}\lesssim \langle l,j\rangle \|\Delta_{12}V_m\|^{q,\kappa}
$$
and
$$
\|a_{l,j}\Delta_{12}A_{l,j}\|^{q,\kappa}\lesssim\kappa^{-\varrho} \langle l,j\rangle^{1+\tau_1} \|\Delta_{12}V_m\|^{q,\kappa}
$$
we find 
\begin{align*}
\|{\bf{I}}_2\|_{s}^{q,\kappa}&\lesssim {\kappa ^{-q(1+\varrho)}}\|\Delta_{12}V_m\|^{q,\kappa}\|\Pi_{N_m}f_m^{[2]}\|_{s+\tau_{1} q+2\tau_{1}+1}^{q,\kappa}\\
&\lesssim\kappa ^{-1}\|\Delta_{12}V_m\|^{q,\kappa}\|\Pi_{N_m}f_m^{[2]}\|_{s+\tau_{1} q+2\tau_{1}+1}^{q,\kappa}.
\end{align*}
It follows that  for any $s\geqslant s_{0}$
\begin{align}\label{estimate delta12 gm}
\|\Delta_{12}g_{m}\|_{s}^{q,\kappa}&\lesssim\kappa^{-1}\|\Pi_{N_{m}}\Delta_{12}f_{m}\|_{s+\tau_{1}q+\tau_{1}}^{q,\kappa}+\kappa^{-1}\|\Delta_{12}V_{m}\|^{q,\kappa}\|\Pi_{N_{m}}f_{m}^{[2]}\|_{s+\tau_{1}q+2\tau_{1}+1}^{q,\kappa}\\
&\lesssim\kappa^{-1}N_{m}^{\tau_{1}q+\tau_{1}}\|\Delta_{12}f_{m}\|_{s}^{q,\kappa}+\kappa^{-1}\|\Delta_{12}V_{m}\|^{q,\kappa}\|\Pi_{N_m}f_{m}^{[2]}\|_{s+\tau_{1}q+2\tau_{1}+1}^{q,\kappa}\nonumber.
\end{align}
Now, we need to estimate $\Delta_{12}f_{m}.$ Then by setting 
\begin{align}\label{um-est}
u_{m}\triangleq \Pi_{N_{m}}^{\perp}f_{m}+f_{m}\partial_{\theta}g_{m},
\end{align}
and making appeal to \eqref{definition Vm+1 and fm+1} we obtain 
$$\Delta_{12}f_{m+1}=({\bf B}_{m}^{-1})^{[1]}\Delta_{12}u_{m}+\left(\Delta_{12}{\bf B}_{m}^{-1}\right)u_{m}^{[2]}$$
with
$$\Delta_{12}u_{m}=\Pi_{N_{m}}^{\perp}\Delta_{12}f_{m}+\Delta_{12}f_{m}\partial_{\theta}g_{m}^{[1]}+f_{m}^{[2]}\partial_{\theta}\Delta_{12}g_{m}.
$$
Thus we get for any  $s\geqslant s_{0}$
\begin{align}\label{IT diff fm+1}
\|\Delta_{12}f_{m+1}\|_{s}^{q,\kappa}&\leqslant\|({\bf B}_{m}^{-1})^{[1]}\Delta_{12}u_{m}\|_{s}^{q,\kappa}+\|(\Delta_{12}{\bf B}_{m}^{-1})u_{m}^{[2]}\|_{s}^{q,\kappa}.
\end{align}
By Lemma \ref{Compos1-lemm}, \eqref{mathscrBB1}, \eqref{control of g by f} and Lemma \ref{Law-prodX1}, we have for all $s\geqslant s_{0}$
\begin{align*}
\|({\bf B}_{m}^{-1})^{[1]}\Delta_{12}u_{m}\|_{s}^{q,\kappa}&\leqslant\|\Delta_{12}u_{m}\|_{s}^{q,\kappa}\left(1+C\|(\widehat{g}_{m}^{[1]}\|_{s_{0}}^{q,\kappa}\right)+C\|(\widehat{g}_{m}^{[1]}\|_{s}^{q,\kappa}\|\Delta_{12}u_{m}\|_{s_{0}}^{q,\kappa}\\
&\leqslant\|\Delta_{12}u_{m}\|_{s}^{q,\kappa}\left(1+C\|g_{m}^{[1]}\|_{s_{0}}^{q,\kappa}\right)+C\|g_{m}^{[1]}\|_{s}^{q,\kappa}\|\Delta_{12}u_{m}\|_{s_{0}}^{q,\kappa}\\
&\leqslant\|\Delta_{12}u_{m}\|_{s}^{q,\kappa}\left(1+C\kappa^{-1}\|f_{m}^{[1]}\|_{s_{0}+\tau_{1}q+\tau_{1}}^{q,\kappa}\right)+C\kappa^{-1}\|f_{m}^{[1]}\|_{s+\tau_{1}q+\tau_{1}}^{q,\kappa}\|\Delta_{12}u_{m}\|_{s_{0}}^{q,\kappa}.
\end{align*}
From \eqref{uniform estimate of deltams} and \eqref{small-C2}, one has
\begin{align}\label{unif fm shb+1}
\kappa^{-1}\sup_{m\in\mathbb{N}}\max_{k\in\{1,2\}}\|f_{m}^{[k]}\|_{\overline{s}_{h}+\mathtt{a}}^{q,\kappa}&\leqslant C\varepsilon\kappa^{-1}\left(1+\max_{k\in\{1,2\}}\|(\mathfrak{I}_{0}^{[k]}\|_{\overline{s}_{h}+\mathtt{a}+\sigma_{1}}^{q,\kappa}\right)\nonumber\\
&\leqslant C.
\end{align}
Thus, by \eqref{hypothesis of induction deltam} and \eqref{unif fm shb+1}, we obtain for all $s\in[s_{0},\overline{s}_{h}+\mathtt{a}]$
\begin{align}\label{Takpal1}
\|({\bf B}_{m}^{-1})^{[1]}\Delta_{12}u_{m}\|_{s}^{q,\kappa}\leqslant\|\Delta_{12}u_{m}\|_{s}^{q,\kappa}\left(1+CN_{0}^{\overline{\mu}_{2}}N_{m}^{-\overline{\mu}_{2}}\right)+CN_{m}^{\tau_{1}q+\tau_{1}}\|\Delta_{12}u_{m}\|_{s_{0}}^{q,\kappa}.
\end{align}
To estimate  $\Delta_{12}u_{m}$ we turn to \eqref{um-est} and use  Lemma \ref{Law-prodX1} in order to get  for all $s\geqslant s_{0}$
\begin{align*}
\|\Delta_{12}u_{m}\|_{s}^{q,\kappa}\leqslant \|\Pi_{N_{m}}^{\perp}\Delta_{12}f_{m}\|_{s}^{q,\kappa}&+C\|\Delta_{12}f_{m}\|_{s}^{q,\kappa}\|g_{m}^{[1]}\|_{s_{0}+1}^{q,\kappa}+C\|\Delta_{12}f_{m}\|_{s_{0}}^{q,\kappa}\|g_{m}^{[1]}\|_{s+1}^{q,\kappa}\\
&+C\|f_{m}^{[2]}\|_{s}^{q,\kappa}\|\Delta_{12}g_{m}\|_{s_{0}+1}^{q,\kappa}+C\|f_{m}^{[2]}\|_{s_{0}}^{q,\kappa}\|\Delta_{12}g_{m}\|_{s+1}^{q,\kappa}.
\end{align*}
Applying  \eqref{control of g by f} yields for all $s\geqslant s_{0}$
\begin{align*}
\|\Delta_{12}u_{m}\|_{s}^{q,\kappa}\leqslant&\, \|\Pi_{N_{m}}^{\perp}\Delta_{12}f_{m}\|_{s}^{q,\kappa}+C\kappa^{-1}\|\Delta_{12}f_{m}\|_{s}^{q,\kappa}\|f_{m}^{[1]}\|_{s_{0}+\tau_{1}q+\tau_{1}+1}^{q,\kappa}\\
&+C\kappa^{-1}\|\Delta_{12}f_{m}\|_{s_{0}}^{q,\kappa}\|\Pi_{N_m}f_{m}^{[1]}\|_{s+\tau_{1}q+\tau_{1}+1}^{q,\kappa}+C\|f_{m}^{[2]}\|_{s}^{q,\kappa}\|\Delta_{12}g_{m}\|_{s_{0}+1}^{q,\kappa}\\
&\quad+C\|f_{m}^{[2]}\|_{s_{0}}^{q,\kappa}\|\Delta_{12}g_{m}\|_{s+1}^{q,\kappa}.
\end{align*}
Combining this estimate with  \eqref{estimate delta12 gm}, we finally obtain for all $s\geqslant s_{0}$
\begin{align*}
\|\Delta_{12}u_{m}\|_{s}^{q,\kappa}&\leqslant \|\Pi_{N_{m}}^{\perp}\Delta_{12}f_{m}\|_{s}^{q,\kappa}+C\kappa^{-1}N_m^{\tau_1 q+\tau_1+1}\|\Delta_{12}f_{m}\|_{s}^{q,\kappa}\max_{k\in\{1,2\}}\|f_{m}^{[k]}\|_{s_{0}}^{q,\kappa}\\
&+C\kappa^{-1}N_{m}^{\tau_{1}q+\tau_{1}+1}\|\Delta_{12}f_{m}\|_{s_{0}}^{q,\kappa}\max_{k\in\{1,2\}}\|f_{m}^{[k]}\|_{s}^{q,\kappa}\\
&+C\kappa^{-1}N_{m}^{\tau_{1}q+2\tau_{1}+1}\max_{k\in\{1,2\}}\|f_{m}^{[k]}\|_{s}^{q,\kappa}\max_{k\in\{1,2\}}\|f_{m}^{[k]}\|_{s_{0}}^{q,\kappa}\|\Delta_{12}V_{m}\|^{q,\kappa}.
\end{align*}
For $s=s_0$ we can actually get a better estimate,
\begin{align*}
\|\Delta_{12}u_{m}\|_{s_0}^{q,\kappa}&\leqslant \|\Pi_{N_{m}}^{\perp}\Delta_{12}f_{m}\|_{s_0}^{q,\kappa}+C\kappa^{-1}N_m^{\tau_1 q+\tau_1+1}\|\Delta_{12}f_{m}\|_{s_0}^{q,\kappa}\max_{k\in\{1,2\}}\|f_{m}^{[k]}\|_{s_{0}}^{q,\kappa}\\
&+C\kappa^{-1}N_m^{\tau_1q+2\tau_1+1}\max_{k\in\{1,2\}}\|f_{m}^{[k]}\|_{s_0}^{q,\kappa}\max_{k\in\{1,2\}}\|f_{m}^{[k]}\|_{s_{0}}^{q,\kappa}\|\Delta_{12}V_{m}\|^{q,\kappa}.
\end{align*}
By \eqref{hypothesis of induction deltam}, Lemma \ref{orthog-Lem1} and \eqref{unif fm shb+1},  using in particular that $s_{0}+\tau_1q+2\tau_1+1\leqslant \overline{s}_h$, we deduce
\begin{align}\label{foufg90}
\nonumber\|\Delta_{12}u_{m}\|_{s_{0}}^{q,\kappa}&\leqslant N_{m}^{s_{0}-\overline{s}_{h}-\mathtt{a}}\|\Delta_{12}f_{m}\|_{\overline{s}_{h}+\mathtt{a}}^{q,\kappa}+CN_{0}^{\overline{\mu}_{2}}N_{m}^{\tau_{1}q+\tau_{1}+1-\overline{\mu}_{2}}\delta_{0}^{[1,2]}(\overline{s}_{h})\|\Delta_{12}f_{m}\|_{s_{0}}^{q,\kappa}\\
&\qquad+CN_m^{\tau_1q+2\tau_1+1}N_{0}^{2\overline{\mu}_{2}}N_{m}^{-2\overline{\mu}_{2}}\delta_{0}^{[1,2]}(\overline{s}_{h})\|\Delta_{12}V_{m}\|^{q,\kappa}
\end{align}
and
\begin{align*}
\|\Delta_{12}u_{m}\|_{\overline{s}_{h}+\mathtt{a}}^{q,\kappa}&\leqslant\|\Delta_{12}f_{m}\|_{\overline{s}_{h}+\mathtt{a}}^{q,\kappa}\left(1+CN_{0}^{\mu_{2}}N_{m}^{\tau_{1}q+\tau_{1}+1-\overline{\mu}_{2}}\delta_{0}^{[1,2]}(\overline{s}_{h})\right)+CN_{m}^{\tau_{1}q+\tau_{1}+1}\delta_{0}^{[1,2]}(\overline{s}_{h})\|\Delta_{12}f_{m}\|_{s_{0}}^{q,\kappa}\\
&\qquad +CN_{0}^{\overline{\mu}_{2}}N_{m}^{\tau_{1}q+2\tau_{1}+1-\overline{\mu}_{2}}\delta_{0}^{[1,2]}(\overline{s}_{h})\|\Delta_{12}V_{m}\|^{q,\kappa}
\end{align*}
where
$$\delta_{0}^{[1,2]}(s)\triangleq\kappa^{-1}\max_{k\in\{1,2\}}\|f_{0}^{[k]}\|_{s}^{q,\kappa}.
$$
Hence, inserting the  estimate \eqref{foufg90} into \eqref{Takpal1} implies
\begin{align}\label{Gm-1 diff s0}
\nonumber\|({\bf B}_{m}^{-1})^{[1]}\Delta_{12}u_{m}\|_{s_{0}}^{q,\kappa}&\leqslant CN_m^{\tau_1 q+\tau_1}\|\Delta_{12}u_{m}\|_{s_{0}}^{q,\kappa}
\\
\nonumber &\leqslant  N_{m}^{s_{0}-\overline{s}_{h}+\tau_1 q+\tau_1}\|\Delta_{12}f_{m}\|_{\overline{s}_{h}}^{q,\kappa}+CN_{0}^{\overline{\mu}_{2}}N_{m}^{2(\tau_{1}q+\tau_{1})+1-\overline{\mu}_{2}}\delta_{0}^{[1,2]}(\overline{s}_{h})\|\Delta_{12}f_{m}\|_{s_{0}}^{q,\kappa}\\
&\qquad+CN_{0}^{2\overline{\mu}_{2}}N_{m}^{-2\overline{\mu}_{2}+\tau_1 q+\tau_1}\delta_{0}^{[1,2]}(\overline{s}_{h})\|\Delta_{12}V_{m}\|^{q,\kappa}.
\end{align}
Similarly and after  straightforward computations we get
\begin{align}\label{Gm-1 diff shb}
\|({\bf B}_{m}^{-1})^{[1]}\Delta_{12}u_{m}\|_{\overline{s}_{h}+\mathtt{a}}^{q,\kappa}\leqslant&\|\Delta_{12}f_{m}\|_{\overline{s}_{h}+\mathtt{a}}^{q,\kappa}\left(1+N_{m}^{s_0+\tau_{1}q+\tau_{1}-\overline{s}_h}+N_{0}^{\overline{\mu}_{2}}N_{m}^{\tau_{1}q+\tau_{1}+1-\mu_{2}}\delta_{0}^{[1,2]}(\overline{s}_{h})\right)\nonumber\\
&+C\big(N_0^{\overline{\mu}_2}N_m^{2(\tau_1 q+\tau_1)+1-\overline{\mu}_2}+N_{m}^{\tau_{1}q+\tau_{1}+1}\big)\delta_{0}^{[1,2]}(\overline{s}_{h})\|\Delta_{12}f_{m}\|_{s_{0}}^{q,\kappa}\nonumber\\
&\quad+CN_{0}^{\overline{\mu}_{2}}N_{m}^{\tau_{1}q+\tau_{1}+1-\overline{\mu}_{2}}\delta_{0}^{[1,2]}(\overline{s}_{h})\|\Delta_{12}V_{m}\|^{q,\kappa}.
\end{align} 
Next we intend to estimate $(\Delta_{12}{\bf B}_{m}^{-1})u_{m}^{[2]}.$ For this aim we write according to 
 Taylor formula, 
$$(\Delta_{12}{\bf B}_{m}^{-1})u_{m}^{[2]}(\theta)=\Delta_{12}\widehat{g}_{m}(\theta)\int_{0}^{1}\big(\partial_{\theta}u_{m}^{[2]}\big)\Big(\theta+\widehat{g}_{m}^{[2]}(\theta)+t\Delta_{12}\widehat{g}_{m}(\theta)\Big)dt.$$
Coming back to \eqref{um-est} and using   Lemma \ref{Law-prodX1}, \eqref{control of g by f} and \eqref{hypothesis of induction deltam}, we have for all $s\geqslant s_{0}$
\begin{align}\label{link um-fm}
\|u_{m}^{[2]}\|_{s}^{q,\kappa}\leqslant&\|\Pi_{N_{m}}^{\perp}f_{m}^{[2]}\|_{s}^{q,\kappa}+C\|f_{m}^{[2]}\|_{s}^{q,\kappa}\|\partial_{\theta}g_{m}^{[2]}\|_{s_{0}}^{q,\kappa}+C\|f_{m}^{[2]}\|_{s_{0}}^{q,\kappa}\|\partial_{\theta}g_{m}^{[2]}\|_{s}^{q,\kappa}\nonumber\\
&\leqslant \|f_{m}^{[2]}\|_{q,s}\left(1+C\kappa^{-1}\|f_{m}^{[2]}\|_{s_{0}+\tau_{1}q+\tau_{1}+1}^{q,\kappa}\right)\nonumber\\
&\quad\leqslant C\|f_{m}^{[2]}\|_{s}^{q,\kappa}.
\end{align}
Applying  Lemma \ref{Law-prodX1} and  Lemma \ref{Compos1-lemm}, we obtain for all $s\geqslant s_{0}$
\begin{align*}
\|(\Delta_{12}{\bf B}_{m}^{-1})u_{m}^{[2]}\|_{s}^{q,\kappa}\leqslant &C\|\Delta_{12}\widehat{g}_{m}\|_{s}^{q,\kappa}\|u_{m}^{[2]}\|_{s_{0}+1}^{q,\kappa}\left(1+\|\widehat{g}_{m}^{[2]}\|_{s_{0}}^{q,\kappa}+\|\Delta_{12}\widehat{g}_{m}\|_{s_{0}}^{q,\kappa}\right)\\
&+C\|\Delta_{12}\widehat{g}_{m}\|_{s_{0}}^{q,\kappa}\|u_{m}^{[2]}\|_{s+1}^{q,\kappa}\left(1+\|\widehat{g}_{m}^{[2]}\|_{s_{0}}^{q,\kappa}+\|\Delta_{12}\widehat{g}_{m}\|_{s_{0}}^{q,\kappa}\right)\\
&\quad+C\|\Delta_{12}\widehat{g}_{m}\|_{s_{0}}^{q,\kappa}\|u_{m}^{[2]}\|_{s_{0}+1}^{q,\kappa}\left(\|\widehat{g}_{m}^{[2]}\|_{s}^{q,\kappa}+\|\Delta_{12}\widehat{g}_{m}\|_{s}^{q,\kappa}\right).
\end{align*}
From \eqref{link diff beta hat and diff beta}, \eqref{deacr gm} and Sobolev embeddings, one obtains for all $s\in[s_{0},\overline{s}_{h}+\mathtt{a}]$
\begin{align*}
\|\Delta_{12}\widehat{g}_{m}\|_{s}^{q,\kappa}&\leqslant C\left(\|\Delta_{12}g_{m}\|_{s}^{q,\kappa}+\|\Delta_{12}g_{m}\|_{s_{0}}^{q,\kappa}\max_{k\in\{1,2\}}\|g_{m}^{[k]}\|_{s+1}^{q,\kappa}\right)\\
&\leqslant C\|\Delta_{12}g_{m}\|_{s}^{q,\kappa}.
\end{align*}
Combining the last two estimates with \eqref{maj gmrk} and Lemma \ref{Compos1-lemm} we deduce that for all $s\in[s_{0},\overline{s}_{h}+\mathtt{a}]$
\begin{align*}
\|(\Delta_{12}{\bf B}_{m}^{-1})u_{m}^{[2]}\|_{s}^{q,\kappa}&\leqslant C\|\Delta_{12}g_{m}\|_{s}^{q,\kappa}\|u_{m}^{[2]}\|_{q,s_{0}+1}+C\|\Delta_{12}g_{m}\|_{s_{0}}^{q,\kappa}\|u_{m}^{[2]}\|_{s+1}^{q,\kappa}.
\end{align*}
Using \eqref{estimate delta12 gm}, \eqref{link um-fm}, we obtain for all $s\in[s_{0},\overline{s}_{h}+\mathtt{a}]$
\begin{align*}
\|(\Delta_{12}{\bf B}_{m}^{-1})u_{m}^{[2]}\|_{s}^{q,\kappa}\leqslant& C\kappa^{-1}N_{m}^{\tau_{1}q+\tau_{1}}\|\Delta_{12}f_{m}\|_{s}^{q,\kappa}\max_{k\in\{1,2\}}\|f_{m}^{[k]}|_{s_{0}+1}^{q,\kappa}\\
&+C\kappa^{-1}N_{m}^{\tau_{1}q+\tau_{1}}\|\Delta_{12}f_{m}\|_{s_{0}}^{q,\kappa}\max_{k\in\{1,2\}}\|f_{m}^{[k]}\|_{s+1}^{q,\kappa}\\
&\quad+C\kappa^{-1}N_{m}^{\tau_{1}q+2\tau_{1}}\|\Delta_{12}V_{m}\|^{q,\kappa}\max_{k\in\{1,2\}}\|f_{m}^{[k]}\|_{s_{0}+1}^{q,\kappa}\max_{k\in\{1,2\}}\|f_{m}^{[k]}\|_{s+1}^{q,\kappa}.
\end{align*}
Hence, by \eqref{hypothesis of induction deltam}, \eqref{uniform estimate of deltams} and \eqref{unif fm shb+1}, we have
\begin{align}\label{diff Gm-1 s0}
\|(\Delta_{12}{\bf B}_{m}^{-1})u_{m}^{[2]}\|_{s_0}^{q,\kappa}&\leqslant CN_{0}^{\overline{\mu}_{2}}N_{m}^{\tau_{1}q+\tau_{1}-\overline{\mu}_{2}}\delta_{0}^{[1,2]}(\overline{s}_{h}+1)\|\Delta_{12}f_{m}\|_{s_{0}}^{q,\kappa}\\
\nonumber&+N_{0}^{2\overline{\mu}_{2}}N_{m}^{\tau_1 q+2\tau_1-2\overline{\mu}_{2}}\delta_{0}^{[1,2]}(\overline{s}_{h}+1)\|\Delta_{12}V_{m}\|^{q,\kappa}
\end{align}
and
\begin{align}\label{diff Gm-1 shb}
\nonumber\|(\Delta_{12}{\bf B}_{m}^{-1})u_{m}^{[2]}\|_{\overline{s}_h+\mathtt{a}}^{q,\kappa}\leqslant &CN_{0}^{\overline{\mu}_{2}}N_{m}^{\tau_{1}q+\tau_{1}-\overline{\mu}_{2}}\delta_{0}^{[1,2]}(\overline{s}_{h}+1)\|\Delta_{12}f_{m}\|_{\overline{s}_{h}+\mathtt{a}}^{q,\kappa}\\
 &+CN_{m}^{\tau_{1}q+\tau_{1}}\delta_{0}^{[1,2]}(\overline{s}_{h}+\mathtt{a}+1)\|\Delta_{12}f_{m}\|_{s_{0}}^{q,\kappa}\\
\nonumber&\quad +N_{0}^{\overline{\mu}_{2}}N_{m}^{\tau_1 q+2\tau_1-\overline{\mu}_{2}}\delta_{0}^{[1,2]}(\overline{s}_{h}+\mathtt{a}+1)\|\Delta_{12}V_{m}\|^{q,\kappa}.
\end{align}
Putting together \eqref{IT diff fm+1}, \eqref{Gm-1 diff s0} and \eqref{diff Gm-1 s0} and using Sobolev embeddings  yield
\begin{align}\label{link diff fm+1 and fm s0}
\nonumber\|\Delta_{12}f_{m+1}\|_{s_{0}}^{q,\kappa} \leqslant &N_{m}^{s_{0}-\overline{s}_{h}+\tau_1 q+\tau_1}\|\Delta_{12}f_{m}\|_{\overline{s}_{h}+\mathtt{a}}^{q,\kappa}\\
&+CN_{0}^{\overline{\mu}_{2}}N_{m}^{2(\tau_{1}q+\tau_{1})+1-\overline{\mu}_{2}}\delta_{0}^{[1,2]}(\overline{s}_{h}+1)\|\Delta_{12}f_{m}\|_{s_{0}}^{q,\kappa}\\
\nonumber&\qquad+CN_{0}^{2\overline{\mu}_{2}}N_{m}^{-2\overline{\mu}_{2}+\tau_1 q+2\tau_1}\delta_{0}^{[1,2]}(\overline{s}_{h}+1)\|\Delta_{12}V_{m}\|^{q,\kappa}.
\end{align}
Similarly, we obtain according to \eqref{IT diff fm+1}, \eqref{Gm-1 diff shb} and \eqref{diff Gm-1 shb} 
\begin{align}\label{link diff fm+1 and fm shb}
\nonumber\|\Delta_{12}f_{m+1}\|_{\overline{s}_{h}+\mathtt{a}}^{q,\kappa}\leqslant& \|\Delta_{12}f_{m}\|_{\overline{s}_{h}+\mathtt{a}}^{q,\kappa}\left(1+N_{m}^{s_0+\tau_{1}q+\tau_{1}-\overline{s}_h}+N_{0}^{\overline{\mu}_{2}}N_{m}^{\tau_{1}q+\tau_{1}+1-\overline{\mu}_{2}}\delta_{0}^{[1,2]}(\overline{s}_{h}+1)\right)\\
 &+C\big(N_0^{\overline{\mu}_2}N_m^{2(\tau_1 q+\tau_1)+1-\overline{\mu}_2}+N_{m}^{\tau_{1}q+\tau_{1}+1}\big)\delta_{0}^{[1,2]}(\overline{s}_{h}+\mathtt{a}+1)\|\Delta_{12}f_{m}\|_{s_{0}}^{q,\kappa}\\
\nonumber&\quad +N_{0}^{\overline{\mu}_{2}}N_{m}^{\tau_1 q+2\tau_1-\overline{\mu}_{2}}\delta_{0}^{[1,2]}(\overline{s}_{h}+\mathtt{a}+1)\|\Delta_{12}V_{m}\|^{q,\kappa}.
\end{align}
Next, we  introduce the quantities
$$\overline{\delta}_{m}(s)\triangleq \kappa^{-1}\|\Delta_{12}f_{m}\|_{s}^{q,\kappa}\quad\textnormal{and}\quad\varkappa_{m}\triangleq\kappa^{-1}\|\Delta_{12}V_{m}\|^{q,\kappa}.$$
According to \eqref{definition Vm+1 and fm+1} one deduces that
$$\Delta_{12}V_{m+1}=\Delta_{12}V_{m}+\langle\Delta_{12}f_{m}\rangle_{\varphi,\theta}\quad\hbox{and}\quad \Delta_{12}V_{0}=0,
$$
which implies in view of   Sobolev embeddings
\begin{equation}\label{link varkappa and overline delta}
\varkappa_{m}\leqslant C \sum_{k=0}^{m-1}\overline{\delta}_{k}(s_{0}).
\end{equation}
The goal now is to  prove by induction that for any $0\leqslant \mathtt{a}\leqslant s_h-\overline{s}_{h}+{\sigma_1-5}$
\begin{equation}\label{hypothesis of induction overline delta}
\forall k\leqslant m,\quad\overline{\delta}_{k}(s_{0})\leqslant N_{0}^{\overline{\mu}_{2}}N_{k}^{-\overline{\mu}_{2}}\nu(\overline{s}_{h}+\mathtt{a})\quad\textnormal{and}\quad\overline{\delta}_{k}(\overline{s}_{h}+\mathtt{a})\leqslant\left(2-\frac{1}{k+1}\right)\nu(\overline{s}_{h}+\mathtt{a})
\end{equation}
with
$$\nu(s)\triangleq\overline{\delta}_{0}(s)+\varepsilon\kappa^{-1}\|\Delta_{12}i\|_{s_{0}+{4}}.
$$
Notice that the property \eqref{hypothesis of induction overline delta} is obvious for $m=0$ due to Sobolev embeddings. Let us now assume that \eqref{hypothesis of induction overline delta} is true at the order $m$ and let us check it at the order $m+1.$ By hypothesis of induction \eqref{hypothesis of induction overline delta} and \eqref{link varkappa and overline delta}, one has
\begin{align}\label{mah00}
\varkappa_{m}\leqslant C\nu(\overline{s}_{h}+\mathtt{a}),
\end{align}
with $C$ independent of $m.$
By \eqref{link diff fm+1 and fm s0} and hypothesis of induction \eqref{hypothesis of induction overline delta}, we can write
\begin{align*}
\overline{\delta}_{m+1}(s_{0})&\leqslant N_{m}^{s_{0}-\overline{s}_{h}+{\tau_1q+\tau_1}}\overline{\delta}_{m}(\overline{s}_{h}+\mathtt{a})+CN_{0}^{\overline{\mu}_{2}}N_{m}^{2(\tau_{1}q+\tau_{1})+1-\overline{\mu}_{2}}\delta_{0}^{[1,2]}(\overline{s}_{h}+1)\overline{\delta}_{m}(s_{0})\\
&\qquad\qquad +CN_{0}^{2\overline{\mu}_{2}}N_{m}^{-2\overline{\mu}_{2}+\tau_1 q+2\tau_1}\delta_{0}^{[1,2]}(\overline{s}_{h}+1)\varkappa_{m}\\
&\leqslant \left[2N_{m}^{s_{0}-\overline{s}_{h}+{\tau_1q+\tau_1}}+CN_{0}^{2\overline{\mu}_{2}}N_{m}^{2(\tau_{1}q+\tau_{1})+1-2\overline{\mu}_{2}}\delta_{0}^{[1,2]}(\overline{s}_{h}+1)\right]\nu(\overline{s}_{h}+\mathtt{a}).
\end{align*}
From the constraints  on $\overline{s}_{h}$ and $\overline{\mu}_{2}$ fixed by  \eqref{Conv-T2N} and \eqref{definition of Nm} we get 
\begin{align*}
2N_{m}^{s_{0}-\overline{s}_{h}+{\tau_1q+\tau_1}}=& 2N_{m}^{-\frac32\overline{\mu}_{2}-2}=2N_m^{-2}N_{m+1}^{-\overline{\mu}_{2}}\\
&\leqslant 2N_0^{-2}N_{m+1}^{-\overline{\mu}_{2}}\\
&\quad\leqslant \frac12 N_0^{\overline{\mu}_2}N_{m+1}^{-\overline{\mu}_{2}}.
\end{align*}
The last inequality occurs provided that 
$$
4\leqslant N_0^{\overline{\mu}_2+2}
$$
which is true since $N_0\geqslant 2$ and $\overline\mu_2\geqslant 2$. Similarly, one gets from \eqref{estimate delta0 and I0} and \eqref{Conv-T2N} 
\begin{align*}
CN_{0}^{2\overline{\mu}_{2}}N_{m}^{2(\tau_{1}q+\tau_{1})+1-2\overline{\mu}_{2}}\delta_{0}^{[1,2]}(\overline{s}_{h}+1)\leqslant &C\varepsilon\kappa^{-1}N_{0}^{\overline{\mu}_{2}}N_{m}^{2(\tau_{1}q+\tau_{1})+1-\frac12\overline{\mu}_{2}}N_{0}^{\overline{\mu}_{2}}N_{m+1}^{-\overline{\mu}_2}\\
&\leqslant C\varepsilon\kappa^{-1}N_{0}^{2(\tau_{1}q+\tau_{1})+1+\frac12\overline{\mu}_{2}}N_{0}^{\overline{\mu}_{2}}N_{m+1}^{-\overline{\mu}_2}\\
 &\quad \leqslant C\varepsilon\kappa^{-1}N_{0}^{\overline{\mu}_{2}}N_{0}^{\overline{\mu}_{2}}N_{m+1}^{-\overline{\mu}_2}.
\end{align*}
Therefore, from \eqref{small-C2} and by taking $\varepsilon_0$ small enough we infer
\begin{align*}
CN_{0}^{2\overline{\mu}_{2}}N_{m}^{2(\tau_{1}q+\tau_{1})+1-2\overline{\mu}_{2}}\delta_{0}^{[1,2]}(\overline{s}_{h}+1)&
 \leqslant \frac12N_{0}^{\overline{\mu}_{2}}N_{m+1}^{-\overline{\mu}_2}.
\end{align*}
Putting together the preceding estimates we find that
\begin{align*}
\overline{\delta}_{m+1}(s_{0})&\leqslant N_{0}^{\overline{\mu}_{2}}N_{m+1}^{-\overline{\mu}_2}\nu(\overline{s}_{h}+\mathtt{a}).
\end{align*}
This achieves the induction of the first statement in \eqref{hypothesis of induction overline delta}. Now,  let us move to the second one. Then  we can write by virtue of \eqref{link diff fm+1 and fm shb}, the  first statement of \eqref{hypothesis of induction overline delta}, \eqref{estimate delta0 and I0} and \eqref{small-C2}
\begin{align*}
\nonumber\overline{\delta}_{m+1}(\overline{s}_{h}+\mathtt{a})\leqslant& \overline{\delta}_{m}(\overline{s}_{h}+\mathtt{a})\left(1+N_{m}^{s_0+\tau_{1}q+\tau_{1}-\overline{s}_h}+N_{0}^{\overline{\mu}_{2}}N_{m}^{\tau_{1}q+\tau_{1}+1-\overline{\mu}_{2}}\right)\\
 &+C\big(N_0^{\overline{\mu}_2}N_m^{2(\tau_1 q+\tau_1)+1-\overline{\mu}_2}+N_{m}^{\tau_{1}q+\tau_{1}+1}\big)N_{0}^{\overline{\mu}_{2}}N_{m}^{-\overline{\mu}_2}\varepsilon\kappa^{-1}\nu(\overline{s}_h+\mathtt{a}).
\end{align*}
Notice that we have used according to \eqref{estimate delta0 and I0} and \eqref{small-C2} that if $\overline{s}_{h}+\mathtt{a}+5\leqslant s_h+\sigma_1$
\begin{align*}
\delta_{0}^{[1,2]}(\overline{s}_{h}+\mathtt{a}+1)&\lesssim \varepsilon\kappa^{-1}\left(1+\max_{k=1,2}\|\mathfrak{I}_{k}\|_{\overline{s}_{h}+\mathtt{a}+{5}}^{q,\kappa}\right)\\
&\lesssim \varepsilon\kappa^{-1}.
\end{align*} 
By hypothesis of induction \eqref{hypothesis of induction overline delta}, we obtain
\begin{align*}
\nonumber\overline{\delta}_{m+1}(\overline{s}_{h}+\mathtt{a})\leqslant& \left(2-\frac{1}{m+1}\right)\left(1+N_{m}^{s_0+\tau_{1}q+\tau_{1}-\overline{s}_h}+N_{0}^{\overline{\mu}_{2}}N_{m}^{\tau_{1}q+\tau_{1}+1-\overline{\mu}_{2}}\right)\nu(\overline{s}_h+\mathtt{a})\\
 &+C\big(N_0^{\overline{\mu}_2}N_m^{2(\tau_1 q+\tau_1)+1-\overline{\mu}_2}+N_{m}^{\tau_{1}q+\tau_{1}+1}\big)N_{0}^{\overline{\mu}_{2}}N_{m}^{-\overline{\mu}_2}\varepsilon\kappa^{-1}\nu(\overline{s}_h+\mathtt{a}).
\end{align*}
In the same way to \eqref{conv-t2}  we conclude in view of \eqref{small-C2}
\begin{align*}
 &\left(2-\frac{1}{m+1}\right)\left(1+N_{m}^{s_0+\tau_{1}q+\tau_{1}-\overline{s}_h}+N_{0}^{\overline{\mu}_{2}}N_{m}^{\tau_{1}q+\tau_{1}+1-\overline{\mu}_{2}}\right)\\
 &+C\big(N_0^{\overline{\mu}_2}N_m^{2(\tau_1 q+\tau_1)+1-\overline{\mu}_2}+N_{m}^{\tau_{1}q+\tau_{1}+1}\big)N_{0}^{\overline{\mu}_{2}}N_{m}^{-\overline{\mu}_2}\varepsilon\kappa^{-1}\\
 &\leqslant \left(2-\frac{1}{m+2}\right).
\end{align*}
Consequently, we find 
$$\overline{\delta}_{m+1}(\overline{s}_{h}+\mathtt{a})\leqslant\left(2-\frac{1}{m+2}\right)\nu(\overline{s}_{h}+\mathtt{a})
$$
which achieves the proof of the second statement in \eqref{hypothesis of induction overline delta} at the order $m+1$.

\smallskip

$\blacktriangleright$ \textbf{Estimate on $\Delta_{12}c_{i_{0}}.$} According to \eqref{struc-ci} and since $V_0$ is independent of $r$  then
$$\Delta_{12}c_{i_{0}}=\sum_{m=0}^{+\infty}\Delta_{12}(V_{m+1}-V_{m}).$$
By \eqref{definition Vm+1 and fm+1}, Sobolev embeddings and implementing  \eqref{hypothesis of induction overline delta}with  $\mathtt{a}=0$ 
\begin{align*}
\|\Delta_{12}(V_{m+1}-V_{m})\|^{q,\kappa}&=\|\langle \Delta_{12}f_{m}\rangle_{\varphi,\theta}\|^{q,\kappa}\\
&\leqslant C\kappa\overline{\delta}_{m}(s_{0})\\
&\leqslant C\kappa N_{0}^{\overline{\mu}_{2}}N_{m}^{-\overline{\mu}_{2}}\nu(\overline{s}_{h}).
\end{align*}
It follows by the preceding estimates combined with   Lemma \ref{lemma sum Nn} and \eqref{ain-z1} 
\begin{align*}
\|\Delta_{12}c_{i_{0}}\|^{q,\kappa}&\leqslant \sum_{m=0}^{\infty}\|\Delta_{12}(V_{m+1}-V_{m})\|^{q,\kappa}\\
& \leqslant  C\kappa\, \nu(\overline{s}_{h})N_{0}^{\underline{\mu_{2}}}\sum_{m=0}^{+\infty}N_{m}^{-\underline{\mu_{2}}}\\
& \leqslant C\varepsilon\|\Delta_{12}i\|_{\overline{s}_{h}+4}^{q,\kappa}.
\end{align*}

\smallskip

$\blacktriangleright$ \textbf{Estimates of $\Delta_{12}{\beta}$ and $\Delta_{12}\widehat{\beta}.$}
Now, inserting the estimates of  \eqref{hypothesis of induction overline delta} into \eqref{estimate delta12 gm}, taken with $s=s_{0},$ we obtain 
\begin{align*}
	\|\Delta_{12}g_{m}\|_{s_{0}}^{q,\kappa}&\lesssim \overline{\delta}_{m}(s_{0}+\tau_{1}q+\tau_{1})+\varkappa_{m}\delta_{m}(s_{0}+\tau_{1}q+2\tau_{1}+1).
\end{align*}
Applying Sobolev embeddings,  Lemma \ref{interpolation-In}, \eqref{hypothesis of induction deltam} and  \eqref{hypothesis of induction overline delta}  we get for some   $\overline{\vartheta}\in(0,1)$
\begin{align*}
	\overline{\delta}_{m}(s_{0}+\tau_{1}q+\tau_{1})&\leqslant  \overline{\delta}_{m}(s_{0}+\tau_{1}q+2\tau_{1}+1)\\
	&\lesssim\overline{\delta}_{m}(s_{0})^{\overline{\vartheta}}\overline{\delta}_{m}(\overline{s}_{h})^{1-\overline{\vartheta}}\\
	&\lesssim N_{0}^{\overline{\vartheta}\mu_{2}}N_{m}^{-\overline{\vartheta}\mu_{2}}\nu(\overline{s}_{h})
\end{align*}
and
\begin{align*}
	\delta_{m}(s_{0}+\tau_{1}q+2\tau_{1}+1)&\lesssim \delta_{m}(s_{0})^{\overline{\vartheta}}\delta_{m}(\overline{s}_{h})^{1-\overline{\vartheta}}\\
	&\lesssim N_{0}^{\overline{\vartheta}\overline{\mu}_{2}}N_{m}^{-\overline{\vartheta}\overline{\mu}_{2}}\delta_{0}(\overline{s}_{h})\\
	&\lesssim N_{0}^{\overline{\vartheta}\overline{\mu}_{2}}N_{m}^{-\overline{\vartheta}\overline{\mu}_{2}}.
\end{align*}
Thus we deduce from the preceding estimates and \eqref{mah00}
\begin{align}\label{ain-b1}
\|\Delta_{12}g_{m}\|_{s_{0}}^{q,\kappa}\lesssim N_{0}^{\overline{\vartheta}\overline{\mu}_{2}}N_{m}^{-\overline{\vartheta}\overline{\mu}_{2}}\nu(\overline{s}_{h}).
\end{align}
Now from \eqref{estimate delta12 gm}  we infer
\begin{align}\label{ain-m1}
\|\Delta_{12}g_{m}\|_{\overline{s}_{h}+\mathtt{a}+1}^{q,\kappa}&\lesssim \overline{\delta}_{m}(\overline{s}_{h}+\mathtt{a}+\tau_{1}q+\tau_{1}+1)+\varkappa_{m}{{\delta}_{m}(\overline{s}_{h}+\mathtt{a}+\tau_{1}q+2\tau_{1}+2)}.
\end{align}
Coming back to  \eqref{hypothesis of induction overline delta} we infer that  under the assumption 
\begin{equation}\label{asump-ain}
\overline{s}_{h}+\mathtt{a}+\tau_{1}q+\tau_{1}+1\leqslant s_h+\sigma_1-5
\end{equation}
one gets 
\begin{align*}
\overline{\delta}_{m}(\overline{s}_{h}+\mathtt{a}+\tau_{1}q+\tau_{1}+1)&\leqslant 2\nu(\overline{s}_{h}+\mathtt{a}+\tau_{1}q+\tau_{1}+1)\\&\leqslant 2\overline{\delta}_{0}(\overline{s}_{h}+\mathtt{a}+\tau_{1}q+\tau_{1}+1)+2\varepsilon\kappa^{-1}\|\Delta_{12}i\|_{s_{0}+{4}}^{q,\kappa}.
\end{align*}
Similarly to \eqref{estimate delta0 and I0}, one gets from Lemma \ref{lemma-reste}-(iii) and Proposition \ref{lemma-GS0}-(i)
\begin{align}\label{ain-z1}
\nonumber\forall s\geqslant s_0,\quad \overline{\delta}_{0}(s)&=\kappa^{-1}\|\Delta_{12}V_{\varepsilon r}\|_{s}^{q,\kappa}\nonumber\\
&\lesssim\varepsilon\kappa^{-1}\Big(\|\Delta_{12} i\|_{s+4}^{q,\kappa}+{ \|\Delta_{12} i\|_{s_0+4}^{q,\kappa}
\max_{\ell=1,2} \|r_\ell\|_{s+4}^{q,\kappa}.}
\Big).
\end{align}
Then under the assumption 
\begin{equation}\label{asump-ain1}
{\overline{s}_{h}+\mathtt{a}+\tau_{1}q+2\tau_{1}+6\leqslant s_h+\sigma_1}
\end{equation}
combined with   \eqref{small-C2}  and \eqref{ain-z1} we find
\begin{align*}
\overline{\delta}_{0}(\overline{s}_{h}+\mathtt{a}+\tau_{1}q+\tau_{1}+1)
&\lesssim\varepsilon\kappa^{-1}\|\Delta_{12} i\|_{\overline{s}_{h}+\mathtt{a}+\tau_{1}q+\tau_{1}+5}^{q,\kappa}.
\end{align*}
From \eqref{uniform estimate of deltams}, \eqref{small-C2} and \eqref{asump-ain1} we obtain the following estimate
\begin{align*}
\delta_{m}((\overline{s}_{h}+\mathtt{a}+\tau_{1}q+2\tau_{1}+2)& \leqslant C\varepsilon\kappa^{-1}\left(1+\|\mathfrak{I}_{0}\|_{\overline{s}_{h}+\mathtt{a}+\tau_{1}q+2\tau_{1}+6}^{q,\kappa}\right)\\
& \leqslant C\varepsilon\kappa^{-1}.
\end{align*}
Plugging the preceding estimates into \eqref{ain-m1} and using \eqref{mah00}  yield
\begin{align}\label{ain-mm1}
\|\Delta_{12}g_{m}\|_{\overline{s}_{h}+\mathtt{a}+1}^{q,\kappa}&\lesssim \varepsilon\kappa^{-1}\|\Delta_{12} i\|_{\overline{s}_{h}+\mathtt{a}+\tau_{1}q+\tau_{1}+5}^{q,\kappa}.
\end{align}
provided that \eqref{asump-ain1} is satisfied.
By virtue of Lemma \ref{interpolation-In}, \eqref{ain-mm1} and \eqref{ain-b1} we obtain for some  $\overline{\theta}\in(0,1)$
\begin{equation}\label{final dlt12 gm shb+1}
\|\Delta_{12}g_{m}\|_{\overline{s}_{h}+\mathtt{a}}^{q,\kappa}\lesssim\varepsilon \kappa^{-1} N_{0}^{\overline{\theta}\underline{\mu_{2}}}N_{m}^{-\overline{\theta}\underline{\mu_{2}}}\|\Delta_{12} i\|_{\overline{s}_{h}+\mathtt{a}+\tau_{1}q+\tau_{1}+5}^{q,\kappa}.
\end{equation}
Consequently, we deduce from \eqref{final dlt12 gm shb+1}   and Lemma \ref{lemma sum Nn}
\begin{align}\label{final sum dlt12 gm shb+1}
\sum_{m=0}^{\infty}\|\Delta_{12}g_{m}\|_{\overline{s}_{h}+\mathtt{a}}^{q,\kappa}&\lesssim \varepsilon \kappa^{-1}\|\Delta_{12} i\|_{\overline{s}_{h}+\mathtt{a}+\tau_{1}q+\tau_{1}+5}^{q,\kappa}N_{0}^{\overline{\theta}\underline{\mu_{2}}}\sum_{m=0}^{\infty}N_{m}^{-\overline{\theta}\underline{\mu_{2}}}\nonumber\\
&\lesssim\varepsilon\kappa^{-1}\|\Delta_{12} i\|_{\overline{s}_{h}+\mathtt{a}+\tau_{1}q+\tau_{1}+5}^{q,\kappa}.
\end{align}
Finally, putting together \eqref{estimate series Delta12 beta}, \eqref{telscopic diff betam}, \eqref{final dlt12 gm shb+1} and \eqref{final sum dlt12 gm shb+1}, we obtain
$$\|\Delta_{12}\beta\|_{\overline{s}_{h}+\mathtt{a}}^{q,\kappa}\lesssim\varepsilon\kappa^{-1}\|\Delta_{12}i\|_{q,\overline{s}_{h}+\mathtt{a}+\tau_{1}q+\tau_{1}+5}.$$
Thus, combining this estimate with  \eqref{link diff beta hat and diff beta}, \eqref{maj betamrk} and Sobolev embeddings  we find
\begin{align*}
\|\Delta_{12}\widehat{\beta}\|_{\overline{s}_{h}+\mathtt{a}}^{q,\kappa}&\lesssim\|\Delta_{12}\beta\|_{\overline{s}_{h}+\mathtt{a}}^{q,\kappa}\\
&\lesssim\varepsilon\kappa^{-1}\|\Delta_{12}i\|_{q,\overline{s}_{h}+\mathtt{a}+\tau_{1}q+\tau_{1}+5}
\end{align*}
The proof of Proposition \ref{QP-change} is now achieved.
\end{proof}

\subsection{First conjugation of the linearized operator}
The main goal of this section  is to explore the conjugation of the operator $ \mathcal{L}_{r,\lambda}$ introduced in Lemma \ref{lemma-reste} using the symplectic change of coordinates $\mathscr{B}$  constructed in Proposition \ref{QP-change}. Before stating our main result, we need to prove the following technical lemma.
\begin{lemma}\label{Lem-T-sing}
Let $\widehat\beta$ the function constructed in Proposition \ref{QP-change},  then we have the decomposition
\begin{align*}
\frac{\left|\sin\left(\frac{\theta-\eta}{2}\right)\right|^\alpha}{\left|\sin\left(\frac{\theta-\eta+\widehat{\beta}(\varphi,\theta)-\widehat{\beta}(\varphi,\eta)}{2}\right)\right|^\alpha}
&=\frac12\left(\big(1+ \partial_\theta\widehat\beta(\varphi,\theta)\big)^{-\alpha}+\big(1+ \partial_\eta\widehat\beta(\varphi,\eta)\big)^{-\alpha}\right)\\
&+\sin^2\left(\frac{\theta-\eta}{2}\right)\mathscr{K}_{r,\lambda}(\varphi,\theta,\eta),
\end{align*}
with the following estimate
$$
\forall \, s\in [s_0,S],\quad  \sup_{\eta\in\T}\|\mathscr{K}_{r,*}(\cdot,\centerdot, \centerdot+\eta)\|_{s}^{q,\kappa}\lesssim \varepsilon \kappa^{-1}\left(1+\| \mathfrak{I}_{0}\|_{s+\sigma_1+4}^{q,\kappa}\right),
$$
where the notation $*$ stands for $\lambda$,  $\cdot$  for $\varphi$ and $\centerdot$ denotes  $\theta$. The number $\sigma_1$ is given in \eqref{Conv-T2N}.

\end{lemma}
\begin{proof}
We shall start with the  the following expression  which follows from standard  trigonometric identities
\begin{align}\label{Trigo1}
\nonumber f(\varphi,\theta,\eta)&\triangleq \frac{\sin\left(\frac{\theta-\eta+\widehat{\beta}(\varphi,\theta)-\widehat{\beta}(\varphi,\eta)}{2}\right)}{\sin\left(\frac{\theta-\eta}{2}\right)}\\
&=\cos\left(\frac{\widehat{\beta}(\varphi,\theta)-\widehat{\beta}(\varphi,\eta)}{2}\right)+\frac{\sin\left(\frac{\widehat{\beta}(\varphi,\theta)-\widehat{\beta}(\varphi,\eta)}{2}\right)}{\sin\left(\frac{\theta-\eta}{2}\right)}\cos\left(\frac{\theta-\eta}{2}\right).
\end{align}
Applying the mean value theorem combined with the periodicity of $\widehat\beta$  we infer
\begin{align*}
\left| \frac{\sin\left(\frac{\widehat{\beta}(\varphi,\theta)-\widehat{\beta}(\varphi,\eta)}{2}\right)}{\sin\left(\frac{\theta-\eta}{2}\right)}\right|&\leqslant \frac{\left|{\widehat{\beta}(\varphi,\theta)-\widehat{\beta}(\varphi,\eta)}\right|}{2\left|\sin\left(\frac{\theta-\eta}{2}\right)\right|}\\
&\leqslant \|\widehat{\beta}\|_{\textnormal{Lip}(\T^{d+1})}.
\end{align*}
Putting together   \eqref{est-beta-r}  and  \eqref{small-C2} combined with  Sobolev embeddings  yield
\begin{align}\label{betasmal}
\nonumber \|\widehat{\beta}\|_{ s_0+1}^{q,\kappa}\leqslant& C\,\varepsilon \kappa^{-1}\left(1+\| \mathfrak{I}_{0}\|_{s_0+1+\sigma_1}^{q,\kappa}\right)\\
 \leqslant& C\varepsilon\kappa^{-1}
 \leqslant C\varepsilon_0,
\end{align}
which implies in turn from Sobolev embeddings that  $\|\widehat{\beta}\|_{\textnormal{Lip}(\T^{d+1})}$ is small enough.  Consequently we deduce from \eqref{Trigo1}  that $f$ is strictly positive. Next, we shall  use the splitting   $f=f_1+f_2 f_3$ with 
\begin{align*}
\nonumber f_1(\varphi,\theta,\eta)=\cos\left(\frac{\widehat{\beta}(\varphi,\theta)-\widehat{\beta}(\varphi,\eta)}{2}\right),&\quad 
f_2(\varphi,\theta,\eta)=\frac{\sin\left(\frac{\widehat{\beta}(\varphi,\theta)-\widehat{\beta}(\varphi,\eta)}{2}\right)}{\frac{\widehat{\beta}(\varphi,\theta)-\widehat{\beta}(\varphi,\eta)}{2}},\\
f_3(\varphi,\theta,\eta)&=\frac{\widehat{\beta}(\varphi,\theta)-\widehat{\beta}(\varphi,\eta)}{2\tan\left(\frac{\theta-\eta}{2}\right)}.
\end{align*}
 Then using Lemma \ref{Compos-lemm} and \eqref{est-beta-r} we obtain for any $s\in[s_0,S],$
\begin{align*}
\sup_{\eta\in\T}\|f_1(\cdot,\centerdot,\centerdot+\eta)-1\|_{s}^{q,\kappa}&\lesssim \|\widehat{\beta}\|_{s}^{q,\kappa}\\
 &\lesssim \varepsilon \kappa^{-1}\left(1+\| \mathfrak{I}_{0}\|_{s+\sigma_1}^{q,\kappa}\right).
\end{align*}
Since $x\mapsto \frac{\sin x}{x}$ is $\mathscr{C}^\infty$ then  we get in a similar way to  $f_1$
\begin{align*}
\sup_{\eta\in\T}\|f_2(\cdot,\centerdot,\centerdot+\eta)-1\|_{s}^{q,\kappa}
 \lesssim\, \varepsilon \kappa^{-1}\left(1+\| \mathfrak{I}_{0}\|_{s+\sigma_1}^{q,\kappa}\right).
\end{align*}
As  to the term $f_3$ we simply apply  Lemma \ref{lem-Reg1}-$($i$)$ in order to get in view of Proposition \ref{QP-change}-$($i$)$ 
\begin{align*}
\sup_{\eta\in\T}\|f_3(\cdot,\centerdot,\centerdot+\eta)\|_{s}^{q,\kappa}\lesssim&\, \|\widehat{\beta}\|_{s+1 }^{q,\kappa}\\
 \lesssim&\, \varepsilon \kappa^{-1}\left(1+\| \mathfrak{I}_{0}\|_{s+\sigma_1+1}^{q,\kappa}\right).
\end{align*}
Using the law products stated in Lemma \ref{Law-prodX1} combined with \eqref{small-C2} we deduce that
 \begin{align}\label{gy1}
\nonumber \sup_{\eta\in\T}\|f(\cdot,\centerdot,\centerdot+\eta)-1\|_{s}^{q,\kappa}\lesssim& \sup_{\eta\in\T}\Big[ \|f_1(\cdot,\centerdot,\centerdot+\eta)-1\|_{s}^{q,\kappa}+ \|f_3(\cdot,\centerdot,\centerdot+\eta)\|_{s}^{q,\kappa}\Big]\\
\nonumber&+ \sup_{\eta\in\T}\Big[\|f_2(\cdot,\centerdot,\centerdot+\eta)-1\|_{s_0}^{q,\kappa}\|f_3(\cdot,\centerdot,\centerdot+\eta)\|_{s}^{q,\kappa}\Big]\\
\nonumber&\quad+ \sup_{\eta\in\T}\Big[\|f_2(\cdot,\centerdot,\centerdot+\eta)-1\|_{s}^{q,\kappa}\|f_3\|_{s_0}^{q,\kappa}\Big]\\
&\qquad\lesssim \varepsilon \kappa^{-1}\left(1+\| \mathfrak{I}_{0}\|_{s+\sigma_1+1}^{q,\kappa}\right).
 \end{align} 
 In particular we get in view of the smallness condition \eqref{small-C2}
\begin{align*}
 \sup_{\eta\in\T}\|f(\cdot,\centerdot,\centerdot+\eta)-1\|_{s_0}^{q,\kappa}&\leqslant C  \varepsilon \kappa^{-1}\\
 &\leqslant C  \varepsilon_0.
 \end{align*}
 This implies that $f$ is  close to $1$ when $\varepsilon_0$ is small enough and therefore we may apply  Lemma \ref{Compos-lemm} to deduce that  $f^{-\alpha}$ satisfies similar  estimates, that is,
\begin{align}\label{gy2bis}
\sup_{\eta\in\T}\|f^{-\alpha}(\cdot,\centerdot,\centerdot+\eta)-1\|_{s}^{q,\kappa}&\lesssim \varepsilon \kappa^{-1}\left(1+\| \mathfrak{I}_{0}\|_{s+\sigma_1+1}^{q,\kappa}\right) \end{align} 
 Arguing similarly  we also get $(\varphi,\theta,\eta)\mapsto f^{-\alpha}(\varphi,\theta,\eta)\in {W}^{q,\infty}(\mathcal{O};H^s(\T^{d+2})$ with
 \begin{align}\label{gy2}
\|f^{-\alpha}-1\|_{s}^{q,\kappa}&\lesssim \varepsilon \kappa^{-1}\left(1+\| \mathfrak{I}_{0}\|_{s+\sigma_1+1}^{q,\kappa}\right).
 \end{align} 
 According to \eqref{Trigo1} one easily gets that the function $f$ is symmetric, that is,  $f(\varphi,\theta,\eta)=f(\varphi,\eta,\theta)$. Then using  Lemma \ref{Lemm-RegZ}-$($ii$)$ combined with \eqref{gy2} we find
\begin{align}\label{mathscrK}
{f^{-\alpha}(\varphi,\theta,\eta)}=\frac12\left({f^{-\alpha}(\varphi,\theta,\theta)}+{f^{-\alpha}(\varphi,\eta,\eta)}\right)+\sin^2\left(\frac{\theta-\eta}{2}\right)\mathscr{K}_{r,\lambda}(\varphi,\theta,\eta),
\end{align}
with the estimate
\begin{align}\label{mathscrK1}
\nonumber \sup_{\eta\in\T}\|\mathscr{K}_{r,\lambda}(\cdot,\centerdot, \centerdot+\eta)\|_{s}^{q,\kappa}\lesssim&\|\Delta_{\theta,\eta}(f^{-\alpha})\|_{s+\frac12+\epsilon}^{q,\kappa} \\
\nonumber &\lesssim\|f^{-\alpha}-1\|_{s+\frac52+\epsilon}^{q,\kappa} \\
&\quad \lesssim\,\varepsilon \kappa^{-1}\left(1+\| \mathfrak{I}_{0}\|_{s+\sigma_1+4}^{q,\kappa}\right)
\end{align}
On the other hand, we infer from  \eqref{Trigo1}
$$
f(\varphi,\theta,\theta)=1+ \partial_\theta\widehat\beta(\varphi,\theta)\quad\hbox{and}\quad  f(\varphi,\eta,\eta)=1+ \partial_\eta\widehat\beta(\varphi,\eta).
$$
Consequently, we find
$$
{f^{-\alpha}(\varphi,\theta,\eta)}=\frac12\left(\big(1+ \partial_\theta\widehat\beta(\varphi,\theta)\big)^{-\alpha}+\big(1+ \partial_\eta\widehat\beta(\varphi,\eta)\big)^{-\alpha}\right)+\sin^2\left(\frac{\theta-\eta}{2}\right)\mathscr{K}_{r,\lambda}(\varphi,\theta,\eta)
$$
and the proof is now achieved.
\end{proof}
Next, we plan to describe the action of the transformation $\mathscr{B}$ on the operator $ \mathcal{L}_{r,\lambda}$ introduced in Lemma \ref{lemma-reste} and derive some useful analytical estimates.
\begin{proposition}\label{prop-chang}
Under the same assumptions and notations  of \mbox{Proposition $\ref{QP-change}$},  one gets
for any $n\in\NN$ and for any $\lambda\in \mathcal{O}_{\infty,n}^{\kappa,\tau_{1}}(i_{0})$
\begin{align}\label{linXZ1}
\mathcal{L}_{r,\lambda}^0\triangleq \mathscr{B}^{-1}\mathcal{L}_{\varepsilon r,\lambda}\mathscr{B}=\omega\cdot\partial_\varphi+c(\lambda,i_0)\partial_\theta-\partial_\theta\Big(\mu_{r,\lambda}|\textnormal{D}|^{\alpha-1}+\textnormal{D}|^{\alpha-1} \mu_{r,\lambda}\Big)+\partial_\theta\mathcal{R}_{r,\lambda}+\mathtt{E}_{n}^{0},
\end{align}
where  $\mathcal{L}_{r,\lambda}$ was defined in \eqref{linXX1}-\eqref{W-eq}-\eqref{V-eq}, $\mathscr{B}$, $ c_{r,\lambda}$ and $\mathtt{E}_{n}^{0}$ satisfy the same estimates of \mbox{Proposition $\ref{QP-change}$}. In addition, 
 $\mu_{r}\in W^{q,\infty}_{\kappa}(\mathcal{O},H_{\textnormal{even}}^{s}) $ for any $ s\in[s_0,S]$   with the estimate
$$
\forall s\in[s_0,S],\quad \|\mu_{r,\lambda}-2^{-\alpha-1}C_\alpha\|_{s}^{q,\kappa}\lesssim\,\varepsilon \kappa^{-1}\left(1+\| \mathfrak{I}_{0}\|_{s+\sigma_1+1}^{q,\kappa}\right)
$$
and the remainder $\partial_\theta\mathcal{R}_{r,\lambda}$ is a reversible integral operator with,
\begin{align*}
\forall s\geqslant s_0, \forall\,\gamma\in\N\quad\hbox{with}\quad s+\gamma\leqslant S,\quad 
 \interleave\partial_\theta \mathcal{R}_{r,\lambda}\interleave_{-1,s,\gamma}^{q,\kappa}\lesssim\,\varepsilon \kappa^{-1}\left(1+\| \mathfrak{I}_{0}\|_{s+\sigma_1+7+\gamma}^{q,\kappa}\right).
\end{align*}
Moreover, there exists $\overline\sigma_1=\sigma(\tau_1,d,q)\geqslant\sigma_1$ such that for any $\gamma\in \N, s\geqslant 0$ with $ s+\gamma\leqslant \overline{s}_h$
$$
\|\Delta_{12}\mu_r\|_{\overline{s}_{h}}^{q,\kappa}+\interleave \Delta_{12}\partial_\theta\mathcal{R}_r\interleave_{-1,{s},\gamma}^{q,\kappa} \lesssim \,\varepsilon\kappa^{-1}\|\Delta_{12}i\|_{\overline{s}_{h}+\overline\sigma_1}^{q,\kappa}.
$$
\end{proposition}
\begin{proof}
Applying Proposition \ref{QP-change} and conjugating  by  $\mathscr{B}$ the   operator $\mathcal{L}_{\varepsilon r,\alpha}$ described in \eqref{linXX1}  we deduce that
 \begin{align}\label{lin1}
\mathscr{B}^{-1}\mathcal{L}_{\varepsilon r,\lambda}\mathscr{B}&=\omega\cdot\partial_\varphi+c(\lambda,i_0)\partial_\theta-\mathscr{B}^{-1} \partial_\theta\mathcal{T}_{ r}\mathscr{B}+\, \mathscr{B}^{-1}\partial_\theta\mathscr{R}_{\varepsilon r,\lambda}\mathscr{B},
\end{align}
where $\mathcal{T}_r$ denotes the self-adjoint operator
$$
\mathcal{T}_r= \mathscr{W}_{\varepsilon r,\alpha} |{\textnormal D}|^{\alpha-1}+|{\textnormal D}|^{\alpha-1}\mathscr{W}_{\varepsilon r,\alpha}
$$
and  $\mathscr{W}_{r,\alpha}$ is defined in \eqref{W-eq}. From the definition of the fractional Laplacian \eqref{fract1} we may write 
\begin{align}\label{KZ2}
\mathcal{T}_rh(\varphi,\theta)=
&\frac{ 1 }{2\pi }\bigintsss_{0}^{2\pi} \frac{\mathscr{W}_{\varepsilon r,\alpha}(\varphi,\theta)+\mathscr{W}_{\varepsilon r,\alpha}(\varphi,\eta)}{|\sin(\frac{\eta-\theta}{2})|^{\alpha}} h(\varphi,\eta)d\eta.
\end{align}
Applying  Lemma \ref{algeb1}$-($i$)$ we infer 
\begin{align}\label{KZ02}
\nonumber\mathscr{B}^{-1}\partial_\theta \mathcal{T}_r\mathscr{B}=&\partial_\theta {\bf{B}}^{-1} \mathcal{T}_r\mathscr{B}\\
&\triangleq \partial_\theta  \mathcal{T}_r^1.
\end{align}
Then using \eqref{mathscrB} and \eqref{KZ2}  together with  the change of variable $\eta\leadsto \eta+\widehat{\beta}(\varphi,\eta)$  we get
\begin{align}\label{KZ3}
\mathcal{T}_r^1h(\varphi,\theta)=
&\frac{ 1 }{2\pi }\bigintss_{0}^{2\pi} \frac{{\mathscr{W}}^1_{ r,\alpha}(\varphi,\theta,\eta)\,h(\varphi,\eta)}{\big|\sin(\frac{\theta-\eta+\widehat{\beta}(\varphi,\theta)-\widehat{\beta}(\varphi,\eta)}{2})\big|^{\alpha}} d\eta
\end{align}
where
\begin{align}\label{widehatW}
{\mathscr{W}}^1_{r,\alpha}(\varphi,\theta,\eta)&\triangleq \mathscr{W}^0_{r,\alpha}(\varphi,\theta)+\mathscr{W}^0_{r,\alpha}(\varphi,\eta)\\
\nonumber \mathscr{W}^0_{r,\alpha}(\varphi,\theta)&\triangleq {\mathscr{W}}_{\varepsilon r,\alpha}\big(\varphi,\theta+\widehat{\beta}(\varphi,\theta)\big)
\end{align}
and $\mathscr{W}_{r,\alpha}$ was introduced in \eqref{W-eq}. To estimate $ \mathscr{W}^0_{r,\alpha}$ we first note that
$$
\mathscr{W}_{0,\alpha}^0(\varphi,\theta)=C_\alpha 2^{-\alpha-1}=\mathscr{W}_{0,\alpha},
$$
which is independent of $\varphi$ and $\theta$. Then applying Lemma \ref{Compos1-lemm}  leads 
\begin{align*}
\|\mathscr{W}^0_{r,\alpha}-\mathscr{W}^0_{0,\alpha}\|_{s}^{q,\kappa}
\lesssim&\,\|\widehat\beta\|_{s_0+1}^{q,\kappa}\|\mathscr{W}_{\varepsilon r,\alpha}-\mathscr{W}_{0,\alpha}\|_{s+1}^{q,\kappa}+\|\widehat{\beta}\|_{s+1}^{q,\kappa}
\|\mathscr{W}_{ \varepsilon r,\alpha}-\mathscr{W}_{0,\alpha}\|_{s_0+1}^{q,\kappa}
\end{align*}
From Lemma \ref{lemma-reste}-$($i$)$  combined with Proposition \ref{QP-change}-$($ii$)$   we infer
\begin{align*}
\nonumber\|\mathscr{W}^0_{r,\alpha}-C_\alpha 2^{-\alpha-1}\|_{s}^{q,\kappa}
\lesssim& \,\varepsilon\,\|\widehat\beta\|_{s_0+1}^{q,\kappa}\|r\|_{s+2}^{q,\kappa}+\varepsilon\|\widehat{\beta}\|_{s+1}^{q,\kappa}
\|r\|_{s_0+2}^{q,\kappa}\\
\nonumber&\lesssim  \,\varepsilon^2\,\kappa^{-1}\left(1+\| \mathfrak{I}_{0}\|_{s_0+\sigma_1+1}^{q,\kappa}\right)\|r\|_{s+2}^{q,\kappa}+ \,\varepsilon^2\,\kappa^{-1}\left(1+\| \mathfrak{I}_{0}\|_{s+\sigma_1+1}^{q,\kappa}\right)
\|r\|_{s_0+2}^{q,\kappa}.
\end{align*}
Therefore, we deduce from Proposition \ref{lemma-GS0}-(i), Sobolev embeddings and \eqref{small-C2}
\begin{align}\label{Wr01}
\|\mathscr{W}^0_{r,\alpha}-C_\alpha 2^{-\alpha-1}\|_{s}^{q,\kappa}
&\lesssim \varepsilon \kappa^{-1}\left(1+\| \mathfrak{I}_{0}\|_{s+\sigma_1+1}^{q,\kappa}\right).
\end{align}
Putting together \eqref{Wr01} with \eqref{widehatW} yields for any $s\in[s_0,S],$
\begin{align}\label{Wr11}
\|\mathscr{W}^1_{ r,\alpha}-C_\alpha 2^{-\alpha}\|_{s}^{q,\kappa}
&\lesssim\, \varepsilon \kappa^{-1}\left(1+\| \mathfrak{I}_{0}\|_{s+\sigma_1+1}^{q,\kappa}\right).
\end{align}
Notice that the Sobolev norm $H^s$ in the preceding quantity  concerns all the variables $\varphi, \theta,\eta.$ By 
inserting the identity of Lemma \ref{Lem-T-sing} into \eqref{KZ3} yields
\begin{align}\label{KZ4}
\nonumber\mathcal{T}_r^1h(\varphi,\theta)=
&\frac{ 1 }{2\pi }\bigintss_{0}^{2\pi} \frac{{\mathscr{W}}^2_{r,\alpha}(\varphi,\theta,\eta)}{\big|\sin(\frac{\theta-\eta}{2})\big|^{\alpha}}h(\varphi,\eta) d\eta\\
&+\frac{ 1 }{2\pi }\bigintss_{0}^{2\pi} \frac{{\mathscr{W}}^1_{r,\alpha}(\varphi,\theta,\eta)\,\mathscr{K}_{r,\lambda}(\varphi,\theta,\eta)}{\big|\sin(\frac{\theta-\eta}{2})\big|^{\alpha-2}}h(\varphi,\eta) d\eta,
\end{align}
with
\begin{align}\label{Wrt2}
\mathscr{W}^2_{r,\alpha}(\varphi,\theta,\eta)=\frac12\mathscr{W}^1_{r,\alpha}(\varphi,\theta,\eta)\left(\big(1+ \partial_\theta\widehat\beta(\varphi,\theta)\big)^{-\alpha}+\big(1+ \partial_\eta\widehat\beta(\varphi,\eta)\big)^{-\alpha}\right).
\end{align}
Observe that the function $\mathscr{W}^2_{r,\alpha}$ is smooth and symmetric, $\mathscr{W}^2_{r,\alpha}(\varphi,\theta,\eta)=\mathscr{W}^2_{r,\alpha}(\varphi,\eta,\theta),$ then we apply Lemma \ref{Lemm-RegZ} leading to
\begin{align}\label{Wrt3}
\mathscr{W}^2_{r,\alpha}(\varphi,\theta,\eta)=\frac12\left( \mathscr{W}^2_{r,\alpha}(\varphi,\theta,\theta)+\mathscr{W}^2_{r,\alpha}(\varphi,\eta,\eta)\right)
+\sin^2\left(\frac{\theta-\eta}{2}\right)\mathscr{K}^1_{r,\lambda}(\varphi,\theta,\eta),
\end{align}
with the estimate
\begin{align}\label{KL1}
\sup_{\eta\in\T}\|(\partial_\theta^{p_1}\partial_\eta^{p_2}\mathscr{K}^1_{r,\lambda})(\cdot,\centerdot,\centerdot+\eta)\|_{s}^{q,\kappa}\lesssim \big\|\nabla_{\theta,\eta}\mathscr{W}^2_{r,\alpha}\big\|_{s+\frac32+p_1+p_2+\epsilon}^{q,\kappa}.
\end{align}
By writing \eqref{Wrt2} in the form
\begin{align*}
\mathscr{W}^2_{r,\alpha}(\varphi,\theta,\eta)&=\frac12\mathscr{W}^1_{r,\alpha}(\varphi,\theta,\eta)\left(\big(1+ \partial_\theta\widehat\beta(\varphi,\theta)\big)^{-\alpha}+\big(1+ \partial_\eta\widehat\beta(\varphi,\eta)\big)^{-\alpha}-2\right)\\
&+\big(\mathscr{W}^1_{r,\alpha}(\varphi,\theta,\eta)-C_\alpha 2^{-\alpha}\big)+C_\alpha 2^{-\alpha}
\end{align*}
and  applying {Lemmata  \ref{Law-prodX1}} and  \ref{Compos-lemm} combined with  \eqref{Wr11}  we find under the condition \eqref{small-C2},
\begin{align*}
\big\|\nabla_{\theta,\eta}\mathscr{W}^2_{ r,\alpha}\big\|_{s+p_1+p_2+\frac32+\epsilon}^{q,\kappa}&\lesssim \big\|\mathscr{W}^2_{r,\alpha}-C_\alpha 2^{-\alpha}\big\|_{s+\frac52+p_1+p_2+\epsilon}^{q,\kappa}
\\
&\lesssim \|\mathscr{W}^1_{ r,\alpha}\big\|_{s+p_1+p_2+\frac52+\epsilon}^{q,\kappa}\|\widehat\beta\|_{s_0+1}^{q,\kappa}+ \|\mathscr{W}^1_{r,\alpha}\big\|_{s_0}^{q,\kappa}\|\widehat\beta\|_{s+p_1+p_2+\frac72+\epsilon}^{q,\kappa}\\
&+\varepsilon\kappa^{-1}\left(1+\| \mathfrak{I}_{0}\|_{s+p_1+p_2+\frac52+\epsilon+\sigma_1}^{q,\kappa}\right).
\end{align*}
From \eqref{Wr11}, Proposition \ref{QP-change}-$($ii$)$ and  \eqref{small-C2}  we infer after straightforward manipulations
\begin{align}\label{WrtM4}
\big\|\nabla_{\theta,\eta}\mathscr{W}^2_{r,\alpha}\big\|_{s+p_1+p_2+\frac32+\epsilon}^{q,\kappa}&
\lesssim\,  \varepsilon \kappa^{-1} \big(1+\|\mathfrak{I}_{0}\|_{s+p_1+p_2+\frac72+\epsilon+\sigma}^{q,\kappa}\big).
\end{align}
Inserting \eqref{WrtM4} into \eqref{KL1} yields in view of Sobolev embeddings 
\begin{align}\label{KL2}
\nonumber\sup_{\eta\in\T}\|(\partial_\theta^{p_1}\partial_\eta^{p_2})\mathscr{K}^1_{r,\lambda}(\cdot,\centerdot,\centerdot+\eta)\|_{s}^{q,\kappa}&\lesssim \varepsilon \kappa^{-1}\left(1+\| \mathfrak{I}_{0}\|_{s+p_1+p_2+\frac72+\epsilon+\sigma_1}^{q,\kappa}\right) \\
&\lesssim \varepsilon\,\kappa^{-1} \left(1+\| \mathfrak{I}_{0}\|_{s+p_1+p_2+4+\sigma_1}^{q,\kappa}\right).
\end{align}
On the other hand, using \eqref{widehatW} and \eqref{Wrt2}   we get
\begin{align}\label{Wrt4}
\nonumber\mu_{r,\lambda}(\varphi,\theta)\triangleq&\frac12 \mathscr{W}^2_{r,\alpha}(\varphi,\theta,\theta)
\\
=&\frac{\mathscr{W}_{\varepsilon r,\alpha}\big(\varphi,\theta+\widehat\beta(\varphi,\theta)\big)}{\big(1+ \partial_\theta\widehat\beta(\varphi,\theta)\big)^{\alpha} }\cdot\end{align}
Plugging \eqref{Wrt3} and \eqref{Wrt4} into \eqref{KZ4} implies
\begin{align}\label{KZ6}
\nonumber \mathcal{T}_r^1h(\varphi,\theta)=
&\frac{ 1 }{2\pi }\bigintss_{0}^{2\pi} \frac{\mu_{r,\lambda}(\varphi,\theta)+\mu_{r,\lambda}(\varphi,\eta)}{\big|\sin(\frac{\theta-\eta}{2})\big|^{\alpha}}h(\varphi,\eta) d\eta\\
&+\frac{ 1 }{2\pi }\bigintss_{0}^{2\pi} \frac{\mathscr{K}^2_{r,\lambda}(\varphi,\theta,\eta)}{\big|\sin(\frac{\theta-\eta}{2})\big|^{\alpha-2}}h(\varphi,\eta) d\eta,
\end{align}
with
\begin{align}\label{Wrt5}
\mathscr{K}^2_{r,\lambda}(\varphi,\theta,\eta)\triangleq \mathscr{K}^1_{r,\lambda}(\varphi,\theta,\eta)+\mathscr{W}^1_{r,\alpha}(\varphi,\theta,\eta)\mathscr{K}_{r,\lambda}(\varphi,\theta,\eta)
\end{align}
and where $\mathscr{K}_{r,\lambda}$ is defined in \eqref{mathscrK} , $\mathscr{K}^1_{r,\lambda}$ in \eqref{Wrt3} and $\mathscr{W}_{r,\alpha}$ in \eqref{widehatW}.  To estimate $\mathscr{K}^2_{r,\lambda}$ we use Lemma \ref{Law-prodX1}
\begin{align}\label{WrtL1}
\nonumber \|\mathscr{K}^2_{r,\lambda}(\cdot,\centerdot,\centerdot+\eta)\|_{s}^{q,\kappa}\lesssim\|\mathscr{K}^1_{r,\lambda}(\cdot,\centerdot,\centerdot+\eta)\|_{s}^{q,\kappa}+&
\|\mathscr{W}^1_{r,\alpha}(\cdot,\centerdot,\centerdot+\eta)\|_{s_0}^{q,\kappa}\|\mathscr{K}_{r,\lambda}(\cdot,\centerdot,\centerdot+\eta)\|_{s}^{q,\kappa}\\
+&\|\mathscr{W}^1_{r,\alpha}(\cdot,\centerdot,\centerdot+\eta)\|_{s}^{q,\kappa}\|\mathscr{K}_{r,\lambda}(\cdot,\centerdot,\centerdot+\eta)\|_{s_0}^{q,\kappa}.
\end{align}
Then  \eqref{WrtL1}, \eqref{KL2}, \eqref{mathscrK1}, \eqref{Wr11} combined with \eqref{small-C2} imply
\begin{align}\label{WrtX}
\sup_{\eta\in\T}\|\mathscr{K}^2_{r,\lambda}(\cdot,\centerdot, \centerdot+\eta)\|_{s}^{q,\kappa}\lesssim&\,\,\varepsilon\kappa^{-1}  \left(1+\| \mathfrak{I}_{0}\|_{s+4+\sigma_1}^{q,\kappa}\right).
\end{align}
Similarly we also get forv any $p_1,o_2\in\NN$
\begin{align}\label{WrtXbis}
\sup_{\eta\in\T}\|(\partial_\theta^{p_1}\partial_\eta^{p_2}\mathscr{K}^2_{r,\lambda})(\cdot,\centerdot, \centerdot+\eta)\|_{s}^{q,\kappa}
\lesssim&\varepsilon\kappa^{-1}  \left(1+\| \mathfrak{I}_{0}\|_{s+p_1+p_2+4+\sigma_1}^{q,\kappa}\right).
\end{align}
Coming back to   \eqref{KZ6} and using the definition \eqref{fract1}  we deduce that 
\begin{align}\label{KZ005}
\nonumber \mathcal{T}_r^1h(\varphi,\theta)=
&\left(\mu_{r,\lambda}|\textnormal{D}|^{\alpha-1}+|\textnormal{D}|^{\alpha-1}\mu_{r,\lambda}\right)h(\varphi,\theta)+\frac{ 1 }{2\pi }\bigintss_{0}^{2\pi} \frac{\mathscr{K}^2_{r,\lambda}(\varphi,\theta,\eta)}{\big|\sin(\frac{\theta-\eta}{2})\big|^{\alpha-2}}h(\varphi,\eta) d\eta\\
\triangleq &\left(\mu_{r,\lambda}|\textnormal{D}|^{\alpha-1}+|\textnormal{D}|^{\alpha-1}\mu_{r,\lambda}\right)h(\varphi,\theta)+\mathcal{R}_r^{0,1}h(\varphi,\theta).
\end{align}
Next, we intend to  estimate the functions $\mu_{r,\lambda}$ defined in \eqref{Wrt4}.  Using  \eqref{W-eq} we get
\begin{align*}
\mu_{0,\lambda}(\varphi,\theta)&= \frac{\mathscr{W}_{0,\alpha}(\varphi,\theta)}{\big(1+ \partial_\theta\widehat\beta(\varphi,\theta)\big)^{\alpha}}\\
&=C_\alpha 2^{-\alpha-1}\big(1+ \partial_\theta\widehat\beta(\varphi,\theta)\big)^{-\alpha}.
\end{align*}
 In addition, we obtain from \eqref{Wrt4} and \eqref{widehatW}
\begin{align*}
\nonumber\mu_{r,\lambda}(\varphi,\theta)-\mathscr{W}_{0,\alpha}=&\frac{\mathscr{W}_{r,\alpha}^0(\varphi,\theta)-\mathscr{W}_{0,\alpha}^0}{\big(1+ \partial_\theta\widehat\beta(\varphi,\theta)\big)^{\alpha} }+\mathscr{W}_{0,\alpha}\left( \big(1+ \partial_\theta\widehat\beta(\varphi,\theta)\big)^{-\alpha}-1\right)\\
\triangleq &\mathcal{I}_1+\mathscr{W}_{0,\alpha}\mathcal{I}_2\cdot
\end{align*}
Now, applying {Lemmata  \ref{Law-prodX1}} and  \ref{Compos-lemm} yields  under the smallness condition \eqref{small-C2}, 
\begin{align*}
\|\mathcal{I}_1\|_{s}^{q,\kappa}&\lesssim \big(1+\|\partial_\theta\widehat\beta\|_{s_0}^{q,\kappa}\big)\,\big\|\mathscr{W}_{r,\alpha}^0-{\mathscr{W}}_{0,\alpha}^0\big\|_{s}^{q,\kappa}\\
&+\big(1+\|\partial_\theta\widehat\beta\|_{s}^{q,\kappa}\big)\,\big\|\mathscr{W}_{r,\alpha}^0-\mathscr{W}_{0,\alpha}^0\big\|_{s_0}^{q,\kappa}.
\end{align*}
Therefore we deduce from \eqref{Wr01}, \eqref{small-C2} and  Proposition \ref{QP-change}-$($ii$)$ 
\begin{align}\label{MW1}
\|\mathcal{I}_1\|_{s}^{q,\kappa}
&\lesssim\varepsilon \kappa^{-1}\left(1+\| \mathfrak{I}_{0}\|_{s+\sigma_1+1}^{q,\kappa}\right).
\end{align}
As to the term $\mathscr{W}_{0,\alpha}\mathcal{I}_2$ we shall combine   {Lemmata  \ref{Law-prodX1}}, \ref{Compos-lemm} with   Proposition \ref{QP-change}-$($ii$)$
\begin{align*}
\|\mathscr{W}_{0,\alpha}\mathcal{I}_2\|_{s}^{q,\kappa}&\lesssim \|\partial_\theta\widehat\beta\|_{s}^{q,\kappa}\\
&\lesssim  \varepsilon \kappa^{-1}\left(1+\| \mathfrak{I}_{0}\|_{s+\sigma_1+1}^{q,\kappa}\right)
.\end{align*}
It follows from the preceding estimates that
\begin{align*}
\|\mu_{r,\lambda}-\mathscr{W}_{0,\alpha}\|_{s}^{q,\kappa}
&\lesssim \varepsilon \kappa^{-1}\left(1+\| \mathfrak{I}_{0}\|_{s+\sigma_1+1}^{q,\kappa}\right).
\end{align*}
Putting together \eqref{lin1}, \eqref{KZ02} and \eqref{KZ005} we get
\begin{align*}
\mathscr{B}^{-1}\mathcal{L}_{\varepsilon r,\lambda}\mathscr{B}&=\omega\cdot\partial_\varphi+c_{r,\lambda}\partial_\theta-\partial_\theta\left(\mu_{r,\lambda}|\textnormal{D}|^{\alpha-1}+|\textnormal{D}|^{\alpha-1}\mu_{r,\lambda}\right)+\partial_\theta\mathcal{R}_r,
\end{align*}
with
\begin{align}\label{remaindX1}
\nonumber \partial_\theta\mathcal{R}_r\triangleq &-\partial_\theta\mathcal{R}_r^{0,1}+\, \mathscr{B}^{-1}\partial_\theta\mathscr{R}_{\varepsilon r,\lambda}\mathscr{B}\\
=&\partial_\theta(-\mathcal{R}_r^{0,1}+{\bf{B}}^{-1}\mathscr{R}_{\varepsilon r,\lambda}\mathscr{B}\big).
\end{align}
We shall start with estimating  $\partial_\theta \mathcal{R}_r^{0,1}$. Coming back to the definition \eqref{KZ005}  and differentiating in the variable  $\theta$ we get after straightforward computations 

\begin{align}\label{RrM1}
\partial_\theta \mathcal{R}_r^{0,1}h(\varphi,\theta)=\frac{ 1 }{2\pi }\bigintss_{0}^{2\pi} \frac{\sin(\frac{\eta-\theta}{2})\mathscr{K}^3_{r,\lambda}(\varphi,\theta,\eta)}{\big|\sin(\frac{\theta-\eta}{2})\big|^{\alpha}}h(\varphi,\eta) d\eta,
\end{align}
with
\begin{align}\label{K3bis0}
\mathscr{K}^3_{r,\lambda}(\varphi,\theta,\eta)=\frac{\alpha-1}{2}\cos\big((\eta-\theta)/2\big)\mathscr{K}^2_{r,\lambda}(\varphi,\theta,\eta)+\sin\big((\eta-\theta)/2\big)\partial_\theta\mathscr{K}^2_{r,\lambda}(\varphi,\theta,\eta).
\end{align}
Making the change of variable $\eta\leadsto \eta+\theta$ allows to get
\begin{align}\label{RrMU}
\partial_\theta \mathcal{R}_r^{0,1}h(\varphi,\theta)=\frac{ 1 }{2\pi }\bigintss_{0}^{2\pi} \frac{\sin(\frac{\eta}{2})\mathscr{K}^3_{r,\lambda}(\varphi,\theta,\theta+\eta)}{\big|\sin(\frac{\eta}{2})\big|^{\alpha}}h(\varphi,\theta+\eta) d\eta.
\end{align}
From Lemma  \ref{Law-prodX1} combined with \eqref{WrtX} we infer
\begin{align}\label{WM1}
\sup_{\eta\in\T}\|(\partial_\theta^{p_1}\partial_\eta^{p_2}\mathscr{K}^3_{r,\lambda})(\cdot,\centerdot,\centerdot+\eta)\|_{s}^{q,\kappa}\lesssim\,\varepsilon\kappa^{-1}  \left(1+\| \mathfrak{I}_{0}\|_{s+5+\sigma_1}^{q,\kappa}\right).
\end{align}
Next, we move to the estimate of the symbol associated to $ \partial_\theta \mathcal{R}_r^{0,1}$ whose expression can be derived  from \eqref{symb-kern} and takes the form
\begin{align*}
\sigma_{\partial_\theta \mathcal{R}_r^{0,1}}(\varphi,\theta,\xi)&=\frac{ 1 }{2\pi }\bigintsss_{0}^{2\pi} \frac{\sin(\frac{\eta}{2})\mathscr{K}^3_{r,\lambda}(\varphi,\theta,\theta+\eta)}{\big|\sin(\frac{\eta}{2})\big|^{\alpha}} e^{\ii \eta\xi}d\eta.
\end{align*}
Then integration by parts combined with \eqref{It-action} yield
\begin{align*}
\xi^{1+\gamma}\Delta_\xi^\gamma\sigma_{\partial_\theta \mathcal{R}_r^{0,1}}(\varphi,\theta,\xi)&=\frac{ i^{1+\gamma} }{2\pi }\bigintsss_{0}^{2\pi} e^{\ii \eta\xi}\partial_\eta^{1+\gamma} \left[\frac{ \sin(\frac{\eta}{2})(e^{\ii \eta}-1)^\gamma}{|\sin(\frac{\eta}{2})|^{\alpha}}\mathscr{K}^3_{r,\lambda}(\varphi,\theta,\theta+\eta)\right]\, d\eta.
\end{align*}
Therefore 
\begin{align*}
\langle \xi\rangle^{1+\gamma}\big\|\Delta_\xi^\gamma\sigma_{\partial_\theta \mathcal{R}_r^{0,1}}(\cdot,\centerdot,\xi)\big\|_{H^s}&\lesssim\sum_{0\leqslant \gamma^\prime\leqslant 1+\gamma}\bigintsss_{0}^{2\pi} \frac{\| \partial_\eta^{\gamma^\prime}\mathscr{K}^3_{r,\lambda}(\cdot,\centerdot,\centerdot+\eta)\|_{H^s}}{|\sin({\eta}/{2})|^{\alpha}}\, d\eta.
\end{align*}
Similarly we get for any $|\beta|\leqslant q$ and $\alpha\in(0,\overline\alpha)$, using Leibniz formula and Sobolev embeddings,
\begin{align*}
\kappa^{|\beta|}\langle \xi\rangle^{1+\gamma}\big\|\partial_\lambda^\beta\Delta_\xi^\gamma\sigma_{\partial_\theta \mathcal{R}_r^{0,1}}(\cdot,\centerdot,\xi)\big\|_{H^{s-|\beta|}}&\lesssim\kappa^{|\beta|}\sum_{0\leqslant \gamma^\prime\leqslant 1+\gamma\atop 
0\leqslant\beta^\prime\leqslant\beta}\bigintsss_{0}^{2\pi} \frac{\| \partial_\lambda^{\beta^\prime}\partial_\eta^{\gamma^\prime}\mathscr{K}^3_{r,\lambda}(\cdot,\centerdot,\centerdot+\eta)\|_{H^{s-|\beta^\prime|}}}{|\sin({\eta}/{2})|^{\frac12}}\, d\eta\\
&\lesssim\sum_{0\leqslant \gamma^\prime\leqslant 1+\gamma}\bigintsss_{0}^{2\pi} \frac{\| \partial_\eta^{\gamma^\prime}\mathscr{K}^3_{r,\lambda}(\cdot,\centerdot,\centerdot+\eta)\|_{s}^{q,\kappa}}{|\sin({\eta}/{2})|^{\frac12}}\, d\eta.
\end{align*}
Plugging  \eqref{WM1} into the preceding estimate and using \eqref{Def-pseud-w} allow to get
\begin{align}\label{oud-1}
\interleave\partial_\theta \mathcal{R}_r^{0,1}\interleave_{-1,s,\gamma}^{q,\kappa} &\lesssim \,\varepsilon \kappa^{-1}\left(1+\| \mathfrak{I}_{0}\|_{s+\gamma+\sigma_1+6}^{q,\kappa}\right).
\end{align}
Let us now move to the estimate of $\mathscr{B}_r^{-1}\partial_\theta\mathscr{R}_{r,\alpha}\mathscr{B}_r$. Then using \eqref{K2} and Lemma \ref{algeb1}-(i) we find
\begin{align}\label{KZL02}
\nonumber\mathscr{B}^{-1}\partial_\theta\mathscr{R}_{\varepsilon r,\lambda}\mathscr{B}=&\partial_\theta {\bf{B}}^{-1} \mathscr{R}_{\varepsilon r,\lambda}\mathscr{B}\\
\triangleq &\partial_\theta  \mathcal{S}_r,
\end{align}
with
\begin{align}\label{KZD3}
\mathcal{S}_rh(\varphi,\theta)=
&\frac{ 1 }{2\pi }\bigintss_{0}^{2\pi} \frac{\mathscr{A}_{r,\alpha}\big(\varphi,\theta+\widehat{\beta}(\varphi,\theta),\eta+\widehat{\beta}(\varphi,\eta)\big)}{\big|\sin(\frac{\theta-\eta+\widehat{\beta}(\varphi,\theta)-\widehat{\beta}(\varphi,\eta)}{2})\big|^{\alpha-2}} h(\varphi,\eta)d\eta.
\end{align}
Applying Lemma \ref{Lem-T-sing}, which  remains true if we replace $\alpha$ with $\alpha-2$,  we get
$$
{\left|\sin\left(\frac{\theta-\eta+\widehat{\beta}(\varphi,\theta)-\widehat{\beta}(\varphi,\eta)}{2}\right)\right|^{2-\alpha}}\triangleq \frac{F_1(\varphi,\theta,\eta)}{\big|\sin(\frac{\theta-\eta}{2})\big|^{\alpha-2}}\cdot
$$
The estimate of $F_1$ is similar to \eqref{Wr11} and one gets
\begin{align}\label{WrS11}
 \|F_1\|_{s}^{q,\kappa}  &\lesssim\, \varepsilon \kappa^{-1}\left(1+\| \mathfrak{I}_{0}\|_{s+\sigma_1+1}^{q,\kappa}\right).
\end{align}
The $H^s$ norm  of $F_1$ concerns  all the variables. Let us introduce  
$$
F_2(\varphi,\theta,\eta)=F_1(\varphi,\theta,\eta) \mathscr{A}_{r,\alpha}\big(\varphi,\theta+\widehat{\beta}(\varphi,\theta),\eta+\widehat{\beta}(\varphi,\eta)\big).
$$
Then similarly to the estimates \eqref{WrtL1} and \eqref{WrtXbis} we deduce that
\begin{align}\label{WrtPX}
 \sup_{\eta\in\T}\|F_2(\cdot,\centerdot, \centerdot+\eta)\|_{s}^{q,\kappa} &\lesssim\, \varepsilon \kappa^{-1}\left(1+\| \mathfrak{I}_{0}\|_{s+\sigma_1+4}^{q,\kappa}\right)
\end{align}
and
\begin{align}\label{WrtXbisMM}
\sup_{\eta\in\T}\|(\partial_\theta^{p_1}\partial_\eta^{p_2}F_2)(\cdot,\centerdot, \centerdot+\eta)\|_{s}^{q,\kappa}
&\lesssim\, \varepsilon \kappa^{-1}\left(1+\| \mathfrak{I}_{0}\|_{s+\sigma_1+p_1+p_2+4}^{q,\kappa}\right).
\end{align}
In addition, \eqref{KZD3} becomes
\begin{align}\label{KZDD3}
\mathcal{S}_rh(\varphi,\theta)=
&\frac{ 1 }{2\pi }\bigintss_{0}^{2\pi} \frac{F_2(\varphi,\theta+\eta)}{\big|\sin(\frac{\eta}{2})\big|^{\alpha-2}} h(\varphi,\theta+\eta)d\eta.
\end{align}
Now, we proceed as  for  \eqref{RrMU} by  combining \eqref{WrtPX}, \eqref{WrtXbisMM} and \eqref{KZDD3} in order to get
\begin{align*}
\interleave\mathcal{S}_r\interleave_{-2,s,\gamma}^{q,\kappa}&\lesssim\, \varepsilon \kappa^{-1}\left(1+\| \mathfrak{I}_{0}\|_{s+\sigma_1+\gamma+6}^{q,\kappa}\right).
\end{align*}
On the other hand one gets from \eqref{kernel-phi-lambda} the following identity
$$
\sigma_{\partial_\theta\mathcal{S}_r}=\partial_\theta\sigma_{\mathcal{S}_r}+i\xi \sigma_{\mathcal{S}_r}.
$$
Consequently, we find
\begin{align*}
\interleave\partial_\theta\mathcal{S}_r\interleave_{-1,s,\gamma}^{q,\kappa}&\lesssim\interleave\mathcal{S}_r\interleave_{-2,s+1,\gamma}^{q,\kappa}\\
&\lesssim\, \varepsilon \kappa^{-1}\left(1+\| \mathfrak{I}_{0}\|_{s+\sigma_1+\gamma+7}^{q,\kappa}\right),
\end{align*}
which gives in view of \eqref{KZL02}
\begin{align*}
\interleave\mathscr{B}^{-1}\partial_\theta\mathscr{R}_{\varepsilon r,\lambda}\mathscr{B}\interleave_{-1,s,\gamma}^{q,\kappa}\lesssim\,  \varepsilon \kappa^{-1}\left(1+\| \mathfrak{I}_{0}\|_{s+\sigma_1+\gamma+7}^{q,\kappa}\right).
\end{align*}
 Next we shall move to the estimate of the differences. Let us start with $\Delta_{12}\mu_r$ where $\mu_r$ is defined in \eqref{Wrt4}. The computations are slightly long but they are classical and resemble  to the preceding estimates. Then  based on \eqref{difference beta}, Lemma \ref{Compos1-lemm} combined with  \eqref{small-C2} we find
\begin{equation*}
\|\Delta_{12}\mu_r\|_{\overline{s}_{h}}^{q,\kappa}\lesssim\varepsilon\kappa^{-1}\|\Delta_{12}i\|_{\overline{s}_{h}+\sigma}^{q,\kappa},
\end{equation*}
for some $\sigma$ depending only on $\tau_1, d$ and $q.$
Similar arguments combined with the estimates used for getting  \eqref{WM1}  allow to get
\begin{align*}
\forall\, s\geqslant 0, \gamma\in\N\quad  \hbox{with}\,\quad  s+\gamma\leqslant \overline{s}_h,\quad \sup_{\eta\in\T}\|\partial_\eta^\gamma\Delta_{12}\mathscr{K}^3_{r,\lambda}(\cdot,\centerdot,\centerdot+\eta)\|_{{s}}^{q,\kappa}\lesssim\,\varepsilon\kappa^{-1}\|\Delta_{12}i\|_{\overline{s}_{h}+\sigma}^{q,\kappa}.
\end{align*}
Consequently, we find in view of \eqref{RrM1} and implementing similar arguments to \eqref{oud-1} 
\begin{align*}
\forall\, s\geqslant 0, \gamma\in\N\quad  \hbox{with}\,\quad  s+\gamma\leqslant \overline{s}_h,\quad\interleave \Delta_{12}\partial_\theta \mathcal{R}_r^{0,1}\interleave_{-1,{s},\gamma}^{q,\kappa} &\lesssim \,\varepsilon\kappa^{-1}\|\Delta_{12}i\|_{\overline{s}_{h}+\sigma}^{q,\kappa}.
\end{align*}
Coming back to \eqref{KZDD3} and implementing the same arguments we find
\begin{align*}
\forall\, s\geqslant 0, \gamma\in\N\quad  \hbox{with}\,\quad  s+\gamma\leqslant \overline{s}_h\quad\interleave \Delta_{12}\partial_\theta\mathcal{S}_r\interleave_{-1,{s},\gamma}^{q,\kappa} &\lesssim \,\varepsilon\kappa^{-1}\|\Delta_{12}i\|_{\overline{s}_{h}+\sigma}^{q,\kappa}.
\end{align*}
Finally, we get from \eqref{remaindX1} and the preceding estimates
\begin{equation}\label{puir1}
\forall\, s\geqslant 0, \gamma\in\N\quad  \hbox{with}\,\quad  s+\gamma\leqslant \overline{s}_h,\quad \interleave \Delta_{12}\partial_\theta\mathcal{R}_r\interleave_{-1,{s},\gamma}^{q,\kappa} \lesssim \,\varepsilon\kappa^{-1}\|\Delta_{12}i\|_{\overline{s}_{h}+\sigma}^{q,\kappa}
\end{equation}
This achieves the proof of Proposition \ref{prop-chang}.

\end{proof}
\subsection{Reduction of the nonlocal part} \label{Sec-Local-P} 
The main goal of this section is to reduce to constant coefficients the   leading  term in \eqref{linXZ1} with positive order. Notice that at this level  the transport part is with constant coefficient and we need to find a judicious conjugation  in order to reduce the nonlocal part to a  Fourier multiplier with the same order, of course without altering the transport part. As we shall see below, this will be done using the infinite dimensional hyperbolic flow $\Phi$ described by \eqref{eqn:omega0} through a specific choice of the density $\rho.$  We emphasize that this kind of reduction was used for water waves  in  \cite{BertiMontalto}. Our main result can be stated as follows. 
\begin{proposition}\label{prop-constant-coe}
There exists $\sigma_2=\sigma(\tau_1,d,q)\geqslant\overline\sigma_1$  such that under \eqref{Conv-T2} and \eqref{Conv-T2N}
 and by assuming 
\begin{align}
\label{small-CC2}\|\mathfrak{I}_{0}\|_{q,s_{h}+\sigma_2}\leqslant1\quad\textnormal{and}\quad N_{0}^{\mu_{2}}\varepsilon{\kappa^{-2}}\leqslant{\varepsilon}_0,
\end{align}
the following assertions hold true. 
 There exists  
 ${\Psi}:{{\mathcal{O}}}\to \mathscr{L}\big(H^{s}(\T^{d+1}),H^{s}(\T^{d+1})\big) $  a family of reversibility preserving invertible linear bounded  operators  and  such  that for any $n\in\NN$ and for any $\lambda$ in the set $ \mathcal{O}_{\infty,n}^{\kappa,\tau_{1}}(i_{0})$ defined in Proposition $\ref{QP-change}$ we have 
\begin{align}\label{linX1}
\Psi^{-1}\mathcal{L}_{\varepsilon r,\lambda}\Psi\triangleq \mathcal{L}_{r,\lambda}^1\triangleq\omega\cdot\partial_\varphi+{c}(\lambda,i_0)\partial_\theta- \textnormal{m}(\lambda,i_0)\partial_\theta|\textnormal{D}|^{\alpha-1}+\mathcal{R}^1_{r,\lambda}+\mathtt{E}_{n}^{1},
\end{align}
with the following properties. \begin{enumerate}
\item  We have $\displaystyle{\textnormal{m}(\lambda,i_0)=\tfrac{2}{(2\pi)^{d+1}}\int_{\T^{d+1}} \mu_{r,\lambda}(\varphi,\theta)d\varphi d\theta}$ and  
$$
\|\textnormal{m}(\cdot,i_0)-2^{-\alpha}C_\alpha\|^{q,\kappa}+\|c(\cdot,i_0)-V_{0,\alpha}\|^{q,\kappa}\lesssim \varepsilon\kappa^{-1}.
$$
Moreover
$$
\|\Delta_{12}\textnormal{m}(\cdot,i)\|^{q,\kappa}+\|\Delta_{12}{c}(\cdot,i)\|^{q,\kappa}\lesssim \varepsilon\kappa^{-1}\|\Delta_{12} i\|_{\overline{s}_h+\sigma_2}^{q,\kappa}.
$$
 
\item The operator $\mathcal{ R}_{r,\lambda}^1$ is reversible and for any  $s\in[s_0,S]$ and for   $\epsilon>0$ small enough 
\begin{align*}
\interleave{\mathcal R}_{r,\lambda}^1\interleave_{2\overline\alpha-1+\epsilon,s,0}^{q,\kappa}+{\interleave{\mathcal R}_{r,\lambda}^{1\,\star}\interleave_{2\overline\alpha-1+\epsilon,s,0}^{q,\kappa}}\lesssim \varepsilon \kappa^{-2}\left(1+\| \mathfrak{I}_{0}\|_{s+\sigma_2}^{q,\kappa}\right).
\end{align*}
In addition,
\begin{align*}
\|\Delta_{12}\mathcal{R}^1_{r,\lambda}h\|_{\overline{s}_h }^{q,\kappa}+\|\Delta_{12}\mathcal{R}^{1,\star}_{r,\lambda}h\|_{\overline{s}_h }^{q,\kappa}&\lesssim
\varepsilon \kappa^{-2}\|\Delta_{12} i\|_{q,\overline{s}_h+\sigma_2}^{q,\kappa} \| h\|_{\overline{s}_h+\sigma_2}^{q,\kappa}\end{align*}
and
\begin{align*}
\interleave \Delta_{12}\mathcal{R}^1_{r,\lambda}\interleave_{0,0,0 }^{q,\kappa}&\lesssim
\varepsilon \kappa^{-2}\|\Delta_{12} i\|_{q,\overline{s}_h+\sigma_2}^{q,\kappa} .
\end{align*}
\item The operators $\Psi^{\pm1}$ satisfy for any $s\geqslant s_0$ 
$$
\big\|\Psi^{\pm1}h\big\|_{s}^{q,\kappa}\lesssim \|h\|_{s}^{q,\kappa}+\varepsilon\kappa^{-2}\big(1+\| \mathfrak{I}_{0}\|_{s+\sigma_2}^{q,\kappa}\big)
\| h\|_{s_0}^{q,\kappa}
$$
and 
\begin{align*}
\big\|(\Psi^{\pm1}-\textnormal{Id})h\big\|_{s}^{q,\kappa}&+\big\|((\Psi^{\pm1})^{\star}-\textnormal{Id})h\big\|_{s}^{q,\kappa}\lesssim \varepsilon\kappa^{-2}
\|h\|_{s+1}^{q,\kappa}+\varepsilon\kappa^{-2}\| \mathfrak{I}_{0}\|_{s+\sigma_2}^{q,\kappa}
\| h\|_{s_0}^{q,\kappa}.
\end{align*}
\begin{align*}
\|\Delta_{12}\Psi^{\pm1} h\|_{\overline{s}_h}^{q,\kappa}+\|\Delta_{12}(\Psi^{\pm1})^{\star} h\|_{\overline{s}_h}^{q,\kappa}
&\lesssim \varepsilon{\kappa^{-2}}\|\Delta_{12}i\|_{\overline{s}_h+\sigma_2}^{q,\kappa} \| h\|_{\overline{s}_h+\sigma_2}^{q,\kappa}.
\end{align*}
\item The operator  $\mathtt{E}_{n}^{1}$ is reversible and satisfies the estimate
$$
\forall s\in[s_0,S],\quad \big\|\mathtt{E}_n^{1}h\big\|_{s_0}^{q,\kappa}\lesssim\varepsilon{\kappa^{-2}}\Big(N_{0}^{\mu_{2}}N_{n+1}^{-\mu_{2}}+  N_{n}^{s_0-s}\big(1+\| \mathfrak{I}_{0}\|_{s+\sigma_2}^{q,\kappa}\big)\Big)\| h\|_{s_0+3}^{q,\kappa}.$$

\end{enumerate}
\end{proposition}
 \begin{proof}
We have seen in Proposition \ref{prop-chang} that for any $\lambda\in \mathcal{O}_{\infty,n}^{\kappa,\tau_{1}}(i_{0})$
\begin{align}\label{etaM1}
\nonumber \mathcal{L}_{r,\lambda}^0=\mathscr{B}^{-1}\mathcal{L}_{\varepsilon r,\lambda}\mathscr{B}&=\omega\cdot\partial_\varphi+c(\lambda,i_0)\partial_\theta-\partial_\theta\mathcal{A}+\partial_\theta\mathcal{R}_{r,\lambda}+\mathtt{E}_n^0,\\
\hbox{with}\quad \mathcal{A}\triangleq &\mu_{r,\lambda}|\textnormal{D}|^{\alpha-1}+|\textnormal{D}|^{\alpha-1} \mu_{r,\lambda}.
\end{align}  
Next, we shall conjugate this operator by  the flow $\Phi(t)$ described in  Proposition \ref{flowmap00}  
  \begin{align}\label{lin2}
\mathcal{L}_{r,\lambda}^t&\triangleq\Phi^{-1}(t)\mathcal{L}_{r,\lambda}^0\Phi(t),\quad \hbox{with}\quad \Phi(t)=e^{t \mathbb{A}},\\
\nonumber&=\Psi^{-1}(t)\mathcal{L}_{\varepsilon r,\lambda}\Psi(t),\quad \hbox{with}\quad \Psi(t)=\mathscr{B}\Phi(t),
\end{align}
 where the generator  $\mathbb{A}$ takes the form
\begin{align}\label{self-ad1}
\mathbb{A}=\partial_\theta\mathcal{T}\quad\hbox{and}\quad \mathcal{T}=\rho|\textnormal{D}|^{\alpha-1}+|\textnormal{D}|^{\alpha-1} \rho,
\end{align}
with $\rho$ being  a smooth function to be chosen later. 
Then Taylor formula implies
\begin{align}\label{Lie01}
\mathcal{L}_{r,\lambda}^1=\mathcal{L}_{r,\lambda}^0+\bigintsss_0^1\,\Phi^{-1}(t)\big[\mathcal{L}_{r,\lambda}^0,\mathbb{A}\big]\Phi(t) dt.
\end{align}
Let us evaluate the commutator $\big[\mathcal{L}_{r,\lambda}^0,\mathbb{A}\big]$ using the precise structure of $\mathbb{A}$. By direct computations we may check  that
$$
\big[\omega\cdot\partial_\varphi,\mathbb{A}\big]=\partial_\theta\Big(\left(\omega\cdot\partial_\varphi\rho\right)|\textnormal{D}|^{\alpha-1}+|\textnormal{D}|^{\alpha-1} \left(\omega\cdot\partial_\varphi\rho\right)\Big)
$$
and
$$
\big[c(\lambda,i_0)\partial_\theta,\mathbb{A}\big]=\partial_\theta\Big(\big(c(\lambda,i_0)\partial_\theta\rho\big)|\textnormal{D}|^{\alpha-1}+|\textnormal{D}|^{\alpha-1} \big(c(\lambda,i_0)\partial_\theta\rho\big)\Big).
$$
Consequently,
\begin{align}\label{Cod-v1}
\big[\omega\cdot\partial_\varphi+c(\lambda,i_0)\partial_\theta,\mathbb{A}\big]=\partial_\theta\Big(\widehat\rho\, |\textnormal{D}|^{\alpha-1}+|\textnormal{D}|^{\alpha-1} \widehat\rho\Big)\quad\hbox{with}\quad \widehat\rho\triangleq \omega\cdot\partial_\varphi\rho+c(\lambda,i_0)\partial_\theta\rho.
\end{align}
Define 
\begin{align}\label{hatA}
 \widehat{\mathbb{A}}&=\partial_\theta\big(\widehat\rho\, |\textnormal{D}|^{\alpha-1}+|\textnormal{D}|^{\alpha-1} \widehat\rho\big).
\end{align} 
It follows that
\begin{align}\label{comP1}
\big[\mathcal{L}_{r,\lambda}^0,\mathbb{A}\big]=& \widehat{\mathbb{A}}+\Big[-\partial_\theta\mathcal{A}+\partial_\theta\mathcal{R}_{r,\lambda},\mathbb{A}\Big]+\big[\mathtt{E}_n^0,\mathbb{A}\big].
\end{align}
Inserting \eqref{comP1} into \eqref{Lie01} yields 
\begin{align}\label{LrY}
\nonumber\mathcal{L}_{r,\lambda}^1=\omega\cdot\partial_\varphi+c(\lambda,i_0)\partial_\theta-\partial_\theta\mathcal{A}+\partial_\theta \mathcal{R}_{r,\lambda}+&\int_0^1\Phi(-t)  \widehat{\mathbb{A}}\Phi(t)dt\\
\nonumber &+\int_0^1\Phi(-t)\Big[\partial_\theta\mathcal{R}_{r,\lambda}-\partial_\theta\mathcal{A},\mathbb{A}\Big]\Phi(t)dt\\
&\quad+\mathtt{E}_n^0+\int_0^1\Phi(-t)\big[\mathtt{E}_n^0,\mathbb{A}\big]\Phi(t)dt.
\end{align}
Then  using once again Taylor formula we find
$$
{\Phi(-t)} \widehat{\mathbb{A}}{\Phi(t)}=\widehat{\mathbb{A}}+\int_0^t\Phi(-t^\prime) \big[\widehat{\mathbb{A}},\mathbb{A}\big]{\Phi(t^\prime)}dt^\prime.
$$
Plugging this identity into \eqref{LrY} allows to get
\begin{align}\label{LrYP}
\mathcal{L}_{r,\lambda}^1&=\omega\cdot\partial_\varphi+c(\lambda,i_0)\partial_\theta+\big(\widehat{\mathbb{A}}-\partial_\theta\mathcal{A}\big)+{\mathcal{R}}^1_{r,\lambda}\\
\nonumber &\quad+\mathtt{E}_n^0+\int_0^1\Phi(-t)\big[\mathtt{E}_n^0,\mathbb{A}\big]\Phi(t)dt,
\end{align}
with
\begin{align}\label{hatR}
\nonumber{\mathcal R}^1_{r,\lambda}\triangleq&\partial_\theta\mathcal R_{r,\lambda}+\int_0^1\Phi(-t)\Big[\partial_\theta\mathcal{R}_{r,\lambda}-\partial_\theta\mathcal{A},\mathbb{A}\Big]\Phi(t)dt+\int_0^1\int_0^t{\Phi(-t^\prime)} \big[\widehat{\mathbb{A}},\mathbb{A}\big]{\Phi(t^\prime)}dt^\prime dt\\
&=\partial_\theta\mathcal R_{r,\lambda}+\int_0^1\Phi(-t)\Big[\partial_\theta\mathcal{R}_{r,\lambda}-\partial_\theta\mathcal{A},\mathbb{A}\Big]\Phi(t)dt+\int_0^1(1-t)\Phi(-t) \big[\widehat{\mathbb{A}},\mathbb{A}\big]{\Phi(t)} dt,
\end{align}
where we used in the last inequality an integration by parts.
Combining   \eqref{etaM1} and \eqref{hatA} we may write
\begin{align}\label{Lup768}
\widehat{\mathbb{A}}-\partial_\theta\mathcal{A}=\partial_\theta\big[(\widehat\rho-\mu_{r,\lambda})|\textnormal{D}|^{\alpha-1}+|\textnormal{D}|^{\alpha-1} (\widehat\rho-\mu_{r,\lambda})\big].
\end{align}
We recall that  $\mu_{r,\lambda}$ was defined in \eqref{Wrt4} and in what follows $\langle \mu_{r,\lambda}\rangle_{\varphi,\theta}$ stands for its average  in both variables. According to   Lemma \ref{L-Invert}, the following functions
\begin{align}\label{form-rep-X1}
 \rho\triangleq(\omega\cdot\partial_\varphi+c(\lambda,i_0)\partial_\theta)_{\textnormal{ext}}^{-1}\big(\mu_{r,\lambda}-\langle \mu_{r,\lambda}\rangle_{\varphi,\theta}\big)\quad\hbox{and}\quad  \rho_{n}\triangleq\Pi_{N_n}\rho
\end{align}
are well defined in the whole set $\mathcal{O}$ and moreover for any $\lambda\in \mathcal{O}_{\infty,n}^{\kappa,\tau_{1}}(i_{0})$
\begin{align}\label{Eq-etaLL}
\big(\omega\cdot\partial_\varphi +c(\lambda,i_0)\partial_\theta\big)\rho_{n}=\Pi_{N_n}\Big(\mu_{r,\lambda}-\langle \mu_{r,\lambda}\rangle_{\varphi,\theta}\Big).
\end{align}
Therefore, we deduce  from  the decomposition
\begin{align}\label{Density-split12}
\rho=\rho_{n}+\Pi_{N_n}^\perp\rho
\end{align}
and \eqref{Eq-etaLL} that for any $\lambda\in \mathcal{O}_{\infty,n}^{\kappa,\tau_{1}}(i_{0})$ 
\begin{align}\label{Eq-eta}
\big(\omega\cdot\partial_\varphi +c(\lambda,i_0)\partial_\theta\big)\rho&=\Pi_{N_n}\Big(\mu_{r,\lambda}-\langle \mu_{r,\lambda}\rangle_{\varphi,\theta}\Big)+\big(\omega\cdot\partial_\varphi +c(\lambda,i_0)\partial_\theta\big)\Pi_{N_n}^\perp\rho\\
\nonumber&\triangleq \Pi_{N_n}\Big(\mu_{r,\lambda}-\langle \mu_{r,\lambda}\rangle_{\varphi,\theta}\Big)+{\bf R}_n.
\end{align}
Putting together \eqref{Lup768}, \eqref{Cod-v1} and \eqref{Eq-eta}  yields  for  $\lambda\in \mathcal{O}_{\infty,n}^{\kappa,\tau_{1}}(i_{0})$
\begin{align}\label{Eq-eta-M}
\nonumber\widehat{\mathbb{A}}-\partial_\theta\mathcal{A}=-2\langle \mu_{r,\lambda}\rangle_{\varphi,\theta}\,\partial_\theta\,|\textnormal{D}|^{\alpha-1}&-\partial_\theta\big[(\Pi_{N_n}^\perp\mu_{r,\lambda})|\textnormal{D}|^{\alpha-1}+|\textnormal{D}|^{\alpha-1} (\Pi_{N_n}^\perp \mu_{r,\lambda})\big]\\
&+{\partial_\theta\big({\bf R}_n|\textnormal{D}|^{\alpha-1}+|\textnormal{D}|^{\alpha-1} {\bf R}_n\big)}.
\end{align}
Inserting \eqref{Eq-eta-M} into \eqref{LrYP} allows to get for any $\lambda\in \mathcal{O}_{\infty,n}^{\kappa,\tau_{1}}(i_{0})$
\begin{align}\label{LrYP0P}
\nonumber \mathcal{L}_{r,\lambda}^1&=\Psi^{-1}\mathcal{L}_{\varepsilon r,\lambda}\Psi,\quad\hbox{with}\quad \Psi\triangleq \Psi(1)=\mathscr{B}\Phi(1)\\
&=\omega\cdot\partial_\varphi+c(\lambda,i_0)\partial_\theta-\textnormal{m}(\lambda,i_0)\,\partial_\theta\,|\textnormal{D}|^{\alpha-1}+{\mathcal{R}}^1_{r,\lambda}+\mathtt{E}_n^1,
\end{align}
with $\textnormal{m}(\lambda,i_0)\triangleq2\langle \mu_{r,\lambda}\rangle_{\varphi,\theta}$ and
\begin{align}\label{En-11}
 \nonumber\mathtt{E}_n^1\triangleq \mathtt{E}_n^0+\int_0^1\Phi(-t)\big[\mathtt{E}_n^0,\mathbb{A}\big]\Phi(t)dt&-\partial_\theta\big[(\Pi_{N_n}^\perp\mu_{r,\lambda})|\textnormal{D}|^{\alpha-1}+|\textnormal{D}|^{\alpha-1} (\Pi_{N_n}^\perp \mu_{r,\lambda})\big]\\
&+\partial_\theta\big({\bf R}_n|\textnormal{D}|^{\alpha-1}+|\textnormal{D}|^{\alpha-1} {\bf R}_n\big).
\end{align}

{\bf{$($i$)$}} The estimate of $c(\lambda,i_0)$  was done before in Proposition \ref{QP-change}-$($i$)$. However the  estimate of $\mu_{r,\lambda}$ is given in  Proposition \ref{prop-chang}. Then the   estimate of the average  $\textnormal{m}_{r,\lambda}$ follows easily from  Sobolev embeddings and \eqref{small-C2}
\begin{align*}
\nonumber\big\|\textnormal{m}(\cdot,i_0)-2^{-\alpha} C_\alpha\big\|^{q,\kappa}=&2\big\|\langle {\mu}_{r,\lambda}-2^{-\alpha-1} C_\alpha\rangle_{\varphi,\theta} \big\|^{q,\kappa}
\\
\nonumber&\lesssim \big\|\mu_{r,\lambda}-2^{-\alpha-1} C_\alpha\big\|_{s_0}^{q,\kappa}\\
&\quad \lesssim  \,\varepsilon \kappa^{-1}\left(1+\| \mathfrak{I}_{0}\|_{s_0+\sigma_1+1}^{q,\kappa}\right)\\
&\quad \quad \lesssim  \,\varepsilon \kappa^{-1}.
\end{align*}
For the estimate of the difference $\textnormal{m}_{r,\lambda}$ we write in view of Proposition \ref{prop-chang} 
\begin{align*}
\nonumber\big\|\Delta_{12}\textnormal{m}(\cdot,i)\big\|^{q,\kappa}=&2\big\|\langle \Delta_{12}{\mu}_{r,\lambda}\rangle_{\varphi,\theta} \big\|^{q,\kappa}
\\
\lesssim& \,\varepsilon\kappa^{-1}\|\Delta_{12}i\|_{\overline{s}_{h}+\overline\sigma_1}^{q,\kappa}
.
\end{align*}
{\bf{$($ii$)$}} We shall start with estimating $\rho$  given in  \eqref{form-rep-X1}. Then applying  
Lemma \ref{L-Invert}   gives 
\begin{align}\label{Plus-11}
\|\rho\|_{s}^{q,\kappa}&\lesssim \kappa^{-1}\big(1+\|c(\cdot,i_0)\|^{q,\kappa}\big)\big\|\mu_{r,\lambda}-\langle \mu_{r,\lambda}\rangle_{\varphi,\theta}\big\|_{s+\tau_1(q+1)}^{q,\kappa}.
\end{align}
Using \eqref{est-r1} combined with \eqref{small-CC2} yield
\begin{align}\label{Tyk6}
 \|c(\cdot,i_0)\|^{q,\kappa}&\lesssim 1.
\end{align}
On the other hand,  we infer from Proposition \ref{prop-chang} 
$$
\big\|\mu_{r,\lambda}-\langle \mu_{r,\lambda}\rangle_{\varphi,\theta}\big\|_{s+\tau_1(q+1)}^{q,\kappa}\lesssim\,\varepsilon \kappa^{-1}\left(1+\| \mathfrak{I}_{0}\|_{s+\sigma_1+\tau_1(1+q)+1}^{q,\kappa}\right).
$$
Hence, combining this estimate with \eqref{Tyk6} and \eqref{Plus-11} we deduce that
\begin{align}\label{Est-eta}
\|\rho\|_{s}^{q,\kappa}&\lesssim \varepsilon \kappa^{-2}\left(1+\| \mathfrak{I}_{0}\|_{s+\sigma_1+\tau_1(1+q)+1}^{q,\kappa}\right).
\end{align}
{Implementing similar estimates used before to get \eqref{estimate delta12 gm},  by  making appeal   in particular to  the estimates of $\Delta_{12}c_r$ and $\Delta_{12}\mu_r$ detailed in Proposition \ref{prop-chang}, one gets after tedious and long computations  
\begin{align}\label{west-eta}
\|\Delta_{12}\rho\|_{\overline{s}_h}^{q,\kappa}&\lesssim \varepsilon\kappa^{-2}\|\Delta_{12} i\|_{\overline{s}_h+\sigma}^{q,\kappa},
\end{align}
}
for some $\sigma=\sigma(\tau_1,d,q)$.
 The next goal is to estimate ${\mathcal R}^1_{r,\lambda}$ whose expression is detailed in \eqref{hatR}. To do that, we shall first use the splitting
\begin{align}\label{Split-P}
{\mathcal R}^1_{r,\lambda}=\partial_\theta\mathcal{R}_{r,\lambda}+\widehat{\mathcal R}_1+\widehat{\mathcal R}_2,
\end{align}
with
\begin{align}\label{split-R}
\widehat{\mathcal R}_1=&\int_0^1\Phi(-t)\Big[\partial_\theta\mathcal{R}_{r,\lambda},\mathbb{A}\Big]\Phi(t)dt,\\
\nonumber\widehat{\mathcal R}_2=&\int_0^1\Phi(-t)\Big[\partial_\theta\mathcal{A}+(t-1)\widehat{\mathbb{A}},\mathbb{A}\Big]\Phi(t)dt.
\end{align}
Recall from  Proposition \ref{prop-chang} that
\begin{align}\label{G-V-L}
\interleave\partial_\theta \mathcal{R}_{r,\lambda}\interleave_{-1,s,\gamma}^{q,\kappa}\lesssim\varepsilon \kappa^{-1}\left(1+\| \mathfrak{I}_{0}\|_{s+\sigma_1+7+\gamma}^{q,\kappa}\right).
\end{align}
Applying  Lemma \ref{Lem-Commutator}-$($iv$)$ and using \eqref{self-ad1}
\begin{align*}
\interleave[\partial_\theta \mathcal{R}_{r,\lambda},\mathbb{A}\interleave_{-1,s,\gamma}^{q,\kappa}&\lesssim
\sum_{0\leqslant\beta\leqslant\gamma}\interleave\partial_\theta \mathcal{R}_{r,\lambda}\interleave_{-1,s_0+2,1+\beta}^{q,\kappa}\|\rho\|_{s+4+\gamma-\beta}^{q,\kappa}+\interleave\partial_\theta \mathcal{R}_{r,\lambda}\interleave_{-1,s+2,1+\beta}^{q,\kappa}\|\rho\|_{s_0+4+\gamma-\beta}^{q,\kappa}\\
&+\sum_{ 0\leqslant\beta\leqslant\gamma}\interleave \partial_\theta \mathcal{R}_{r,\lambda}\interleave_{-1,s_0+2+\beta,\gamma-\beta}^{q,\kappa}\|\rho\|_{s+1}^{q,\kappa} +\interleave \partial_\theta \mathcal{R}_{r,\lambda}\interleave_{-1,s+2+\beta,\gamma-\beta}^{q,\kappa}\|\rho\|_{s_0+1}^{q,\kappa}.
\end{align*}
Then using  \eqref{Est-eta} and   \eqref{G-V-L}  combined with Sobolev embeddings and \eqref{small-CC2}, we may find an explicit $\sigma=\sigma(\tau_1,d,q)$   such that
\begin{align*}
\interleave[\partial_\theta \mathcal{R}_{r,\lambda},\mathbb{A}\interleave_{-1,s,\gamma}^{q,\kappa}&\lesssim  \varepsilon \kappa^{-1}\big(1+\| \mathfrak{I}_{0}\|_{s+\sigma+\gamma}^{q,\kappa}\big)+
  \varepsilon \kappa^{-1}\sum_{0\leqslant\beta\leqslant\gamma}\|\mathfrak{I}_{0}\|_{s_0+\sigma+\beta}^{q,\kappa}\|\mathfrak{I}_{0}\|_{s+\sigma+\gamma-\beta}^{q,\kappa}\\
&+ \varepsilon \kappa^{-1}\sum_{ 0\leqslant\beta\leqslant\gamma}\|\mathfrak{I}_{0}\|_{s_0+\sigma+\gamma}^{q,\kappa}\|\mathfrak{I}_{0}\|_{s+\sigma}^{q,\kappa} \\
&\lesssim
  \varepsilon \kappa^{-1}\big(1+\| \mathfrak{I}_{0}\|_{s+\sigma+\gamma}^{q,\kappa}\big)+
  \varepsilon \kappa^{-1}\sum_{0\leqslant\beta\leqslant\gamma}\|\mathfrak{I}_{0}\|_{s_0+\sigma+\beta}^{q,\kappa}\|\mathfrak{I}_{0}\|_{s+\sigma+\gamma-\beta}^{q,\kappa}.
  \end{align*}
  We point out that along the computations below the value of  $\sigma$   may change from line to line but still only depends on $\tau_1,d$ and $q.$
Then using interpolation inequalities we find for any $0\leqslant\beta\leqslant\gamma$ and $s\geqslant s_0$ \begin{align*}
\|\mathfrak{I}_{0}\|_{s_0+\sigma+\beta}^{q,\kappa}\|\mathfrak{I}_{0}\|_{s+\sigma+\gamma-\beta}^{q,\kappa}
\lesssim\|\mathfrak{I}_{0}\|_{s_0+\sigma}^{q,\kappa}\|\mathfrak{I}_{0}\|_{s+\sigma+\gamma}^{q,\kappa}.
\end{align*}
Thus from the smallness condition \eqref{small-CC2}   we infer  for any $\gamma\in\NN, s\geqslant   s_0$ 
\begin{align}\label{YiL0}
 \interleave[\partial_\theta \mathcal{R}_{r,\lambda},\mathbb{A}]\interleave_{-1,s,\gamma}^{q,\kappa}&\lesssim  \varepsilon \kappa^{-1}\big(1+\| \mathfrak{I}_{0}\|_{s+\sigma+\gamma}^{q,\kappa}\big).
\end{align}
Applying Theorem \ref{Prop-EgorV}-(ii)
we get 
 for any $s\geqslant s_0$ and for any $t\in[0,1]$ 
\begin{align}\label{Tah0}
\interleave \Phi(-t)[\partial_\theta \mathcal{R}_{r,\lambda},\mathbb{A}]\Phi(t)\interleave_{-1,s,0}^{q,\kappa}\leqslant C e^{C\overline\mu^{s+1}(0)}\overline\mu(s)
\end{align}
with
$$
\overline\mu(s)\triangleq \overline{F}_0(s+27+s_0)+\|\rho\|_{{s+27+s_0}}^{q,\kappa}
$$
and 
$$
\overline{F}_0(s)\triangleq C\varepsilon \kappa^{-1}\big(1+\| \mathfrak{I}_{0}\|_{s+\sigma+4}^{q,\kappa}\big).
$$
Notice that $s\mapsto \overline{F}_0(s)$ is log-convex since Sobolev norms enjoy this property.  Using \eqref{Est-eta} combined with  Sobolev embeddings we deduce
$$ \overline\mu(s) \lesssim \varepsilon \kappa^{-1}\big(1+\| \mathfrak{I}_{0}\|_{s+\sigma}^{q,\kappa}\big)
$$
and consequently we get from \eqref{Tah0} and   \eqref{small-CC2} that for any $t\in[0,1]$
\begin{align*}
\interleave \Phi(-t)[\partial_\theta \mathcal{R}_{r,\lambda},\mathbb{A}]\Phi(t)\interleave_{-1,s,0}^{q,\kappa}\lesssim\, \varepsilon \kappa^{-1}\big(1+\| \mathfrak{I}_{0}\|_{s+\sigma}^{q,\kappa}\big).
\end{align*}
We remind that the value of  $\sigma$   may change from line to line but still depends only on $\tau_1,d$ and $q.$
Combining this estimate with the definition of $\widehat{\mathcal{R}}_1$ in \eqref{split-R} we obtain 
\begin{align}\label{Tah-m1}
\interleave\widehat{\mathcal{R}}_1\interleave_{-1,s,0}^{q,\kappa}\lesssim\,\varepsilon\, \kappa^{-1}\|r\|_{s+\sigma}^{q,\kappa}.
\end{align}
Let us move to the estimate of the first member of $\widehat{\mathcal R}_2$ in \eqref{split-R}. Then using Lemma \ref{Lem-Commutator}-(iii) combined with \eqref{etaM1} and \eqref{self-ad1}
\begin{eqnarray*}
\nonumber&&\interleave [\partial_\theta\mathcal{A},\mathbb{A}]\interleave_{2\overline\alpha-1+\epsilon,s,\gamma}^{q,\kappa}\leqslant C\big(\|\mu_{r,\lambda}\|_{s_0+3}^{q,\kappa}+\|\rho\|_{s_0+3}^{q,\kappa}\big)\big(\|\mu_{r,\lambda}\|_{s+\gamma+3}^{q,\kappa}+\|\rho\|_{s+\gamma+3}^{q,\kappa}\big).
\end{eqnarray*}
Applying \eqref{Est-eta} and the estimate of $\mu_{r,\lambda}$ stated in  Proposition \ref{prop-chang} together  with Sobolev embeddings and \eqref{small-CC2},  one gets for any $s\geqslant s_0, \gamma\in\N$
\begin{eqnarray*}
\interleave [\partial_\theta\mathcal{A},\mathbb{A}]\interleave_{2\overline\alpha-1+\epsilon,s,\gamma}^{q,\kappa}&\lesssim& \varepsilon \kappa^{-2}\big(1+\| \mathfrak{I}_{0}\|_{s_0+\sigma}^{q,\kappa}\big)\big(1+\| \mathfrak{I}_{0}\|_{s+\sigma+\gamma}^{q,\kappa}\big)\\
&\lesssim&  \varepsilon \kappa^{-2}\big(1+\| \mathfrak{I}_{0}\|_{s+\sigma+\gamma}^{q,\kappa}\big).
\end{eqnarray*}
Then according to  Theorem \ref{Prop-EgorV}-(ii),
we deduce that 
 for any $s\geqslant s_0$ and for any $t\in[0,1]$ 
\begin{align}\label{Tah0D}
\interleave \Phi(-t)[\partial_\theta\mathcal{A},\mathbb{A}]\Phi(t)\interleave_{2\overline\alpha-1+\epsilon,s,0}^{q,\kappa}\leqslant C e^{C\overline\mu^{s+1}(0)}\overline\mu(s)
\end{align}
with
$$
\overline\mu(s)=\overline{F}_0(s+27+s_0)+\|\rho\|_{{s+27+s_0}}^{q,\kappa}
$$
and 
$$
\overline{F}_0(s)=C\varepsilon \kappa^{-2}\big(1+\| \mathfrak{I}_{0}\|_{s+\sigma}^{q,\kappa}\big).
$$
Notice that $\overline{F}_0(s)$ is log-convex. Using \eqref{Est-eta} we find
$$
 \overline\mu(s) \lesssim\varepsilon \kappa^{-2}\big(1+\| \mathfrak{I}_{0}\|_{s+\sigma}^{q,\kappa}\big).
$$
Therefore we get from \eqref{Tah0D} and \eqref{small-CC2} that for any $t\in[0,1]$ and $s\geqslant s_0$
\begin{align}\label{Tah-m2}
\interleave \Phi(-t)[\partial_\theta\mathcal{A},\mathbb{A}]\Phi(t)\interleave_{2\overline\alpha-1+\epsilon,s,0}^{q,\kappa}\lesssim \varepsilon \kappa^{-2}\big(1+\| \mathfrak{I}_{0}\|_{s+\sigma}^{q,\kappa}\big).
\end{align}
For the second member of  $\widehat{\mathcal R}_2$ in \eqref{split-R}, we use  Lemma \ref{Lem-Commutator}-(iii) combined with \eqref{hatA} and \eqref{self-ad1}
\begin{eqnarray*}
\nonumber&&\interleave [\widehat{\mathbb{A}},\mathbb{A}]\interleave_{2\overline\alpha-1+\epsilon,s,\gamma}^{q,\kappa}\leqslant C\big(\|\rho\|_{s_0+3}^{q,\kappa}+\|\widehat\rho\|_{s_0+3}^{q,\kappa}\big)\big(\|\rho\|_{s+\gamma+3}^{q,\kappa}+\|\widehat\rho\|_{s+\gamma+3}^{q,\kappa}\big).
\end{eqnarray*}
Applying \eqref{Cod-v1} combined with the Lemma \ref{Law-prodX1}, \eqref{Tyk6} and \eqref{Est-eta} we infer 
\begin{align*}
\|\widehat\rho\|_{s}^{q,\kappa}&\lesssim \|\rho\|_{s+1}^{q,\kappa}\\
&\lesssim \varepsilon\kappa^{-2} \big(1+\| \mathfrak{I}_{0}\|_{s+\sigma}^{q,\kappa}\big).
\end{align*}
It follows from  \eqref{small-CC2} that
\begin{eqnarray*}
\nonumber&&\interleave [\widehat{\mathbb{A}},\mathbb{A}]\interleave_{2\overline\alpha-1+\epsilon,s,\gamma}^{q,\kappa}\leqslant \varepsilon\kappa^{-2} \big(1+\| \mathfrak{I}_{0}\|_{s+\sigma+\gamma}^{q,\kappa}\big).
\end{eqnarray*}
As before we can use  Theorem \ref{Prop-EgorV}-(ii)
and get 
 for any $s\geqslant s_0$ and $t\in[0,1]$ 
\begin{align*}
\interleave \Phi(-t) [\widehat{\mathbb{A}},\mathbb{A}]\Phi(t)\interleave_{2\overline\alpha-1+\epsilon,s,0}^{q,\kappa}\leqslant C e^{C\overline\mu^{s+1}(0)}\overline\mu(s)
\end{align*}
with
$$
\overline\mu(s)=\overline{F}_0(s+27+s_0)+\|\rho\|_{{s+27+s_0}}^{q,\kappa}
$$
and 
$$
\overline{F}_0(s)=C\varepsilon\kappa^{-2} \big(1+\| \mathfrak{I}_{0}\|_{s+\sigma}^{q,\kappa}\big).
$$
Therefore  we deduce  for any $t\in[0,1]$ and $s\geqslant s_0$
\begin{align}\label{Tah-m03}
\interleave \Phi(-t)[\widehat{\mathbb{A}},\mathbb{A}]\Phi(t)\interleave_{2\overline\alpha-1+\epsilon,s,0}^{q,\kappa}\leqslant C\varepsilon\kappa^{-2} \big(1+\| \mathfrak{I}_{0}\|_{s+\sigma}^{q,\kappa}\big).
\end{align}
Combining \eqref{split-R} with  \eqref{Tah-m2} and \eqref{Tah-m03} we find
\begin{align}\label{Tah-m3}
\interleave\widehat{\mathcal R}_2\interleave_{2\overline\alpha-1+\epsilon,s,0}^{q,\kappa}\leqslant C\varepsilon\kappa^{-2} \big(1+\| \mathfrak{I}_{0}\|_{s+\sigma}^{q,\kappa}\big).
\end{align}
Putting together \eqref{Split-P},\eqref{G-V-L},\eqref{Tah-m1} and \eqref{Tah-m3} yields
\begin{align*}
\interleave{\mathcal R}_{r,\lambda}^1\interleave_{2\overline\alpha-1+\epsilon,s,0}^{q,\kappa}\leqslant C\varepsilon\kappa^{-2} \big(1+\| \mathfrak{I}_{0}\|_{s+\sigma}^{q,\kappa}\big).
\end{align*}
 The estimate for the $L^2$ adjoint ${\mathcal R}_{r,\lambda}^{1\,\star}$ can be done using Lemma \ref{Lem-Rgv1}-(vi)
 \begin{align*}
\interleave{\mathcal R}_{r,\lambda}^{1\,\star}\interleave_{2\overline\alpha-1+\epsilon,s,0}^{q,\kappa}&\lesssim \interleave{\mathcal R}_{r,\lambda}^1\interleave_{2\overline\alpha-1+\epsilon,s+s_0+1,0}^{q,\kappa}\\
&\leqslant C\varepsilon\kappa^{-2} \big(1+\| \mathfrak{I}_{0}\|_{s+\sigma}^{q,\kappa}\big).
\end{align*}
Hence  we proved the first part of the second point (ii).\\
 For the difference operator $\Delta_{12}\mathcal {R}^1_{r,\lambda}$ we write according to \eqref{Split-P}
\begin{align*}
\Delta_{12}\mathcal{ R}^1_{r,\lambda}=&\Delta_{12}\partial_\theta\mathcal{R}_{r,\lambda}+\Delta_{12}\widehat{\mathcal R}_1+\interleave\Delta_{12}\widehat{\mathcal R}_2,
\end{align*}
 The estimate of $\Delta_{12}\partial_\theta\mathcal{R}_{r,\lambda}$ is given in  Proposition \ref{prop-chang}. Thus it remains to estimate the last two terms. The arguments are similar for both terms and we shall only discuss how to get the suitable  estimate for $\Delta_{12}\widehat{\mathcal R}_1$. We start with writing
 \begin{align*}
\nonumber\Delta_{12}\widehat{\mathcal R}_1=&\int_0^1\Phi_1(-t)\Big\{\Delta_{12}\big[\partial_\theta\mathcal{R}_{r,\lambda},\mathbb{A}\big]\Big\}\Phi_1(t)dt+\int_0^1\Delta_{12}\Big\{\Phi(-t)\big[\partial_\theta\mathcal{R}_{r_2,\lambda},\mathbb{A}_2\big]\Phi(t)\Big\}dt,
\end{align*}
 Let us estimate the first term of the right hand side. Then similarly to \eqref{YiL0},  using in particular the estimate \eqref{puir1} combined with the law products we obtain 
 \begin{align}\label{YiL000}
\nonumber\forall\, s+\gamma\leqslant \overline{s}_h,\quad \interleave\Delta_{12}[\partial_\theta \mathcal{R}_{r,\lambda},\mathbb{A}]\interleave_{-1,s,\gamma}^{q,\kappa}&\lesssim \varepsilon\kappa^{-1}\|\Delta_{12}i\|_{\overline{s}_{h}+\sigma}^{q,\kappa}\\
&\lesssim \varepsilon\kappa^{-1}\|\Delta_{12}i\|_{\overline{s}_{h}+\sigma+s+\gamma}^{q,\kappa}\end{align}
 Therefore applying Theorem \ref{Prop-EgorV}-(ii) and using \eqref{YiL000} lead in view of \eqref{small-CC2} to 
 \begin{align}\label{YiLP1}
\quad\int_0^1\interleave\Phi_1(-t) \Delta_{12}[\partial_\theta \mathcal{R}_{r,\lambda},\mathbb{A}]\Phi_1(t) \interleave_{-1,\overline{s}_h-s_0-27,0}^{q,\kappa}&\lesssim \varepsilon\kappa^{-1}\|\Delta_{12}i\|_{2\overline{s}_{h}+\sigma}^{q,\kappa}.
\end{align}
On the other hand, we apply Theorem \ref{Prop-EgorV}-(iii) combined with \eqref{YiL0} and \eqref{west-eta}
\begin{align}\label{YiLP2}
\nonumber\quad\int_0^1\interleave\Delta_{12}\Phi(-t)\big[\partial_\theta\mathcal{R}_{r_2,\lambda},\mathbb{A}_2\big]\Phi(t)\interleave_{0,0,0}^{q,\kappa}&\lesssim \|\Delta_{12}\rho\|_{s_0+5}^{q,\kappa}
\\
&\lesssim \varepsilon\kappa^{-2}\|\Delta_{12}i\|_{\overline{s}_{h}+\sigma}^{q,\kappa}.
\end{align}
Then combining \eqref{YiLP1} and \eqref{YiLP2} yields
\begin{align}\label{YiLP23}
\interleave\Delta_{12}\widehat{\mathcal R}_1\interleave_{0,0,0}^{q,\kappa}&\lesssim \varepsilon\kappa^{-2}\|\Delta_{12}i\|_{\overline{s}_{h}+\sigma}^{q,\kappa}.
\end{align}
Remind once again that the value of $\sigma$ may change from line to line and it depends only on $\tau_1, d$ and $q.$
The estimates of $\Delta_{12}\mathcal{ R}^1_{r,\lambda}h$ and $\Delta_{12}\mathcal{ R}^{1,\star}_{r,\lambda}h$ can be implemented in a straightforward way where we use in particular the preceding estimates combined with Lemma \ref{Lem-Rgv1},  Proposition \ref{flowmap00}-(iii) \mbox{and \eqref{es-phi-m1}.}

\smallskip

{\bf{(iii)}} Recall that  $\Psi=\mathscr{B}\Phi(1)$, then  using Proposition \ref{QP-change}-(ii) yields
\begin{equation*}
\|\Psi h\|_{s}^{q,\kappa}\lesssim\| \Phi(1)h\|_{s}^{q,\kappa}+\varepsilon\kappa^{-1}\| \mathfrak{I}_{0}\|_{s+\sigma}^{q,\kappa}\| \Phi(1)h\|_{s_{0}}^{q,\kappa}.
\end{equation*}
Applying Proposition \ref{flowmap00}-(ii) combined with \eqref{Est-eta},  the smallness condition \eqref{small-CC2} and Sobolev embeddings  we obtain
\begin{align*}
\| \Phi(1)h\|_{s}^{q,\kappa}&\lesssim\| h\|_{s}^{q,\kappa}+\varepsilon\kappa^{-2}\big(1+\| \mathfrak{I}_{0}\|_{s+\sigma}^{q,\kappa}\big)\| h\|_{s_{0}}^{q,\kappa}\\
&\lesssim\| h\|_{s}^{q,\kappa}+\varepsilon\kappa^{-2}\| \mathfrak{I}_{0}\|_{s+\sigma}^{q,\kappa}\| h\|_{s_{0}}^{q,\kappa},
\end{align*}
which implies 
 \begin{align}\label{Est-etaY01}
 \|\Psi h\|_{s}^{q,\kappa}&\lesssim\| h\|_{s}^{q,\kappa}+\varepsilon\kappa^{-2}\| \mathfrak{I}_{0}\|_{s+\sigma}^{q,\kappa}\| h\|_{s_{0}}^{q,\kappa}.
\end{align}
We point out that one would expect in the preceding estimate to get $\| h\|_{s_{0}+q}^{q,\kappa}$ instead of $\| h\|_{s_{0}}^{q,\kappa}$.  However according to Proposition \ref{flowmap00}-(ii), $s_0>\frac{d+5}{2}$ and then  $s_0+q>\frac{d+5}{2}+q$. Therefore one may replace $s_0+q$ by $s_0$ by virtue of the assumption \eqref{Conv-T2}. \\
Similar arguments give also 
\begin{align}\label{Est-etaY02}
\|\Psi^{-1} h\|_{s}^{q,\kappa}\lesssim \|h\|_{s}^{q,\kappa}+\varepsilon\kappa^{-2}\| \mathfrak{I}_{0}\|_{s+\sigma}^{q,\kappa}\| h\|_{s_{0}}^{q,\kappa}.
\end{align}
On the other hand one has
$$
\Psi-\textnormal{Id}=\mathscr{B}\big(\Phi(1)-\textnormal{Id}\big)+\big(\mathscr{B}-\textnormal{Id}\big).
$$
For the last part we use Proposition \ref{QP-change} 
\begin{equation*}
\big\|\left(\mathscr{B}-\textnormal{Id}\right)h\big\|_{s}^{q,\kappa}\lesssim \varepsilon\kappa^{-1}\left(\| h\|_{s+1}^{q,\kappa}+\| \mathfrak{I}_{0}\|_{s+\sigma}^{q,\kappa}\| h\|_{s_{0}}^{q,\kappa}\right).
\end{equation*}
The same proposition gives
\begin{equation*}
\big\|\mathscr{B}\big(\Phi(1)-\textnormal{Id}\big)h\big\|_{s}^{q,\kappa}\lesssim\| (\Phi(1)-\textnormal{Id})h\|_{s}^{q,\kappa}+\varepsilon\kappa^{-1}\|  \mathfrak{I}_{0}\|_{s+\sigma}^{q,\kappa}\| (\Phi(1)-\textnormal{Id})h\|_{s_{0}}^{q,\kappa}.
\end{equation*}
Applying once again Proposition \ref{flowmap00} combined  with \eqref{Est-eta} and  \eqref{small-CC2} yields
$$
\|\mathscr{B}\big(\Phi(1)-\textnormal{Id}\big)h\|_{s}^{q,\kappa}\lesssim \varepsilon {\kappa^{-2}}\|h\|_{s+1}^{q,\kappa}+\varepsilon{\kappa^{-2}}\| \mathfrak{I}_{0}\|_{s+\sigma}^{q,\kappa}\| h\|_{s_0}^{q,\kappa}.$$
Putting together the preceding estimates  gives for any $s\geqslant s_0$
\begin{equation*}
\big\|\big(\Psi-\textnormal{Id}\big)h\big\|_{s}^{q,\kappa}\lesssim\varepsilon {\kappa^{-2}}\|h\|_{s+1}^{q,\kappa}+\varepsilon{\kappa^{-2}}\| \mathfrak{I}_{0}\|_{s+\sigma}^{q,\kappa}\| h\|_{s_0}^{q,\kappa}.
\end{equation*}
In a similar way we obtain
$$
\|(\Psi^{-1}-\textnormal{Id})h\|_{s}^{q,\kappa}+\|((\Psi^{\pm1})^{\star}-\textnormal{Id})h\|_{s}^{q,\kappa}\lesssim \varepsilon {\kappa^{-2}}\|h\|_{s+1}^{q,\kappa}+\varepsilon{\kappa^{-2}}\| \mathfrak{I}_{0}\|_{s+\sigma}^{q,\kappa}\| h\|_{s_0}^{q,\kappa}.
$$
Let us now move to the estimate of the difference $\Delta_{12}\Psi$ and sketch the main arguments. According to Proposition \ref{flowmap00}-(iii), \eqref{Est-eta}, \eqref{west-eta} and \eqref{small-CC2} one gets
\begin{align}\label{es-phi-m1}
\| \Delta_{12}\Phi(1)h\|_{\overline{s}_h}^{q,\kappa}&\lesssim \varepsilon{\kappa^{-2}}\|\Delta_{12}i\|_{\overline{s}_h+\sigma}^{q,\kappa} \| h\|_{\overline{s}_h+\sigma}^{q,\kappa}.
\end{align}
To estimate $\Delta_{12}\mathscr{B}$ where $\mathscr{B}$ is defined in \eqref{mathscrB} we first write
\begin{align*}
\Delta_{12}\mathscr{B}=(\partial_{\theta}\Delta_{12}\beta)\textbf{B}_1+(1+\partial_{\theta}\beta_2)\Delta_{12}\textbf{B}.
\end{align*}
For the first term of the right hand side we use the law products combined with Lemma \ref{Compos1-lemm},   Proposition \ref{QP-change}-(iv) and \eqref{small-CC2}
\begin{align*}
\|(\partial_{\theta}\Delta_{12}\beta)\textbf{B}_1h\|_{\overline{s}_h}^{q,\kappa}\lesssim \varepsilon{\kappa^{-1}}\|\Delta_{12}i\|_{\overline{s}_h+\sigma}^{q,\kappa} \| h\|_{\overline{s}_h}^{q,\kappa}.
\end{align*}
Using the fundamental theorem of calculus we infer
\begin{align*}
\Delta_{12}\textbf{B}h=(\beta_1-\beta_2)\int_0^1(\partial_\theta h)\big(\cdot,(1-\tau)\beta_1+\tau\beta_2\big) d\tau.
\end{align*}
Then the same arguments as for the preceding estimate allow to get
\begin{align*}
\|\Delta_{12}\textbf{B}h\|_{\overline{s}_h}^{q,\kappa}\lesssim \varepsilon{\kappa^{-1}}\|\Delta_{12}i\|_{\overline{s}_h+\sigma}^{q,\kappa} \| h\|_{\overline{s}_h+1}^{q,\kappa}.
\end{align*}
Putting together the foregoing estimates imply
\begin{align}\label{es-bet-m1}
\|\Delta_{12}\mathscr{B}h\|_{\overline{s}_h}^{q,\kappa}\lesssim \varepsilon{\kappa^{-1}}\|\Delta_{12}i\|_{\overline{s}_h+\sigma}^{q,\kappa} \| h\|_{\overline{s}_h+1}^{q,\kappa}.
\end{align}
Next we use the identity
$$
\Delta_{12}\Psi=(\Delta_{12}\mathscr{B})\Phi_1(1)+\mathscr{B}_2\Delta_{12}\Phi(1)\quad\hbox{with}\quad \Psi=\mathscr{B}\Phi(1).
$$
Applying \eqref{es-bet-m1} with Proposition \ref{flowmap00}-(ii), \eqref{Est-eta} and \eqref{small-CC2} leads to
\begin{align*}
\|(\Delta_{12}\mathscr{B})\Phi_1(1)h\|_{\overline{s}_h}^{q,\kappa}&\lesssim \varepsilon{\kappa^{-1}}\|\Delta_{12}i\|_{\overline{s}_h+\sigma}^{q,\kappa} \| \Phi_1(1)h\|_{\overline{s}_h+1}^{q,\kappa}\\
&\lesssim \varepsilon{\kappa^{-1}}\|\Delta_{12}i\|_{\overline{s}_h+\sigma}^{q,\kappa} \| h\|_{\overline{s}_h+1}^{q,\kappa}.
\end{align*}
Similarly, \eqref{es-phi-m1} combined with \eqref{estimate on the first reduction operator and its inverse} and \eqref{small-CC2} give
\begin{align*}
\|\mathscr{B}_2\Delta_{12}\Phi(1)h\|_{\overline{s}_h}^{q,\kappa}&\lesssim \|\Delta_{12}\Phi(1)h\|_{\overline{s}_h}^{q,\kappa}\\
&\lesssim \varepsilon{\kappa^{-2}}\|\Delta_{12}i\|_{\overline{s}_h+\sigma}^{q,\kappa} \| h\|_{\overline{s}_h+\sigma}^{q,\kappa}.
\end{align*}
Consequently, we obtain
\begin{align*}
\|\Delta_{12}\Psi h\|_{\overline{s}_h}^{q,\kappa}
&\lesssim \varepsilon{\kappa^{-2}}\|\Delta_{12}i\|_{\overline{s}_h+\sigma}^{q,\kappa} \| h\|_{\overline{s}_h+\sigma}^{q,\kappa}.
\end{align*}
The estimates of $\Delta_{12}\Psi^{-1},\Delta_{12}\Psi^{\star}$ and $\Delta_{12}(\Psi^{\star})^{-1}$ are quite similar to $\Delta_{12}\Psi.$ We mainly use the same arguments supplemented in particular  with Lemma \ref{algeb1}-(iii) and Proposition \ref{flowmap00}-(ii)-(iii).

\smallskip

{\bf{(iv)}} Recall from \eqref{En-11} that
\begin{align}\label{En-101}
 \nonumber\mathtt{E}_n^1&=\mathtt{E}_n^0+\int_0^1\Phi(-t)\big[\mathtt{E}_n^0,\mathbb{A}\big]\Phi(t)dt-\partial_\theta\big[(\Pi_{N_n}^\perp\mu_{r,\lambda})|\textnormal{D}|^{\alpha-1}+|\textnormal{D}|^{\alpha-1} (\Pi_{N_n}^\perp \mu_{r,\lambda})\big]\\
 \nonumber&\qquad\qquad\qquad +\partial_\theta\big({\bf R}_n|\textnormal{D}|^{\alpha-1}+|\textnormal{D}|^{\alpha-1} {\bf R}_n\big)\\
&\triangleq\mathtt{E}_n^0+\mathtt{E}_n^{0,1}-\mathtt{E}_n^{0,2}+\mathtt{E}_n^{0,3}.
\end{align}
The estimate of $\mathtt{E}_n^0$ is stated in Proposition \ref{QP-change}. For the second one we may use the law products combined with \eqref{Est-eta},  Proposition \ref{QP-change} and \eqref{small-CC2}  in order to get 
\begin{align*}
\big\|\big[\mathtt{E}_n^0,\mathbb{A}\big]h\big\|_{s_0}^{q,\kappa}&\lesssim \|\mathtt{E}_n^0\mathbb{A}h\big\|_{s_0}^{q,\kappa}+\varepsilon\kappa^{-2}\|\mathtt{E}_n^0h\|_{s_0+1}^{q,\kappa}
\\
&\lesssim\varepsilon N_{0}^{\mu_{2}}N_{n+1}^{-\mu_{2}}\|h\|_{q,s_{0}+3}^{q,\kappa}.
\end{align*}
Combining this estimate with  Proposition \ref{flowmap00} and \eqref{small-C2} implies
$$\big\|\mathtt{E}_n^{0,1}h\big\|_{s_0}^{q,\kappa}\lesssim\varepsilon N_{0}^{\mu_{2}}N_{n+1}^{-\mu_{2}}\|h\|_{q,s_{0}+3}^{q,\kappa}.$$
For the estimate of $\mathtt{E}_n^{0,2},$ we have
from the law products, the projector properties and \mbox{Proposition \ref{prop-chang}}
\begin{align*}
\forall\, s\geqslant s_0,\,\quad\big\|\mathtt{E}_n^{0,2}h\big\|_{s_0}^{q,\kappa}\lesssim&\|(\Pi_{N_n}^\perp(\mu_{r,\lambda}-\langle\mu_{r,\lambda}\rangle )\|_{s_0+1}^{q,\kappa}\|h\|_{s_0+1}^{q,\kappa}\\
&\lesssim N_{n}^{s-s_0} \|\mu_{r,\lambda}-\langle\mu_{r,\lambda}\rangle\|_{s+1}^{q,\kappa} \| h\|_{s_0+1}^{q,\kappa}
\\
&\quad \lesssim  \varepsilon{\kappa^{-1}}  N_{n}^{s-s_0}\big(1+\| \mathfrak{I}_{0}\|_{s+\sigma}^{q,\kappa}\big)\| h\|_{s_0+1}^{q,\kappa}.
\end{align*}
As to the estimate of $\mathtt{E}_n^{0,3}$ in \eqref{En-101}, it is quite  similar to $\mathtt{E}_n^{0,2},$ and one gets first
\begin{align}\label{Fr-Bel-0-2}
\big\|\mathtt{E}_n^{0,3}h\big\|_{s_0}^{q,\kappa}&\lesssim\|{\bf R}_n\|_{s_0+1}^{q,\kappa}\|h\|_{s_0+1}^{q,\kappa}.
\end{align}
Then, coming back to the definition of ${\bf R}_n$ in \eqref{Eq-eta}, we obtain by virtue of  the estimate $c(\lambda,i_0)$   in Proposition \ref{QP-change}-$($i$)$,  \eqref{small-C2} and \eqref{Est-eta}
\begin{align*}
\forall\, s\geqslant s_0,\,\quad\|{\bf R}_n\|_{s_0+1}^{q,\kappa}&\lesssim \|\Pi_{N_n}^\perp\rho\|_{s_0+2}^{q,\kappa}\\
&\lesssim N_n^{s_0-s}\|\Pi_{N_n}^\perp\rho\|_{s+2}^{q,\kappa}\\
&\lesssim\varepsilon{\kappa^{-2}}  N_{n}^{s-s_0}\big(1+\| \mathfrak{I}_{0}\|_{s+\sigma}^{q,\kappa}\big).
\end{align*}
Plugging this estimate into \eqref{Fr-Bel-0-2} implies
\begin{align*}
\forall\, s\geqslant s_0,\,\quad\big\|\mathtt{E}_n^{0,3}h\big\|_{s_0}^{q,\kappa}&\lesssim \varepsilon{\kappa^{-2}}  N_{n}^{s-s_0}\big(1+\| \mathfrak{I}_{0}\|_{s+\sigma}^{q,\kappa}\big)\|h\|_{s_0+1}^{q,\kappa}.
\end{align*}
Combining the preceding estimates yields in view of \eqref{small-CC2}
$$\big\|\mathtt{E}_n^{1}h\big\|_{s_0}^{q,\kappa}\lesssim\varepsilon N_{0}^{\mu_{2}}N_{n+1}^{-\mu_{2}}\|h\|_{q,s_{0}+3}^{q,\kappa}+\varepsilon{\kappa^{-2}}  N_{n}^{s-s_0}\big(1+\| \mathfrak{I}_{0}\|_{s+\sigma}^{q,\kappa}\big)\| h\|_{s_0+1}^{q,\kappa}.$$
This completes the proof of Proposition \ref{prop-constant-coe}.
 \end{proof}

\subsection{Complete reduction up to small errors}\label{Norm-sec-1}
The main concern of this section is to reduce to a Fourier multiplier  the operator $\widehat{\mathcal{L}}_\omega$ described by  \eqref{Norm-proj} and \mbox{Proposition \ref{lemma-GS0}} and which takes the form (To alleviate the notation we replace therein $r_0$ by $r$)
\begin{align}\label{Norm-local-z}
\widehat{\mathcal{L}}_{\omega}=\Pi_{\mathbb{S}_0}^{\perp}\left(\mathcal{L}_{\varepsilon r,\lambda}-\varepsilon\partial_{\theta}\mathcal{R}\right)| \Pi_{\mathbb{S}_0}^\bot.
\end{align}
Notice that the set $\mathbb{S}_0$ is defined in \eqref{tangent-set2}. 
 According to  Proposition \ref{prop-constant-coe}, we have seen that when the parameter $\lambda$  belongs to 
 the Cantor like set  $ \mathcal{O}_{\infty,n}^{\kappa,\tau_{1}}(i_{0})$ then 
\begin{align}\label{dida-lip1}
\Psi^{-1}\mathcal{L}_{\varepsilon r,\lambda}\Psi=\mathcal{L}_{r,\lambda}^1=\omega\cdot\partial_\varphi+{c}(\lambda,i_0)\partial_\theta-\textnormal{m}(\lambda,i_0)\partial_\theta|\textnormal{D}|^{\alpha-1}+\mathcal{R}^1_{ r,\lambda}+\mathtt{E}_n^1
\end{align}
where ${\Psi}:{{\mathcal{O}}}\to \mathscr{L}\big(H^{s}(\T^{d+1}),H^{s}(\T^{d+1})\big) $ is  a family of invertible linear operators  preserving the symmetry and satisfying  suitable tame estimates.\\
In \eqref{tangent-set2} we introduced the set $\mathbb{S}_{0}=\mathbb{S}\cup(-\mathbb{S})\cup\{0\}$ and the space $ H_{\perp}^{s}$  is described through \eqref{decoacca}. Now, we consider the orthogonal projector  $\Pi_{\mathbb{S}_{0}}:L^{2}(\mathbb{T}^{d+1},\mathbb{C})\to L^{2}(\mathbb{T}^{d+1},\mathbb{C})$ given by 
\begin{align}\label{orth-proj-1}
h(\varphi,\theta)=\sum_{j\in\mathbb{Z}}h_{j}(\varphi){\bf{e}}_j(\theta)\Longrightarrow\Pi_{\mathbb{S}_{0}}h(\varphi,\theta)=\sum_{j\in\mathbb{S}_{0}}h_{j}(\varphi)\,{\bf{e}}_j(\theta),\, {\bf{e}}_j(\theta)=e^{\ii j\theta}
\end{align}
and $\Pi_{\mathbb{S}_{0}}^{\perp}=\textnormal{Id}-\Pi_{\mathbb{S}_{0}}.$
Define 
$$\Psi_{\perp}\triangleq\Pi_{\mathbb{S}_{0}}^{\perp}\Psi\Pi_{\mathbb{S}_{0}}^{\perp}:{{\mathcal{O}}}\longrightarrow \mathscr{L}\big(H_{\perp}^{s}(\T^{d+1}),H_{\perp}^{s}(\T^{d+1})\big).
$$
We shall discuss two points. The first one is related to the invertibility of $\Psi_{\perp}$ with suitable tame estimates. We will see that the operators $\Psi_{\perp}^{\pm1}$  can be viewed as finite rank perturbations of $\Psi^{\pm1}$. This allows to get a precise description of  the truncated operator
\begin{align}\label{Lin-orthog}
\Psi_{\perp}^{-1} \widehat{\mathcal{L}}_{\omega}  \Psi_{\perp}=\Psi_{\perp}^{-1}\mathcal{L}_{\varepsilon r,\lambda}^0\Psi_{\perp}-\varepsilon\Psi_{\perp}^{-1}\partial_{\theta}\mathcal{R}\Psi_{\perp}.
\end{align}
 Actually, one deduces that $\Psi_{\perp}^{-1} \widehat{\mathcal{L}}_{\omega}  \Psi_{\perp}$ is a small finite rank perturbation of  the operator $\mathcal{L}_{\varepsilon r,\lambda}$.  Then the second main goal is to construct  an approximate  right inverse of  this latter  operator with tame estimates.  For this aim we first use the reductions implemented before throughout Sections \ref{sec-transport-1} and \ref{Sec-Local-P}. By this way, the positive order part of the new operator is a Fourier multiplier. Thus, at this stage it remains to reduce the remainder which of order zero and enjoying  a good behavior with respect to suitable norms seen in \eqref{Def-pseud-w}. The remainder reduction follows standard approach  as in \cite{BertiMontalto, BCP}.
\subsubsection{Frequency localization  of operators}
We shall be concerned with some analytical properties of operators generated by  frequency localization of given operators. Before stating the result, we shall first fix the framework. Let ${\Psi}:{{\mathcal{O}}}\to \mathscr{L}\big(H^{s}(\T^{d+1}),H^{s}(\T^{d+1})\big) $ be  an arbitrary  family  of invertible linear operators (not necessary those of \eqref{dida-lip1}) and let $\mathbb{S}_0\subset\mathbb{Z}$ be a finite set as \mbox{in \eqref{tangent-set2}.} Consider the orthogonal projectors $\Pi_{\mathbb{S}_0}$ and $\Pi_{\mathbb{S}_0}^\perp$defined \mbox{by \eqref{orth-proj-1}} and let us introduce the restricted  transformation 
$$\Psi_{\perp}\triangleq\Pi_{\mathbb{S}_{0}}^{\perp}\Psi\Pi_{\mathbb{S}_{0}}^{\perp}:{{\mathcal{O}}}\longrightarrow \mathscr{L}\big(H_{\perp}^{s}(\T^{d+1}),H_{\perp}^{s}(\T^{d+1})\big).
$$
Consider for $\lambda\in {\mathcal{O}}$  the finite-dimensional matrix
\begin{align}\label{M-atrix}
{\bf{M}} (\lambda,\varphi)=\left(\Big\langle {\bf e}_m,{\big(\Psi^\star(\lambda)\big)}^{-1}{\bf{e}}_k\Big\rangle_{L^2_\theta(\T)}\right)_{m,k\in\mathbb{S}_0},\quad {\bf{e}}_k(\theta)=e^{\ii k\theta}
\end{align}
where $\Psi^\star$ is the $L^2_\theta-$adjoint of $\Psi.$ The scalar product $\langle\cdot\rangle_{L^2_\theta(\T)}$ is defined by
$$
\langle f,g\rangle_{L^2_\theta(\T)}=\frac{1}{2\pi}\int_{\T} f(\theta)\overline{g(\theta)}d\theta.
$$
Remark that the matrix entries depend only  on $\lambda$ and $\varphi$ but not on $\theta$ because we are taking the average  over $\theta.$
Our main task is to establish  the following result.
\begin{lemma}\label{lemm-decompY}
The following assertions hold true.
\begin{enumerate}
\item 
Assume that  $\big[\Psi,\rho\big]=0$ for any smooth real function $\rho: \varphi\in \T^d\mapsto \rho(\varphi)\in\RR,$ and 
$$
 \forall (\lambda,\varphi)\in {{\mathcal{O}}}\times \T^{d},\quad \textnormal{det } {\bf{M}} (\lambda,\varphi)\neq0.
$$
Then,  for any  $ h\in H_{\perp}^{s}$ 
$$
\Psi_{\perp}^{-1}h=\Psi^{-1}h-\sum_{m\in\mathbb{S}_0}\big\langle h, \big({(\Psi^\star)}^{-1}-\textnormal{Id}\big)g_m  \big\rangle_{L^2_\theta(\T)}\Psi^{-1}{\bf e}_m
$$
where the functions $g_m $ are defined via
$$
   {\bf {M}}^{-1}(\lambda,\varphi)=\Big(\alpha_{m,k}(\lambda,\varphi)\Big)_{(m,k)\in\mathbb{S}_0^2}\quad\hbox{and}\quad g_m(\lambda,\varphi,\theta)\triangleq \sum_{
k\in\mathbb{S}_0}\overline{\alpha_{m,k}}(\lambda,\varphi){\bf e}_k(\theta).
$$
\item Let $\Psi$ be the  operator constructed in Proposition $\ref{prop-constant-coe}$ and assume \eqref{small-CC2},  then 
\begin{align*}
\forall\, s\in[s_0,S],\quad \|\Psi_{\perp}^{\pm1}g\|_{s}^{q,\kappa}&\lesssim \| g\|_{s}^{q,\kappa}+\varepsilon\kappa^{-2}\| \mathfrak{I}_{0}\|_{s+\sigma}^{q,\kappa}
\| g\|_{s_0}^{q,\kappa}.
\end{align*}
In addition, we have
\begin{align*}
\| g_m\|_{s}^{q,\kappa}\lesssim\big(1+\varepsilon\kappa^{-2}\| \mathfrak{I}_{0}\|_{s+\sigma}^{q,\kappa}\big)
\end{align*}
and 
\begin{align*}
\max_{m\in\mathbb{S}_0}\| \Delta_{12}g_m\|_{\overline{s}_h}^{q,\kappa}\lesssim \varepsilon {\kappa^{-2}}\|\Delta_{12}i\|_{\overline{s}_h+\sigma}^{q,\kappa}.
\end{align*}

\end{enumerate}
\end{lemma}
\begin{proof}

{\bf{(i)}}
 Given $g\in W^{q,\infty}_{\kappa}(\mathcal{O}, H_{\perp}^{s})$, we want  to solve the equation
$$
f\in W^{q,\infty}_{\kappa}(\mathcal{O}, H_{\perp}^{s}),\quad  \Psi_{\perp}f=\Pi_{\mathbb{S}_{0}}^{\perp}\Psi\Pi_{\mathbb{S}_{0}}^{\perp}f=g.
$$ 
This is equivalent to
$$
\Psi f=g+h, \quad \hbox{with}\quad \Pi_{\mathbb{S}_{0}}h=h\quad\hbox{and}\quad  \Pi_{\mathbb{S}_{0}}f=0.
$$
Then we get 
\begin{align}\label{f-Id}
f= \Psi^{-1}\big(g+h\big), \quad \hbox{with}\quad \Pi_{\mathbb{S}_{0}}h=h,
\end{align}
provided that  $\Pi_{\mathbb{S}_{0}}f=0$, that is,
$$
\big\langle\Psi^{-1}\big(g+h\big), {\bf e}_k   \big\rangle_{L^2_\theta(\T)}=0,\quad \forall k\in \mathbb{S}_{0}.
$$
Making appeal to the identity $(\Psi^\star)^{-1}=(\Psi^{-1})^\star$, this  equation is equivalent to  
\begin{align}\label{Def-tra1}
\big\langle g+h, \widehat{{\bf e}}_k   \big\rangle_{L^2_\theta(\T)}=0,\quad \forall k\in \mathbb{S}_{0},\quad\hbox{with}\quad \widehat{{\bf e}}_k\triangleq (\Psi^\star(\lambda))^{-1}{\bf e}_k.
\end{align}
These  constraints will uniquely determine  $h$. Indeed,  by expanding  $h$ in the form 
$$ h(\varphi,\theta)=\sum_{m\in \mathbb{S}_{0}} a_m(\lambda,\varphi){\bf e}_m(\theta),$$we can transform the preceding system  into
\begin{align}\label{matrix-inv}
\sum_{m\in \mathbb{S}_{0}} a_m(\varphi)\big\langle {\bf e}_m, \widehat{{\bf e}}_k \big\rangle_{L^2_\theta(\T)}=-\big\langle g, \widehat{{\bf e}}_k   \big\rangle_{L^2_\theta(\T)},\quad \forall k\in \mathbb{S}_{0}.
\end{align}
By assumption, the matrix  ${\bf {M}}(\lambda,\varphi)$ defined in \eqref{M-atrix} is invertible. Therefore the system \eqref{matrix-inv} is invertible and one gets a unique solution given by
$$
a_m(\lambda,\varphi)=-\sum_{k\in\mathbb{S}_0}\alpha_{m,k}(\lambda,\varphi)\big\langle g, \widehat{{\bf e}}_k   \big\rangle_{L^2_\theta(\T)}\quad\hbox{with}\quad {\bf {M}}^{-1}(\lambda,\varphi)\triangleq \Big(\alpha_{m,k}(\lambda,\varphi)\Big)_{(m,k)\in\mathbb{S}_0^2}.
$$
This implies
\begin{align*}
h(\varphi,\theta)=-\sum_{m\in \mathbb{S}_{0}} \big\langle g, \widehat{g}_m   \big\rangle_{L^2_\theta(\T)}{\bf e}_m(\theta)\quad\hbox{with}\quad \widehat{g}_m\triangleq  \sum_{k\in\mathbb{S}_0}\overline{\alpha_{m,k}}\,\widehat{{\bf e}}_k,   
\end{align*}
where we denote by $\overline{\alpha_{m,k}}$ the complex conjugate of ${\alpha_{m,k}}$.
Using the commutation assumption $[\Psi,\rho]=0$ we also get  $[\Psi^{-1},\rho]=0$ and therefore
\begin{align}\label{matrix-inv01}
\Psi^{-1}h=-\sum_{m\in \mathbb{S}_{0}} \big\langle g, \widehat{g}_m   \big\rangle_{L^2_\theta(\T)}\Psi^{-1}{\bf e}_m.   
\end{align}
Notice that one also has  $\big[(\Psi^\star)^{-1},\rho\big]=0$, which implies in view of \eqref{Def-tra1}
\begin{align}\label{Mer-09}
\widehat{g}_m=(\Psi^\star(\lambda))^{-1}g_m\quad\hbox{with}\quad g_m\triangleq  \sum_{k\in\mathbb{S}_0}\overline{\alpha_{m,k}}\,{{\bf e}}_k.
\end{align}
Thus from \eqref{matrix-inv01} and \eqref{f-Id} we deduce the formula
\begin{align}\label{Mer-11}
 \Psi_{\perp}^{-1}g&=\Psi^{-1}g-\sum_{m\in \mathbb{S}_{0}} \big\langle g, (\Psi^\star)^{-1}g_m   \big\rangle_{L^2_\theta(\T)}\Psi^{-1}{\bf e}_m.
 \end{align}
Since $\Pi_{\mathbb{S}_{0}}^{\perp}g=g $ and  $\Pi_{\mathbb{S}_{0}}^{\perp}g_m=0 $ then $\big\langle g,g_m \big\rangle_{L^2_\theta(\T)}=0$ and therefore
$$
\big\langle g,(\Psi^\star)^{-1}g_m \big\rangle_{L^2_\theta(\T)}=\big\langle g,\big((\Psi^\star(\lambda))^{-1}-\textnormal{Id}\big)g_m \big\rangle_{L^2_\theta(\T)}.
$$
Plugging this identity into \eqref{Mer-11} yields
\begin{align}\label{Decomp-ZZ1}
\Psi_{\perp}^{-1}g&=\Psi^{-1}g-\sum_{m\in\mathbb{S}_0}\big\langle g,\big((\Psi^\star)^{-1}-\textnormal{Id}\big)g_m \big\rangle_{L^2_\theta(\T)}\Psi^{-1}{\bf e}_m.
\end{align}
This ends the proof of the first point.

\smallskip

{\bf{(ii)}} The estimate of $\Psi_{\perp}$ follows from the continuity of the orthogonal projection $\Pi_{\mathbb{S}_0}^\perp$ on  Sobolev spaces combined with Proposition \ref{prop-constant-coe}-(iii). Let us now move to the estimate of $\Psi_{\perp}^{-1}$ which is more delicate.  We shall  first estimate in $H^s_\varphi$ the term $\mathtt{I}_m\triangleq \big\langle g,\big((\Psi^*)^{-1}-\textnormal{Id}\big)g_m \big\rangle_{L^2_\theta(\T)}$ that appears in the identity \eqref{Decomp-ZZ1}. Using the law products in  Lemma \ref{Law-prodX1} gives
\begin{align*}
\big\|\mathtt{I}_m\big\|_{H^{s}_\varphi}^{q,\kappa}&\lesssim \bigintssss_{\T}\|g(\cdot,\theta)\|_{H^{s}_\varphi}^{q,\kappa}\big\|\big((\Psi^\star)^{-1}-\textnormal{Id}\big)g_m(\cdot,\theta)\big\|_{H^{s_0}_{\varphi}}^{q,\kappa}d\theta\\
&+\bigintssss_{\T}\|g(\cdot,\theta)\|_{H^{s_0}_\varphi}^{q,\kappa}\big\|\big((\Psi^\star)^{-1}-\textnormal{Id}\big)g_m(\cdot,\theta)\big\|_{H^{s}_{\varphi}}^{q,\kappa}d\theta.
\end{align*}
According to    Cauchy-Schwarz inequality we find
\begin{align}\label{ES-Im}
\nonumber \big\|\mathtt{I}_m\big\|_{H^{s}_\varphi}^{q,\kappa}&\lesssim \| g\|_{L^2_\theta H^{s}_\varphi}^{q,\kappa}\big\|\big((\Psi^\star)^{-1}-\textnormal{Id}\big)g_m\big\|_{L^2_\theta H^{s_0}_{\varphi}}^{q,\kappa}+\| g\|_{L^2_\theta H^{s_0}_\varphi}^{q,\kappa}\big\|\big((\Psi^\star)^{-1}-\textnormal{Id}\big)g_m\big\|_{L^2_\theta H^{s}_{\varphi}}^{q,\kappa}\\
&\lesssim \| g\|_{s}^{q,\kappa}\big\|\big((\Psi^*)^{-1}-\textnormal{Id}\big)g_m\big\|_{{s_0}}^{q,\kappa}+\| g\|_{s_0}^{q,\kappa}\big\|\big((\Psi^\star)^{-1}-\textnormal{Id}\big)g_m\big\|_{s}^{q,\kappa}.
\end{align}
Thus using Proposition \ref{prop-constant-coe}-(iii)
 we get under the smallness condition \eqref{small-CC2}
\begin{align}\label{Yu-Tu1}
\| ((\Psi^\star)^{-1}-\textnormal{Id})g_m\|_{s}^{q,\kappa}&\lesssim \varepsilon\kappa^{-2}
\|g_m\|_{s+1}^{q,\kappa}+\varepsilon\kappa^{-2}\| \mathfrak{I}_{0}\|_{s+\sigma}^{q,\kappa}
\| g_m\|_{s_0+q+1}^{q,\kappa}.
\end{align}
At this level we need to establish an estimate for $\|g_m\|_{s}^{q,\kappa}$. From \eqref{Mer-09} we infer
\begin{align}\label{Mer-009}
\| g_m\|_{s}^{q,\kappa}\lesssim& \sum_{k\in\mathbb{S}_0}\|{\alpha_{m,k}}\|_{H^s_\varphi}^{q,\kappa}.
\end{align}
Recall that the coefficients $\alpha_{m,k}$ are the entries of the matrix ${\bf{M}}^{-1},$ where ${\bf{M}}=\big(c_{m,k}\big)_{m,k\in\mathbb{S}_0}$ is defined \mbox{in \eqref{M-atrix}.}
Then it follows that
\begin{align}\label{sig-stt}
c_{m,k}(\varphi)-\delta_{km}=\big\langle {\bf{e}}_k, \big((\Psi^\star)^{-1}-\textnormal{Id}\big){\bf{e}}_m\big\rangle_{L^2_\theta(\T)},
\end{align}
where $\delta_{km}$ is the Kronecker symbol. Consequently, we obtain
\begin{align*}
\nonumber\|c_{m,k}-\delta_{km}\|_{H^s_\varphi}^{q,\kappa}&\leqslant  \big\|(\Psi^\star)^{-1}-\textnormal{Id}\big){e}_m\big\|_{L^2_\theta H^s_\varphi }^{q,\kappa}\\
 &\lesssim \sup_{m\in\mathbb{S}_0}\big\|(\Psi^\star)^{-1}-\textnormal{Id}\big){e}_m\big\|_{s}^{q,\kappa}.
\end{align*}
Applying \eqref{Yu-Tu1} by replacing  $g_m$ with $e_m$ gives for any $s\in[ s_0,S]$,
\begin{align}\label{refm}
\max_{k,m\in\mathbb{S}_0}\|c_{m,k}-\delta_{km}\|_{H^s_\varphi}^{q,\kappa}&\lesssim \,\varepsilon\kappa^{-2}\left(1+\| \mathfrak{I}_{0}\|_{s+\sigma}^{q,\kappa}\right).
\end{align}
Finally, we deduce that
\begin{align}\label{vdo1}
{\bf{M}}(\lambda,\varphi)=\textnormal{Id}-{\bf R}(\lambda,\varphi)\quad\hbox{with} \quad \|{\bf R}\|_{H^s_\varphi}^{q,\kappa}\lesssim \varepsilon\kappa^{-2}\left(1+\| \mathfrak{I}_{0}\|_{s+\sigma}^{q,\kappa}\right).
\end{align}
Hence under the smallness condition \eqref{small-CC2} combined with the law products in  Lemma \ref{Law-prodX1}   we get 
\begin{align*}
\nonumber\|{\bf M}^{-1}-\textnormal{Id}\|_{H^s_\varphi}^{q,\kappa}&\leqslant \sum_{n\geqslant 1}\| {\bf R}^n\|_{s}^{q,\kappa}
\\
&\leqslant C\varepsilon\kappa^{-2}\left(1+\| \mathfrak{I}_{0}\|_{s+\sigma}^{q,\kappa}\right)
.
\end{align*}
Therefore  we deduce that
\begin{align}\label{H-D-P0}
\max_{m.k\in\mathbb{S}_0}\|\alpha_{m,k}\|_{H^s_\varphi}^{q,\kappa}&\leqslant C\big(1+\varepsilon\kappa^{-2}\| \mathfrak{I}_{0}\|_{s+\sigma}^{q,\kappa}\big).
\end{align}
Inserting \eqref{H-D-P0} into \eqref{Mer-009} allows to get 
\begin{align}\label{Mer-008}
\| g_m\|_{s}^{q,\kappa}\leqslant C\big(1+\varepsilon\kappa^{-2}\| \mathfrak{I}_{0}\|_{s+\sigma}^{q,\kappa}\big).
\end{align}
Putting together \eqref{Mer-008}, \eqref{Yu-Tu1} and \eqref{small-CC2} we find
\begin{align}\label{Yu-Tu2}
\| ((\Psi^\star)^{-1}-\textnormal{Id})g_m\|_{s}^{q,\kappa}\leqslant C\big(1+\varepsilon\kappa^{-2}\| \mathfrak{I}_{0}\|_{s+\sigma}^{q,\kappa}\big).
\end{align}
Plugging \eqref{Yu-Tu2} into \eqref{ES-Im} and using \eqref{small-CC2} yield
\begin{align}\label{ES-Im1}
 \big\|\mathtt{I}_m\big\|_{H^{s}_\varphi}^{q,\kappa}&\lesssim \varepsilon\kappa^{-2}\| g\|_{s}^{q,\kappa}+\varepsilon\kappa^{-2}\| \mathfrak{I}_{0}\|_{s+\sigma}^{q,\kappa}\| g\|_{s_0}^{q,\kappa}.
  \end{align}
Inserting this into \eqref{Decomp-ZZ1} and using the law products in Lemma \ref{Law-prodX1} we obtain
\begin{align}\label{Decomp-ZZe}
\|\Psi_{\perp}^{-1}g\|_{s}^{q,\kappa}&\lesssim \|\Psi^{-1}g\|_{s}^{q,\kappa}+\sum_{m\in\mathbb{S}_0}\big(\|\mathtt{I}_m\|_{s}^{q,\kappa}\|\Psi^{-1}{\bf e}_m\|_{s_0}^{q,\kappa}+\big\|\mathtt{I}_m\|_{s_0}^{q,\kappa}\|\Psi^{-1}{\bf e}_m\|_{s}^{q,\kappa}\big).
\end{align}
Thus using Proposition \ref{prop-constant-coe}-(iii) and \eqref{small-CC2}, we get successively
$$
\|\Psi^{\pm1}g\|_{s}^{q,\kappa}\lesssim \|g\|_{s}^{q,\kappa}+\varepsilon\kappa^{-2}\| \mathfrak{I}_{0}\|_{s+\sigma}^{q,\kappa}\| g\|_{s_0}^{q,\kappa}.
$$
and 
\begin{align}\label{tspr1}
\nonumber\max_{m\in\mathbb{S}_0}\|\Psi^{-1}{\bf e}_m\|_{s}^{q,\kappa}&\lesssim \|{\bf e}_m\|_{s}^{q,\kappa}+\varepsilon\kappa^{-2}\| \mathfrak{I}_{0}\|_{s+\sigma}^{q,\kappa}\| {\bf e}_m\|_{s_0}^{q,\kappa}\\
&\lesssim 1+\varepsilon\kappa^{-2}\| \mathfrak{I}_{0}\|_{s+\sigma}^{q,\kappa}.
\end{align}
Combining the preceding two estimates  with \eqref{Decomp-ZZe}, \eqref{ES-Im1} and \eqref{small-CC2} we deduce for $s\in[ s_0,S]$
\begin{align*}
\|\Psi_{\perp}^{-1}g\|_{s}^{q,\kappa}&\lesssim \| g\|_{s}^{q,\kappa}+\varepsilon\kappa^{-2}\| \mathfrak{I}_{0}\|_{s+\sigma}^{q,\kappa}\| g\|_{s_0}^{q,\kappa}.
\end{align*}
This completes  the proof of the desired estimate.  It remains to check the estimate of $\Delta_{12}g_m$. For this purpose we come back to 
\eqref{sig-stt} which gives 
\begin{align}\label{sig-stt1}
\Delta_{12}c_{m,k}(\varphi)=\big\langle {\bf{e}}_k, \Delta_{12}(\Psi^\star)^{-1}{\bf{e}}_m\big\rangle_{L^2_\theta(\T)}.
\end{align}
Thus we get from  Proposition \eqref{prop-constant-coe}-(iii)
\begin{align*}
\sup_{m,k\in\mathbb{S}_0}\|\Delta_{12}c_{m,k}\|_{H^{\overline{s}_h}_\varphi}^{q,\kappa}&\lesssim \varepsilon{\kappa^{-2}}\|\Delta_{12}i\|_{\overline{s}_h+\sigma}^{q,\kappa}.
\end{align*}
This implies from \eqref{vdo1}
\begin{align*}
\|\Delta_{12}{\bf R}\|_{H^{\overline{s}_h}_\varphi}^{q,\kappa}&\lesssim \varepsilon{\kappa^{-2}}\|\Delta_{12}i\|_{\overline{s}_h+\sigma}^{q,\kappa}.
\end{align*}
Now we remark that 
$$
\Delta_{12}{\bf M}^{-1}=\sum_{n\geqslant 1}\Delta_{12} {\bf R}^n.
$$
Therefore, the law products in  Lemma \ref{Law-prodX1}  combined with \eqref{vdo1} and \eqref{small-CC2} allow to get
\begin{align*}
 \|\Delta_{12}{\bf M}^{-1}\|_{H^{\overline{s}_h}_\varphi}^{q,\kappa}&\lesssim\|\Delta_{12}{\bf R}\|_{H^{\overline{s}_h}_\varphi}^{q,\kappa}\\
 &\lesssim \varepsilon{\kappa^{-2}}\|\Delta_{12}i\|_{\overline{s}_h+\sigma}^{q,\kappa}.
 \end{align*}
 From this we deduce that
\begin{align*}
\max_{m,k\in\mathbb{S}_0}\|\Delta_{12}\alpha_{m,k}\|_{H^{\overline{s}_h}_\varphi}^{q,\kappa}&\leqslant \varepsilon {\kappa^{-2}}\|\Delta_{12}i\|_{\overline{s}_h+\sigma}^{q,\kappa}.
\end{align*}
Hence  combining this estimate with \eqref{Mer-09} yields
\begin{align*}
\max_{m\in\mathbb{S}_0}\| \Delta_{12}g_m\|_{\overline{s}_h}^{q,\kappa}\leqslant \varepsilon {\kappa^{-2}}\|\Delta_{12}i\|_{\overline{s}_h+\sigma}^{q,\kappa}.
\end{align*}
This achieves the proof of Lemma \ref{lemm-decompY}.
\end{proof}
\subsubsection{Localization on the normal direction}
The main concern of this section is to describe the localization effects on the normal direction of the operator  $\widehat{\mathcal{L}}_{\omega}$ described in \eqref{Norm-local-z}. Our main result reads as follows.
\begin{proposition}\label{projection in the normal directions}
With the same notations and assumptions of Propositions $\ref{lemma-GS0}-\ref{QP-change}-\ref{prop-chang}-\ref{prop-constant-coe}$, there exists $\sigma_3=\sigma(\tau_1,d,q)\geqslant \sigma_2$ such that the following results hold true.
 \begin{enumerate}
\item  On the Cantor set $\mathcal{O}_{\infty,n}^{\kappa,\tau_{1}}(i_{0})$ introduced in Proposition $\ref{QP-change},$  we have
\begin{align*}
{\mathcal{L}_{\omega}^2}\triangleq\Psi_{\perp}^{-1}\widehat{\mathcal{L}}_{\omega}\Psi_{\perp}&=\big(\omega\cdot\partial_\varphi+{c}_{r,\lambda}\partial_\theta-\textnormal{m}_{r,\lambda}\partial_\theta|\textnormal{D}|^{\alpha-1}\big)\Pi_{\mathbb{S}_0}^{\perp}+\mathcal{R}^2_{r,\lambda}+\mathtt{E}_n^2\\
\nonumber&\triangleq\big(\omega\cdot\partial_{\varphi}+\mathscr{D}_{0}\big)\Pi_{\mathbb{S}_0}^{\perp}+\mathcal{R}^2_{r,\lambda}+\mathtt{E}_n^2
\end{align*}
where  $\mathcal{R}^2_{r,\lambda}=\Pi_{\mathbb{S}_0}^\perp \mathcal{R}^2_{r,\lambda}\Pi_{\mathbb{S}_0}^\perp$ is reversible  and  $\mathscr{D}_{0}$ is   reversible Fourier multiplier with 
$$
\forall (l,j)\in \mathbb{Z}^{d}\times\mathbb{S}_{0}^{c},\quad  \mathscr{D}_{0}\mathbf{e}_{l,j}=\ii\,\mu_{j}^{0}(\lambda,i_{0})\mathbf{e}_{l,j},$$
 with 
$$\mu_{j}^{0}(\lambda,i_{0})\triangleq \Omega_{j}(\alpha)+jr^{1}(\lambda,i_{0})-{
\tfrac{j\, \Gamma(j+\frac{\alpha}{2})}{\Gamma(1+j-\frac{\alpha}{2})}r^{2}(\lambda,i_{0})}.
$$
In addition, for $k\in\{1,2\}$
\begin{equation*}
\|r^k\|^{q,\kappa}\lesssim {\varepsilon \kappa^{-1}\left(1+\| \mathfrak{I}_{0}\|_{\overline{s}_h}^{q,\kappa}\right)},\quad \|{\Delta_{12}r^{k}}\|^{q,\kappa}\lesssim {\varepsilon \kappa^{-1}}\| \Delta_{12}i\|_{\overline{s}_h+\sigma_3}^{q,\kappa}.
\end{equation*}

\item The operator $\mathtt{E}_n^2$ is real reversible and satisfies the estimate
$$\big\|\mathtt{E}_n^{2}h\big\|_{s_0}^{q,\kappa}\lesssim\varepsilon{\kappa^{-2}}\Big(N_{0}^{\mu_{2}}N_{n}^{-\mu_{2}}+  N_{n}^{s-s_0}\big(1+\| \mathfrak{I}_{0}\|_{s+\sigma_3}^{q,\kappa}\big)\Big)\| h\|_{s_0+3}^{q,\kappa}.$$
\item The operator   $\mathcal{R}^2_{r,\lambda}$ is a real and reversible T\"oplitz in time operator satisfying 
\begin{equation*}
\forall\,\, s\in[s_0,S],\mbox{ }\interleave\mathcal{R}^2_{r,\lambda}\interleave_{2\overline\alpha-1+\epsilon,s,0}^{q,\kappa}\lesssim\varepsilon{\kappa^{-2}}\left(1+\| \mathfrak{I}_{0}\|_{s+\sigma_3}^{q,\kappa}\right)
\end{equation*}
and
\begin{equation*}
\interleave \Delta_{12}\mathcal{R}^2_{r,\lambda}\interleave_{0,0,0}^{q,\kappa}\lesssim \varepsilon{\kappa^{-2}}\|\Delta_{12}i\|_{\overline{s}_h+\sigma_3}^{q,\kappa}.
\end{equation*}

\end{enumerate}

\end{proposition}
\begin{proof}
{{\bf{(i)}} } According to \eqref{Lin-orthog} and  the decomposition $\textnormal{Id}=\Pi_{\mathbb{S}_{0}}+\Pi_{\mathbb{S}_{0}}^{\perp}$  we may write
\begin{align*}
\Psi_{\perp}^{-1} \widehat{\mathcal{L}}_{\omega}  \Psi_{\perp}&=\Psi_{\perp}^{-1}\mathcal{L}_{\varepsilon r,\lambda}\Psi_{\perp}-\varepsilon\Psi_{\perp}^{-1}\partial_{\theta}\mathcal{R}\Psi_{\perp}
\\
&=\Psi_{\perp}^{-1}\Pi_{\mathbb{S}_{0}}^{\perp}\mathcal{L}_{\varepsilon r,\lambda}\Psi\Pi_{\mathbb{S}_{0}}^{\perp}-\Psi_{\perp}^{-1}\Pi_{\mathbb{S}_{0}}^{\perp}\mathcal{L}_{\varepsilon r,\lambda}\Pi_{\mathbb{S}_{0}}\Psi\Pi_{\mathbb{S}_{0}}^{\perp}-\varepsilon\Psi_{\perp}^{-1}\Pi_{\mathbb{S}_{0}}^{\perp}\partial_{\theta}\mathcal{R}\Psi_{\perp}.
\end{align*}
Recall from Lemma \ref{lemma-reste} that 
\begin{align*}
 \mathcal{L}_{\varepsilon r,\lambda}&=\omega\cdot\partial_\varphi+\partial_\theta\Big[ V_{\varepsilon r,\alpha}-\big(\mathscr{W}_{\varepsilon r,\alpha}\, |{\textnormal D}|^{\alpha-1}+|{\textnormal D}|^{\alpha-1}\mathscr{W}_{\varepsilon r,\alpha}\big)+\mathscr{R}_{\varepsilon r,\lambda}\Big]\\
&\triangleq \omega\cdot\partial_\varphi+Y_{ r,\lambda}.
\end{align*}
Applying \eqref{dida-lip1} yields
$$
\mathcal{L}_{\varepsilon r,\lambda}\Psi=\Psi \mathcal{L}_{ r,\lambda}^1
$$
with
$$
 \mathcal{L}_{ r,\lambda}^1=\omega\cdot\partial_\varphi+{c}(\lambda,i_0)\partial_\theta- \textnormal{m}(\lambda,i_0)\partial_\theta|\textnormal{D}|^{\alpha-1}+\mathcal{R}^1_{ r,\lambda}+\mathtt{E}_n^1.
$$
Therefore, we find
\begin{align*}
\Psi_{\perp}^{-1}\widehat{\mathcal{L}}_{\omega}\Psi_{\perp}=&\Psi_{\perp}^{-1}\Pi_{\mathbb{S}_{0}}^{\perp}\Psi\mathcal{L}_{ r,\lambda}^1\Pi_{\mathbb{S}_{0}}^{\perp}-\Psi_{\perp}^{-1}\Pi_{\mathbb{S}_{0}}^{\perp}\left(Y_{ r,\lambda}-Y_{0,\lambda}\right)\Pi_{\mathbb{S}_{0}}\Psi\Pi_{\mathbb{S}_{0}}^{\perp}-\varepsilon\Psi_{\perp}^{-1}\partial_{\theta}\mathcal{R}\Psi_{\perp}
\end{align*}
where we have used the identities
$$
\Psi_{\perp}^{-1}\Pi_{\mathbb{S}_{0}}^{\perp}=\Psi_{\perp}^{-1}\quad\hbox{and}\quad [\Pi_{\mathbb{S}_{0}}^{\perp},T]=0=[\Pi_{\mathbb{S}_{0}},T], 
$$
 for any  Fourier multiplier $T$. In particular, we have  used tha $Y_{0,\lambda}$ is a Fourier multiplier and  thus
 $$
 \Pi_{\mathbb{S}_{0}}^{\perp} Y_{0,\lambda}\Pi_{\mathbb{S}_{0}}=0.
 $$ 
  From the structure of $\mathcal{L}_{ r,\lambda}^1$  stated before and the preceding commutation relations we infer \begin{align*}\Pi_{\mathbb{S}_{0}}^{\perp}\Psi \mathcal{L}_{ r,\lambda}^1\Pi_{\mathbb{S}_{0}}^{\perp}&=\Pi_{\mathbb{S}_{0}}^{\perp}\Psi\big(\omega\cdot\partial_\varphi+c(\lambda,i_0)\partial_\theta- \textnormal{m}(\lambda,i_0)\partial_\theta|\textnormal{D}|^{\alpha-1}+\mathcal{R}^1_{r,\lambda}\big)\Pi_{\mathbb{S}_{0}}^{\perp}+\Pi_{\mathbb{S}_{0}}^{\perp}\Psi \mathtt{E}_n^1\Pi_{\mathbb{S}_{0}}^{\perp}\\
&=\Psi_\perp\big(\omega\cdot\partial_\varphi+c(\lambda,i_0)\partial_\theta-\textnormal{m}(\lambda,i_0)\partial_\theta|\textnormal{D}|^{\alpha-1}\big)+\Pi_{\mathbb{S}_{0}}^{\perp}\Psi\mathcal{R}^1_{r,\lambda}\Pi_{\mathbb{S}_{0}}^{\perp}+\Pi_{\mathbb{S}_{0}}^{\perp}\Psi \mathtt{E}_n^1\Pi_{\mathbb{S}_{0}}^{\perp}.
 \end{align*}
It follows that
\begin{align*}
\Psi_{\perp}^{-1}\Pi_{\mathbb{S}_{0}}^{\perp}\Psi\mathcal{L}_{ r,\lambda}^1\Pi_{\mathbb{S}_{0}}^{\perp}
=&\,\omega\cdot\partial_\varphi+{c}(\lambda,i_0)\partial_\theta- \textnormal{m}(\lambda,i_0)\partial_\theta|\textnormal{D}|^{\alpha-1}+\Psi_{\perp}^{-1}\Pi_{\mathbb{S}_{0}}^{\perp}\Psi\mathcal{R}^1_{ r,\lambda}\Pi_{\mathbb{S}_{0}}^{\perp}+\Psi_{\perp}^{-1}\Pi_{\mathbb{S}_{0}}^{\perp}\Psi \mathtt{E}_n^1\Pi_{\mathbb{S}_{0}}^{\perp}\\
&=\omega\cdot\partial_\varphi+{c}(\lambda,i_0)\partial_\theta-\textnormal{m}(\lambda,i_0)\partial_\theta|\textnormal{D}|^{\alpha-1}+\Pi_{\mathbb{S}_{0}}^{\perp}\mathcal{R}^1_{ r,\lambda}\Pi_{\mathbb{S}_{0}}^{\perp}+\Psi_{\perp}^{-1}\Psi\Pi_{\mathbb{S}_{0}}\mathcal{R}^1_{ r,\lambda}\Pi_{\mathbb{S}_{0}}^{\perp}\\
&\quad+\Psi_{\perp}^{-1}\Pi_{\mathbb{S}_{0}}^{\perp}\Psi \mathtt{E}_n^1\Pi_{\mathbb{S}_{0}}^{\perp}.
 \end{align*}
 Consequently, we find
\begin{align}\label{Form1A}
\nonumber\Psi_{\perp}^{-1}\widehat{\mathcal{L}}_{\omega}\Psi_{\perp}=&\omega\cdot\partial_\varphi+{c}(\lambda,i_0)\partial_\theta-\textnormal{m}(\lambda,i_0)\partial_\theta|\textnormal{D}|^{\alpha-1}+\Pi_{\mathbb{S}_{0}}^{\perp}\mathcal{R}^1_{ r,\lambda}\Pi_{\mathbb{S}_{0}}^{\perp}+\Psi_{\perp}^{-1}\Psi\Pi_{\mathbb{S}_{0}}\mathcal{R}^1_{ r,\lambda}\Pi_{\mathbb{S}_{0}}^{\perp}\\
\nonumber&-\Psi_{\perp}^{-1}\Pi_{\mathbb{S}_{0}}^{\perp}\left(Y_{ r,\lambda}-Y_{0,\lambda}\right)\Pi_{\mathbb{S}_{0}}\Psi\Pi_{\mathbb{S}_{0}}^{\perp}-\varepsilon\Psi_{\perp}^{-1}\partial_{\theta}\mathcal{R}\Psi_{\perp}+\Psi_{\perp}^{-1}\Psi \mathtt{E}_n^1\Pi_{\mathbb{S}_{0}}^{\perp}\\
&\quad\triangleq\big(\omega\cdot\partial_{\varphi}+\mathscr{D}_{0}\big)\Pi_{\mathbb{S}_0}^{\perp}+\mathcal{R}^2_{r,\lambda}+\mathtt{E}_n^2,
\end{align}
with
$$
\mathtt{E}_n^2\triangleq \Psi_{\perp}^{-1}\Psi \mathtt{E}_n^1\Pi_{\mathbb{S}_{0}}^{\perp}\qquad\hbox{and}\qquad \mathscr{D}_{0}\triangleq {c}(\lambda,i_0)\partial_\theta-\textnormal{m}(\lambda,i_0)\partial_\theta|\textnormal{D}|^{\alpha-1} .
$$
Notice that $$
\forall (l,j)\in \mathbb{Z}^{d}\times\mathbb{S}_{0}^{c},\quad  \mathscr{D}_{0}\mathbf{e}_{l,j}=\ii\,\mu_{j}^{0}(\lambda,i_{0})\mathbf{e}_{l,j},$$
 with 
$$\mu_{j}^{0}(\lambda,i_{0})\triangleq \Omega_{j}(\alpha)+j\,r^{1}(\lambda,i_{0})-{
\tfrac{j\, \Gamma(j+\frac{\alpha}{2})}{\Gamma(1+j-\frac{\alpha}{2})}r^{2}(\lambda,i_{0})},
$$
where $W_{0,\alpha}$ is defined in Lemma \ref{lem-asym} and
$$
r^{1}(\lambda,i_{0})\triangleq {c}(\lambda,i_0)-V_{0,\alpha}\quad\hbox{and}\quad r^{2}(\lambda,i_{0})\triangleq \frac{\textnormal{m}(\lambda,i_0)-2^{-\alpha}C_\alpha}{W_{0,\alpha}}.
$$
Hence, the estimates of $r^1$ and $r^2$ follow from Proposition \ref{prop-constant-coe}-(i).

\smallskip

{\bf{(ii)}}  To get the suitable  estimates for  $\mathtt{E}_n^2$ it suffices to combine  Lemma \ref{lemm-decompY}-(ii) with \mbox{Proposition \ref{prop-constant-coe}-(iii)-(iv)} and \eqref{small-CC2}. Then one gets first the estimates
 \begin{align*}
 \big\|\Psi_{\perp}^{-1}\Pi_{\mathbb{S}_{0}}^{\perp}\Psi \mathtt{E}_n^1\Pi_{\mathbb{S}_{0}}^{\perp}\big\|_{s_0}^{q,\kappa}\lesssim \big\|\mathtt{E}_n^1\Pi_{\mathbb{S}_{0}}^{\perp}h\big\|_{s_0}^{q,\kappa}
 \end{align*}
 and 
 \begin{equation*}
\big\|\mathtt{E}_n^1\Pi_{\mathbb{S}_{0}}^{\perp}h\big\|_{s_0}^{q,\kappa}
\lesssim\varepsilon\kappa^{-2}\Big(N_{0}^{\mu_{2}}N_{n}^{-\mu_{2}}+  N_{n}^{s-s_0}\big(1+\| \mathfrak{I}_{0}\|_{s+\sigma}^{q,\kappa}\big)\Big)\| h\|_{s_0+3}^{q,\kappa}.
\end{equation*}
\\
{\bf{(iii)}} As to the term $\mathcal{R}^2_{ r,\lambda}$ we shall prove the following estimate,
\begin{equation}\label{vpn-1}
\interleave\mathcal{R}^2_{ r,\lambda}\interleave_{2\overline\alpha-1+\epsilon,s}^{q,\kappa}\lesssim\varepsilon\kappa^{-2}\left(1+\| \mathfrak{I}_{0}\|_{s+\sigma}^{q,\kappa}\right).
\end{equation}
 $\blacktriangleright$ {\it Estimate of} $\Pi_{\mathbb{S}_{0}}^{\perp}\mathcal{R}^1_{ r,\lambda}\Pi_{\mathbb{S}_{0}}^{\perp}$. One gets  according to Proposition \ref{prop-constant-coe}-(ii) and the continuity of the projectors
\begin{align}\label{DT-URTT}
\nonumber\interleave\Pi_{\mathbb{S}_{0}}^{\perp}\mathcal{R}^1_{ r,\lambda}\Pi_{\mathbb{S}_{0}}^{\perp}\interleave_{2\overline\alpha-1+\epsilon,s,0}^{q,\kappa}&\lesssim\interleave\mathcal{R}^1_{ r,\lambda}\interleave_{2\overline\alpha-1+\epsilon,s,0}^{q,\kappa}\\
 &\lesssim\varepsilon\kappa^{-2}\left(1+\| \mathfrak{I}_{0}\|_{s+\sigma}^{q,\kappa}\right).
\end{align}
$\blacktriangleright$ {\it Estimate of} $\Psi_{\perp}^{-1}\Psi \Pi_{\mathbb{S}_{0}}\mathcal{R}^1_{\varepsilon r,\lambda}\Pi_{\mathbb{S}_{0}}^{\perp}$. Using the identity of Lemma \ref{lemm-decompY}-(i) yields 
\begin{align}\label{DT-U}
\Psi_{\perp}^{-1}\Psi\Pi_{\mathbb{S}_{0}}\mathcal{R}^1_{ r,\lambda}\Pi_{\mathbb{S}_{0}}^{\perp}&=\Pi_{\mathbb{S}_{0}}\mathcal{R}^1_{ r,\lambda}\Pi_{\mathbb{S}_{0}}^{\perp}-\mathcal{T}_0 \mathcal{S}_1
\end{align}
where
\begin{align*}
\mathcal{T}_0h=\sum_{m\in\mathbb{S}_0}\big\langle h, \big({(\Psi^\star)}^{-1}-\textnormal{Id}\big)g_m  \big\rangle_{L^2_\theta(\T)}\Psi^{-1}e_m\quad\hbox{and}\quad \mathcal{S}_1\triangleq \Psi\Pi_{\mathbb{S}_{0}}\mathcal{R}^1_{ r,\lambda}\Pi_{\mathbb{S}_{0}}^{\perp}.
\end{align*}
To estimate the first term,  we use the same arguments as for  \eqref{DT-URTT}
 \begin{align*}
\nonumber\interleave\Pi_{\mathbb{S}_{0}}\mathcal{R}^1_{r,\lambda}\Pi_{\mathbb{S}_{0}}^{\perp}\interleave_{2\overline\alpha-1+\epsilon,s,0}^{q,\kappa}&\lesssim\interleave\mathcal{R}^1_{ r,\lambda}\interleave_{2\overline\alpha-1+\epsilon,s,0}^{q,\kappa}\\
 &\lesssim\varepsilon\kappa^{-2}\left(1+\| \mathfrak{I}_{0}\|_{s+\sigma}^{q,\kappa}\right).
\end{align*}
As to the second term, we write 
\begin{align*}
\mathcal{T}_0\mathcal{S}_1h&=\sum_{m\in\mathbb{S}_0}\big\langle \mathcal{S}_1 h, \big({(\Psi^\star)}^{-1}-\textnormal{Id}\big)g_m  \big\rangle_{L^2_\theta(\T)}\Psi^{-1}e_m\\
&=\sum_{m\in\mathbb{S}_0}\big\langle  h, \mathcal{S}_1^\star\big({(\Psi^\star)}^{-1}-\textnormal{Id}\big)g_m  \big\rangle_{L^2_\theta(\T)}\Psi^{-1}e_m,
\end{align*}
where $\mathcal{S}_1^\star$ is the adjoint of $\mathcal{S}_1.$ Observe that this operator  is an integral operator taking the form 
\begin{align}\label{Tya-re1}
\mathcal{T}_0\mathcal{S}_1h(\varphi,\theta)&=\int_{\T} \mathcal{K}_1(\varphi,\theta,\eta)h(\varphi,\eta) d\eta
\end{align}
with
$$
\mathcal{K}_1(\varphi,\theta,\eta)\triangleq\sum_{m\in\mathbb{S}_0}\mathcal{S}_1^\star\big({(\Psi^\star)}^{-1}-\textnormal{Id}\big)g_m(\varphi,\eta)  \big(\Psi^{-1}e_m\big)(\varphi,\theta).
$$
Applying Lemma \ref{lemma-Sym-R}-(ii)-(b), there $(-\Delta_\eta)^{\frac12}$ can be replaced by $\partial_\eta$, combined with Lemma \ref{Law-prodX1}we get for any $s\geqslant s_0$
\begin{align}\label{Z01}
\interleave \mathcal{T}_0\mathcal{S}_1\interleave_{-1,s,0}^{q,\kappa}\lesssim& \bigintssss_{\T}\big\|\partial_\eta\mathcal{K}_1(\cdot,\centerdot,\centerdot+\eta)\big\|_{s}^{q,\kappa} d\eta\\
\nonumber&\lesssim \sum_{m\in\mathbb{S}_0}\big\|\mathcal{S}_1^\star\big((\Psi^\star)^{-1}-\textnormal{Id}\big)g_m\big\|_{s+1}^{q,\kappa}\big\|\Psi^{-1}e_m\big\|_{s_0}^{q,\kappa}\\
\nonumber&\quad+\sum_{m\in\mathbb{S}_0}\big\|\mathcal{S}_1^\star\big((\Psi^\star)^{-1}-\textnormal{Id}\big)g_m\big\|_{s_0+1}^{q,\kappa}\big\|\Psi^{-1}e_m\big\|_{s}^{q,\kappa}.
\end{align}
Since  
$$
\mathcal{S}_1^\star=-\Pi_{\mathbb{S}_{0}}^{\perp}({\mathcal R}_{r,\lambda}^1)^\star\Pi_{\mathbb{S}_{0}}\Psi^\star.
$$
This implies
\begin{align}\label{Iddd}
\mathcal{S}_1^\star\big((\Psi^\star)^{-1}-\textnormal{Id}\big)g_m=-\Pi_{\mathbb{S}_{0}}^{\perp}({\mathcal R}_{r,\lambda}^1)^\star\Pi_{\mathbb{S}_{0}}\big(\textnormal{Id}-\Psi^\star\big)g_m.
\end{align}
Therefore, by applying  Proposition \ref{prop-constant-coe}-(ii), one gets
\begin{align*}
\interleave{(\mathcal R}_{ r,\lambda}^1)^\star\interleave_{0,s,0}^{q,\kappa}\lesssim & {\interleave{(\mathcal R}_{ r,\lambda}^1)^\star\interleave_{2\overline\alpha-1+\epsilon,s,0}^{q,\kappa}}\\
&\quad\lesssim  \varepsilon\kappa^{-2}\left(1+\|\mathfrak{I}_{0}\|_{s+\sigma}^{q,\kappa}\right).
\end{align*}
Then combining this estimate with  Lemma \ref{Lem-Rgv1}-$($iii$)$ and using the continuity of the orthogonal projectors yield
\begin{align}\label{D-L-O}
\nonumber \|\mathcal{S}_1^\star\big((\Psi^\star)^{-1}-\textnormal{Id}\big)g_m\|_{s}^{q,\kappa}
 &\lesssim \varepsilon\kappa^{-2}\left(1+\|\mathfrak{I}_{0}\|_{s+\sigma}^{q,\kappa}\right)\|\big(\textnormal{Id}-\Psi^\star\big)g_m\|_{s_0}^{q,\kappa}\\
&+ \varepsilon\kappa^{-2}\left(1+\|\mathfrak{I}_{0}\|_{s_0+\sigma}^{q,\kappa}\right)\|\big(\textnormal{Id}-\Psi^\star\big)g_m\|_{s}^{q,\kappa}.
\end{align}
Putting together    Proposition \ref{prop-constant-coe},  the smallness condition \eqref{small-CC2} and \eqref{Mer-008} lead to 
\begin{align}\label{prod-pid1}
\nonumber\big\|\big((\Psi^{\star})^{-1}-\textnormal{Id}\big)g_m\big\|_{s}^{q,\kappa}&\lesssim\|g_m\|_{s}^{q,\kappa}+{\varepsilon\kappa ^{-2}}\| \mathfrak{I}_{0}\|_{s+\sigma}^{q,\kappa}\|g_m\|_{s_{0}}^{q,\kappa}\\
&\lesssim 1+{\varepsilon\kappa ^{-2}}\| \mathfrak{I}_{0}\|_{s+\sigma}^{q,\kappa}.
\end{align}
Inserting \eqref{tspr1} and \eqref{prod-pid1} into \eqref{D-L-O} and using  \eqref{small-CC2} allow to get
\begin{align}\label{mima-l1}
\|\mathcal{S}_1^\star\big((\Psi^{\star})^{-1}-\textnormal{Id}\big)g_m\|_{s}^{q,\kappa}&\lesssim \varepsilon\kappa^{-2}\left(1+\|\mathfrak{I}_{0}\|_{s+\sigma}^{q,\kappa}\right).
\end{align}
Plugging this estimate into \eqref{Z01} yields for $s\geqslant s_0$
\begin{align}\label{Z02}
\interleave \mathcal{T}_0\mathcal{S}_1\interleave_{-1,s,0}^{q,\kappa}&\lesssim \varepsilon\kappa ^{-2}\left(1+\|\mathfrak{I}_{0}\|_{s+\sigma}^{q,\kappa}\right).
\end{align}

$\blacktriangleright$ {\it  Estimate of $ \Psi_{\perp}^{-1}\Pi_{\mathbb{S}_{0}}^{\perp} \left(Y_{ r,\lambda}-Y_{0,\lambda}\right)\Pi_{\mathbb{S}_{0}}\Psi\Pi_{\mathbb{S}_{0}}^{\perp}$.} 
We first write,
$$
\Psi_{\perp}^{-1}\Pi_{\mathbb{S}_{0}}^{\perp}\left(Y_{ r,\lambda}-Y_{0,\lambda}\right)\Pi_{\mathbb{S}_{0}}\Psi\Pi_{\mathbb{S}_{0}}^{\perp}\triangleq \Psi_{\perp}^{-1}\mathcal{S}_2\Psi\Pi_{\mathbb{S}_{0}}^{\perp}
$$
with
\begin{align*}
\mathcal{S}_2&\triangleq \left(Y_{ r,\lambda}-Y_{0,\lambda}\right)\Pi_{\mathbb{S}_{0}}\\
&\triangleq \partial_\theta\Big({\overline{V}}_{\varepsilon r,\alpha}-\big(\overline{\mathscr{W}}_{\varepsilon r,\alpha}\, |{\textnormal D}|^{\alpha-1}+|{\textnormal D}|^{\alpha-1}\overline{\mathscr{W}}_{\varepsilon r,\alpha}\big)+\mathscr{R}_{\varepsilon r,\lambda}\Big)\Pi_{\mathbb{S}_{0}}
\end{align*}
and 
$$
{\overline{V}}_{\varepsilon r,\alpha}\triangleq {{V}}_{\varepsilon r,\alpha}-{{V}}_{0,\alpha}\quad\hbox{and}\quad \overline{\mathscr{W}}_{\varepsilon r,\lambda}\triangleq {\mathscr{W}}_{\varepsilon r,\lambda}-{\mathscr{W}}_{0,\alpha}.
$$
Recall from \eqref{K2} that
$$
\left(Y_{ r,\lambda}-Y_{0,\lambda}\right) h(\varphi,\theta)=\partial_\theta \big({\overline{V}}_{\varepsilon r,\alpha}(\varphi,\theta)h(\varphi,\theta)\big)-\partial_\theta\int_{\T}\mathcal{K}_2(\varphi,\theta,\eta) \rho(\varphi,\eta) d\eta,
$$
with
$$
\mathcal{K}_2(\varphi,\theta,\eta)\triangleq\frac{\overline{\mathscr{W}}_{\varepsilon r,\lambda}(\varphi,\theta)+\overline{\mathscr{W}}_{\varepsilon r,\lambda}(\varphi,\eta)}{|\sin(\frac{\theta-\eta}{2})|^\alpha}-\frac{\mathscr{A}_{\varepsilon r,\alpha}(\varphi,\theta,\eta)}{|\sin(\frac{\eta-\theta}{2})|^{\alpha-2}}\cdot
$$
Then from elementary computations we find
\begin{align}\label{huita1}
\mathcal{S}_2\rho(\varphi,\theta)=\partial_\theta  \int_{\T}\mathcal{K}_3(\varphi,\theta,\eta) \rho(\varphi,\eta) d\eta,
\end{align}
with 
\begin{align*}
\mathcal{K}_3(\varphi,\theta,\eta) \triangleq{\overline{V}}_{\varepsilon r,\alpha}(\varphi,\theta)\,D_{\mathbb{S}_0}(\theta-\eta)&-\int_{\T}\mathcal{K}_2(\varphi,\theta,\theta+\eta^\prime)D_{\mathbb{S}_0}(\theta+\eta^\prime-\eta)d\eta^\prime\\
D_{\mathbb{S}_0}(\theta)&=\sum_{n\in\mathbb{S}_0} e^{\ii  n\theta}.
\end{align*}
Consequently, we get from \eqref{mathscrB} and \eqref{mathscrB1}
\begin{align}\label{huita2-1}
(\mathscr{B}^{-1}\mathcal{S}_2\mathscr{B}) \rho(\varphi,\theta)=\partial_\theta  \int_{\T}\mathcal{K}_4(\varphi,\theta,\eta) \rho(\varphi,\eta) d\eta,
\end{align}
with 
\begin{align*}
\mathcal{K}_4(\varphi,\theta,\eta) \triangleq \mathcal{K}_3\big(\varphi,\theta+\widehat{\beta}(\varphi,\theta),\eta+\widehat{\beta}(\varphi,\eta)\big).
\end{align*}
Straightforward long computations similar to that used to get the estimates of \eqref{KZ02} and \eqref{K3bis0} yield for any $\gamma\in\N$
\begin{align}\label{KKK-1}
\nonumber \sup_{\eta\in\T}\big\|\partial_\eta^\gamma\big((\partial_\theta \mathcal{K}_4)(\cdot,\centerdot,\centerdot+\eta)\big)\big\|_{s}^{q,\kappa}&\lesssim \varepsilon \kappa^{-1} \|r\|_{s+\sigma+\gamma}^{q,\kappa}\\
&\lesssim  \varepsilon\kappa ^{-1}\left(1+\|\mathfrak{I}_{0}\|_{s+\sigma+\gamma}^{q,\kappa}\right)
.
\end{align}
Applying Lemma \ref{lemma-Sym-R}-$($ii$)$ combined with \eqref{KKK-1} we infer 
\begin{align}\label{KKK-2}
\interleave \mathscr{B}^{-1}\mathcal{S}_2\mathscr{B}\interleave_{-1,s,\gamma}^{q,\kappa}\lesssim \varepsilon \kappa^{-1}\left(1+\|\mathfrak{I}_{0}\|_{s+\sigma+\gamma}^{q,\kappa}\right).
\end{align}
Thus applying Theorem \ref{Prop-EgorV} combined with \eqref{KKK-2}, \eqref{Est-eta}  and  \eqref{small-CC2} yields 
\begin{align}\label{KKK-3}
\nonumber \interleave \Psi^{-1}\mathcal{S}_2\Psi\interleave_{-1,s,0}^{q,\kappa}&\lesssim \varepsilon \kappa^{-1}\left(1+\|\mathfrak{I}_{0}\|_{s+\sigma}^{q,\kappa}\right)+ \varepsilon\kappa^{-2} \|r\|_{s+\sigma}^{q,\kappa}\\
&\lesssim  \varepsilon\kappa ^{-2}\left(1+\|\mathfrak{I}_{0}\|_{s+\sigma}^{q,\kappa}\right).
\end{align}
Now, proceeding as in \eqref{DT-U} we obtain
\begin{align}\label{mod-T-V}
\Psi_{\perp}^{-1}\mathcal{S}_2\Psi\Pi_{\mathbb{S}_{0}}^{\perp}&=\Psi^{-1}\mathcal{S}_2\Psi\Pi_{\mathbb{S}_{0}}^{\perp}-\mathcal{T}_0\mathcal{S}_2\Psi\Pi_{\mathbb{S}_{0}}^{\perp}.
\end{align}
It follows that 
\begin{align}\label{WL1}
\interleave\Psi_{\perp}^{-1}\mathcal{S}_2\Psi\Pi_{\mathbb{S}_{0}}^{\perp}\interleave_{-1,s,0}^{q,\kappa}&\lesssim \varepsilon\kappa ^{-2}\left(1+\|\mathfrak{I}_{0}\|_{s+\sigma}^{q,\kappa}\right)
+\interleave\mathcal{T}_0\mathcal{S}_2\Psi\Pi_{\mathbb{S}_{0}}^{\perp}\interleave_{-1,s,0}^{q,\kappa}.
\end{align}
 The estimate of the second term in \eqref{WL1} is similar to \eqref{Z02} with slight modifications and one gets
 \begin{align}\label{Z05}
 \interleave\mathcal{T}_0\mathcal{S}_2\Psi\Pi_{\mathbb{S}_{0}}^{\perp}\interleave_{-1,s,0}^{q,\kappa} &\lesssim \varepsilon\kappa^{-2}\left(1+\|\mathfrak{I}_{0}\|_{s+\sigma}^{q,\kappa}\right).
\end{align}
Plugging \eqref{Z05} into \eqref{WL1} we find
  \begin{align*}
\interleave\Psi_{\perp}^{-1}\mathcal{S}_2\Psi\Pi_{\mathbb{S}_{0}}^{\perp}\interleave_{-1,s,0}^{q,\kappa}&\lesssim\varepsilon\kappa^{-2}\left(1+\|\mathfrak{I}_{0}\|_{s+\sigma}^{q,\kappa}\right).
\end{align*}
$\blacktriangleright$ {\it Estimate of} $\varepsilon\Psi_{\perp}^{-1}\partial_{\theta}\mathcal{R}\Psi_{\perp}.$ The estimate of this term is quite   similar to the estimate of $\Psi_{\perp}^{-1}\mathcal{S}_2\Psi$, since $\partial_{\theta}\mathcal{R}$ admits an inetgral representation as in  \eqref{huita1}. Then using   the estimates of the kernel of $\partial_{\theta}\mathcal{R}$ stated in Proposition \ref{lemma-GS0} we obtain 
 \begin{align}
 \nonumber\varepsilon\interleave \Psi_{\perp}^{-1}\partial_{\theta}\mathcal{R}\Psi_{\perp}\interleave_{-1,s,0}^{q,\kappa} &\lesssim \varepsilon\kappa^{-2}\left(1+\|\mathfrak{I}_{0}\|_{s+\sigma}^{q,\kappa}\right).
\end{align}
Putting together the preceding estimates gives \eqref{vpn-1}.

\smallskip

$\blacktriangleright$
{{Estimate of $\Delta_{12}\mathcal{R}^2_{r,\lambda}$}}. Recall from \eqref{Form1A} that the operator $\mathcal{R}^2_{r,\lambda}$ takes the form
\begin{align}\label{Rr-np}
\nonumber \mathcal{R}^2_{r,\lambda}=&\Pi_{\mathbb{S}_{0}}^{\perp}\mathcal{R}^1_{ r,\lambda}\Pi_{\mathbb{S}_{0}}^{\perp}+\Psi_{\perp}^{-1}\Psi\Pi_{\mathbb{S}_{0}}\mathcal{R}^1_{ r,\lambda}\Pi_{\mathbb{S}_{0}}^{\perp}\\
&-\Psi_{\perp}^{-1}\Pi_{\mathbb{S}_{0}}^{\perp}\left(Y_{ r,\lambda}-Y_{0,\lambda}\right)\Pi_{\mathbb{S}_{0}}\Psi\Pi_{\mathbb{S}_{0}}^{\perp}-\varepsilon\Psi_{\perp}^{-1}\partial_{\theta}\mathcal{R}\Psi_{\perp}.
\end{align}

 The estimates of the different parts are slightly long and tedious sharing lot of similarities with the preceding discussion. For this reason we shall  explain how to proceed for the first  two terms and the remaining ones  can be treated in a similar way. Concerning the first one we simply use
$$
 \Delta_{12}\Pi_{\mathbb{S}_{0}}^{\perp}\mathcal{R}^1_{ r,\lambda}\Pi_{\mathbb{S}_{0}}^{\perp}=\Pi_{\mathbb{S}_{0}}^{\perp} \Delta_{12}\mathcal{R}^1_{ r,\lambda}\Pi_{\mathbb{S}_{0}}^{\perp}.
$$
 Then by the continuity of the projectors combined with Proposition \ref{prop-constant-coe}-(ii)
 \begin{align}\label{psk1}
\interleave \Delta_{12}\Pi_{\mathbb{S}_{0}}^{\perp}\mathcal{R}^1_{ r,\lambda}\Pi_{\mathbb{S}_{0}}^{\perp}\interleave_{0,0,0 }^{q,\kappa}\lesssim
\varepsilon \kappa^{-2}\|\Delta_{12} i\|_{q,\overline{s}_h+\sigma}^{q,\kappa}.
\end{align}
 As to the second term of $\mathcal{R}^2_{r,\lambda}$, we use the following representation detailed in \eqref{DT-U} and \eqref{Tya-re1}
 \begin{align*}
\Psi_{\perp}^{-1}\Psi\Pi_{\mathbb{S}_{0}}\mathcal{R}^1_{ r,\lambda}\Pi_{\mathbb{S}_{0}}^{\perp}&=\Pi_{\mathbb{S}_{0}}\mathcal{R}^1_{ r,\lambda}\Pi_{\mathbb{S}_{0}}^{\perp}-\mathcal{T}_0 \mathcal{S}_1
\end{align*}
where
\begin{align*}
\mathcal{T}_0\mathcal{S}_1h(\varphi,\theta)&=\int_{\T} \mathcal{K}_1(\varphi,\theta,\eta)h(\varphi,\eta) d\eta
\end{align*}
and
$$
\mathcal{K}_1(\varphi,\theta,\eta)\triangleq\sum_{m\in\mathbb{S}_0}\mathcal{S}_1^\star\big({(\Psi^\star)}^{-1}-\textnormal{Id}\big)g_m(\varphi,\eta)  \big(\Psi^{-1}e_m\big)(\varphi,\theta).
$$
The estimate of the first term is similar to \eqref{psk1}.  For the operator $\mathcal{T}_0 \mathcal{S}_1$ we use
\begin{align}\label{Z0DD1}
\interleave \Delta_{12}\mathcal{T}_0\mathcal{S}_1\interleave_{0,s_0,0}^{q,\kappa}\lesssim& \bigintssss_{\T}\big\|\Delta_{12}\mathcal{K}_1(\cdot,\centerdot,\centerdot+\eta)\big\|_{s_0}^{q,\kappa} d\eta\triangleq {\bf{I}}.
\end{align}
From the law products we get 
\begin{align*}
{\bf{I}}\lesssim& \sum_{m\in\mathbb{S}_0}\big\| \Delta_{12}\mathcal{S}_1^\star\big((\Psi^\star)^{-1}-\textnormal{Id}\big)g_m\big\|_{s_0+1}^{q,\kappa}\big\|\Psi_1^{-1}\textbf{e}_m\big\|_{s_0}^{q,\kappa}\\
&+\sum_{m\in\mathbb{S}_0}\big\| \mathcal{S}_{1,2}^\star\big((\Psi_2^\star)^{-1}-\textnormal{Id}\big)g_m\big\|_{s_0+1}^{q,\kappa}\big\|\Delta_{12}\Psi^{-1}\textbf{e}_m\big\|_{s_0}^{q,\kappa}\\
&\quad\triangleq{\bf{I}}_1+{\bf{I}}_2
\end{align*}
where the notation $f_j$ with $j=1,2$ denotes the value of  $f$ at the state $i_j$. We shall start with the estimate  of ${\bf{I}}_2$. For this purpose, we use Proposition \ref{prop-constant-coe}-(iii) combined with Sobolev embeddings
\begin{align*}
\max_{m\in\mathbb{S}_0}\|\Delta_{12}\Psi^{-1} \textbf{e}_m\|_{{s}_0}^{q,\kappa}
&\lesssim \varepsilon{\kappa^{-2}}\|\Delta_{12}i\|_{\overline{s}_h+\sigma}^{q,\kappa}.
\end{align*}
Thus we obtain from \eqref{mima-l1} and \eqref{small-CC2}
\begin{align}\label{mima-l2}
{\bf{I}}_2&\lesssim\varepsilon{\kappa^{-2}}\|\Delta_{12}i\|_{\overline{s}_h+\sigma}^{q,\kappa}.
\end{align}
In a similar way, coming back to the identity \eqref{Iddd} and applying   Proposition \ref{prop-constant-coe}-(ii)-(iii) combined with \eqref{tspr1},  \eqref{prod-pid1} and \eqref{small-CC2}
\begin{align*}
{\bf{I}}_1&\lesssim \max_{m\in\mathbb{S}_0}\big\| \Delta_{12}({\mathcal R}_{r,\lambda}^{1,\star}\Pi_{\mathbb{S}_{0}}\big(\textnormal{Id}-\Psi^\star\big)g_m)\big\|_{s_0+1}^{q,\kappa}\\
&\lesssim\varepsilon{\kappa^{-2}}\|\Delta_{12}i\|_{\overline{s}_h+\sigma}^{q,\kappa}.
\end{align*}
Therefore inserting the preceding estimate and \eqref{mima-l2} into \eqref{Z0DD1} gives
\begin{align*}
\interleave \Delta_{12}\mathcal{T}_0\mathcal{S}_1\interleave_{0,s_0,0}^{q,\kappa}\lesssim& \varepsilon{\kappa^{-2}}\|\Delta_{12}i\|_{\overline{s}_h+\sigma}^{q,\kappa}.
\end{align*}
The remaining terms of $\mathcal{R}^2_{r,\lambda} $ in  \eqref{Rr-np} can be analyzed following the same arguments. Then the proof of Proposition \ref{projection in the normal directions} is achieved.%

\end{proof}

\subsubsection{KAM reduction of the remainder term}\label{section-KAM-Red}
According to Proposition \ref{projection in the normal directions}, the operator  $\mathcal{L}_{\omega}^2=\Psi_{\perp}^{-1}\widehat{\mathcal{L}}_{\omega}\Psi_{\perp}$, which is well defined in the whole set $\mathcal{O}$, decomposes in the Cantor set $\mathcal{O}_{\infty,n}^{\kappa,\tau_{1}}(i_{0})$ described in Proposition $\ref{QP-change},$ as
\begin{align}\label{def-L0}
{\mathcal{L}_{\omega}^2}&=\mathscr{L}_0+\mathtt{E}_n^2\\
&\triangleq\big(\omega\cdot\partial_{\varphi}+\mathscr{D}_{0}\big)\Pi_{\mathbb{S}_0}^{\perp}+\mathcal{R}^2_{r,\lambda}+\mathtt{E}_n^2.\nonumber\end{align}
This section is devoted to the reduction   of the operator  $\mathscr{L}_{0}$. It will be implemented in a classical way using the KAM reduction provided that the exterior parameters belong to a suitable Cantor set.  Our main result in this section reads as follows.
\begin{proposition}\label{reduction of the remainder term}
With the constraints \eqref{Conv-T2} and using the same notations  of the preceding Propositions $\ref{lemma-GS0}-\ref{QP-change}-\ref{prop-chang}-\ref{prop-constant-coe}-\ref{projection in the normal directions}$, the following results hold true. Let  
\begin{align}\label{cond-diman1}
\mu_{2}\geqslant\overline{\mu}_{2}+{2\tau_2(1+q)}\quad\hbox{and}\quad s_{h}{=}\tfrac{3}{2}\mu_{2}+s_l+{\tau_2(1+q)}+1\triangleq \tfrac{3}{2}\mu_{2}+\overline{s}_{l}+1, 
\end{align}
where $\overline{\mu}_{2} $ and $ s_l$ are defined in \eqref{Conv-T2N}.
 There exists ${\varepsilon}_0>0$ and $\sigma_4=\sigma(\tau_1,d,q)\geqslant \sigma_3$ such that if
\begin{align}
\label{small-C3}\|\mathfrak{I}_{0}\|_{q,s_{h}+\sigma_4}\leqslant1\quad\textnormal{and}\quad N_{0}^{\mu_{2}}\varepsilon{\kappa^{-3-q}}\leqslant{\varepsilon}_0,
\end{align}
then the following assertions hold true.
\begin{enumerate}
\item There exists a  family of reversibility preserving invertible linear bounded operator $\Phi_{\infty}:\mathcal{O}\to \mathcal{L}\big(H_{\mathbb{S}_{0}}^{\perp},H_{\mathbb{S}_{0}}^{\perp}\big)$ satisfying the estimates
\begin{equation}\label{estimate on Phiinfty and its inverse}
\forall s\in[s_{0},S],\quad \mbox{ }\|\Phi_{\infty}^{\pm 1}h\|_{s}^{q,\kappa}\leqslant C\|h\|_{s}^{q,\kappa}+C{\varepsilon\overline \gamma^{-3}}\| \mathfrak{I}_{0}\|_{s+\sigma_4}^{q,\kappa}\|h\|_{s_0}^{q,\kappa}.
\end{equation}
There exists a family of diagonal operators $\mathscr{L}_{\infty}(\lambda,i_{0})$ taking the form
\begin{align*}\mathscr{L}_{\infty}(\lambda,i_{0})&=\big(\omega\cdot\partial_{\varphi}+\mathscr{D}_{\infty}(\lambda,i_{0})\big)\Pi_{\mathbb{S}_0}^{\perp}
\end{align*}
where $\mathscr{D}_{\infty}(\lambda,i_{0})$ is a reversible Fourier multiplier  operator given  by,
$$
\forall (l,j)\in \mathbb{Z}^{d}\times\mathbb{S}_{0}^{c},\quad  \mathscr{D}_{\infty}(\lambda,i_{0})\mathbf{e}_{l,j}={\ii}\,\mu_{j}^{\infty}(\lambda,i_0)\mathbf{e}_{l,j},$$
with
$$\forall j\in\mathbb{S}_{0}^{c},\quad\mu_{j}^{\infty}(\lambda,i_0)=\mu_{j}^{0}(\lambda,i_0)+r_{j}^{\infty}(\lambda,i_0)$$
and for any $\epsilon>0$
\begin{align}\label{estimate rjinfty}
\sup_{j\in\mathbb{S}_{0}^{c}}\| r_{j}^{\infty}\|^{q,\kappa}\leqslant C \varepsilon\kappa^{-2}\quad\quad \mbox{ and }\quad\quad{\sup_{j\in\mathbb{S}_{0}^{c}}|j|^{1-\epsilon-2\overline\alpha}\| r_{j}^{\infty}\|^{q,\kappa}\leqslant C\varepsilon\kappa^{-2}
}
\end{align}
such that on the Cantor set
\begin{align*}
\mathcal{O}_{{\infty,n}}^{\kappa,\tau_1,\tau_{2}}(i_{0})\triangleq \bigg\{\lambda=(\omega,\alpha)\in\mathcal{O}_{{\infty,n}}^{\kappa,\tau_{1}}(i_{0});\;\, &\forall |l|\leqslant {N_n},\,j,j_{0}\in\mathbb{S}_{0}^{c},\,(l,j)\neq(0,j_{0}),\\
&\big|\omega\cdot l+\mu_{j}^{\infty}(\lambda,i_0)-\mu_{j_{0}}^{\infty}(\lambda,i_0)\big|>2\tfrac{\kappa\langle j-j_{0}\rangle }{\langle l\rangle^{\tau_{2}}}\bigg\}
\end{align*}
we have 
\begin{align*}
\Phi_{\infty}^{-1}(\lambda,\omega,i_0)\mathscr{L}_{0}\Phi_{\infty}(\lambda,)&=\mathscr{L}_{\infty}(\lambda,i_{0})+{\mathtt{E}^3_n(\lambda,i_0)},
\end{align*}
with ${\mathtt{E}^3_n(\lambda,i_0)}$ a linear operator satisfying 
\begin{equation}\label{Error-Est-2D}
\|{\mathtt{E}^3_n(\lambda,i_0)}h\|_{s_0}^{q,\kappa}\leqslant C\varepsilon\kappa^{-3}N_{0}^{{\mu}_{2}}{N_{n+1}^{-\mu_{2}}} \|h\|_{s_0+3}^{q,\kappa}.
\end{equation}
Notice that the Cantor set $\mathcal{O}_{\infty,n}^{\kappa,\tau_{1}}(i_{0})$ was introduced in Proposition $\ref{QP-change}$, the operator $\mathscr{L}_{0}$ and the frequencies $\{\mu_{j}^{0}(\lambda)\}$ were described in  Proposition $ \ref{projection in the normal directions}$ and \eqref{def-L0}.
\item Given two tori $i_{1}$ and $i_{2}$  satisfying \eqref{small-C3} (replacing $\mathfrak{I}_{0}$ by $\mathfrak{I}_{1}$ or $\mathfrak{I}_{2}$), then if 
$${s_h}\geqslant 2\overline{s}_h+8\tau_2(q+1)
$$
   we get
\begin{equation*}
\forall j\in \mathbb{S}_0^c,\quad  \|\Delta_{12}r_{j}^{\infty}\|^{q,\kappa}
\lesssim  \varepsilon {\kappa^{-2}}\big(\|\Delta_{12}i\|_{\overline{s}_h+\sigma_4}^{q,\kappa}\big)^{\frac12}
\end{equation*}
and
\begin{equation*}
 \forall j\in \mathbb{S}_0^c,\quad \|\Delta_{12}\mu_{j}^{\infty}\|^{q,\kappa} \lesssim \varepsilon {\kappa^{-2}}\big(\|\Delta_{12}i\|_{\overline{s}_h+\sigma_4}^{q,\kappa}\big)^{\frac12}|j|,
\end{equation*}
where $\overline{s}_h$ is defined in \eqref{Conv-T2N}.
\end{enumerate}

\end{proposition}
\begin{proof}
{\bf{(i)}} 
Set 
$$
\delta_{0}(s)\triangleq\kappa ^{-1}\interleave\mathscr{R}_{0}\interleave_{s}^{q,\kappa},
$$
where $\mathscr{R}_0\triangleq \mathcal{R}_{r,\lambda}^2$ being  the remainder introduced  in  Proposition  \ref{projection in the normal directions} and the norm $\interleave\cdot\interleave_{s}^{q,\kappa} $ is described  in \eqref{Top-NormX}. According to  Proposition  \ref{projection in the normal directions}-(iii) one gets 
\begin{align}\label{whab1}\delta_{0}(s)\leqslant C\varepsilon\kappa^{-3}\left(1+\|\mathfrak{I}_{0}\|_{s+\sigma}^{q,\kappa}\right).
\end{align}
Then  with the notation of \eqref{Conv-T2}, the smallness condition \eqref{small-C3} implies that 
\begin{align}\label{Conv-P3}
 \nonumber N_{0}^{\mu_{2}}\delta_{0}(s_{h}) &\leqslant  C N_{0}^{\mu_{2}}\varepsilon\kappa^{-3}\\
 &\leqslant {C}{\varepsilon}_0.
\end{align}
$\blacktriangleright$ \textbf{KAM step.} In view of \eqref{def-L0} we write
\begin{equation}\label{mouka1}
\mathscr{L}_0=\big(\omega\cdot\partial_{\varphi}+\mathscr{D}_0\big)\Pi_{\mathbb{S}_0}^{\perp}+\mathscr{R}_0,\,\quad  \mathscr{R}_0\triangleq  \mathcal{R}_{r,\lambda}^2
\end{equation}
 and  on the Cantor set $\mathcal{O}_{\infty,n}^{\gamma,\tau_1}(i_0)$ one has the structure 
\begin{equation*}
\mathcal{L}_{\omega}^2=\mathscr{L}_0+\mathtt{E}_n^2
\end{equation*}
 with $\mathscr{D}_0$ a diagonal T\"oplitz operator and $\mathscr{R}_0$  a real and reversible T\"oplitz in time operator  of zero order and satisfies $\Pi_{\mathbb{S}_0}^\perp \mathscr{R}_0\Pi_{\mathbb{S}_0}^\perp =\mathscr{R}_0.$   We intend to explain the KAM step that will be implemented later during the KAM scheme in order to reduce completely  $\mathscr{R}_0$ into a diagonal operator. The scheme can  cover  more general linear operators. Actually, assume that we have a linear operator $\mathscr{L}$ such that on some Cantor set  $\mathscr{O}$ one has 
$$
\mathscr{L}=\big(\omega\cdot\partial_{\varphi}+\mathscr{D}\big)\Pi_{\mathbb{S}_0}^{\perp}+\mathscr{R},
$$
where $\mathscr{D}$ is real and reversible diagonal T\"oplitz in time  operator, that is,
\begin{align}\label{TUma1}
\mathscr{D}\mathbf{e}_{l,j}=i\mu_j(\lambda) \,\mathbf{e}_{l,j}\quad\hbox{and}\quad \mu_{-j}(\lambda)=-\mu_{j}(\lambda).
\end{align}
We shall also assume that the operator  $\mathscr{R}$ is  real and reversible T\"oplitz in time operator  of zero order and satisfies $\Pi_{\mathbb{S}_0}^\perp \mathscr{R}\Pi_{\mathbb{S}_0}^\perp =\mathscr{R}.$
Take  a linear invertible transformation close to the identity $$\Phi=\Pi_{\mathbb{S}_0}^{\perp}+\Psi:\mathcal{O}\rightarrow\mathcal{L}(H_{\mathbb{S}_{0}}^{\perp})$$
with $\Psi$ being small and will depend on $\mathscr{R}$.
Then 
$$\begin{array}{rcl}
\Phi^{-1}\mathscr{L}\Phi & = & \Phi^{-1}\Big(\Phi\left(\omega\cdot\partial_{\varphi}+\mathscr{D}\right)\Pi_{\mathbb{S}_0}^{\perp}+\left[\omega\cdot\partial_{\varphi}\Pi_{\mathbb{S}_0}^{\perp}+\mathscr{D},\Psi\right]+\mathscr{R}+\mathscr{R}\Psi\Big)\\
& = & \big(\omega\cdot\partial_{\varphi}+\mathscr{D}\big)\Pi_{\mathbb{S}_0}^{\perp}+\Phi^{-1}\Big(\big[\big(\omega\cdot\partial_{\varphi}+\mathscr{D}\big)\Pi_{\mathbb{S}_0}^{\perp},\Psi\big]+P_{N}\mathscr{R}+P_{N}^{\perp}\mathscr{R}+\mathscr{R}\Psi\Big),
\end{array}$$
where the projector $P_{N}$ was defined in \eqref{definition of projections for operators}.
The basic idea consists in replacing the remainder $\mathscr{R}$ with another quadratic one  up to  a diagonal part and provided that the parameters $\lambda=(\omega,\alpha)$ belongs to a Cantor set connected with non-resonance conditions associated to the {\it homological equation}. Iterating this scheme  will generate new remainders  which  become smaller and smaller  up to  new contributions on  the diagonal part and with  more excision on the parameters. Then by passing   to the limit we expect to diagonalize completely the operators provided that the parameters belong to Cantor set limit.    Now the first step is to impose the following  homological equation,
\begin{equation}\label{equation Psi}
\big[\big(\omega\cdot\partial_{\varphi}+\mathscr{D}\big)\Pi_{\mathbb{S}_0}^{\perp},\Psi\big]+P_{N}\mathscr{R}=\lfloor P_{N}\mathscr{R}\rfloor
\end{equation}
where $\lfloor P_{N}\mathscr{R}\rfloor$ is the diagonal part of the operator $\mathscr{R}$. We point out that the notation   $\lfloor  \mathscr{R}\rfloor$ with a general  operator $ \mathscr{R}$ is defined as follows, {for all $(l_{0},j_{0})\in\mathbb{Z}^{d}\times\mathbb{S}_{0}^{c}$},
\begin{equation}\label{Dp1X}
 \mathscr{R}{\bf e}_{l_0,j_0}=\sum_{{(}l,j{)\in\mathbb{Z}^{d}\times\mathbb{S}_{0}^{c}}}\mathscr{R}_{l_0,j_0}^{l,j}{\bf e}_{l,j}\Longrightarrow  \lfloor\mathscr{R}\rfloor {\bf e}_{l_0,j_0}=\mathscr{R}_{l_0,j_0}^{l_0,j_0}{\bf e}_{l_0,j_0}=\big \langle \mathscr{R}{\bf e}_{\ell_0,j_0}, {\bf e}_{l_0,j_0}\big\rangle_{L^2(\T^{d+1})}\,{\bf e}_{l_0,j_0}.
\end{equation}
Recall that ${\bf e}_{l_0,j_0}(\varphi,\theta)=e^{\ii (l_0\cdot\varphi+j_0\theta)}$. 
 We  define the Fourier coefficients {expansion} of $\Psi $ in a standard way as follows
$$
\Psi {\bf e}_{l_0,j_0}=\sum_{{(}l,j{)\in\mathbb{Z}^{d}\times\mathbb{S}_{0}^{c}}}\Psi_{l_0,j_0}^{l,j}{\bf e}_{l,j}.
$$
After straightforward computations, using Fourier expansion we find 
$$
\big[\omega\cdot\partial_{\varphi}\Pi_{\mathbb{S}_0}^{\perp},\Psi\big]\mathbf{e}_{l_{0},j_{0}}=\ii\sum_{(l,j)\in\mathbb{Z}^{d }\times\mathbb{S}_{0}^{c}}\Psi_{l_{0},j_{0}}^{l,j}\,\,\omega\cdot(l-l_{0})\,\,\mathbf{e}_{l,j}
$$
and similarly, since $\mathscr{D}$ is diagonal,
$$
[\mathscr{D},\Psi]\mathbf{e}_{l_{0},j_{0}}=\ii\sum_{(l,j)\in\mathbb{Z}^{d }\times\mathbb{S}_{0}^{c}}\Psi_{l_{0},j_{0}}^{l,j}\left({\mu_{j}(\lambda)-\mu_{j_0}(\lambda)}\right)\mathbf{e}_{l,j}.
$$
By assumption, $\mathscr{R}$ is a real and reversible T\"oplitz in time operator. Thus its Fourier coefficients  satisfy by Proposition \ref{characterization of real operator by its Fourier coefficients},
\begin{equation}\label{coefficients of the remainder operator R}
\mathscr{R}_{l_{0},j_{0}}^{l,j}\triangleq \ii\,r_{j_{0}}^{j}(\lambda,l_{0}-l)\in \ii\,\mathbb{R}\quad\mbox{ and }\quad \mathscr{R}_{-l_{0},-j_{0}}^{-l,-j}=-\mathscr{R}_{l_{0},j_{0}}^{l,j}.
\end{equation}
Hence $\Psi$ is a solution of \eqref{equation Psi} if and only if 
$$
\Psi {\bf e}_{l_0,j_0}=\sum_{|\ell-\ell_0|\leqslant N\atop |j-j_{0}|\leqslant N}\Psi_{l_0,j_0}^{l,j}{\bf e}_{l,j}
$$
and
$$\Psi_{l_{0},j_{0}}^{l,j}\Big(\omega\cdot (l-l_{0})+{\mu_{j}(\lambda)-\mu_{j_0}(\lambda)}\Big)=\left\lbrace\begin{array}{ll}
-r_{j_{0}}^{j}(\lambda,l_{0}-l) & \mbox{if }(l,j)\neq(l_{0},j_{0})\\
0 & \mbox{if }(l,j)=(l_{0},j_{0}).
\end{array}\right.$$
This implies that $\Psi$ is a T\"oplitz in time operator with $\Psi_{l_{0},j_{0}}^{l,j}\triangleq \Psi_{j_{0}}^{j}(l_{0}-l)$. In addition for  $(l,j,j_{0})\in\mathbb{Z}^{d }\times(\mathbb{S}_{0}^{c})^{2}$ with $|l|,|j-j_{0}|\leqslant N,$ one gets
\begin{equation}\label{Psi-gh}
\Psi_{j_{0}}^{j}(\lambda,l)=\left\lbrace\begin{array}{ll}
\frac{-r_{j_{0}}^{j}(\lambda,l)}{\omega\cdot l+\mu_{j}(\lambda)-\mu_{j_{0}}(\lambda)} & \mbox{if }(l,j)\neq(0,j_{0})\\
0 & \mbox{if }(l,j)=(0,j_{0}),
\end{array}\right.
\end{equation}
provided that the denominator is not vanishing. Moreover since $\Pi_{\mathbb{S}_0}^\perp \mathscr{R}\Pi_{\mathbb{S}_0}^\perp =\mathscr{R},$ then one finds that
$$
\forall \,l\,\in\mathbb{Z}^d,\,\forall\, j\,\,\hbox{or}\,\,\,j_0\in\mathbb{S}_0,\quad r_{j_{0}}^{j}(\lambda,l)=0.
$$
Consequently one may impose
$$
\forall \,l\,\in\mathbb{Z}^d,\,\forall\, j\,\,\hbox{or}\,\,\,j_0\in\mathbb{S}_0,\quad \Psi_{j_{0}}^{j}(\lambda,l)=0,
$$
which implies in particular that   $\Pi_{\mathbb{S}_0}^\perp \Psi\Pi_{\mathbb{S}_0}^\perp=\Psi.$
Next,  to justify the expression given by \eqref{Psi-gh} we need to avoid resonances and   restrict the parameters to the following  set 
$$
{\mathscr{O}_{+}^{\gamma }=\bigcap_{\underset{|l|\leqslant N}{(l,j,j_{0})\in\mathbb{Z}^{d }\times({\mathbb{S}}_{0}^{c})^{2}}\atop (l,j)\neq(0,j_0)}\left\lbrace\lambda\in\mathscr{O}
, |\omega\cdot l+\mu_{j}(\lambda)-\mu_{j_{0}}(\lambda)|>\frac{\kappa\langle j-j_{0}\rangle }{\langle l\rangle^{\tau_{2}}}\right\rbrace}.
$$ 
Then with this restriction the identity  \eqref{Psi-gh} is well defined and to extend  $\Psi$ to the whole set of parameters $\mathcal{O}$ we shall   use the cut-off function $\chi$ of  \eqref{chi-def-1}   by setting 
\begin{equation}\label{Ext-psi-op}
\Psi_{j_{0}}^{j}(\lambda,l)=\left\lbrace\begin{array}{ll}
-\varrho_{j_{0}}^{j}(\lambda,l)\,\, r_{j_{0}}^{j}(\lambda,l),& \mbox{if }\quad (l,j)\neq(0,j_{0})\\
0, & \mbox{if }\quad (l,j)=(0,j_{0}),
\end{array}\right.
\end{equation}
with
\begin{align}\label{varr-d}
\varrho_{j_{0}}^{j}(\lambda,l)\triangleq \frac{\chi\left((\omega\cdot l+\mu_{j}(\lambda)-\mu_{j_{0}}(\lambda))(\kappa\langle j-j_{0}\rangle)^{-1}\langle l\rangle^{\tau_{2}}\right)}{\omega\cdot l+\mu_{j}(\lambda)-\mu_{j_{0}}(\lambda)}\cdot
\end{align}
For the sake of simplicity  we shall  keep the same notation and identify $\Psi$ to its extension. Notice that the extension \eqref{Ext-psi-op} is smooth and coincides with $\Psi$ on the set $\mathscr{O}_{+}^{\gamma }$. 
In particular, we get from \eqref{coefficients of the remainder operator R} and \eqref{Ext-psi-op}  that $\Psi_{j_{0}}^{j}(l)\in\mathbb{R}.$ Moreover, combining \eqref{varr-d} and   \eqref{TUma1} yields 
$$\Psi_{-j_{0}}^{-j}(-l)=\Psi_{j_{0}}^{j}(l).$$
This implies by  Proposition \ref{characterization of real operator by its Fourier coefficients} that $\Psi$ is a real and reversibility preserving operator.
 { Let us define,
\begin{equation}\label{PU-RT}
\mathscr{D}_{+}=\mathscr{D}+\lfloor P_{N}\mathscr{R}\rfloor,\quad \quad \mathscr{R}_{+}=\Phi^{-1}\big(-\Psi\,\lfloor P_{N}\mathscr{R}\rfloor +P_{N}^{\perp}\mathscr{R}+\mathscr{R}\Psi\big)
\end{equation}
and
$$
\mathscr{L}_{+}\triangleq\big(\omega\cdot\partial_{\varphi}+\mathscr{D}_{+}+\mathscr{R}_{+}\big)\Pi_{\mathbb{S}_0}^{\perp}.$$
Then on the Cantor set $\mathscr{O}_{+}^{\gamma }$, we have
$$
\mathscr{L}_{+}=\Phi^{-1}\mathscr{L}\Phi.
$$}
 Notice that these operators are defined in the ambient set $\mathcal{O}.$
{Next we intend to  estimate $\varrho_{j_{0}}^{j}$  introduced in \eqref{varr-d}, which can be written in the form
\begin{align}\label{varr-dX1}
\varrho_{j_{0}}^{j}(\lambda,l)= a \chi_1(a A_{j,j_0}^l(\lambda)\big), \quad A_{j,j_0}^l(\lambda)=\omega\cdot l+\mu_{j}(\lambda)-\mu_{j_{0}}(\lambda), \quad a= (\kappa\langle j-j_{0}\rangle )^{-1}\langle l\rangle^{\tau_{2}},
\end{align}
where the function  ${\chi}_1(x)\triangleq\frac{\chi(x)}{x}$  is $\mathcal{C}^\infty$ with bounded derivatives and $\chi_1(0)=0$. We shall impose the following condition 
\begin{align}\label{reg-G-10}\exists C>0,\,\, \forall\,j,j_{0}\in\mathbb{S}_0^c,\quad\max_{0\leqslant|\gamma| \leqslant q}\sup_{\lambda\in\mathcal{O}}\left|\partial_{\lambda}^{\gamma}\left(\mu_{j}(\lambda)-\mu_{j_{0}}(\lambda)\right)\right|\leqslant C\,|j-j_{0}|.
\end{align}
Hence we obtain
\begin{align}\label{reg-G-11}\exists C>0,\,\, \forall\, l\in\mathbb{Z}^{d}, j,j_{0}\in\mathbb{S}_0^c,\quad\max_{0\leqslant|\gamma| \leqslant q}\sup_{\lambda\in\mathcal{O}}\left|\partial_{\lambda}^{\gamma}A^l_{j,j_0}(\lambda)\right|\leqslant C\,\big(| l |+|j-j_{0}|\big).
\end{align}
Similarly to \eqref{Dent-oist},  we may use Lemma \ref{Compos-lemm-VM} in order to get 

\begin{align}\label{eg-k-1}
 \forall |\gamma|\leqslant q,\quad \sup_{\lambda\in\mathcal{O}}\big|\partial_{\lambda}^{\gamma}\varrho_{j_{0}}^{j}(\lambda,l)\big|\leqslant C\kappa^{-(|\gamma|+1)} \langle l,j-j_{0}\rangle^{\tau_2(1+|\gamma|)+|\gamma|}.
\end{align}
}
Using the characterization \eqref{Top-NormX} combined with \eqref{Ext-psi-op}, \eqref{eg-k-1}  and Leibniz formula  yield 
\begin{align}\label{link Psi and R}
\nonumber  \interleave \Psi\interleave_{s}^{q,\kappa}=&  \Big(\max_{|\gamma|\leqslant q}\sup_{\lambda\in\mathcal{O}}\sup_{n\in\Z}\sum_{|l|\leqslant N\atop |m|\leqslant N}\kappa^{2|\gamma|}\langle l,m\rangle^{2(s-|\gamma|)}\big|\partial_{\lambda}^{\gamma}\Psi_{n}^{n+m}(l,\lambda)\big|^2\Big)^{\frac12}\\
\nonumber &\leqslant C  \kappa ^{-1}\interleave  P_{N}\mathscr{R}\interleave_{s+\tau_{2} q+\tau_{2}}^{q,\kappa}\\
&\quad \leqslant C  \kappa ^{-1}\interleave  \mathscr{R}\interleave_{s+\tau_{2} q+\tau_{2}}^{q,\kappa}.
\end{align}
Assume that 
\begin{align}\label{assum-Zi1}
 C\kappa ^{-1}\interleave \mathscr{R}\interleave_{{s_0+\tau_{2} q+\tau_{2}}}^{q,\kappa}\leqslant \varepsilon_0,
\end{align}
then applying \eqref{link Psi and R} yields 
\begin{align}\label{ZZ-KL}
\nonumber \|\Psi\interleave_{s_0}^{q,\kappa}&\leqslant C \kappa ^{-1}\interleave  \mathscr{R}\interleave_{s_0+\tau_{2} q+\tau_{2}}^{q,\kappa}\\
&\leqslant \varepsilon_0.
\end{align}
{This implies that  $\Phi^{-1}$ is invertible for small $\varepsilon_0$ and 
$$\Phi^{-1}=\displaystyle\sum_{n=0}^{+\infty}(-1)^{n}\Psi^{n}\triangleq \textnormal{Id}+\Sigma.
$$
}
Applying the law products from Lemma \ref{comm-pseudo1}, Lemma \ref{Lem-Rgv1}-(iv), \eqref{link Psi and R} and \eqref{ZZ-KL} allow to get
\begin{equation}\label{QWX0}\begin{array}{rcl}
\interleave\Sigma\interleave_{s}^{q,\kappa} & \leqslant & \displaystyle\|\Psi\interleave_{s}^{q,\kappa}\left(1+\sum_{n=1}^{+\infty}\left(C\interleave\Psi\interleave_{s_0}^{q,\kappa}\right)^{n}\right)\\
&{\leqslant} & \displaystyle C\,\kappa ^{-1}N^{\tau_{2} q+\tau_{2}}\interleave\mathscr{R}\interleave_{s}^{q,\kappa}.
\end{array}
\end{equation}
 Consequently, we conclude under the smallness condition \eqref{assum-Zi1} that  $\Phi^{-1}$ is invertible with
 \begin{equation}\label{QWX1}
\interleave\Phi^{-1}-\hbox{Id}\interleave_{s}^{q,\kappa}  {\leqslant}  C\kappa ^{-1}N^{\tau_{2} q+\tau_{2}}\interleave\mathscr{R}\interleave_{s}^{q,\kappa}.
\end{equation}
Coming back to \eqref{PU-RT}, we may write
$$
\mathscr{R}_{+}=P_{N}^{\perp}\mathscr{R}+\Phi^{-1}\mathscr{R}\Psi-\Psi\,\lfloor P_{N}\mathscr{R}\rfloor +\Sigma\big(P_{N}^{\perp}\mathscr{R}-\Psi\,\lfloor P_{N}\mathscr{R}\rfloor \big).
$$
Hence, by virtue of Lemma \ref{Lem-Rgv1}-(iv) and \eqref{QWX1}, we deduce that
\begin{align}\label{k-12}
\nonumber \interleave\mathscr{R}_{+}\interleave_{s}^{q,\kappa} & \leqslant  \interleave P_{N}^{\perp}\mathscr{R}\interleave_{s}^{q,\kappa}+C\interleave \Sigma\interleave_{s}^{q,\kappa}\left(\interleave P_{N}^{\perp}\mathscr{R}\interleave_{s_0}^{q,\kappa}+\interleave\Psi\interleave_{s_0}^{q,\kappa}\interleave\mathscr{R}\interleave_{s_0}^{q,\kappa}\right)\\
&+C\left(1+\interleave\Sigma\interleave_{s_0}^{q,\kappa}\right)
\Big(\interleave\Psi\interleave_{s}^{q,\kappa}\interleave\mathscr{R}\interleave_{s_0}^{q,\kappa}+\interleave\Psi\interleave_{s_0}^{q,\kappa}\interleave\mathscr{R}\interleave_{s}^{q,\kappa}\Big).
\end{align}
Combining  Lemma \ref{Lem-Rgv1}-(iv), \eqref{link Psi and R},\eqref{ZZ-KL} and  \eqref{QWX1}, we obtain for all $S\geqslant \overline{s}\geqslant s\geqslant s_{0}$,
\begin{equation}\label{KAM step remainder term}
\interleave\mathscr{R}_{+}\interleave_{s}^{q,\kappa}\leqslant N^{s-\overline{s}}\interleave\mathscr{R}\interleave_{\overline s}^{q,\kappa}+C\kappa ^{-1}N^{\tau_{2} q+\tau_{2}}\interleave\mathscr{R}\interleave_{s_0}^{q,\kappa}\interleave \mathscr{R}\interleave_{s}^{q,\kappa}.
\end{equation}
Recall  that  $S$ is a fixed  arbitrary large number. 

\smallskip

$\blacktriangleright$ \textbf{Initialization} We shall check that the assumptions \eqref{reg-G-10} and \eqref{assum-Zi1} required along the KAM step to reach   the final form \eqref{KAM step remainder term} are satisfied with the operator $\mathscr{L}=\mathscr{L}_0$ in \eqref{mouka1}. Indeed, for \eqref{reg-G-10} we shall use Lemma \ref{lem-asym}-(iv) in order to get
\begin{align}\label{reg-G-1}\exists C>0,\,\, \forall(j,j_{0})\in\mathbb{Z}^{2},\quad\max_{0\leqslant|\beta| \leqslant q}\sup_{\alpha\in(0,\overline\alpha)}\left|\partial_{\alpha}^{\beta}\left(\Omega_{j}(\alpha)-\Omega_{j_{0}}(\alpha)\right)\right|\leqslant C\,|j-j_{0}|.
\end{align}
Then applying Proposition \ref{projection in the normal directions}-(i) and \eqref{small-C3} we find 
\begin{align}\label{reg-GT-10}
\exists C>0,\,\, \forall(j,j_{0})\in\mathbb{Z}^{2},\quad\max_{0\leqslant|\beta| \leqslant q}\sup_{\lambda\in\mathcal{O}}\left|\partial_{\lambda}^{\beta}\left(\mu_{j}^{0}(\lambda)-\mu_{j_{0}}^{0}(\lambda)\right)\right|&\leqslant C\,|j-j_{0}|\big(1+\varepsilon\kappa^{-q-1}\big)\\
\nonumber &\leqslant C\,|j-j_{0}|.
\end{align}
As to the second assumption \eqref{assum-Zi1}, we apply  Proposition \ref{projection in the normal directions}-(iii) with $ \mathscr{R}_{0}=\mathscr{R}_{r,\lambda}^2$ combined with \eqref{small-C3}  leading to 
\begin{align*}
\nonumber C\kappa ^{-1}\interleave \mathscr{R}_0\interleave_{{s_0+\tau_{2} q+\tau_{2}}}^{q,\kappa}\leqslant&C\kappa ^{-1}\interleave\mathcal{R}^2_{r,\lambda}\interleave_{2\overline\alpha-1+\epsilon,{s_0+\tau_{2} q+\tau_{2}},0}^{q,\kappa}\\
&\lesssim\varepsilon {\kappa^{-3}}\left(1+\| \mathfrak{I}_{0}\|_{s_h+\sigma}^{q,\kappa}\right)\\
&\quad \leqslant C \varepsilon_0.
\end{align*}
$\blacktriangleright$ \textbf{KAM iteration}. For given  $m\in\mathbb{N}$ assume that we have  a linear  operator 
\begin{align}\label{Op-Lm}
\mathscr{L}_{m}\triangleq\big(\omega\cdot\partial_{\varphi}+\mathscr{D}_{m}+\mathscr{R}_{m}\big)\Pi_{\mathbb{S}_0}^\perp
\end{align}
with $\mathscr{D}_m$ a diagonal real reversible  T\"oplitz operator and $\mathscr{R}_m$ is a  real and reversible T\"oplitz in time operator  of zero order and satisfies $\Pi_{\mathbb{S}_0}^\perp \mathscr{R}_m\Pi_{\mathbb{S}_0}^\perp =\mathscr{R}_m.$ We assume that both assumptions \eqref{reg-G-10} and \eqref{assum-Zi1} are satisfied for $\mathscr{D}_m$ and $\mathscr{R}_m$. Notice that for $m=0$ we take the operator $\mathscr{L}_0$ defined in \eqref{mouka1}.
{ Let  $\Phi_{m}=\hbox{Id}+\Psi_{m}$ be  a linear invertible operator  such that
\begin{align}\label{Op-Lm1}
\Phi_{m}^{-1}\mathscr{L}_{m}\Phi_{m}\triangleq\big(\omega\cdot\partial_{\varphi}+\mathscr{D}_{m+1}+\mathscr{R}_{m+1}\big)\Pi_{\mathbb{S}_0}^\perp,
\end{align}
}
with $\Psi_{m}$ satisfying the homological equation
$$
\big[\big(\omega\cdot\partial_{\varphi}+\mathscr{D}_{m}\big)\Pi_{\mathbb{S}_0}^\perp,\Psi_{m}\big]+P_{N_{m}}\mathscr{R}_{m}=\lfloor P_{N_{m}}\mathscr{R}_{m}\rfloor.
$$
Remind that $N_{m}$ was introduced in \eqref{definition of Nm}. The diagonal parts   $(\mathscr{D}_{m})_{m\in\mathbb{N}}$ and the remainders $(\mathscr{R}_{m})_{m\in\mathbb{N}}$ are defined similarly to \eqref{PU-RT}  by the recursive formulas,
\begin{align}\label{Deco-T1}
\mathscr{D}_{m+1}=\mathscr{D}_{m}+\lfloor P_{N_{m}}\mathscr{R}_{m}\rfloor\quad \mbox{ and }\quad \mathscr{R}_{m+1}=\Phi_{m}^{-1}\left(-\Psi_{m}\,\lfloor P_{N_{m}}\mathscr{R}_{m}\rfloor +P_{N_{m}}^{\perp}\mathscr{R}_{m}+\mathscr{R}_{m}\Psi_{m}\right).
\end{align}
Notice that $\mathscr{D}_{m}$ and $\lfloor P_{N_{m}}\mathscr{R}_{m}\rfloor$ are T\"oplitz  Fourier multiplier operators that can be identified to their spectra  $(\ii\,\mu_{j}^{m})_{j\in\mathbb{S}_{0}^{c}}$ and $(\ii\,r_{j}^{m})_{j\in\mathbb{S}_{0}^{c}}$ in the following sense 
\begin{align}\label{Spect-T1}
\forall (l,j)\in\mathbb{Z}^d\times\mathbb{S}_0^c,\quad \mathscr{D}_{m} {\bf e}_{l,j}=\ii \,\mu_{j}^{m}\,{\bf e}_{l,j}\quad\hbox{and}\quad \lfloor P_{N_{m}}\mathscr{R}_{m}\rfloor {\bf e}_{l,j}=\ii\,r_{j}^{m}\,{\bf e}_{l,j}.
\end{align}
By construction, we deduce that
\begin{align}\label{Spect-T2}
\mu_{j}^{m+1}=\mu_{j}^{m}+r_{j}^{m}.
\end{align}
Similarly to \eqref{Psi-gh} we obtain 
\begin{equation}\label{Psi-gh1}
(\Psi_{m})_{j_{0}}^{j}(\lambda,l)=\left\lbrace\begin{array}{ll}
\frac{-r_{j_{0},m}^{j}(\lambda,l)}{\omega\cdot l+\mu_{j}^{m}(\lambda)-\mu_{j_{0}}^{m}(\lambda)} & \mbox{if }(l,j)\neq(0,j_{0})\\
0 & \mbox{if }(l,j)=(0,j_{0}),
\end{array}\right.
\end{equation}
where  the sequence  $\big\{r_{j_{0},m}^{j}(\lambda,l)\big\}$ describes  the Fourier coefficients of $\mathscr{R}_{m}$, that is,
$$
\mathscr{R}_{m} {\bf e}_{l_0,j_0}=\ii\,\sum_{(l,j)\in\mathbb{Z}^{d+1}}{r_{j_{0},m}^{j}(\lambda, l_0-l){\bf e}_{l,j}}.
$$
We shall introduce the open Cantor set where the preceding formula has a meaning,
\begin{equation}\label{Cantor-SX}
\mathscr{O}_{m+1}^{\kappa }=\bigcap_{\underset{{(l,j)\neq(0,j_0)}}{\underset{|l|\leqslant N_{m}}{(l,j,j_{0})\in\mathbb{Z}^{d }\times(\mathbb{S}_{0}^{c})^{2}}}}\left\lbrace\lambda\in\mathscr{O}_{m}^{\kappa };\;\, |\omega\cdot l+\mu_{j}^{m}(\lambda)-\mu_{j_{0}}^{m}(\lambda)|>\kappa\tfrac{\langle j-j_{0}\rangle }{\langle l\rangle^{\tau_{2}}}\right\rbrace.
\end{equation}
 As in \eqref{Ext-psi-op} and \eqref{varr-d} we may extend \eqref{Psi-gh1} as follows 
$$
(\Psi_{m})_{j_{0}}^{j}(\lambda,\omega,l)=\left\lbrace\begin{array}{ll}
-\frac{\chi\left((\omega\cdot l+\mu_{j}^{m}(\lambda)-\mu_{j_{0}}^{m}(\lambda))(\kappa|j-j_{0}|)^{-1}\langle l\rangle^{\tau_{2}}\right){r}_{j_{0},m}^{j}(\lambda,l)}{\omega\cdot l+\mu_{j}^{m}(\lambda)-\mu_{j_{0}}^{m}(\lambda)} & \mbox{if }(l,j)\neq(0,j_{0})\\
0 & \mbox{if }(l,j)=(0,j_{0}).
\end{array}\right.
$$
{We emphasize that working with this extension for $\Psi_m$ allows to extend naturally both $\mathscr{D}_{m+1}$ and   the remainder $\mathscr{R}_{m+1}$ provided that the  operators $\mathscr{D}_{m}$ and $\mathscr{R}_{m}$ are defined in the ambient set of parameters $\mathcal{O}.$ Thus the  operator defined by   the right-hand side in \eqref{Op-Lm1} can be  extended to  the whole set of parameters $\mathcal{O}$ and for simplicity we  still denote this extension by $\mathscr{L}_{m+1}$, that is,
\begin{align}\label{Op-Lm2}
\big(\omega\cdot\partial_{\varphi}+\mathscr{D}_{m+1}+\mathscr{R}_{m+1}\big)\Pi_{\mathbb{S}_0}^\perp\triangleq\mathscr{L}_{m+1}.
\end{align}
This allows to construct by induction the sequence of operators $\left(\mathscr{L}_{m+1}\right)$ in the full set  $\mathcal{O}$. Similarly the operator $\Phi_{m}^{-1}\mathscr{L}_{m}\Phi_{m}$ admits an extension in $\mathcal{O}$ by simply extending $\Phi_m^{\pm1}$ . Nevertheless, by   construction the identity $\mathscr{L}_{m+1}=\Phi_{m}^{-1}\mathscr{L}_{m}\Phi_{m}$ in \eqref{Op-Lm1} occurs on the  Cantor set $\mathscr{O}_{m+1}^{\kappa }$ and may fail outside this set.
}We set 
\begin{equation}\label{Def-Taou}
\delta_{m}(s)\triangleq \kappa ^{-1}\interleave \mathscr{R}_{m}\interleave_{s}^{q,\kappa}
\end{equation}
and we intend to   prove by induction in $m\in\mathbb{N}$ that 
\begin{equation}\label{hypothesis of induction for deltamprime}
\forall\, m\in\mathbb{N},\, \forall s\in[s_{0},\overline{s}_{l}],\,\delta_{m}(s)\leqslant \delta_{0}(s_{h})N_{0}^{\mu_{2}}N_{m}^{-\mu_{2}}\quad \mbox{ and }\quad \delta_{m}(s_{h})\leqslant\left(2-\frac{1}{m+1}\right)\delta_{0}(s_{h}).
\end{equation}
In addition, we should check the validity of the  assumptions \eqref{reg-G-10} and \eqref{assum-Zi1}  for $\mathscr{D}_{m+1}$ and $\mathscr{R}_{m+1}$.
 Observe that by Sobolev embeddings, it is sufficient to prove the first inequality with $s=\overline{s}_{l}.$ The  property is trivial  for $m=0$. Now, assume that  the property \eqref{hypothesis of induction for deltamprime} is  true for $m\in\mathbb{N}$ and let us check it at the order $m+1.$
We write 
\begin{align}\label{phi-inverse1}
\Phi_{m}^{-1}=\hbox{Id}+\Sigma_{m}\quad\hbox{with}\quad \Sigma_{m}=\displaystyle\sum_{n=1}^{+\infty}(-1)^{n}\Psi_{m}^{n}.
\end{align} 
Then proceeding as for \eqref{QWX0}, using in particular \eqref{link Psi and R} and \eqref{ZZ-KL} allows to deduce successively
 \begin{align*}
\interleave\Sigma_{m}\interleave_{s_0}^{q,\kappa}  & \leqslant  \interleave\Psi_{m}\interleave_{s_0}^{q,\kappa}\left(1+\sum_{n=0}^{+\infty}\big(C\interleave\Psi_{m}\interleave_{s_0}^{q,\kappa}\big)^{n}\right)\\
& \leqslant  \displaystyle \delta_{m}(s_0+\tau_{2} q+\tau_{2})\left(1+\sum_{n=0}^{+\infty}\big(C\delta_{m}(s_{0}+\tau_{2} q+\tau_{2})\big)^{n}\right)
\end{align*}
and for any $s\in[s_0,S],$
$$\begin{array}{rcl}
\interleave\Sigma_{m}\interleave_{s}^{q,\kappa} & \leqslant & \displaystyle\|\Psi_{m}\interleave_{s}^{q,\kappa}\left(1+\sum_{n=0}^{+\infty}\big(C\interleave\Psi_{m}\interleave_{s_0}^{q,\kappa}\big)^{n}\right)\\
& \leqslant & \displaystyle N_{m}^{\tau_{2} q+\tau_{2}}\delta_{m}(s)\left(1+\sum_{n=0}^{+\infty}\big(C\delta_{m}(s_{0}+\tau_{2} q+\tau_{2})\big)^{n}\right).
\end{array}$$
Consequently, we get by the induction assumption, since ${s_{0}+\tau_{2} q+\tau_{2}\leqslant \overline{s}_l},\,N_{m}\geqslant N_{0}$,
$$\begin{array}{rcl}
\|\Sigma_{m}\interleave_{s_0}^{q,\kappa} & \leqslant & \displaystyle CN_{0}^{\mu_{2}}N_{m}^{-\mu_{2}}\delta_{0}(s_{h})\left(1+\sum_{n=0}^{+\infty}\big(CN_{0}^{\mu_{2}}N_{m}^{-\mu_{2}}\delta_{0}(s_{h})\big)^{n}\right)\\
& \leqslant & \displaystyle CN_{0}^{\mu_{2}}N_{m}^{-\mu_{2}}\delta_{0}(s_{h})\left(1+\sum_{n=0}^{+\infty}\big(C\delta_{0}(s_{h})\big)^{n}\right)
\end{array}$$
and
\begin{align*}
\|\Sigma_{m}\interleave_{s}^{q,\kappa}  &\leqslant  \displaystyle N_{m}^{\tau_{2} q+\tau_{2}}\delta_{m}(s)\left(1+\sum_{n=0}^{+\infty}\big(CN_{0}^{\mu_{2}}N_{m}^{-\mu_{2}}\delta_{0}(s_{h})\big)^{n}\right)\\
 &\leqslant  \displaystyle N_{m}^{\tau_{2} q+\tau_{2}}\delta_{m}(s)\left(1+\sum_{n=0}^{+\infty}\big(C\delta_{0}(s_{h})\big)^{n}\right).
\end{align*}
Thus under the smallness condition \eqref{Conv-P3} we obtain for any $s\in[s_0,S],$
\begin{align}\label{AP-W34}
\|\Sigma_{m}\interleave_{s_0}^{q,\kappa}\leqslant CN_{0}^{\mu_{2}}N_{m}^{-\mu_{2}}\delta_{0}(s_{h})\quad\hbox{and}\quad \|\Sigma_{m}\interleave_{s}^{q,\kappa} & \leqslant  CN_{m}^{\tau_{2} q+\tau_{2}}\delta_{m}(s).
\end{align}
Notice that one also may obtain 
\begin{align}\label{Bir-11}
\|\Sigma_{m}\interleave_{s}^{q,\kappa} & \leqslant  C\delta_{m}(s+\tau_2(1+q)).
\end{align}
Applying  KAM step \eqref{KAM step remainder term} and using Sobolev embeddings,  we deduce that
$$\begin{array}{rcl}
\delta_{m+1}(\overline{s}_{l}) & \leqslant & N_{m}^{\overline{s}_{l}-s_{h}}\delta_{m}(s_{h})+CN_{m}^{\tau_{2}q+\tau_{2}}\left(\delta_{m}(\overline{s}_{l})\right)^{2}.
\end{array}$$
Thus using the induction assumption \eqref{hypothesis of induction for deltamprime} yields
$$\begin{array}{rcl}
\delta_{m+1}(\overline{s}_{l}) & \leqslant & N_{m}^{\overline{s}_{l}-s_{h}}\left(2-\frac{1}{m+1}\right)\delta_{0}(s_{h})+CN_{m}^{\tau_{2}q+\tau_{2}}\delta_{0}^2(s_{h})N_{0}^{2\mu_{2}}N_{m}^{-2\mu_{2}}\\
& \leqslant & 2 N_{m}^{\overline{s}_{l}-s_{h}}\delta_{0}(s_{h})+CN_{m}^{\tau_{2}q+\tau_{2}}\delta_{0}^2(s_{h})N_{0}^{2\mu_{2}}N_{m}^{-2\mu_{2}}.
\end{array}$$
If we select the parameters $\overline{s}_{l},s_h$ and $\mu_2$ such that
\begin{equation}\label{Conv-T1}
N_{m}^{\overline{s}_{l}-s_{h}} \leqslant  \frac{1}{4}N_{0}^{\mu_{2}}N_{m+1}^{-\mu_{2}}\qquad\hbox{and}\qquad 
CN_{m}^{\tau_{2}q+\tau_{2}}\delta_{0}(s_{h})N_{0}^{2\mu_{2}}N_{m}^{-2\mu_{2}}  \leqslant  \frac{1}{2}N_{0}^{\mu_{2}}N_{m+1}^{-\mu_{2}}
\end{equation}
then we find
$$
\delta_{m+1}(\overline{s}_{l}) \leqslant \delta_{0}(s_{h})N_{0}^{\mu_{2}}N_{m+1}^{-\mu_{2}}.
$$
By elementary arguments based on \eqref{definition of Nm} and \eqref{cond-diman1} we can check   that the assumptions  of \eqref{Conv-T1} hold true provided that 
\begin{equation}\label{Cond1}
4N_{0}^{-\mu_2}\leqslant 1\qquad\hbox{and}\qquad 
2C N_{0}^{\mu_{2}}\delta_{0}(s_{h}) \leqslant  1.
\end{equation}
The second constraint  follows from  \eqref{Conv-P3}, however  the first one   is automatically satisfied due to  $N_0\geqslant 2$ and $\mu_2\geqslant 2,$ according to \eqref{Conv-T2}.
This achieves  the first statement of the induction in \eqref{hypothesis of induction for deltamprime}. Let us now move to the second estimate  in \eqref{hypothesis of induction for deltamprime}.
Applying  the KAM step \eqref{KAM step remainder term} combined with  the induction assumptions \eqref{hypothesis of induction for deltamprime}
\begin{align*}
\delta_{m+1}(s_{h}) & \leqslant  \delta_{m}(s_{h})+C N_{m}^{\tau_{2}q+\tau_{2}}\delta_{m}(s_{0}) \delta_{m}(s_{h})\\
&{\leqslant}  \Big(2-\frac{1}{m+1}\Big)\delta_{0}(s_{h})\Big(1+CN_{0}^{\mu_{2}}N_{m}^{\tau_{2}q+\tau_{2}-\mu_{2}}\delta_{0}(s_{h})\Big).
\end{align*}
Therefore if one has
\begin{equation}\label{Conv-4}
\Big(2-\frac{1}{m+1}\Big)\Big(1+CN_{0}^{\mu_{2}}N_{m}^{\tau_{2}q+\tau_{2}-\mu_{2}}\delta_{0}(s_{h})\Big)\leqslant 2-\frac{1}{m+2}
\end{equation}
then we find
$$
\delta_{m+1}(s_{h})  \leqslant  \left(2-\frac{1}{m+2}\right)\delta_{0}(s_{h})
$$
which achieves the induction argument of \eqref{hypothesis of induction for deltamprime}. Observe from \eqref{cond-diman1} that  $\mu_2\geqslant 2\tau_2(1+q)$, then the  condition  \eqref{Conv-4} holds true if  
\begin{equation}\label{Conv-5}
CN_{0}^{\mu_{2}}N_{m}^{-\tau_{2}q-\tau_{2}}\delta_{0}(s_{h})\leqslant\frac{1}{(2m+1)(m+2)}\cdot
\end{equation}
Then as $N_{0}\geqslant 2$ we may find a constant $c_0>0$ small enough  such that
$$
\forall\, m\in\mathbb{N},\quad c_0N_m^{-1}\leqslant \frac{1}{(2m+1)(m+2)}\cdot
$$
Therefore \eqref{Conv-5} is satisfied  provided that
\begin{equation}\label{Conv-6}
CN_{0}^{\mu_{2}}N_{m}^{-\tau_{2}q-\tau_{2}+1}\delta_{0}(s_{h})\leqslant c_0.
\end{equation}
From the assumption   \eqref{Conv-T2} we get in particular
$
\tau_{2}q+\tau_{2}-1\geqslant0.
$
Then \eqref{Conv-6} is satisfied according to \eqref{Conv-P3} when $\varepsilon_0$ is small enough. To achieve the induction   proof of \eqref{hypothesis of induction for deltamprime} it remains to check that the  assumptions \eqref{reg-G-10} and \eqref{assum-Zi1} are satisfied for $\mathscr{D}_{m+1}$ and $\mathscr{R}_{m+1}$. First, the assumption  \eqref{assum-Zi1} follows from the first inequality of  \eqref{hypothesis of induction for deltamprime}  applied at the order $m+1$  with $s=s_0+\tau_2(1+q)\leqslant \overline{s}_l$ supplemented with \eqref{Conv-P3}. Second,  to check the  validity of  \eqref{reg-G-10} for the eigenvalues of  $\mathscr{D}_{m+1}$, we combine \eqref{Spect-T1} and \eqref{Spect-T2} and  \eqref{Dp1X}, in order to find
\begin{align*}
\|\mu_{j}^{m+1}-\mu_{j}^{m}\|^{q,\kappa}
&=\big\|\big\langle  P_{N_{m}}\mathscr{R}_{m}\mathbf{e}_{l,j},\mathbf{e}_{l,j}\big\rangle_{L^{2}(\mathbb{T}^{d +1})}\big\|^{q,\kappa}.
\end{align*}
Since $\mathscr{R}_{m}$ is T\"oplitz then 
\begin{align*}
\|\mu_{j}^{m+1}-\mu_{j}^{m}\|^{q,\kappa}&=\big\|\big\langle P_{N_{m}}\mathscr{R}_{m}\mathbf{e}_{0,j},\mathbf{e}_{0,j}\big\rangle_{L^{2}(\mathbb{T}^{d +1})}\big\|^{q,\kappa}.
\end{align*}
A duality argument combined with  Lemma \ref{Lem-Rgv1}-(iii) and \eqref{Def-Taou} yield
\begin{align}\label{dual-1X}
\nonumber \|\mu_{j}^{m+1}-\mu_{j}^{m}\|^{q,\kappa}\lesssim& \big\|\mathscr{R}_{m}\mathbf{e}_{0,j}\|_{s_0}^{q,\kappa}\,\,\langle j\rangle^{-s_0}\\
\nonumber &\lesssim\interleave\mathscr{R}_{m}\interleave_{s_0}^{q,\kappa}\|\mathbf{e}_{0,j}\|_{H^{s_0}}\,\,\langle j\rangle^{-s_0}\\
&\quad\lesssim\interleave\mathscr{R}_{m}\interleave_{s_0}^{q,\kappa}=\kappa\, \delta_{m}(s_{0}).
\end{align}
Consequently  we deduce from \eqref{hypothesis of induction for deltamprime},\eqref{whab1} and \eqref{small-C3}
\begin{align}\label{Mahma1-R}
\|\mu_{j}^{m+1}-\mu_{j}^{m}\|^{q,\kappa}
 &\leqslant C \kappa\,\delta_{0}(s_{h})N_{0}^{\mu_{2}}N_{m}^{-\mu_{2}}\\
 \nonumber&\leqslant C \varepsilon \kappa^{-2}\,N_{0}^{\mu_{2}}N_{m}^{-\mu_{2}}.
\end{align}
Since the assumption \eqref{reg-G-10} holds true with  $\mathscr{D}_{m}$, that is, 
\begin{align}\label{Poki1} \forall\,j,j_{0}\in\mathbb{S}_0^c,\quad\max_{0\leqslant|\beta| \leqslant q}\sup_{\lambda\in\mathcal{O}}\left|\partial_{\lambda}^{\beta}\left(\mu_{j}^m(\lambda)-\mu_{j_{0}}^m(\lambda)\right)\right|\leqslant C\,|j-j_{0}|
\end{align}
then we get by the triangle inequality, \eqref{Mahma1-R} and \eqref{small-C3}
\begin{align*} \forall\,j, j_{0}\in\mathbb{S}_0^c,\quad\max_{0\leqslant|\beta| \leqslant q}\sup_{\lambda\in\mathcal{O}}\left|\partial_{\lambda}^{\beta}\left(\mu_{j}^{m+1}(\lambda)-\mu_{j_{0}}^{m+1}(\lambda)\right)\right|&\leqslant C\big(1+\varepsilon \kappa^{-2-q}\,N_{0}^{\mu_{2}}N_{m}^{-\mu_{2}}\big)\,|j-j_{0}|\\
&\leqslant C\big(1+\varepsilon^{\frac{1}{3+q}} N_{0}^{\mu_{2}}N_{m}^{-\mu_{2}}\big)\,|j-j_{0}|
\end{align*}
Therefore using    the convergence of the series $\sum_m N_m^{-\mu_2}$ and \eqref{Mahma1-R} allows to guarantee the required  assumption with the  same constant $C$ independently of $m$ and this completes the induction principle. \\
Next, we shall provide some estimates for $\Psi_m$ that will be used later to study the strong convergence.  
Applying \eqref{link Psi and R} combined with Lemma  \ref{Lem-Rgv1} and $s_0+\tau_2(1+q)+1\leqslant \overline{s}_l$ yield
\begin{align}\label{link-ZD3}
\nonumber  \|\Psi_m\interleave_{s_0+1}^{q,\kappa}&\leqslant C \kappa ^{-1}\interleave P_{N_{m}}\mathscr{R}_{m}\interleave_{{s_0+\tau_2(q+1)+1}}^{q,\kappa}\\
 &\leqslant\, C\, \delta_m(\overline{s}_l).
\end{align}
Then we infer from \eqref{hypothesis of induction for deltamprime} and \eqref{Conv-P3} 
\begin{align}\label{link-Z3}
\nonumber  \|\Psi_m\interleave_{s_0+1}^{q,\kappa}\leqslant&\, C\, \delta_{0}(s_{h})N_{0}^{\mu_{2}}N_{m}^{-\mu_{2}}\\
 &\leqslant\, C\,\varepsilon \kappa^{-3}N_0^{\mu_2} \,N_{m}^{-\mu_{2}}.
\end{align}
Let us now discuss the persistence of higher regularity. Take   ${s\in[s_{0},S],}$ then from \eqref{KAM step remainder term} and \eqref{hypothesis of induction for deltamprime} we find 
\begin{align*}
\delta_{m+1}(s)  \leqslant &\, \delta_{m}(s)\Big(1+CN_{m}^{\tau_{2} q+\tau_{2}}\delta_{m}(s_{0})\Big)\\
& \leqslant \, \delta_{m}(s)\Big(1+CN_{0}^{{\mu}_{2}}N_{m}^{\tau_{2} q+\tau_{2}-\mu_{2}}\delta_{0}({s}_{h})\Big).
\end{align*}
Putting together this estimate with  Lemma \ref{lemma sum Nn},  \eqref{whab1} and \eqref{Conv-P3}  we infer \begin{align}\label{uniform estimate of deltamsprime}
\nonumber \forall\, s\in[s_0,S],\,\, \forall\, m\in\mathbb{N},\quad \delta_{m}(s)\leqslant&\,\delta_{0}(s)\prod_{n=0}^{\infty}\Big(1+CN_{0}^{{\mu}_{2}}N_{m}^{\tau_{2} q+\tau_{2}-\mu_{2}}\delta_{0}({s}_{h})\Big)\\
\nonumber&\leqslant C\delta_{0}(s)e^{\sum_{n\in\N}CN_{0}^{{\mu}_{2}}N_{n}^{\tau_{2} q+\tau_{2}-\mu_{2}}\delta_{0}({s}_{h})}\\
\nonumber&\quad\leqslant C\delta_{0}(s)e^{CN_{0}^{\tau_{2} q+\tau_{2}}\delta_{0}({s}_{h})}\\
&\qquad \leqslant C\varepsilon\kappa^{-3}\left(1+\|\mathfrak{I}_{0}\|_{s+\sigma}^{q,\kappa}\right).
\end{align}
Thus  \eqref{link Psi and R} together  with Lemma  \ref{Lem-Rgv1} and  interpolation inequalities   lead to
 \begin{align}\label{link-Z1}
\nonumber \|\Psi_m\interleave_{s}^{q,\kappa}\leqslant&\, C \kappa ^{-1}\interleave P_{N_{m}}\mathscr{R}_{m}\interleave_{{s+\tau_{2} q+\tau_{2}}}^{q,\kappa}\\
\nonumber&\leqslant\, C\, \delta_m(s+\tau_{2} q+\tau_{2})\\
&\quad \leqslant\, C\, \delta_m^{\overline\theta}(s_0)\delta_m^{1-\overline\theta}(s+\tau_{2} q+\tau_{2}+1).
\end{align}
with $\overline\theta=\frac{1}{s-s_0+\tau_2(1+q)}.$ Plugging \eqref{hypothesis of induction for deltamprime} and  \eqref{uniform estimate of deltamsprime}  into \eqref{link-Z1} and using    \eqref{Conv-P3} we obtain
\begin{align}\label{link-Z2}
\nonumber \|\Psi_m\interleave_{s}^{q,\kappa}&\leqslant\, C\, \delta_0^{\overline\theta}(s_h)\delta_0^{1-\overline\theta}(s+\tau_{2} q+\tau_{2}+1)N_0^{\mu_2\overline\theta} N_m^{-\mu_2\overline\theta}\\
&\leqslant\, C\,\varepsilon_0^{\overline\theta} \delta_0^{1-\overline\theta}(s+\tau_{2} q+\tau_{2}+1) N_m^{-\mu_2\overline\theta}.
\end{align}
We observe that one also deduces from \eqref{Bir-11},  the second inequality of  \eqref{link-Z1} and \eqref{uniform estimate of deltamsprime} that
\begin{align}\label{link-ZP2}
 \forall\, s\in[ s_0,S],\quad \sup_{m\in\N}\left(\interleave\Sigma_m\interleave_{s}^{q,\kappa}+\interleave\Psi_m\interleave_{s}^{q,\kappa}\right)&\leqslant\,C\varepsilon\kappa^{-3}\left(1+\|\mathfrak{I}_{0}\interleave_{s+\sigma+\tau_2(1+q)}^{q,\kappa}\right).
\end{align}
Hence we find
\begin{align}\label{link-ZPDD2}
 \forall\, s\in[ s_0,S],\quad \sup_{m\in\N}\interleave\Phi_m^{\pm1}-\textnormal{Id}\interleave_{s}^{q,\kappa}&\leqslant\,C\varepsilon\kappa^{-3}\left(1+\|\mathfrak{I}_{0}\interleave_{s+\sigma+\tau_2(1+q)}^{q,\kappa}\right).
\end{align}
$\blacktriangleright$ \textbf{KAM conclusion}. Consider the sequence of operators 
$\left(\widehat{\Phi}_m\right)_{m\in\mathbb{N}},$ 
\begin{equation}\label{Def-Phi}\widehat\Phi_0\triangleq \Phi_0\quad\hbox{and}\quad \quad  \forall m\geqslant1,\,\, \widehat\Phi_m\triangleq \Phi_0\circ\Phi_1\circ...\circ\Phi_m.
\end{equation} 
Then it is clear from the identity ${\Phi}_{m}=\hbox{Id}+\Psi_m$ that $\widehat\Phi_{m+1}=\widehat\Phi_{m}+\widehat\Phi_{m}\Psi_{m+1}$.  Using  the law products of Lemma \ref{comm-pseudo1} we obtain
$$\begin{array}{rcl}
\interleave\widehat{\Phi}_{m+1}\interleave_{s_{0}+1}^{q,\kappa} & \leqslant & \interleave \widehat{\Phi}_{m}\interleave_{s_{0}+1}^{q,\kappa} \Big(1+C\interleave\Psi_{m+1}\interleave_{s_{0}+1}^{q,\kappa} \Big).
\end{array}$$
Then iterating this inequality and using \eqref{Conv-P3} and \eqref{link-Z3} yield
$$\begin{array}{rcl}
\interleave\widehat{\Phi}_{m+1}\interleave_{s_{0}+1}^{q,\kappa}
& \leqslant & \displaystyle\interleave\widehat{\Phi}_{0}\interleave_{s_{0}+1}^{q,\kappa}\prod_{n=1}^{m+1}\Big(1+C\interleave\Psi_{n}\interleave_{s_{0}+1}^{q,\kappa}\Big)\\
& \leqslant & \displaystyle\prod_{n=0}^{\infty}\Big(1+C\,{\varepsilon}_0N_{n}^{-\mu_{2}}\Big).
\end{array}$$
From the first condition of \eqref{Cond1} and  \eqref{definition of Nm} and one gets
$$\begin{array}{rcl}
\|\widehat{\Phi}_{m+1}\|_{\textnormal{{O-d}},q,s_{0}+1} 
& \leqslant & \displaystyle\prod_{n=0}^{\infty}\Big(1+C\,{\varepsilon}_04^{-(\frac32)^n}\Big).
\end{array}$$
Since  the infinite product converges and for ${\varepsilon}_0$ small enough  then  we obtain 
\begin{align}\label{phim-T}
\sup_{m\in\mathbb{N}}\interleave\widehat{\Phi}_{m}\interleave_{s_{0}+1}^{q,\kappa}\leqslant 2.
\end{align}
Next,  we intend to estimate the difference $\widehat{\Phi}_{m+1}-\widehat{\Phi}_{m}$ and for this aim we use the law products of \mbox{Lemma \ref{comm-pseudo1}} combined with \eqref{link-Z3} and \eqref{phim-T}
{
\begin{align}\label{Trank-P1}
 \interleave\widehat{\Phi}_{m+1}-\widehat{\Phi}_{m}\interleave_{s_{0}+1}^{q,\kappa}
\nonumber& \leqslant   C \interleave \widehat{\Phi}_{m}\interleave_{s_{0}+1}^{q,\kappa}\interleave\Psi_{m+1}\interleave_{s_{0}+1}^{q,\kappa}\\
& \leqslant C\,\delta_{0}(s_{h})N_{0}^{\mu_{2}} \,N_{m+1}^{-\mu_2}.
\end{align}
Consequently, Lemma \ref{lemma sum Nn} gives
\begin{align}\label{Conv-Z1}
\sum_{m=0}^\infty\interleave\widehat{\Phi}_{m+1}-\widehat{\Phi}_{m}\interleave_{s_{0}+1}^{q,\kappa}\leqslant C\,\delta_{0}(s_{h}).
\end{align}
Then by a completeness argument we deduce that the series $\displaystyle \widehat\Phi_{0}+\sum_{m\in\mathbb{N}}(\widehat\Phi_{m+1}-\widehat\Phi_{m})$ converges to an element $\Phi_\infty$. Furthermore, we get in view of \eqref{Trank-P1}
\begin{align}\label{Conv-op-od}
\nonumber\interleave\widehat\Phi_{m}-\Phi_\infty\interleave_{s_{0}+1}^{q,\kappa}\leqslant & \sum_{j\geqslant m}\interleave\widehat\Phi_{j+1}-\widehat\Phi_{j}\interleave_{s_{0}+1}^{q,\kappa}\\
\nonumber&\leqslant  C\,\delta_{0}(s_{h})N_{0}^{\mu_{2}} \, \sum_{j\geqslant m}N_{j+1}^{-\mu_2}\\
&\quad \leqslant  C\,\delta_{0}(s_{h})N_{0}^{\mu_{2}}N_{m+1}^{-\mu_2}.
\end{align} 
}
Notice that one also  deduces from \eqref{phim-T}
\begin{equation}\label{Conv-op-ES}
\interleave \Phi_\infty\interleave_{s_{0}+1}^{q,\kappa}\leqslant 2.
\end{equation} 
On the other hand, using \eqref{Conv-Z1} combined with  \eqref{link-Z1} for $m=0$ and \eqref{Conv-T2}
\begin{align}\label{Conv-UL1}
\nonumber \interleave \Phi_\infty-\textnormal{Id}\interleave_{s_0+1}^{q,\kappa}\leqslant&\, \sum_{m=0}^\infty\interleave \widehat{\Phi}_{m+1}-\widehat{\Phi}_{m}\interleave_{s_0+1}^{q,\kappa}+\interleave {\Psi}_{0}\interleave_{s_0+1}^{q,\kappa}\\
 \leqslant&\, C\,\delta_{0}(s_{h}).
\end{align} 
Let us now analyze  the convergence in higher norms. Take  $s\in[ s_{0},S]$, then using the law products of of \mbox{Lemma \ref{comm-pseudo1}}, \eqref{link-Z3},  \eqref{link-Z2} and \eqref{phim-T}  we obtain
\begin{align}\label{Rec-YP1}
 \interleave\widehat{\Phi}_{m+1}\interleave_{s}^{q,\kappa} & \leqslant  \interleave \widehat{\Phi}_{m}\interleave_{s}^{q,\kappa}\left(1+C\interleave \Psi_{m+1}\interleave_{s_0}^{q,\kappa}\right)+C\interleave\widehat{\Phi}_{m}\interleave_{s_0}^{q,\kappa}\interleave\Psi_{m+1}\interleave_{s}^{q,\kappa}\\
\nonumber& \leqslant  \interleave\widehat{\Phi}_{m}\interleave_{s}^{q,\kappa}\left(1+C\,{\varepsilon_0} \,N_{m+1}^{-\mu_{2}}\right)+C\,\delta_{0}^{\overline\theta}(s_{h})N_{0}^{\mu_{2} \overline\theta} \delta_0^{1-\overline\theta}(s+\tau_{2} q+\tau_{2}+1) N_m^{-\mu_2\overline\theta}.
\end{align}
From the first condition of \eqref{Cond1} and  \eqref{definition of Nm} and one gets
\begin{align*}
\prod_{n=0}^{\infty}\left(1+C\,{\varepsilon}_0 \,N_{n}^{-\mu_{2}}\right)
& \leqslant  \prod_{n=0}^{\infty}\Big(1+C\,{\varepsilon}_04^{-(\frac32)^n}\Big)\\
& \leqslant 2
\end{align*}
where the last inequality holds if ${\varepsilon}_0$ is chosen small enough. Similarly we get
$$
\sum_{n=0}^{\infty}\,N_{n}^{-\mu_{2}\overline\theta}<\infty.
$$
Then applying \eqref{Ind-res1} to  \eqref{Rec-YP1} and using \eqref{link Psi and R} yields
\begin{align*}
\nonumber \sup_{m\in\mathbb{N}}\interleave\widehat{\Phi}_{m}\interleave_{s}^{q,\kappa} 
& \leqslant  C\Big( \interleave{\Phi}_{0}\interleave_{s}^{q,\kappa}+\delta_{0}^{\overline\theta}(s_{h})N_{0}^{\mu_{2} \overline\theta} \delta_0^{1-\overline\theta}(s+\tau_{2} q+\tau_{2}+1) \Big)\\
& \leqslant  C\Big(1+ \delta_0(s+\tau_{2} q+\tau_{2})+ \delta_{0}^{\overline\theta}(s_{h})N_{0}^{\mu_{2} \overline\theta}\delta_0^{1-\overline\theta}\big(s+\tau_{2} q+\tau_{2}+1\big) \Big).
\end{align*}
Using interpolation inequalities and \eqref{Conv-P3} we find
\begin{align}\label{Rec-YP2}
 \sup_{m\in\mathbb{N}}\interleave\widehat{\Phi}_{m}\interleave_{s}^{q,\kappa} 
& \leqslant  C\Big(  1+\delta_0\big(s+\tau_{2} q+\tau_{2}+1\big) \Big).
\end{align}
The next goal is to estimate the difference $\interleave\widehat{\Phi}_{m+1}-\widehat{\Phi}_{m}\interleave_{s}^{q,\kappa}  $. Then by the law products of Lemma \ref{comm-pseudo1} combined with the first inequality in \eqref{link-Z3}, \eqref{link-Z2}, \eqref{phim-T} and \eqref{Rec-YP2} one gets
\begin{align*}
\interleave\widehat{\Phi}_{m+1}-\widehat{\Phi}_{m}\interleave_{s}^{q,\kappa}
\leqslant&\, C\Big(\interleave\widehat{\Phi}_{m}\interleave_{s}^{q,\kappa}\interleave\Psi_{m+1}\interleave_{s_0}^{q,\kappa}+\interleave\widehat{\Phi}_{m}\interleave_{s_0}^{q,\kappa}\interleave \Psi_{m+1}\interleave_{s}^{q,\kappa}\Big)\\
&\leqslant C\,\delta_{0}(s_{h})N_{0}^{\mu_{2}} \,N_{m+1}^{-\mu_2}\Big(1+\delta_0\big(s+\tau_{2} q+\tau_{2}+1\big) \Big)\\
&\quad+C\,\delta^{\overline\theta}_{0}(s_{h})N_{0}^{\mu_{2} \overline\theta}\,\delta_0^{1-\overline\theta}\big(s+\tau_{2} q+\tau_{2}+1\big) N_{m+1}^{-\mu_2\overline\theta}.
\end{align*}
Consequently, we obtain in view of Lemma \ref{lemma sum Nn}
\begin{align*}
\sum_{m=0}^\infty\interleave \widehat{\Phi}_{m+1}-\widehat{\Phi}_{m}\interleave_{s}^{q,\kappa}
&\leqslant C\,\delta_{0}(s_{h})\,\Big(1+\delta_0\big(s+\tau_{2} q+\tau_{2}+1\big) \Big)\\
&+C\,\delta^{\overline\theta}_{0}(s_{h})\,\delta_0^{1-\overline\theta}\big(s+\tau_{2} q+\tau_{2}+1\big).
\end{align*}
Using interpolation inequalities and the second condition in \eqref{Cond1} we get
\begin{align}\label{P-R-S1}
\sum_{m=0}^\infty\interleave \widehat{\Phi}_{m+1}-\widehat{\Phi}_{m}\interleave_{s}^{q,\kappa}
&\leqslant C\Big(\delta_{0}(s_{h})+\delta_0\big(s+\tau_{2} q+\tau_{2}+1\big) \Big).
\end{align}
Therefore we find from this latter inequality  combined with \eqref{Conv-P3} and \eqref{Rec-YP2} 
\begin{align}\label{Conv-opV1}
\nonumber \interleave\Phi_\infty\interleave_{s}^{q,\kappa}&\leqslant \sum_{m=0}^\infty\interleave \widehat{\Phi}_{m+1}-\widehat{\Phi}_{m}\interleave_{s}^{q,\kappa}+\interleave \widehat{\Phi}_{0}\interleave_{s}^{q,\kappa}\\
&\leqslant C\Big(1+\delta_0\big(s+\tau_{2} q+\tau_{2}+1\big) \Big).
\end{align} 
On the other hand, one can easily check that, using \eqref{P-R-S1} and the second inequality in  \eqref{link-Z1} with $m=0,$
\begin{align}\label{Conv-U1}
\nonumber \|\Phi_\infty-\textnormal{Id}\interleave_{s}^{q,\kappa}&\leqslant \sum_{m=0}^\infty\|\widehat{\Phi}_{m+1}-\widehat{\Phi}_{m}\interleave_{s}^{q,\kappa}+\|{\Psi}_{0}\interleave_{s}^{q,\kappa}\\
&\leqslant C\Big(\delta_{0}(s_{h})+\delta_0\big(s+\tau_{2} q+\tau_{2}+1\big) \Big).
\end{align} 
Applying Lemma \ref{Lem-Rgv1}-(iii), we obtain by virtue of \eqref{Conv-op-ES}, \eqref{Conv-opV1} and Sobolev embeddings 
\begin{align}\label{Ti-Ti}
\nonumber \|\Phi_{\infty}\rho\|_{s}^{q,\kappa}&\lesssim\interleave  \Phi_{\infty}\interleave_{s_0}^{q,\kappa}\|\rho\|_{s}^{q,\kappa}+\interleave  \Phi_{\infty}\interleave_{s}^{q,\kappa}\|\rho\|_{s_0}^{q,\kappa}\\
\nonumber &\lesssim \|\rho\|_{s}^{q,\kappa}+\Big(1+\delta_0\big(s+\tau_{2} q+\tau_{2}+1\big) \Big)\|\rho\|_{s_0}^{q,\kappa}\\
&\lesssim \|\rho\|_{s}^{q,\kappa}+\delta_0\big(s+\tau_{2} q+\tau_{2}+1\big) \|\rho\|_{s_0}^{q,\kappa}.
\end{align}
According to  \eqref{whab1} we infer 
\begin{align}\label{fgh1} 
 \delta_0\big(s+\tau_{2} q+\tau_{2}+1\big)
&\leqslant C\varepsilon\kappa^{-3}\left(1+\| \mathfrak{I}_{0}\|_{s+\sigma+\tau_{2} q+\tau_{2}+1}^{q,\kappa}\right).
\end{align}
Inserting \eqref{fgh1} into \eqref{Ti-Ti} and using \eqref{small-C3} combined with Sobolev embeddings and \eqref{Conv-T2} gives
\begin{align}\label{Ti-TiX}
\nonumber \|\Phi_{\infty}\rho\|_{s}^{q,\kappa}&\lesssim \|\rho\|_{s}^{q,\kappa}+\varepsilon\kappa^{-3}\left(1+\| \mathfrak{I}_{0}\|_{s+\sigma+\tau_{2} q+\tau_{2}+1}^{q,\kappa}\right) \|\rho\|_{s_0}^{q,\kappa}\\
&\lesssim \|\rho\|_{s}^{q,\kappa}+\varepsilon\kappa^{-3}\| \mathfrak{I}_{0}\|_{s+\sigma}^{q,\kappa} \|\rho\|_{s_0}^{q,\kappa}.
\end{align}
Notice that in the last line we have replaced $\sigma+\tau_{2} q+\tau_{2}+1$ by $\sigma$ which can be justified  by taking $\sigma$ sufficient large.\\
Similarly to \eqref{Ti-Ti} we get  by applying   Lemma \ref{Lem-Rgv1}-(iii)  combined with \eqref{Conv-U1} and \eqref{Conv-UL1}
\begin{align*}
\nonumber \big\|\big(\Phi_{\infty}-\textnormal{Id}\big)\rho\big\|_{s}^{q,\kappa}&\lesssim\interleave \Phi_{\infty}-\textnormal{Id}\interleave_{s_{0}}^{q,\kappa}\|\rho\|_{s}^{q,\kappa}+\interleave \Phi_{\infty}-\textnormal{Id}\interleave_{s}^{q,\kappa}\|\rho\|_{s_0}^{q,\kappa}\\
\nonumber &\lesssim \delta_0(s_h)\|\rho\|_{s}^{q,\kappa}+\Big(\delta_{0}(s_{h})+\delta_0\big(s+\tau_{2} q+\tau_{2}+1\big) \Big)\|\rho\|_{s_0}^{q,\kappa}\\
&\lesssim \delta_0(s_h)\|\rho\|_{s}^{q,\kappa}+\delta_0\big(s+\tau_{2} q+\tau_{2}+1\big) \|\rho\|_{s_0}^{q,\kappa}.
\end{align*}
Then we find from  \eqref{fgh1} and \eqref{Conv-P3} combined with Sobolev embeddings and \eqref{whab1}
\begin{align}\label{Ti-TiP}
\nonumber \big\|\big(\Phi_{\infty}-\textnormal{Id}\big)\rho\big\|_{s}^{q,\kappa}
&\lesssim \big(\varepsilon\kappa^{-3}+\delta_0(s_h)\big)\|\rho\|_{s}^{q,\kappa}+\varepsilon\kappa^{-3}\| \mathfrak{I}_{0}\|_{s+\sigma+\tau_{2} q+\tau_{2}+1}^{q,\kappa}  \|\rho\|_{s_0}^{q,\kappa}\\
&\lesssim \varepsilon\kappa^{-3}\|\rho\|_{s}^{q,\kappa}+\varepsilon\kappa^{-3}\| \mathfrak{I}_{0}\|_{s+\sigma+\tau_{2} q+\tau_{2}+1}^{s,\kappa}  \|\rho\|_{s_0}^{q,\kappa}.
\end{align}
{The estimates of $\Phi_{\infty}^{-1}, \Phi_{\infty}^{-1}-\widehat{\Phi}_{n}^{-1}$ can be checked using similar arguments.}

\smallskip

\ding{70} Next we shall study the asymptotic expansion of the eigenvalues. Summing up  in $m$  the estimates  \eqref{Mahma1-R} and using Lemma \ref{lemma sum Nn} we find
\begin{align}\label{sum-R1}
\nonumber \sum_{m=0}^\infty\|\mu_{j}^{m+1}-\mu_{j}^{m}\|^{q,\kappa}
&\leqslant C\kappa\,\delta_{0}(s_{h})N_{0}^{\mu_{2}}\sum_{m=0}^\infty N_{m}^{-\mu_{2}}\\
& \leqslant C\kappa\delta_{0}(s_{h}).
\end{align}
This shows that for each $j\in\mathbb{S}_0^c$  the sequence $(\mu_{j}^{m})_{m\in\mathbb{N}}$ converges in the space $W^{q,\infty,\gamma }(\mathcal{O},\mathbb{C})$ to an element denoted by  $\mu_{j}^{\infty}\in W^{q,\infty,\gamma }(\mathcal{O},\mathbb{C})$. 
In addition, for any $m\in\mathbb{N},$ we find in view of \eqref{Mahma1-R}
\begin{align*}
\|\mu_{j}^{\infty}-\mu_{j}^{m}\|^{q,\kappa}&\leqslant  \sum_{n=m}^\infty\|\mu_{j}^{n+1}-\mu_{j}^{n}\|^{q,\kappa}\\
&\leqslant \kappa\,\delta_{0}(s_{h})N_{0}^{\mu_{2}}\sum_{n=m}^\infty N_{n}^{-\mu_{2}}.
\end{align*}
Therefore we deduce according ton Lemma \ref{lemma sum Nn}  
\begin{align}\label{Spect-TD1}
\sup_{j\in \mathbb{S}_0^c}\|\mu_{j}^{\infty}-\mu_{j}^{m}\|^{q,\kappa}
&\leqslant C \kappa\delta_{0}(s_{h})N_{0}^{\mu_{2}} {N_{m}^{-\mu_{2}}}.
\end{align}
One also gets
\begin{align}\label{dekomp}
\nonumber \mu_{j}^{\infty}&=\mu_{j}^{0}+\sum_{m=0}^\infty\big(\mu_{j}^{m+1}-\mu_{j}^{m}\big)\\
&\triangleq \mu_{j}^{0}+ r_{j}^{\infty}
\end{align}
where $(\mu_{j}^{0})$ is  described in Proposition \ref{projection in the normal directions}-(i) and given by
$$
\mu_{j}^{0}(\lambda,i_0)= \Omega_{j}(\alpha)+jr^{1}(\lambda,i_0)+j|j|^{\alpha-1}r^{2}(\lambda,i_0).
$$
Therefore   \eqref{sum-R1} and \eqref{Conv-P3} yield
\begin{align*}
\|r_{j}^{\infty}\|^{q,\kappa}
&\leqslant C\, \kappa\,\delta_{0}(s_{h})\\
&\leqslant C\, \varepsilon \kappa^{-2}
\end{align*}
and this gives the first  result in \eqref{estimate rjinfty}.

\smallskip

Now, let us 
 consider the diagonal  operator $\mathscr{D}_{\infty}$ defined on the normal modes  by
 \begin{align}\label{Dinfty-op}
\forall (l,j)\in\mathbb{Z}^d\times\mathbb{S}_0^c,\quad \mathscr{D}_{\infty} {\bf e}_{l,j}=\ii\mu_{j}^{\infty}{\bf e}_{l,j}.
\end{align}
By virtue of \eqref{Top-NormX}  we obtain
\[
\|\mathscr{D}_{m}-\mathscr{D}_{\infty}\interleave_{s_0}^{q,\kappa}=\displaystyle\sup_{j\in\mathbb{S}_{0}^{c}}\|\mu_{j}^{m}-\mu_{j}^{\infty}\|^{q,\kappa}
\]
which gives by virtue  of \eqref{Spect-TD1}
\begin{align}\label{Conv-Dinf}
\nonumber  \interleave\mathscr{D}_{m}-\mathscr{D}_{\infty}\interleave_{s_{0}}^{q,\kappa}
&\leqslant C\, \kappa\,\delta_{0}(s_{h})N_{0}^{\mu_{2}} {N_{m}^{-\mu_{2}}}\\
 &\leqslant C\, \varepsilon\kappa^{-2}\,N_{0}^{\mu_{2}} {N_{m}^{-\mu_{2}}}.
 \end{align}
\ding{70} In what follows we shall prove that the Cantor set $\mathcal{O}_{\infty,n}^{\kappa,\tau_1,\tau_{2}}(i_{0})$ defined in Proposition \ref{reduction of the remainder term} satisfies
$$\mathcal{O}_{\infty,n}^{\kappa,\tau_1,\tau_{2}}(i_{0})\subset\bigcap_{m=0}^{n+1}\mathscr{O}_{m}^{\kappa}=\mathscr{O}_{n+1}^{\kappa}.$$
whereas  the intermediate Cantor sets are defined in \eqref{Cantor-SX}. For this goal we shall proceed by finite  induction on $m$ with $n$ fixed. First, we notice that by construction 
$\mathcal{O}_{\infty,n}^{\kappa,\tau_1,\tau_{2}}(i_{0})\subset\mathcal{O}\triangleq \mathscr{O}_{0}^{\kappa }.$ 
Now  assume  that $\mathscr{O}_{\infty,n}^{\kappa,\tau_1,\tau_{2}}(i_{0})\subset\mathscr{O}_{m}^{\kappa }$ for all $m\leqslant n$ and let us check that 
\begin{align}\label{Inc-H}
\mathcal{O}_{\infty,n}^{\kappa,\tau_1,\tau_{2}}(i_{0})\subset\mathscr{O}_{m+1}^{\kappa}.
\end{align}
Take $\lambda\in\mathcal{O}_{\infty,n}^{\kappa,\tau_1,\tau_{2}}(i_{0})$ and give  $(l,j,j_{0})\in\mathbb{Z}^{d }\times(\mathbb{S}_{0}^{c})^{2}$ such that $0\leqslant |l|\leqslant N_{m}$ and ${(l,j)\neq(0,j_0)}.$ Then we  may write by the triangle inequality, \eqref{Spect-TD1}, \eqref{Conv-T2} and \eqref{Conv-P3}
\begin{align*}
|\omega\cdot l+\mu_{j}^{m}(\lambda)-\mu_{j_{0}}^{m}(\lambda)| & \geqslant  \displaystyle|\omega\cdot l+\mu_{j}^{\infty}(\lambda)-\mu_{j_{0}}^{\infty}(\lambda)|-2\sup_{j\in\mathbb{S}_{0}^{c}}\|\mu_{j}^{m}-\mu_{j}^{\infty}\|^{q,\kappa}\\
& \geqslant  \displaystyle\frac{2\kappa\langle j-j_{0}\rangle }{\langle l\rangle^{\tau_{2}}}-C\kappa\delta_{0}(s_{h})N_{0}^{\mu_{2}} {N_{m}^{-\mu_{2}}}\\
& \geqslant  \displaystyle 2\frac{\kappa \langle j-j_{0}\rangle}{\langle l\rangle^{\tau_{2}}}-C\kappa\varepsilon_0 \langle l\rangle^{-\mu_2}\langle j-j_{0}\rangle.
\end{align*}
Then for $ \varepsilon_0$ small enough and using that  $\mu_2\geqslant\mu_c \geqslant \tau_2$, which follows from \eqref{cond-diman1}, we get
$$
\big|\omega\cdot l+\mu_{j}^{m}(\lambda)-\mu_{j_{0}}^{m}(\lambda)\big| > \frac{\kappa\langle j-j_{0}\rangle }{\langle l\rangle^{\tau_{2}}}.
$$
This  shows that $\lambda\in\mathscr{O}_{m+1}^{\kappa}$ and therefore the inclusion \eqref{Inc-H} holds.

\smallskip

\ding{70} In what follows  the {\it convergence} of the sequence $\left(\mathscr{L}_{m}\right)_{m\in\mathbb{N}}$ introduced in \eqref{Op-Lm} towards  the diagonal operator  $\mathscr{L}_{\infty}\triangleq \omega\cdot\partial_{\varphi}+\mathscr{D}_{\infty},$ where $\mathscr{D}_\infty$ is given by \eqref{Dinfty-op}. Making use of \eqref{Conv-Dinf} and \eqref{hypothesis of induction for deltamprime} we find
\begin{align}\label{FD-TU-1} 
\nonumber \interleave\mathscr{L}_{m}-\mathscr{L}_{\infty}\interleave_{s_{0}}^{q,\kappa }&\leqslant\interleave\mathscr{D}_{m}-\mathscr{D}_{\infty}\interleave_{s_{0}}^{q,\kappa }+\interleave\mathscr{R}_{m}\interleave_{s_{0}}^{q,\kappa }\\
&\leqslant C\, \kappa\,\delta_{0}(s_{h})N_{0}^{\mu_{2}} {N_{m}^{-\mu_{2}}},
\end{align}
which implies in particular that
\begin{align}\label{strong-Cv}
\lim_{m\rightarrow\infty} \interleave\mathscr{L}_{m}-\mathscr{L}_{\infty}\interleave_{s_{0}}^{q,\kappa }=0.
\end{align}
{On the other hand, we remark from \eqref{Def-Phi} and \eqref{Op-Lm1} that  
\begin{align*}
\forall \,\lambda\in {\mathscr{O}_{n+1}^{\kappa }},\quad\widehat{\Phi}_{n}^{-1}\mathscr{L}_{0}\widehat{\Phi}_{n}&=\big(\omega\cdot\partial_{\varphi}+\mathscr{D}_{n+1}+\mathscr{R}_{n+1}\big)\Pi_{\mathbb{S}_0}^\perp\\
&=\mathscr{L}_{\infty}+\big(\mathscr{D}_{n+1}-\mathscr{D}_{\infty}+\mathscr{R}_{n+1}\big)\Pi_{\mathbb{S}_0}^\perp
\end{align*} 
and therefore for any $\lambda\in  {\mathscr{O}_{n+1}^{\kappa }}$ we may write the decomposition
\begin{align}\label{wahda-mall-1}
\nonumber\Phi_{\infty}^{-1}\mathscr{L}_{0}\Phi_{\infty}&=\mathscr{L}_{\infty}+\big(\mathscr{D}_{n+1}-\mathscr{D}_{\infty}+\mathscr{R}_{n+1}\big)\Pi_{\mathbb{S}_0}^\perp\\
\nonumber&\quad+\Phi_{\infty}^{-1}\mathscr{L}_{0}\left(\Phi_{\infty}-\widehat{\Phi}_{n}\right)+\left(\Phi_{\infty}^{-1}-\widehat{\Phi}_{n}^{-1}\right)\mathscr{L}_{0}\widehat{\Phi}_{n}\\
&\triangleq\mathscr{L}_{\infty}+\mathtt{E}_{n,1}^2+\mathtt{E}_{n,2}^2+\mathtt{E}_{n,3}^2\triangleq\mathscr{L}_{\infty}+\mathtt{E}_{n}^3.
\end{align}
To  estimate $\mathtt{E}_{n,1}^2$ we use \eqref{Conv-Dinf} combined with \eqref{Def-Taou},
 \eqref{hypothesis of induction for deltamprime}, \eqref{whab1} and  \eqref{small-C3}
\begin{align}\label{bard1}
\nonumber \interleave\mathtt{E}_{n,1}^2\interleave_{s_0}^{q,\kappa}
&\leqslant C\, \kappa\,\delta_{0}(s_{h})N_{0}^{\mu_{2}} {N_{n+1}^{-\mu_{2}}}\\
&\leqslant C\varepsilon\kappa^{-2}N_{0}^{{\mu}_{2}}{N_{n+1}^{-\mu_{2}}}.
\end{align}
Putting together  Lemma \ref{Lem-Rgv1}-(iii) with \eqref{bard1} yields
\begin{align*}
\|\mathtt{E}_{n,1}^2 h\|_{s_0}^{q,\kappa}
&\leqslant C\varepsilon\kappa^{-2}N_{0}^{{\mu}_{2}}{N_{n+1}^{-\mu_{2}}} \|h\|_{s_0}^{q,\kappa}.
\end{align*}
Next let us prove  the estimates of  $\mathtt{E}_{n,2}^2$ and $ \mathtt{E}_{n,3}^2$ defined in \eqref{wahda-mall-1}. They can be implemented in a similar way and therefore we shall restrict the discussion to the term $\mathtt{E}_{n,2}^2$. Using \eqref{estimate on Phiinfty and its inverse} yields
\begin{align}\label{Taptap2}
\|\mathtt{E}_{n,2}^2 h\|_{s_0}^{q,\kappa}\lesssim\big\|\mathscr{L}_{0}\big(\Phi_{\infty}-\widehat{\Phi}_{n}\big) h\big\|_{s_0}^{q,\kappa}+{\varepsilon\kappa^{-3}}\| \mathfrak{I}_{0}\|_{s_0+\sigma}^{q,\kappa}\big\|\mathscr{L}_{0}\big(\Phi_{\infty}-\widehat{\Phi}_{n}\big) h\big\|_{s_0}^{q,\kappa}.
\end{align}
Then combining this estimate with  \eqref{Taptap2} and  \eqref{small-C3} we obtain
 \begin{align*}
\|\mathtt{E}_{n,2}^2h\|_{s_0}^{q,\kappa}&\lesssim\big\|\mathscr{L}_{0}\big(\Phi_{\infty}-\widehat{\Phi}_{n}\big)h\big\|_{s_0}^{q,\kappa}\\
&\lesssim \big\|\big(\Phi_{\infty}-\widehat{\Phi}_{n}\big)h\big\|_{s_0+3}^{q,\kappa}.
\end{align*}
Hence \eqref{Conv-op-od} together with   Lemma \ref{Lem-Rgv1},  \eqref{small-C3} and \eqref{whab1} allow to find
\begin{align*}
\|\mathtt{E}_{n,2}^2h\|_{s_0}^{q,\kappa}&\lesssim\interleave\Phi_{\infty}-\widehat{\Phi}_{n}\interleave_{s_0+3}^{q,\kappa}\|h\|_{s_0+3}^{q,\kappa}\\
&\leqslant  C\,\delta_{0}(s_{h})N_{0}^{\mu_{2}}N_{m+1}^{-\mu_2}\|h\|_{s_0+3}^{q,\kappa}\\
&\leqslant  C\,\varepsilon\kappa^{-3}N_0^{\mu_{2}}N_{m+1}^{-\mu_2}\|h\|_{s_0+3}^{q,\kappa}.
\end{align*}
The estimate of  $\mathtt{E}_{n,3}^2$ can be done in a similar way to $\mathtt{E}_{n,1}^2$  and we get the same estimate. Putting together the foregoing estimates yields  \eqref{Error-Est-2D}.\\
}\\
\ding{70} The next task is to prove the second estimate of \eqref{estimate rjinfty}. Define
$$
\widehat{\delta}_{m}(s)\triangleq \kappa ^{-1}\interleave\mathscr{R}_{m}\interleave_{2\overline\alpha-1+\epsilon,s,0}^{q,\kappa},
$$
where the corresponding norm is introduced in \eqref{Def-pseud-w}.
Then we shall prove by induction on $m\in\mathbb{N}$ that
\begin{align}\label{Ind-Ty1}
\widehat{\delta}_{m}(s_{0})\leqslant C\varepsilon\kappa^{-3}N_{0}^{\mu_{2}}N_{m}^{-\mu_{2}}\quad \mbox{ and }\quad \widehat{\delta}_{m}(s_{h})\leqslant C\varepsilon\kappa^{-3}\left(2-\frac{1}{m+1}\right).
\end{align}
 Using   Proposition  \ref{projection in the normal directions}-(iii) combined with \eqref{small-C3} one gets 
\begin{align}\label{tip-op1}
\nonumber\widehat{\delta}_{0}(s_{h}) 
&\leqslant  C\varepsilon\kappa^{-3}\left(1+\|\mathfrak{I}_{0}\|_{s_h+\sigma}^{q,\kappa}\right)
\\
& \leqslant C\varepsilon\kappa^{-3}.
\end{align}
which shows that \eqref{Ind-Ty1} holds true for $m=0$.
 Now assume that \eqref{Ind-Ty1} is satisfied at the order $m$ and let us check it at the order $m+1$. Applying Lemma \ref{comm-pseudo1} to \eqref{Deco-T1} and using \eqref{phi-inverse1} we obtain the expression
 \begin{align}\label{Deco-TP1}
\interleave \mathscr{R}_{m+1}\interleave_{2\overline\alpha-1+\epsilon,s,0}^{q,\kappa}\leqslant  \interleave \Phi_{m}^{-1}P_{N_{m}}^{\perp}\mathscr{R}_{m}\interleave_{2\overline\alpha-1+\epsilon,s,0}^{q,\kappa}&+C
\interleave \Phi_{m}^{-1}\interleave_{0,s,0}^{q,\kappa}\interleave \mathcal{U}_m\interleave_{2\overline\alpha-1+\epsilon,s_0,0}^{q,\kappa}\\
\nonumber&
+C\interleave \Phi_{m}^{-1}\interleave_{0,s_0,0}^{q,\kappa}\interleave \mathcal{U}_m\interleave_{2\overline\alpha-1+\epsilon,s,0}^{q,\kappa},
\end{align}
with
$$
\mathcal{U}_m\triangleq -\Psi_{m}\,\lfloor P_{N_{m}}\mathscr{R}_{m}\rfloor +\mathscr{R}_{m}\Psi_{m}.
$$

Then putting together \eqref{Deco-TP1}, \eqref{link-ZPDD2} and \eqref{small-C3} gives 
 \begin{align}\label{Deco-TPS1}
\nonumber \interleave \mathscr{R}_{m+1}\interleave_{2\overline\alpha-1+\epsilon,s,0}^{q,\kappa}\leqslant  \interleave \Phi_{m}^{-1}P_{N_{m}}^{\perp}\mathscr{R}_{m}\interleave_{2\overline\alpha-1+\epsilon,s,0}^{q,\kappa}&+C\varepsilon\kappa^{-3}\left(1+\|\mathfrak{I}_{0}\interleave_{s+\sigma}^{q,\kappa}\right)\interleave \mathcal{U}_m\interleave_{2\overline\alpha-1+\epsilon,s_0,0}^{q,\kappa}\\
&
+C\interleave \mathcal{U}_m\interleave_{2\overline\alpha-1+\epsilon,s,0}^{q,\kappa}.
\end{align}

Using Lemma \ref{comm-pseudo1} we deduce that
 \begin{align*}
  \interleave \Phi_{m}^{-1}P_{N_{m}}^{\perp}\mathscr{R}_{m}\interleave_{2\overline\alpha-1+\epsilon,s,0}^{q,\kappa}\leqslant
  \interleave P_{N_{m}}^{\perp}\mathscr{R}_{m}\interleave_{2\overline\alpha-1+\epsilon,s,0}^{q,\kappa} &+C
\interleave \Sigma_{m}\interleave_{0,s,0}^{q,\kappa}\interleave \mathscr{R}_{m}\interleave_{2\overline\alpha-1+\epsilon,s_0,0}^{q,\kappa}\\
\nonumber&
+C\interleave \Sigma_{m}\interleave_{0,s_0,0}^{q,\kappa}\interleave \mathscr{R}_{m}\interleave_{2\overline\alpha-1+\epsilon,s,0}^{q,\kappa}.
\end{align*}
Consequently we find from Lemma \ref{Lem-Rgv1}-(iv), \eqref{AP-W34} and \eqref{whab1}
\begin{align}\label{Deco-THP1}
\nonumber  \interleave \Phi_{m}^{-1}P_{N_{m}}^{\perp}\mathscr{R}_{m}\interleave_{2\overline\alpha-1+\epsilon,s,0}^{q,\kappa}&\leqslant
  N_{m}^{s-\overline{s}}\interleave\mathscr{R}_{m}\interleave_{2\overline\alpha-1+\epsilon,\overline{s},0}^{q,\kappa} +C
N_{m}^{\tau_{2} q+\tau_{2}}\delta_{m}(s)\interleave \mathscr{R}_{m}\interleave_{2\overline\alpha-1+\epsilon,s_0,0}^{q,\kappa}\\
&
\quad+C\varepsilon \kappa^{-3}N_{0}^{\mu_{2}}N_{m}^{-\mu_{2}}\interleave \mathscr{R}_{m}\interleave_{2\overline\alpha-1+\epsilon,s,0}^{q,\kappa},
\end{align}
and for $s=s_0$, we get in view of \eqref{whab1} and \eqref{small-C3}
\begin{align}\label{Deco-THP3}
  \interleave \Phi_{m}^{-1}P_{N_{m}}^{\perp}\mathscr{R}_{m}\interleave_{2\overline\alpha-1+\epsilon,s_0,0}^{q,\kappa}&\leqslant
  N_{m}^{s_0-\overline{s}}\interleave\mathscr{R}_{m}\interleave_{2\overline\alpha-1+\epsilon,\overline{s},0}^{q,\kappa}\\
 \nonumber  & \quad+C\varepsilon \kappa^{-3}N_{0}^{\mu_{2}}N_{m}^{-\mu_{2}}\interleave \mathscr{R}_{m}\interleave_{2\overline\alpha-1+\epsilon,s,0}^{q,\kappa}.
\end{align}
Applying once again Lemma \ref{comm-pseudo1} yields
\begin{align*}
\interleave\mathcal{U}_m\interleave_{2\overline\alpha-1+\epsilon,s,0}^{q,\kappa}&\lesssim\interleave \Psi_m\interleave_{0,s,0}^{q,\kappa}\interleave \mathscr{R}_{m}\interleave_{2\overline\alpha-1+\epsilon,s_0,0}^{q,\kappa}
+\interleave \Psi_m\interleave_{0,s_0,0}^{q,\kappa}\interleave \mathscr{R}_{m}\interleave_{2\overline\alpha-1+\epsilon,s,0}^{q,\kappa}.
\end{align*}
Hence we derive from  \eqref{link-Z3} and \eqref{link-Z1}
\begin{align*}
\interleave\mathcal{U}_m\interleave_{2\overline\alpha-1+\epsilon,s_h,0}^{q,\kappa}&\lesssim\,\varepsilon\kappa^{-3}\interleave \mathscr{R}_{m}\interleave_{2\overline\alpha-1+\epsilon,s_0,0}^{q,\kappa}
+\varepsilon \kappa^{-3}N_0^{\mu_2} \,N_{m}^{-\mu_{2}}\interleave \mathscr{R}_{m}\interleave_{2\overline\alpha-1+\epsilon,s_h,0}^{q,\kappa}.
\end{align*}
In addition, for $s=s_0$, we obtain 
\begin{align}\label{Deco-ZZPP1}
\interleave\mathcal{U}_m\interleave_{2\overline\alpha-1+\epsilon,s_0,0}^{q,\kappa}&\lesssim \,\varepsilon \kappa^{-3}N_0^{\mu_2} \,N_{m}^{-\mu_{2}}\interleave \mathscr{R}_{m}\interleave_{2\overline\alpha-1+\epsilon,s_0,0}^{q,\kappa}.
\end{align}
Then plugging \eqref{Deco-ZZPP1} and \eqref{Deco-THP3} into \eqref{Deco-TPS1} and using \eqref{small-C3} we find
\begin{align*}
\widehat{\delta}_{m+1}(s_0)\leqslant \, C  \,\varepsilon \kappa^{-3}N_0^{\mu_2} \,N_{m}^{-\mu_{2}}\widehat{\delta}_{m}(s_0)
+C N_{m}^{s_0-{s}_h}\widehat{\delta}_{m}({s}_h).
\end{align*}
In a similar way, we find
$$
\widehat{\delta}_{m+1}(s_h)\leqslant  \,\widehat{\delta}_{m}(s_h)+C\varepsilon \kappa^{-3}N_0^{\mu_2} \,N_{m}^{-\mu_{2}}\widehat{\delta}_{m}(s_h)\
+C N_{m}^{\tau_{2} q+\tau_{2}}\delta_{m}(s_h)\widehat{\delta}_{m}({s}_0)+\varepsilon\kappa^{-3} \widehat{\delta}_{m}(s_0).
$$
Applying \eqref{hypothesis of induction for deltamprime} yields
\begin{align*}
\widehat{\delta}_{m+1}(s_h)\leqslant  \,\widehat{\delta}_{m}(s_h)+C\varepsilon \kappa^{-3}N_0^{\mu_2} \,N_{m}^{-\mu_{2}}\widehat{\delta}_{m}(s_h)\
+C \varepsilon \kappa^{-3}N_{m}^{\tau_{2} q+\tau_{2}}\widehat{\delta}_{m}({s}_0)+\varepsilon\kappa^{-3} \widehat{\delta}_{m}(s_0).
\end{align*}
Therefore with  the constraints on $\mu_2$ and $s_h$ imposed in \eqref{cond-diman1} we achieve the validity of the induction \eqref{Ind-Ty1}  at the order $m+1$ in a similar way to that of \eqref{hypothesis of induction for deltamprime}.\\
{
 Next, we will see how to  deduce   the second estimate of \eqref{estimate rjinfty}.
Recall from \eqref{Spect-T2} and \eqref{dekomp} that
$$r_{j}^{\infty}=\sum_{m=0}^{+\infty}r_{j}^{m}\quad\quad \textnormal{with}\quad\quad r_{j}^{m}=\big\langle P_{N_{m}}\mathscr{R}_{m}\mathbf{e}_{0,j},\mathbf{e}_{0,j}\big\rangle_{L^{2}(\mathbb{T}^{d +1},\mathbb{C})}.$$
Then for any $a\in\RR$ and $j\in\Z^\star$ it is clear from integration by parts that
$$\big\langle P_{N_{m}}\mathscr{R}_{m}\mathbf{e}_{0,j},\mathbf{e}_{0,j}\big\rangle_{L^{2}(\mathbb{T}^{d +1},\mathbb{C})}=\frac{1}{|j|^a}\big\langle P_{N_{m}}\Lambda^a\mathscr{R}_{m}\mathbf{e}_{0,j},\mathbf{e}_{0,j}\big\rangle_{L^{2}(\mathbb{T}^{d +1},\mathbb{C})},
$$
with $\Lambda=\sqrt{-\Delta}.$ 
Using a duality argument  combined with   Lemma \ref{Lem-Rgv1}-(iii) we obtain
\begin{align*}
\|\langle P_{N_{m}}\Lambda^a\mathscr{R}_{m}\mathbf{e}_{0,j},\mathbf{e}_{0,j}\rangle_{L^{2}(\mathbb{T}^{d +1},\mathbb{C})}\|^{q,\kappa}&\leqslant C \interleave\Lambda^a\mathscr{R}_{m}\interleave_{0,s_0-1,0}^{q,\kappa}.
\end{align*}
Applying Lemma \ref{comm-pseudo1} with $a=1-2\overline\alpha-\epsilon$ combined with \eqref{Ind-Ty1} yields
\begin{align*}
 \interleave\Lambda^a\mathscr{R}_{m}\interleave_{0,s_0-1,0}^{q,\kappa} &\leqslant C \interleave\mathscr{R}_{m}\interleave_{2\overline\alpha-1+\epsilon,s_0,0}^{q,\kappa} \\
 &\leqslant  C \kappa\,\widehat{\delta}_{m}(s_{0})\\
& \leqslant C\varepsilon\kappa^{-2}\,N_{0}^{\mu_{2}}N_{m}^{-\mu_{2}}.
\end{align*}
It follows that
$$
\| r_{j}^{m}\|^{q,\kappa}\leqslant  C|j|^{2\overline\alpha-1+\epsilon}\varepsilon\kappa^{-2}\,N_{0}^{\mu_{2}}N_{m}^{-\mu_{2}}.
$$
Consequently, we obtain in view of Lemma \ref{lemma sum Nn}
\begin{align*}
\| r_{j}^{\infty}\|^{q,\kappa}&\lesssim|j|^{2\overline\alpha-1+\epsilon}\varepsilon\kappa^{-2}\,N_{0}^{\mu_{2}}\sum_{m=0}^{+\infty}N_{m}^{-\mu_{2}}\\
& \lesssim|j|^{2\overline\alpha-1+\epsilon}\varepsilon\kappa^{-2}.
\end{align*}
This completes the proof of \eqref{estimate rjinfty}.

\smallskip

{\bf{(ii)}} From \eqref{Deco-T1} and \eqref{phi-inverse1} we can write
$$
\mathscr{R}_{m+1}=(\hbox{Id}+\Sigma_m)U_{m}
$$
with
\begin{align}\label{bkl1}
U_{m}\triangleq P_{N_{m}}^{\perp}\mathscr{R}_{m}+\mathscr{R}_{m}\Psi_{m}-\Psi_{m}\big\lfloor P_{N_{m}}\mathscr{R}_{m}\big\rfloor.
\end{align}
Thus we get  from  straightforward computations
\begin{align}\label{Gall-01}
\nonumber\Delta_{12}U_{m}&=P_{N_{m}}^{\perp}\Delta_{12}\mathscr{R}_{m}+(\Delta_{12}\mathscr{R}_{m})\Psi_{m}^{[1]}+\mathscr{R}_{m}^{[2]}(\Delta_{12}\Psi_{m})\\
&\quad -(\Delta_{12}\Psi_{m})\big\lfloor P_{N_{m}}\mathscr{R}_{m}^{[1]}\big\rfloor-\Psi_{m}^{[2]}\big\lfloor P_{N_{m}}\Delta_{12}\mathscr{R}_{m}\big\rfloor
\end{align}
and
\begin{align}\label{gui1}
\Delta_{12}\mathscr{R}_{m+1}=\Delta_{12}U_{m}+(\Delta_{12}\Sigma_{m})U_{m}^{[1]}+\Sigma_{m}^{[2]}\Delta_{12}U_{m}.
\end{align}
Where we have used the notation $f^{[k]}$ to denote  the value of $f$ at the torus $i_k$ for $ k=1,2.$ 
Elementary calculations based on \eqref{phi-inverse1} give
$$\Delta_{12}\Sigma_{m}=\Delta_{12}\Phi_{m}^{-1}=-(\Phi_{m}^{[2]})^{-1}(\Delta_{12}\Psi_{m})(\Phi_{m}^{[1]})^{-1}.$$
 Applying  Lemma \ref{comm-pseudo1} and using \eqref{link-ZP2} combined with \eqref{small-C3}we obtain
 \begin{align}\label{link-ZP3}
\forall \, s\in[s_0,s_h],\quad   \interleave\Delta_{12}\Sigma_{m}\interleave_{s}^{q,\kappa}&\lesssim 
\interleave \Delta_{12}\Psi_{m}\interleave_{s}^{q,\kappa}.
\end{align}
By virtue of Lemma \ref{comm-pseudo1}, \eqref{link-ZP3} and  \eqref{gui1} we find
\begin{align}\label{gui2}
\nonumber \interleave \Delta_{12}\mathscr{R}_{m+1}\interleave_{s_0}^{q,\kappa}&\leqslant  \interleave\Delta_{12}U_{m}\interleave_{s_0}^{q,\kappa}+\interleave\Delta_{12}\Psi_{m}\interleave_{s_0}^{q,\kappa}\|U_{m}^{[1]}\interleave_{s_0}^{q,\kappa}\\
 &\quad+\interleave\Sigma_{m}^{[2]}\interleave_{s_0}^{q,\kappa}\interleave \Delta_{12}U_{m}\interleave_{s_0}^{q,\kappa}
\end{align}
and
\begin{align}\label{gui3}
\nonumber \interleave\Delta_{12}\mathscr{R}_{m+1}\interleave_{s_h}^{q,\kappa}&\leqslant  \interleave\Delta_{12}U_{m}\interleave_{s_h}^{q,\kappa}+\interleave\Delta_{12}\Psi_{m}\interleave_{s_0}^{q,\kappa}\interleave U_{m}^{[1]}\interleave_{s_h}^{q,\kappa}\\
\nonumber &\quad+\interleave\Delta_{12}\Psi_{m}\interleave_{s_h}^{q,\kappa}\interleave U_{m}^{[1]}\interleave_{s_0}^{q,\kappa}+\interleave\Sigma_{m}^{[2]}\interleave_{s_0}^{q,\kappa}\interleave\Delta_{12}U_{m}\interleave_{s_h}^{q,\kappa}\\
&\quad+ \interleave \Sigma_{m}^{[2]}\interleave_{s_h}^{q,\kappa}\interleave\Delta_{12}U_{m}\interleave_{s_0}^{q,\kappa}.
\end{align}
To estimate $U_m^{[1]}$ (to alleviate the notation we  remove the subscript $[1]$) described by \eqref{bkl1} we use Lemma \ref{comm-pseudo1} in order to get 
\begin{align}\label{bkl3}
\|U_{m}\interleave_{s_0}^{q,\kappa}&\leqslant \|\mathscr{R}_{m}\interleave_{s_0}^{q,\kappa}+\|\mathscr{R}_{m}\interleave_{s_0}^{q,\kappa} \|\Psi_m\interleave_{s_0}^{q,\kappa}
\end{align}
and
\begin{align}\label{bkl4}
\nonumber \interleave U_{m}\interleave_{s_h}^{q,\kappa}&\leqslant \interleave\mathscr{R}_{m}\interleave_{s_h}^{q,\kappa}+\interleave\mathscr{R}_{m}\interleave_{s_0}^{q,\kappa} \interleave\Psi_m\interleave_{s_h}^{q,\kappa}\\
&\quad+\interleave\mathscr{R}_{m}\interleave_{s_h}^{q,\kappa} \interleave\Psi_m\interleave_{s_0}^{q,\kappa}.
\end{align}
Applying \eqref{hypothesis of induction for deltamprime},\eqref{whab1}  and \eqref{link-Z3} together with \eqref{bkl3} we obtain
\begin{align}\label{bkl03}
\interleave U_{m}\interleave_{s_0}^{q,\kappa}&\leqslant C\varepsilon \kappa^{-2} N_0^{\mu_2} N_m^{-\mu_2}.\end{align}
On the other hand, 
from the first estimate of \eqref{link-Z1}, \eqref{hypothesis of induction for deltamprime} and \eqref{whab1} we deduce that
\begin{align}\label{bahya-11}
\nonumber \max_{k=1,2}\interleave\Psi_{m}^{[k]}\interleave_{s_h}^{q,\kappa}&\lesssim N_m^{\tau_2q+q}\delta_m(s_h)\\
&\lesssim N_m^{\tau_2q+q}\varepsilon \overline \gamma^{-3}.
\end{align}
Thus we deduce from \eqref{bkl4},  \eqref{hypothesis of induction for deltamprime} and \eqref{small-C3}
\begin{align}\label{bkl04}
 \sup_{m\in\NN} \interleave U_{m}\interleave_{s_h}^{q,\kappa}&\leqslant C\varepsilon \kappa^{-2}.
\end{align}
Inserting \eqref{bkl03} and \eqref{bkl04} into \eqref{gui3} yields
\begin{align}\label{guiM3}
\nonumber \interleave \Delta_{12}\mathscr{R}_{m+1}\interleave_{s_h}^{q,\kappa} &\leqslant  \interleave \Delta_{12}U_{m}\interleave_{s_h}^{q,\kappa}+C\varepsilon \kappa^{-2}\|\Delta_{12}\Psi_{m}\interleave_{s_0}^{q,\kappa}\\
\nonumber &\quad+C\varepsilon \kappa^{-2} N_0^{\mu_2} N_m^{-\mu_2}\interleave \Delta_{12}\Psi_{m}\interleave_{s_h}^{q,\kappa}+\interleave \Sigma_{m}^{[2]}\interleave_{s_0}^{q,\kappa}\interleave \Delta_{12}U_{m}\interleave_{s_h}^{q,\kappa}\\
&\quad+ \interleave\Sigma_{m}^{[2]}\interleave_{s_h}^{q,\kappa}\interleave \Delta_{12}U_{m}\interleave_{s_0}^{q,\kappa}.
\end{align}
From \eqref{AP-W34} and \eqref{whab1} we get
\begin{align}\label{tahya-b1}
\interleave \Sigma_{m}^{[2]}\interleave_{s_0}^{q,\kappa}\leqslant C\varepsilon \kappa^{-3}N_0^{\mu_2} N_m^{-\mu_2}\quad\hbox{and}\quad \interleave \Sigma_{m}^{[2]}\interleave_{s_h}^{q,\kappa}\leqslant C\varepsilon \kappa^{-3}N_m^{\tau_2 q+q}.
\end{align}
Plugging \eqref{tahya-b1} into \eqref{guiM3} implies
\begin{align}\label{guiM4}
\nonumber \interleave \Delta_{12}\mathscr{R}_{m+1}\interleave_{s_h}^{q,\kappa}&\leqslant \big(1+C\varepsilon \kappa^{-3}N_0^{\mu_2} N_m^{-\mu_2}\big) \interleave \Delta_{12}U_{m}\interleave_{s_h}^{q,\kappa}+C\varepsilon \kappa^{-3}N_m^{\tau_2 q+q}\interleave \Delta_{12}U_{m}\interleave_{s_0}^{q,\kappa}\\
&\quad+C\varepsilon \kappa^{-2} N_0^{\mu_2} N_m^{-\mu_2}\interleave \Delta_{12}\Psi_{m}\interleave_{s_h}^{q,\kappa}+C\varepsilon \kappa^{-2}\interleave \Delta_{12}\Psi_{m}\interleave_{s_0}^{q,\kappa}.
\end{align}
Similarly, by putting together \eqref{bkl03}, \eqref{tahya-b1} and \eqref{gui2} we infer
\begin{align}\label{guiP2}
\nonumber \interleave\Delta_{12}\mathscr{R}_{m+1}\interleave_{s_0}^{q,\kappa}&\leqslant \big(1+C\varepsilon \kappa^{-3}N_0^{\mu_2} N_m^{-\mu_2}\big) \interleave\Delta_{12}U_{m}\interleave_{s_0}^{q,\kappa}\\
&\quad+C\varepsilon \kappa^{-2} N_0^{\mu_2} N_m^{-\mu_2}\interleave\Delta_{12}\Psi_{m}\interleave_{s_0}^{q,\kappa}.
\end{align}
Coming back to \eqref{Gall-01}  and using Lemma \ref{comm-pseudo1} together with Lemma \ref{Lem-Rgv1} we get $ \forall \, s\in[s_0,s_h],$
\begin{align*}
  \interleave \Delta_{12}U_{m}\interleave_{s}^{q,\kappa}& \leqslant N_{m}^{s-s_h}   
\interleave \Delta_{12}\mathscr{R}_{m}\interleave_{s_h}^{q,\kappa}
+C\interleave \Delta_{12}\mathscr{R}_{m}\interleave_{s}^{q,\kappa}\max_{k=1,2}\interleave \Psi_{m}^{[k]}\interleave_{s_0}^{q,\kappa}\\
&\quad+C\interleave \Delta_{12}\mathscr{R}_{m}\interleave_{s_0}^{q,\kappa}\max_{k=1,2}\interleave \Psi_{m}^{[k]}\interleave_{s}^{q,\kappa}+C\interleave \Delta_{12}\Psi_{m}\interleave_{s}^{q,\kappa}\max_{k=1,2}\|\mathscr{R}_m^{[k]}\interleave_{s_0}^{q,\kappa}\\
&\quad+C\|\Delta_{12}\Psi_{m}\interleave_{s_0}^{q,\kappa}\max_{k=1,2}\|\mathscr{R}_m^{[k]}\interleave_{s}^{q,\kappa}.
\end{align*}
Putting together the preceding estimate with \eqref{bahya-11},  \eqref{hypothesis of induction for deltamprime} and \eqref{link-Z3}
\begin{align}\label{jima1}
 \nonumber \interleave \Delta_{12}U_{m}\interleave_{s_0}^{q,\kappa}&\leqslant N_{m}^{s_0-s_h}   
\interleave\Delta_{12}\mathscr{R}_{m}\interleave_{s_h}^{q,\kappa}
+C\varepsilon \kappa^{-3} N_0^{\mu_2} N_m^{-\mu_2}\interleave \Delta_{12}\mathscr{R}_{m}\interleave_{s_0}^{q,\kappa}\\
&\quad+C\varepsilon \kappa^{-2} N_0^{\mu_2} N_m^{-\mu_2}\interleave \Delta_{12}\Psi_{m}\interleave_{s_0}^{q,\kappa}
\end{align}
and
\begin{align*}
  \interleave\Delta_{12}U_{m}\interleave_{s_h}^{q,\kappa}&\leqslant   
\big(1+C\varepsilon \kappa^{-3} N_0^{\mu_2} N_m^{-\mu_2}\big)\interleave \Delta_{12}\mathscr{R}_{m}\interleave_{s_h}^{q,\kappa}+CN_m^{\tau_2q+q}\varepsilon \kappa^{-3}\interleave\Delta_{12}\mathscr{R}_{m}\interleave_{s_0}^{q,\kappa}\\
&\quad+C\varepsilon \kappa^{-2} N_0^{\mu_2} N_m^{-\mu_2}\interleave \Delta_{12}\Psi_{m}\interleave_{s_h}^{q,\kappa}+C\varepsilon \kappa^{-2}\interleave\Delta_{12}\Psi_{m}\interleave_{s_0}^{q,\kappa}.\end{align*}
Combining the foregoing estimate with \eqref{guiM4}, \eqref{guiP2}, \eqref{Conv-T2} and  \eqref{small-C3} allows to get
\begin{align}\label{guiPMS}
\nonumber \interleave \Delta_{12}\mathscr{R}_{m+1}\interleave_{s_h}^{q,\kappa}&\leqslant \Big(1+C\varepsilon \kappa^{-3} N_0^{\mu_2} N_m^{-\mu_2}+C\varepsilon \kappa^{-3}N_{m}^{s_0-s_h+\tau_2(1+q)}\Big)\interleave\Delta_{12}\mathscr{R}_{m}\interleave_{s_h}^{q,\kappa}\\
\nonumber &\quad+CN_m^{\tau_2q+q}\varepsilon \kappa^{-3}\interleave\Delta_{12}\mathscr{R}_{m}\interleave_{s_0}^{q,\kappa}+C\varepsilon \kappa^{-2} N_0^{\mu_2} N_m^{-\mu_2}\interleave\Delta_{12}\Psi_{m}\interleave_{s_h}^{q,\kappa}\\
&\quad +C\varepsilon \kappa^{-2}\interleave\Delta_{12}\Psi_{m}\interleave_{s_0}^{q,\kappa}.
\end{align}
Similarly, by putting together \eqref{guiP2}, \eqref{jima1} and \eqref{small-C3} we infer
\begin{align}\label{guiPM2}
\nonumber \interleave \Delta_{12}\mathscr{R}_{m+1}\interleave_{s_0}^{q,\kappa}&\leqslant N_{m}^{s_0-s_h}   
\interleave\Delta_{12}\mathscr{R}_{m}\interleave_{s_h}^{q,\kappa}
+C\varepsilon \kappa^{-3} N_0^{\mu_2} N_m^{-\mu_2}\interleave\Delta_{12}\mathscr{R}_{m}\interleave_{s_0}^{q,\kappa}\\
&\quad+C\varepsilon \kappa^{-2} N_0^{\mu_2} N_m^{-\mu_2}\interleave \Delta_{12}\Psi_{m}\interleave_{s_0}^{q,\kappa}.
\end{align}
In what follows we intend to estimate  $\Delta_{12}\Psi_{m}$ which is slightly delicate. According to \eqref{link Psi and R}, \eqref{Ext-psi-op} and \eqref{varr-d}  we may write 
\begin{align}\label{delta12psi}
  \interleave \Delta_{12}\Psi_m\interleave_{s}^{q,\kappa}=  \Big(\max_{|\gamma|\leqslant q}\sup_{\lambda\in\mathcal{O}}\sup_{n\in\Z}\sum_{|l|\leqslant N_m\atop |k|\leqslant N_m}\kappa^{2|\gamma|}\langle l,k\rangle^{2(s-|\gamma|)}\big|\partial_{\lambda}^{\gamma}\Delta_{12}\Psi_{n}^{n+k}(\lambda,l)\big|^2\Big)^{\frac12}
\end{align}
with
\begin{equation}\label{Ext-psi-opg}
\Psi_{j_{0}}^{j}(\lambda,l)=\left\lbrace\begin{array}{ll}
-\varrho_{j_{0}}^{j}(\lambda,l)\,\, r_{j_{0}}^{j}(\lambda,l),& \mbox{if }\quad (l,j)\neq(0,j_{0})\\
0, & \mbox{if }\quad (l,j)=(0,j_{0}),
\end{array}\right.
\end{equation}
and
\begin{align}\label{varr-dg}
\varrho_{j_{0}}^{j}(\lambda,l)\triangleq \frac{\chi\left((\omega\cdot l+\mu_{j}^m(\lambda)-\mu_{j_{0}}^m(\lambda))(\kappa\langle j-j_{0}\rangle)^{-1}\langle l\rangle^{\tau_{2}}\right)}{\omega\cdot l+\mu_{j}^m(\lambda)-\mu_{j_{0}}^m(\lambda)}\cdot
\end{align}
Remind from \eqref{coefficients of the remainder operator R} that $( r_{j_{0}}^{j}(\lambda,l)$ are the Fourier coefficients of $P_{N_m}\mathscr{R}_m$, that is,
$$
\ii\, r_{j_{0}}^{j}(\lambda,l)=\big \langle P_{N_m}\mathscr{R}_m{\bf e}_{0,j_0}, {\bf e}_{l,j}\big\rangle_{L^2(\T^{d+1})}.
$$
From this latter identity we get
$$
\ii\, \Delta_{12}r_{j_{0}}^{j}(\lambda,l)=\big \langle \Delta_{12}P_{N_m}\mathscr{R}_m{\bf e}_{0,j_0}, {\bf e}_{l,j}\big\rangle_{L^2(\T^{d+1})}.
$$
Now we need to estimate $\Delta_{12}\varrho_{j_{0}}^{j}(\lambda,l)$. For this aim, we write 
$$
\varrho_{j_{0}}^{j}(\lambda,l)=b_{l}^j\,\chi_1\big(b_{l}^jB_{l,j_0}^j(\lambda)\big)
$$
with
\begin{align}\label{g-lj-8}
B_{l,j_0}^j(\lambda)=\omega\cdot l+\mu_{j}^m(\lambda)-\mu_{j_{0}}^m(\lambda),\quad b_{l}^j=(\kappa\langle j\rangle)^{-1}\langle l\rangle^{\tau_{2}}, \quad \chi_1(x)=\tfrac{\chi(x)}{x}\cdot
\end{align}
According to \eqref{eg-k-1} one gets 
\begin{align}\label{eg-k-100}
 \forall |\gamma|\leqslant q,\quad \sup_{\lambda\in\mathcal{O}}\big|\partial_{\lambda}^{\gamma}\varrho_{j_{0}}^{j}(\lambda,l)\big|\leqslant C\kappa^{-(|\gamma|+1)} \langle l,j-j_{0}\rangle^{\tau_2(1+|\gamma|)+|\gamma|}.
\end{align}
Then using Taylor formula in a similar way to \eqref{II22} we find
\begin{align}\label{II23}
\Delta_{12}\varrho_{j_{0}}^{j}(\lambda,l)=(\Delta_{12}B_{l,j_0}^j)\,(b_{l}^j)^2\,\int_0^1\chi_1^\prime\left(b_{l}^j\left[(1-\tau)B_{l,j_0}^{j,[1]}+\tau B^{j,[2]}_{l,j_0}\right]\right)d\tau.
\end{align}
By \eqref{Poki1} combined with Lemma \ref{Compos-lemm-VM} we deduce from straightforward computations that
\begin{align}\label{sanin11}
(b_{l}^j)^2\left\|\chi_1^\prime\left(b_{l}^j\left[(1-\tau)B^{j,[1]}_{l,j_0}+\tau B^{j,[2]}_{l,j_0}\right]\right)\right\|^{q,\kappa}\lesssim \kappa^{-2}\langle l,j-j_0\rangle^{q+\tau_2(q+2)}.
\end{align}
Next, we shall move to the estimate of $\Delta_{12}B_{l,j_0}^j$. For this purpose we use \eqref{Spect-T2} leading to 
 \begin{align}\label{Spect-TMs2}
\mu_{j}^{m}=\mu_{j}^{0}+\sum_{n=0}^{m-1}\big \langle P_{N_n}\mathscr{R}_n{\bf e}_{0,j}, {\bf e}_{0,j}\big\rangle_{L^2(\T^{d+1})}.
\end{align}
We recall from Proposition \ref{projection in the normal directions} and the Gamma quotient introduced in Lemma \ref{lem-asym} 
$$\mu_{j}^{0}(\lambda,i_{0})\triangleq \Omega_{j}(\alpha)+jr^{1}(\lambda,i_{0})-\mathtt{W}(j,\alpha)r^{2}(\lambda,i_{0}).
$$
Therefore
\begin{align}\label{Spect-TMp2} 
\Delta_{12}\mu_{j}^{m}=\Delta_{12}\mu_{j}^{0}+\sum_{n=0}^{m-1}\big \langle P_{N_n}\Delta_{12}\mathscr{R}_n{\bf e}_{0,j}, {\bf e}_{0,j}\big\rangle_{L^2(\T^{d+1})}.
\end{align}
Combining \eqref{g-lj-8} with \eqref{Spect-TMp2} and \eqref{Spect-TMs2} allows to get
\begin{align*}
\Delta_{12}B_{l,j_0}^j&=\Delta_{12}(\mu_{j}^m-\mu_{j_{0}}^m)\\
&=(j-j_0)\Delta_{12} r^1-\big(\mathtt{W}(j,\alpha)-\mathtt{W}(j_0,\alpha)\big)\Delta_{12}r^2\\
&\quad+\sum_{n=0}^{m-1}\big \langle P_{N_n}\Delta_{12}\mathscr{R}_n{\bf e}_{0,j}, {\bf e}_{0,j}\big\rangle_{L^2(\T^{d+1})}-\sum_{n=0}^{m-1}\big \langle P_{N_n}\Delta_{12}\mathscr{R}_n{\bf e}_{0,j_0}, {\bf e}_{0,j_0}\big\rangle_{L^2(\T^{d+1})}.
\end{align*}
It follows from Proposition \ref{projection in the normal directions}-(i), Lemma \ref{Lem-Rgv1}-(i) and Lemma \ref{Stirling-formula}
\begin{align}\label{visitpl1}
\left\|\Delta_{12}B_{l,j_0}^j\right\|^{q,\kappa}
&\lesssim |j-j_0|  {\varepsilon \kappa^{-1}}\| \Delta_{12}i\|_{\overline{s}_h+\sigma}^{q,\kappa}+\sum_{n=0}^{m-1}\interleave P_{N_{n}}\Delta_{12}\mathscr{R}_{n}\interleave_{s_0}^{q,\kappa}.
\end{align}
Combining  \eqref{II23}, \eqref{sanin11}, \eqref{visitpl1} together with the law products in  Lemma \ref{Law-prodX1} we find
\begin{align}\label{IIM23}
\|\Delta_{12}\varrho_{j_{0}}^{j}(\cdot,l)\|^{q,\kappa}\lesssim &\varepsilon\kappa^{-3}\langle l,j-j_0\rangle^{q+\tau_2(q+2)+1} \| \Delta_{12}i\|_{\overline{s}_h+\sigma}^{q,\kappa}\\
\nonumber &+\kappa^{-2}\langle l,j-j_0\rangle^{q+\tau_2(q+2)}\sum_{n=0}^{m-1}\interleave P_{N_{n}}\Delta_{12}\mathscr{R}_{n}\interleave_{s_0}^{q,\kappa}. 
\end{align}
By the law products(or Leibniz rule) we find from \eqref{Ext-psi-opg}
\begin{align*}
\|\Delta_{12}\Psi_{n}^{n+k}(\cdot,l)\|_{W^{q,\infty}}&\lesssim  \|\Delta_{12}\varrho_{n}^{n+k}(\cdot,l)\|_{W^{q,\infty}}\max_{\ell=1,2}\|(r_{n}^{n+k})^{[\ell]}(\cdot,l)\|_{W^{q,\infty}}\\
&\quad+ \|\Delta_{12}r_{n}^{n+k}(\cdot,l)\|_{W^{q,\infty}}\max_{\ell=1,2}\|(\varrho_{n}^{n+k})^{[\ell]}(\cdot,l)\|_{W^{q,\infty}}.
\end{align*}
Putting together \eqref{IIM23}, \eqref{Tippal2} and using \eqref{small-C3}
\begin{align*}
\|\Delta_{12}\Psi_{n}^{n+k}(\cdot,l)\|_{W^{q,\infty}}&\lesssim \varepsilon \kappa^{-3-q}\,\langle l,k\rangle^{q+\tau_2(q+2)+1} \| \Delta_{12}i\|_{\overline{s}_h+\sigma}^{q,\kappa}\max_{\ell=1,2}\|(r_{n}^{n+k})^{[\ell]}(\cdot,l)\|_{W^{q,\infty}}\\
\nonumber &\quad +\kappa^{-q-2}\langle l,k\rangle^{q+\tau_2(q+2)}\sum_{n=0}^{m-1}\interleave P_{N_{n}}\Delta_{12}\mathscr{R}_{n}\interleave_{s_0}^{q,\kappa}\max_{\ell=1,2}\|(r_{n}^{n+k})^{[\ell]}(\cdot,l)\|_{W^{q,\infty}}\\
&\quad+\kappa^{-(q+1)} \langle l,k\rangle^{\tau_2(1+q)+q} \|\Delta_{12}r_{n}^{n+k}(\cdot,l)\|_{W^{q,\infty}}. 
\end{align*}
Notice that in the preceding estimate we can replace $q$ by any $|\gamma|\leqslant q$ and get after a slight modification the estimate 
\begin{align*}
 |\partial_\lambda^\gamma\Delta_{12}\Psi_{n}^{n+k}(\lambda,l)|&\lesssim \varepsilon \kappa^{-3-|\gamma|}\,\langle l,k\rangle^{|\gamma|+\tau_2(|\gamma|+2)+1} \| \Delta_{12}i\|_{\overline{s}_h+\sigma}^{q,\kappa}\max_{\ell=1,2}\max_{|\beta|\leqslant |\gamma|}|\partial_\lambda^\beta(r_{n}^{n+k})^{[\ell]}(\lambda,l)|\\
\nonumber &\quad+\kappa^{-|\gamma|-2}\langle l,k\rangle^{|\gamma|+\tau_2(|\gamma|+2)}\sum_{n=0}^{m-1}\interleave P_{N_{n}}\Delta_{12}\mathscr{R}_{n}\interleave_{s_0}^{q,\kappa}\max_{\ell=1,2}\max_{|\beta|\leqslant |\gamma|}|\partial_\lambda^\beta(r_{n}^{n+k})^{[\ell]}(\lambda,l)|\\
&\quad+\kappa^{-(|\gamma|+1)} \langle l,k\rangle^{\tau_2(1+|\gamma|)+|\gamma|} \max_{|\beta|\leqslant |\gamma|}|\partial_\lambda^\beta\Delta_{12}r_{n}^{n+k}(\lambda,l)|. \end{align*}
Inserting this estimate into \eqref{delta12psi} yields
\begin{align*}
  \interleave \Delta_{12}\Psi_m\interleave_{s}^{q,\kappa}&\lesssim  \varepsilon \kappa^{-3-q}\, \| \Delta_{12}i\|_{\overline{s}_h+\sigma}^{q,\kappa}\interleave P_{N_m}\mathscr{R}_m\interleave_{s+\tau_2(q+2)+1}^{q,\kappa}\\
  &\quad+\kappa^{-q-2}\interleave P_{N_m}\mathscr{R}_m\interleave_{s+\tau_2(q+2)}^{q,\kappa}\sum_{n=0}^{m-1}\interleave P_{N_{n}}\Delta_{12}\mathscr{R}_{n}\interleave_{s_0}^{q,\kappa}\\
  &\quad +\kappa^{-q-1}\interleave \Delta_{12}P_{N_m}\mathscr{R}_m\interleave_{s+\tau_2(1+q)}^{q,\kappa}
   .
\end{align*}
By setting
\begin{equation}\label{fouf-M1}
\varkappa_m(s)\triangleq \kappa^{-1}\sum_{n=0}^{m-1}\interleave \Delta_{12}\mathscr{R}_{n}\interleave_{s}^{q,\kappa}
\end{equation}
and using \eqref{Def-Taou}
we get successively  
\begin{align}\label{galeria-1}
\nonumber\|\Delta_{12}\Psi_{m}\interleave_{s_0}^{q,\kappa}&\lesssim  \varepsilon \kappa^{-2-q}\,N_m^{\tau_2+1} \| \Delta_{12}i\|_{\overline{s}_h+\sigma}^{q,\kappa}\delta_m(\overline{s}_l)+\kappa^{-q-1}N_m^{\tau_2}\delta_m(\overline{s}_l)\varkappa_m(s_0)\\
&\quad+\kappa^{-q-1}N_m^{\tau_2(1+q)}\interleave \Delta_{12}\mathscr{R}_m\interleave_{s_0}^{q,\kappa}
   \end{align}
and 
\begin{align}\label{Galeria-2}
\nonumber\interleave\Delta_{12}\Psi_{m}\interleave_{s_h}^{q,\kappa}&\lesssim \varepsilon \kappa^{-2-q}\,N_m^{\tau_2(q+2)+1}\delta_m(s_h) \| \Delta_{12}i\|_{\overline{s}_h+\sigma}^{q,\kappa}+\kappa^{-q-1}N_m^{\tau_2(q+2)}\delta_m({s}_h)\varkappa_m(s_0)\\
&\quad+\kappa^{-q-1}N_m^{\tau_2(1+q)}\interleave \Delta_{12}\mathscr{R}_m\interleave_{s_h}^{q,\kappa}.
\end{align}
   According to  \eqref{hypothesis of induction for deltamprime}, \eqref{whab1} and \eqref{small-C3} one obtains 
\begin{align}\label{Tippal2}
\delta_m(\overline{s}_l)\leqslant C \varepsilon \kappa^{-3}\,N_{0}^{\mu_{2}}N_{m}^{-\mu_{2}}\quad\hbox{and}\quad \delta_m(s_h)\leqslant C \varepsilon \kappa^{-3}\lesssim 1
\end{align}
Putting together \eqref{galeria-1} and \eqref{Tippal2} and using \eqref{small-C3} yields
\begin{align}\label{galeria-MM1}
\nonumber\|\Delta_{12}\Psi_{m}\interleave_{s_0}^{q,\kappa}&\lesssim  \varepsilon \kappa^{-2}\,N_m^{\tau_2+1-\mu_2} \| \Delta_{12}i\|_{\overline{s}_h+\sigma}^{q,\kappa}+\kappa^{-1}N_m^{\tau_2-\mu_2}\varkappa_m(s_0)\\
&\quad+\kappa^{-q-1}N_m^{\tau_2(1+q)}\interleave \Delta_{12}\mathscr{R}_m\interleave_{s_0}^{q,\kappa}.
   \end{align}
In a similar way, one gets by plugging \eqref{Tippal2} into \eqref{Galeria-2} 
\begin{align}\label{Galeria-22}
\nonumber\interleave\Delta_{12}\Psi_{m}\interleave_{s_h}^{q,\kappa}&\lesssim \varepsilon \kappa^{-2-q}\,N_m^{\tau_2(q+2)+1} \| \Delta_{12}i\|_{\overline{s}_h+\sigma}^{q,\kappa}+\varepsilon\kappa^{-q-4}N_m^{\tau_2(q+2)}\varkappa_m(s_0)\\
&\quad+\kappa^{-q-1}N_m^{\tau_2(1+q)}\interleave \Delta_{12}\mathscr{R}_m\interleave_{s_h}^{q,\kappa}\end{align}

Inserting \eqref{galeria-MM1} into \eqref{guiPM2} yields by virtue of  \eqref{hypothesis of induction for deltamprime} and \eqref{small-C3}
\begin{align}\label{guiPM3}
\nonumber \interleave\Delta_{12}\mathscr{R}_{m+1}\interleave_{s_0}^{q,\kappa}&\leqslant N_{m}^{s_0-s_h}   
\interleave\Delta_{12}\mathscr{R}_{m}\interleave_{s_h}^{q,\kappa}+C\varepsilon\kappa^{-1}N_m^{\tau_2(q+1)-2\mu_2}\,\|\Delta_{12}i\|_{{\overline{s}_h+\sigma}}^{q,\kappa}
\\
&\quad+C N_m^{\tau_2(q+1)-\mu_2}\|\Delta_{12}\mathscr{R}_{m}\interleave_{s_0}^{q,\kappa}+C\, N_m^{\tau_2-2\mu_2}\varkappa_m(s_0).
\end{align}
Similarly, by inserting \eqref{galeria-1}, \eqref{Galeria-2} into \eqref{guiPMS} and using \eqref{small-C3} we find
\begin{align}\label{guiPMSY}
\nonumber \interleave\Delta_{12}\mathscr{R}_{m+1}\interleave_{s_h}^{q,\kappa}&\leqslant \Big(1+C N_m^{-\mu_2+\tau_2(q+1)}+CN_{m}^{s_0-s_h+\tau_2(q+1)}\Big)\interleave\Delta_{12}\mathscr{R}_{m}\interleave_{s_h}^{q,\kappa}\\
\nonumber &\quad+C\varepsilon\kappa^{-3}N_m^{\tau_2(1+q)+q}\interleave\Delta_{12}\mathscr{R}_{m}\interleave_{s_0}^{q,\kappa}+C\varepsilon\kappa^{-3}\, N_m^{\tau_2(q+2)-\mu_2}\varkappa_m(s_0).\\
&\quad +C\varepsilon\kappa^{-1}\|\Delta_{12}i\|_{\overline{s}_h+\sigma}^{q,\kappa}N_m^{\tau_2(q+2)+1-\mu_2}.
\end{align}

Introduce 
$$\overline{\delta}_{m}(s)\triangleq \kappa^{-1}\interleave\Delta_{12}\mathscr{R}_{m}\interleave_{s}^{q,\kappa}.
$$
Then \eqref{guiPM3} and \eqref{guiPMSY} become
\begin{align}\label{guiPM4}
\nonumber \overline{\delta}_{m+1}(s_0)&\leqslant N_{m}^{s_0-s_h}   
\overline{\delta}_{m}(s_h)
+C\varepsilon\kappa^{-2}\interleave \Delta_{12}i\interleave_{{\overline{s}_h+\sigma}}^{q,\kappa}N_m^{\tau_2(q+1)-2\mu_2}\\
&\quad+C N_m^{\tau_2(q+1)-\mu_2}\overline{\delta}_{m}(s_0)+CN_m^{\tau_2-2\mu_2}\varkappa_{m}(s_0)
\end{align}
and
\begin{align}\label{guiPMSZ}
 \nonumber\overline{\delta}_{m+1}(s_h)&\leqslant \Big(1+C N_m^{-\mu_2+\tau_2(q+1)}+CN_{m}^{s_0-s_h+\tau_2(q+1)}\Big)\overline{\delta}_{m}(s_h)+C\varepsilon\kappa^{-3}N_m^{\tau_2(q+1)+q}\overline{\delta}_{m}(s_0)\\
&\quad+CN_m^{\tau_2(q+2)-\mu_2}\varkappa_{m}(s_0)+C\varepsilon\kappa^{-2}\|\Delta_{12}i\|_{{\overline{s}_h}+\sigma}^{q,\kappa}N_m^{\tau_2 (q+2)+1-\mu_2}.
\end{align}
It follows from  \eqref{fouf-M1} that for any $s\geqslant s_0$
\begin{align}\label{rec-dp-1}
\varkappa_{m}(s)\leqslant  \sum_{k=0}^{m-1} \overline{\delta}_{k}(s).
\end{align}
We shall prove by induction in $m\in\N$ that 
\begin{equation}\label{hypothesis of induction differences of mathscr Rm}
\forall\, k\leqslant m,\quad \overline{\delta}_{k}(s_{0})\leqslant N_{0}^{\mu_{2}}N_{k}^{-\mu_{2}}\nu(s_h)\quad \mbox{ and }\quad \overline{\delta}_{k}(s_{h})\leqslant\left(2-\frac{1}{k+1}\right)\nu(s_h),
\end{equation}
with
\begin{equation}\label{ITTL12}
\nu(s)\triangleq \overline\delta_0(s)+\varepsilon\kappa^{-2}\interleave\Delta_{12}i\interleave_{{\overline{s}_h+\sigma}}^{q,\kappa}.
\end{equation}
The validity of \eqref{hypothesis of induction differences of mathscr Rm} for $m=0$ is obvious from Sobolev embeddings. Now let us 
assume that the  property \eqref{hypothesis of induction differences of mathscr Rm} holds true at the order $m$ and let us check it at the order $m+1$. Then from \eqref{rec-dp-1} and Lemma \ref{lemma sum Nn} we find an absolute constant $C>0$ such that
$$
\varkappa_{m}(s_0)\leqslant C \nu(s_h).
$$
Combining this estimate with the induction assumption,  \eqref{guiPM4} and \eqref{guiPMSZ} implies
 \begin{align}
\nonumber\overline{\delta}_{m+1}(s_0)&\leqslant 2N_{m}^{s_0-s_h}   
\nu(s_{h})
+C  N_m^{\tau_2(q+1)-2\mu_2}\nu(s_h)+CN_0^{\mu_2} N_m^{\tau_2(q+1)-2\mu_2}\nu(s_h)\\
&\leqslant 2N_{m}^{s_0-s_h}   
\nu(s_{h})
+CN_0^{\mu_2} N_m^{\tau_2(q+1)-2\mu_2}\nu(s_h)
\end{align}
and
\begin{align*}
 \overline{\delta}_{m+1}(s_h)&\leqslant \Big(1+C N_m^{-\mu_2+\tau_2(q+1)}+CN_{m}^{s_0-s_h+\tau_2(q+1)}\Big)\left(2-\frac{1}{m+1}\right)\nu(s_{h})\\
&\quad+CN_m^{\tau_2(q+1)+q-\mu_2}\nu(s_h)+C\,N_m^{\tau_2(q+2)-\mu_2}\nu(s_h)+CN_m^{\tau_2(q+2)+1-\mu_2}\nu(s_h)\\
&\leqslant \Big(1+C N_m^{-\mu_2+\tau_2(q+1)}\Big)\left(2-\frac{1}{m+1}\right)\nu(s_{h})+CN_m^{\tau_2(q+2)+q-\mu_2}\nu(s_h)
\end{align*}
where in the last inequality we have used that $s_h\geqslant s_0+\mu_2$ accordng to \eqref{cond-diman1}.
Then with the choice of $\mu_2$  made in \eqref{cond-diman1}  we deduce in a standard way  that
\begin{align*}
\overline{\delta}_{m+1}(s_0)&\leqslant N_0^{\mu_2}N_{m+1}^{-\mu_2}\nu(s_h)\quad\hbox{and}\quad \overline{\delta}_{m+1}(s_h)\leqslant \left(2-\frac{1}{m+2}\right)\nu(s_h)
\end{align*}
which achieves the induction. Notice that one can fix in \eqref{hypothesis of induction differences of mathscr Rm} the free numbers  $\mu_2$ and $s_h$ at their lower bounds given in \eqref{cond-diman1}, that is, in view of the notations \eqref{Conv-T2N},
\begin{align*}
\mu_2&=\overline{\mu}_2+2\tau_2(1+q)\triangleq \mu_c\\
s_h&=\tfrac32\overline{\mu}_2+s_l+1+4\tau_2(1+q)=\overline{s}_h+4\tau_2(1+q)\triangleq s_c.
\end{align*} 
Therefore, one  finds  in particular
\begin{align}\label{low-es-3}
\overline{\delta}_{m}(s_{0})\leqslant N_{0}^{{\mu}_{c}}N_{m}^{-\mu_{c}}\nu({s}_c).
\end{align}
The next task is to estimate $\Delta_{12}r_{j}^{\infty}$. Then proceeding as for \eqref{dual-1X} we get by  using  
a duality argument, Lemma \ref{Lem-Rgv1}, \eqref{low-es-3} and Lemma \ref{lemma sum Nn} 
\begin{align*}
\|\Delta_{12}r_{j}^{\infty}\|^{q,\kappa}&\leqslant\sum_{m=0}^{\infty}\|\langle (P_{N_{m}}\Delta_{12}\mathscr{R}_{m})\mathbf{e}_{0,j},\mathbf{e}_{0,j}\rangle\|^{q,\kappa}\\
&\leqslant C \sum_{m=0}^{\infty}\interleave\Delta_{12}\mathscr{R}_{m}\interleave_{s_0}^{q,\kappa}\leqslant C \kappa \nu({s}_c)\sum_{m=0}^{\infty} N_{0}^{{\mu}_{c}}N_{m}^{-\mu_{c}}\\
& \leqslant C \kappa\,\nu({s}_c).
\end{align*}
Thus  we get from \eqref{ITTL12}
 \begin{align}\label{oppa1}
 \|\Delta_{12}r_{j}^{\infty}\|^{q,\kappa}
&\lesssim  \interleave\Delta_{12}\mathscr{R}_{0}\interleave_{{s}_c}^{q,\kappa}+\varepsilon\kappa^{-1}\interleave\Delta_{12}i\interleave_{{\overline{s}_h+\sigma}}^{q,\kappa}.
\end{align}
Combining   \eqref{dekomp} with   Proposition \ref{projection in the normal directions}-(i) and  using Sobolev embeddings we obtain
 \begin{align}\label{oppa2}
 \nonumber\forall j\in \mathbb{S}_0^c,\quad \|\Delta_{12}\mu_{j}^{\infty}\|^{q,\kappa}&\leqslant\ \|\Delta_{12}\mu_{j}^{0}\|^{q,\kappa}+ \|\Delta_{12}r_{j}^{\infty}\|^{q,\kappa}\\
&\lesssim\interleave\Delta_{12}\mathscr{R}_{0}\interleave_{{s}_c}^{q,\kappa}+\varepsilon\kappa^{-1}|j|\|\Delta_{12}i\|_{{\overline{s}_h}+\sigma}^{q,\kappa}.
\end{align}
Recall that the operator $\mathscr{R}_{0}$ coincides with $\mathcal{R}^2_{r,\lambda}$ defined in Proposition \ref{projection in the normal directions}. Then from Proposition \ref{projection in the normal directions}-(iii) and \eqref{small-C3} we infer
\begin{equation*}
\interleave \Delta_{12}\mathscr{R}_{0}\interleave_{0,0,0}^{q,\kappa}\lesssim \varepsilon{\kappa^{-2}}\|\Delta_{12}i\|_{\overline{s}_h+\sigma}^{q,\kappa}
\end{equation*}
and
\begin{equation*}
\interleave\mathscr{R}_{0}\interleave_{0,s_h,0}^{q,\kappa}\lesssim\varepsilon{\kappa^{-2}}.
\end{equation*}
From the latter inequality we easily  get by  applying the triangle inequality, since $s_h\geqslant 2 {s}_c$
\begin{equation*}
\interleave\Delta_{12}\mathscr{R}_{0}\interleave_{0,2{s}_c,0}^{q,\kappa}\lesssim\varepsilon{\kappa^{-2}}.
\end{equation*}
Applying the interpolation inequality of Lemma \ref{Lem-Rgv1}-(v) we get
\begin{align*}
\interleave\Delta_{12}\mathscr{R}_{0}\interleave_{0,{s}_c,0}^{q,\kappa}&\lesssim\big( \interleave\Delta_{12}\mathscr{R}_{0}\interleave_{0,0,0}^{q,\kappa}\big)^{\frac12}\big(\interleave\Delta_{12}\mathscr{R}_{0}\interleave_{0,2{s}_c,0}^{q,\kappa}\big)^{\frac12}\\
&\lesssim \varepsilon {\kappa^{-2}}\big(\|\Delta_{12}i\|_{\overline{s}_h+\sigma}^{q,\kappa}\big)^{\frac12}.
\end{align*}
Plugging this estimate into \eqref{oppa1} and \eqref{oppa2} yields successively, 
\begin{align*}
\forall j\in \mathbb{S}_0^c,\quad  \|\Delta_{12}r_{j}^{\infty}\|^{q,\kappa}
&\lesssim  \varepsilon {\kappa^{-2}}\big(\|\Delta_{12}i\|_{\overline{s}_h+\sigma}^{q,\kappa}\big)^{\frac12}+\varepsilon\kappa^{-1}\interleave\Delta_{12}i\interleave_{{\overline{s}_h+\sigma}}^{q,\kappa}\\
&\lesssim  \varepsilon {\kappa^{-2}}\big(\|\Delta_{12}i\|_{\overline{s}_h+\sigma}^{q,\kappa}\big)^{\frac12}
\end{align*}
and
 \begin{align*}
 \nonumber\forall j\in \mathbb{S}_0^c,\quad \|\Delta_{12}\mu_{j}^{\infty}\|^{q,\kappa}&\lesssim\varepsilon {\kappa^{-2}}\big(\|\Delta_{12}i\|_{\overline{s}_h+\sigma}^{q,\kappa}\big)^{\frac12}+\varepsilon\kappa^{-1}|j|\interleave\Delta_{12}i\interleave_{{\overline{s}_h+\sigma}}^{q,\kappa}\\
 &\lesssim \varepsilon {\kappa^{-2}}\big(\|\Delta_{12}i\|_{\overline{s}_h+\sigma}^{q,\kappa}\big)^{\frac12}|j|.
\end{align*}
This ends the proof of Proposition \ref{reduction of the remainder term}.}
\end{proof}

\subsubsection{Approximate  inverse in the normal direction}\label{Inverse-est-Op}
The main concern of this section is to find an approximate right inverse of the linearized operator $\widehat{\mathcal{L}}_{\omega}$ defined \mbox{in \eqref{Norm-proj}} and detailed in \mbox{Proposition \ref{lemma-GS0}} when the set of parameters is restricted to a suitable Cantor like set obtained through successive excisions.  The main result reads as follows.
\begin{theorem}\label{inversion of the linearized operator in the normal directions}
Let $(\gamma,q,d,\tau_{1},\tau_2,s_{0},s_{h},S)$ satisfy  \eqref{Conv-T2}, \eqref{Conv-T2N}, \eqref{cond-diman1} and assume the smallness condition \eqref{small-C3}. Then the following assertions hold true.
\begin{enumerate}
\item Consider the operator  $\mathscr{L}_{\infty}$ defined in Proposition $\ref{reduction of the remainder term},$ then there exists a family of  linear operators $\big(\mathtt{T}_n\big)_{n\in\mathbb{N}}$  defined in the whole set $\mathcal{O}$ with  the estimate
$$
\forall \, s\in[s_0,S],\quad \sup_{n\in\mathbb{N}}\|\mathtt{T}_{n}\rho\|_{s}^{q,\kappa }\lesssim \kappa ^{-1}\|\rho\|_{s+\tau_{1}q+\tau_{1}}^{q,\kappa }$$
and such that  for any $n\in\mathbb{N}$ and  
 in the Cantor set
$$\Lambda_{\infty,{n}}^{\kappa,\tau_{1}}(i_{0})=\bigcap_{(l,j)\in\mathbb{Z}^{d }\times\mathbb{S}_{0}^{c}}\left\lbrace\lambda\in\mathcal{O};\;\,\forall\, |l|\leqslant N_n,\,j\in\mathbb{S}_0^c,\,|\omega\cdot l+\mu_{j}^{\infty}(\lambda,i_0)|>\tfrac{\kappa {\langle j\rangle } }{\langle l\rangle^{\tau_{1}}}\right\rbrace
$$
we have
$$
\mathscr{L}_{\infty}\mathtt{T}_n=\textnormal{Id}+{\mathtt{E}^4_n}
$$
with 
$$
\forall \, s\in\,[ s_0, S],\quad \quad \|\mathtt{E}^4_{n}h\|_{s_0}^{q,\kappa }\lesssim 
 N_n^{s_0-s}\,\kappa^{-1}\,\|h\|_{s+\tau_{1}q+\tau_{1}+1}^{q,\kappa }.
 $$
\item 
There exist $\sigma_5=\sigma(\tau_1,\tau_2,d,q)\geqslant \sigma_4$ and  a family of linear  operators $(\mathtt{T}_{\omega,n})_n$ satisfying
\begin{equation}\label{estimate mathcalTomega}
\forall \, s\in\,[ s_0, S] ,\quad\sup_{n\in\mathbb{N}}\|\mathtt{T}_{\omega,n}h\|_{s}^{q,\kappa }\lesssim\kappa^{-1}\left(\|h\|_{s+{{\sigma_5}}}^{q,\kappa }+\| \mathfrak{I}_{0}\|_{s+{{\sigma_5}}}^{q,\kappa }\|h\|_{s_{0}+{{\sigma_5}}}^{q,\kappa }\right)
\end{equation}
and such that  
 in the Cantor set
$$\mathtt{G}_n(\kappa,\tau_{1},\tau_{2},i_{0})\triangleq\mathcal{O}_{\infty,n}^{\kappa,\tau_1,\tau_{2}}(i_{0})\cap\Lambda_{\infty,n}^{\kappa,\tau_{1}}(i_{0})$$
we have
$$
\widehat{\mathcal{L}}_{\omega}\mathtt{T}_{\omega,n}=\textnormal{Id}+{\mathtt{E}_n},
$$
with the estimates
\begin{align*}
\nonumber  \forall\, s\in [s_0,S],\quad\|\mathtt{E}_nh\|_{s_0}^{q,\kappa}
 \nonumber\lesssim N_n^{s_0-s}\kappa^{-1}\Big( \|h\|_{s+\sigma_5}^{q,\kappa}&+{\varepsilon\kappa^{-3}}\| \mathfrak{I}_{0}\|_{s+\sigma_5}^{q,\kappa}\|h\|_{s_0+{\sigma_5}}^{q,\kappa} \Big)\\
 &+ \varepsilon\kappa^{-4}N_{0}^{{\mu}_{2}}{N_{n}^{-\mu_{2}}} \|h\|_{s_0+\sigma_5}^{q,\kappa},
\end{align*}
where  $\widehat{\mathcal{L}}_{\omega},  \mathcal{O}_{\infty,n}^{\kappa,\tau_{1}}(i_{0})$ and $\mathcal{O}_{\infty,n}^{\kappa,\tau_1,\tau_{2}}(i_{0})$ are defined in \eqref{Norm-proj} and Proposition $\ref{QP-change}$ and Proposition $\ref{reduction of the remainder term}$, respectively.
\item  On the Cantor set $\mathtt{G}_n(\kappa,\tau_{1},\tau_{2},i_{0})$, we have the splitting
$$
\widehat{\mathcal{L}}_{\omega}=\widehat{\mathtt{L}}_{\omega,n}+\widehat{\mathtt{R}}_n,\quad\hbox{with}\quad \widehat{\mathtt{L}}_{\omega,n}\mathtt{T}_{\omega,n}=\textnormal{Id}\quad\hbox{and}\quad \widehat{\mathtt{R}}_n={\mathtt{E}_n}\widehat{\mathtt{L}}_{\omega}
$$
where the operators $\widehat{\mathtt{L}}_{\omega,n}$ and $\widehat{\mathtt{R}}_n$ are defined in the whole set $\mathcal{O}$ with the estimates
$$
 \forall\, s\in [s_0,S],\quad\sup_{n\in\mathbb{N}}\|\widehat{\mathtt{L}}_{\omega,n} h\|_{s}^{q,\kappa}\lesssim \|h\|_{s+1}^{q,\kappa}{+{{\varepsilon\kappa^{-3}}\| \mathfrak{I}_{0}\|_{q,s+\sigma_5}^{q,\kappa }\|h\|_{s_{0}+1}^{q,\kappa }}}
$$
and
\begin{align*}
\nonumber  \|\widehat{\mathtt{R}}_nh\|_{s_0}^{q,\kappa}
 \nonumber\lesssim N_n^{s_0-s}\kappa^{-1}\Big( \|h\|_{s+\sigma_5}^{q,\kappa_5}&+{\varepsilon\kappa^{-3}}\| \mathfrak{I}_{0}\|_{s+\sigma_5}^{q,\kappa}\|h\|_{s_0+\sigma_5}^{q,\kappa} \Big)+ \varepsilon\kappa^{-4}N_{0}^{{\mu}_{2}}{N_{n}^{-\mu_{2}}} \|h\|_{s_0+\sigma_5}^{q,\kappa}.
\end{align*}

\end{enumerate}
\end{theorem}
\begin{proof}
{\bf{(i)}} From   Proposition \ref{reduction of the remainder term} we recall that
\begin{align*}
\mathscr{L}_{\infty}(\lambda,\omega,i_{0})&=\big(\omega\cdot\partial_{\varphi}+\mathscr{D}_{\infty}(\lambda,i_{0})\big)\Pi_{\mathbb{S}_0}^{\perp}
.
\end{align*}
Then we may split this operator as follows
\begin{align}\label{Tikl1}
\nonumber\mathscr{L}_{\infty}(\lambda,\omega,i_{0})&=\Pi_{N_n}\omega\cdot\partial_{\varphi}\Pi_{N_n}\Pi_{\mathbb{S}_0}^{\perp}+\mathscr{D}_{\infty}(\lambda,i_{0})\Pi_{\mathbb{S}_0}^{\perp}-\Pi_{N_n}^\perp\omega\cdot\partial_{\varphi}\Pi_{N_n}^\perp\Pi_{\mathbb{S}_0}^{\perp}\\
&\triangleq\mathtt{L}_n-\mathtt{R}_n,
\end{align}
with $\mathtt{R}_n=\Pi_{N_n}^\perp\omega\cdot\partial_{\varphi}\Pi_{N_n}^\perp\Pi_{\mathbb{S}_0}^{\perp}$ and the projector $\Pi_{N_n}$ is defined by
$$
\Pi_{N}\sum_{l,j}h_{l,j}{\bf e}_{l,j}=\sum_{|l|\leqslant N\atop j\in\Z}h_{l,j}{\bf e}_{l,j}
$$
 From this definition and the structure of $\mathscr{D}_{\infty}(\lambda,i_{0})$ in Proposition \ref{reduction of the remainder term} we deduce that
$${\bf e}_{-l,-j}\mathtt{L}_n{\bf e}_{l,j}=\left\lbrace\begin{array}{rcl}
 \ii\big(\omega\cdot l+\mu_j^\infty\big);& \hbox{if}& \quad |l|\leqslant N_n,\quad j\in\mathbb{S}_0^c\\
  \ii\,\mu_j^\infty;& \hbox{if} & \quad |l|> N_n,\quad j\in\mathbb{S}_0^c.
\end{array}\right.$$
Define the diagonal  operator  $\mathtt{T}_n$ by 
\begin{eqnarray*}
\mathtt{T}_{n}h(\lambda,\varphi,\theta)&\triangleq&
-\ii\sum_{|l|\leqslant N_n \atop 
  j\in\mathbb{S}_{0}^{c}}\tfrac{\chi\left((\omega\cdot l+\mu_{j}^{\infty}(\lambda))\kappa ^{-1}\langle l\rangle^{\tau_{1}}\right)}{\omega\cdot l+\mu_{j}^{\infty}(\lambda)}h_{l,j}(\lambda)\,{\bf e}_{l,j}(\varphi,\theta)\\
&&-\ii\sum_{|l|>N_n\atop   j\in\mathbb{S}_{0}^{c}}\frac{h_{l,j}(\lambda)}{\mu_{j}^{\infty}(\lambda)}\,{\bf e}_{l,j}(\varphi,\theta)
\end{eqnarray*}
where $\chi$ is the cut-off function defined in \eqref{chi-def-1} and $\left(h_{l,j}(\lambda)\right)$ are  the Fourier coefficients of $h$. We recall from Proposition \ref{reduction of the remainder term} that
$$\mu_{j}^{\infty}(\lambda)=\mu_{j}^{0}(\lambda)+r_{j}^{\infty}(\lambda)\quad \mbox{ with }\quad 
\sup_{j\in\mathbb{S}_{0}^{c}}\| r_{j}^{\infty}\|^{q,\kappa}\leqslant C \varepsilon\kappa^{-2}.
$$
Combining this estimate with  Proposition \ref{projection in the normal directions}-(i) we find 
$$
\forall\, j\in \mathbb{S}_{0}^{c},\, \|\mu^\infty_{j}\|^{q,\kappa}\leqslant C |j|.
$$
Similar arguments based on Lemma \ref{lem-asym}-(iv) give under the smallness condition \eqref{small-C3}
$$|j|\lesssim \|\mu^\infty_{j}\|^{0,\kappa }\leqslant  \|\mu^\infty_{j}\|^{q,\kappa}.
$$
Then implementing in part Lemma \ref{L-Invert} yields
\begin{align}\label{Es-Mu1}
\forall s\geqslant s_0,\quad  \|\mathtt{T}_{n}h\|_{s}^{q,\kappa}\lesssim \kappa^{-1}\|h\|_{s+\tau_{1}q+\tau_{1}}^{q,\kappa}.
\end{align}
Moreover, by construction we get
 \begin{align}\label{Inv-Ty1}
 \forall \lambda\in \Lambda_{\infty,n}^{\kappa,\tau_{1}}(i_{0}),\quad\mathtt{L}_{n}\mathtt{T}_n=\textnormal{Id},
 \end{align}
 since $\chi(\cdot)=1$ in the set $ \Lambda_{\infty,n}^{\kappa,\tau_{1}}(i_{0})$. It follows according to  \eqref{Tikl1} that 
\begin{align}\label{Inv-op-1}
\nonumber\forall\, \lambda\in \Lambda_{\infty,n}^{\kappa,\tau_{1}}(i_{0}),\quad \mathscr{L}_\infty \mathtt{T}_n&=\textnormal{Id}-\mathtt{R}_n\mathtt{T}_n\\
&\triangleq \textnormal{Id}+\mathtt{E}_n^4.
\end{align}
Notice that one gets from direct computations that
$$
 \forall \, s_0\leqslant s\leqslant \overline{s},\quad \|\mathtt{R}_{n}h\|_{s}^{q,\kappa}\lesssim N_n^{s-\overline{s}}\|h\|_{\overline{s}+1}^{q,\kappa}.
 $$
 Combining this estimate with \eqref{Es-Mu1} yields
\begin{align}\label{Mila-1}
\nonumber  \forall \, s_0\leqslant s\leqslant \overline{s},\quad \|\mathtt{E}^4_{n}h\|_{s}^{q,\kappa}&\lesssim N_n^{s-\overline{s}}\|\mathtt{T}_nh\|_{\overline{s}+1}^{q,\kappa}\\
 &\lesssim 
 N_n^{s-\overline{s}}\kappa^{-1}\|h\|_{\overline{s}+\tau_{1}q+\tau_{1}+1}^{q,\kappa}.
 \end{align}
 This achieves the proof of the first point.
 
 \smallskip
 
{{\bf{(ii)}} }Let us define  
\begin{align}\label{tomega}
\mathtt{T}_{\omega,n}\triangleq \Psi_{\perp}\Phi_{\infty}\mathtt{T}_{n}\Phi_{\infty}^{-1}\Psi_{\perp}^{-1}.
\end{align}
where the operators $ \Psi_{\perp}$ and $\Phi_{\infty}$ are defined in Proposition \ref{projection in the normal directions} and Proposition \ref{reduction of the remainder term}, respectively. Notice that $\mathtt{T}_{\omega,n}$ is defined in the whole range of parameters $\mathcal{O}$. 
{
Applying Lemma \ref{lemm-decompY}-(ii) combined with \eqref{small-C3} yields 
$$\forall s\geqslant s_{0},\quad\|\mathtt{T}_{\omega,n}h\|_{s}^{q,\kappa}\lesssim\|\Phi_{\infty}\mathtt{T}_{n}\Phi_{\infty}^{-1}\Psi_{\perp}^{-1}h\|_{s}^{q,\kappa}+\|\mathfrak{I}_{0}\|_{s+\sigma}^{q,\kappa}\|\Phi_{\infty}\mathtt{T}_{n}\Phi_{\infty}^{-1}\Psi_{\perp}^{-1}h\|_{s_0}^{q,\kappa}.$$
By using \eqref{estimate on Phiinfty and its inverse} and \eqref{small-C3}  one gets
$$\forall s\in[s_{0},S],\quad\|\Phi_{\infty}\mathtt{T}_{n}\Phi_{\infty}^{-1}\Psi_{\perp}^{-1}h\|_{s}^{q,\kappa}\lesssim\|\mathtt{T}_{n}\Phi_{\infty}^{-1}\Psi_{\perp}^{-1}h\|_{s}^{q,\kappa}+\|\mathfrak{I}_{0}\|_{q,s+{\sigma}}^{q,\kappa}\|\mathtt{T}_{n}\Phi_{\infty}^{-1}\Psi_{\perp}^{-1}h\|_{s_0}^{q,\kappa}.$$
From \eqref{Es-Mu1}  we find
$$\forall s\geqslant s_{0},\quad\|\mathtt{T}_{n}\Phi_{\infty}^{-1}\Psi_{\perp}^{-1}h\|_{s}^{q,\kappa}\lesssim\kappa^{-1}\|\Phi_{\infty}^{-1}\Psi_{\perp}^{-1}h\|_{s+\tau_{1}q+\tau_{1}}^{q,\kappa}.$$
Applying  \eqref{estimate on Phiinfty and its inverse} and \eqref{small-C3} yields
$$\forall s\in[s_{0},S],\quad\|\Phi_{\infty}^{-1}\Psi_{\perp}^{-1}h\|_{s}^{q,\kappa}\lesssim\|\Psi_{\perp}^{-1}h\|_{s}^{q,\kappa}+\|\mathfrak{I}_{0}\|_{s+\overline{\sigma}}^{q,\kappa}\|\Psi_{\perp}^{-1}h\|_{s_0}^{q,\kappa}.$$
Putting  together the preceding  three  estimates with Proposition \ref{lemm-decompY}-(ii) then  we get \eqref{estimate mathcalTomega}.}\\
Now combining Propositions \ref{projection in the normal directions}-\ref{reduction of the remainder term} we find that  that in the Cantor set  $\mathcal{O}_{\infty,n}^{\kappa,\tau_{1}}(i_{0})\cap\mathscr{O}_{\infty,n}^{\kappa,\tau_{2}}(i_{0})$
\begin{align*}
\Phi^{-1}_{\infty}\Psi_{\perp}^{-1}\widehat{\mathcal{L}}_\omega\Psi_{\perp}\Phi_{\infty}&=\Phi^{-1}_{\infty}\mathscr{L}_{0}\Phi_{\infty}+\Phi^{-1}_{\infty}\mathtt{E}_{n}^2\Phi_{\infty}\\
&=\mathscr{L}_{\infty}+\mathtt{E}_n^3+\Phi^{-1}_{\infty}\mathtt{E}_{n}^2\Phi_{\infty}.
\end{align*}
It follows that in the Cantor set  $\mathtt{G}_n(\kappa,\tau_{1},\tau_{2},i_{0})=\mathcal{O}_{\infty,n}^{\kappa,\tau_{1},\tau_2}(i_{0})\cap \Lambda_{\infty,n}^{\kappa,\tau_{1}}(i_{0})$ one has  by virtue of  \eqref{Inv-op-1} and the identity \eqref{tomega}
\begin{align*}
\Phi^{-1}_{\infty}\Psi_{\perp}^{-1}\widehat{\mathcal{L}}_\omega\Psi_{\perp}\Phi_{\infty}\mathtt{T}_n
&=\textnormal{Id}+\mathtt{E}_n^4+\mathtt{E}_n^3\mathtt{T}_n+\Phi^{-1}_{\infty}\mathtt{E}_{n}^2\Phi_{\infty}\mathtt{T}_n\,\end{align*}
which gives in view of \eqref{tomega}  the identity 
\begin{align}\label{Id-Moh1}
\nonumber\widehat{\mathcal{L}}_\omega\mathtt{T}_{\omega,n}&=\textnormal{Id}+\Psi_{\perp}\Phi_{\infty}\big(\mathtt{E}_n^4+\mathtt{E}_n^3\mathtt{T}_n+\Phi^{-1}_{\infty}\mathtt{E}_{n}^2\Phi_{\infty}\mathtt{T}_n\big)\Phi^{-1}_{\infty}\Psi_{\perp}^{-1}
\\
\nonumber&\triangleq\textnormal{Id}+\Psi_{\perp}\Phi_{\infty}\mathtt{E}_n^5\Phi^{-1}_{\infty}\Psi_{\perp}^{-1}\\
&\triangleq\textnormal{Id}+\mathtt{E}_n
\end{align}
 provided that  $\lambda\in \mathtt{G}_n(\kappa,\tau_{1},\tau_{2},i_{0})$. The  estimate of the first term of  $\mathtt{E}_n^5$ is given in \eqref{Mila-1}.
  For the second term of $\mathtt{E}_n^5$ we use the estimate  \eqref{Error-Est-2D} combined with \eqref{Es-Mu1} leading to
\begin{align}\label{Error-Est-2LD}
\nonumber \|{\mathtt{E}^3_n\mathtt{T}_n}h\|_{s_0}^{q,\kappa}&\lesssim\,\varepsilon\kappa^{-3}N_{0}^{{\mu}_{2}}{N_{n}^{-\mu_{2}}} \|\mathtt{T}_nh\|_{s_0+1}^{q,\kappa}\\
&\lesssim \,\varepsilon\kappa^{-4}N_{0}^{{\mu}_{2}}{N_{n}^{-\mu_{2}}} \|h\|_{s_0+1+\tau_1(1+q)}^{q,\kappa}.
\end{align}
To estimate  $\Phi^{-1}_{\infty}\mathtt{E}_{n}^2\Phi_{\infty}\mathtt{T}_n$ we combine Proposition \ref{projection in the normal directions}-(ii) with \eqref{estimate on Phiinfty and its inverse}, \eqref{Es-Mu1}  and \eqref{small-C3} 
\begin{align}\label{Error-Est-2LLD}
\nonumber \|\Phi^{-1}_{\infty}\mathtt{E}_{n}^2\Phi_{\infty}\mathtt{T}_nh\|_{s_0}^{q,\kappa}&\leqslant C\varepsilon\kappa^{-2}\Big(N_{0}^{\mu_{2}}N_{n}^{-\mu_{2}}+  N_{n}^{s-s_0}\big(1+\| \mathfrak{I}_{0}\|_{s+\sigma}^{q,\kappa}\big)\Big)\| \mathtt{T}_n h\|_{s_0+{3}}^{q,\kappa}\\
&\leqslant C\varepsilon\kappa^{-3}\Big(N_{0}^{\mu_{2}}N_{n}^{-\mu_{2}}+  N_{n}^{s-s_0}\big(1+\| \mathfrak{I}_{0}\|_{s+\sigma}^{q,\kappa}\big)\Big) \|h\|_{s_0+\tau_1(1+q)+{3}}^{q,\kappa}.
\end{align}
Similarly we get from \eqref{Mila-1}, \eqref{Error-Est-2LD} and \eqref{Error-Est-2LLD}
\begin{align}\label{Fif-1}
 \|{\mathtt{E}^5_n}h\|_{s_0}^{q,\kappa}&\lesssim  \varepsilon\kappa^{-4}\Big(N_{0}^{\mu_{2}}N_{n}^{-\mu_{2}}+  N_{n}^{s-s_0}\big(1+\| \mathfrak{I}_{0}\|_{s+\sigma}^{q,\kappa}\big)\Big) \|h\|_{s_0+3+\tau_1(1+q)}^{q,\kappa}\\
 \nonumber&\quad +N_n^{s_0-s}\kappa^{-1}\|h\|_{{s}+1+\tau_{1}q+\tau_{1}}^{q,\kappa}.
\end{align}
Set $\overline\Psi\triangleq\Psi_{\perp}\Phi_{\infty}$ then from  \eqref{estimate on Phiinfty and its inverse}, Propostion \ref{lemm-decompY}-(ii) and using \eqref{small-C3}   we deduce that
\begin{align}\label{MLLD001}
\forall s\in[ s_{0},S],\quad \mbox{ }\|\overline\Psi^{\pm1}h\|_{s}^{q,\kappa}\lesssim \|h\|_{s}^{q,\kappa}+{\varepsilon\kappa^{-3}}\| \mathfrak{I}_{0}\|_{s+\sigma}^{q,\kappa}\|h\|_{s_0}^{q,\kappa}.
\end{align}
In particular, we get
\begin{align}\label{MLLD01}
\nonumber \mbox{ }\|\overline\Psi^{\pm1}h\|_{s_0}^{q,\kappa}&\lesssim \|h\|_{s_0}^{q,\kappa}+{\varepsilon\kappa^{-3}}\| \mathfrak{I}_{0}\|_{s_0+\sigma}^{q,\kappa}\|h\|_{s_0}^{q,\kappa}\\
&\lesssim \|h\|_{s_0}^{q,\kappa}.
\end{align}
Thus we deduce from  \eqref{MLLD01} combined with  \eqref{Fif-1} and  \eqref{small-C3} we deduce 
\begin{align*}
\nonumber \|\overline\Psi{\mathtt{E}^5_n}\overline\Psi^{-1}h\|_{s_0}^{q,\kappa}&\lesssim \,\,\|\mathtt{E}^5_n\overline\Psi^{-1}h\|_{s_0}^{q,\kappa}\\
 \nonumber&\lesssim  \varepsilon\kappa^{-4}\Big(N_{0}^{\mu_{2}}N_{n}^{-\mu_{2}}+  N_{n}^{s_0-s}\big(1+\| \mathfrak{I}_{0}\|_{s+\sigma}^{q,\kappa}\big)\Big) \|\overline\Psi^{-1}h\|_{s_0+3+\tau_1(1+q)}^{q,\kappa}\\
 \nonumber&\quad+N_n^{s_0-s}\kappa^{-1}\|\overline\Psi^{-1}h\|_{{s}+1+\tau_{1}q+\tau_{1}}^{q,\kappa}\\
 \nonumber& \lesssim\varepsilon\kappa^{-4}\Big(N_{0}^{\mu_{2}}N_{n}^{-\mu_{2}}+  N_{n}^{s_0-s}\big(1+\| \mathfrak{I}_{0}\|_{s+\sigma}^{q,\kappa}\big)\Big) \|h\|_{s_0+3+\tau_1(1+q)}^{q,\kappa}\\
 &\quad + N_n^{s_0-s}\kappa^{-1}\|h\|_{s+1+\tau_1(1+q)}^{q,\kappa}.
\end{align*}
Consequently we obtain from \eqref{small-C3}
\begin{align}\label{MLLDRP1}
  \sup_{n\in\mathbb{N}}\|\mathtt{E}_nh\|_{s_0}^{q,\kappa}
\lesssim N_n^{s_0-s}\kappa^{-1}\Big( \|h\|_{s+\sigma}^{q,\kappa}&+{\varepsilon\kappa^{-3}}\| \mathfrak{I}_{0}\|_{s+\sigma}^{q,\kappa}\|h\|_{s_0+{\sigma}}^{q,\kappa} \Big)+ \varepsilon\kappa^{-4}N_{0}^{{\mu}_{2}}{N_{n}^{-\mu_{2}}} \|h\|_{s_0+\sigma}^{q,\kappa}.
\end{align}
This achieves the proof of the second point.\\

{\bf{(iii)}} This result follows easily from the preceding point (ii). Indeed, according to \eqref{Id-Moh1} one may write on the Cantor set 
$\mathtt{G}_n(\kappa,\tau_{1},\tau_{2},i_{0})$
\begin{align}\label{Id-Moh11}
\widehat{\mathcal{L}}_\omega=&\mathtt{T}_{\omega,n}^{-1}+\mathtt{E}_n\mathtt{T}_{\omega,n}^{-1}
\end{align}
In addition, we get from \eqref{tomega} and \eqref{Inv-Ty1}
\begin{align}\label{tomega11}
\nonumber\mathtt{T}_{\omega,n}^{-1}&=\Psi_{\perp}\Phi_{\infty}\mathtt{L}_{n}\Phi_{\infty}^{-1}\Psi_{\perp}^{-1}\\
&\triangleq \widehat{\mathtt{L}}_{\omega,n}.
\end{align}
This provides the splitting  
$$
\widehat{\mathcal{L}}_\omega=\widehat{\mathtt{L}}_{\omega,n}+\widehat{\mathtt{R}}_{n}\quad\hbox{and}\quad \widehat{\mathtt{R}}_{n}\triangleq \mathtt{E}_n \widehat{\mathtt{L}}_{\omega,n} .
$$
From \eqref{Tikl1} and \eqref{MLLD001} combined with Proposition \ref{reduction of the remainder term} we obtain
$$
\|\widehat{\mathtt{L}}_{\omega,n} h\|_{s}^{q,\kappa}\lesssim \|h\|_{s+1}^{q,\kappa}+{\varepsilon\kappa^{-3}}\| \mathfrak{I}_{0}\|_{s+\sigma}^{q,\kappa}\|h\|_{s_0+1}^{q,\kappa}.
$$
Inserting this estimate  into \eqref{MLLDRP1} yields
\begin{align*}
 \nonumber \forall\, s\in [s_0,S],\quad \sup_{n\in\mathbb{N}}\|\widehat{\mathtt{R}}_{n}h\|_{q,s}^{q,\kappa } &\lesssim \kappa^{-1}\|h\|_{{s}+\sigma+1}^{q,\kappa }+{\varepsilon\kappa^{-4}}\| \mathfrak{I}_{0}\|_{q,s+\sigma}^{q,\kappa }\|h\|_{s_{0}+\sigma+1}^{q,\kappa }.
\end{align*}
and
\begin{align*}
  \sup_{n\in\mathbb{N}}\|\widehat{\mathtt{R}}_{n}h\|_{s_0}^{q,\kappa}
&\lesssim N_n^{s_0-s}\kappa^{-1}\Big( \|h\|_{s+\sigma+1}^{q,\kappa}+{\varepsilon\kappa^{-3}}\| \mathfrak{I}_{0}\|_{s+\sigma+1}^{q,\kappa}\|h\|_{s_0+\sigma+1}^{q,\kappa} \Big)\\
&\quad+ \varepsilon\kappa^{-4}N_{0}^{{\mu}_{2}}{N_{n}^{-\mu_{2}}} \|h\|_{s_0+\sigma+1}^{q,\kappa}.
\end{align*}
Thus changing  $\sigma+1$ with $\sigma$  gives the suitable result.
\end{proof}

\section{Proof of the main result}\label{N-M-section12}
The main concern of this section is to implement  Nash-Moser scheme in order to construct  zeros  for the nonlinear functional
${\mathcal F} $ defined in \eqref{operatorF} when $\varepsilon>0$ and  small enough. We shall prove that solutions do  exist provided that the  parameters $\lambda=(\omega,\alpha)$ belong to a final Cantor set $\mathtt{G}_{\infty}^{\kappa}.$ More precisely, we are able to construct  smooth  functions $\lambda\in\mathcal{O}\mapsto i_\infty( \lambda),\, \mathtt{c}_\infty( \lambda) $ in $W^{q,\infty}(\mathcal{O})$  such that
\begin{align}\label{eq-NM}
\forall\, \lambda\in \mathtt{G}_{\infty}^{\kappa},\quad \mathcal{F} (i_\infty( \lambda), \mathtt{c}_\infty( \lambda), \lambda,\varepsilon )=0.
\end{align}
The next stage is to find solutions to the origininal Hamiltonian equation \eqref{NLeq}. For this aim,  we should adjust the parameters in such a way that   $\mathtt{c}_\infty( \omega,  \alpha)=-{\omega}_{\textnormal{Eq}}(\alpha)$, where the latter quantity corresponds to the equilibrium frequency defined in \eqref{tan-nor}. This equation is invertible in $\omega$ close to the equilibrium frequencies and $\omega$ becomes an implicit function of $\alpha$. Consequently we generate solutions when the parameter $\alpha$ belongs to the  Cantor set
$$
\mathtt{C}_{\infty}^{\kappa,\varepsilon}=\Big\{\alpha\in(\underline\alpha,\overline\alpha);\;\, \mathtt{c}_\infty( \omega(\alpha),  \alpha)=-{\omega}_{\textnormal{Eq}}\quad\hbox{and}\quad  ( \omega,  \alpha) \in \mathtt{G}_{\infty}^{\kappa}\Big\}.
$$
The ultimate goal is to  measure this last Cantor and show that it is asymptotically with full  Lebesgue measure when $\varepsilon\to0$.  This important point  will be discussed in Section \ref{Section 6.2} and based on the rigidity of the equilibrium frequencies combined with perturbative arguments.
\subsection{Nash-Moser scheme}\label{N-M-S1}
The construction of the solutions  to the nonlinear equation \eqref{eq-NM} stems from the  modified  Nash-Moser scheme as in \cite{BB13,BB10}  and reproduced later in several papers as for instance in \cite{BertiMontalto,Baldi-berti}. The scheme is based on the construction of  successive approximations belonging to a finite-dimensional space given by 
$$
E_{m}\triangleq\Big\{\mathfrak{I};\;\, \Pi_m\mathfrak{I}=\mathfrak{I}\Big\}\quad\hbox{with}\quad  \mathfrak{I}:\varphi\in\T^d\mapsto \mathfrak{I}(\varphi)=\big(\Theta(\varphi),I(\varphi),z(\varphi)\big)
$$
where $\Pi_{m}$ is the projector defined   
$$
h(\varphi,\theta)=\sum_{(l,j)\in\mathbb{Z}^{d+1}}h_{l,j}e^{\ii(l\cdot\varphi+j\theta)},\quad\Pi_{m}h(\varphi,\theta)=\sum_{|l|+|j|\leqslant N_{m}}h_{l,j}e^{\ii(l\cdot\varphi+j\theta)},
$$
and when $h$ depends only on $\varphi$ it is given by
$$
h(\varphi)=\sum_{l\in\mathbb{Z}^{d}}h_{l}e^{\ii\,l\cdot\varphi},\quad\Pi_{m}h(\varphi)=\sum_{|l |\leqslant N_{m}}h_{l}e^{\ii\,l\cdot\varphi},
$$
{where the sequence of numbers  $(N_m)_{m}$ is defined  in \eqref{definition of Nm}. We point out that in the preceding sections we have used the same  notation $\Pi_{m}$ to denote the orthogonal projector localizing only on the time frequency set $\{(l,j)\in\mathbb{Z}^{d+1},\,|l|\leqslant N_m\}$. This latter projector will not be used throughout this section and therefore  there is no confusion to fear from this notation in this section. \\ 
During the Nash-Moser scheme we need a list of significant parameters $\{a_1,a_2,\overline{a},\mu_1,\mu_2, s_h\}$  that will be fixed with respect to the geometry of Cantor sets through the parameters $\tau_1,\tau_2$ and $d$. Below, we  opt for a particular choice, which is not optimal but fulfills all the required  constraints in the Nash-Moser scheme,
\begin{equation}\label{Assump-DRX}
\left\lbrace\begin{array}{rcl}
\overline{a} & = & \tau_{2}+2\\
\mu_1 & = & 6\overline{\sigma}+6+3q(\tau_2+2)\\
a_{1} & = & 6q(2+\tau_2)+12\overline{\sigma}+15\\
a_{2} & = & 3q(2+\tau_2)+6\overline{\sigma}+9\\
s_h & = & s_0+9\overline{\sigma}+\tau_2(5q+2)+9q+12\\
b_1&=&2s_h-s_0.
\end{array}\right.
\end{equation}
We will also need   the parameter $\mu_2$ introduced in  Proposition \ref{reduction of the remainder term} and given by
\begin{equation}\label{Assump-mu2}
\mu_2=\tfrac23\big(s_h-(s_0+\tau_1(q+1)+{\tau_2(q+1)}+3)\big).
\end{equation}
{Let us now give an insight  about the roles played by these numbers in the Nash-Moser scheme, see below Proposition \ref{Nash-Moser}. \\
 \ding{70} The number $\overline\sigma=\sigma(\tau_1,\tau_2,d)$ is arbitrary large and measures the total loss of regularity in the construction of the approximate inverse according to  Theorem \ref{thm:stima inverso approssimato}.\\
\qquad \ding{71} The list  $\{a_1,a_2\}$ is introduced   to basically describe the convergence rates at different  lower regularity indexes $s_0, s_0+\overline\sigma$.\\
 \qquad \ding{70} The number $\mu_1$ measures  the  norm inflation  at a higher regularity level related to $s_h$.\\
  \qquad \ding{71} The number $\mu_2$  is related to the rate of convergence of the errors in the approximate inverse, see  Theorem \ref{thm:stima inverso approssimato}.\\
  \qquad \ding{70}  The number  $\overline{a}$ is associated to the enlargement of the intermediate Cantor sets by open sets whose thickness is proportional to  $N_m^{-\overline a}.$ This procedure is needed to extend in a classical way the approximate solutions to the whole set of parameters $\mathcal{O}.$}\\}
We shall also impose the following conditions
\begin{equation}\label{choice of gamma and N0 in the Nash-Moser}
\begin{array}{cccc}
N_{0}\triangleq\kappa^{-1}, & \kappa\triangleq\varepsilon^{a} &\quad\hbox{with}\qquad   0<a<\frac{1}{\mu_2+q+3}\cdot
\end{array}
\end{equation}
The main goal is to prove by induction the following result.
\begin{proposition}[Nash-Moser]\label{Nash-Moser}
{Let $(\tau_{1},\tau_{2},q,d,s_{0},s_h)$ satisfy \eqref{Conv-T2},  \eqref{Assump-DRX} and  assume    \eqref{choice of gamma and N0 in the Nash-Moser}. There exist $C_{\ast}>0$ and ${\varepsilon}_0>0$ such that for any $\varepsilon\in[0,\varepsilon_0]$ and under  the second condition in \eqref{small-C3}  we get  for all $m\in\mathbb{N}$ the following properties,
\begin{itemize}
\item  $(\mathcal{P}1)_{m}$ There exists a $q$-times differentiable function 
$${W}_{m}:\begin{array}[t]{rcl}
\mathcal{O} & \rightarrow &  E_{m-1}\times\mathbb{R}^{d}\\
\lambda & \mapsto & ({\mathfrak{I}}_{m},\mathtt{c}_{m}{-}\omega)
\end{array}$$
satisfying 
$$
{W}_{0}=0\quad\mbox{ and }\quad\,\| {W}_{m}\|_{s_{0}+\overline{\sigma}}^{q,\kappa}\leqslant C_{\ast}\varepsilon\kappa^{-1}{N_0^{q\overline{a}}}\quad \hbox{for}\quad m\geqslant1.
$$
By setting 
$$
U_{0}=\big((\varphi,0,0),\omega\big),\quad {U}_{m}=U_{0}+{W}_{m}\quad \hbox{and} \quad {H}_{m} ={U}_{m}-{U}_{m-1} \quad \hbox{for}\quad m\geqslant1,
$$
 then we have $\forall s\in[s_0,S],$
$$
\| {H}_{1}\|_{s}^{q,\kappa}\leqslant \tfrac12C_{\ast}\varepsilon\kappa^{-1}{N_0}^{q\overline{a}}\quad\mbox{ and }\quad\| {H}_{k}\|_{{{s}_0+\overline{\sigma}}}^{q,\kappa}\leqslant C_{\ast}\varepsilon\kappa^{-1}N_{k-1}^{-a_{2}}\quad  \forall\,\, 2\leqslant k\leqslant m.
$$
\item $(\mathcal{P}2)_{m}$ Define, ${i}_{m}=(\varphi,0,0)+{\mathfrak{I}}_{m}$, $ \kappa_{m}=\kappa(1+2^{-m})\in[\kappa,2\kappa]\,$ and 
$$
{\mathcal{A}_{0}^{\kappa}=\mathcal{O}\quad \mbox{ and }\quad \mathcal{A}_{m+1}^{\kappa}=\mathcal{A}_{m}^{\kappa}\cap\mathtt{G}_{m}(\kappa_{m+1},\tau_{1},\tau_{2},{i}_{m})} \quad\forall m\in\mathbb{N},
$$
where $\mathtt{G}(\kappa_{m+1},\tau_{1},\tau_{2},i_{m})$ is defined in Theorem $\ref{inversion of the linearized operator in the normal directions}$ and $\mathtt {DC}_{N_{0}}$ is described by \eqref{DC tau gamma N}.\\ Consider the open sets 

$$ \forall\, r>0,\quad  \mathrm{O}_{m}^{r}\triangleq \Big\{\lambda\in\mathcal{O};\;\, {\mathtt{dist}}\big(\lambda,\mathcal{A}_{m}^{\kappa}\big)< r\,N_{{m }}^{-\overline{a}}\Big\}
$$
where $\displaystyle {\mathtt{dist}}(x,A)=\inf_{y\in A}\|x-y\|$.  Then we have the following estimate 
$$\|\mathcal{F}({U}_{m})\|_{s_{0}}^{q,\kappa,\mathrm{O}_{m}^{\kappa}}\leqslant C_{\ast}\varepsilon N_{m-1}^{-a_{1}}.
$$
\item $(\mathcal{P}3)_{m}$ $\| {W}_{m}\|_{b_1+\overline\sigma}^{q,\kappa,\mathcal{O}}\leqslant C_{\ast}\varepsilon\kappa^{-1}N_{m-1}^{\mu_1}.$
\end{itemize}}
\end{proposition}
\begin{remark}
For any open set $O\subset \mathcal{O},$ the norm $\|\cdot\|_{s}^{q,\kappa,O}$ is defined according to the  \mbox{Definition $\ref{Def-WS}$} by simply changing the set $\mathcal{O}$ with $O$.
\end{remark}
\begin{proof}

The proof will be done using an induction principle.

\smallskip

\ding{202} \textit{ Initialization.} By construction, $U_{0}=\Big((\varphi,0,0),\omega\Big)$ and thus using  \eqref{operatorF} we obtain
$$\mathcal{F}(U_{0})=\varepsilon\left(\begin{array}{c}
-\partial_{I}\mathcal{P}_{\varepsilon}((\varphi,0,0))\\
\partial_{\vartheta}\mathcal{P}_{\varepsilon}((\varphi,0,0))\\
-\partial_{\theta}\nabla_{z}\mathcal{P}_{\varepsilon}((\varphi,0,0)).
\end{array}\right)$$
Thus, using Lemma \ref{tame estimates for the vector field XmathcalPvarepsilon}-{(i)}, we find 
\begin{equation}\label{estimate mathcalF(U0)}
\forall s\geqslant 0,\quad \|\mathcal{F}(U_{0})\|_{s}^{q,\kappa,\mathcal{O}}\leqslant C_*\varepsilon
\end{equation}
for some positive  constant $C_{\ast}$. The properties $(\mathcal{P}1)_{0},$ $(\mathcal{P}2)_{0}$ and $(\mathcal{P}3)_{0}$ then follow immediately.

\smallskip

\ding{203} \textit{ Induction step.} Given $m\in\N^\star$ and assume that $(\mathcal{P}1)_{k},$ $(\mathcal{P}2)_{k}$ and $(\mathcal{P}3)_{k}$ are true for any $0\leqslant k\leqslant m$. The goal is to  check the validity of theses properties  at the  order $m+1$. To do so, we start first with the  construction of   the approximation $U_{m+1}$ by using  a modified Nash-Moser scheme. We shall introduce the linearized operator of $\mathcal{F}$ at the state $(i_m,\mathtt{c}_m)$
\begin{align*}
L_{m}&\triangleq L_{m}(\lambda)\\
&=d_{i,\mathtt{c}}\mathcal{F}({i}_{m}(\lambda,\omega){,\mathtt{c}_{m}(\lambda)}).
\end{align*}
As we shall see later, the constriction of the next approximation  $U_{m+1}$ requires an approximate right inverse of $L_m$. This task is discussed in the preceding sections and a precise statement is stated in Theorem \ref{thm:stima inverso approssimato}. To apply this result and get some bounds on $U_{m+1}$  we need to establish first some intermediate results connected  with the smallness condition and some Cantor set inclusions.

\smallskip

\foreach \x in {\bf a} {%
  \textcircled{\x}
}   \textbf{Smallness/boundedness conditions.} We remark that with  the conditions \eqref{Assump-DRX} and \eqref{Assump-mu2} on the parameters $s_h$ and $\mu_2$ the conditions listed in \eqref{cond-diman1} are automatically satisfied and therefore the results of Proposition \ref{reduction of the remainder term} hold true provided that the boundedness/smallness conditions \eqref{small-C3} are verified. For the smallness condition, it is satisfied  provided that
\begin{align*}
\quad N_{0}^{\mu_{2}}\varepsilon{\kappa^{-3-q}}&=\varepsilon^{1-a(\mu_2+3+q)}\\
&\leqslant {\varepsilon}_0,
\end{align*} 
with $\varepsilon_0$ being a small positive  number. This holds true for small $\varepsilon$ by virtue of \eqref{choice of gamma and N0 in the Nash-Moser}. Concerning the boundedness condition in \eqref{small-C3}, we recall that $\sigma\geqslant\sigma_4$ and use 
 Lemma \ref{interpolation-In} leading to 
\begin{align}\label{modon1}
 \|{H}_{m}\|_{s_{h}+\overline\sigma}^{q,\kappa}\lesssim\left(\|{H}_{m}\|_{s_{0}+\overline\sigma}^{q,\kappa}\right)^{{\frac12}}\left(\|{H}_{m}\|_{b_{1}+\overline\sigma}^{q,\kappa}\right)^{\frac12},\quad  b_1= 2s_h-s_0.
 \end{align}
On the other hand, by using $(\mathcal{P}1)_{m}$ we find
\begin{align}\label{rent-P76}
\forall s\in[s_0,S],\,\|{H}_{1}\|_{s}^{q,\kappa}\leqslant\tfrac12 C_{\ast}\varepsilon\kappa^{-1}N_0^{q\overline{a}}\quad\hbox{and}\quad \quad \|{H}_{m}\|_{s_{0}+\overline\sigma}^{q,\kappa}\leqslant C_{\ast}\varepsilon\kappa^{-1}N_{m-1}^{-a_{2}}.
\end{align}
Now applying  $(\mathcal{P}3)_{m}$ and $(\mathcal{P}3)_{m-1}$ yields
\begin{align*}
\|{H}_{m}\|_{b_{1}+\overline\sigma}^{q,\kappa}&\leqslant\|{W}_{m}\|_{b_{1}+\overline\sigma}^{q,\kappa}+\|{W}_{m-1}\|_{b_{1}+\overline\sigma}^{q,\kappa}\\
& \leqslant 2C_{\ast}\varepsilon\kappa^{-1}N_{m-1}^{\mu}.
\end{align*}
Inserting the foregoing estimates into \eqref{modon1} allows to get
$$\|{H}_{m}\|_{s_{h}+\overline\sigma}^{q,\kappa}\leqslant CC_*\varepsilon\kappa^{-1}N_{m-1}^{\frac12\mu-\frac12a_{2}}.$$
However for $m=1$ one has
$$\|{H}_{1}\|_{s_{h}+\overline\sigma}^{q,\kappa}\leqslant \tfrac12 C_{\ast}\varepsilon\kappa^{-1}N_0^{q\overline{a}}.
$$
From \eqref{Assump-DRX}, one gets  
\begin{equation}\label{cond_b1_1}
\mu-a_2\leqslant -2.
\end{equation}
Thus, using \eqref{choice of gamma and N0 in the Nash-Moser} and taking  $\varepsilon$ small enough, we obtain
\begin{align*}
\|{W}_{m}\|_{s_{h}+\overline\sigma}^{q,\kappa}&\leqslant\|{H}_{1}\|_{s_{h}+\overline\sigma}^{q,\kappa}+\sum_{k=2}^{m}\|{H}_{k}\|_{s_{h}+\overline\sigma}^{q,\kappa}\\
&\leqslant \tfrac12C_*\varepsilon\kappa^{-1}N_0^{q\overline{a}}+CC_*\varepsilon\kappa^{-1}\sum_{k=0}^{\infty}N_{k}^{-1}\\
&\leqslant C_* \varepsilon^{1-a(1+q\overline{a})}.
\end{align*}
Hence  by virtue of \eqref{choice of gamma and N0 in the Nash-Moser} and \eqref{Assump-DRX} we find
\begin{equation}\label{cond-Lu1}
a\leqslant \frac{1}{2(1+q\overline{a})}
\end{equation}
which implies in turn  that 
\begin{align}\label{imp-yi1}
\|{W}_{m}\|_{s_{h}+\sigma}^{q,\kappa}&\leqslant C\varepsilon^{\frac12}\\
\nonumber&\leqslant 1,
\end{align}
provided that $\varepsilon$ is small enough.  Another important observation is that $\overline\sigma$ can be chosen large enough in order to get $s_0+\overline\sigma\geqslant \overline{s}_h+\sigma_4$ where $\overline{s}_h$ is defined in \eqref{Conv-T2N}. By this way we get from the second estimate of  \eqref{rent-P76} and Sobolev embeddings
\begin{align}\label{rent-P78}
  \|{H}_{m}\|_{\overline{s}_{h}+\sigma_4}^{q,\kappa}\leqslant C_* \varepsilon\kappa^{-1}N_{m-1}^{-a_{2}}.
\end{align}

\smallskip

\foreach \x in {\bf b} {%
  \textcircled{\x}
}   \textbf{Set inclusions.} 
Notice that by the previous point and the conditions \eqref{Assump-DRX} and \eqref{choice of gamma and N0 in the Nash-Moser} we can perform the reducibility  of the linearized operator in the normal directions at the step $m$. In particular, Propositions \ref{QP-change} and \ref{reduction of the remainder term} apply. Thus the sets $\mathcal{A}_{k}^{\kappa}$  are well-defined for all $k\leqslant m+1$ and to develop  later suitable estimates  we need to establish  first  the following inclusions
\begin{align}\label{Inclu-z0}
\mathcal{A}_{m+1}^{\kappa}\subset\mathrm{O}_{m+1}^{2\kappa}\subset\left(\mathcal{A}_{m+1}^{\frac{\kappa}{2}}\cap\mathrm{O}_{m}^{\kappa}\right).
\end{align}
{We draw the reader's attention not to confuse $\kappa$ and $k$ in this part.} The first left inclusion in \eqref{Inclu-z0}   is obvious by construction of $\mathrm{O}_{m+1}^{2\kappa}$ which is an enlargement of $\mathcal{A}_{m+1}^{\kappa}.$ As to the second inclusion, we first claim that 
			\begin{equation}\label{inco}
				\forall k\in \llbracket 0,m\rrbracket,\quad  \mathrm{O}_{k+1}^{2\kappa}\subset \mathrm{O}_{k}^{\kappa}.
			\end{equation}
			Indeed, since by construction $\mathcal{A}_{k+1}^{\kappa}\subset\mathcal{A}_{k}^{\kappa}$ then taking $\lambda\in \mathrm{O}_{k+1}^{2\kappa}$ we have the following estimates
			\begin{align*}
				{\mathtt{dist}}\big(\lambda,\mathcal{A}_{k}^{\kappa}\big)&\leqslant  {\mathtt{dist}}\big(\lambda,\mathcal{A}_{k+1}^{\kappa}\big)\\
				&< 2\kappa N_{k+1}^{-\overline{a}}= 2\kappa N_{k}^{-\overline{a}}N_{0}^{-\frac{1}{2}\overline{a}}\\
				& < \kappa N_{k}^{-\overline{a}}
			\end{align*}
			provided that $2N_0^{-\frac12\overline a}<1,$ which is satisfied for $N_0$ large enough, that is in view of \eqref{choice of gamma and N0 in the Nash-Moser} for $\varepsilon$ small enough. The next task is to prove by induction in $k$ that
			\begin{equation}\label{hyprec O in A}
				{\forall k\in \llbracket 0,m+1\rrbracket,\quad  \mathrm{O}_{k}^{2\kappa}\subset\mathcal{A}_{k}^{\frac{\kappa}{2}}.}
			\end{equation}
			Remark that for $k=0$, one has $\mathrm{O}_{0}^{2\kappa}=\mathcal{O}=\mathcal{A}_{0}^{\frac{\kappa}{2}}.$ Now assume that \eqref{hyprec O in A} occurs at the order $k\in \llbracket 0,n\rrbracket$ and let us check  it  at the order $k+1.$ Using \eqref{inco} and the induction assumption \eqref{hyprec O in A}, we get
			$$\mathrm{O}_{k+1}^{2\kappa}\subset\mathrm{O}_{k}^{\kappa}\subset\mathrm{O}_{k}^{2\kappa}\subset\mathcal{A}_{k}^{\frac{\kappa}{2}}.$$
			Hence, it remains to check that 
			$$\mathrm{O}_{k+1}^{2\kappa}\subset\mathtt{G}_{k}\Big(\tfrac{\kappa_{k+1}}{2},\tau_{1},\tau_{2},i_{k}\Big).$$
Given  $\lambda\in\mathrm{O}_{k+1}^{2\kappa}$
there exists $\lambda'=(\omega',\alpha')\in\mathcal{A}_{k+1}^{\kappa}$ such that $\mathtt{dist}\left(\lambda,\lambda'\right)<2\kappa N_{k+1}^{-\overline{a}}.$  Then for all $(l,j)\in\mathbb{Z}^{d}\times\mathbb{S}_{0}^{c}$ with $|l|\leqslant N_{k},$ we get using  the triangle inequality and   Cauchy-Schwarz inequalities
\begin{align*}
\big|\omega\cdot l+\mu_{j}^{\infty}(\lambda,{i}_{k})\big|&\geqslant\big|\omega'\cdot l+\mu_{j}^{\infty}(\lambda',{i}_{k})\big|-|\omega-\omega'||l|-\big|\mu_{j}^{\infty}(\lambda,{i}_{k})-\mu_{j}^{\infty}(\lambda',{i}_{k})\big|\\
&>\tfrac{\kappa_{k+1}\langle j\rangle}{\langle l\rangle^{\tau_{1}}}-2\kappa N_{k+1}^{1-\overline{a}}-\big|\mu_{j}^{\infty}(\lambda,{i}_{k})-\mu_{j}^{\infty}(\lambda',{i}_{k})\big|.
\end{align*}
Using the Mean Value Theorem yields
 \begin{align}\label{Tb-X1}
\left|\mu_{j}^{\infty}(\lambda,{i}_{k})-\mu_{j}^{\infty}(\lambda',{i}_{k})\right|&\leqslant |\lambda-\lambda^\prime|\kappa^{-1}\|\mu_{j}^{\infty}({i}_{k})\|^{q,\kappa}\\
 \nonumber &\leqslant 2 N_{k+1}^{-\overline{a}}\|\mu_{j}^{\infty}({i}_{k})\|^{q,\kappa}.
\end{align}
Therefore combining  \eqref{Tb-X1} with the assumption $(\mathcal{P}1)_{m}$ and the eigenvalues estimates given in  Propositions \ref{reduction of the remainder term}-\ref{projection in the normal directions} we get in view of \eqref{small-C3}
\begin{align}\label{Tb-X2}
\big|\mu_{j}^{\infty}(\lambda,{i}_{k})-\mu_{j}^{\infty}(\lambda',{i}_{k})\big| &\leqslant C\kappa\,\langle j\rangle N_{k+1}^{-\overline{a}}.
\end{align}
Hence, by \eqref{small-C3} we get  since $|l|\leqslant N_{k}$ and $\kappa_{k+1}\geqslant \kappa$, 
\begin{align*}
\big|\omega\cdot l+\mu_{j}^{\infty}(\lambda,{i}_{k})\big|&>\tfrac{\kappa_{k+1}\langle j\rangle}{\langle l\rangle^{\tau_{1}}}-C\kappa\langle j\rangle N_{k+1}^{1-\overline{a}}\\
&>\tfrac{\kappa_{k+1}\langle j\rangle}{\langle l\rangle^{\tau_{1}}}\left(1-CN_{k+1}^{\tau_{1}+1-\overline{a}}\right).
\end{align*}
From \eqref{Assump-DRX} and \eqref{Conv-T2} we find 
$
\tau_{1}+2\leqslant\overline{a}
$ and taking  $N_{0}$ sufficiently large we get
$$CN_{k+1}^{\tau_{1}+1-\overline{a}}\leqslant CN_{0}^{-1}\leqslant\tfrac12$$
so that
$$
\big|\omega\cdot l+\mu_{j}^{\infty}(\lambda,{i}_{k})\big|>\tfrac{\kappa_{k+1}\langle j\rangle}{2\langle l\rangle^{\tau_{1}}}\cdot
$$
As a consequence, we deduce that  $\lambda\in\Lambda_{\infty,k}^{\frac{\kappa_{k+1}}{2},\tau_{1}}({i}_{k}).$ Let us now check that $\lambda\in \mathcal{O}_{\infty,k}^{\frac{\kappa_{k+1}}{2},\tau_{1}}({i}_{k})$ defined in Proposition \ref{QP-change} (using the new notation $c_{i}(\lambda)=c(\lambda,i)$),
\begin{align*}
\big|\omega\cdot l+jc_{{i}_{k}}(\lambda)\big|&\geqslant\big|\omega'\cdot l+jc_{{i}_{k}}(\lambda')\big|-|\omega-\omega'||l|-|j|\big|c_{{i}_{k}}(\lambda)-c_{{i}_{k}}(\lambda')\big|\\
&>4\,\kappa_{k+1}^{\varrho}\tfrac{\langle j\rangle}{\langle l\rangle^{\tau_{1}}}-2\kappa N_{k+1}^{1-\overline{a}}-\langle j\rangle\big|c_{{i}_{k}}(\lambda)-c_{{i}_{k}}(\lambda')\big|.
\end{align*}
Using the Mean Value Theorem and \eqref{est-r1} combined with \eqref{small-C3} and the  assumption $(\mathcal{P}1)_{m}$
\begin{align*}
\big|c_{{i}_{k}}(\lambda)-c_{{i}_{k}}(\lambda')\big|&\leqslant |\lambda-\lambda^\prime|\big(C+\kappa^{-1
}\|c_{r,\lambda}-V_{0,\alpha}\|^{q,\kappa}\big)+ 2\kappa\,N_{k+1}^{-\overline{a}}\big(C+C\varepsilon\kappa^{-2
}\big)\\
&\leqslant C\kappa\,N_{k+1}^{-\overline{a}}.
\end{align*}
Hence as before  we get from the definition of $\kappa_{m}$ and $\varrho\in(0,1)$
\begin{align*}
\big|\omega\cdot l+c_{{i}_{k}}(\lambda)\big|&> 4\,\kappa_{k+1}^{\varrho}\tfrac{\langle j\rangle}{\langle l\rangle^{\tau_{1}}}-C\kappa\langle j\rangle N_{k+1}^{1-\overline{a}}\\
&\geqslant 4\tfrac{(\tfrac{\kappa_{k+1}}{2})^{\varrho}\langle j\rangle}{\langle l\rangle^{\tau_{1}}}\left(2^{\varrho}-C2^{\varrho}N_{k+1}^{\tau_{1}+1-\overline{a}}\right).
\end{align*}
Since by assumption  $\overline{a}\geqslant \tau_{1}+2$ then taking  $N_{0}$ sufficiently large we get
$$
CN_{k+1}^{\tau_{1}+1-\overline{a}}\leqslant CN_{0}^{-1}<1-2^{-\varrho}
$$
implying  that
$$
\left|\omega\cdot l+jc_{{i}_{k}}(\lambda)\right|>4\left(\tfrac{\kappa_{k+1}}{2}\right)^{\varrho}\tfrac{\langle j\rangle}{\langle l\rangle^{\tau_{1}}}\cdot
$$
This shows that  $\lambda\in\mathcal{O}_{\infty,k}^{\frac{\kappa_{k+1}}{2},\tau_{1}}({i}_{k}).$ It remains to check that  $\lambda\in \mathcal{O}_{\infty,k}^{\frac{\kappa_{k+1}}{2},\tau_{1},\tau_2}({i}_{k})$ defined in \mbox{Proposition \ref{reduction of the remainder term}.} We write by the triangle inequality and the condition $\lambda'\in\mathcal{A}_{m+1}^{\kappa}$
\begin{align*}
\left|\omega\cdot l+\mu_{j}^{\infty}(\lambda,{i}_{k})-\mu_{j_{0}}^{\infty}(\lambda,{i}_{k})\right|&\geqslant\left|\omega'\cdot l+\mu_{j}^{\infty}(\lambda',{i}_{k})-\mu_{j_{0}}^{\infty}(\lambda',{i}_{k})\right|-|\omega-\omega'||l|\\
&\quad-\left|\mu_{j}^{\infty}(\lambda,{i}_{k})-\mu_{j_{0}}^{\infty}(\lambda,{i}_{k})+\mu_{j_{0}}^{\infty}(\lambda',{i}_{k})-\mu_{j}^{\infty}(\lambda',{i}_{k})\right|\\
&>\tfrac{2\kappa_{k+1}\langle j-j_{0}\rangle}{\langle l\rangle^{\tau_{2}}}-2\kappa N_{k+1}^{1-\overline{a}}\\
&\quad -\left|\mu_{j}^{\infty}(\lambda,{i}_{k})-\mu_{j_{0}}^{\infty}(\lambda,{i}_{k})+\mu_{j_{0}}^{\infty}(\lambda',{i}_{k})-\mu_{j}^{\infty}(\lambda',{i}_{k})\right|.
\end{align*}
Now recall that from Proposition \ref{reduction of the remainder term}, we can write
$$\mu_{j}^{\infty}(\lambda,{i}_{k})=\mu_{j}^{0}(\lambda,{i}_{k})+r_{j}^{\infty}(\lambda,{i}_{k})
$$
and then
\begin{align*}
&\left|\mu_{j}^{\infty}(\lambda,{i}_{k})-\mu_{j_{0}}^{\infty}(\lambda,{i}_{k})+\mu_{j_{0}}^{\infty}(\lambda',{i}_{k})-\mu_{j}^{\infty}(\lambda',{i}_{k})\right|\\
&\qquad\leqslant\left|\mu_{j}^{0}(\lambda,{i}_{k})-\mu_{j_{0}}^{0}(\lambda,{i}_{k})+\mu_{j_{0}}^{0}(\lambda',{i}_{k})-\mu_{j}^{0}(\lambda',{i}_{k})\right|\\
&\qquad+\left|r_{j}^{\infty}(\lambda,{i}_{k})-r_{j}^{\infty}(\lambda',{i}_{k})\right|+\left|r_{j_{0}}^{\infty}(\lambda,{i}_{k})-r_{j_{0}}^{\infty}(\lambda',{i}_{k})\right|.
\end{align*}
Using the Mean Value Theorem, combined with Proposition \ref{projection in the normal directions},  Lemma \ref{lem-asym}-(vi) and \eqref{small-C3} yield
$$\left|\mu_{j}^{0}(\lambda,{i}_{k})-\mu_{j_{0}}^{0}(\lambda,{i}_{k})+\mu_{j_{0}}^{0}(\lambda',{i}_{k})-\mu_{j}^{0}(\lambda',{i}_{k})\right|\leqslant C\kappa \,N_{k+1}^{-\overline{a}}\langle j-j_{0}\rangle.$$
Similarly we get in view of  \eqref{estimate rjinfty} 
\begin{align*}\left|r_{j}^{\infty}(\lambda,{i}_{k})-r_{j}^{\infty}(\lambda',{i}_{k})\right|&\leqslant C\kappa N_{k+1}^{-\overline{a}}\varepsilon \kappa^{-3}\\
&\leqslant C\kappa N_{k+1}^{-\overline{a}}\langle j-j_{0}\rangle.
\end{align*}
Putting together the foregoing estimates and  the facts that $|l|\leqslant N_{k}$ and $\kappa_{k+1}\geqslant\kappa$ yield
\begin{align*}
\left|\omega\cdot l+\mu_{j}^{\infty}(\lambda,{i}_{k})-\mu_{j_{0}}^{\infty}(\lambda,{i}_{k})\right|&\geqslant\kappa_{k+1}\tfrac{\langle j-j_{0}\rangle}{\langle l\rangle^{\tau_{2}}}\left(2-CN_{k+1}^{\tau_{2}+1-\overline{a}}\right).
\end{align*}
Since  $\overline{a}= \tau_{2}+2$, see \eqref{Assump-DRX}, 
then  taking $N_{0}$ sufficiently large we find
$$CN_{k+1}^{\tau_{2}+1-\overline{a}}\leqslant CN_{0}^{-1}<1$$
we deduce that
$$\left|\omega\cdot l+\mu_{j}^{\infty}(\lambda,{i}_{k})-\mu_{j_{0}}^{\infty}(\lambda,{i}_{k})\right|>\kappa_{k+1}\tfrac{\langle j-j_{0}\rangle}{\langle l\rangle^{\tau_{2}}}\cdot
$$
It follows that  $\lambda\in\mathcal{O}_{\infty,k}^{\frac{\kappa_{k+1}}{2},\tau_{1},\tau_{2}}({i}_{k}).$ Consequently, we deduce that   $\lambda\in\mathtt{G}_{k}(\tfrac12\kappa_{k+1},\tau_{1},\tau_{2},{i}_{k})$ and thus we conclude that $\lambda\in\mathcal{A}_{k+1}^{\frac{\kappa}{2}}.$ This ends the proof of \eqref{hyprec O in A} and therefore the proof of \eqref{Inclu-z0} is now complete. 

\smallskip

\foreach \x in {\bf c} {%
  \textcircled{\x} {\bf Construction of the next approximation.}
}  
Since the assumption \eqref{small-C3} us satisfied at the order $m$ we can  perform the  reduction of the linearized operator stated in  Theorem  \ref{inversion of the linearized operator in the normal directions}. Then  coming back to  Theorem  \ref{thm:stima inverso approssimato} applied with  $L_{m}$ we find  an operator $\mathbf{T}_{m}\triangleq\mathbf{T}_{m}(\lambda)$ well-defined in the whole set of parameters $\lambda\in\mathcal{O}$ with the estimates
\begin{equation}\label{estimate Tm}
\forall s\in[s_{0},S ],\quad\|\mathbf{T}_{m}h\|_{s}^{q,\kappa}\lesssim\kappa^{-1}\left(\|h\|_{s+\overline{\sigma}}^{q,\kappa}+\|{\mathfrak{I}}_{m}\|_{s+\overline{\sigma}}^{q,\kappa}\|h\|_{s_{0}+\overline{\sigma}}^{q,\kappa}\right)
\end{equation}
and
\begin{equation}\label{estimate Tm in norm s0}
\|\mathbf{T}_{m}h\|_{s_{0}}^{q,\kappa}\lesssim\kappa^{-1}\|h\|_{s_{0}+\overline{\sigma}}^{q,\kappa}.
\end{equation}
In addition, when the parameters are restricted in  the Cantor set $\mathtt{G}_{n}(\kappa_{m+1},\tau_{1},\tau_{2},{i}_{m}),$ $\mathbf{T}_{m}$ is  an approximate right inverse \mbox{of $L_{m}$} in the sense of \eqref{splitting per approximate inverse} which will be useful later.

Next, we  define the  function $\widetilde{U}_{m+1}$ as follows
$$
\widetilde{U}_{m+1}\triangleq {U}_{m}+\widetilde{H}_{m+1}\quad\mbox{ with }\quad \widetilde{H}_{m+1}\triangleq(\widehat{\mathfrak{I}}_{m+1},\widehat{\mathtt{c}}_{m+1})\triangleq-\mathbf{\Pi}_{m}\mathbf{T}_{m}\Pi_{m}\mathcal{F}({U}_{m})\in E_{m}\times\mathbb{R}^{d}$$
where $\mathbf{\Pi}_{m}$ is defined by
\begin{align}\label{and-bp1}
\mathbf{\Pi}_{m}(\mathfrak{I},\alpha)=(\Pi_{m}\mathfrak{I},\alpha)\quad \mbox{ and }\quad\mathbf{\Pi}_{m}^{\perp}(\mathfrak{I},\alpha)=(\Pi_{m}^\perp\mathfrak{I},0).
\end{align}
Observe that ${U}_{m}$ is defined in the full set $\mathcal{O}$ which implies that $\widetilde{U}_{m+1}$ is defined in  $\mathcal{O}$ too. However we shall not work with this natural extension  but we will localize it around the Cantor set $\mathcal{A}_{m+1}^{\kappa}$ in order to get a good decay.
Let us introduce the quadratic function 
\begin{align}\label{Def-Qm}
Q_{m}=\mathcal{F}(\widetilde{U}_{m+1})-\mathcal{F}({U}_{m})-L_{m}\widetilde{H}_{m+1}.
\end{align}
Then it is easy to check that
\begin{align}\label{Decom-RTT1}
 \mathcal{F}(\widetilde{U}_{m+1})
& =  \mathcal{F}({U}_{m})-L_{m}\mathbf{\Pi}_{m}\mathbf{T}_{m}\Pi_{m}\mathcal{F}({U}_{m})+Q_{m}\\
\nonumber& =  \mathcal{F}({U}_{m})-L_{m}\mathbf{T}_{m}\Pi_{m}\mathcal{F}({U}_{m})+L_{m}\mathbf{\Pi}_{m}^{\perp}\mathbf{T}_{m}\Pi_{m}\mathcal{F}({U}_{m})+Q_{m}\\
\nonumber& =  \mathcal{F}({U}_{m})-\Pi_{m}L_{m}\mathbf{T}_{m}\Pi_{m}\mathcal{F}({U}_{m})+(L_{m}\mathbf{\Pi}_{m}^{\perp}-\Pi_{m}^{\perp}L_{m})\mathbf{T}_{m}\Pi_{m}\mathcal{F}({U}_{m})+Q_{m}\\
\nonumber& =  \Pi_{m}^{\perp}\mathcal{F}({U}_{m})-\Pi_{m}(L_{m}\mathbf{T}_{m}-\textnormal{Id})\Pi_{m}\mathcal{F}({U}_{m})+(L_{m}\mathbf{\Pi}_{m}^{\perp}-\Pi_{m}^{\perp}L_{m})\mathbf{T}_{m}\Pi_{m}\mathcal{F}({U}_{m})+Q_{m}.
\end{align}

\smallskip

$\blacktriangleright$ \textbf{Estimates of $\mathcal{F}(\widetilde{U}_{m+1}), \, m\geqslant 1$.}
We shall  estimate $\mathcal{F}(U_{m+1})$ in norm $\|\cdot\|_{s}^{q,\kappa,\mathrm{O}_{m+1}^{2\kappa}}$ through estimating the right-hand side terms in \eqref{Decom-RTT1}. We note that the parameters $\lambda=(\omega,\alpha)$ are located in the small open set $\mathrm{O}_{m+1}^{2\kappa}$ containing the Cantor set $\mathcal{A}_{m+1}^{\kappa}$ as indicated in \eqref{Inclu-z0}.
\\
\quad $\quad \bullet$ \textit{Estimate of $\Pi_{m}^{\perp}\mathcal{F}({U}_{m}).$} 
By Taylor formula applied with \eqref{operatorF} and Lemma \ref{tame estimates for the vector field XmathcalPvarepsilon} combined with  \eqref{estimate mathcalF(U0)} and $(\mathcal{P}1)_{m}$, we get 
\begin{align}\label{link mathcalF(Um) and Wm}
\nonumber \forall s\geqslant s_{0},\quad\|\mathcal{F}({U}_{m})\|_{s}^{q,\kappa,\mathrm{O}_{m}^{\kappa}}&\leqslant\|\mathcal{F}(U_{0})\|_{s}^{q,\kappa,\mathrm{O}_{0}^{\kappa}}+\|\mathcal{F}({U}_{m})-\mathcal{F}(U_{0})\|_{s}^{q,\kappa,\mathrm{O}_{m}^{\kappa}}\\
&\lesssim\varepsilon+\|{W}_{m}\|_{s+\overline{\sigma}}^{q,\kappa,{\mathcal{O}}}.
\end{align}
Therefore we deduce from the smallness condition \eqref{choice of gamma and N0 in the Nash-Moser} and $(\mathcal{P}1)_{m}$,
\begin{equation}\label{estimate mathcalF(U0) in norm s0}
\kappa^{-1}\|\mathcal{F}({U}_{m})\|_{s_{0}}^{q,\kappa,\mathrm{O}_{m}^{\kappa}}\leqslant 1.
\end{equation}
By Lemma \ref{orthog-Lem1} combined with  \eqref{link mathcalF(Um) and Wm}, we infer since $b_1=2s_h-s_0$
\begin{align*}
\|\Pi_{m}^{\perp}\mathcal{F}({U}_{m})\|_{s_{0}}^{q,\kappa,\mathrm{O}_{m}^{\kappa}}&\leqslant N_{m}^{s_0-b_{1}}\|\mathcal{F}({U}_{m})\|_{b_{1}}^{q,\kappa,\mathrm{O}_{m}^{\kappa}}\\
&\lesssim N_{m}^{s_0-b_{1}}\left(\varepsilon+\|{W}_{m}\|_{b_{1}+\overline{\sigma}}^{q,\kappa,\mathcal{O}}\right).
\end{align*}
Now, by using $(\mathcal{P}3)_{m}$ combined with \eqref{definition of Nm} and \eqref{choice of gamma and N0 in the Nash-Moser} 
we obtain 
\begin{align}\label{Wm in high norm}
\nonumber\varepsilon+\|{W}_{m}\|_{b_{1}+\overline{\sigma}}^{q,\kappa,\mathcal{O}}&\leqslant\varepsilon\left(1+C_{\ast}\kappa^{-1}N_{m-1}^{\mu_1}\right)\\
&\leqslant 2C_{\ast}\varepsilon N_{m}^{1+\frac{2}{3}\mu_1}.
\end{align}
Inserting this in the previous estimate and making appeal to \eqref{Inclu-z0} we finally get 

\begin{align}\label{final estimate PiperpF(Um)}
\nonumber \|\Pi_{m}^{\perp}\mathcal{F}({U}_{m})\|_{s_{0}}^{q,\kappa,\mathrm{O}_{m}^{\kappa}} & \leqslant  \|\Pi_{m}^{\perp}\mathcal{F}({U}_{m})\|_{s_{0}}^{q,\kappa,\mathrm{O}_{m}^{\kappa}}\\
&\lesssim \varepsilon N_{m}^{s_0+1+\frac23\mu_1-b_{1}}.
\end{align}

Notice that one also gets from \eqref{link mathcalF(Um) and Wm} and \eqref{Wm in high norm}
\begin{equation}\label{HDP10}
{\|\mathcal{F}({U}_{m})\|_{b_1+\overline\sigma}^{q,\kappa,\mathrm{O}_{m}^{\kappa}}}  \leqslant  C_*\varepsilon N_{m}^{\overline{\sigma}+\frac{2}{3}\mu_1+1}.
\end{equation}

\smallskip

$\bullet$ \textit{Estimate of $\Pi_{m}(L_{m}\mathbf{T}_{m}-\textnormal{Id})\Pi_{m}\mathcal{F}({U}_{m})$.} According to  \eqref{Inclu-z0} we may write
$$
\mathrm{O}_{m+1}^{2\kappa}\subset \mathcal{A}_{m+1}^{\frac{\kappa}{2}}\subset\mathtt{G}_{m}(\tfrac12\kappa_{m+1},\tau_{1},\tau_{2},{i}_{m})\cap \mathrm{O}_{m}^{\kappa}.
$$
Then we can apply   Theorem \ref{thm:stima inverso approssimato} 
 \begin{align}\label{Rent-M1}
\nonumber\|\Pi_{m}(L_{m}\mathbf{T}_{m}-\textnormal{Id})\Pi_{m}\mathcal{F}({U}_{m})\|_{s}^{q,\kappa,\mathrm{O}_{m+1}^{2\kappa}}&\leqslant \|\mathcal{E}^{(m)}_1 \Pi_{m}\mathcal{F}({U}_{m}) \|_{s}^{q,\kappa,\mathrm{O}_{m}^{\kappa}}+\|\mathcal{E}^{(m)}_2 \Pi_{m}\mathcal{F}({U}_{m})\|_{s}^{q,\kappa,\mathrm{O}_{m}^{\kappa}}\\
\nonumber &\quad+\|\mathcal{E}^{(m)}_3 \Pi_{m}\mathcal{F}({U}_{m}) \|_{s}^{q,\kappa,\mathrm{O}_{m}^{\kappa}}\\
& \triangleq \mathtt{I}^{(m)}_1(s)+ \mathtt{I}^{(m)}_2(s)+ \mathtt{I}^{(m)}_3(s).
\end{align}
Let us start with the estimate of the three terms on the right hand side for $m\geqslant 1$.
Applying the estimates of Theorem \ref{thm:stima inverso approssimato} and using    $(\mathcal{P}2)_{m}$ 
 \begin{align}\label{IJKL1}
\nonumber \mathtt{I}^{(m)}_1(s_0)  &\lesssim \kappa^{-1 } \|\mathcal{F}({U}_{m}) \|_{s_0 +\overline\sigma}^{q,\kappa,\mathrm{O}_{m}^{\kappa}} \| \Pi_{m}\mathcal{F}({U}_{m}) \|_{s_0 + \overline\sigma}^{q,\kappa,\mathrm{O}_{m}^{\kappa}}\\
 \nonumber& \lesssim\kappa^{-1 } N_m^{\overline\sigma}\|\mathcal{F}({U}_{m}) \|_{s_0 +\overline\sigma}^{q,\kappa,\mathrm{O}_{m}^{\kappa}} \| \mathcal{F}({U}_{m}) \|_{s_0 }^{q,\kappa,\mathrm{O}_{m}^{\kappa}}\\
 &\lesssim \varepsilon\kappa^{-1} N_{m}^{\overline\sigma-\frac23a_{1}} \| \mathcal{F}({U}_{m}) \|_{s_0+\overline\sigma }^{q,\kappa,\mathrm{O}_{m}^{\kappa}}.
  \end{align}
   On the other hand, one may write in view of Lemma \ref{orthog-Lem1},  $(\mathcal{P}2)_{m}$, $(\mathcal{P}3)_{m}$ and \eqref{HDP10}
 \begin{align*}
  \| \mathcal{F}({U}_{m}) \|_{s_0+\sigma }^{q,\kappa,\mathrm{O}_{m}^{\kappa}}&\lesssim N_m^{\overline\sigma} \| \mathcal{F}({U}_{m}) \|_{s_0 }^{q,\kappa,\mathrm{O}_{m}^{\kappa}}+ \|\Pi_m^{\perp} \mathcal{F}({U}_{m}) \|_{s_0+\overline\sigma }^{q,\kappa,\mathrm{O}_{m}^{\kappa}}\\
  &\lesssim N_m^{\overline\sigma} \| \mathcal{F}({U}_{m}) \|_{s_0 }^{q,\kappa,\mathrm{O}_{m}^{\kappa}}+N_m^{s_0-b_1} \|\mathcal{F}({U}_{m}) \|_{b_1+\overline\sigma }^{q,\kappa,\mathrm{O}_{m}^{\kappa}}\\
  &\lesssim \varepsilon N_{m}^{\overline\sigma-\frac23a_1} +\varepsilon   N_{m}^{s_0+\overline{\sigma}+\frac{2}{3}\mu_1+1-b_1}. \end{align*}
 Putting together the preceding two estimates and using $\varepsilon\kappa^{-1}\leqslant 1$,  which follows from  \eqref{choice of gamma and N0 in the Nash-Moser}, give
 \begin{align}\label{Es-I10}
 \mathtt{I}^{(m)}_1(s_0) & \lesssim  \varepsilon\Big(N_{m}^{2\sigma-\frac43a_{1}}+ N_{m}^{s_0+2\overline{\sigma}+\frac{2}{3}\mu_1+1-\frac23a_1-b_1}\Big).
  \end{align}
Coming back once again  to Theorem \ref{thm:stima inverso approssimato} 
$$
   \mathtt{I}^{(m)}_2(s_0)  \lesssim 
\kappa^{- 1} N_m^{s_0- b_1 } \Big( \| \mathcal{F}({U}_{m}) \|_{b_1+ \overline\sigma }^{q,\kappa,\mathrm{O}_{m}^{\kappa}}+\eps
N_m^{\overline\sigma}\| W_{m} \|_{b_1
+ \overline\sigma    }^{q,\kappa} \big \|  \mathcal{F}({U}_{m}) \|_{s_0 }^{q,\kappa,\mathrm{O}_{m}^{\kappa}} \Big).
$$
It follows from $(\mathcal{P}2)_{m}$, $(\mathcal{P}3)_{m}$ and \eqref{HDP10}
  \begin{align*}
   \mathtt{I}^{(m)}_2(s_0) & \lesssim 
\varepsilon\kappa^{- 1} N_m^{s_0- b_1 } \Big( N_{m}^{\overline{\sigma}+\frac{2}{3}\mu+1}+
N_m^{\overline\sigma}\kappa^{-1}N_{m-1}^{\mu_1}\varepsilon N_{m-1}^{-a_{1}} \Big)\\
&\lesssim 
\varepsilon N_m^{\overline{\sigma}+\frac{2}{3}\mu_1+2+s_0- b_1 }.
  \end{align*}
 For the third term, we shall  use once again  Theorem \ref{thm:stima inverso approssimato}
 allowing to get in view of $(\mathcal{P}1)_{m}$,  $(\mathcal{P}2)_{m}$, $(\mathcal{P}3)_{m}$, \eqref{HDP10} and \eqref{small-C3}
\begin{align*}
\mathtt{I}^{(m)}_3(s_0) & \lesssim N_m^{s_0-b_1}\kappa^{-2}\Big( \|  \Pi_{m}\mathcal{F}({U}_{m})\|_{b_1+\overline\sigma}^{q,\kappa,\mathrm{O}_{m}^{\kappa}}+{\varepsilon\kappa^{-3}}\| W_{m}\|_{b_1+\overline\sigma}^{q,\kappa}\|  \Pi_{m}\mathcal{F}({U}_{m})\|_{s_0+{\overline\sigma}}^{q,\kappa,\mathrm{O}_{m}^{\kappa}} \Big)\\
&\quad+ \varepsilon\kappa^{-5}N_{0}^{{\mu}_{2}}{N_{m}^{-\mu_{2}}} \|  \Pi_{m}\mathcal{F}({U}_{m})\|_{s_0+\overline\sigma}^{q,\kappa}\\
& \lesssim \varepsilon\Big( N_m^{s_0-b_1+\overline{\sigma}+\frac{2}{3}\mu_1+3}
+N_m^{s_0-b_1-a_1+\overline{\sigma}+\frac{2}{3}\mu_1}+ {N_{m}^{\overline\sigma+2-\mu_{2}-\frac23a_1}}\Big).
\end{align*}
Therefore
\begin{align*}
\mathtt{I}^{(m)}_3(s_0) & \lesssim \varepsilon\Big( N_m^{s_0-b_1+\overline{\sigma}+\frac{2}{3}\mu_1+3}
+ {N_{m}^{\overline\sigma+2-\mu_{2}-\frac23a_1}}\Big).
\end{align*}
 Plugging the preceding estimates into \eqref{Rent-M1} yields for $m\geqslant 1$
 \begin{equation}\label{Rent-M01}
\|\Pi_{m}(L_{m}\mathbf{T}_{m}-\textnormal{Id})\Pi_{m}\mathcal{F}({U}_{m})\|_{s_0}^{q,\kappa,\mathrm{O}_{m+1}^{2\kappa}}\lesssim \varepsilon\Big(N_{m}^{2\overline\sigma-\frac43a_{1}}+ N_{m}^{s_0+2\overline{\sigma}+\frac{2}{3}\mu_1+1-b_1}
+ {N_{m}^{\overline\sigma+2-\mu_{2}-\frac23a_1}}\Big).
\end{equation}
For $m=0$ we obtain from slight modifications of \eqref{IJKL1} and \eqref{estimate mathcalF(U0)} 
$$
 \mathtt{I}^{(0)}_1(s_0) 
 \lesssim \varepsilon^2 \kappa^{-1 }.
 $$
  Similarly, we get
$$
   \mathtt{I}^{(0)}_2(s_0)  \lesssim 
\varepsilon \kappa^{- 1}
$$
and since $\varepsilon\kappa^{-3}\lesssim1$ in view of \eqref{small-C3}, then 
  \begin{align*}
\mathtt{I}^{(0)}_3(s_0) & \lesssim\, \varepsilon \kappa^{-2}N_0^{s_0-b_1}+ \varepsilon^2\kappa^{-5} \\
&\lesssim  \varepsilon\kappa^{-2}.
\end{align*}
Hence we find
  
\begin{align}\label{final estimate term L0T0-Id}
 \|\Pi_{0}(L_{0}\mathbf{T}_{0}-\textnormal{Id})\Pi_{0}\mathcal{F}({U}_{0})\|_{s_{0}}^{q,\kappa,\mathrm{O}_{1}^{2\kappa}}&\lesssim\varepsilon\kappa^{-2}.
\end{align}

\smallskip

$\bullet$  \textit{Estimate of $\big(L_{m}\mathbf{\Pi}_{m}^{\perp}-\Pi_{m}^{\perp}L_{m}\big)\mathbf{T}_{m}\Pi_{m}\mathcal{F}({U}_{m}).$}
Using \eqref{NLeq-eps} one finds  for $H=(\widehat{\mathfrak{I}},\widehat{\mathtt{c}}),$
\begin{align}\label{MTU1}
L_{m}H=\omega\cdot\partial_{\varphi}\widehat{\mathfrak{I}}-(0,0,\partial_{\theta}\mathrm{L}(\alpha)\widehat{w})-\varepsilon d_{i}X_{\mathcal{P}_{\varepsilon}}(i_{n})\widehat{\mathfrak{I}}-(\widehat{\mathtt{c}},0,0).
\end{align}
As $\mathrm{L}(\alpha)$ is a Fourier multiplier  then  $[\Pi_{m}^{\perp},\mathrm{L}(\alpha)]=0$. Therefore we infer from \eqref{and-bp1} 
$$
\left(L_{m}\mathbf{\Pi}_{m}^{\perp}-\Pi_{m}^{\perp}L_{m}\right)H=-\varepsilon\,\big[d_{i}X_{{P}_{\varepsilon}}\big(i_{m},\mathtt{c}_m\big),\Pi_{m}^{\perp}\big]\widehat{\mathfrak{I}}.
$$
According to  Lemma \ref{tame estimates for the vector field XmathcalPvarepsilon}-{(ii)}, Lemma \ref{Lem-Rgv1}-(iv) and $(\mathcal{P}1)_{m}$ we obtain 
$$
\|\left(L_{m}\mathbf{\Pi}_{m}^{\perp}-\Pi_{m}^{\perp}L_{m}\right)H\|_{s_{0}}^{q,\kappa,\mathrm{O}_{m+1}^{2\kappa}}\lesssim\varepsilon N_{m}^{s_0-b_{1}}\left(\|\widehat{\mathfrak{I}}\|_{b_{1}+1}^{q,\kappa,\mathrm{O}_{m}^{\kappa}}+\|{\mathfrak{I}}_{m}\|_{b_{1}+\overline\sigma}^{q,\kappa,\mathcal{O}}\|\widehat{\mathfrak{I}}\|_{s_{0}+1}^{q,\kappa,\mathrm{O}_{m}^{\kappa}}\right).
$$
It follows that
\begin{align*}
\mathtt{I}_4^{(m)}\triangleq \big\|\big(L_{m}\mathbf{\Pi}_{m}^{\perp}-\Pi_{m}^{\perp}L_{m}\big)\mathbf{T}_{m}\Pi_{m}\mathcal{F}({U}_{m})\big\|_{s_{0}}^{q,\kappa,\mathrm{O}_{m}^{\kappa}}&\lesssim\varepsilon N_{m}^{s_0-b_{1}}\|\mathbf{T}_{m}\Pi_{m}\mathcal{F}({U}_{m})\|_{b_{1}+1}^{\gamma,\mathrm{O}_{m}^{\kappa}}\\
&\quad+\varepsilon N_{m}^{s_0-b_{1}}\|{\mathfrak{I}}_{m}\|_{b_{1}+\overline\sigma}^{q,\kappa,\mathcal{O}}\|\mathbf{T}_{m}\Pi_{m}\mathcal{F}({U}_{m})\|_{s_{0}+1}^{q,\kappa,\mathrm{O}_{m}^{\kappa}}.
\end{align*}
Therefore putting together \eqref{estimate Tm}  with  Lemma \ref{orthog-Lem1}, Sobolev embeddings, \eqref{choice of gamma and N0 in the Nash-Moser}  and $(\mathcal{P}1)_{m}$ allows to get 
\begin{align*}
\mathtt{I}_4^{(m)}&
\lesssim \varepsilon\kappa^{-1}N_{m}^{s_0-b_{1}}\left[\|\Pi_{m}\mathcal{F}({U}_{m})\|_{b_{1}+\overline{\sigma}+1}^{q,\kappa,\mathrm{O}_{m}^{\kappa}}+\|{\mathfrak{I}}_{m}\|_{b_{1}+\overline{\sigma}+1}^{q,\kappa,\mathcal{O}}\|\Pi_{m}\mathcal{F}({U}_{m})\|_{s_{0}+\overline{\sigma}}^{q,\kappa,\mathrm{O}_{m}^{\kappa}}\right]\\
&\quad+\varepsilon\kappa^{-1}N_{m}^{s_0-b_{1}}\|{\mathfrak{I}}_{m}\|_{b_{1}+\overline\sigma}\left(\|\Pi_{m}\mathcal{F}({U}_{m})\|_{s_{0}+\overline{\sigma}+1}^{q,\kappa,\mathrm{O}_{m}^{\kappa}}+\|{\mathfrak{I}}_{m}\|_{s_{0}+\overline{\sigma}+1}^{q,\kappa,\mathcal{O}}\|\Pi_{m}\mathcal{F}({U}_{m})\|_{s_{0}+\overline{\sigma}}^{q,\kappa,\mathrm{O}_{m}^{\kappa}}\right)\\
&\lesssim \varepsilon\, N_{m}^{s_0+2-b_{1}}\left(\|\mathcal{F}({U}_{m})\|_{b_{1}+\sigma}^{q,\kappa,\mathrm{O}_{m}^{\kappa}}+\| {W}_{m}\|_{q,b_{1}+\overline\sigma}^{\gamma,\mathcal{O}}\|\Pi_{m}\mathcal{F}({U}_{m})\|_{s_{0}+\overline{\sigma}}^{q,\kappa,\mathrm{O}_{m}^{\kappa}}\right).
\end{align*}

Applying Lemma \ref{orthog-Lem1} and $(\mathcal{P}2)_{m}$, we find for $m\in\N$
\begin{align*}
\|\Pi_{m}\mathcal{F}(\widetilde{U}_{m})\|_{s_{0}+{\sigma}}^{q,\kappa,\mathrm{O}_{m}^{\kappa}}&\lesssim\varepsilon N_{m}^{\overline{\sigma}}N_{m-1}^{-a_{1}}\\
&\lesssim \varepsilon N_{m}^{\overline{\sigma}-\frac23 a_1}.
\end{align*}
Combining this estimate with  \eqref{HDP10} and $(\mathcal{P}3)_{m}$, we obtain for $m\in\N$
\begin{equation}\label{final estimate commutator}
 \| (L_{m}{\Pi}_{m}^{\perp}-\Pi_{m}^{\perp}L_{m})\mathbf{T}_{m}\Pi_{m}\mathcal{F}({U}_{m})\|_{s_{0}}^{q,\kappa,\mathrm{O}_{m+1}^{2\kappa}}
 \leqslant  C_*\varepsilon N_{m}^{s_0+\overline{\sigma}+\frac{2}{3}\mu_1+3-b_{1}}.
\end{equation}

\smallskip

$\bullet$  \textit{Estimates if  $Q_{m}$.} Coming back to \eqref{Def-Qm} and using 
 Taylor Formula at the second order, one finds
$$Q_{m}=\int_{0}^{1}(1-t)d_{i,\alpha}^{2}\mathcal{F}({U}_{m}+t\widetilde{H}_{m+1})[\widetilde{H}_{m+1},\widetilde{H}_{m+1}]dt.
$$
Hence  \eqref{MTU1} and Lemma \ref{tame estimates for the vector field XmathcalPvarepsilon}-{(iii)} allow to get
\begin{align}\label{mahma-YDa1}
\| Q_{m}\|_{s_{0}}^{q,\kappa,\mathrm{O}_{m+1}^{2\kappa}}\lesssim\varepsilon\left(1+\|{W}_{m}\|_{s_{0}+2}^{q,\kappa,\mathcal{O}}+\| \widetilde{H}_{m+1}\|_{s_{0}+2}^{q,\kappa,\mathrm{O}_{m}^{\kappa}}\right)\left(\| \widetilde{H}_{m+1}\|_{s_{0}+2}^{q,\kappa,\mathrm{O}_{m}^{\kappa}}\right)^{2}.
\end{align}
Combining \eqref{estimate Tm}, \eqref{link mathcalF(Um) and Wm} and \eqref{estimate mathcalF(U0) in norm s0}, we have for all $s\in[s_{0},S],$
\begin{equation}\label{link Hm+1 and Wm}
\begin{aligned}
\|\widetilde{H}_{m+1}\|_{s}^{q,\kappa,\mathrm{O}_{m+1}^{2\kappa}} & =  \|{\Pi}_{m}\mathbf{T}_{m}\Pi_{m}\mathcal{F}({U}_{m})\|_{s}^{q,\kappa,\mathrm{O}_{m+1}^{2\kappa}}\\
& \lesssim  \kappa^{-1}\left(\|\Pi_{m}\mathcal{F}({U}_{m})\|_{s+\overline{\sigma}}^{q,\kappa,\mathrm{O}_{m}^{\kappa}}+\|{\mathfrak{I}}_{m}\|_{s+\overline{\sigma}}^{q,\kappa,\mathcal{O}}\|\Pi_{m}\mathcal{F}({U}_{m})\|_{s_{0}+\overline{\sigma}}^{q,\kappa,\mathrm{O}_{m}^{\kappa}}\right)\\
& \lesssim   \kappa^{-1}\left(N_{m}^{\overline{\sigma}}\|\mathcal{F}({U}_{m})\|_{s}^{q,\kappa,\mathrm{O}_{m}^{\kappa}}+N_{m}^{2\overline{\sigma}}\|{\mathfrak{I}}_{m}\|_{s}^{q,\kappa,\mathcal{O}}\|\mathcal{F}({U}_{m})\|_{s_{0}}^{q,\kappa,\mathrm{O}_{m}^{\kappa}}\right)\\
& \lesssim   \kappa^{-1}N_{m}^{2\overline{\sigma}}\left(\varepsilon+\| {W}_{m}\|_{s}^{q,\kappa,\mathcal{O}}\right).
\end{aligned}
\end{equation}
Similarly, in view of  \eqref{Inclu-z0}, \eqref{estimate Tm in norm s0}, $(\mathcal{P}1)_{m}$ and $(\mathcal{P}2)_{m}$,  we get for $m\in\N$
\begin{align}\label{link Hm+1 and mathcalF(Um) in norm s0}
\nonumber \| \widetilde{H}_{m+1}\|_{s_{0}}^{q,\kappa,\mathrm{O}_{m+1}^{2\kappa}}&\lesssim \kappa^{-1}N_{m}^{\overline{\sigma}}\|\mathcal{F}({U}_{m})\|_{s_{0}}^{q,\kappa,\mathrm{O}_{m}^{\kappa}}\\
&\lesssim \varepsilon\kappa^{-1} N_{m}^{\overline{\sigma}}N_{m-1}^{-a_{1}}.
\end{align}
From  \eqref{imp-yi1},  $(\mathcal{P}1)_{m}$ and \eqref{link Hm+1 and mathcalF(Um) in norm s0},  we find for  $m\in\N$
\begin{align*}
\|{W}_{m}\|_{s_{0}+2}^{q,\kappa,\mathcal{O}}+\| \widetilde{H}_{m+1}\|_{s_{0}+2}^{q,\kappa,\mathrm{O}_{m+1}^{2\kappa}} & \leqslant  1+N_{m}^{2}\| H_{m+1}\|_{s_{0}}^{q,\kappa,\mathrm{O}_{m+1}^{2\kappa}}\\
& \leqslant 1+C\varepsilon\kappa^{-1}N_{m}^{\overline{\sigma}+2}N_{m-1}^{-a_{1}}\\
& \leqslant  1+C\,\varepsilon N_{m-1}^{3+\frac32\overline{\sigma}-a_{1}}.
\end{align*}
By virtue of \eqref{Assump-DRX} one has
\begin{align}\label{kita1}
a_1\geqslant 3+\tfrac32{\sigma}
\end{align}
and using  \eqref{choice of gamma and N0 in the Nash-Moser}   we may write for small $\varepsilon$
$$\|{W}_{m}\|_{s_{0}+2}^{q,\kappa,\mathcal{O}}+\| \widetilde{H}_{m+1}\|_{s_{0}+2}^{q,\kappa,\mathrm{O}_{m+1}^{2\kappa}}\leqslant 2.$$
Thus, by  inserting this estimate and \eqref{link Hm+1 and mathcalF(Um) in norm s0} into \eqref{mahma-YDa1} we deduce that
\begin{align*}
\| Q_{m}\|_{s_{0}}^{q,\kappa,\mathrm{O}_{m+1}^{2\kappa}} & \lesssim  \varepsilon\left(\| H_{m+1}\|_{s_{0}+2}^{q,\kappa,\mathrm{O}_{m+1}^{2\kappa}}\right)^{2}\\
& \lesssim  \varepsilon N_{m}^{4}\left(\| H_{m+1}\|_{s_{0}}^{q,\kappa,\mathrm{O}_{m+1}^{2\kappa}}\right)^{2}\\
& \lesssim  \varepsilon^3\kappa^{-2} N_{m}^{2\overline{\sigma}+4}N_{m-1}^{-2a_1}.
\end{align*}
Consequently,  we find according to  \eqref{definition of Nm} and $\varepsilon\kappa^{-1}\lesssim1$ that \mbox{for $m\geqslant 1$}
\begin{align}\label{final estimate for Qm}
\nonumber \|Q_{m}\|_{s_{0}}^{q,\kappa,\mathrm{O}_{m}^{\kappa}}&\lesssim\,\varepsilon^3\kappa^{-2} N_{m}^{2\overline{\sigma}+4-\frac{4}{3}a_{1}}\\
&\lesssim\,\varepsilon N_{m}^{2\overline{\sigma}+4-\frac{4}{3}a_{1}}.
\end{align}
As to the case  $m=0$, we come back to \eqref{link Hm+1 and Wm} and make suitable adjustments using in \mbox{particular \eqref{estimate mathcalF(U0)}}
\begin{equation}\label{link-H0}
\begin{array}{rcl}
\| \widetilde{H}_{1}\|_{s}^{q,\kappa,\mathrm{O}_{0}^{\kappa}} & \lesssim & \kappa^{-1}\|\Pi_{0}\mathcal{F}({U}_{0})\|_{s+\overline{\sigma}}^{q,\kappa,\mathrm{O}_{0}^{\kappa}}\\
& \lesssim &  \varepsilon\kappa^{-1}.
\end{array}
\end{equation}
Consequently, the inequality \eqref{final estimate for Qm} becomes for $m=0$,
 \begin{equation}\label{final-Q0}
\|Q_{0}\|_{s_{0}}^{q,\kappa,\mathrm{O}_{0}^{\kappa}}\lesssim\, \varepsilon^3\kappa^{-2}.
\end{equation}

\smallskip

$\bullet$ \textit{Parameters constraints.}
	Plugging  \eqref{final estimate PiperpF(Um)}, \eqref{Rent-M01}, \eqref{final estimate commutator} and \eqref{final estimate for Qm}, into \eqref{Decom-RTT1} yields \mbox{for $m\geqslant 1$}
\begin{align*}
\|\mathcal{F}(\widetilde{U}_{m+1})\|_{q,s_{0}}^{\gamma,\mathrm{O}_{m+1}^{2\kappa}}&\leqslant C C_*  \varepsilon\Big(N_{m}^{s_0+2\overline{\sigma}+3+\frac23\mu_1-b_{1}}+N_{m}^{\overline{\sigma}+2-\frac{2}{3}a_{1}-\mu_{2}}+N_{m}^{2\overline{\sigma}+4-\frac{4}{3}a_{1}}\Big).
\end{align*}
We want to check that with the assumptions fixed in  \eqref{Assump-DRX}  one gets the conditions,
$$\left\lbrace\begin{array}{rcl}
CN_{m}^{s_0+2\overline{\sigma}+\frac{2}{3}\mu_1+3-b_{1}} & \leqslant & \frac{1}{3}N_{m}^{-a_{1}}\\
CN_{m}^{\overline{\sigma}+2-\frac{2}{3}a_{1}-\mu_{2}}& \leqslant & \frac{1}{3}N_{m}^{-a_{1}}\\
CN_{m}^{2{\overline\sigma}+4-\frac{4}{3}a_{1}} & \leqslant & \frac{1}{3}N_{m}^{-a_{1}}.
\end{array}\right.$$
Actually, they are  satisfied  by  taking $N_{0}$ large enough, that is $\varepsilon$ small enough, provided that
\begin{equation}\label{Assump-DR1}
\left\lbrace\begin{array}{rcl}
s_0+2\overline{\sigma}+\frac{2}{3}\mu_1+4+a_{1} & \leqslant & b_{1}\\
\overline{\sigma}+\frac{1}{3}a_{1}{+3} & \leqslant & \mu_{2}\\
2\overline{\sigma}+5& \leqslant & \frac{1}{3}a_{1}.
\end{array}\right.
\end{equation}
Hence it is easy to verify that \eqref{Assump-DR1} follows  immediately from \eqref{Assump-DRX}.
Thus, we deduce for  $m\geqslant 1$
\begin{align}\label{Dabdoub1}
\nonumber\|\mathcal{F}(\widetilde{U}_{m+1})\|_{q,s_{0}}^{\gamma,\mathrm{O}_{m+1}^{2\kappa}}&\leqslant C\,C_{\ast}N_0^{-1}\varepsilon N_{m}^{-a_{1}}\\
&\leqslant C_{\ast}\varepsilon N_{m}^{-a_{1}}.
\end{align}
Concerning the case $m=0$, we insert \eqref{final estimate PiperpF(Um)}, \eqref{final estimate term L0T0-Id}, \eqref{final estimate commutator} and \eqref{final-Q0} into \eqref{Decom-RTT1} 
\begin{align*}
\|\mathcal{F}(\widetilde{U}_{1})\|_{q,s_{0}}^{\kappa,\mathrm{O}_{1}^{2\kappa}}&\leqslant CC_*\varepsilon\Big(N_{0}^{s_0+\overline{\sigma}+3+\frac23\mu_1-b_{1}}+\varepsilon\kappa^{-2}  +\varepsilon^3\kappa^{-2}\Big).\end{align*}
Then by virtue of \eqref{choice of gamma and N0 in the Nash-Moser} and the assumption on $b_1$ in \eqref{Assump-DR1} we get for $\varepsilon$ small enough
\begin{align*}
 C\Big(N_{0}^{s_0+\overline{\sigma}+3+\frac23\mu_1-b_{1}}+\varepsilon\kappa^{-2}  +\varepsilon^3\kappa^{-2}\Big) \leqslant 1
 \end{align*}
 leading to
 \begin{align*}
\|\mathcal{F}({U}_{1})\|_{q,s_{0}}^{\gamma,\mathrm{O}_{1}^{2\kappa}}&\leqslant C_*\varepsilon.
\end{align*}
$\blacktriangleright$ \textbf{Extension and achievement  of $(\mathcal{P}1)_{m+1}$-$(\mathcal{P}2)_{m+1}$-$(\mathcal{P}3)_{m+1}.$} We shall now extend $\widetilde{H}_{m+1}$ to the whole set of parameters $\mathcal{O}$ using a suitable cut-off function. For this purpose, we consider a $C^{\infty}$ cut-off function $\chi_{m+1}:\mathcal{O}\rightarrow[0,1]$ defined by
$$\chi_{m+1}(\lambda)=\left\lbrace\begin{array}{ll}
1 & \textnormal{in }\quad\mathrm{O}_{m+1}^{\kappa}\\
0 & \textnormal{in }\quad\mathcal{O}\backslash\mathrm{O}_{m+1}^{2\kappa}
\end{array}\right.$$
and satisfying the additional conditions
\begin{align}\label{Blow-loc1}
\forall \alpha\in\mathbb{N}^{d},\,|\alpha|\leqslant q,\quad \|\partial_{\lambda}^{\alpha}\chi_{m+1}\|_{L^{\infty}(\mathcal{O})}\lesssim\left(\kappa^{-1}N_{m}^{\overline{a}}\right)^{|\alpha|}.
\end{align}
Then we  define the extension ${H}_{m+1}$ of $\widetilde{H}_{m+1}$ by
$$
{H}_{m+1}(\lambda)=\left\lbrace\begin{array}{ll}
\chi_{m+1}(\lambda)\widetilde{H}_{m+1}(\lambda) &\quad  \textnormal{in }\quad \mathrm{O}_{m+1}^{2\kappa}\\
0 & \quad \textnormal{in }\quad\mathcal{O}\backslash\mathrm{O}_{m+1}^{2\kappa}.
\end{array}\right.
$$
and the extension $U_{m+1}$ is given by
$$
U_{m+1}\triangleq U_m+H_{m+1}.
$$
Notice that
$$
H_{m+1}=\widetilde{H}_{m+1},\quad\mathcal{F}(U_{m+1})=\mathcal{F}(\widetilde{U}_{m+1}) \quad\textnormal{in}\quad \mathrm{O}_{m+1}^{\kappa}.
$$
Therefoe we infer from \eqref{Dabdoub1}
\begin{align*}
\|\mathcal{F}({U}_{m+1})\|_{q,s_{0}}^{\gamma,\mathrm{O}_{m+1}^{\kappa}}&\leqslant C_{\ast}\varepsilon N_{m}^{-a_{1}},
\end{align*}
which achieves the induction of $(\mathcal{P}2)_{m+1}$.\\
 Next, using the law product in Lemma \ref{Law-prodX1} together with \eqref{Blow-loc1} we find
\begin{equation}\label{link Hm+1tilde Hm+1}
\forall s\in[s_{0},S],\quad\|{H}_{m+1}\|_{s}^{q,\kappa,\mathcal{O}}\lesssim N_{m}^{q\overline{a}}\|\widetilde{H}_{m+1}\|_{q,s}^{\gamma,\mathrm{O}_{m+1}^{2\kappa}}.
\end{equation}
Combining  \eqref{link Hm+1tilde Hm+1} and \eqref{link Hm+1 and mathcalF(Um) in norm s0} we deduce for $m\geqslant 1$
\begin{align*}
 \| {H}_{m+1}\|_{s_{0}+\overline{\sigma}}^{q,\kappa,\mathcal{O}}&\leqslant N_{m}^{q\overline{a}+\overline\sigma}\|\widetilde{H}_{m+1}\|_{s_{0}}^{q,\kappa,\mathrm{O}_{m+1}^{2\kappa}}\\
&\leqslant C C_{\ast}\varepsilon\kappa^{-1} N_{m}^{q\overline{a}+2\overline{\sigma}-\frac{2}{3}a_{1}}.
\end{align*}
From the choice done in \eqref{Assump-DRX} one gets  
\begin{align}\label{ConDDR4}
a_{2}= \tfrac{2}{3}a_{1}-q\overline{a}-2\overline{\sigma}-1\geqslant 1
\end{align}
and by choosing $\varepsilon$ small enough we obtain
\begin{align}\label{douza1}
\nonumber \|{H}_{m+1}\|_{s_{0}+\overline{\sigma}}^{q,\kappa,\mathcal{O}}&\leqslant CN_0^{-1}C_{\ast}\varepsilon\kappa^{-1}N_{m}^{-a_{2}}\\
&\leqslant C_{\ast}\varepsilon\kappa^{-1}N_{m}^{-a_{2}}.
\end{align}
For $m=0$ we may combine \eqref{link Hm+1tilde Hm+1} and \eqref{link-H0} in order to get
\begin{align}\label{douzami1}\|{H}_{1}\|_{s_{0}+\overline{\sigma}}^{q,\kappa}\leqslant \tfrac12C_{\ast}\varepsilon\kappa^{-1}N_{0}^{q\overline{a}}.
\end{align}
Let us define 
\begin{align}\label{Constru-1}
{W}_{m+1}\triangleq{W}_{m}+{H}_{m+1},
\end{align} 
then it is obvious  that $U_{m+1}=U_0+W_{m+1}$. In addition,  applying  $(\mathcal{P}1)_{m}$ with \eqref{douza1} and \eqref{douzami1} we find in view of Lemma \ref{lemma sum Nn} and for small $\varepsilon$
\begin{align*}
\|{W}_{m+1}\|_{s_{0}+\overline\sigma}^{q,\kappa}&\leqslant\|{H}_{1}\|_{s_{0}+\overline\sigma}^{q,\kappa}+\sum_{k=2}^{m+1}\|{H}_{k}\|_{s_{0}+\overline\sigma}^{q,\kappa}\\
&\leqslant \tfrac12C_*\varepsilon\kappa^{-1}N_0^{q\overline{a}}+C_*\varepsilon\kappa^{-1}\sum_{k=0}^{\infty}N_{k}^{-1}\\
&\leqslant \tfrac12C_* \varepsilon\kappa^{-1}N_0^{q\overline{a}}+CC_*\varepsilon\kappa^{-1}N_0^{-1}\leqslant C_* \varepsilon\kappa^{-1}N_0^{q\overline{a}}.
\end{align*}
Now, applying 
 \eqref{link Hm+1 and Wm}, \eqref{link Hm+1tilde Hm+1} and $(\mathcal{P}3)_{m}$, we deduce that
\begin{align*}
\|{W}_{m+1}\|_{b_{1}+\overline\sigma}^{q,\kappa} & \leqslant  \| {W}_{m}\|_{b_{1}+\overline\sigma}^{q,\kappa}+C N_m^{q\overline{a}}\|{H}_{m+1}\|_{b_{1}+\overline\sigma}^{q,\kappa}\\
& \leqslant   C_*\varepsilon\kappa^{-1}N_{m-1}^{\mu_1}+\kappa^{-1}N_{m}^{q\overline{a}+2\overline{\sigma}}\left(\varepsilon+\| {W}_{m}\|_{b_{1}+\overline\sigma}^{q,\kappa,\mathcal{O}_{m}^{\gamma}}\right)\\
& \leqslant  CC_{\ast}\varepsilon\kappa^{-1}N_{m}^{q\overline{a}+2\overline{\sigma}+1+\frac{2}{3}\mu_1}.
\end{align*}
According to \eqref{Assump-DRX} one may check that
\begin{align}\label{ConDDR2}
q\overline{a}+2\overline{\sigma}+2= \tfrac{\mu_1}{3}
\end{align}
and therefore we deduce
\begin{align*}
\|{W}_{m+1}\|_{b_{1}+\overline\sigma}^{q,\kappa,\mathcal{O}}&\leqslant CN_0^{-1}C_{\ast}\varepsilon\kappa^{-1}N_{m}^{\mu_1}\\
&\leqslant C_{\ast}\varepsilon\kappa^{-1}N_{m}^{\mu_1}
\end{align*}
for $\varepsilon$ small enough and \eqref{modon1} and 
this achieves the proof of $(\mathcal{P}3)_{m+1}.$ 
This achieves  the proof of $(\mathcal{P}1)_{m+1}$ and therefore  the proof of Proposition \ref{Nash-Moser} is now achieved.
\end{proof}
As a by-product of Proposition \ref{Nash-Moser} we shall construct solutions  to  the nonlinear equation \eqref{eq-NM} provided that the parameters belong to a suitable Cantor set. Later, we shall investigate in Section \ref{Section 6.2}  the Lebesgue measure of this set. 
\begin{coro}\label{cor-ima12}
There exists $\varepsilon_0>0$ such that for any $\varepsilon\in(0,\varepsilon_0)$ the following assertions hold true.  Consider the Cantor set (depending implicitly in $\varepsilon$)
$$\mathtt{G}_{\infty}^{\kappa}\triangleq\bigcap_{m\in\mathbb{N}}\mathcal{A}_{m}^{\kappa}.$$
 There exists  a function (depending implicitly on $\varepsilon$)
$$U_{\infty}:\begin{array}[t]{rcl}
\mathcal{O} & \rightarrow & \left(\mathbb{T}^{d}\times\mathbb{R}^{d}\times H_{\mathbb{S}}^{\perp}\right)\times\mathbb{R}^{d}\\
\lambda=(\omega,\alpha) & \mapsto & (i_{\infty}(\lambda),\mathtt{c}_{\infty}(\lambda))
\end{array}$$
such that 
$$\forall\,\lambda\in\mathtt{G}_{\infty}^{\kappa}\quad\mathcal{F}\big(U_{\infty}(\lambda),\lambda,\varepsilon\big)=0.$$
In addition, the function  $\lambda\in\mathcal{O}\mapsto \mathtt{c}_{\infty}(\lambda,\varepsilon)$ belongs to $W^{q,\infty}(\mathcal{O})$ with  
 \begin{equation*}
\mathtt{c}_{\infty}(\lambda,\varepsilon)=\omega+\mathrm{r}_{\varepsilon}(\lambda)\quad\mbox{ and }\quad \|\mathrm{r}_{\varepsilon}\|^{q,\kappa}\lesssim\varepsilon\kappa^{-1}{N_0^{q\overline{a}}}.
\end{equation*}
Moreover, there exists a function $\alpha\in(\underline\alpha,\overline\alpha)\mapsto  \lambda(\alpha,\varepsilon)\triangleq \big({\omega}(\alpha,\varepsilon),\alpha\big)$ in  $W^{q,\infty}$ with 
\begin{equation}\label{estimate repsilon1}
{\omega}(\alpha,\varepsilon)=-{\omega}_{\textnormal{Eq}}(\alpha)+\overline{\mathrm{r}}_{\varepsilon}(\alpha),\quad \quad \|\overline{\mathrm{r}}_{\varepsilon}\|^{q,\kappa}\lesssim\varepsilon \kappa^{-1}{N_0^{q\overline{a}}}
\end{equation}
and
$$\forall\,\alpha\in\mathtt{C}_{\infty}^{\kappa,\varepsilon},\quad\mathcal{F}\big(U_{\infty}(\lambda(\alpha,\varepsilon),\varepsilon),\lambda(\alpha,\varepsilon),\varepsilon\big)=0,\quad \mathtt{c}_{\infty}(\lambda(\alpha,\varepsilon),\varepsilon)=-{\omega}_{\textnormal{Eq}}(\alpha),
$$
where the Cantor set $\mathtt{C}_{\infty}^{\kappa,\varepsilon}$ is defined by
$$
\mathtt{C}_{\infty}^{\kappa,\varepsilon}=\Big\{\alpha\in(\underline\alpha,\overline\alpha),   \big( \omega(\alpha,\varepsilon),  \alpha\big) \in \mathtt{G}_{\infty}^{\kappa}\Big\}.
$$

\end{coro}
\begin{proof} Combining  \eqref{Constru-1} and \eqref{douza1}, we obtain 
$$\|{W}_{m+1}-{W}_{m}\|_{s_{0}}^{q,\kappa}=\|{H}_{m+1}\|_{s_{0}}^{q,\kappa}\leqslant\|{H}_{m+1}\|_{s_{0}+\overline{\sigma}}^{q,\kappa}\leqslant C_{\ast}\varepsilon\kappa^{-1}N_{m}^{-a_{2}}.$$
Hence, the sequence $\left({W}_{m}\right)_{m\in\mathbb{N}}$ is convergent and then we can define
\begin{align*}
W_{\infty}&\triangleq \lim_{m\rightarrow\infty}{W}_{m}\\
&\triangleq (\mathfrak{I}_{\infty},\mathtt{c}_{\infty}-\omega)
\end{align*}
and
$$U_{\infty}= (i_{\infty},\mathtt{c}_{\infty})\triangleq U_0+W_\infty=\Big((\varphi,0,0),\omega\Big)+W_{\infty}.
$$
According to  the point $(\mathcal{P}2)_{m}$ of Proposition \ref{Nash-Moser}, we have for small $\varepsilon$
$$\forall\,\lambda\in\mathtt{G}_{\infty}^{\gamma},\quad\mathcal{F}\big(U_{\infty}(\lambda,\varepsilon),\lambda,\varepsilon\big)=0.$$
We recall that $\mathcal{F}$ is given in \eqref{operatorF}  and the Cantor set depends also in $\varepsilon.$ Now, applying  the point $(\mathcal{P}1)_{m}$ of Proposition \ref{Nash-Moser}, we obtain,
 \begin{equation}\label{alpha infty}
\mathtt{c}_{\infty}(\lambda,\varepsilon)=\omega+\mathrm{r}_{\varepsilon}(\lambda)\quad\mbox{ with }\quad \|\mathrm{r}_{\varepsilon}\|^{q,\kappa}\lesssim\varepsilon\kappa^{-1}{N_0^{q\overline{a}}}.
\end{equation}
Let us move to the second result and check the existence of solutions to the original Hamiltonian equation. First recall that the open set $\mathcal{O}$ is given by
$$
\mathcal{O}=\mathcal{U}\times (\underline\alpha,\overline\alpha),\quad\hbox{with}\quad \mathcal{U}\triangleq B(0,R),
$$
with $R$ large enough such that the open  ball $B(0,R)$ contains  the equilibrium frequency vector $\big\{{\omega}_{\textnormal{Eq}}(\alpha),\alpha\in[\underline\alpha,\overline\alpha]\big\}$
Applying  \eqref{alpha infty}, we find that for any  $\alpha\in(\underline\alpha,\overline\alpha)$ the partial mapping  $\omega\mapsto \mathtt{c}_{\infty}(\omega,\alpha)$ is invertible form $\mathcal{U}$ into its image $\mathtt{c}_{\infty}(\mathcal{U},\alpha)$ and we have similarly to \eqref{IInv-77}  
$$\widehat\omega=\mathtt{c}_{\infty}(\omega,\alpha)=\omega+\mathrm{r}_{\varepsilon}(\omega,\alpha)\Longleftrightarrow\omega=\mathtt{c}_{\infty}^{-1}(\widehat\omega,\alpha)=\widehat\omega+\widehat{\mathrm{r}}_{\varepsilon}(\widehat\omega,\alpha).$$
Applying Lemma \ref{algeb1} we deduce that  $\widehat{\mathrm{r}}_{\varepsilon}$ satisfies the estimate
\begin{align}\label{estimate mathrm repsilon}
\nonumber\|\widehat{\mathrm{r}}_{\varepsilon}\|^{q,\kappa}&\lesssim\|{\mathrm{r}}_{\varepsilon}\|^{q,\kappa}\\
&\lesssim\varepsilon\kappa^{-1} {N_0^{q\overline{a}}}.
\end{align}
We denote 
$${\omega}(\alpha,\varepsilon)\triangleq\mathtt{c}_{\infty}^{-1}\big(-{\omega}_{\textnormal{Eq}}(\alpha),\alpha\big)=-{\omega}_{\textnormal{Eq}}(\alpha)+\overline{\mathrm{r}}_{\varepsilon}(\alpha)\quad \mbox{ with }\quad \overline{\mathrm{r}}_{\varepsilon}(\alpha)\triangleq\widehat{\mathrm{r}}_{\varepsilon}\big(-{\omega}_{\textnormal{Eq}}(\alpha),\alpha\big)
$$
and
\begin{align}\label{lambda-m1}
\lambda(\alpha,\varepsilon)\triangleq \big({\omega}(\alpha,\varepsilon),\alpha\big).
\end{align}
It follows that
$$\forall\,\alpha\in\mathtt{C}_{\infty}^{\kappa,\varepsilon},\quad\mathcal{F}\big(U_{\infty}(\lambda(\alpha,\varepsilon)),\lambda(\alpha,\varepsilon),\varepsilon\big)=0,
$$
where the Cantor set $\mathtt{C}_{\infty}^{\kappa,\varepsilon}$ was defined in Corollary \ref{cor-ima12}.
This gives nontrivial solutions for the original Hamiltonian equation provided that $\alpha\in\mathtt{C}_{\infty}^{\kappa,\varepsilon}.$
It remains to check the suitable estimates for the function $\overline{\mathrm{r}}_{\varepsilon}$. According to Lemma \ref{lem-asym}-(vi) we know that all  the derivatives $\partial_\alpha^j\omega_{\textnormal{Eq}}$ are uniformly bounded on the interval $(\underline\alpha,\overline\alpha),$ then by the chain rule,  \eqref{estimate mathrm repsilon} and  \eqref{choice of gamma and N0 in the Nash-Moser} we find
\begin{equation}\label{ouz-end}
\|\overline{\mathrm{r}}_{\varepsilon}\|^{q,\kappa}\lesssim\varepsilon\kappa^{-1} {N_0^{q\overline{a}}}\quad\hbox{and}\quad \|\omega(\cdot,\varepsilon)\|^{q,\kappa}\lesssim 1+\varepsilon\kappa^{-1} {N_0^{q\overline{a}}}\lesssim 1.
\end{equation}
The proof of Corollary \ref{cor-ima12} is now complete.
\end{proof}

\subsection{Final Cantor set  estimates}\label{Section 6.2}
The goal of this section is to give a lower bound Lebesgue measure for  the Cantor set $\mathtt{C}_{\infty}^{\kappa,\varepsilon}$ constructed in  Corollary \ref{cor-ima12} and show that  when $\varepsilon\to0$   it will  get a full Lebesgue measure.
From Corollary \ref{cor-ima12}, the Cantor set $\mathtt{C}_{\infty}^{\kappa,\varepsilon}$ can be written in the form 
\begin{equation}\label{definition of the final Cantor set}
\mathtt{C}_{\infty}^{\kappa,\varepsilon}=\bigcap_{m\in\mathbb{N}}\mathtt{C}_{m}^{\kappa}\quad \mbox{ where }\quad \mathtt{C}_{m}^{\kappa}=\Big\{\alpha\in(\underline\alpha,\overline\alpha)\quad\hbox{s.t.}\quad\lambda(\alpha,\varepsilon)\in\mathcal{A}_m^{\kappa}\Big\},
\end{equation}
where the sets $\mathcal{A}_m^{\kappa}$ are defined in Proposition \ref{Nash-Moser}.
Our result reads as follows
\begin{proposition}\label{lem-meas-es1}
Let $q_0$ as in Proposition $\ref{lemma transversality}$ and assume \eqref{Assump-DRX} and \eqref{choice of gamma and N0 in the Nash-Moser} with $q=q_0+1$. Assume in addition that  $\tau_{1},\tau_{2}$ and $\varrho$ satisfy
$$\left\lbrace\begin{array}{l}
\tau_{1}> d\,q_{0},
\tau_2>\frac{ \tau_1}{1-2\overline\alpha}+ d\,q_0,
\\
0<\varrho< \min\left(\tfrac{1}{2}-\overline{\alpha},\tfrac{1}{q_0+2}, 2\tfrac{1-a}{a}\right).
\end{array}\right.$$
Then, there exists $C>0$ and $\varepsilon_0$   such that  for any $\varepsilon\in(0,\varepsilon_0]$
$$
\big|\mathtt{C}_{\infty}^{\kappa,\varepsilon}\big|\geqslant \overline\alpha-\underline\alpha-C\varepsilon^{\frac{a\varrho}{q_0}}.
$$
In particular, $\displaystyle \lim_{\varepsilon\to0}\big|\mathtt{C}_{\infty}^{\kappa,\varepsilon}\big|=\overline\alpha-\underline\alpha.$
\end{proposition}
\begin{proof}
By \eqref{definition of the final Cantor set}, we can write 
\begin{align}\label{Mir-1}
\nonumber \left|(\underline\alpha,\overline\alpha)\backslash\mathtt{C}_{\infty}^{\kappa,\varepsilon}\right|&\leqslant\big|(\underline\alpha,\overline\alpha)\backslash\mathtt{C}_{0}^{\kappa}\big|+\sum_{m=0}^{\infty}\big|\mathtt{C}_{m}^{\kappa}\backslash\mathtt{C}_{m+1}^{\kappa}\big|\\
&\triangleq\sum_{m=0}^\infty\mathcal{S}_{m},
\end{align}
where we use the notation $|A|$ for the Lebesgue measure of a given measurable  set $A$ and the fact that by by construction  $\mathtt{C}_{0}^{\kappa}=(\underline\alpha,\overline\alpha)$. 
According to the notations introduced in  Proposition \ref{reduction of the remainder term} and  Proposition \ref{projection in the normal directions} combined with \eqref{lambda-m1} one may write 
\begin{align}\label{asy-z1}
\nonumber \mu_{j}^{\infty,m}(\alpha,\varepsilon)&\triangleq\mu_{j}^{\infty}\big(\lambda(\alpha,\varepsilon),i_{m}\big)\\
&=\Omega_{j}(\alpha)+j\,r^{1,m}(\alpha,\varepsilon)-
j\mathtt{W}(j,\alpha)\,r^{2,m}(\alpha,\varepsilon)+r_{j}^{\infty,m}(\alpha,\varepsilon)
\end{align}
with $\mathtt{W}(j,\alpha)$ is defined in Lemma \ref{lem-asym}, 
\begin{align*}
\ell\in\{1,2\},\quad r^{\ell,m}(\alpha,\varepsilon)&\triangleq r^{\ell}\big(\lambda(\alpha,\varepsilon),i_m\big)
\end{align*}
and 
$$r_{j}^{\infty,m}(\alpha,\varepsilon)\triangleq r_{j}^{\infty}\big(\lambda(\alpha,\varepsilon),i_{m}\big).$$
Applying  Proposition \ref{projection in the normal directions}-(i) combined with Proposition \ref{Nash-Moser}-$(\mathcal{P}1)_{m}$, \eqref{rent-P78}, the chain rule and \eqref{ouz-end} yield
\begin{align}\label{uniform estimate r1}
 \nonumber 0\leqslant k\leqslant q,\quad \sup_{m\in\mathbb{N}}\sup_{\alpha\in(\underline\alpha,\overline\alpha)}|\partial_\alpha^k r^{\ell,m}(\alpha,\varepsilon)|&\lesssim \sup_{m\in\mathbb{N}} {\varepsilon \kappa^{-1-k}\left(1+\| \mathfrak{I}_{m}\|_{\overline{s}_h}^{q,\kappa}\right)}\\
 &\lesssim\varepsilon\kappa^{-1-k}.
\end{align}
In a similar way, using  \eqref{estimate rjinfty}, \eqref{ouz-end} and Proposition \ref{Nash-Moser}-$(\mathcal{P}1)_{m}$, we find
\begin{align}\label{uniform estimate rjinfty}
\nonumber0\leqslant k\leqslant q,\quad\sup_{m\in\mathbb{N}}\sup_{j\in\mathbb{S}_{0}^{c}}\sup_{\alpha\in(\underline\alpha,\overline\alpha)}|j|^{1-\epsilon-2\overline\alpha}|\partial_\alpha^kr_{j}^{\infty,m}(\alpha,\varepsilon)|&\leqslant\kappa^{-k}\sup_{m\in\mathbb{N}}\sup_{j\in\mathbb{S}_{0}^{c}}|j|^{1-\epsilon-2\overline\alpha}\| r_{j}^{\infty}(i_{m})\|^{q,\kappa}\\&\lesssim\varepsilon\kappa^{-2-k}.
\end{align}
Coming back to \eqref{definition of the final Cantor set} and using the Cantor sets introduced in Proposition \ref{reduction of the remainder term}, Theorem \ref{inversion of the linearized operator in the normal directions} and Proposition \ref{QP-change} one obtains by construction that   for any $m\in\mathbb{N}$ 
\begin{equation}\label{set-U0}
\mathtt{C}_{m}^{\kappa}\backslash\mathtt{C}_{m+1}^{\kappa}=\bigcup_{{|l|\leqslant N_m},\atop j\in\mathbb{Z}\backslash\{0\}}\mathcal{R}_{l,j}^{(0)}(i_{m})\bigcup_{{|l|\leqslant N_m}\atop j,j_{0}\in\mathbb{S}_{0}^{c}}\mathcal{R}_{l,j,j_{0}}(i_{m})\bigcup_{{|l|\leqslant N_m}\atop j\in\mathbb{S}_{0}^{c}}\mathcal{R}_{l,j}^{(1)}(i_{m}),
\end{equation}
where
\begin{align*}
\mathcal{R}_{l,j}^{(0)}(i_{m})&\triangleq\left\lbrace\alpha\in\mathtt{C}_{m}^{\kappa};\;\,|{\omega}(\alpha,\varepsilon)\cdot l+jc_{m}(\alpha,\varepsilon)|\leqslant\tfrac{4\kappa_{m+1}^{\varrho}\langle j\rangle}{\langle l\rangle^{\tau_{1}}}\right\rbrace,\\
\mathcal{R}_{l,j,j_{0}}(i_{m})&\triangleq\left\lbrace\alpha\in\mathtt{C}_{m}^{\kappa};\;\,|{\omega}(\alpha,\varepsilon)\cdot l+\mu_{j}^{\infty,m}(\lambda,\varepsilon)-\mu_{j_{0}}^{\infty,m}(\lambda,\varepsilon)|\leqslant\tfrac{2\kappa_{m+1}\langle j-j_{0}\rangle}{\langle l\rangle^{\tau_{2}}}\right\rbrace
\\
\mathcal{R}_{l,j}^{(1)}(i_{m})&\triangleq\left\lbrace\alpha\in\mathtt{C}_{m}^{\kappa};\;\,|{\omega}(\alpha,\varepsilon)\cdot l+\mu_{j}^{\infty,m}(\alpha,\varepsilon)|\leqslant\tfrac{\kappa_{m+1}\langle j\rangle}{\langle l\rangle^{\tau_{1}}}\right\rbrace,
\end{align*}
with the following definition $c_{m}(\alpha,\varepsilon)\triangleq c\big(\lambda(\alpha,\varepsilon),i_m(\lambda(\alpha,\varepsilon),\varepsilon)\big)$, where the coefficient $c(\lambda,i_0)$ was introduced in \mbox{Proposition \ref{QP-change}.} 
Let us now move to the estimate of $\mathcal{S}_{m}$ defined in \eqref{Mir-1}. Then using \eqref{set-U0} together with  Lemma \ref{lemm-dix1}-(iv) and Lemma  \ref{some cantor set are empty}, we find
\begin{align*}
\mathcal{S}_{m}&\leqslant \sum_{\underset{|j|\leqslant C_{0}\langle l\rangle,|l|>N_{m-1}}{(l,j)\in\mathbb{Z}^{d+1}\backslash\{0\}}}\left|\mathcal{R}_{l,j}^{(0)}(i_{m})\right|+\sum_{\underset{\underset{\min(|j|,|j_{0}|)^{1-2\overline\alpha-\epsilon}\leqslant c_{2}\kappa^{-\varrho}\langle l\rangle^{\tau_{1}}}{|j-j_{0}|\leqslant C_{0}\langle l\rangle,|l|>N_{m-1}}}{(l,j,j_{0})\in\mathbb{Z}^{d}\times(\mathbb{S}_{0}^{c})^{2}}}\left|\mathcal{R}_{l,j,j_{0}}(i_{m})\right|\\
&\qquad +\sum_{\underset{|j|\leqslant C_{0}\langle l\rangle,|l|>N_{m-1}}{(l,j)\in\mathbb{Z}^{d}\times\mathbb{S}_{0}^{c}}}\left|\mathcal{R}_{l,j}^{(1)}(i_{m})\right|.
\end{align*}
By a direct application of Lemma \ref{Piralt} combined with  Lemma \ref{lemma R\"ussmann condition for the perturbed frequencies} we obain for any  $m\in\mathbb{N}$ 
\begin{equation}\begin{array}{l}\label{kio1}
|\mathcal{R}_{l,j}^{(0)}(i_{m})|\lesssim\kappa^{\frac{\varrho}{q_0}}\langle j\rangle^{\frac{1}{q_0}}\langle l\rangle^{-1-\frac{\tau_{1}+1}{q_0}}\\
|\mathcal{R}_{l,j}^{(1)}(i_{m})|\lesssim\kappa^{\frac{1}{q_0}}\langle j\rangle^{\frac{1}{q_0}}\langle l\rangle^{-1-\frac{\tau_{1}+1}{q_0}}\\
|\mathcal{R}_{l,j,j_{0}}(i_{m})|\lesssim\kappa^{\frac{1}{q_0}}\langle j-j_{0}\rangle ^{\frac{1}{q_0}}\langle l\rangle^{-1-\frac{\tau_{2}+1}{q_0}}.
\end{array}
\end{equation}
Notice that if $|j-j_{0}|\leqslant C_{0}\langle l\rangle$ and $\min(|j|,|j_{0}|)^{1-2\overline\alpha-\epsilon}\lesssim\kappa^{-\varrho}\langle l\rangle^{\tau_{1}}$, then 
$$
\max(|j|,|j_{0}|)=\min(|j|,|j_{0}|)+|j-j_{0}|\lesssim \left(\kappa^{-\varrho}\langle l\rangle^{\tau_{1}}\right)^{\frac{1}{1-2\overline\alpha-\epsilon}}+C_{0}\langle l\rangle\lesssim\kappa^{-\frac{\varrho}{1-2\overline\alpha-\epsilon}}\langle l\rangle^{\frac{\tau_{1}}{1-2\overline\alpha-\epsilon}}.
$$
Hence using \eqref{kio1} allows to get
\begin{align}\label{dtu-c1}
\nonumber	\sum_{m\in\mathbb{N}}{\mathcal{S}_{m}}&\lesssim  \kappa^{\frac{1}{q_0}}\Big(\sum_{m\in\NN\atop |l|>N_{m-1}}\langle l\rangle^{-\frac{\tau_{1}}{q_0}}+\kappa^{-\frac{\varrho}{q(1-2\overline\alpha-\epsilon)}}\sum_{m\in\N\atop |l|>N_{m-1}}\langle l\rangle^{\frac{\tau_1}{q_0(1-2\overline\alpha-\epsilon)}-\frac{\tau_{2}}{q_0}}\Big)
+\kappa^{\frac{\varrho}{q_0}}\sum_{m\in\N\atop|l|>N_{m-1}}\langle l\rangle^{-\frac{\tau_{1}}{q_0}}\\
 & \lesssim  \kappa^{\min\left(\frac{\varrho}{q_0},\frac{1}{q_0}-\frac{\varrho}{q_0(1-2\overline\alpha-\epsilon)}\right)},
\end{align}
provided that
\begin{align}\label{assum2-d}
\tau_1>d \,q_0\quad\hbox{and}\quad \tau_2>\tfrac{ \tau_1}{1-2\overline\alpha}+ d\,q_0,
\end{align}
and  $\epsilon$ is chosen small enough. Thus, by imposing
\begin{align}\label{pourq1}
0<\varrho<\tfrac{1}{2}-\overline{\alpha}
\end{align}
we obtain from \eqref{dtu-c1}
\begin{align}\label{set-U07}
\sum_{m\in\mathbb{N}}{\mathcal{S}_{m}}&\lesssim  \kappa^{\frac{\varrho}{q_0}}.
\end{align}
Inserting  \eqref{set-U07} into  \eqref{Mir-1} implies
$$
 \left|(\underline\alpha,\overline\alpha)\backslash\mathtt{C}_{\infty}^{\kappa,\varepsilon}\right|  \lesssim \kappa^{\tfrac{\varrho}{q_0}},
 $$
provided that the conditions \eqref{assum2-d} and \eqref{pourq1} are satisfied. 
By making the choice  
\begin{equation}\label{choice of upsilon}
\quad \varrho< \min\left(\tfrac{1}{2}-\overline{\alpha},\tfrac{1}{q_0+2}, 2\tfrac{1-a}{a}\right)
\end{equation}
which is compatible with the condition on $\varrho$ given by Proposition \ref{QP-change}
 with $q=q_0+1$, \eqref{pourq1} and Lemma \ref{lemm-dix1}, we find,  since $\kappa=\varepsilon^a$ according to \eqref{choice of gamma and N0 in the Nash-Moser},
$$ \left|(\underline\alpha,\overline\alpha)\backslash\mathtt{C}_{\infty}^{\kappa,\varepsilon}\right|\lesssim\varepsilon^{\frac{a\varrho}{q_0}}.$$
This completes the proof of the proposition. 
\end{proof}
It remains now to prove  Lemma \ref{lemm-dix1} and Lemma \ref{some cantor set are empty} below,   used in the proof of Proposition \ref{lem-meas-es1}.
\begin{lemma}\label{lemm-dix1}
Assume  \eqref{Assump-DRX},  \eqref{choice of gamma and N0 in the Nash-Moser} and $\varrho\in\left(0,\tfrac{2(1-a)}{a}\right)$. 
{Let $m\in\mathbb{N}^{*}$ and  $l\in\mathbb{Z}^{d}$ such that $|l|\leqslant N_{m-1}.$ Then the following assertions hold true.
\begin{enumerate}
\item For   $ j\in\mathbb{Z}$ with $(l,j)\neq(0,0)$, we get  $\,\,\mathcal{R}_{l,j}^{(0)}(i_{m})=\varnothing.$
\item For  $ (j,j_{0})\in(\mathbb{S}_{0}^{c})^{2}$ with $(l,j)\neq(0,j_0),$ we get $\,\,\mathcal{R}_{l,j,j_{0}}(i_{m})=\varnothing.$
\item For  $j\in\mathbb{S}_{0}^{c}$ with $(l,j)\neq(0,0)$, we get $\,\,\mathcal{R}_{l,j}^{(1)}(i_{m})=\varnothing.$
\item For any $m\in\mathbb{N},$
\begin{equation*}
\mathtt{C}_{m}^{\kappa}\backslash\mathtt{C}_{m+1}^{\kappa}=\bigcup_{\underset{|l|>N_{m-1}}{(l,j)\in\mathbb{Z}^{d}\times\mathbb{Z}\backslash\{(0,0)\}}}\mathcal{R}_{l,j}^{(0)}(i_{m})\cup\bigcup_{\underset{|l|>N_{m-1}}{(l,j,j_{0})\in\mathbb{Z}^{d}\times(\mathbb{S}_{0}^{c})^{2}}}\mathcal{R}_{l,j,j_{0}}(i_{m})\cup\bigcup_{\underset{|l|>N_{m-1}}{(l,j)\in\mathbb{Z}^{d}\times\mathbb{S}_{0}^{c}}}\mathcal{R}_{l,j}^{(1)}(i_{m}).
\end{equation*}
\end{enumerate}}
\end{lemma}
\begin{proof}
$\blacktriangleright$ The key estimate that will be used for all the points  comes from  \eqref{rent-P78}
\begin{align}\label{Bio-X1}
\nonumber \| i_{m}-i_{m-1}\|_{\overline{s}_{h}+\sigma_4}^{q,\kappa} &\leqslant\| U_{m}-U_{m-1}\|_{\overline{s}_{h}+\sigma_4}^{q,\kappa}=\| H_{m}\|_{\overline{s}_h+\sigma_4}^{q,\kappa}\\
&\leqslant C_{\ast}\varepsilon\kappa^{-1}N_{m-1}^{-a_{2}}.
\end{align}
\textbf{(i)} We shall first prove that if $|l|\leqslant N_{m-1}$ and $(l,j)\neq (0,0),$  then $\mathcal{R}_{l,j}^{(0)}(i_{m}) \subset\mathcal{R}_{l,j}^{(0)}(i_{m-1}).$ Assume for a while this inclusion and let us check how this implies that $\mathcal{R}_{l,j}^{(0)}(i_{m})=\varnothing.$ According to \eqref{set-U0} one deduces that 
$$
\mathcal{R}_{l,j}^{(0)}(i_{m}) \subset\mathcal{R}_{l,j}^{(0)}(i_{m-1})\subset \mathtt{C}_{m-1}^{\kappa}\backslash\mathtt{C}_{m}^{\kappa}.
$$
Using  once again \eqref{set-U0}  gives  o $\mathcal{R}_{l,j}^{(0)}(i_{m}) \subset\mathtt{C}_{m}^{\kappa}\backslash\mathtt{C}_{m+1}^{\kappa}$, and therefore  we get
$$
\mathcal{R}_{l,j}^{(0)}(i_{m}) \subset(\mathtt{C}_{m}^{\kappa}\backslash\mathtt{C}_{m+1}^{\kappa})\cap  (\mathtt{C}_{m-1}^{\kappa}\backslash\mathtt{C}_{m}^{\kappa})=\varnothing.
$$
Let us now turn to the proof of the claim $\mathcal{R}_{l,j}^{(0)}(i_{m}) \subset\mathcal{R}_{l,j}^{(0)}(i_{m-1}).$  To do that, consider a point  $\alpha\in\mathcal{R}_{l,j}^{(0)}(i_{m}).$ Then we get first by construction that  $\alpha\in\mathtt{C}_{m}^{\kappa}\subset \mathtt{C}_{m-1}^{\kappa}.$ In addition, by the triangle inequality and Sobolev embeddings we deduce 
\begin{align*}
\big|{\omega}(\alpha,\varepsilon)\cdot l+jc_{m-1}(\alpha,\varepsilon)\big| & \leqslant  \big|{\omega}(\alpha,\varepsilon)\cdot l+jc_{m}(\alpha,\varepsilon)\big|+|j|\big|c_{m}(\alpha,\varepsilon)-c_{m-1}(\alpha,\varepsilon)\big|\\
&\leqslant  4\tfrac{\kappa_{m+1}^{\varrho}|j|}{\langle l\rangle^{\tau_{1}}}+C|j|\|c(\cdot,{ i_{m}})-c(\cdot,i_{m-1})\|^{q,\kappa}.
\end{align*}
Thus applying \eqref{difference ci}, \eqref{Bio-X1} and \eqref{choice of gamma and N0 in the Nash-Moser} we obtain with the choice $\kappa=\varepsilon^a$
\begin{align*}
\big|{\omega}(\alpha,\varepsilon)\cdot l+jc_{m-1}(\alpha,\varepsilon)\big| & \leqslant   \tfrac{4\kappa_{m+1}^{\varrho}|j|}{\langle l\rangle^{\tau_{1}}}+C\varepsilon\kappa^{-1}|j|\| i_{m}-i_{m-1}\|_{\overline{s}_{h}+4}^{q,\kappa}\\
& \leqslant \tfrac{4\kappa_{m+1}^{\varrho}|j|}{\langle l\rangle^{\tau_{1}}}+C\varepsilon^{2(1-a)}|j|N_{m-1}^{-a_{2}}.
\end{align*}
According to the definition of $\kappa_m$ in Proposition \ref{Nash-Moser}-$(\mathcal{P}2)_m$ one gets
$$
\exists c_0>0,\quad \mbox{s.t}\quad \forall m\in \mathbb{N},\quad \kappa_{m+1}^{\varrho}-\kappa_{m}^{\varrho}\leqslant - c_0\,{\kappa^{\varrho}} 2^{-m}.
$$
Then under the assumption
\begin{align}\label{assum-co1}
a_2>\tau_1\quad\hbox{and}\quad 2(1-a)-{a \varrho>0}
\end{align}
one finds successively that $\displaystyle \sup_{m\in\mathbb{N}}2^{m}N_{m-1}^{-a_2+\tau_1}<\infty$ and when $\varepsilon$ is small enough and  $|l|\leqslant N_{m-1}$  we deduce that
\begin{align*}
\big|{\omega}(\alpha,\varepsilon)\cdot l+jc_{m-1}(\alpha,\varepsilon)\big| & \leqslant \tfrac{4\kappa_{m}^{\varrho}|j|}{\langle l\rangle^{\tau_{1}}}+\frac{|j|{\kappa^{\varrho}}}{2^m\langle l\rangle^{\tau_{1}}}\Big(-4c_0+C\varepsilon^{2(1-a)-{a \varrho}}2^{m}N_{m-1}^{-a_2+\tau_1}\Big)\\ & <  \tfrac{4\kappa_{m}^{\varrho}|j|^{\varrho}}{\langle l\rangle^{\tau_{1}}}.
\end{align*}

Hence we obtain  $\alpha\in \mathcal{R}_{l,j}^{(0)}(i_{m-1})$ and this achieves the proof of the first point.
Observe that the first assumption in \eqref{assum-co1} is automatically satisfied from \eqref{Assump-DRX} since $q\geqslant1$.

\smallskip

\textbf{(ii)} Let $j,j_{0}\in\mathbb{S}_{0}^{c}$ and  $(l,j)\neq(0,j_0).$ 
Imitating the proof of  the first point \textbf{(i)}, then to get the desired  result it is enough to  check that $\mathcal{R}_{l,j,j_{0}}(i_{m})\subset \mathcal{R}_{l,j,j_{0}}(i_{m-1}).$
Let  $\alpha\in\mathcal{R}_{l,j,j_{0}}(i_{m})$ then from the definition of this set introduced  in \eqref{set-U0} we deduce that  $\alpha\in\mathtt{C}_{m}^{\kappa}\subset \mathtt{C}_{m-1}^{\kappa}$ and 
\begin{equation}\label{poiH1}\begin{array}{rcl}
&&\big|{\omega}(\alpha,\varepsilon)\cdot l+\mu_{j}^{\infty,m-1}(\alpha,\varepsilon)-\mu_{j_{0}}^{\infty,m-1}(\alpha,\varepsilon)\big|  \leqslant\tfrac{2\kappa_{m+1}\langle j-j_{0}\rangle}{\langle l\rangle^{\tau_{2}}}\\
& & \qquad +\big|\mu_{j}^{\infty,m}(\alpha,\varepsilon)-\mu_{j_0}^{\infty,m}(\alpha,\varepsilon)-\mu_{j}^{\infty,m-1}(\alpha,\varepsilon)+\mu_{j_{0}}^{\infty,m-1}(\alpha,\varepsilon)\big|.
\end{array}
\end{equation}
Set 
$$ \varrho^m_{j,j_0}(\alpha,\varepsilon)\triangleq\big|\mu_{j}^{\infty,m}(\alpha,\varepsilon)-\mu_{j_0}^{\infty,m}(\alpha,\varepsilon)-\mu_{j}^{\infty,m-1}(\alpha,\varepsilon)+\mu_{j_{0}}^{\infty,m-1}(\alpha,\varepsilon)\big|.
$$ 
Then coming back to \eqref{asy-z1}  we deduce that 
\begin{align}\label{asy-z2}
\nonumber \varrho^m_{j,j_0}(\alpha,\varepsilon)&\leqslant |j-j_0|\big|r^{1,m}(\alpha,\varepsilon)-r^{1,m-1}(\alpha,\varepsilon)\big|+\big|j\mathtt{W}(j,\alpha)-j_0\mathtt{W}(j_0,\alpha)\big|\big|r^{2,m}(\alpha,\varepsilon)-r^{2,m-1}(\alpha,\varepsilon)\big|\\
&\quad+\big|r_{j}^{\infty,m}(\alpha,\varepsilon)-r_{j}^{\infty,m-1}(\alpha,\varepsilon)\big|
+\big|r_{j_0}^{\infty,m}(\alpha,\varepsilon)-r_{j_0}^{\infty,m-1}(\alpha,\varepsilon)\big|.
\end{align}
According to Proposition \ref{projection in the normal directions}-(i), \eqref{Bio-X1} and \eqref{choice of gamma and N0 in the Nash-Moser} one gets for $k=1,2$
\begin{align*}
\big|r^{k,m}(\alpha,\varepsilon)-r^{k,m-1}(\alpha,\varepsilon)\big|&\lesssim\, \varepsilon\kappa^{-1}\| i_{m}-i_{m-1}\|_{\overline{s}_{h}+\sigma_3}^{q,\kappa}\\
&\lesssim \varepsilon^{2}\kappa^{-2}N_{m-1}^{-a_{2}}\\
&\lesssim \varepsilon^{2-2a}\langle j-j_0\rangle N_{m-1}^{-a_{2}}.
\end{align*}
Next,
using  Proposition \ref{reduction of the remainder term}-(ii), \eqref{Bio-X1} and \eqref{choice of gamma and N0 in the Nash-Moser}  and using the fact (which follows from the last assumtion in \eqref{Assump-DRX} since $\overline\sigma$ can be taken large enough)
$${s_h}\geqslant 2\overline{s}_h+8\tau_2(q+1)
$$
yield
\begin{align*}
\big|r_{j}^{\infty,m}(\alpha,\varepsilon)-r_{j}^{\infty,m-1}(\alpha,\varepsilon)\big| &\lesssim \varepsilon\kappa^{-2}\left(\| i_{m}-i_{m-1}\|_{\overline{s}_{h}+\sigma_4}^{q,\kappa}\right)^{\frac12}\\
&\lesssim \varepsilon^{\frac32}\kappa^{-\frac52}N_{m-1}^{-\frac12a_{2}}\\
&\lesssim \varepsilon^{\frac{3-5a}{2}}\langle j-j_0\rangle N_{m-1}^{-\frac12a_{2}}.
\end{align*}
Plugging the two preceding  estimates into \eqref{asy-z2} and using the inequality \eqref{ineq:fj} 
 we find
\begin{align}\label{asy-z3}
 \varrho^m_{j,j_0}(\alpha,\varepsilon)\lesssim  \varepsilon^{\frac{3-5a}{2}}\langle j-j_0\rangle N_{m-1}^{-\frac12a_{2}}.
\end{align}
Combining \eqref{asy-z3} with \eqref{poiH1} and using  $\kappa_{m+1}=\kappa_{m}-\varepsilon^a 2^{-m-1}$ allow to get
\begin{align*}
\left|{\omega}(\alpha,\varepsilon)\cdot l+\mu_{j}^{\infty,m-1}(\alpha,\varepsilon)-\mu_{j_{0}}^{\infty,m-1}(\alpha,\varepsilon)\right| 
&\leqslant   \displaystyle\tfrac{2\kappa_{m}\langle j-j_{0}\rangle }{\langle l\rangle^{\tau_{2}}}-{\varepsilon^a \langle j-j_{0}\rangle}2^{-m}\langle l\rangle ^{-\tau_2}\\
&\quad+C \varepsilon^{\frac{3-5a}{2}}\langle j-j_0\rangle N_{m-1}^{-\frac12a_{2}}.
\end{align*} 
Now we write, since $|l|\leqslant N_{m-1}$
\begin{align*}
-{\varepsilon^a }2^{-m}\langle l\rangle ^{-\tau_2}+C\varepsilon^{\frac{3-5a}{2}} N_{m-1}^{-\frac12a_{2}}\leqslant {\varepsilon^a }2^{-m}\langle l\rangle ^{-\tau_2}\Big(-1+C\varepsilon^{\frac{3-7a}{2}}2^{m}N_{m-1}^{\tau_2-\frac12a_{2}}\Big).
\end{align*}
From \eqref{Assump-DRX} and \eqref{choice of gamma and N0 in the Nash-Moser} it is obvious that  
\begin{align}\label{MNKL}
2\tau_2<a_2\quad\hbox{and}\quad a<\tfrac37
\end{align}
implying in turn for  $\varepsilon$ small enough
\begin{align*}
\forall\, m\in\mathbb{N},\quad -1+C\varepsilon^{\frac{3-7a}{2}}2^{m}N_{m-1}^{\tau_2-\frac12a_{2}}\leqslant 0.
\end{align*}
Therefore we obtain 
$$
\big|{\omega}(\alpha,\varepsilon)\cdot l+\mu_{j}^{\infty,m-1}(\alpha,\varepsilon)-\mu_{j_{0}}^{\infty,m-1}(\alpha,\varepsilon)\big| 
\leqslant   \tfrac{2\kappa_{m}\langle j-j_{0}\rangle}{\langle l\rangle^{\tau_{2}}}\cdot
$$ 
This shows that $\alpha \in  \mathcal{R}_{l,j,j_{0}}(i_{m-1})$, which  completes the proof of the second point.

\smallskip

\textbf{(iii)} Let $j\in\mathbb{S}_{0}^{c}, (l,j)\neq (0,0),$ we shall first prove that if $|l|\leqslant N_{m-1}$   then $\mathcal{R}_{l,j}^{(1)}(i_{m}) \subset\mathcal{R}_{l,j}^{(1)}(i_{m-1}).$ As in the point \textbf{(i)} this implies that  $\mathcal{R}_{l,j}^{(1)}(i_{m})=\varnothing.$ Remind that the set $\mathcal{R}_{l,j}^{(1)}(i_{m})$ is defined below the form  \eqref{set-U0}.  Let $\alpha\in\mathcal{R}_{l,j}^{(1)}(i_{m})$ then by construction $\alpha\in \mathtt{C}_{m}^{\kappa}\subset \mathtt{C}_{m-1}^{\kappa}.$\\
 Now by the triangle inequality we may write  in view of Proposition \ref{reduction of the remainder term}-(ii), \eqref{Bio-X1} and  $\kappa=\varepsilon^a$
\begin{align*}
\big|{\omega}(\alpha,\varepsilon)\cdot l+\mu_{j}^{\infty,m-1}(\alpha,\varepsilon)\big| & \leqslant  \big|{\omega}(\alpha,\varepsilon)\cdot l+\mu_{j}^{\infty,m}(\alpha,\varepsilon)\big|+|\mu_{j}^{\infty,m}(\alpha,\varepsilon)-\mu_{j}^{\infty,m-1}(\alpha,\varepsilon)|\\
& {\leqslant} \tfrac{\kappa_{m+1}\langle j\rangle}{\langle l\rangle^{\tau_{1}}}+C\varepsilon\kappa^{-2}|j|\left(\| i_{m}-i_{m-1}\|_{\overline{s}_{h}+\sigma_4}^{q,\kappa}\right)^{\frac12}\\
& \leqslant  \tfrac{\kappa_{m+1}\langle j\rangle}{\langle l\rangle^{\tau_{1}}}+C\varepsilon^{\frac{3-5a}{2}}\langle j\rangle N_{m-1}^{-\frac12a_{2}}.
\end{align*}
Since   $\kappa_{m+1}=\kappa_{m}-\varepsilon^a 2^{-m-1}$ and $|l|\leqslant N_{m-1}$  then 
$$
\big|{\omega}(\alpha,\varepsilon)\cdot l+\mu_{j}^{\infty,m-1}(\alpha,\varepsilon)\big|  \leqslant  \tfrac{\kappa_{m}\langle j\rangle}{\langle l\rangle^{\tau_{1}}}+\tfrac{\langle j\rangle\varepsilon^a}{2^{m+1}\langle l\rangle^{\tau_1}}\Big(-1+\varepsilon^{\frac{3-7a}{2}} 2^{m+1}N_{m-1}^{\tau_1-\frac12a_{2}}\Big).
$$
Putting together \eqref{Conv-T2} and \eqref{MNKL} yields
\begin{align*}
a<\tfrac37 \quad\hbox{and}\quad 2\tau_1<a_2
\end{align*}
and by taking  $\varepsilon$ small enough we find that 
$$
\forall\,m\in \mathbb{N},\quad -1+\varepsilon^{2-4a} 2^{m+1}N_{m-1}^{\tau_1-a_{2}}\leqslant 0.
$$
This  implies in turn that 
$$
\big|{\omega}(\alpha,\varepsilon)\cdot l+\mu_{j}^{\infty,m-1}(\alpha,\varepsilon)\big|  \leqslant \tfrac{\kappa_{m}\langle j\rangle}{\langle l\rangle^{\tau_{1}}}\cdot
$$
It follows   that $\lambda\in \mathcal{R}_{l,j}^{(1)}(i_{m-1})$, which completes the proof. 

\smallskip

\textbf{(iv)} The proof follows easily from \eqref{set-U0} and the points Lemma \ref{lemm-dix1}-{(i)}-{(ii)}-{(iii)}.
\end{proof}
The next result deals with necessary  conditions such that the sets in \eqref{set-U0} are not empty.
\begin{lemma}\label{some cantor set are empty}
Let $a\in(0,\frac12)$ and $\tau_2>\tau_1$, there exists $\varepsilon_0$ such that for any $\varepsilon\in[0,\varepsilon_0]$ and $m\in\mathbb{N},$ the following assertions hold true. 
\begin{enumerate}
\item Let $(l,j)\in\mathbb{Z}^{d}\times\mathbb{Z}\backslash\{(0,0)\}.$ If $\,\displaystyle\mathcal{R}_{l,j}^{(0)}(i_{m})\neq\varnothing,$ then $|j|\leqslant C_{0}\langle l\rangle.$
\item Let $(l,j,j_{0})\in\mathbb{Z}^{d}\times(\mathbb{S}_{0}^{c})^{2}.$ If $\,\displaystyle\mathcal{R}_{l,j,j_{0}}(i_{m})\neq\varnothing,$ then $|j-j_{0}|\leqslant C_{0}\langle l\rangle.$
\item  Let $(l,j)\in\mathbb{Z}^{d}\times\mathbb{S}_{0}^{c}.$ If $\,\displaystyle \mathcal{R}_{l,j}^{(1)}(i_{m})\neq\varnothing,$ then $|j|\leqslant C_{0}\langle l\rangle.$
\item Let $(l,j,j_{0})\in\mathbb{Z}^{d}\times(\mathbb{S}_{0}^{c})^{2}.$ If $\displaystyle \left(\min(|j|,|j_{0}|)\right)^{1-2\overline\alpha-\epsilon}\geqslant c_{2}\kappa^{\varrho}\langle l\rangle^{\tau_{1}},$ then 
$$\mathcal{R}_{l,j,j_{0}}(i_{m})\subset\mathcal{R}_{l,j-j_{0}}^{(0)}(i_{m}),$$
\end{enumerate}
with $c_2$ being a fixed constant and $\epsilon>0$ is an arbitrary small  number.
\end{lemma}
\begin{proof}
\textbf{(i)} Assume $\mathcal{R}_{l,j}^{(0)}(i_{m})\neq\varnothing$ then by construction  we can find   $\alpha\in(\underline\alpha,\overline\alpha)$ such that,
\begin{align}\label{Rad-M1}
\nonumber |c_{m}(\alpha,\varepsilon)||j|&\leqslant 4\langle j\rangle \kappa_{m+1}^{\varrho}\langle l\rangle^{-\tau_{1}}+|{\omega}(\alpha,\varepsilon)\cdot l|\\
\nonumber &\leqslant 4\langle j\rangle \kappa_{m+1}^{\varrho}+C\langle l\rangle\\
&\leqslant C\varepsilon^{a \varrho}\langle j\rangle+C\langle l\rangle,
\end{align}
where we have used $\kappa=\varepsilon^a$ and  the fact that $(\alpha,\varepsilon)\mapsto \omega(\alpha, \varepsilon)$ is bounded. \\
Applying \eqref{est-r1}  combined with Proposition \ref{Nash-Moser}-$(\mathcal{P}1)_{m}$ we write ( by taking $\overline\sigma$ large enough)
\begin{align}\label{c_m-dev}
\nonumber |c_{m}(\alpha,\varepsilon)-V_{\alpha,0}|&\leqslant \|c_{m}-V_{0,\alpha}\|^{q,\kappa}
\\
\nonumber&\lesssim \varepsilon\kappa^{-1}\left(1+\|\mathfrak{I}_{m}\|_{\overline{s}_{h}+\sigma_1}^{q,\kappa}\right)\\
&\lesssim \varepsilon^{1-a}.
\end{align}
On the other hand, one deduces from \eqref{defTheta} that
$$\inf_{\alpha\in(\underline\alpha,\overline\alpha)}V_0(\alpha)\triangleq c_{1}>0,
$$ 
and then taking $\varepsilon$ small enough, we  obtain 
$$\inf_{m\in\mathbb{N}}\inf_{\alpha\in(\underline\alpha,\overline\alpha)}|c_{m}(\alpha,\varepsilon)|\geqslant \tfrac{c_{1}}{2}\cdot
$$
Thus, coming back to \eqref{Rad-M1} and   taking $\varepsilon$ small enough we deduce that  $|j|\leqslant C_{0}\langle l\rangle$ for some constant~$C_{0}>0.$

\smallskip

\textbf{(ii)} First, observe that the  case $j=j_0$ is trivial and the conclusion is true, then we shall assume that  $j\neq j_0$. If  $\mathcal{R}_{l,j,j_{0}}(i_{m})\neq\varnothing$ then there exists $\alpha\in(\underline\alpha,\overline\alpha)$ such that 
\begin{align*}
|\mu_{j}^{\infty,m}(\alpha,\varepsilon)-\mu_{j_{0}}^{\infty,m}(\alpha,\varepsilon)|&\leqslant 2\kappa_{m+1}\langle j-j_{0}\rangle \langle l\rangle^{-\tau_{2}}+|{\omega}(\alpha,\varepsilon)\cdot l|\\
&\leqslant 2\kappa_{m+1}\langle j-j_{0}\rangle+C\langle l\rangle\\
&\leqslant 4\varepsilon^a\langle j-j_{0}\rangle+C\langle l\rangle.
\end{align*}
By virtue of the triangle inequality, Lemma \ref{lem-asym}-{(iv)}, \eqref{asy-z1}, \eqref{ineq:fj},  \eqref{uniform estimate r1} and \eqref{uniform estimate rjinfty} we find 
\begin{align*}
|\mu_{j}^{\infty,m}(\alpha,\varepsilon)-\mu_{j_{0}}^{\infty,m}(\alpha,\varepsilon)| & \geqslant  |\Omega_{j}(\alpha)-\Omega_{j_{0}}(\alpha)|-\big(|r^{1,m}(\alpha,\varepsilon)|+|r^{2,m}(\alpha,\varepsilon)|\big)|j-j_{0}|\\
&\quad -|r_{j}^{\infty,m}(\alpha,\varepsilon)|-|r_{j_{0}}^{\infty,m}(\alpha,\varepsilon)|\\
& \geqslant  \big(C_0-C\varepsilon^{1-2a}\big)|j-j_{0}|\geqslant  \tfrac{C_{0}}{2}|j-j_{0}|
\end{align*}
provided that $\varepsilon$ is small enough. Then combining the preceding inequalities implies  for $\varepsilon $ small enough that  $|j-j_{0}|\leqslant C_{0}\langle l\rangle,$ for some $C_{0}>0.$

\smallskip

\textbf{(iii)} Notice that for  $j=0$ the conclusion  is always true. Now assume for $j\neq 0$ that  $\mathcal{R}_{l,j}^{(1)}(i_{m})\neq\varnothing$ then there exists   $\alpha\in(\underline\alpha,\overline\alpha)$ such that 
\begin{align*}
|\mu_{j}^{\infty,m}(\alpha,\varepsilon)|&\leqslant \kappa_{m+1}\langle j\rangle \langle l\rangle^{-\tau_{1}}+|{\omega}(\alpha,\varepsilon)\cdot l|\\
&\leqslant  2\varepsilon^a|j|+C\langle l\rangle.
\end{align*}
Using the definition \eqref{asy-z1} combined with 
 the triangle  inequality, Lemma \ref{lem-asym}-{(iii)}, \eqref{uniform estimate r1} and \eqref{uniform estimate rjinfty}, we obtain
\begin{align*}
|\mu_{j}^{\infty,m}(\alpha,\varepsilon)|&\geqslant C_0|j|-|j|\big(|r^{1,m}(\alpha,\varepsilon)|+|r^{2,m}(\alpha,\varepsilon)|\big)-|r_{j}^{\infty,m}(\alpha,\varepsilon)|\\
&
\geqslant C_0|j|-C\varepsilon^{1-2a}|j|.
\end{align*}
Putting together these inequalities yields, provided that $a\in(0,\frac12)$,
\begin{align*}
\big( C_{0}-C\varepsilon^{1-2a}-2\varepsilon^a\big)|j|
&\leqslant  C\langle l\rangle.
\end{align*}
Therefore by  taking $\varepsilon$ small enough, we find $|j|\leqslant C_0\langle  l\rangle,$ for some $C_{0}>0.$

\smallskip

\textbf{(iv)} The case $j=j_0$ is trivial and follows from the definition \eqref{set-U0} and the fact $\tau_2>\tau_1$. Consider  $j\neq j_0$ and $\alpha\in\mathcal{R}_{l,j,j_{0}}(i_{m})$ then by definition
$$\big|{\omega}(\alpha,\varepsilon)\cdot l+\mu_{j}^{\infty,m}(\alpha,\varepsilon)-\mu_{j_{0}}^{\infty,m}(\alpha,\varepsilon)\big|\leqslant\tfrac{2\kappa_{m+1}\langle j-j_{0}\rangle}{\langle l\rangle^{\tau_{2}}}\cdot
$$
Combining \eqref{asy-z1} with Lemma \ref{lem-asym}-(iii) gives the asymptotic 
\begin{align*}
 \mu_{j}^{\infty,m}(\alpha,\varepsilon)&=j\,c_m(\alpha,\varepsilon)-j|j|^{\alpha-1}\,\big(W_{0,\alpha}+r^{2,m}(\alpha,\varepsilon)\big)+r_{j}^{\infty,m}(\alpha,\varepsilon)+O\left(\tfrac{1}{|j|^{1-\alpha}}\right).
\end{align*}
Therefore using   this expansion   with the triangle inequality allow to get
\begin{align*}
\big|{\omega}(\alpha,\varepsilon)\cdot l+(j-j_{0})c_{m}(\alpha,\varepsilon)\big|  &\leqslant  \big|{\omega}(\alpha,\varepsilon)\cdot l+\mu_{j}^{\infty,m}(\alpha,\varepsilon)-\mu_{j_{0}}^{\infty,m}(\alpha,\varepsilon)\big|\\
&\quad+\big(W_{0,\alpha}+|r^{2,m}(\alpha,\varepsilon)|\big)\big|j|j|^{\alpha-1}-j_0|j_0|^{\alpha-1}\big|\\
&\quad+\big|r_{j}^{\infty,m}(\lambda,\varepsilon)|+|r_{j_{0}}^{\infty,m}(\lambda,\varepsilon)\big|\\
&\quad +O\left(\tfrac{1}{|j|^{1-\alpha}}\right)+O\left(\tfrac{1}{|j_0|^{1-\alpha}}\right).
\end{align*}
From  \eqref{uniform estimate r1} and \eqref{ineq:fj} we obtain
$$
\big(W_{0,\alpha}+|r^{2,m}(\alpha,\varepsilon)|\big)\big|j|j|^{\alpha-1}-j_0|j_0|^{\alpha-1}\big|\leqslant C |j-j_0|\left(\frac{1}{|j|^{1-\alpha}}+\frac{1}{|j_0|^{1-\alpha}}\right).
$$
Thus combining the last two estimates and \eqref{uniform estimate rjinfty} we find for $j\neq j_0\in\mathbb{S}_0^c$
\begin{align}\label{ZaraX1}
\nonumber \big|{\omega}(\alpha,\varepsilon)\cdot l+(j-j_{0})c_{m}(\alpha,\varepsilon)\big| & \leqslant  \frac{2\kappa_{m+1}\langle j-j_{0}\rangle}{\langle l\rangle^{\tau_{2}}}+C |j-j_0|\left(\frac{1}{|j|^{1-\alpha}}+\frac{1}{|j_0|^{1-\alpha}}\right)
\\
\nonumber&\quad+C\varepsilon^{1-2a}\left(\frac{1}{|j|^{1-\epsilon-2\overline\alpha}}+\frac{1}{|j_0|^{1-\epsilon-2\overline\alpha}}\right)\\
 & \leqslant\frac{2\kappa_{m+1}\langle j-j_{0}\rangle}{\langle l\rangle^{\tau_{2}}}+C \frac{|j-j_0|}{\left(\min(|j|,|j_{0}|)\right)^{1-2\overline\alpha-\epsilon}}\cdot
\end{align}
Since $\varrho, \kappa_m\in(0,1)$ then $\kappa_{m+1}\leqslant \kappa_{m+1}^\varrho$. Therefore if we assume (we know that $\tau_{2}>\tau_{1}$)
$$ \left(\min(|j|,|j_{0}|)\right)^{1-2\overline\alpha-\epsilon}\geqslant \tfrac12 C\kappa_{m+1}^{-\varrho}\langle l\rangle^{\tau_{1}},
$$ 
 we deduce from  \eqref{ZaraX1} that
$$\begin{array}{rcl}
\big|{\omega}(\alpha,\varepsilon)\cdot l+(j-j_{0})c_{m}(\alpha,\varepsilon)\big| & \leqslant & \frac{4\kappa_{m+1}^{\varrho}|j-j_{0}|}{\langle l\rangle^{\tau_{1}}}\cdot
\end{array}$$
To get the result it suffices to use  $\kappa_{m+1}\in[\kappa,2\kappa].$
This ends the proof of Lemma \ref{some cantor set are empty}.
\end{proof}

The last result to prove deals with the uniform transversality during Nash-Moser scheme.
\begin{lemma}\label{lemma R\"ussmann condition for the perturbed frequencies}
{Let $q=q_0$ and $\rho_{0}$ as in Proposition  $\ref{lemma transversality},$ and consider the function  $(\alpha,\varepsilon)\mapsto \omega(\alpha,\varepsilon)$ stated in  \eqref{estimate repsilon1}. There exist $\varepsilon_{0}>0$ small enough such that for any   $\varepsilon\in[0,\varepsilon_{0}]$  we have the following assertions.
\begin{enumerate}
\item For all $m\in\N,\,(l,j)\in\mathbb{Z}^{d+1}\backslash\{0\}$ such that $|j|\leqslant C_{0}\langle l\rangle,$ we have
$$\inf_{\alpha\in(\underline\alpha,\overline\alpha)}\max_{0\leqslant k\leqslant q}|\partial_{\alpha}^{k}\big(\omega(\alpha,\varepsilon)\cdot l+jc_{m}(\alpha,\varepsilon)\big)|\geqslant\frac{\rho_{0}\langle l\rangle}{2}\cdot
$$
\item For all $m\in\N,\,(l,j)\in\mathbb{Z}^{d}\times\mathbb{S}_{0}^{c}$ such that $|j|\leqslant C_{0}\langle l\rangle,$ we have
$$\inf_{\alpha\in(\underline\alpha,\overline\alpha)}\max_{0\leqslant k\leqslant q}\big|\partial_{\alpha}^{k}\big({\omega}(\alpha,\varepsilon)\cdot l+\mu_{j}^{\infty,m}(\alpha,\varepsilon)\big)\big|\geqslant\frac{\rho_{0}\langle l\rangle}{2}\cdot
$$
\item For all $(l,j,j_{0})\in\mathbb{Z}^{d}\times(\mathbb{S}_{0}^{c})^{2}$ such that $|j-j_{0}|\leqslant C_{0}\langle l\rangle,$ we have
$$\inf_{\alpha\in(\underline\alpha,\overline\alpha)}\max_{0\leqslant k\leqslant q}\big|\partial_{\alpha}^{k}\big({\omega}(\alpha,\varepsilon)\cdot l+\mu_{j}^{\infty,m}(\alpha,\varepsilon)-\mu_{j_{0}}^{\infty,m}(\alpha,\varepsilon)\big)\big|\geqslant\frac{\rho_{0}\langle l\rangle}{2}\cdot$$
\end{enumerate}}
\end{lemma}
\begin{proof}
\textbf{(i)}
Using the triangle  inequality combined with  \eqref{estimate repsilon1},  \eqref{c_m-dev},  \eqref{choice of gamma and N0 in the Nash-Moser}, Proposition \ref{lemma transversality}-\textbf{(ii)} and the fact that $|j|\leqslant C_{0}\langle l\rangle$ we get
\begin{align*}
\displaystyle\max_{0\leqslant k\leqslant q}|\partial_{\alpha}^{k}\left({\omega}(\alpha,\varepsilon)\cdot l+jc_{m}(\alpha,\varepsilon)\right)| & \geqslant  \displaystyle\max_{0\leqslant k\leqslant q}|\partial_{\alpha}^{k}\left({\omega}_{\textnormal{Eq}}(\alpha)\cdot l+jV_0(\alpha)\right)|\\
&\quad-\max_{k\leqslant q}|\partial_{\alpha}^{k}\overline{\mathrm{r}}_{\varepsilon}(\alpha)\cdot l|-C|j|\varepsilon ^{1-a(1+q)}\\
&  \geqslant  \displaystyle\rho_{0}\langle l\rangle-C\varepsilon^{1-a(1+q+q\overline{a})}\langle l\rangle-C\langle l\rangle\varepsilon ^{1-a(1+q)}.\end{align*}
It follows that
\begin{align*}
\displaystyle\max_{0\leqslant k\leqslant q}|\partial_{\alpha}^{k}\left({\omega}(\alpha,\varepsilon)\cdot l+jc_{m}(\alpha,\varepsilon)\right)| 
 \geqslant & \displaystyle\left(\rho_{0}-C\varepsilon^{1-a(1+q+q\overline{a})}\right)\langle l\rangle\\
 \geqslant & \displaystyle\frac{\rho_{0}\langle l\rangle}{2}
\end{align*}
provided that
\begin{align}\label{est-cond-a}
{a<\frac{1}{1+q(1+\overline{a})}}
\end{align}
and $\varepsilon$ is taken small enough. We observe that the condition \eqref{est-cond-a} follows from \eqref{choice of gamma and N0 in the Nash-Moser} and \eqref{Assump-DRX}.

\smallskip

\textbf{(ii)} This proof is similar to  that of the first point. Indeed,  we combine the  triangle inequality, \eqref{estimate repsilon1}, \eqref{asy-z1},\eqref{uniform estimate r1}, \eqref{uniform estimate rjinfty},  \eqref{choice of gamma and N0 in the Nash-Moser} and Proposition \ref{lemma transversality} with  $|j|\leqslant C_{0}\langle l\rangle$ in order to get
\begin{align*}
\max_{
0\leqslant k\leqslant q}\big|\partial_{\alpha}^{k}\big({\omega}(\alpha,\varepsilon)\cdot l+\mu_{j}^{\infty,m}(\alpha,\varepsilon)\big)\big| & \geqslant  \displaystyle\max_{0\leqslant k\leqslant q}|\partial_{\alpha}^{k}\left({\omega}_{\textnormal{Eq}}(\alpha)\cdot l+\Omega_{j}(\alpha)\right)|\\
&\quad \displaystyle-\max_{0\leqslant k\leqslant q}\big|\partial_{\alpha}^{k}\big(\overline{\mathrm{r}}_{\varepsilon}(\alpha)\cdot l+jr^{1,m}(\lambda,\varepsilon)+j|j|^{\alpha-1}\,r^{2,m}(\alpha,\varepsilon)\big)\big|\\
&\quad\displaystyle-\max_{0\leqslant k\leqslant q}\big|\partial_{\alpha}^{k}r_{j}^{\infty,m}(\alpha,\varepsilon)\big)\big|\\
& \geqslant  \displaystyle\rho_{0}\langle l\rangle-C\varepsilon^{1-a(1+q+q\overline{a})}\langle l\rangle-C\varepsilon^{1-a(2+q)}|j|.
\end{align*}
Therefore 
\begin{align*}
\max_{
0\leqslant k\leqslant q}\big|\partial_{\alpha}^{k}\big({\omega}(\alpha,\varepsilon)\cdot l+\mu_{j}^{\infty,m}(\alpha,\varepsilon)\big)\big| & \geqslant   \displaystyle\frac{\rho_{0}\langle l\rangle}{2}
\end{align*}
provided that  $\varepsilon$ is small enough and the condition \eqref{est-cond-a} is satisfied.

\smallskip

\textbf{(iii)} The case $j=j_0$ is trivial. Arguing as in the preceding point, using  \eqref{asy-z1} combined with  Proposition \ref{lemma transversality}, \eqref{ouz-end},   \eqref{uniform estimate r1}, \eqref{uniform estimate rjinfty}, \eqref{choice of gamma and N0 in the Nash-Moser}   and the fact that $0<|j-j_{0}|\leqslant C_{0}\langle l\rangle$ we have 
\begin{align*}
\max_{0\leqslant k\leqslant q}\big|\partial_{\alpha}^{k}\big({\omega}(\lambda,\varepsilon)\cdot l&+\mu_{j}^{\infty,m}(\alpha,\varepsilon)-\mu_{j_{0}}^{\infty,m}(\lambda,\varepsilon)\big)\big|
\geqslant\displaystyle\max_{0\leqslant k\leqslant q}\big|\partial_{\alpha}^{k}\big({\omega}_{\textnormal{Eq}}(\alpha)\cdot l+\Omega_{j}(\alpha)-\Omega_{j_{0}}(\alpha)\big)\big|\displaystyle\\
&-\max_{0\leqslant k\leqslant q}\big|\partial_{\alpha}^{k}\big(\overline{\mathrm{r}}_{\varepsilon}(\alpha)\cdot l+(j-j_{0})r^{1,m}(\alpha,\varepsilon)+r_{j}^{\infty,m}(\alpha,\varepsilon)-r_{j_{0}}^{\infty,m}(\alpha,\varepsilon)\big)\big|\\
&-\max_{0\leqslant k\leqslant q}\big|\partial_{\alpha}^{k}\big[(j\mathtt{W}(j,\alpha)-j_0\mathtt{W}(j_0,\alpha))r^{2,m}(\alpha,\varepsilon)\big]\big|\\
&\geqslant\displaystyle\rho_{0}\langle l\rangle-C\varepsilon^{1-a(1+q+q\overline{a})}\langle l\rangle-C\varepsilon^{1-a(2+q)}|j-j_{0}|-C\varepsilon^{1-a(1+q+q\overline{a})}|j-j_0|\\
&\geqslant \displaystyle\frac{\rho_{0}\langle l\rangle}{2}
\end{align*}
for $\varepsilon$ small enough.  Notice that we have used the following inequality \eqref{ineq:fj}.
This ends the proof of Lemma \ref{lemma R\"ussmann condition for the perturbed frequencies}.
\end{proof}

\appendix

 \section{Appendix}
We intend in this section to recall and establish some results related to Gamma function and  used before in  some proofs. 
The  function  $\Gamma:\CC\backslash(-\NN)\to \CC$  refers to the gamma function which is the analytic continuation to the negative half plane of the usual gamma function defined on the positive half-plane $\big\{\hbox{Re}z>0\big\}$ by the integral representation
$$
\Gamma(z)=\int_{0}^{\infty} t^{z-1}\, e^{-t}dt.
$$
It satisfies some algebraic identities such as  the relation
\begin{equation}\label{Gamma1c1}
\Gamma(z+1)=z\,\Gamma(z), \quad \forall z\in \CC \backslash(-\NN).
\end{equation}
or the Legendre duplication formula,
\begin{equation}\label{leg-dup}
\Gamma(z)\Gamma(z+\frac12)=2^{1-2z}\sqrt{\pi}\,\Gamma(2z).
\end{equation}

 It is known  that $\Gamma$ function does not vanish and admits simple  poles at  $\{-n, n\in \NN\}$ and the reciprocal gamma function $\frac{1}{\Gamma}$ is an entire function. Moreover, the real function $x\in(0,\infty)\mapsto\Gamma(x)$ reaches its absolute minimum at a point  $x_0\in(7/5,3/2),$ that is, 
  $$\forall\, x>0,\quad \Gamma(x)\geq\Gamma(x_0)>0.$$ 
For $x\in\RR,$  we denote by $(x)_n$ the Pokhhammer's symbol defined by
\begin{equation}\label{Poch}
(x)_n=\left\{ \begin{array}{ll}
x(x+1)...(x+n-1),\quad \hbox{if}\quad n\geq1, &\\
1,\quad \hbox{if}\quad n=0.
\end{array} \right.
\end{equation} 
From the identity \eqref{Gamma1c1} we deduce the relations
\begin{equation}\label{Pocc1}
(x)_n=\frac{\Gamma(x+n)}{\Gamma(x)}\quad\hbox{and}\quad (x)_n=(-1)^n\frac{\Gamma(1-x)}{\Gamma(1-x-n)},
\end{equation}
provided all the quantities   are well-defined. 
The next result is very important and has been frequently used. It describes some asymptotic behavior of Gamma quotient (also named Wallis quotient) given by, 
$$
 \alpha\in[0,1],\,j\in[0,\infty),\quad \mathtt{W}(j,\alpha)\triangleq \frac{\Gamma\big(j+\frac{\alpha}{2}\big)}{\Gamma\big(j+1-\frac{\alpha}{2}\big)}\cdot
$$
Our result reads as follows. 
\begin{lemma}\label{Stirling-formula}

 For positive large real number  $j$, we have 
$$
\mathtt{W}(j,\alpha)=\frac{1}{j^{1-\alpha}}+O\big(\tfrac{1}{j^{3-\alpha}}\big).
$$
Moreover, the following estimates hold true.
\begin{enumerate}
\item For any $m,k \in\NN,$ there exists a constant $C(k,m)$ such that   \begin{align*}
\forall \alpha\in [0,1], \; \forall j\in[1,\infty),\quad  |\partial_{\alpha}^k\partial_j^m\mathtt{W}(j,\alpha)| &\leqslant  C(k,m)j^{\alpha-1-{m}}\left(1+\ln^{k}(j)\right)\nonumber \\&
\leqslant  C(\epsilon,k,m) j^{\alpha+\epsilon-1-{m}}.
\end{align*}
\item For any $\gamma,k \in\NN,$ there exists a constant $C(k,m)$ such that  \begin{align*}
\forall \alpha\in [0,1],\;\forall\, j\in\NN,\quad |\partial_{\alpha}^k\Delta_j^\gamma\mathtt{W}(j,\alpha)|&\leqslant C(\gamma,q,\epsilon)\,\langle j\rangle^{\alpha+\epsilon-1-{\gamma}}.
\end{align*}

We recall that the difference operator $\Delta_j$ was introduced in \eqref{Difference-op78}.
\end{enumerate}
\end{lemma}
\begin{proof}
Recall Hermite's formula  of the $\Gamma$ function, see formula (1.10) in \cite{Nemes},
\begin{equation}\label{gamma-expan}
\begin{aligned}
\forall\, \textnormal{Re}\{z\}>0,\quad \forall a\in[0,1],\quad \log\Gamma(z+a)&=\Big(z+a-\frac12\Big)\log (z)-z+\frac12\log(2\pi)\\ &+\int_0^{\infty}\bigg(\frac{e^{(1-a)t}}{e^{t}-1}+a-\frac12-\frac1t\bigg)\frac{e^{-zt}}{t}dt.
\end{aligned}
\end{equation}
This  gives the following expression for Wallis quotient, 
\begin{align}\label{Exact-wallis}
\forall j\in [1,\infty),\quad \mathtt{W}(j,\alpha) = \frac{1}{j^{1-\alpha}}\exp\Big( \mathtt{w}_1(j,\alpha)+ \mathtt{w}_2(j,\alpha)\Big)
\end{align}
with
\begin{align*}
\nonumber
 \mathtt{w}_1(j,\alpha)\triangleq \bigintsss_0^{1} \mathtt{K}(\alpha,t)\frac{e^{-jt}}{t}dt,&\quad 
   \mathtt{w}_2(j,\alpha)\triangleq \bigintsss_1^{\infty} \mathtt{K}(\alpha,t)\frac{e^{-jt}}{t}dt
 \\ 
 \textnormal{and}\qquad  \mathtt{K}(\alpha,t)&\triangleq \frac{e^{(1-\alpha)t}-1}{e^{t}-1}e^{\frac\alpha2 t}+\alpha-1 .
\end{align*}
It is straightforward that, for all $k\geq 0$ there exist $C=C(k)$ such that
\begin{align*}
\forall \alpha\in [0,1],\quad \forall t\in [1,\infty),\quad   |\partial_{\alpha}^k\mathtt{K}(\alpha,t)|&  \leqslant C(1+t^k).
\end{align*}
It follows that $\mathtt{w}_2(j,\alpha)$ satisfies the estimate
\begin{align}\label{estimate:w2}
\forall j\in[1,\infty),\quad\forall \alpha\in [0,1],\quad   |\partial_{\alpha}^k\partial_j^m\mathtt{w}_2(j,\alpha)|& \leqslant  \int_1^{\infty} |\partial_{\alpha}^k \mathtt{K}(\alpha,t)|t^{m-1} e^{-jt}dt\nonumber\\
& \leqslant  C \int_1^{\infty} \big(1+t^{k+m-1}\big) e^{-jt}dt\nonumber\\
& \leqslant  C(k,m) {e^{-j}}.
\end{align}
On the other hand, using Taylor formula we get
$$
\forall\, X>0,\quad \frac{X^{1-\alpha}-1}{X-1}=(1-\alpha)\int_0^1\big(1-\tau+\tau X)^{-\alpha} d\tau.
$$
Thus by the substitution $X=e^t$ we infer 
$$
 \mathtt{K}(\alpha,t)=(1-\alpha)\int_0^{1}\Big(\mathtt{F}(\alpha,\tau,t )-1\Big) d\tau \quad \textnormal{with }\quad \mathtt{F}(\alpha,\tau,t )\triangleq  \frac{e^{\frac\alpha2 t}}{\big(1-\tau+\tau e^{t}\big)^\alpha}\cdot
$$
 Applying Taylor expansion of order two in the variable $t$, combined with 
 $$\mathtt{F}(\alpha,\tau,0)=1\quad\hbox{and}\quad \partial_t\mathtt{F}(\alpha,\tau,0)=\alpha\big(\tfrac12-\tau\big)
 $$
  lead  to
\begin{align*}
\int_0^{1}\Big(\mathtt{F}(\alpha,\tau,t )-1\Big) d\tau&=\alpha\, t\int_0^{1}(\tfrac12-\tau) d\tau+
t^2\int_0^1\int_0^1(\partial_t^2 \mathtt{F})(\alpha,\tau,t s) (1-s)d\tau ds\\
&=t^2\,\int_0^1\int_0^1(\partial_t^2 \mathtt{F})(\alpha,\tau,t s) (1-s)d\tau ds.
\end{align*}
Consequently, we obtain
\begin{align*}
 \mathtt{K}(\alpha,t) &= t^2 \mathtt{G}(\alpha, t )\quad \textnormal{with }\quad  \mathtt{G}(\alpha, t )\triangleq (1-\alpha)\int_0^1\int_0^1(\partial_t^2 \mathtt{F})(\alpha,\tau,t s) (1-s)d\tau ds.
\end{align*}
Since $F$ is $\mathcal{C}^\infty$ in $[0,1]^3$, then $G$  should be  $\mathcal{C}^\infty$  in $[0,1]^2$. Thus, for any  $k\geq 0$ there exists $C>0$ such that
$$
\forall \alpha, t\in [0,1],\quad |\partial^k_\alpha\, \mathtt{G}(\alpha, t )|\leqslant C.
$$
It follows that  the expression of $\mathtt{w}_1(j,\alpha) $ becomes
\begin{align}
 \mathtt{w}_1(j,\alpha)&=\int_0^{1}\mathtt{G}(\alpha, t ) t e^{-jt}  dt
\end{align}
and one may easily check that
\begin{align}\label{estimate:w3}
\forall j\in[1,\infty),\quad\forall \alpha\in [0,1],\quad   |\partial_{\alpha}^k\partial_j^m\mathtt{w}_1(j,\alpha)|&\leqslant \int_0^{1}|\partial_{\alpha}^k\mathtt{G}(\alpha, t )|  t^{m+1} e^{-jt}  dt\nonumber \\ &\leqslant   C(k,m)\, j^{-2-{m}}.
\end{align}
Putting together \eqref{Exact-wallis}, \eqref{estimate:w2} and  \eqref{estimate:w3} and using Leibniz rule give the desired asymptotic of $\mathtt{W}(j,\alpha)$ with the suitable estimates as stated in  {\rm (i)}. The same  arguments used to establish  (i) combined with the identity \eqref{It-action} allow to  get  (ii)  . Therefore, the proof of Lemma \ref{Stirling-formula} is achieved.

\end{proof}

{The following elementary result has been used before at several points.}
\begin{lemma}\label{lemma sum Nn}
	Let $N_{0}\geqslant 2.$ Consider the sequence $(N_{m})_{m\in\mathbb{N}}$ defined by \eqref{definition of Nm}. Then for all $\alpha>0$, we have
	$$\sum_{k=m}^{\infty}N_{k}^{-\alpha}\underset{m\rightarrow\infty}{\sim}N_{m}^{-\alpha}.$$

\end{lemma}
\begin{proof}
Let us consider the decreasing  continuous function
$$t\in\mathbb{R}_{+}^{*}\mapsto N_{0}^{-\alpha\left(\frac{3}{2}\right)^{t}}={e^{-\alpha\ln(N_{0})e^{t\ln\frac{3}{2}}}}.$$
Then applying  the series-integral comparison test  and making a change of variables yield, 
\begin{align}\label{Rem-estP}
\nonumber \sum_{k=m+1}^{\infty}N_{k}^{-\alpha}&\leqslant\int_{m}^{\infty}e^{-\alpha\ln(N_{0})e^{t\ln\frac{3}{2}}}{dt}\\
&\leqslant \int_{0}^{\infty} e^{-\alpha\ln(N_{0})e^{u\ln\frac{3}{2}}e^{m\ln\frac{3}{2}}}du.
\end{align}
Now remark that
$${N_{m}^{\alpha}}e^{-\alpha\ln(N_{0})e^{u\ln\frac{3}{2}}e^{m\ln\frac{3}{2}}}=\exp\left(\alpha\ln(N_{0})\left(1-e^{u\ln\left(\frac{3}{2}\right)}\right)e^{m\ln\left(\frac{3}{2}\right)}\right).$$
Then using 
$$\forall u\geqslant 0,\quad 1-e^{u\ln\left(\tfrac{3}{2}\right)}\leqslant -u\ln\left(\tfrac{3}{2}\right)
$$
we find 
$$
\forall u\geqslant 0,\quad{N_{m}^{\alpha}}e^{-\alpha\ln(N_{0})e^{u\ln\frac{3}{2}}e^{m\ln\frac{3}{2}}}\leqslant \exp\left(-\alpha\ln(N_{0})u\ln\left(\tfrac{3}{2}\right)\left(\tfrac{3}{2}\right)^m\right)
.
$$
Inserting this estimate into \eqref{Rem-estP} and using  $N_0\geqslant 2$ we get an absolute constant $C>0$ such 
\begin{align*}\forall u\in\mathbb{R}_{+}^{*},\forall m\in\mathbb{N},\quad{N_{m}^{\alpha}}e^{-\alpha\ln(N_{0})e^{u\ln\frac{3}{2}}e^{m\ln\frac{3}{2}}}&\leqslant {e^{-C\alpha\,(\frac32)^m u}}\in L^{1}(\mathbb{R}_{+}).
\end{align*}
It follows that  
$$\sum_{k=m+1}^{\infty}N_{k}^{-\alpha}\underset{m\rightarrow\infty}{=}o\left(N_{m}^{-\alpha}\right).$$
As a consequence, we find
$$\sum_{k=m}^{\infty}N_{k}^{-\alpha}=N_{m}^{-\alpha}+\sum_{k=m+1}^{\infty}N_{k}^{-\alpha}\underset{m\rightarrow\infty}{\sim}N_{m}^{-\alpha},
$$	
which achieves the proof of Lemma \ref{lemma sum Nn}.
\end{proof}

We shall end this section with recalling the following  lemma due to R\"usseman \cite[Theorem 17.1]{Russ} and used in a crucial way to establish \eqref{kio1}. 
\begin{lemma}\label{Piralt}
Let $ q\in\NN$  $(\alpha,\beta)\in(\mathbb{R}_{+}^{*})^{2}$ and $f:[a,b]\to \RR$ be  a function of class $\mathscr{C}^q$  such that 
$$\inf_{x\in[a,b]}\max_{0\le j\leqslant q}|f^{(j)}(x)|\geqslant\beta.$$
Then  there exists $C>0$ depending only on $b-a, q$ and $ \|g\|_{\mathscr{C}^k}$ such that 
$$\Big|\big\lbrace x\in [a,b];\;\, |f(x)|\leqslant\alpha\big\rbrace\Big|\leqslant C\frac{\alpha^{\frac{1}{q}}}{\beta^{1+\frac{1}{q}}}\cdot $$
\end{lemma}

\section*{Acknowledgements:}
 The work of Z. Hassainia and  N. Masmoudi  is supported   by Tamkeen under the NYU Abu Dhabi Research Institute grant of the center SITE.    The work of N. M is supported by  NSF grant DMS-1716466.

Z. Hassainia and N. Masmoudi  would like to thank M. Berti for many discussions about KAM techniques.

%
%
%
%
%


\begin{thebibliography}{99}

  \bibitem{Alazard-baldi} T. Alazard and  P. Baldi. 
  {\it Gravity capillary standing water waves}. 
Arch. Ration. Mech. Anal. 217 (2015), no. 3, 741--830.
\bibitem{ADDMW} W. Ao, J. Davila, M. del Pino, M. Musso, and J. Wei.
{\it Travelling and rotating solutions to the generalized inviscid surface quasi-geostrophic equation}.
Trans. Amer. Math. Soc. 374 (2021), 6665--6689.
 \bibitem{Arnold} V.  Arnold.
{\it  Small denominators and problems of stability of motion in
classical mechanics and celestial mechanics}. 
Uspekhi Mat. Nauk 18 (1963), 91--192.
\bibitem{Baldi-berti} P. Baldi, M. Berti, E. Haus and R. Montalto. 
{\it Time quasi-periodic gravity water waves in finite depth.} 
Invent. Math., 214 ( 2018) no.2, 739--911.
 \bibitem{Baldi-Montalto21} P. Baldi and  R. Montalto.
  {\it Quasi-periodic incompresible Euler fows in 3D}. 
   Adv. Math. 384 (2021), Paper No. 107730, 74 pp.
 \bibitem{Baldi-Berti-Montalto14} P. Baldi, M. Berti and R.  Montalto.
   {\it KAM for quasi-linear and fully nonlinear forced perturbations of Airy equation}. 
   Math. Ann. 359 (2014), no. 1-2, 471--536.
   \bibitem{BBM-auto} P. Baldi, M. Berti and  R. Montalto.
    {\it KAM for autonomous quasi-linear perturbations of KdV}. 
Ann. Inst. H. Poincar\'e Analyse Non. Lin. 33 (2016), no. 6, 1589--1638. 
\bibitem{bambusi-Berti} D. Bambusi, M. Berti and E.  Magistrelli.
 {\it Degenerate KAM theory for partial differential equations}.
  Journal Diff. Equations, 250 (2011), no. 8, 3379--3397.
  \bibitem{BB13} M. Berti and   P. Bolle,  
{\it A Nash-Moser approach to KAM theory}.
 Fields Institute Communications, special volume ``Hamiltonian PDEs and Applications'', 
255--284,  2015. 
\bibitem{Berti-Bolle2013} M. Berti, P.  Bolle.
{\it Quasi-periodic Solutions of Nonlinear Wave Equations on the d-Dimensional Torus}. 
European Mathematical Society, ISSN 3037192119, 9783037192115, 2020.
 \bibitem{BB10} M. Berti,  P. Bolle and  M. Procesi.
  {\it An abstract Nash-Moser theorem with parameters and applications to PDEs}. 
  Ann. Inst. H. Poincar\'e Anal. Non Lin\'eaire 27 (2010), no. 1, 377--399.
  \bibitem{BCP} M. Berti, L. Corsi and M. Procesi.
   {\it An Abstract Nash Moser Theorem and Quasi-Periodic Solutions for NLW and NLS on Compact Lie Groups and Homogeneous Manifolds.} 
   Comm. Math. Phys. 334, no. 3, 1413--1454, 2015.
\bibitem{BFM21}M. Berti, L. Franzoi and  A. Maspero.
  {\it Traveling quasi-periodic water waves with constant vorticity}.
Archive for Rational Mechanics and Analysis, 240 (2021), 99--202.
\bibitem{BFM} M. Berti, L. Franzoi and  A. Maspero.
 {\it Pure gravity traveling quasi-periodic water waves with constant vorticity}. 
 	arXiv:2101.12006.
\bibitem{BHM} M. Berti, Z. Hassainia and N. Masmoudi.	{\it Time quasi-periodic vortex patches}. Preprint.
\bibitem{BertiMontalto} 
M. Berti and  R. Montalto. 
{\it Quasi-periodic standing wave solutions of gravity-capillary water waves}.
 MEMO, Volume 263, 1273, Memoires AMS, ISSN 0065-9266, 2020.
 \bibitem{B-C} A. L. Bertozzi and P. Constantin.
  {\it  Global regularity for vortex patches}.
   Comm. Math. Phys., 152 (1993), no. 1, 9--28.
\bibitem{Bourgain}  J.  Bourgain, {\it Construction of quasi-periodic solutions for Hamiltonian perturbations of linear equations and applications to nonlinear PDE}. Internat. Math. Res. Notices (1994), no. 11, 475ff., approx. 21 pp
\bibitem{BSV} T. Buckmaster, S. Shkoller and V. Vicol.
 {\it Nonuniqueness of weak solutions to the SQG equation}. 
 Comm. Pure Appl. Math. 72 (2019), no. 9, 1809--1874.
\bibitem{B} J.\ Burbea.
 {\it Motions of vortex patches.} 
Lett. Math. Phys. 6 (1982), no. 1,  1--16.
\bibitem{CQZZ} D. Cao, G. Qin, W. Zhan, C. Zou.
{\it Existence and regularity of co-rotating and travelling global solutions for the generalized SQG equation}.  
 arXiv:2103.03992
 \bibitem{CCG} A. Castro, D. C\'ordoba and  J. G\'omez-Serrano.
  {\it  Existence and regularity of rotating global solutions for the generalized surface quasi-geostrophic equations}. 
  Duke Math. J. 165 (2016), no. 5, 935--984.
 \bibitem{CCG4} A. Castro, D. C\'ordoba and J. G\'omez-Serrano.
  {\it   Uniformly rotating analytic global patch solutions for active scalars}.
   Ann. PDE 2  (2016) no.1,  1--34.
 \bibitem{CCG2} A. Castro, D. C\'ordoba and  J. G\'omez-Serrano.
  {\it  Uniformly rotating smooth solutions for the incompressible 2D Euler equations}.
   Arch. Ration. Mech. Anal. 231 (2019), no. 2, 719--785.
\bibitem{C-C-C-G-W} D. Chae, P. Constantin, D. C\'ordoba, F. Gancedo, and J. Wu.
 {\it Generalized surface quasi-geostrophic equations with singular velocities}.
    Comm. Pure Appl. Math., 65 (2012), no. 8, 1037--1066.
\bibitem{C-C-W} D. Chae, P. Constantin and J. Wu. 
{\it Inviscid models generalizing the two-dimensional Euler and the surface quasi- geostrophic equations}.
 Arch. Ration. Mech. Anal., 202 (2011), no. 1, 35--62.
\bibitem{Ch} J.Y.\ Chemin.
 {\it Fluides parfaits incompressibles}.
  Ast\'{e}risque {230}, Soci\'{e}t\'{e} Math\'{e}matique de France (1995).
\bibitem{C-M-T} P. Constantin, A. J. Majda, and E. Tabak.
 {\it Formation of strong fronts in the 2-D quasigeostrophic thermal active scalar}.
   Nonlinearity, 7 (1994), no. 6, 1495--1533.
\bibitem{Cord} D. C\'ordoba. 
{\it Nonexistence of simple hyperbolic blow-up for the quasi-geostrophic equation}.
 Ann. of Math., 148 (1998), no. 2, 1135--1152.
\bibitem{C-F} D. C\'ordoba and C. Fefferman. 
{\it Growth of solutions for QG and 2D Euler equations}. 
J. Amer. Math. Soc., 15 (2002), no. 3, 665--670.
\bibitem{C-F-M-R} D. C\'ordoba, M. A. Fontelos, A. M. Mancho and J. L. Rodrigo.
 {\it Evidence of singularities for a family of contour dynamics equations.}  
 Proc. Natl. Acad. Sci. USA 102 (2005),  5949--5952.
\bibitem{Cord-Cord-Gan} A. C\'ordoba, D. C\'ordoba and F. Gancedo. 
{\it Uniqueness for SQG patch solutions.} 
Trans. Amer. Math. Soc. Ser. B, 5 (2018), 1--31.
\bibitem{DZ} G.S. Deem and N. J. Zabusky,
 {\it Vortex waves : Stationary "V-states", Interactions, Recurrence, and Breaking.}
Phys. Rev. Lett.  40  (1978), no. 13, 859--862.
\bibitem{DHH} F. de la Hoz, Z. Hassainia and T. Hmidi. 
 {\it Doubly connected V-states for the generalized surface quasi-geostrophic equations}.
  Arch. Ration. Mech. Anal., 220 (2016), no. 3, 1209--1281.
\bibitem{DHHM} F. de la Hoz, Z. Hassainia, T. Hmidi and J. Mateu.
 {\it An analytical and numerical study of steady patches in the disc}. 
 Anal. PDE, 9 (2016), no. 7, 1609--1670.
  \bibitem{HFMV} F. de la Hoz,  T. Hmidi,  J.  Mateu and  J.  Verdera. 
  {\it  Doubly connected V-states for the planar Euler equations}.
    SIAM J. Math. Anal. 48 (2016), no. 3, 1892--1928. 
\bibitem{DHR} D. G. Dritschel, T. Hmidi, C. Renault. {\it  Imperfect bifurcation for the quasi-geostrophic shallow-water equations.} Arch. Ration. Mech. Anal. 231 (2019), no. 3, 1853--1915
\bibitem{Ell-Kuk } L.H. Eliasson and S.B. Kuksin.
{\it  KAM for the nonlinear Schr\"odinger equation}.
 Ann. Math 172 (2010), 371--435.
\bibitem{FGMP19}
R.~Feola, F.~Giuliani, R.~Montalto, and M.~Procesi.
{\it Reducibility of first order linear operators on tori via {M}oser's theorem}.
 {J. Funct. Anal.}, 276 (2019), no.3, 932--970.
\bibitem{Fraenkel} L. E. Fraenkel.
{  \it  An introduction to maximum principles and symmetry in elliptic problems}. 
Number 128. Cambridge University Press, 2000.
\bibitem{Gancedo} F. Gancedo.
{  \it Existence for the $\alpha$-patch model and the QG sharp front in Sobolev spaces}.
  Adv. Math., 217 (2008), no. 6, 2569--2598.
\bibitem{GP} F. Gancedo and N. Patel.
 {\it On the local existence and blow-up for generalized SQG patches}.
  Ann. PDE 7 (2021), no. 1, 63 pp.
   \bibitem{G-Kar}
C. Garc\'{\i}a.
{\it K\'{a}rm\'{a}n vortex street in incompressible fluid models}.
 { Nonlinearity}, 33(4):1625--1676, 2020.
\bibitem{Garcia}
C. Garc{\'\i}a.
{\it Vortex patches choreography for active scalar equations}.
 {J Nonlinear Sci 31, 75 (2021)}.
 \bibitem{GHS} C. Garc\'ia, T. Hmidi and  J. Soler, {\it Non uniform rotating vortices and periodic orbits for the two-dimensional Euler Equations.} Arch. for Ration. Mech. and Anal.  238, 929--1085 (2020).
\bibitem{GGS}
L. Godard-Cadillac, P. Gravejat, and D. Smets.
{\it Co-rotating vortices with n fold symmetry for the inviscid surface quasi-geostrophic equation}.
{ arXiv preprint arXiv:2010.08194}, 2020.
\bibitem{Gomez-Serrano}  J. G\'omez-Serrano.
 {\it On the existence of stationary patches}.
  Advances in Mathematics 343 (2019),  110--140.
\bibitem{gomez2019symmetry}  J. G\'omez-Serrano,  J. Park, J. Shi and  Y.  Yao.
 {\it Symmetry in stationary and uniformly-rotating solutions of active scalar equations.}
  Duke Math. J. 170 (2021), no. 13, 2957--3038.
 \bibitem{Grav-Smets} 
 P. Gravejat and  D.  Smets.
{\it Smooth travelling-wave solutions to the inviscid surface quasi-geostrophic equation}.
Int. Math. Res. Not. IMRN 2019, no. 6, 1744--1757
\bibitem{Grebert-Kappeler} B. Gr\'ebert and  T. Kappeler.
{\it The defocusing NLS equation and its normal form}.
 EMS Series of Lectures in Mathematics. European Mathematical Society (EMS), Z\"urich, 2014. x+166 pp. ISBN: 978-3-03719-131-6.
 \bibitem{Grebert-Paturel} B.  Gr\'ebert and  E. Paturel.
 {\it KAM for the Klein Gordon equation on $\mathbb{S}^d$}. 
 Boll. Unione Mat. Ital. 9 (2016), no. 2, 237--288. 
 \bibitem{HH}  Z. Hassainia and  T. Hmidi.
  {\it On the V-States for the generalized quasi-geostrophic equations}. 
  Comm. Math. Phys. 337
(2015), no. 1, 321--377.
\bibitem{HH2} Z. Hassainia and T. Hmidi.
 {\it Steady asymmetric vortex pairs for Euler equations}.
 { Discrete Contin. Dyn. Syst.}, 41(2021), no. 4, 1939--1969.
 \bibitem{Hass-Mass-Wheel} Z. Hassainia, N. Masmoudi and M. H. Wheeler. 
 {\it Global bifurcation of rotating vortex patches}.
 Comm. Pure Appl Math., 73 (2020), no. 9, 1933--1980.
 \bibitem{HW} Z. Hassainia and M. H. Wheeler.
 {\it Multipole vortex patch equilibria for active scalar equations}. 
 arXiv:2103.06839.
 \bibitem{H-K} S. He and  A. Kiselev.
  {\it Small scale creation for solutions of the SQG equation}.
   Duke Math. J. 170 (2021), no. 5, 1027--1041.
 \bibitem{Hm}  T. Hmidi.
   {\it On the trivial solutions for the rotating patch model}.
    J. Evol. Equ., 15 (2015 ), no. 4, 801--816.
    \bibitem{HM2} T. Hmidi and J. Mateu.
     {\it Bifurcation of rotating patches from Kirchhoff vortices}. 
     Discrete Contin. Dyn. Syst. 36 (2016), no. 10, 5401--5422. 
 \bibitem{HM3} T. Hmidi and  J. Mateu.
  {\it Degenerate bifurcation of the rotating patches}. 
  Adv. Math. 302 (2016), 799--850. 
  \bibitem{HM} T. Hmidi and  J.  Mateu.
   {\it  Existence of corotating and counter-rotating vortex pairs for active scalar equations}.
    Comm. Math. Phys. 350 (2017), no. 2, 699--747.   
 \bibitem{HMV} T. Hmidi, J.  Mateu and J.  Verdera.
  {\it Boundary Regularity of Rotating Vortex Patches}. 
  Arch. Ration. Mech. Anal. 209 (2013), no. 1, 171--208.
  \bibitem{HR}  T. Hmidi and  C. Renault.
  {\it Existence of small loops in a bifurcation diagram near degenerate eigenvalues}. 
  Nonlinearity, 30 (2017), no. 10, 3821--3852.
  \bibitem{Holder} E. H\"older.
   {\it\"Uber die unbeschr\"ankte Fortsetzbarkeit einer stetigen ebenen Bewegung in einer unbegrenzten inkompressiblen Fl\"ussigkeit}.
    Math. Z. 37 (1933), 727--738
\bibitem{Ioos} G. Iooss, P. I. Plotnikov and J. F. Toland. {\it Standing waves on an infinitely deep perfect fluid under gravity.} Arch. Ration. Mech. Anal., 177(3):367--478, 2005
\bibitem{Kirc} G. Kirchhoff,.
{\it Vorlesungen uber mathematische Physik}.
 (Leipzig, 1874).
\bibitem{KYZ} A. Kiselev, Y. Yao and A.  Zlato$\check{\textnormal{s}}$.
 {\it Local regularity for the modified SQG patch equation}.
  Comm. Pure Appl. Math. 70 (2017), no. 7, 1253--1315.
\bibitem{KRYZ} A. Kiselev, L. Ryzhik, Y. Yao and A.  Zlato$\check{\textnormal{s}}$.
{\it Finite time singularity for the modified SQG patch equation}. 
Ann. of Math. (2) 184 (2016), no. 3, 909--948. 
\bibitem{Kolmogorov} A.N. Kolmogorov.
 {\it On the persistence of conditionally periodic motions under a small change of the hamiltonian function}.
 Doklady Akad. Nauk SSSR 98 (1954), 527--530.
\bibitem{Kuksin1} S.  B. Kuksin.
{\it  Nearly integrable infinite-dimensional Hamiltonian systems}.
Lecture Notes in Mathematics, vol. 1556, Springer-Verlag, Berlin, 1993.
\bibitem{Kuksin-Poschel} S. Kuksin and J. P\"oschel.
 {\it  Invariant cantor manifolds of quasi-periodic oscillations for a nonlinear schr\"odinger equation}.
  Annals of Math 2 (1996), no. 143, 149--179.
\bibitem{Marchand} F. Marchand. 
{\it Existence and regularity of weak solutions to the quasi-geostrophic equations in the spaces $L^p$ or $\dot{H}^{-1/2}$}. Comm. Math. Phys., 277(1), (2008), 45--67.
\bibitem{moser} Moser.
{\it  On invariant curves of area-preserving mappings of an annulus}.
 Nachr. Akad. Wiss., G\"ottingen, Math. Phys. Kl. (1962), 1--20.
\bibitem{Nemes} G. Nemes.
 {\it Generalization of Binet’s Gamma function formulas}. 
 Integral Transforms Spec. Funct. 24 (8),  597--606.
\bibitem{Coralie} C. Renault.
 {\it Relative equilibria with holes for the surface quasi-geostrophic equations}.
  J. Differential Equations, 263 (2017), no. 1567--614.
\bibitem{Resnick} S. G. Resnick.
 {\it Dynamical problems in non-linear advective partial differential equations}. 
 PhD thesis, University of Chicago, Department of Mathematics, 1995.
\bibitem{R} J. L. Rodrigo,  {\it On the evolution of sharp fronts for the quasi-geostrophic equation}.
   Comm. Pure Appl. Math., 58, no.6  (2005), 821--866.
\bibitem{Russ} H.  R\"ussmann.
 {\it  Invariant tori in non-degenerate nearly integrable Hamiltonian systems}.
  Regul. Chaotic Dyn. 6 (2001), no. 2, 119--204.
\bibitem{Ruzhansky} M. V. Ruzhansky and  V.  Turunen.
 {\it Pseudo-Differential Operators and Symmetries. Background analysis and advanced topics. Pseudo-Differential Operators. Theory and Applications 2}.
  Birkh\"auser Verlag, Basel, 2010. 
\bibitem{sevryuk} M. B. Sevryuk.
 {\it Invariant tori in quasi-periodic non-autonomous dynamical
systems via Herman's method}.
 Discrete Contin. Dyn. Syst. 18 (2007), no. 2--3, 569--595.
\bibitem{T}
B. Turkington.
{\it Corotating steady vortex flows with {$N$}-fold symmetry}.
 { Nonlinear Anal.}, 9 (1985), no. 4, 351--369.
\bibitem{W} Watson.
 {\it A Treatise on the Theory of Bessel Functions}.
   Cambridge University Press, (1944).
   \bibitem{wayne} C. E. Wayne.
 {\it Periodic and quasi-periodic solutions of nonlinear wave equations
via KAM theory}. 
Communications in Mathematical Physics 127 (1990), no. 3, 479--528.
\bibitem{Wolibner} W. Wolibner.
{\it Un theor\`eme sur l’existence du mouvement plan d’un fluide parfait, homog\`ene, incompressible, pendant un temps infiniment long (French)}.
 Mat. Z., 37 (1933), 698--726.
\bibitem{Y1} Y. Yudovich.
 {\it Nonstationary flow of an ideal incompressible liquid}. 
 Zh. Vych. Mat.,  3, (1963), 1032--1066.
\end{thebibliography}
\end{document}